\documentclass[10pt,leqno]{amsart}
\usepackage{amssymb}
\usepackage{mathrsfs}
\usepackage{xcolor}
\usepackage{soul}
\usepackage{dsfont}
\usepackage{enumitem}
\usepackage{fixltx2e}

\usepackage[pdfborder={0 0 0}]{hyperref}

\usepackage{vmargin}

\newcommand\N{{\mathbb N}}
\newcommand\R{{\mathbb R}}
\newcommand\T{{\mathbb T}}
\newcommand\C{{\mathbb C}}

\newcommand\Ss{{\mathbb S}}
\newcommand\Z{{\mathbb Z}}

\def\AA{{\mathcal A}}
\def\BB{{\mathcal B}}
\def\CC{{\mathcal C}}
\def\DD{{\mathcal D}}
\def\EE{{\mathcal E}}
\def\FF{{\mathcal F}}
\def\GG{{\mathcal G}}
\def\HH{{\mathcal H}}
\def\II{{\mathcal I}}
\def\JJ{{\mathcal J}}
\def\KK{{\mathcal K}}
\def\LL{{\mathcal L}}
\def\MM{{\mathcal M}}

\def\OO{{\mathcal O}}
\def\PP{{\mathcal P}}
\def\QQ{{\mathcal Q}}
\def\RR{{\mathcal R}}
\def\SS{{\mathcal S}}
\def\TT{{\mathcal T}}
\def\UU{{\mathcal U}}
\def\VV{{\mathcal V}}
\def\WW{{\mathcal W}}
\def\VV{{\mathcal V}}
\def\XX{{\mathcal X}}
\def\YY{{\mathcal Y}}

\def\BBB{{\mathscr B}}
\def\CCC{{\mathscr C}}
\def\DDD{{\mathscr D}}
\def\EEE{{\mathscr E}}

\def\GGG{{\mathscr G}}

\def\KKK{{\mathscr K}}

\def\MMM{{\mathscr M}}

\def\RRR{{\mathscr R}}

\def\TTT{{\mathscr T}}

\def\kappaa{a}
\def\frakg{{\mathfrak{g}}}

\def\VVb{{{\mathcal V}_{\bf b}}}

\def\eps{{\varepsilon}}

\newcommand{\gpf}{\gamma_{_{\!+}\!} f}
\newcommand{\gp}{\gamma_{_{\!+}\!} }

\newcommand{\gpm}{\gamma_{_{\!\pm}\!} }
\newcommand{\gmp}{\gamma_{_{\!\mp}\!} }

\def\Wloc{W_{\rm loc}} 
\def\Mloc{M_{\rm  loc}} 
 
\def\Lloc{L_{\rm loc}} 
\def\Lip{\rm Lip} 
 
\def\Hloc{H_{\rm  loc}} 

\newcommand{\Sp}{{\mathbb S}}

\def\esssup{\mathop{\rm ess\,sup\,}}
\def\essinf{\mathop{\rm ess\,inf\,}}
\newcommand{\wto}{\rightharpoonup}
\newcommand{\rto}{\mathop{\rightharpoonup}^r}

\def\1{\mathbf 1}

\DeclareMathOperator{\sign}{sign}

\DeclareMathOperator{\supp}{supp}
\DeclareMathOperator{\Span}{Span}
\DeclareMathOperator{\tb}{{t_{\bf b}}} 

\newtheorem{theo}{Theorem}
\newtheorem{prop}[theo]{Proposition}
\newtheorem{lem}[theo]{Lemma}
\newtheorem{cor}[theo]{Corollary}
\newtheorem{rem}[theo]{Remark}

\newtheorem{defin}[theo]{Definition}
\newtheorem{exple}[theo]{Example}
\newtheorem{exples}[theo]{Examples}

\newcommand{\beqn}{\begin{equation}}
\newcommand{\eeqn}{\end{equation}}
\newcommand{\bear}{\begin{eqnarray}}
\newcommand{\eear}{\end{eqnarray}}
\newcommand{\bean}{\begin{eqnarray*}}
\newcommand{\eean}{\end{eqnarray*}}

\newcommand{\ds}{\displaystyle}

\newcommand{\intot}{\ds\int_0^t}

\def\Nt{|\hskip-0.04cm|\hskip-0.04cm|}

\newcommand{\Black}{\color{black}}
\newcommand{\Red}{\color{red}}
\newcommand{\Blue}{\color{black}}
\newcommand{\Cyan}{\color{black}}

\def\eps{{\varepsilon}}
\def\Div{{\rm div}}
\def\Trace{{\rm Tr}}

\def\bb{\mathbf{b}}
\def\ii{\mathbf{i}}

\def\NormT{|\hskip-0.04cm|\hskip-0.04cm|}

\begin{document}

\makeatletter
\@addtoreset{equation}{section}
\@addtoreset{theo}{section}
\makeatother

\author[C. Fonte Sanchez]{Claudia Fonte Sanchez}

\address[C. Fonte Sanchez]{Universit\'e Grenoble Alpes, Inria, 38000 Grenoble, France}

\email{claudia.fonte-sanchez@inria.fr}

\author[P. Gabriel]{Pierre Gabriel}
\address[P. Gabriel]{Institut Denis Poisson, Université de Tours, Université d'Orléans, CNRS,
Parc Grandmont, 37000 Tours, France}
\email{pierre.gabriel@univ-tours.fr}

 \author[S.~Mischler]{St\'ephane Mischler}

\address[S.~Mischler]{IUF \& CEREMADE (CNRS UMR n$^\circ$ 7534), PSL
  university, Universit\'e Paris-Dauphine, Place de Lattre de
  Tassigny, 75775 Paris 16, France}

\email{mischler@ceremade.dauphine.fr}


\bigskip
\bigskip
\bigskip

\parindent 0 cm
\title{
On the Krein-Rutman theorem and beyond
}
\maketitle

\bigskip
\bigskip
 
\tableofcontents

%
%
%
%
%
%
%
%
%
%
%


\bigskip 
\section{Introduction}
\label{sec:KRpossible}


\subsection{Framework and the main result}\label{ssec:main-result}

In this work, we revisit the Krein-Rutman theory for semigroups of positive operators in a Banach lattice framework and we provide some very general, efficient and handy results
with possible constructive estimates about 

- the existence of a solution to the first eigentriplet problem; 

- the geometry of the principal eigenvalue problem; 

- the asymptotic stability of the first eigenvector.  

\smallskip
This abstract theory is motivated and illustrated by several examples of (partial) differential, integro-differential and integral operators. In particular, we revisit the first eigenvalue problem and the asymptotic stability of the first eigenvector for 

- some parabolic equations in a bounded domain and in the whole space;  

- some transport equations in a bounded or unbounded domain, including some growth-fragmentation models and some kinetic models; 

- the kinetic Fokker-Planck equation in a bounded domain; 

- some mutation-selection models. 

\smallskip 
The results we establish on these examples are more general and more accurate that what we can find in the literature. 
Our approach is in the same time able to tackle some critical cases, but also it is very natural and makes possible to bring out the main important properties for each example
and   to get rid of many technical issues.

\smallskip
The present work is motivated by new problems and ideas  presented in the lectures on the Krein-Rutman theorem by P.-L. Lions at Collège de France \cite{PL2} and by the recent contributions by Bansaye et al \cite{Bansaye2022}  and  by Cañizo and Mischler \cite{MR4534707} developing Doblin-Harris techniques.  
Bringing and developing these ideas and techniques together with the more classical  {\it spectral analysis} approach developed or synthesized in  previous contributions by Krein and Rutman \cite{MR0027128}, by Arendt et al \cite{MR839450}, by Mischler and Scher \cite{MR3489637}, by B\'{a}tkai et al \cite{MR3616245} and many others, we are then able to significantly 
 generalize and improve the Krein-Rutman theory for positive semigroups. 


\smallskip

\Blue
The abstract results are developed in a Banach lattice framework, that is in a Banach space $(X,\| \cdot \|)$ endowed with a compatible order relation $\ge$, and thus with an associated positive cone $X_+ := \{ f \in X; \, f \ge 0 \}$. We further assume  either $X=Y'$ or $X' = Y$ for another dual Banach lattice $Y$ as well as, when needed,  some additional structure properties. 
The (standard) Banach framework will be presented in Section~\ref{subsec:BanachLattice} and the  additional structure properties will be introduced in Sections~\ref{subsec:geo1} and \ref{subsec:SP-defStrongP}. 
We are  not seeking generality for its own sake, but rather in order to have a unified framework for all the applications we have in mind, in order to fit with the standard framework and thus to be  able to compare our results with the abundant literature, and finally in order to understand what is really needed at each step for proving our several results. 
We emphasize that 
all the above requirements  hold in the case of the usual Banach lattices used in PDE and stochastic processes theory, in particular for the following examples of usual Banach lattices with which we will work in the applications sections:
 
 \Black


\smallskip
$\bullet$ $X := C_0(E)$, the space of continuous functions which tend to $0$ at infinity (when $E$ is not a compact set) endowed with the uniform norm,  or $X := C_{0,m}(E)$ its weighted variant; 

\smallskip
$\bullet$ $X := L^p(E) = L^p(E,\EE,\mu)$, the Lebesgue space of functions associated to the Borel $\sigma$-algebra $\EE$, 
a positive $\sigma$-finite measure $\mu$ and an exponent $p \in [1,\infty]$,  or $X := L^p_m(E)$ its weighted variant;  

\smallskip
$\bullet$ $X := M^1(E) = (C_0(E))'$, the space of Radon measures defined as the dual space of $C_0(E)$,   or $X := M^1_m(E)$ its weighted variant. 

 \smallskip 
 In all the above examples, $E$ denotes a $\sigma$-compact metric space, and we write $E = \cup E_R$, with $E_R \subset E_{R+1}$, $E_R$ compact.

\smallskip\smallskip
We next consider a positive one-parameter semigroup of operators $S = S_\LL$ on $X$ (we will indifferently writes $S_t = S(t) = S_\LL(t)$ for $t\ge0$), and we denote by $\LL$ its generator, by $D(\LL) \subset X$ the domain of $\LL$, by $\rho(\LL) \subset \C$ the resolvent set of $\LL$  and by $\Sigma(\LL) = \C \backslash \rho(\LL)$ the spectrum of $\LL$. 
We also denote by $S^*$ and $\LL^*$ the corresponding semigroup and generator on the dual space $Y$, and we refer to Section~\ref{subsec:BanachLattice} for more notations. 

\smallskip 
As announced, we may split the issue into several pieces concerning the stationary and the evolution associated problems. 

\smallskip\smallskip
$\bullet$ {\bf Existence.} 
We are first interested in the existence part of the first or principal eigentriplet problem, namely we wish to bring out very general conditions under which 

\begin{enumerate}[label={\bf(S1)},itemindent=12mm,leftmargin=0mm]
\item\label{S1} there exists a solution $(\lambda_1, f_1, \phi_1) \in \R \times X \times Y$ to the eigentriplet problem 
\bear
\label{eq:triplet1}&&\LL f_1 = \lambda_1 f_1, \quad f_1 \ge 0, \quad   f_1 \not= 0, 
\\
\label{eq:triplet2}&&\LL^* \phi_1 = \lambda_1 \phi_1, \quad \phi_1 \ge 0, \quad \phi_1 \not = 0, 
\eear
and furthermore $\lambda_1$ coincides with the spectral bound, namely
\beqn
\label{eq:triplet3} \lambda_1 =s(\LL) := \sup  \{\Re e(\lambda); \lambda \in \Sigma(\LL) \} = \inf \{ \kappa \in \R; \, \Delta_\kappa \subset \rho(\LL)\},
\eeqn
where $\Delta_\alpha$ is the open half plan $\Delta_\alpha := \{ z \in \C; \ \Re e (z) > \alpha \}$.
\end{enumerate}

{\Cyan Defining $\Sigma_P(\LL)$ as the point spectrum (or set of eigenvalues)},  we emphasize on the fact that this problem is named as the principal eigenvalue problem because 
$$
\lambda_1 \in \Sigma_P(\LL)\subset\{z\in\C,\ \Re e(z)\leq\lambda_1\}.
$$

\smallskip\smallskip
$\bullet$ {\bf Geometry.}  A second issue is about an accurate analysis of the principal eigentriplet solution and of the geometry of the (principal part of the) spectrum. 

\smallskip
On the one hand, concerning the eigentriplet solution, we investigate conditions such that  
 
\begin{enumerate}[label={\bf(S2)},itemindent=12mm,leftmargin=0mm]
\item\label{S2} $f_1$ is strictly positive (we refer to Section~\ref{subsec:MorePositive} for a definition) and $f_1$ is the unique (up to normalization) positive eigenvector  
for $\LL$, $\phi_1$ is strictly positive and $\phi_1$ is the unique (up to normalization) positive eigenvector for $\LL^*$, and finally $\lambda_1$ is geometrically  and algebraically simple for both $\LL$ and $\LL^*$. We then may make the (usual) normalization choice
\beqn\label{intro-normalizationf1phi1}
\bigl( \|f_1\|=1, \ \langle f_1,\phi_1 \rangle = 1 \bigr) \ \hbox{ or } \ 
\bigl( \|\phi_1\|=1, \ \langle f_1,\phi_1 \rangle = 1 \bigr). 
\eeqn
\end{enumerate}

\smallskip
We are next interested by describing the boundary point spectrum
$$
\Sigma_P^+ (\LL):= \Sigma_P (\LL) \cap \Sigma_+ (\LL), 
$$
where we define the boundary spectrum  
$\Sigma_+(\LL) := \Sigma(\LL) \cap \{ z \in \C; \, \Re e z = s(\LL) \}$.
More precisely, we exhibit some conditions   such that 

\begin{enumerate}[label={\bf(S3$_\arabic*$)},itemindent=13mm,leftmargin=0mm,itemsep=1mm]
\item\label{S31}  $\Sigma^+_P(\LL) - \lambda_1$ is a (discrete) additive subgroup of $i\R$; 
\item\label{S32} $\Sigma_P^+ (\LL)$ is trivial, namely  
\beqn\label{eq:intro-Sigma+=lambda1}
  \Sigma_{P}^{+} (\LL)  = \{ \lambda_1 \}; 
\eeqn
\noindent  or  
\item\label{S33}  $\Sigma_P^+ (\LL)$ is trivial and $\Sigma(\LL)$ enjoys  a spectral gap property (on its principal part), namely 
\beqn\label{eq:intro-Sigma=gap}
\exists \, \kappa < \lambda_1; \quad \Sigma (\LL) \cap \Delta_\kappa =  \{ \lambda_1 \}. 
\eeqn
\end{enumerate}

In the last situation \eqref{eq:intro-Sigma=gap}, a band separates the spectral value $\lambda_1$ to the remainder of the spectrum, while there is no spectral gap when \eqref{eq:intro-Sigma+=lambda1} holds but \eqref{eq:intro-Sigma=gap} does not.

%
%
%

\smallskip

The importance of such a eigentriplet comes from the fact that we may associate the Malthusian function 
$$
F_1(t) := e^{\lambda_1 t} f_1, 
$$
which is a particular solution to the evolution equation (with maximal growth) and a natural candidate to capture the main asymptotic feature of generic trajectories.

\smallskip\smallskip
$\bullet$ {\bf Asymptotic stability.}  In order to formulate our third main issue, namely the asymptotic stability of  $F_1$, 
we introduce the rescaled operators
$\widetilde \LL = \LL - \lambda_1$ and $\widetilde \LL^* = \LL^* - \lambda_1$, so that 
\bean
 \widetilde \LL f_1 = 0, \quad \widetilde \LL^* \phi_1 = 0,
\eean
or in other words, $f_1$ is a stationary state of the semigroup $\widetilde S = S_{\widetilde \LL}$ and $\phi_1$ is a stationary state of the semigroup $\widetilde S^* = S_{\widetilde \LL^*}$, and thus a conservation law for $\widetilde S$:
\bean
\widetilde S(t) f_1 = f_1, \quad \widetilde S^*(t) \phi_1 = \phi_1, \quad \langle S(t) f,\phi_1 \rangle = \langle  f,\phi_1 \rangle, 
\eean
for any $t \ge 0$ and any $f \in X$. Because of the property of the eigentriplet and of the normalization assumption \eqref{intro-normalizationf1phi1}, we may reduce the issue to considering the case $f \in X$ satisfies $ \langle \phi_1,f \rangle= 0$ when \ref{S32} or  \ref{S33} holds and more generally $f \in Y_0^\perp$ when \ref{S31} holds, 
where $Y_0 \subset Y$ stands for the eigenspace associated to the eigenvalues belonging to $\Sigma_P^+(\LL^*)$. 
Depending of the hypotheses we made on $\LL$ and $S_\LL$, we are able to establish for such a   $f$ some

\begin{enumerate}[label={\bf(E\arabic*)},itemindent=13mm,leftmargin=0mm,itemsep=1mm]
\item\label{E1intro}{\bf mean ergodic property}, namely
$$
\frac{1}{T} \int_0^T \widetilde S_t f dt \to 0 \quad\hbox{as}\quad T \to \infty; 
$$

\item\label{E2intro}{\bf ergodic property}, namely
$$
 \widetilde S_t f \to 0  \quad\hbox{as}\quad t\to \infty; 
$$

\item\label{E3intro}{\bf quantitative asymptotic stability}, which may be  geometric (or exponential) in the spectral gap \eqref{eq:intro-Sigma=gap} case, namely
\bigskip
\begin{enumerate}[label={\bf(E3$_1$)},itemindent=22mm,leftmargin=0mm]
\item\label{E31intro} \qquad$\| \widetilde S(t) f \| \le  C \, e^{-\eps t} \| f\|, \quad \forall \, t \ge 0, $
\end{enumerate}
\bigskip
for possible constructive constants $\eps > 0$ and $C \ge 1$, or under the weaker condition \eqref{eq:intro-Sigma+=lambda1} only subgeometric, namely 
\bigskip
\begin{enumerate}[label={\bf(E3$_2$)},itemindent=22mm,leftmargin=0mm]
\item\label{E32intro} \qquad$\| \widetilde S(t) f \|_1 \le \Theta(t) \| f  \|_2, \quad \forall \, t \ge 0, $
\end{enumerate}
\bigskip
where $\| \cdot \|_2 = \| \cdot \|_X$,  
$\| \cdot \|_1$ is a weaker norm and  $\Theta : \R_+ \to \R_+$ is a constructive decay function satisfying $\Theta(t) \searrow 0$ when $t \nearrow \infty$. 

\end{enumerate}

%
%

\smallskip\smallskip\indent
 We aim now to allude some general hypotheses on the semigroup $S_\LL$ or its generator $\LL$ such that the above three main issues may be tackled. 
 Additionally to the yet mentioned fact that $S_\LL$ is positive 
 (which is almost equivalent to the fact that its resolvent is a positive operator,  that $\LL$ enjoys a weak maximum principle or that $\LL$ enjoys Kato's inequality)
 our hypotheses are mainly of two kinds : 
 
\quad  - strict positivity conditions; 
 
\quad  - regularity conditions  on the dominant part of the semigroup; 
 
 and these ones may be formulated at the stationary level directly on the generator $\LL$ or its resolvent $\RR_\LL$ or they may be formulated at the evolution  level on the semigroup of operators $S_\LL$. 
 Of course, in order to establish constructive results these hypotheses will have to be formulated in a quantitative way. 
 
 \smallskip
The strict positivity we will introduce and use are of different kinds:  
 
\quad  - strong maximum principle on the generator, or equivalently  irreducibility of the semigroup; 
 
\quad  - reverse Kato's inequality for the generator or aperiodicity condition of the semigroup; 
 
\quad  - Doblin-Harris condition on the semigroup, which may be formulated as 
\beqn\label{eq:intro-DHcondition}
S_T f \ge g_0 \langle \psi_0, f \rangle, \quad \forall \, f \in X_+, 
\eeqn
for some for some $T > 0$ and convenient $g_0 \in X_+ \backslash \{0 \}$, $\psi_0 \in Y_+ \backslash \{0 \}$. 
 
 \smallskip
Less systematically but in a crucial way, we will make use of somehow related 

\quad  - barrier functions and positive subeigenfunctions, which for the last one typically writes 
 \beqn\label{eq:intro-subeigenfunction}
 \exists \, \kappa_0 \in \R, \ \exists \, \phi_0 \in Y_+ \backslash \{0 \}, \quad
 \LL^* \phi_0 \ge \kappa_0 \phi_0.
 \eeqn

 \smallskip
 On the other hand, some regularity is needed on the dominant part of the semigroup. In order to briefly explain the issue, we assume that $\LL = \AA + \BB$ with $\AA \in \BBB(X)$ and $\BB$ is 
the generator of a   semigroup $S_\BB$. In such a context, we may write the resolvent  factorization identity  
$$
\RR_\LL  = \RR_\BB  + \RR_\BB  \AA \RR_\LL 
$$
on the resolvent $\RR_\LL$ of $\LL$ and $\RR_\BB$ of $\BB$, and its  iterated version 
\beqn\label{eq:intro:RLRB-RLARBN}
\RR_\LL  = \VV + \WW \RR_\LL, \quad \VV := \RR_\BB + \dots + \RR_\BB (\AA  \RR_\BB)^{N-1} , \quad \WW :=   (\RR_\BB \AA)^N. 
\eeqn
At the level of the generator, our regularity assumption then  typically  writes
\beqn\label{eq:intro-regRL}
\sup_{z \in \Delta_\kappa} \| \VV(z) \|_{\BBB(X)} < \infty,  \quad 
\sup_{z \in \Delta_\kappa} \| \WW(z) \|_{\Blue\BBB(\XX_0,\XX_1)} < \infty, 
\eeqn
for some $\kappa \in \R$ and $\XX_1 \subset X \subset \XX_0$, which is nothing but the   classical Voigt's power compact condition when $\XX_0 = X$ and $\XX_1 \subset X$ with compact embedding. 
Similarly, at the level of the semigroup, we may write the associated Duhamel formula 
$$
S_\LL =S_\BB + (S_\BB \AA) * S_\LL,
$$
(we refer to Section~\ref{subsect-Exist2-Dissip} for a precise definition) and its  iterated version 
\beqn\label{eq:intro:stopDyson}
S_\LL = V + W * S_\LL, 
\quad V := 
  \sum_{\ell=0}^{N-1} S_\BB *  (\AA S_\BB)^{(*\ell)}, \quad W := (S_\BB\AA )^{(*N)},
\eeqn
with $N \ge 1$.
At the level of the semigroup, our regularity assumption then  typically  writes
\beqn\label{eq:intro-regSL}
\sup_{t \ge 0} \| V(t) e^{-\kappa t} \|_{\BBB(X)} < \infty, \quad  
\sup_{t \ge 0} \| W(t) e^{-\kappa t} \|_{\Blue\BBB(\XX_0,\XX_1)} < \infty, 
\eeqn
for some $\kappa \in \R$ and $\XX_1 \subset X \subset \XX_0$ in the dissipative framework and a variant of these estimates in a weak dissipative framework. The crucial information is $\kappa < \kappa_0$ (dissipative framework) or $\kappa = \kappa_0$ (more involved weak dissipative framework).

%
%
%
%
%
%
%
%
%
%
%
%
%
%
%
%
%
%

We are now in position to state in a very informal way our main result in an abstract Banach lattice framework.

\begin{theo}[rough version]\label{theo:main-intro} 
Let us consider a  Banach lattice $X$ picked up in the examples listed above and a positive semigroup $S_\LL$ on $X$ which enjoys the above splitting structure \eqref{eq:intro:RLRB-RLARBN}, \eqref{eq:intro-regRL}, \eqref{eq:intro:stopDyson}, \eqref{eq:intro-regSL}.  

\smallskip
(1) Conclusion \ref{S1} holds under the localization of the principal spectrum assumption $\kappa < \kappa_0$ and a weak compactness assumption on the regular part $\WW$ or $W$ in the splitting.

\smallskip
{\Blue (2)  Under an additional strong positivity assumption (irreducibility property) the conclusions \ref{S2}, \ref{S31} and  \ref{E1intro} hold. }


\smallskip
In order to make one step further, we have the three next possibilities 

\smallskip
(3)  {\Blue Under an even stronger strict positivity condition (reverse Kato's inequality or aperiodicity property),} the  conclusion \ref{S32} holds, as well as \ref{E2intro} when $X \subset \Lloc^1$.

\smallskip
(4) Alternatively, under an additional strong compactness assumption on the regular part $W$  of the splitting, the quantitative exponential asymptotic stability \ref{E31intro} holds (without constructive constants), 
and thus also the spectral gap conclusion \ref{S33} holds  (in a not constructive way).

\smallskip
(5) Alternatively, under an additional Doblin-Harris{\Cyan-type condition  similar to}  
 \eqref{eq:intro-DHcondition} and an appropriate regularity estimate on the regular part  $W$ of the splitting,  the quantitative  asymptotic stability \ref{E3intro} holds for both the geometric and subgeometric framework with now constructive constants.

 \end{theo}

\smallskip
{\Blue
In Section~\ref{sec:ExistenceKR}, we establish  the conclusion \ref{S1} about the existence of an eigentriplet  (and in particular we establish Theorem~\ref{theo:main-intro}-(1)) 
under a family of fundamental conditions denoted by \ref{H1}, \ref{H2}, \ref{H3}, 
and we discuss several handy conditions involving these ones, in particular some based on a splitting condition  \ref{HS1} related to \eqref{eq:intro:RLRB-RLARBN}-\eqref{eq:intro-regRL}. 
An alternative approach for proving at least part of the conclusion \ref{S1} is presented in Section~\ref{sec:DynamicalExistenceKR}, it is based on a dynamical argument and some splitting conditions   \ref{HS2} and \ref{HS3}
related to \eqref{eq:intro:stopDyson}-\eqref{eq:intro-regSL}. 
In Section~\ref{sec:Irreducibility}, we establish  conclusions \ref{S2}  and   \ref{E1}  under an additional structure condition on $X$ denoted by \ref{X1} and an additional positivity condition denoted by \ref{H4} (which is nothing but the  irreducibility for a semigroup). 
In Section~\ref{sec:geo2}, we first introduce some additional structure conditions \ref{X2} and \ref{X3} on $X$ and we establish \ref{S31}, thus concluding the proof of Theorem~\ref{theo:main-intro}-(2). 
We next establish \ref{S32} and \ref{E2} (and thus Theorem~\ref{theo:main-intro}-(3)) under two additional and somehow stronger positivity assumptions  \ref{H5} (reverse Kato's condition) and \ref{H5'} (aperiodicity condition). Introducing an additional alternative and more classical strong compactness assumption on the regular part $W$  of the splitting, we prove \ref{E31intro}  and \ref{S33} without constructive constants (and thus Theorem~\ref{theo:main-intro}-(4)). In Section~\ref{sec:QuantitativeStabilityKR} finally,  we introduce some relevant and somewhat stronger conditions and establish  Theorem~\ref{theo:main-intro}-(5). 
 We emphasize that all the concrete Banach lattices previously introduced satisfy the structure assumptions \ref{X1}, \ref{X2} and \ref{X3}.}
It is also worth mentioning that the assumptions in (4) and (5) may be optimal in the sense that reciprocal implications are likely to be true. We do not follow that line of investigation but rather refer to \cite{MR3489637,Bansaye2022} where such kind of results are established.


 \medskip
\subsection{Discussion about Theorem~\ref{theo:main-intro}}
We discuss several works related to the main Theorem~\ref{theo:main-intro} as well as the hypotheses and the techniques used during the proof.

%
%
%
%
%


\subsubsection{\bf The Krein-Rutman work and related approaches} 


For a strictly positive matrix in a finite dimensional space, Perron \cite{MR1511438}  and Frobenius  \cite{Frobenius} establish at the beginning of the 20th Century that the eigenvalue with largest real part is unique, real and simple. 
%
%
%
%
%
%
%
%
In their pioneer work, Krein and Rutman establish in \cite{MR0027128} for the very first time
{\Blue a Perron-Frobenius type theorem in an  infinite dimensional Banach space.}

\begin{theo}[Krein-Rutman] \label{theo:intro-KRresult} 
Consider a Banach lattice with positive cone $X_+$ and strictly positive cone $X_{++} := \hbox{\rm int}X_+ \not=\emptyset$. Consider a linear and compact operator $\RR : X \to X$ such that $\RR : X_+ \to X_+$ and $\RR : X_+\backslash\{0\} \to X_{++}$. Then there exists a unique eigentriplet $(\mu_1,f_1,\phi_1)$ such that 
$\mu_1 > 0$, $f_1 \in X_{++}$, $f_1 = \mu_1 \RR f_1$, $\phi_1 \in X'_{++}$, $\phi_1 = \mu_1 \RR^* \phi_1$. 
\end{theo}

The non-emptiness of $X_{++}$ and the strict positivity assumption $\RR : X_+\backslash\{0\} \to X_{++}$ can be relaxed, to the price of loosing the uniqueness and strict positivity properties of the eigenvectors.
For a bounded operator $\RR$ on $X$, we denote by $r(\RR)$ the spectral radius
\[r(\RR) := \sup  \{|\lambda|;\, \lambda \in \Sigma(\RR) \} \leq \|\RR\|.\]

\begin{theo}[Krein-Rutman] \label{theo:intro-KRresult2} 
Consider a Banach lattice with positive cone $X_+$ and a linear and compact operator $\RR : X \to X$ such that $\RR : X_+ \to X_+$ and $r(\RR)>0$.
Then there exists an eigentriplet $(\mu_1,f_1,\phi_1)$ with
$\mu_1 = r(\RR)$, $f_1 \in X_{+}\setminus\{0\}$, $f_1 = \mu_1 \RR f_1$, $\phi_1 \in X'_{+}\setminus\{0\}$, $\phi_1 = \mu_1 \RR^* \phi_1$. 
\end{theo}

In Theorems~\ref{theo:intro-KRresult} and~\ref{theo:intro-KRresult2},  the operator $\RR$ corresponds to a resolvent operator $\RR := (\kappa-\LL)^{-1}$ for $\kappa > 0$ large enough, so that when it applies, we deduce in particular that the first eigenvalue problem \eqref{eq:triplet1}-\eqref{eq:triplet2} has a solution with $\lambda_1 = \kappa-\mu_1$. The two conditions $\hbox{\rm int}X_+ \not=\emptyset$ and $\RR : X_+\backslash\{0\} \to X_{++}$ are very strong. The first one essentially imposes to work in the space of continuous functions and the second one to work in a bounded domain. 
The result is however suitable and directly applicable (and somehow restricted) to an elliptic operator with smooth coefficients set in a bounded domain  with suitable boundary conditions or to a Fredholm integral operator with positive kernel also set in a bounded domain. 
In the elliptic context, the property $\RR : X_+\backslash\{0\} \to X_{++}$ is nothing but the strong maximum principle while the compactness property of $\RR$ comes from the elliptic regularity. 
We refer to Section~\ref{subsect:Exist1-discussion} for further discussions. 
The weaker condition $r(\RR)>0$ is less restrictive and is in particular always satisfied for irreducible operators, by virtue of de Pagter's theorem~\cite{Pagter1986}.
In the same framework, Theorems~\ref{theo:intro-KRresult} and~\ref{theo:intro-KRresult2} have been next slightly extended by Bonsall~\cite{MR92938}, Schaefer~\cite{MR106402}, Karlin \cite{MR0114138} or Nussbaum \cite{MR643014}
for instance. We also refer to the book by Dautray and Lions \cite{MR1064315} for a clear and comprehensible presentation and several possible versions.

In his paper~\cite{Birkhoff1957}, G.~Birkhoff derived the Perron-Frobenius theorem by proving a contraction principle in Hilbert's projective metric for positive matrices. His result actually applies to any ``uniformly positive bounded'' linear operators of a Banach lattice, such as integral operators with positive kernels, and also provides geometric stability estimates. A closely related result was proved by E. Hopf~\cite{Hopf1963}, and this Birkhoff-Hopf contraction theorem was subsequently generalized and sharpened, and its proof simplified, by several authors, see in particular~\cite{Bauer1965,Brooks2009,Bushell1973,Eveson1995b,Eveson1995a,Kohlberg1982,Nussbaum1988,Ostrowski1964}. This approach of the Krein-Rutman theorem requires some ``uniform positivity and boundedness'' of the operator, which is quite restrictive, but it nevertheless allows to recover, through an approximation procedure, the standard result of Theorem~\ref{theo:intro-KRresult}, see~\cite[Thm.~6.18]{Brooks2009}. The contraction in Hilbert's projective metric has the advantage to be applicable in partially order linear vector spaces without any topological structure~\cite{Eveson1995a}, and to nonlinear maps~\cite{Nussbaum1988}.

\subsubsection{\bf Spectral analysis approach}  
In his paper \cite{MR146675}, R.S.~Phillips formalized the notion of {\it positive semigroup} acting on a Banach lattice paving the way to a new field of research. 
In the precursory work \cite{MR230531} by  Vidav and next in a series of papers by Greiner and co-authors \cite{MR763356,MR617977,MR839450}, Webb \cite{MR772205,MR902796} and B\"urger~\cite{MR923493} (see also \cite[C-III, Cor.~2.12, Thm.~3.12]{MR839450}, \cite[Thm. VI.1.12, Cor. VI.1.13]{MR1721989} or more recently Theorem~14.17 in the very pedagogical book \cite{MR3616245}) significant generalizations of  the Krein-Rutman theory were established leading to, roughly speaking, the following result.

\begin{theo}\label{theo:Greiner&co}
Consider a  positive  semigroup $S_\LL$ on a  (suitable) Banach lattice $X$ which is irreducible and such that $s(\LL) > - \infty$ is a pole, then 

$\bullet$ $s(\LL)$ is a first-order pole with one-dimensional and strictly positive residue,   so that in particular there exists a solution $(\lambda_1,f_1,\phi_1)$ to the eigentriplet problem;  

$\bullet$ There exists $\alpha \in \R$ such that $\Sigma_+(\LL) =   \{ s(\LL) + i  \alpha \Z  \}$ 
 consists of first-order poles with one-dimensional residue.

$\bullet$ A practical way for verifying that $s(\LL) > - \infty$ is a pole consists in assuming that $\LL$ enjoys the splitting structure $\LL = \AA + \BB$, as described above, with 
  $s(\BB) < s(\LL)$ and $\AA$ is $\BB$ power compact,   that is to say $\WW$ is compact, on $\Delta_{s(\BB)}$. 

Assuming furthermore that  $\widetilde S_\LL$ is quasi-compact then 

$\bullet$ $\widetilde S_\LL$ is exponential asymptotically stable in $\Span \{ \phi_1 \}^\perp$  
(without constructive constants). 

%
%

\end{theo}

The most important improvements in the above result are the fact that the condition $\hbox{\rm int}X_+ \not=\emptyset$ and the strong compactness of the resolvent operator $\RR_\LL$ are removed, 
and also that the exponential asymptotically stability is established.
The hypotheses seem stronger to those stated in Theorem~\ref{theo:main-intro}-(1), where only weak compactness is required what is not the case here. It is however worth emphasizing that in an AL-space 
and an AM-space (what includes the examples $C_0(E)$ and $L^1(E)$)  a power weak compactness implies a
power strong compactness (see \cite[Rk.~2.1]{MR923493} and \cite[Cor.~1 of Thm.~II.9.9]{MR0423039}).  
The hypotheses and conclusions are  similar to those stated in Theorem~\ref{theo:main-intro}-(4). 
The proof is based on the one hand on the Banach lattices theory as formalized for instance by Schaefer~\cite{MR0423039} (see also \cite{MR683043,MR763347,MR790308,MR839450} for significant developments) using notions as ideals and quasi-interior points. On the other hand, it takes advantage on the perturbation techniques initiated by Phillips in \cite{MR54167}  and developed further by 
J\"{o}rgens \cite{MR103419}, Vidav \cite{MR236741,MR259662} and Voigt \cite{MR595321} leading to the notions of power compact resolvent and quasi-compact semigroup, essential spectrum and Calkin algebra.

The above theorem in particular applies to  a positive and  irreducible semigroup which  is eventually norm continuous and its generator has  compact resolvent 
(see for instance Corollary~{VI.1.13} in \cite{MR1721989}).
In that case indeed, one can show that $s(\LL) > - \infty$, $\Sigma_+(\LL)$ is bounded and consists of poles, so that $\Sigma_+(\LL) =   \{ s(\LL) \}$ and the essential growth bound $\omega_{\rm ess}(S)$ associated to the essential spectrum (see for instance \cite[Sec~14.1]{MR3616245} for a definition) satisfies $\omega_{\rm ess}(S) < \omega(S) = s(\LL)$. 
The theorem was motivated and successfully applied to Boltzmann like transport operator \cite{MR230531},   cell division operator \cite{Diekmann1984},  age structured equation \cite{MR902796} and selection-mutation dynamics~\cite{MR923493}. We also refer to \cite[Ch.~VI]{MR1721989} and \cite{MR3616245} for other numerous applications. Although very general and quite efficient, we formulate several criticisms about the above result.

 - The exponential convergence result is definitively not constructive and that approach is not able to say anything about the weak dissipative case (a framework we will introduce latter, see in particular Section~\ref{subsect:AboutWeakDissip}).
 
 - We may observe that Theorem~\ref{theo:Greiner&co} is not so popular in the probability and the PDE communities and still many works in these domains refer to the original 
 Krein and Rutman theorem even when some additional (approximation) arguments are needed rather than applying directly  Theorem~\ref{theo:Greiner&co}. 
 By the way, we did not find in the literature where Theorem~\ref{theo:Greiner&co}  is stated in such an handy way (the closer formulation is probably \cite[Thm.~VI.1.12]{MR1721989} which is given without proof). 

- The proof of Theorem~\ref{theo:Greiner&co} that we may find in the above quoted references is written in a very specific 
and abstract  language which make it quite obscure.




\medskip
In \cite{MR3489637,MischErratum}, one of the authors proposes the following variant.

\begin{theo}\label{theo:MS-KR} Consider a positive semigroup $S_\LL$ which satisfies \eqref{eq:intro-subeigenfunction} with $\kappa_0 \in \R$,  it is irreducible and its generator enjoys the   splitting structure \eqref{eq:intro:RLRB-RLARBN}-\eqref{eq:intro-regRL} where $\kappa < \kappa_0$, $\XX_0 = X$ and $\XX_1\subset X$ with compact imbedding. Assuming furthermore that 
\beqn\label{eq:MS-KRVestim2}
\exists \, \alpha > 0, \quad \sup_{z \in \Delta_\kappa}   \langle z \rangle^{\alpha}  \| \WW(z)  \|_{\BBB(X)}  < \infty, 
\eeqn
the quantitative exponential asymptotic stability \ref{E31intro} holds (without constructive constants). 
\end{theo}

The proof of Theorem~\ref{theo:MS-KR} is based on a partial (but principal) spectral mapping and Weyl’s theorem (in the spirit of Voigt~\cite{MR595321}) coupled with a simple 
analysis of the first eigenelement problem based on the  irreducibility of the semigroup, but which is really simpler than the deep result on irreducible  semigroup stated in Theorem~\ref{theo:Greiner&co}. 
On the other hand, that approach is unable to tackle the situation when $\Sigma_P^+(\LL)$ is not a singleton. 
One of the main features in Theorem~\ref{theo:MS-KR} and the other results established in \cite{MR3489637,MischErratum} 
is the clear identification of the simple localization of the principal spectrum condition with \eqref{eq:intro-subeigenfunction}.  

\subsubsection{\bf Dynamical and probabilistic approach.} 
\label{subsubsec:intro-dynamical&proba}
\

 \Black 

\smallskip

It is well known from the  mean ergodicity theory of Von {N}eumann and Birkhoff introduced in  the 1930s in \cite{vonNeumann,Birkhoff} that for a bounded semigroup a possible stationary state (and thus a first eigenvector associated to the first eigenvalue $\lambda_1 = 0$) can be obtain through a dynamical approach by establishing that the Ces\`aro mean of the semigroup appropriately converges. 
A classical reference is \cite{MR797411}, see also \cite[Sec.~V.4]{MR1721989} for a short presentation.

\medskip

The existence of invariant measures for Markov chains/processes can be derived through a contraction approach by using coupling arguments reminiscent from the ideas of Doblin~\cite{MR0004409} and Harris~\cite{MR0084889}.
This yields a simplified Krein-Rutman theorem in the Banach lattice of finite measures for Markov operators, providing the existence of $f_1$ whilst $\lambda_1=0$ and $\phi_1=1$ are known by definition.
Doblin's condition is a handy criterion which ensures contraction in total variation norm, and hence existence, uniqueness, and geometric stability of the invariant measure, see for instance~\cite{MR3893207,MR4534707} for this very classical and easy result.
It turns out that this contraction is related to the contraction in Hilbert's metric, see~\cite{Gaubert2015}.
The drawback of Doblin's condition is that it is quite demanding and typically requires the state space to be bounded.
Harris's idea allows an extension to the unbounded setting by localizing Doblin's condition in a ``small set'' which is visited infinitely often.
The return to small sets can be obtained by using a Lyapunov function.
When the Lyapunov function is strong enough for ensuring exponential return, contraction in weighted total variation norm can be established and geometric stability of the invariant measure is inferred~\cite{MeynTweedieI,MeynTweedieII,MeynTweedieIII,MeynTweedie,MR2857021,MR4534707},
leading to the following result (which is made constructive in the two last references).

\begin{theo}\label{theo:Harris-famous}
Consider a positive semigroup $S$ on the Banach space $X=M^1_m(E)$
for some weight function $m:E\to[1,\infty)$.
Suppose that $S$ is conservative, in the sense that

\smallskip
(1) $S_t^*\1=\1$ for all $t\geq0$,

\smallskip
and assume that, for some subset $K\subset E$ on which $m$ is bounded and some time $T>0$,

\smallskip
(2) $S_T^*m\leq \alpha m+\theta\1_K$, \ for some $\alpha\in(0,1)$ and $\theta>0$;

\smallskip
(3) $S_T f\geq \langle f,\1_K\rangle g_0$, \ for all $f\in X_+$ and some $g_0\in X_+$ such that $\langle g_0,\1_K\rangle>0$.

\smallskip
Then there exists a unique probability measure $f_1\in M^1_m$ such that $(\lambda_1=0,f_1,\phi_1=\1)$ is solution to the first eigentriplet problem, and the quantitative exponential stability~\ref{E31intro} holds with constructive constants.
Moreover, some reciprocal implication holds true.
\end{theo}

When only a weak version of the above {\it Lyapunov condition (2)} is available, an extension of the theory to a {\it weakly dissipative} framework is possible and has been developed in 
\cite{TT94,MR2071426,MR2499863,Hairer2016notes,MR4534707} leading to 
existence, uniqueness, but only sub-geometric stability of the invariant measure. 
 %
We also mention that ergodicity properties of Feynman-Kac semigroups were investigated in~\cite{DelMoral2001,MR1988460} and~\cite{Kontoyiannis2003,Kontoyiannis2005}.

\medskip
Using a condition proposed in~\cite[Condition~$\mathcal Z$]{MR1988460}, the Doblin-Harris method was extended to non-conservative semigroups in~\cite{Bansaye2022,MR3449390,Champagnat2023,MR4066299,MR4157907}.
In~\cite{Bansaye2022} necessary and sufficient conditions for the geometric stability of $(\lambda_1,f_1,\phi_1)$ in weighted total variation norm are obtained.
To our knowledge, no extension to the above mentioned weakly dissipative setting is available.

\smallskip\Black

The following result is an immediate consequence of~\cite[Thm.~2.1]{Bansaye2022}.

\begin{theo}\label{theo:BCGM-Harris}
Consider the same situation as in Theorem~\ref{theo:Harris-famous} but relax the conservativeness assumption (1) by the assumption that there exists a function $\phi_0:E\to(0,\infty)$, bounded from above and below by positive constants on $K$, such that $\phi_0\leq m$ on $E$, and satisfying

\smallskip
(1a) $S_T^*\phi_0\geq \beta\phi_0$, \  for some $\beta>0$;

\smallskip
(1b) $\1_K S_t^*\phi_0\leq C \langle g_0,\1_KS_t^*\phi_0\rangle$, \ for all $t\geq0$ and some $C>0$;

\smallskip
and replace the condition $\alpha\in(0,1)$ by $\alpha\in(0,\beta)$ in the assumption (2).

\smallskip
Then, there exists a unique solution $(\lambda_1,f_1,\phi_1)$ to the first eigentriplet problem and the quantitative exponential stability~\ref{E31intro} holds with constructive constants.
Moreover, some reciprocal implication holds true.
\end{theo}

%

\Black

Positivity conditions required for the Doblin-Harris approach are less restrictive than for Birkhoff contraction.
Conversely, unlike contraction in Hilbert's metric, Doblin-Harris method strongly uses the linearity of the operators, and may thus not be easily extendable to nonlinear operator.
However, since it is based on contraction arguments, it can be extended to time-inhomogeneous semigroups~\cite{BCG2019}.
Finally,  the existence of a first eigenmeasure in a non-conservative setting were established in~\cite{MR2834717,MR3102473}
through a Lyapunov function property, a suitable renormalization and a fixed point argument.

\smallskip
 
The key point in this approach is that it provides a constructive rate of convergence while its drawback is that it is somehow restricted to a $M^1_m$ (or $L^1_m$) framework and that some of the conditions (typically (1b) in Theorem~\ref{theo:BCGM-Harris})
are not fully intuitive and may be hard to verify in the applications.

\subsubsection{\bf PDE approaches}

\

\smallskip

 At least as far as the existence issue is concerned, one of the most common way in PDE papers  in order to tackle the existence part of the first eigentriplet problem consists in approximating (by regularization of the coefficients, add of a small viscosity, discretization) the eigentriplet problem, then use the most classical Perron-Frobenius Theorem \cite{MR1511438,Frobenius} or Krein-Rutman Theorem \cite{MR0027128,MR1064315} 
and next to derive appropriate estimates and pass to the limit through  a {\it ``stability argument"}. 

\medskip
 Recently, in order to circumvent the above approximation step, a new abstract and general version of the existence part of the Krein-Rutman theory has been developed by Lions in \cite{PL2} which, as for the early works   \cite{MR839040,MR688774}, 
is also adapted to nonlinear operators and it includes the following statement (in the linear operators framework). 

\begin{theo}\label{theo:PLL}
Consider a Banach lattice $X$  and a linear and bounded operator $\RR : X \to X$ such that 

(i) $\RR : X_+ \to X_+$;

(ii)  $\exists \, g_2 \in X_+ \backslash\{0\}$, $\exists \, C_2 > 0$ such that 
$\RR g_2 \le C_2 g_2$;

(iii) $\mu_1 := \sup \JJ < +\infty$, where 
$$
\JJ := \{ \mu \ge 0; \ \exists h \in K_2, h \ge \mu \RR h + g_2 \}, \quad 
 K_2 := \{ g \in X_+; \, \exists a, \, g \le a g_2 \}; 
$$

(iv) any sequence $(g^n)$ of almost first eigenvectors is relatively (possibly weakly) compact, where we say that $(g^n)$ is a sequence of almost first eigenvectors
if $g^n = \mu^n \RR g^n+ \eps^n$, $(g^n)$ is bounded, $\mu^n \nearrow \mu_1$ and $\eps^n \to 0$.

\smallskip
Then there exists $f_1 \in K_2$ such that   $f_1 = \mu_1 \RR f_1$ and $\| f_1\| = 1$. 

\end{theo}


The statement and proof of Theorem~\ref{theo:PLL} somehow generalize the existence part of the Krein-Rutman theorem  presented in  Theorem~\ref{theo:Greiner&co} because the required  splitting structure and associated power compactness are replaced by the very natural stability principle {\it (iv)}. Applications to elliptic operator with strong or critical confinement property in the whole space $\R^d$ setting are also presented in \cite{PL2}.

\medskip
Let us  mention the huge literature on the characterization of the first eigenvalue by a minmax formula. As explained with more details below, this approach has first been introduced 
in the Courant-Fischer min-max theorem \cite{MR1547416,MR1544417,MR0065391} providing a variational characterization of eigenvalues  in an abstract Hilbert setting for self-adjoint elliptic operators. 
Inspired next by pointwise minmax formula established for simple self-adjoint operators \cite{MR20698,MR111923,MR125319} using a technique which goes back to Picard \cite{Picard1890}, it has been then generalized to non self-adjoint elliptic operators in \cite{MR190527,MR1258192} among others. More recently, the same approach has been generalized to non elliptic operators, see for instance \cite{Coville2020} and the references therein. 

\medskip
On the other hand, and beyond the eigentriplet problem, the convergence towards the first eigenfunction may be proved using the {\it general relative entropy (GRE) method} which has been applied to a large class of evolution PDE in 
\cite{MR2162224} which principle is as follows. Assume that $(\lambda_1,f_1,\phi_1) \in \R \times X \times X'$ is a solution to the first eigenvalue problem, that $\lambda_1 = 0$ (a case to which one can always reduces from the general case by a mere change of operator and unknown), that $X,X'   \subset \Lloc^1(E)$ and then define the  generalized relative entropy  
$$
\JJ(f) := \int_E j(f/f_1) \, f_1 \, \phi_1 \, d\mu
$$
for any given convex function $j : \R \to \R_+$. For any solution $f(t) \in X$ to the (appropriate) evolution PDE, one may establish (at least formal) the identity
\beqn\label{eq:GRE}
\JJ(f(t)) + \int_0^t \DD_\JJ(f(s)) \, ds = \JJ(f(0)), \qquad \forall \, t \ge 0, 
\eeqn
where $ \DD_\JJ \ge 0$ is the  associated generalized dissipation of relative entropy, so that $\JJ$  is a Lyapunov functional (it is decreasing along the flow associated to the evolution PDE). Under suitable positivity hypothesis, one has
$\DD_\JJ(f) = 0$ if and only if $f \in \hbox{\rm Vect}(f_1)$, and then one may deduce from \eqref{eq:GRE} and some lower semicontinuity assumption on the operator $\DD_\JJ$ that $f(t) \to c f_1$ as $t\to\infty$ (without rate and with $c \in \R$). The GRE method is of course connected to $j$-divergence in information theory and statistics \cite{MR164374,MR335152,MR672884,MR1430403} and to $j$-entropy in probability and PDE theory \cite{MR2081075,MR2450356}, where however here it is crucial to identify the associated operator $\DD_\JJ$ and that this last one enjoys suitable properties.

\subsubsection{\bf Hypotheses and proof}
\label{subsubsec:intro-Hyp&Proof} 

\

\smallskip
We now briefly discuss the strategy of the proof of Theorem~\ref{theo:main-intro} and how it is connected to the above material. Additional comments will be made in the corresponding Sections~\ref{sec:ExistenceKR} to \ref{sec:QuantitativeStabilityKR}. 
As already said, the first eigenvalue problem  is mainly split into three steps: existence, geometry and asymptotic stability. From a general point of view, our approach is more general than the initial  Krein-Rutman theorem  
as well as less abstract than the usual semigroup school approach. 
We believe it is more intuitive and handy for the possible applications since it is presented as a series of estimates to be checked and  the necessary assumptions are made clearer at each step. 

\medskip
$\bullet$ Concerning the existence of a solution to the first eigentriplet problem, our result improves the previous known results because (1) only weak compactness property is needed (while Theorems~\ref{theo:Greiner&co} \& \ref{theo:MS-KR} require strong compactness assumptions), (2) it is more flexible than Theorems~\ref{theo:Greiner&co}, \ref{theo:MS-KR}, \ref{theo:BCGM-Harris} \& \ref{theo:PLL} (the two first ones being restricted to the generator of a strongly continuous semigroup, the third one being restricted to a $M^1_m$ framework and involving the tricky condition (1b) and the last one being somehow restricted to a weighted $L^\infty$ framework), (3) it applies to weakly dissipative cases (so that no spectral gap is needed). We present two different proofs: one based on a stationary problem approach and another one based on a dynamical problem approach (with which we are able to tackle the weakly dissipative case). 

\smallskip
Our  stationary problem approach mixes in a first step the (clearly formulated) approximation argument of \cite[proof of Thm.~12.15]{MR3616245} together with the stability argument of \cite{PL2}, where it is worth emphasizing that  
the condition $ \kappa<\kappa_0$ in Theorem~\ref{theo:main-intro} is nothing but a practical (and possibly constructive) condition ensuring that assumption $s(\BB) < s(\LL)$ holds in Theorem~\ref{theo:Greiner&co}. On a second step, 
we exhibit several practical situations where the required stability condition is fulfilled recovering as a particular case the existence part in Theorems~\ref{theo:Greiner&co} \& \ref{theo:BCGM-Harris}. We would like to point out here that the splitting hypotheses \eqref{eq:intro:stopDyson}-\eqref{eq:intro-regSL} on the semigroup is a generalization of the Lyapunov condition (2) in Theorem~\ref{theo:Harris-famous} on the semigroup which in turn generalizes the classical 
Lyapunov condition on the generator, namely for instance
$$
\LL^* \psi_2  \le \kappa \psi_1 + K \psi_0
$$
with $\psi_i \in X'$,  $\psi_1, \psi_2 \ge \psi_0$ together with $\psi_2 \le \psi_1$ (super Lyapunov condition), $\psi_2 = \psi_1$ (standard Lyapunov condition), 
 $\psi_2 \ge \psi_1$ (weak Lyapunov condition). We refer to \cite{MR4534707} and to Sections~\ref{sec:ExistenceKR} and \ref{sec:DynamicalExistenceKR} for further discussions on that question. 
 
 \smallskip
On the other hand, our dynamical approach mixes the splitting method yet alluded above together with some argument picked up from Von Neumann \& Birkoff mean ergodic theory in the spirit of but in a more elaborate way than in \cite[Sec.~6]{MR4534707}. 


%
%

\medskip


$\bullet$ The proof about the geometry of the principal eigenvalue problem in Theorem~\ref{theo:main-intro} is a refinement of many arguments already developed in the literature.
More precisely, the uniqueness of the first eigentriplet $(\lambda_1,f_1,\phi_1)$ and the strict positivity of the eigenvectors is established by taking up again in a more general setting some arguments developed in  \cite{MR2114128,MR3489637,MR4265692}. 
The subgroup structure of the boundary point spectrum $\Sigma_P^+(\LL)$ is next established under suitable (but not very restrictive) geometrical properties on the Banach lattice $X$,
these ones being always true for the usual examples we have in mind and that we have already listed above. The proof mainly mimics the usual proof (as for instance presented in \cite[Sec.~14.3]{MR3616245}) but it is less abstract and more general. Especially, the proof does not refer to the notions of ideals, quasi-interior points or Calkin algebra nor uses the Kakutani lattice isomorphism theorem but rather uses the simpler notion of strict positivity (defined by duality) and some convenient structural properties of the signum operator. 
In order to go one step further and to prove the triviality property  $\Sigma_P^+(\LL) = \{\lambda_1\}$, we propose one quite original approach (which we believe to be new at this level of generality) based on an {\it reverse Kato's inequality condition} of $\LL$ (by refining some arguments picked up from \cite{MR3489637,MR4265692})  and some more standard ones based on an {\it aperiodicity condition} on the semigroup $S_\LL$, on a localization of the point spectrum condition or on a {\it quasi-compactness condition} on the semigroup $S_\LL$. 
  
 \medskip

$\bullet$ Finally, the proof on the asymptotic stability of the first eigenvector picks up and mixes some 
spectral analysis, dynamical system, entropy method and Doblin-Harris coupling arguments. 
On a first step, we mainly rewrite some very classical dynamical system results mixed together with some arguments coming from the General Relative Entropy method in order to get our mean ergodicity and ergodicity results which are really general and very little demanding about the trajectories. 
We also rewrite the most classical result about the exponential asymptotic stability (without constructive constants) of the first eigenfunction proposing a very simple (and self-contained) proof which does not make any references to abstract notions as Calkin algebra, essential spectrum or essential growth bound.
Last, we  adapt the Doblin-Harris approach as qualitatively formulated in \cite{MR2857021,MR4534707,Bansaye2022} in order to get the quantitative asymptotic stability of the first eigenfunction with constructive constants.

\medskip

%
%
%
%

\subsection{Some examples of applications}
\label{sec:applicationse}
 
\ \smallskip

The abstract Krein-Rutman theory developed in these notes and alluded above have been cooked up in order to answer to the first eigenvalue problem for PDE. 
We show its efficiency by applying it to several examples of evolution PDE. These examples must be thus considered both as a motivation and an  illustration of simultaneously developed abstract theory. 


\subsubsection{Parabolic equations} In Part~\ref{sec:application1:diffusion}, we are interested in parabolic equations  in divergence form 
$$
\partial_t f =   \partial_i (a_{ij} \partial_j f) +  b_j \partial_j f + \partial_i (\beta_i f)  + c f  \quad\hbox{in}\quad (0,\infty) \times \Omega, 
$$
on the function $f = f(t,x)$, $t \ge 0$, $x \in \Omega$, with general conditions on the coefficients $a_{ij}$,  $b_j$, $\beta_i$,  $c$ and in both  the case of a bounded domain $\Omega \subset \R^d$ (and we then complement the equation with a Dirichlet boundary condition) and the case $\Omega = \R^d$. 
The importance of parabolic equations for Physics, Chemistry, Biology and Economy  modeling is well known and we do not discuss it here. 
We consider the four following cases. 

\smallskip
$\bullet$ For a bounded domain $\Omega \subset \R^d$, we consider a general elliptic operator in divergence form 
$$
\LL f := \partial_i(a_{ij} \partial_j f) + b_i \partial_i f +  \partial_i(\beta_i f) +  cf, \quad f \in H^1_0(\Omega),
$$
under the very general assumption about  the regularity of the coefficients  $a_{ij}  \in L^\infty(\Omega)$,  $a_{ij}  \ge \nu \delta_{ij}$, for some $\nu >0$,
$b_i, \beta_j \in L^r(\Omega)$, $c  \in L^{r/2}(\Omega)$, $r > d$. 
%

\smallskip
$\bullet$ In the case when $\Omega = \R^d$,  we focus first our analysis by considering 
$$
\LL f := \Delta f + b \cdot \nabla f + c f, \quad f \in H^1(\R^d),
$$
with drift $b \in \Lloc^\infty(\R^d)$, potential $c \in \Lloc^2(\R^d)$ and a confinement condition that (roughly speaking) we impose through the properties  $c \to - \infty$ as $|x| \to \infty$ and $b$ is dominated by $c$ at the infinity. 
A typical case is given by $c \sim - |x|^\gamma$ and $b \sim x |x|^{\beta-1}$ as $ |x| \to \infty$, with  $\gamma > \max(0,\beta-1)$.

\smallskip
$\bullet$ Still in the case when $\Omega = \R^d$,  we next consider the similar problem
$$
\LL f := \Delta f + b \cdot \nabla f + r c f, \quad f \in H^1(\R^d),
$$
with now $c \in C_0(\R^d)$, $b \in C_0(\R^d)$ and $r \in \R_+$ a parameter.  
That hypotheses correspond to a critical confinement case and we further assume that $r > 0$ is large enough. 

\smallskip
$\bullet$ In the case when $\Omega = \R^d$ again, we finally consider the elliptic operator 
$$
\LL f := \Delta f + b \cdot \nabla  f + c f,
$$
with the drift confinement 
$$
b = \nabla U, \quad U(x) = \frac1\beta\langle x \rangle^\gamma, \quad \gamma > 0, 
$$
and with $c$ dominated by $b$ at the infinity. We further assume $c \ge \Div b$ when $\gamma \in (0,1]$. It is worth emphasizing that this corresponds to a perturbation of the classical
Fokker-Planck operator associated to the potential $U$.

\medskip
For each of these operators  we are able to complete the existence, geometric and stability program as stated in Theorem~\ref{theo:main-intro}, with constructive estimates on the first eigentriplet solution and more or less explicit rate of convergence to the first eigenfunction.
Few suitable additional assumptions on the coefficients and on the regularity of $\Omega$ as well as  the precise results will be discussed in the corresponding sections. 

\medskip
 The first eigenvalue problem in the three first situations has been studied in \cite[8th and 9th courses]{PL2} which inspired our study and to which we refer for motivations and possible extensions. 
 Since mainly the existence issue is considered in \cite{PL2}, our results supplement the previous analysis by tackling the geometry of the principal spectrum and the exponential asymptotic stability of the first eigenfunction. 
 On the other hand, the fourth situation in the conservative case {\Blue (namely, when we additionally take $c := \Div b$)} is very classical and we refer to \cite{MR2386063,MR2381160,MR3779780,MR3488535,MR4265692} and the references therein. We believe that the extension to a non conservative case as considered here is new.

\smallskip
Of course, when the operator $\LL$ is the Laplace operator or more generally is a self-adjoint elliptic operator, there exists a huge literature about the analysis of its spectrum and in particular about its first eigenvalue problem because among other things this is related to the ground state problem in quantum mechanic. We do not have the precise historical reference where similar results to the ones developed here are established for the first time. 
We may for instance refer to the contributions by Poincaré \cite{MR1505534} and by Courant and Hilbert \cite{MR1544417,MR0065391}. We also refer to the textbook \cite[Thm~8.38]{MR0473443} for the quite general and modern proof which mixes minimization technique, strong maximum principle and Hilbert structure arguments. 
It is worth mentioning that in earlier works,  the Krein-Rutman theorem has been proved using elementary ODE method when considering the Sturm-Liouville operator (in dimension $d=1$), see for instance \cite{Bocher}.
Still for a self-adjoint elliptic operator, the Courant-Fischer min-max theorem \cite{MR1547416,MR1544417} gives a variational characterization of eigenvalues through  Rayleigh quotient \cite{MR0016009}
and the Weyl theorem  \cite{MR1511560,MR1511670,MR0057422,MR133586} provides some information about the distribution of the eigenvalues.
More specifically, some constructive lower bound on the best constant in Poincaré inequality and thus on the first eigenvalue may be obtain through the Faber-Krahn \cite{zbMATH02598582,MR1512244} isoperimetric inequality as presented  in \cite{MR1427759}, see also
 Polya-Svzego \cite{MR0043486,MR133059} and Payne-Weinberger \cite{MR104047,MR149735}. Other results on that direction but based on the Lyapunov condition are obtained in \cite{MR2386063,MR2381160} and we also refer   to \cite{MR3155209} and the references therein. 

\smallskip
 On the other hand, in the case of an elliptic operator which is not self-adjoint the first result on the principal eigenvalue problem seems to be Protter, Weinberger \cite[Rk.~2]{MR190527} who consider the case of smooth domain and coefficients (without precise statement about the regularity) and use minmax formula and the Krein-Rutman Theorem~\ref{theo:intro-KRresult}, see also \cite{MR0236516}. 
Next,  Chicco \cite{MR280858,MR0310462} establishes the existence, uniqueness and some monotony properties of the first eigenvalue-eigenfunction in the weak solutions framework of Stampacchia \cite{MR192177,MR0251373} with mild regularity assumptions on the coefficients and which corresponds to the framework we will consider here (when we will consider the case of a bounded domain). 
These works have been followed by several papers by Donsker and Varadhan \cite{MR425380,MR361998} and next by the famous work of Berestycki,  Nirenberg, Varadhan \cite{MR1258192} opening a new field of research. These last works are mainly based on strong maximum principle technique, see \cite{MR762825}. We also mention the recent works by  Champagnat and Villemonais \cite{MR4066299,Champagnat2023} where similar results to ours for smooth enough coefficients are established using a variant of the probabilistic Doblin-Harris argument as already mentioned in Section~\ref{subsubsec:intro-dynamical&proba}. 
We also emphasize that in the conservative case, the long time behavior problem has been widely studied and some constructive estimates has been obtained in 
\cite{MR889476,MR954373,MR893137,MR1704435,MR1751701} by the mean of log-Sobolev inequality, in  \cite{MR2386063,MR2381160,MR1856277,MR4265692}
by the mean of Poincaré inequality and in \cite{MR3779780,MR4265692} by the mean of semigroup arguments.

\subsubsection{Transport equation}
In Part~\ref{sec:Transport}, we are interested in the general transport equation
\beqn\label{eq:Transport-introEq} 
\partial_t f  + a  \cdot \nabla_y f = \KKK[f] - K f   \quad\hbox{in}\quad (0,\infty) \times \OO, 
\eeqn
on the function $f = f(t,y)$, $t \ge 0$, $y \in \OO$, with $\OO \subset \R^D$, $D \ge 1$, a smooth open connected set.
We assume that $a : \OO \to \R^D$, 
$K : \OO \to \R_+$, 
and that the collision operator $\KKK$ is linear and defined by 
$$
\KKK [g](y) := \int_{\OO} k(y,y_*) \,  g(y_*)  \, dy_* ,
$$
for some kernel $k : \OO \times \OO \to \R_+$.  When $ \OO \not=\R^D$, we complement the equation with a boundary condition on 
the  incoming boundary $\Sigma_-$ which writes 
$$
(\gamma_- f )(t,y) = 
 \RR_\OO [f(t,\cdot)](y) + \RR_\Sigma[\gpf(t,.) ] (y) \,\, \hbox{ on }  \,\,(0,\infty) \times \Sigma_-, 
$$
where $\gamma_\pm f$ are the trace functions on the incoming and out going set $\Sigma_\pm$ and 
$$
\RR_\OO [g](y) := \int_\OO  r_\OO(y,y_*) \, g(y_*) \, dy_*, \quad
\RR_\Sigma [h](y) := \int_{\Sigma_+}   r_\Sigma (y,y_*) \, h(y_*)  \, d\sigma_{\! y_*}, 
$$
for some kernels $r_\OO : \Sigma_- \times \OO \to \R_+$, $ r_\Sigma : \Sigma_- \times \Sigma_+ \to \R_+$.   All the (quite usual) notations will be explained at the begin of Part~\ref{sec:Transport}. 
It is worth emphasizing here that this framework in particular covers the cases of the renewal equation, the growth-fragmentation equation and the kinetic linear Boltzmann equation on which we will come back below.  
This framework is motivated by and generalizes the transport theory developed in \cite{MR274925,MR872231,MR1022305,MR2150445,MR3212249}.

\smallskip  
In a first step, we consider a very general vector field $a$ by assuming that it satisfies the usual Sobolev regularity condition of  DiPerna-Lions transport theory \cite{MR1022305}. 
We also make general assumptions on $\RR_\OO$ and $\RR_\Sigma$, but a very strong and somehow restrictive positivity condition on $\KKK$. Such an equation can be motivated by the abstract transport theory developed
 \cite{MR872231} as well as non-local reaction-diffusion models
\cite{MR4157907,Coville2020,Li2017} and  selection-mutation models in changing environment \cite{Forien2022,Henry2023}.
Under these general conditions and additional ones we will detail later,  we are able to solve the existence and geometrical part of the first eigenvalue problem and to prove an ergodicity result (without rate of convergence)
generalizing some similar results obtained in \cite{MR4157907,Coville2020,Li2017}.

\smallskip  
%
Because of the strong positivity condition made on $\KKK$, the above mentioned result does not apply to the growth-fragmentation equation and the kinetic linear Boltzmann equation.
We thus consider separately these important particular cases in the two next parts.
Other singular jump kernels lacking strong positivity can appear in other models, for instance in neurosciences~\cite{Dumont2020}, and must also be treated through a specific study.


\medskip
Another  related model is the age-structured renewal equation 
\bean
&&\partial_t f  + \partial_y f =  - K f   \quad\hbox{in}\quad (0,\infty) \times (0,\infty),
 \\
 && f(t,0)= \RR [f(t,\cdot) ](y) =  \int_0^\infty r(y_*) f(t,y_*) dy_*. 
 \eean
It corresponds to the case  $D=1$, $\OO = (0,\infty)$, $a=1$, $\RR_\Sigma = 0$, $\Sigma_- = \{ 0 \}$ and $\KKK = 0$ in the transport equation
\eqref{eq:Transport-introEq}. The age structured equation is very popular because it is useful for describing dynamic of populations  \cite{Sharpe1977,Arking,MR1636703,MR772205,MR860959}
and simple neuronal dynamic \cite{MR1239892,MR2576373}.
The long time behaviour can be analyzed though Laplace transform technique \cite{Feller1941,MR0210154,Iannelli}, relative entropy method \cite{MR1946722,MR2162224,Gwiazda2006}, spectral analysis tool \cite{Webb1984,MR763356,MR3489637,MR3850019,MR3859527} and
Doblin approach \cite{BCG2019,MR3900224,MR3893207}.
Because $\KKK = 0$, our previous result on the first eigenvalue problem does not apply. 
We just briefly observe that the method can be applied on the dual equation, thus guaranteeing the existence of $(\lambda_1,\phi_1)$, and then that the validity of Doblin's condition ensures the existence and uniqueness of the triplet $(\lambda_1,f_1,\phi_1)$, its positivity, and the exponential ergodicity.



\smallskip

\Black

\subsubsection{Growth-fragmentation equation}

In Section~\ref{part:application3:GF}, we consider the growth-fragmentation equation
$$
\partial_t f = \LL f = \GG f + \FF f   \quad\hbox{in}\quad (0,\infty) \times (0,\infty),
$$
with the growth operator $\GG f = - \partial_x (af)$ and the fragmentation operator 
$$
(\FF f)(x) = \int_x^\infty k(y,x) f(y) \, dy - K(x) f(x),
\quad K(x) := \int_0^x k(y,x) \frac{y}{x} dy.
$$
Since the work of Diekmann, Heijmans and Thieme~\cite{Diekmann1984}, many authors studied this equation by using various methods.
We can mention, among many others, {\cite{MR2162224,MR2114128,MR2536450,DoumicJauffret2010,MR2832638,MR3030711,MR3928121}} for studies based on  suitable weak distance, entropy and functional inequalities, \cite{MR860959,Pichor2000,MR2821681,MR2915575,BernardGabriel17,BernardGabriel20,MR3489637,MKB2022} in the framework of positive semigroups, 
\cite{BertoinWatson18,Bertoin2020} for a probabilistic approach via the Feynman-Kac formula, \cite{Bansaye2022,MR4312814,GabrielMartin,MR4897760} for Doblin-Harris's method,
 and also~\cite{Franco2023} for a recent new approach based on the reformulation of the equation as an abstract renewal problem.
Our aim here is not to treat the most general cases of coefficients, but rather to illustrate the variety of the possible behaviors of the equation together with the efficiency and flexibility of the method developed in the first sections.
We thus focus on a specific case of fragmentation operator, namely the equal mitosis kernel
$$
k(x,y) = 2 K(x) \delta_{x/2} (dy) = 4 K(x) \delta_{2y} (dx),
$$
so that the equation writes 
\[
\partial_t f (t,x)= - \partial_x \big(a(x)f (t,x)\big) -  K(x) f(t,x) + 4 K(2x) f(t,2x).  
\]

In particular, we are interested in the case when the growth rate $a$ is such that $a(2x)=2a(x)$ for all sizes $x > 0$, for which the boundary spectrum is not trivial and the solutions then exhibit persistent oscillations in time.
When this condition is not satisfied, we recover the more usual exponential convergence to the first eigenfunction.

\smallskip

We also aim at studying the variant of this equation where a variability $v$ is introduced as a growth speed parameter which is inherent to any individual,
in the spirit of~\cite{MR1946722,Rotenberg1983} where such a variable is added in the renewal equation.
More precisely we consider the growth-fragmentation equation with variability $v\in[1,2]$ and the equal mitosis division kernel which reads
\[
\partial_t f (t,x,v)= - v \partial_x \big(a(x)f (t,x,v)\big) -  K(x) f(t,x,v)  + 4 \int_1^2 K(2x) \wp(v,v_*) f(t,2x,v_*) dv_*.
\]
This model was introduced in~\cite{DHKR}, and then also considered in~\cite{Olivier2017}.
We prove that, unlike the case without variability, it exhibits exponential relaxation to the first eigenfunction even when $a(2x)=2a(x)$ for all $x > 0$.

\smallskip

\subsubsection{Kinetic linear Boltzmann equation} 
In Part~\ref{part:application4:Kequation}, we are interested in another important subclass of transport equations, namely in the  kinetic linear Boltzmann equation
\beqn\label{eq:Boltzmann_lin-introEq}
\partial_t f + v \cdot \nabla_x f - \nabla_x \Phi(x) \cdot  \nabla_v f = \KKK[f]- K f \quad\hbox{in}\quad (0,\infty) \times \OO, 
\eeqn
on the function $f = f(t,x,v)$, $t \ge 0$, $(x,v) \in \OO = \Omega \times \VV$, $\Omega \subset \R^d$, $\VV\subset \R^d$, $d \ge 1$. We assume that $K : \OO \to \R_+$, that $\KKK$ is a linear integral operator defined by
$$
 \KKK [g] := \int_{\R^d} r k(x,v,v_*) \,  g(v_*)  \, dv_* ,
 $$
 for some real number $r > 0$ and some kernel $k : \Omega \times \VV \times \VV \to \R_+$, and that $\Phi$ is a space confining potential $\Phi : \Omega \to \R$. We restrict our analysis to the case $\VV := \R^d$ and $\Omega$ 
is either the torus $\Omega := \T^d$ (and we assume $\Phi = 0$) or it is the whole space $\Omega := \R^d$ (and we assume that $\Phi$ is a power function). 
This equation is very famous because it provides a model for neutron transport theory in nuclear reactors \cite{MR0225547,osti_4074688} and for cells migration in a chemotactic gradient \cite{MR609371}. 
We refer to \cite{MR743736,MR855310,MR933978,MR910950,MR1612403}  for a mathematical analysis of the neutron transport equation and its diffusive approximation and to \cite{MR1843975,MR2065025} for the same concerning kinetic models for chemotaxis.
Because the  linear integral operator $\KKK$ is local in the position variable, this problem does not fall in the class of transport equation covered by the Krein-Rutman theorem established in Part~\ref{sec:Transport}
and a specific analysis is necessary. Under suitable positivity and regularity  conditions on the kernel, we are able to complete the existence, geometric and stability program as stated in Theorem~\ref{theo:main-intro}, with constructive estimates in the torus case, generalizing and improving previous works \cite{MR0120814,MR230531,MR259662,MR731337,MR2063095,MR1151539,MR2158174,MR2237499,MR3192457}
where spectral analysis arguments are used and \cite{MR3449390} based on a probability approach. It is worth emphasizing that these works are concerning the same equation in a bounded domain with no-flow boundary condition. 
Most of the literature is about the conservative case (when $\lambda_1 = 0$ and $\phi_1 = 1$) which has been tackled by the mean of spectral analysis method \cite{MR3084493,MR3192457,MR4317820}, of entropy method \cite{MR2243948,MR3101201}, of geometric control method \cite{MR3479064,dietert2023quantitative}, of hypocoercivity method \cite{MR2215889,MR3324910,MR4226238} or of Doblin-Harris coupling approach \cite{MR4063917}.

\medskip

\subsubsection{Kinetic Fokker-Planck equation} In Part~\ref{part:application5:KFPequation}, we consider a kinetic Fokker-Planck equation
\beqn\label{eq:KFP-introEq} 
\partial_t f  + v \cdot \nabla_y f =  \Delta_v f + b \cdot \nabla_v f + cf  \quad\hbox{in}\quad (0,\infty) \times \OO, 
\eeqn
on the function $f = f(t,x,v)$, $t \ge 0$, $(x,v) \in \OO   := \Omega \times \R^d$, $\Omega \subset \R^d$ is a bounded domain,  $b : \OO \to \R^d$ is a given vector field
and  $c : \OO \to \R$ is a given function. 
In contrast with the previous part, collisions are typically modeled by a Fokker-Planck operator
$\Delta_v f + \Div_v (vf)$ which takes into account a  thermal bath of (Gaussian) white-noise, see Kolmogorov \cite{MR1503147},   instead of the integral collisional operator $\KKK [f] - Kf $ in the linear Boltzmann equation \eqref{eq:Boltzmann_lin-introEq}. The above equation is  complemented with the Maxwell boundary condition$$
\gamma_- f  =   \alpha(x) \DD_x \gamma_+ f   + \beta(x)  \Gamma_x \gamma_+ f, 
$$
where $\gamma_\pm f$ stand for the outgoing and incoming trace functions, $\alpha$ and $\beta$ are accommodation coefficients, $\DD_x$ is a boundary diffusive operator and 
$ \Gamma_x$ is the  specular reflection operator. All these classical objects will be precisely defined in Part~\ref{part:application5:KFPequation}. We refer to \cite{MR222474,MR771811,MR875086,MR1200643,MR1634851,MR1949176,MR2721875,MR2562709,MR4253803}  for a mathematical analysis of the kinetic Fokker-Planck equation or related problems. 
Under suitable boundedness and regularity  conditions on the coefficients we are able to complete the existence, geometric and stability (without constructive estimates ) program as stated in Theorem~\ref{theo:main-intro}, generalizing the previous works \cite{MR4347490,MR4756950} (partially based on  \cite{MR2248986,MR3382587,MR4412380}) where similar results are established for the same kind of equation 
 in a bounded domain with no-flow boundary condition. 
From a technical point of view, our proof is based on trace results as those developed in   \cite{MR2721875}, 
boundary estimates picked up from \cite{MR1301931,MR2721875,MR4581432} and regularity estimates recently obtained in \cite{MR2786222,MR3923847,MR4453413}.
We also emphasize that in the conservative case, many works have been done related to hypocoercivity and constructive rate of convergence to the steady state in  \cite{MR1787105,MR1946444,MR1969727,MR2034753,MR2130405,MR2562709} or more recently in 
\cite{MR3324910,MR3488535,MR4069622,MR4113786,MR4776290}.

\subsubsection{Mutation-selection equation}

Last, in Section~\ref{sec:appl:selec-mut}, we consider the mutation-selection evolution equation
\begin{equation*}
    \partial_tf=\LL f=J*f-W(x)f   \quad\hbox{in}\quad (0,\infty) \times \R^d, 
\end{equation*}
This nonlocal-diffusion equation appears for instance in the modeling of genetic variability in evolutionary biology. 
In this context, $f=f(t,x)$ represents the density of a population, at time $t \ge 0$, of phenotypical trait $x$ on the
multi-dimensional phenotypic trait space $\R^d$.
The rate of change in $f$ per generation is  given by the convolution term with kernel $J$ which models the mutations, and the fitness function $-W$ which stands for the difference between birth and death.
This model has been widely used in the literature; we refer, for example, to the works of Kimura~\cite{Kimura1965}, Lande~\cite{Lande1975}, Fleming~\cite{Fleming1979} and B\"urger~\cite{Burger-book} as examples of biological applications.

On the mathematical analysis point of view, the Krein-Rutman problem was investigated by Bürger in~\cite{MR923493,BurgerBomze} and more recently by Coville and co-authors~\cite{Coville2010,Li2017}, as well as by Alfaro et al. in~\cite{MR4622852} where a quantified spectral gap is obtained for symmetric kernels $J$.
A main difference of this equation compared to more classical ``local'' diffusion models, where the convolution is replaced by a Laplacian, is that the first eigenvector $f_1$ might be a measure with atoms~\cite{MR923493,BurgerBomze,Coville2013}.
Some conditions are then needed relating $W$ and $J$ for guaranteeing that the first eigenvector is an eigenfunction~\cite{MR4622852,MR923493,Li2017}.

All the above mentioned results deal with kernels $J$ which are densities, namely absolutely continuous with respect to the Lebesgue measure.
In our study, we allow the convolution kernel to have a singular part.
In Section~\ref{ssec:appl5:regular} we extend the results of the literature to the case of a small enough singular part.
In Section~\ref{ssec:appl5:singular} we consider a specific kernel which is purely singular, supported by the canonical axes of $\R^d$, and we extend the recent result of Velleret~\cite{Velleret2023} to more general confining functions~$W$.

 \medskip
\subsection{Organization of the paper}

The paper is organized in two main parts: the sections~\ref{sec:ExistenceKR} to~\ref{sec:QuantitativeStabilityKR} are dedicated to the development of the abstract results about the Krein-Rutman problem, and the last sections~\ref{sec:application1:diffusion} to~\ref{sec:appl:selec-mut} aim at illustrating the applicability of these results to various linear positivity preserving PDE.

More precisely, we start with the existence part of the Krein-Rutman theorem, namely the conclusion~\ref{S1}.
This question is addressed through a stationary approach in Section~\ref{sec:ExistenceKR} and through a dynamical approach in Section~\ref{sec:DynamicalExistenceKR}.
Section~\ref{sec:Irreducibility} is devoted to the stronger conclusion of uniqueness of the first eigentriplet in the sense of~\ref{S2}, as well as to the mean ergodic property~\ref{E1intro}.
In Section~\ref{sec:geo2}, we are interested in the geometry of the boundary point spectrum, deriving conditions that guarantee~\ref{S31}, \ref{S32} or~\ref{S33}, as well as in the ergodic properties~\ref{E2intro} and~\ref{E31intro}.
Finally, in Section~\ref{sec:QuantitativeStabilityKR}, we tackle the problem of quantifying the conclusions~\ref{S33} and~\ref{E3intro} by using constructive contraction arguments of the Doblin-Harris type.

The purpose of the last six sections is to apply the theory developed in the first sections to the examples of PDE presented in Section~\ref{sec:applicationse}:  parabolic equations (Section~\ref{sec:application1:diffusion}), transport equations with integral terms (Section~\ref{sec:Transport}) and in particular growth-fragmentation equations (Section~\ref{part:application3:GF}) and kinetic equations (Section~\ref{part:application4:Kequation}), kinetic Fokker-Planck equations (Section~\ref{part:application5:KFPequation}), and purely integral mutation-selection equations (Section~\ref{sec:appl:selec-mut}).

\bigskip


\bigskip

\section{Existence through a stationary problem approach}  
\label{sec:ExistenceKR}


In this part we provide a general existence result for the first eigentriplet problem by considering a family of approximating  stationary problems and using a 
stability argument, in particular we prove \ref{S1}.
{\Blue  For that purpose we introduce some conditions \ref{H1}, \ref{H2}, \ref{H3}, as well as a splitting structure \ref{HS1} useful for verifying \ref{H3}.} 
We start by presenting the basic material about the Banach lattice framework and conclude with a comparison with several previous works. 

 \medskip
\subsection{The Banach lattice framework}
\label{subsec:BanachLattice}
We start recalling the Banach lattice framework by stating (most of the time without proof) some well-known facts that one can find in reference monographs as \cite[Chapitre II: Espaces de Riesz]{MR0219684} or \cite{MR0423039,MR839450,MR2178970,MR3616245}.

\medskip
{\bf Banach lattice. } A real Banach lattice is a real Banach space $(X,\| \cdot \|)$ endowed with a partial order denoted by $\ge $ (or $\le$)
such that the following holds: 

(1) The set $ X_+ := \{ f \in X; \,\, f \ge 0 \}$ is a nonempty convex cone (compatibility of the order with the vector space structure). 

(2)  For any  $f \in X$, there exist some unique positive part $f_+ \in X_+$ and negative part $f_- \in X_+$ such that $f = f_+ - f_-$ which are minimal: $f = g - h$, $g,h \ge 0$ imply $g\ge f_+$ and $h \ge f_-$ (generation and properness of the positive cone). We set $|f| := f_+ + f_- \in X_+$ the absolute value of $f \in X$. 

(3) For any $f,g \in X$, $|f| \le |g|$ implies  $\| f \| \le \| g \|$ (compatibility of norm and order structures).

\medskip
Under these assumptions, one can show that 

\smallskip
- The convex cone $X_+$ is closed, pointed $X_+ \cap (-X_+) = \{ 0 \}$ and generating $X = X_+ - X_+$.

\smallskip
- The lattice operations $f \mapsto f_+$,  $f \mapsto f_-$ and $f \mapsto |f|$ are continuous ($1$-Lipschitz). 

\smallskip
- The order intervals $\{ h \in X; \ g \le h \le f \}$ are closed and bounded for any given $f,g \in X$, $f \ge g$. 

\smallskip
It is worth emphasizing that one commonly defines the supremum and infimum operations by 
$$
f \vee g := g + (f-g)_+ \ge f,g, \quad f \wedge g := g - (g-f)_+  \le f,g, 
$$
for any $f,g \in X$, and these operations can be used as an alternative way for defining a Banach lattice (the lattice structure refers indeed to these   supremum and infimum operations). 
 We may note the following elementary formulas
\beqn\label{eq:propprieteDUcone2}
f_+ \wedge f_- = 0, \quad \| |f| \| = \| f \|, \quad \forall \, f \in X.  
\eeqn
We write $f \perp g$ when $|f| \wedge |g| = 0$ or equivalently when $|f|+|g| = |f| \vee |g|$. 

\medskip
{\bf Dual Banach lattice. } 
On the dual space $X'$, we may naturally associate a dual order $\ge$ (or $\le$) by writing for $\varphi \in X'$
$$
\varphi \ge 0 \,\,\,  (\hbox{or } \varphi \in X'_+) \quad\hbox{iff}\quad
\forall \, f \in X_+ \,\,\, \langle \varphi,f \rangle \ge 0.
$$
For $\varphi \in X'$, there exist some unique $\varphi_\pm \in X'_+$ such that $\varphi = \varphi_+ - \varphi_-$
which also satisfy (and are defined by)
$$
\forall \, f \in X_+, \quad 
\langle \varphi_\pm , f \rangle = \sup_{0 \le g \le f} \langle \pm\varphi, g \rangle.
$$
One can show that the above conditions (1), (2) and (3) of a Banach lattice are fulfilled, and thus $X' = (X', \| \cdot \|_{X'},\ge)$ is a Banach lattice.  
We observe that for any $f \in X_+$ there exists $f^* \in X'_+$ such that 
\beqn\label{eq:f*fge0}
\langle f^*, f \rangle = \| f \|^2 = \| f^* \|^2_{X'}, 
\eeqn
as a classical corollary of the Hahn-Banach dominated extension theorem.
 Moreover, for any $f \in X$, 
\beqn\label{eq:fge0}
f \ge 0 \quad\hbox{iff}\quad \langle \varphi,f \rangle \ge 0, \  \forall \, \varphi \in X'_+,
\eeqn
as an immediate application of the Hahn-Banach separation theorem.
In other words, the restriction to $X$ of the dual order in $X''$ associated to the order defined (by duality) on $X'$ is nothing but the initial order, in particular
  the positive cone $X'_+$ is weakly $*$ closed.
 
\medskip
{\bf The functional framework : The duality bracket. } We consider two Banach lattices $X,Y$ such that $X=Y'$ with $Y$ separable or such that $Y=X'$. 
We emphasize on the facts that 
\bear\label{eq:exist1-positive=dualpositive1}
 \hbox{for} \ f \in X:&& f \in X_+ \hbox{ iff } \langle f,\varphi \rangle \ge 0, \ \forall \, \varphi \in Y_+, 
\\ \label{eq:exist1-positive=dualpositive2}
\hbox{for} \ \varphi \in Y: &&  \varphi \in Y_+ \hbox{ iff } \langle f,\varphi \rangle \ge 0, \ \forall \, f  \in X_+,
\eear
which are immediate consequences of \eqref{eq:fge0} and of the definition of the dual order. 

\medskip
{\bf Examples. } For the space $C_0(E)$, the order   is defined by $f \ge 0$ iff $f(x) \ge 0$ for any $x \in E$. 
For a space $L^p(E,\EE,\mu)$, $1 \le p \le \infty$, the order  is defined by $f \ge 0$ iff $f(x) \ge 0$ for $\mu$-a.e. $x \in E$. 
For the space $M^1(E)$, the order is defined by $f \ge 0$ iff in the Hahn decomposition $f = f_+ - f_-$, there holds $f_- = 0$, or equivalently, by duality: $f \ge 0$ iff 
$\langle f, \varphi \rangle \ge 0$ for any $\varphi \in C_0(E)$, $\varphi \ge 0$.

 \smallskip 
 Because confinement will play a major role in our analysis, we will use some weighted version of the above space associated to a weight (continuous or Borel measurable) function 
 $ m : E \to (0,\infty)$  that we introduce now. We recall that $E$ always denotes a $\sigma$-compact metric space, and we write $E = \cup E_R$, with $E_R \subset E_{R+1}$, $E_R$ compact. In that context, we 
 write $x_n \to \infty$ if for any $R \ge 1$ there exists $n_R$ such that $x_n \notin E_R$ for any $n \ge n_R$.

\smallskip
$\bullet$ We denote by $C_{m,0}(E)$ the space 
$$
C_{m,0} (E):= \{ \varphi \in C(E); \ |\varphi(x)|m(x) \to 0 \ \hbox{as} \ x \to \infty \}
$$
endowed with the norm $\| \varphi \|_{C_{m,0}} := \| \varphi m \|_{L^\infty}$. 

\smallskip
$\bullet$ We denote by $M^1_m(E) := (C_{m^{-1},0} (E))'$ the  associated space of Radon measures. 

\smallskip
$\bullet$ We denote by $L^p_m(E) = L^p_m(E,\EE,\mu)$ the space 
$$
L^p_m(E) := \{ f \in \Lloc^1(E); \ \| f \|_{L^p_m} := \| fm \|_{L^p} < \infty \}. 
$$
It is worth emphasizing that $L^p_m(E,\EE,\mu) = L^p(E,\EE,m^{p}\mu)$ when $p \in [1,\infty)$.

\medskip
{\bf Positive operator.} 
We denote by $\BBB(X)$ the set of linear and bounded operators on $X$. 
We also denote by $\KK(X)$ the subspace of compact operators.
We say that a bounded operator  $A \in \BBB(X)$ is positive, and we write $A \ge 0$,  if 
$$
A f \in X_+, \quad \forall \, f \in X_+.
$$  
We will also sometimes abuse notations by writing $A \in \BBB(X_+)$  for meaning that $A \ge 0$.
For a positive operator $A \in \BBB(X)$, we have 
\beqn\label{eq;absAf<Aabsf}
|Af| \le A|f|, \quad \forall \, f \in X,
\quad\hbox{and}\quad
\| A \| =  \sup_{0 \le f \in B_X} \| A f \|, 
\eeqn
where $B_X$  is the unit closed ball.  
More generally, we have 
\beqn\label{eq;AfveeAgLEAfg}
(Af ) \vee (Ag)  \le A(f \vee g), \quad \forall \, f,g \in X.
\eeqn

\smallskip
%

For $X$ and $Y$ in duality, and $A\in\BBB(X)$ and $A^*\in\BBB(Y)$ in duality, in the sense that
$$
\langle A f, \phi \rangle = \langle f, A^* \phi \rangle, \quad 
\forall \, f \in X, \ \phi \in Y,
$$
there holds
\beqn\label{eq:exist1Age0iffA*ge0}
A \ge 0 \quad\hbox{iff}\quad A^* \ge 0.
\eeqn
 Let us present the elementary and classical but instructive proof of the direct implication, the reciprocal way being similar. 
We assume thus $A \ge 0$. We take $\varphi \in Y_+$ and we define $\psi := A^* \varphi$. We then take $f \in X_+$ and we define 
$g := A f$, so that $g \ge 0$ by assumption. We compute 
\bean
\langle \psi, f \rangle 
= \langle A^* \varphi, f \rangle 
= \langle \varphi, A f \rangle =  \langle \varphi, g \rangle \ge 0. 
\eean
Since $f \in X_+$ is arbitrary, we get $\psi \in Y_+$, and thus $A^* \ge 0$.

\medskip
{\bf Semigroup, generator and spectrum. } 
In this work, a semigroup $S = S(t) = (S_t)$ on $X$ will always denote a semigroup of linear and bounded operators on a Banach lattice $X$ which trajectories are 

- either strongly continuous, namely, the mapping $t \mapsto S_t f$ is continuous for the norm of $X$ for any fixed  $ f \in X$;

- either weakly $*$ continuous, namely $X = Y'$ for some separable Banach lattice $Y$ such that  $\forall \, f \in X$, $\forall \, \phi \in Y$, $t \mapsto \langle S_t f, \phi \rangle_{X,Y}$ is continuous
and $\forall \, t \ge 0$, $\forall \, \phi \in Y$, $f \mapsto \langle S_t f, \phi \rangle_{X,Y}$ is continuous.
That is in particular the case when  there exists a strongly continuous semigroup $P$ on $Y$ such that $S_t = P^*_t$ for any $t \ge 0$. 

\smallskip
For a semigroup $S$, we denote by $\LL$ its generator and $D(\LL)$ the associated domain, and thus we sometimes write $S= S_\LL$. 
 We also denote the iterated domain defined recursively by $D(\LL^k) := \{ f \in D(\LL^{k-1}), \, \LL f \in D(\LL^{k-1}) \}$ for any $k \ge 2$ and 
$D(\LL^\infty) := \bigcap_{k \ge 1} D(\LL^k)$.  
We recall that $D(\LL)$ is dense in $X$ and the graph of $\LL$ is closed in $X \times X$. 
We define the growth bound 
\beqn\label{def:GrowthBoundLL}
\omega = \omega(S) :=   \limsup_{t \to \infty} \frac{1}{t} \log \| S(t) \|  \in \R \cup \{-\infty\}, 
\eeqn
so that 
\beqn\label{eq:SMwt}
\forall \, \omega' > \omega, \quad \exists \, M \ge 1, \quad \| S(t) \|_{\BBB(X)} \le M \, e^{\omega' t}, \quad \forall \, t \ge 0, 
\eeqn
and $\omega$ is the infimum of $\omega' \in \R$ such that  \eqref{eq:SMwt} holds. We say that $S$ is a semigroup of contractions when $S$ satisfies  \eqref{eq:SMwt} with $M=1$ and $\omega' = 0$. 

The resolvent set $\rho(\LL)$ is the set of $z \in \C$ such that  if $z- \LL  : D(\LL) \to X$ is bijective and its inverse belongs to $\BBB(X)$.  
We define the resolvent operator by 
\beqn\label{eq:Exist1-DefRR}
\RR(z) = \RR_\LL(z)  := (z-\LL)^{-1}, \quad \forall \, z \in \rho(\LL), 
\eeqn 
and the spectrum by 
$\Sigma(\LL) := \C \backslash \rho(\LL)$. Denoting the  half complex  plane of abscissa $\alpha \in \R$ by 
\beqn\label{eq:Exist1-DefDeltaalpha}
\Delta_\alpha := \{ z \in \C; \,\, \Re e (z) > \alpha \},
\eeqn
we have $\rho(\LL) \supset \Delta_\omega$ and, for any $z \in \Delta_\omega$, there holds 
\beqn\label{eq:Exist1-DefRepresentationRR}
 \RR(z) = \int_0^\infty S(t) e^{-z t} \, dt.
\eeqn
 
\medskip
{\bf Positive semigroup.} 
We say that a semigroup $(S_t)$ 
on a Banach lattice $X$ is positive if
$$
S_t \ge 0, \quad \forall \, t \ge 0.
$$

\begin{lem}\label{lem:sec2-EquivPositiveS}
For a semigroup $S$  on a Banach lattice $X$, there is equivalence between  

\smallskip
(a) $S$ is positive;

\smallskip
(b) the associate resolvent  operator $\RR$ is positive:  $\RR (\kappa) \ge 0$ for all $\kappa > \omega$ (or for all sufficiently large $\kappa$). 

\end{lem}

It is immediate from Hille's identity \eqref{eq:Exist1-DefRepresentationRR} that (a) implies (b). 
The reciprocal implication comes from the relation  $S(t) = \lim_{n\to\infty} [n/t \RR(n/t)]^n$ at the foundation of the Hille-Yosida theory, see for instance  \cite[Thm.~I.8.3]{MR710486}.

 \medskip
\subsection{Existence part of the Krein-Rutman theorem}
\label{subsec-Exist1-KRtheorem}
From now on in this section, we consider a Banach lattice $X$ and an operator $\LL$ with dense domain and closed graph.
Our goal is mainly to prove the existence part for the primal problem in  the Krein-Rutman theorem, namely 
\beqn\label{eq:1stEVP}
\exists \, \lambda_1 \in \R, \ \exists \, f_1 \in X_+ \backslash \{0 \}, \quad \LL f_1 = \lambda_1 f_1.
\eeqn
We will also discuss the existence part for the dual problem at the end of the section. 

\smallskip
We first assume   

\begin{enumerate}[label={\bf(H1)},itemindent=12mm,leftmargin=0mm,itemsep=1mm]
\item\label{H1} $\exists \, \kappa_1 \in \R$ such that $\lambda-\LL$ is invertible and $(\lambda-\LL)^{-1} : X_+ \to X_+$ for any $\lambda \ge \kappa_1$.
\end{enumerate}

Note that an operator $\LL$ satisfying~\ref{H1} is sometimes called a {\it resolvent positive} operator after the paper of Arendt~\cite{Arendt1987}.

\smallskip
We then set  
\beqn\label{eq:exist1-defI}
\II := \{ \kappa \in \R; \ \lambda-\LL \hbox{ is invertible}, \ (\lambda-\LL)^{-1} \ge 0 \ \hbox{ for any } \ \lambda \ge \kappa\},  
\eeqn
which is a non empty and non upper bounded interval due to \ref{H1}. We finally set 
\beqn\label{eq:exist1-deflambda1}
\lambda_1 := \inf \II \in [-\infty, \kappa_1].
\eeqn

\smallskip
For the sake of completeness, we recall now some  
general facts about $\II$ and $\lambda_1$ when $\LL$ is the generator of a positive semigroup.    
We also refer to \cite[Sec.~1.b, Chap.~VI]{MR1721989} or \cite[Chapter~12]{MR3616245} and the references therein for more details.  

\begin{lem}\label{lem:Exist1-RkSG}  When $\LL$ is the generator of a positive semigroup $S=S_\LL$, then 

\ {\bf (i)} \ref{H1} automatically holds with any $\kappa_1 > \omega(S)$, so that $\lambda_1 \le \omega(S)$; 

\smallskip
\ {\bf (ii)} $\Sigma(\LL) \cap \Delta_{\lambda_1} = \emptyset$ and the representation formula \eqref{eq:Exist1-DefRepresentationRR} holds true for any $z \in \Delta_{\lambda_1}$; 

\smallskip
\ {\bf (iii)} it may happen that $\lambda_1 = - \infty$.

\end{lem}

The important property {\bf(ii)} is probably due to \cite{MR617977}. 
 
\begin{proof}[Proof of Lemma~\ref{lem:Exist1-RkSG}.]
The claim {\bf (i)} is an immediate consequence of the representation formula \eqref{eq:Exist1-DefRepresentationRR} for any $\kappa_1 > \omega(S)$ and the positivity of $S(t)$ for any $t \ge 0$  (that is nothing but Lemma~\ref{lem:sec2-EquivPositiveS}). 

\smallskip 
We prove {\bf (ii)}. Take $\lambda > \lambda_1$. From the classical identity 
 $$
S(t)e^{-\lambda t} - I = (\LL -\lambda) \int_0^t S(s)e^{-\lambda s} \, ds, \quad \forall \, t \ge 0, 
$$
and the positivity property of $S$, we have 
$$
0 \le V(t) :=  \int_0^t S(s)e^{-\lambda s} \, ds = \RR(\lambda ) - \RR(\lambda ) S(t)e^{-\lambda t}  \le \RR(\lambda ), 
$$
for any $t \ge 0$. From that estimate, we get $\| V(t) \| \le   \| \RR(\lambda )\|$. 
For any $z \in \Delta_\lambda$, an integration by part yields
$$
\int_0^t e^{-z s} S(s) \, ds = e^{-(z-\lambda )t}V(t) + 
(z-\lambda) \int_0^t e^{-(z-\lambda) s } V(s) \, ds.
$$
The estimate on $V$ makes possible to pass to the limit $t\to\infty$ in the above identity, and  we deduce 
$$
\UU (z) := \int_0^\infty e^{-z s} S(s) \, ds  = 
(z-\lambda) \int_0^\infty e^{-(z-\lambda) s } V(s) \, ds \in \BBB(X). 
$$
In that situation, one classically knows that $z \in \rho(\LL)$ and $ (z-\LL)^{-1} = \UU(z)$. We have thus established $\Sigma(\LL) \cap \Delta_\lambda = \emptyset$ and 
we conclude the proof of {\bf (ii)} by observing that \eqref{eq:Exist1-DefRepresentationRR}  is then nothing but the above formula.

\smallskip
{\bf (iii)} On $L^p(0,1)$, $1 \le p < \infty$, the translation semigroup defined for $a > 0$ by 
$$
S(t) f (x) := f(x+at) {\bf 1}_{x+at \le 1}, \quad \forall \, t \ge 0, \, x \in (0,1), 
$$
is strongly continuous and positive. Since $S(t) \equiv 0$ for any $t \ge 1/a$, we have $\omega (S) = - \infty$, and thus $\lambda_1 = - \infty$ because of {\bf (i)}. 
\end{proof}

For further discussion, we give some probably classical results about the condition \ref{H1} and some equivalent definitions of the set $\II$.

%

\begin{lem}\label{lem:Existe1-H2bis} The operator $\LL$ satisfies \ref{H1} if and only if the operator $\LL^*$ satisfies \ref{H1}.
Furthermore, under condition \ref{H1} for $\LL$ (or  $\LL^*$), we have  
\beqn\label{eq:exist1-defIBIS}
\II = \II_i, \quad \forall \, i=2,3,4, 
\eeqn
with 
\bean
\II_2 &:=&  \{ \kappa \in \R; \ \lambda-\LL \hbox{ is invertible}  \hbox{ for any } \ \lambda \ge \kappa\},
\\
\II_3 &:=&  \{ \kappa \in \R;  \ \lambda-\LL^* \hbox{ is invertible},     \ (\lambda-\LL^*)^{-1} \ge 0 \ \hbox{ for any } \ \lambda \ge \kappa\}, 
\\
\II_4 &:=& \{ \kappa \in \R; \ \lambda-\LL^* \hbox{ is invertible}  \hbox{ for any } \ \lambda \ge \kappa\}.
\eean
\end{lem} 

\begin{proof}[Proof of Lemma~\ref{lem:Existe1-H2bis}.]  
The equivalence of condition \ref{H1} for the operators $\LL$ and $\LL^*$ is an immediate consequence of the identity  $\rho(\LL) = \rho(\LL^*)$ (see for instance \cite[Thm.~III.6.22]{MR0203473}) and 
 the fact that  $(\lambda - \LL)^{-1}  : X_+ \to X_+$  iff $(\lambda - \LL^*)^{-1}  : Y_+ \to Y_+$  as recalled in \eqref{eq:exist1Age0iffA*ge0}. 
 As a consequence, we have  $\II=\II_3$ and  $\II_2 = \II_4$. 

 \smallskip
 We obviously have $\II_2 \subset \II$ and let us show the reverse inclusion.  
 We denote $\RR = \RR_\LL$. On the one hand, for any $z_0 \in \rho(\LL)$ and any $z \in \C$, $|z-z_0| < \|\RR(z_0)\|^{-1}$, we know that
\beqn\label{eq:Uz=sumUz0}
\RR(z) = \RR(z_0) \sum_{k=0}^\infty (z_0-z)^k \RR(z_0)^k,
\eeqn
which gives a proof of the fact that resolvent set $\rho(\LL)$ is open and that $\RR$ is an holomorphic function on $\rho(\LL)$. 
Formula \eqref{eq:Uz=sumUz0} also ensures that for $\lambda_0,\lambda\in \R$, the condition $\RR(\lambda_0) \ge 0$ implies that $\RR(\lambda) \ge 0$ provided that  $\lambda_0-\lambda > 0$ is small enough and thus $\RR(\lambda) \ge 0$ for any $\lambda$ in the non upper bounded connected component of the set $\rho (\LL) \cap \R$  thanks to a continuation argument.  In particular, $\II$ is an open set and $\II = \II_2$.
\end{proof}

We next assume 
 
\begin{enumerate}[label={\bf(H2)},itemindent=12mm,leftmargin=0mm,itemsep=1mm]
\item\label{H2} $\exists \, \kappa_0 \in \R$  
such that 
$\inf \II \ge \kappa_0$.
\end{enumerate}

We do not further consider in these notes the  case when $\inf \II  = - \infty$
 and moreover we will particularly focus on the possibility to exhibit constructive  lower bound $\kappa_0$. 

\smallskip
We point out several conditions under which \ref{H2} is satisfied. {\Blue 
In the applications part and depending of the case, 
we will be able to establish (mostly constructively) one of these different criteria.}

\begin{lem}\label{lem:Existe1-Spectre2bis} Condition \ref{H2} holds under one of the four following conditions

\quad {\bf (i)} $\exists \, \kappa_0 \in \R$,   
$\exists \, \phi_0 \in Y_+ \backslash \{0\}$ such that $  \LL^*\phi_0  \ge \kappa_0 \phi_0$, which means 
$$
\forall \, f   \in D(\LL) \cap  X_+ , \quad \langle \phi_0, (\kappa_0-\LL) f \rangle \le 0;
$$

\quad {\bf (ii)} $\exists \, \kappa_0 \in \R$,  
$\exists \, f_0 \in  X_+ \backslash \{0\}$ such that $  \LL f_0  \ge \kappa_0 f_0$,   
which means 
$$
\forall \, \phi   \in D(\LL^*) \cap  Y_+  , \quad \langle (\kappa_0-\LL^*) \phi,  f_0 \rangle \le 0;
$$

\quad {\bf (iii)} $\LL^*$ is the generator of a positive semigroup $S^* = (S^*_t)$ and
$$
\exists \, \kappa_0 \in \R, \ \exists \, \phi_0 \in Y_+ \backslash \{ 0 \}, \ \exists \, T >0 \  \hbox{ such that } \  S^*_T \phi_0   \ge e^{\kappa_0 T}  \phi_0 ;
$$
\quad {\bf (iv)} $\LL$ is the generator of a positive semigroup $S = (S_t)$ and
$$
\exists \, \kappa_0 \in \R, \ \exists \, f_0 \in X_+ \backslash \{ 0 \}, \ \exists \, T >0 \  \hbox{ such that } \  S_T f_0   \ge e^{\kappa_0 T}  f_0.
$$

 \end{lem} 

\begin{proof}[Proof of Lemma~\ref{lem:Existe1-Spectre2bis}.]
In the  four cases, we claim that $\kappa_0 \notin \II$ and thus $\inf \II \ge \kappa_0$. We argue  by contradiction, 
assuming  $\lambda_1 < \kappa_0$, so that $\kappa_0 \in \II=\II_i$ for any $i=2,3,4$. 
 
\smallskip
We assume  {\bf (i)}.  
For any $g \in X_+$, we define $f := (\kappa_0-\LL)^{-1} g \in X_+$ and we compute
$$
0 \le \langle \phi_0, g \rangle =  \langle \phi_0, (\kappa_0-\LL) f \rangle \le 0. 
$$
That implies $\langle \phi_0, g \rangle = 0$ for any $g \ge 0$, so that $\phi_0 = 0$ and a contradiction. 


\smallskip
We assume first that {\bf (iii)} holds for any $T > 0$.  
For any $f \in  D(\LL) \cap  X_+ \backslash \{0\}$, we compute 
$$
 \langle \phi_0, (\kappa_0-\LL) f \rangle = - \frac{d}{dt} \langle \phi_0,  e^{-\kappa_0 t} S_t f \rangle \le 0, 
$$
which is precisely {\bf (i)}. 
We assume now that {\bf (iii)} holds. If $\kappa_0 \in \II$, for any $g \in X_+$, we may   define $f=(\kappa_0-\LL)^{-1}g\in X_+ \cap D(\LL)$ and from condition {\bf (iii)},  we have
\[
0\leq\langle e^{-n\kappa_0 T}S_{nT}f-f,\phi_0\rangle=\Big\langle(\LL-\kappa_0)\int_0^{nT} e^{-\kappa_0 t}S_t f\,dt,\phi_0\Big\rangle,
\]
for any   $n\in\N$. From the very definition of $f$, we also have 
\[
(\LL-\kappa_0)\int_0^{nT} e^{-\kappa_0 t}S_t f\,dt = \int_0^{nT} e^{-\kappa_0 t}S_t (\LL-\kappa_0) f \,dt = -\int_0^{nT} e^{-\kappa_0 t}S_t g \,dt\leq0.
\]
The two pieces of information together imply 
\[
\Big\langle\int_0^{nT} e^{-\kappa_0 t}S_t g \,dt,\phi_0\Big\rangle=0.
\]
Passing to the limit $n\to\infty$ thanks to Lemma~\ref{lem:Exist1-RkSG}-{\bf (ii)} and using \eqref{eq:Exist1-DefRR}-\eqref{eq:Exist1-DefRepresentationRR}, we obtain 
\[
0 = \Big\langle\int_0^{\infty} e^{-\kappa_0 t}S_t g \,dt,\phi_0\Big\rangle = \langle f,\phi_0 \rangle = \langle g,(\kappa_0-\LL^*)^{-1}\phi_0\rangle.
\]
That implies $(\kappa_0-\LL^*)^{-1}\phi_0=0$ since $g$ is arbitrary, what is not possible since $\phi_0 \not=0$. 
The proofs of \ref{H2} under assumption {\bf (ii)} or {\bf (iv)} are similar to the previous ones and are thus skipped.
 \end{proof}

\begin{rem}\label{rem:Existe1-Spectre2bis} (1) In practice, we may build $f_0$ or $\phi_0$ through an explicit computation or use a barrier function and strong maximum principle techniques. 
We refer to  Lemma~\ref{lem:Irred-barrier+W->H2} for a possible general result in that direction. {\Blue More precisely, in all the applications we will use the criteria {\bf (i)} or {\bf (ii)}, except in Section~\ref{ssec:Transport:KR} where the criteria {\bf (iv)} is used.} 

(2) When {\bf (ii)} holds with $f_0 \in X_+ \backslash \{0 \} \cap D(\LL)$ and $\LL$ is the generator of a positive semigroup $S$, then {\bf (iv)} holds for any $T > 0$. In that case, we may indeed compute 
$$
e^{-\kappa_0 T} S_T   f_0 - f_0 = \int_0^T e^{-\kappa_0 t}  S_t (\LL-\kappa_0) f_0 ds \ge0. 
$$
\end{rem}

\begin{lem}\label{lem:Existe1-Spectre2} Under conditions \ref{H1} and \ref{H2}, there hold
\beqn\label{eq:lambda1inab}
\lambda_1 \in [\kappa_0,\kappa_1]
\eeqn
and 
\beqn\label{eq:lambda1approx}
\exists \, \lambda_n \searrow \lambda_1, \ \exists \, \hat f_n \in D(\LL) \cap X_+,\  \eps_n := \lambda_n \hat f_n - \LL \hat f_n \ge 0, \ \| \hat f_n \| = 1, \ \| \eps_n \| \to 0.
\eeqn
\end{lem} 


\begin{proof}[Proof of Lemma~\ref{lem:Existe1-Spectre2}.]
We obviously have $\lambda_1 \le \kappa_1$ from assumption \ref{H1} and  $\lambda_1 \ge \kappa_0$ by assumption \ref{H2}, so that 
\eqref{eq:lambda1inab} is proved. 

\smallskip
Consider now a sequence $(\lambda_n)_{n \ge 2}$ such that $\lambda_n \searrow \lambda_1$ as $n\to\infty$. We eventually have $\| \RR(\lambda_n) \| \to \infty$ as $n\to\infty$, where we denote by $\RR = \RR_\LL$ the resolvent of $\LL$.
On the contrary, we would have $\| \RR(\lambda_{n'}) \| \le M$ for some subsequence  $\lambda_{n'} \searrow \lambda_1$ and some constant $M >0$. 
Because of \eqref{eq:Uz=sumUz0} this implies that $(\lambda_{n'}-\eps,\lambda_{n'}) \subset \II$ for any $n'$ and some $\eps >0$, and this is
 in contradiction with the definition of $\lambda_1$. 
 The blow up $\| \RR(\lambda_n) \| \to \infty$ means that 
$$
\exists \, f_n \in D(\LL), \, \exists \, g_n \in X, \quad \RR(\lambda_n) g_n = f_n, \ \| f_n \| \to \infty, \ \| g_n \| \le 1.
$$
By splitting $g_n = g_n^+-  g_n^-$, we get
$$
f_n = \RR(\lambda_n) g_n^+ - \RR(\lambda_n) g_n^-
$$
with 
$$
\| g_n^\pm \| \le 1 \quad \hbox{and}\quad (\| \RR(\lambda_n) g_n^+ \| \to \infty \hbox{ or } \| \RR(\lambda_n) g_n^- \| \to \infty).
$$
Changing notations, we have thus 
$$
\exists \, f_n \ge 0, \exists \, g_n \ge 0, \quad \RR(\lambda_n) g_n = f_n, \ \| f_n \| \to \infty, \ \| g_n \| \le 1.
$$
We get \eqref{eq:lambda1approx} by defining $\hat f_n := f_n / \| f_n \|$ and $\eps_n := g_n/ \| f_n \|$.
\end{proof}

We learn a very similar proof in \cite{PL2}, from which our own proof is adapted. The same type of arguments can also be found in  \cite[proof of Thm~12.15]{MR3616245}. 
 
\medskip
We finally assume that 

\begin{enumerate}[label={\bf(H3)},itemindent=12mm,leftmargin=0mm,itemsep=1mm]
\item\label{H3}  for any sequence $(\hat f_n)$ of $X$ such that \eqref{eq:lambda1approx} holds, there exist $f_1 \in X_+ \backslash \{0\}$ and a subsequence $(\hat f_{n'})$ 
such that $\hat f_{n'} \wto f_1$ for the weak convergence or the weak $*$ convergence. 
\end{enumerate}

 \medskip
 We discuss several situations in which assumption \ref{H3} holds. 
We start with a very classical framework formalized for instance by Voigt   \cite{MR595321},  see also  Karlin~\cite[Cor.~1]{MR0114138} or Sasser~\cite{Sasser1964}  for earlier similar situations and results, which is however somehow restrictive since it is based on a strong compactness property assumed at the level of the associated semigroup of operators.

\begin{lem}\label{lem:H3abstract-StrongCbis}
We assume that $\LL$ generates a positive semigroup $S$,  
that \ref{H2} holds for a constant $\kappa_0 \in \R$ and that 
there exists $T> 0$  such that the splitting 
\[
S_T = V_T + K_T
\]
holds with $\| V_T \|_{\BBB(X)} \le e^{\kappa T}$, $\kappa < \kappa_0$, and $K_T \in \KK(X)$.  Then condition \ref{H3} holds for the primal and the dual problems.
\end{lem}
 
\begin{proof}[Proof of   Lemma~\ref{lem:H3abstract-StrongCbis}.]
The condition \ref{H1} holds because of Lemma~\ref{lem:Exist1-RkSG}-{\bf (i)}.
  Let us then consider three sequences $(\lambda_n)$, $(\hat f_n)$ and $(\eps_n)$ satisfying  \eqref{eq:lambda1approx}. 
Integrating along the rescaled flow, this yields
\bean
e^{-\lambda_n T}S_T\hat f_n-\hat f_n &=& \int_0^T e^{-\lambda_n t} S_t (\LL - \lambda_n) \hat f_n dt 
\\
&=& -  \int_0^T e^{-\lambda_n t} S_t\eps_n\,dt =: \tilde\eps_n, 
\eean
which also reads
\[
V_T \hat f_n + K_T \hat f_n - e^{\lambda_n T} \hat f_n = e^{\lambda_n T} \tilde\eps_n.
\]
Since $e^{\lambda_n T} \ge e^{\kappa_0 T} > e^{\kappa T} $, the operator $e^{\lambda_n T}-V_T$ is invertible with inverse $Q(\lambda_n) :=   (e^{\lambda_n T}-V_T)^{-1}$ uniformly bounded 
and converging in $\BBB(X)$ to $Q(\lambda_1) = (e^{\lambda_1 T}-V_T)^{-1}$. We thus have
\[
 \hat f_n = w_n + v_n, \quad w_n := Q(\lambda_n)  K_T  \hat f_n, \quad v_n :=  - Q(\lambda_n) e^{\lambda_n T} \tilde\eps_n ,
 \]
with $\| v_n \|_X \to 0$ and $(w_n)$ relatively compact in $X$. There exist thus a subsequence $( \hat f_{n_k})$ and $g \in X$ such that $K_T  \hat  f_{n_k} \to g$ and next
$$
w_{n_k} - Q(\lambda_1) g = (Q(\lambda_{n_k})-Q(\lambda_1)) K_T  \hat  f_{n_k} + Q(\lambda_1) ( K_T  \hat  f_{n_k} - g) \to 0.
$$
We deduce that $ \hat f_{n_k} \to f_1 := Q(\lambda_1) g$ strongly in $X$. 
Because of the positivity and normalized properties of $\hat f_{n}$, we get $f_1 \in X_+$, $\| f_1 \|_X =1$, and we conclude that  \ref{H3} holds for the primal problem.
Observing that the dual semigroup $S^*$ satisfies $S^*_T = V^*_T + K^*_T$ with 
 $\| V^*_T \|_{\BBB(Y)} \le e^{\kappa T}$ and $K^*_T \in \KK(Y)$, the same proof implies that \ref{H3} holds for the dual problem.
 \end{proof}

 {\Blue Although this classic result can probably be  used in many of the applications to PDE that we present in the second part of this paper, we develop alternative arguments below  because: (1)~some of the applications presented in the sequel certainly do not fall within the scope of the above result (we are thinking in particular to the parabolic equation under weakly dissipative conditions in Section~\ref{subsec-diffusionwith drift}  and   to the mitosis equation with non-mixing growth rate in Section~\ref{ssec:GF:singular})  and (2)~the framework developed below is much simpler and direct to use (it requires no information on the semigroup itself and little information about the properties of the operator).}   
In the six next lemmas, we will assume that  \ref{H1}-\ref{H2} holds associated to some constants $\kappa_i \in \R$, $\kappa_0 < \kappa_1$, and we always make the following splitting structure hypothesis   

\begin{enumerate}[label={\bf(HS1)},itemindent=14mm,leftmargin=0mm,itemsep=1mm] 
\item\label{HS1} there exists a splitting  $\LL = \AA + \BB$ such that $\BB-\alpha$ is invertible for any   $\alpha \ge \kappa_0$ and   
\beqn\label{eq:exist1-defVVWW}
\VV (\alpha) :=  \sum_{i=0}^{N-1}  (\RR_\BB(\alpha) \AA)^i  \RR_\BB(\alpha), \quad \WW(\alpha) := (\RR_\BB(\alpha) \AA)^N,
\eeqn
 are  bounded in $\BBB(X)$ uniformly  with respect to $\alpha \ge \kappa_0$ and for some $N \ge 1$, where we recall that $\RR_\BB(\alpha) := (\alpha-\BB)^{-1}$ is the resolvent of $\BB$. 
 \end{enumerate}

\medskip
We first present a result also based on a strong compactness property which is assumed to hold however at the level of the resolvent operator. 
It turns out that we will be able to use that result in  most of the applications. 

\begin{lem}\label{lem:H3abstract-StrongC} (1) We assume \ref{H1}-\ref{H2}-\ref{HS1} 
and there exists  $N \ge 1$ such that
\beqn\label{lem:H3abstract-StrongCompact-c1}
\WW(\alpha)  \hbox{ is strongly compact locally uniformly on } \alpha \ge \kappa_0,
\eeqn
 in the sense that if $ \alpha_n \to \alpha$,  $\alpha_n \ge \kappa_0$,  and $(g_n)$ is a bounded sequence in $X$, then there exist    $f \in X$ and a subsequence $(g_{n_k})$ such that   $\WW(\alpha_{n_k}) g_{n_k} \to f$ strongly in $X$.
 Then condition \ref{H3} holds.
 
 (2) We assume \ref{H1}-\ref{H2} and \ref{HS1} where  $\RR_\BB(\alpha)$ is  bounded   uniformly in $\alpha \ge \kappa_0$, $\AA \in \BBB(X)$ 
and $\WW(\alpha) \in \KK(X)$ for any fixed $\alpha \ge \kappa_0$ and some $N \ge 1$. Then condition \ref{H3} holds both for the primal and the dual problems.
\end{lem}

 \begin{rem}\label{rem:H3abstract-compactStrong}  
 
(1)  The property \eqref{lem:H3abstract-StrongCompact-c1} holds if we assume $\WW(\alpha) : X \to \XX_1$ is bounded uniformly in  $\alpha \ge \kappa_0$ and $\XX_1 \subset X$ with strong compact embedding. 

\smallskip
(2) The property \eqref{lem:H3abstract-StrongCompact-c1} holds if we assume \ref{H1}-\ref{H2}-\ref{HS1} together with the facts that $\RR_\BB(\alpha)$ and   $ \RR_\BB(\alpha)\AA$  are  bounded   uniformly in $\alpha \ge \kappa_0$
and 
$\WW(\alpha) \in \KK(X)$ for any fixed $\alpha \ge \kappa_0$. Consider indeed $ \alpha_n \to \alpha$,  $\alpha_n \ge \kappa_0$,  and $(g_n)$ a bounded sequence in $X$. On the one hand,  there  exist $f \in X$ and a subsequence $(g_{n_k})$ such that   $\WW(\alpha ) g_{n_k} \to f$ strongly in $X$, because $\WW(\alpha) \in \KK(X)$. On the other hand, using the resolvent identity $\RR_\BB(\lambda) - \RR_\BB(\mu) = (\mu-\lambda) \RR_\BB(\lambda)  \RR_\BB(\mu)$, we have 
\bean
\WW(\alpha) - \WW(\alpha_n) 
= (\alpha_n - \alpha) \sum_{j=1}^{N} (\RR_\BB(\alpha)\AA )^{N-j} \RR_\BB(\alpha)  ( \RR_\BB(\alpha_n) \AA)^{j} \to 0, 
\eean
so that $\WW(\alpha_{n_k} ) g_{n_k} \to f$ strongly in $X$, and \eqref{lem:H3abstract-StrongCompact-c1} holds true. 


\smallskip
{\Blue (3) In the applications part, we will mainly use Lemma~\ref{lem:H3abstract-StrongC} for establishing \ref{H3}.
It is very flexible and allows for generalizations such as using two or more different splittings, see Section~\ref{ssec:GF:variability} for an example.}

\smallskip \Cyan
 (4) We observe that the same conclusion holds when  $\RR_\LL = \VV + \WW \RR_\LL$ for two family of operators $\VV(\alpha)$ and $\WW(\alpha)$ such that  $\VV(\alpha)$ is  bounded in $\BBB(X)$ uniformly on $\alpha \ge \kappa_0$ and $\WW$ satisfies \eqref{lem:H3abstract-StrongCompact-c1}. That follows from the fact that  \eqref{eq:lambda1approx} writes equivalently 
$$
\hat f_n = \VV(\lambda_n) \eps_n +  \WW(\lambda_n) \hat f_n, \quad \lambda_n \searrow \lambda_1, \ \eps_n \to 0, \ \| \hat f_n \| = 1, 
$$
and we may then proceed exactly in the same way as below. 
 \end{rem}

\begin{proof}[Proof of Lemma~\ref{lem:H3abstract-StrongC}.] 
We first assume (1). Taking advantage of the splitting structure \ref{HS1}, we  write equation \eqref{eq:lambda1approx} as 
\beqn\label{eq:existenceStat-H3abstract-StrongC}
(\lambda_n - \BB) \hat f_n = \AA \hat f_n + \eps_n, 
\eeqn
or equivalently
$$
 \hat f_n = \RR_\BB(\lambda_n) \AA \hat f_n +  \RR_\BB(\lambda_n) \eps_n.
$$
Iterating that last identity and using the notations \eqref{eq:exist1-defVVWW}, we get 
\beqn\label{eq:1hatfn=iterate}
\hat f_n = w_n + v_n,  \quad w_n :=  \WW(\lambda_n) \hat f_n ,  \quad  v_n := \VV (\lambda_n) \eps_n. 
\eeqn
We observe that $(w_n)$ is strongly relatively compact from \eqref{lem:H3abstract-StrongCompact-c1} and  $\| \hat f_{n}  \|_{X} = 1$, so that there exist a subsequence $(w_{n_k})$ 
and $f_1 \in X$ such that $ w_{n_k} \to f_1$ strongly in $X$. 
Since $v_n \to 0$ strongly in $X$, we deduce that  $ \hat f_{n_k} \to f_1$ strongly in $X$, {\Cyan in particular $f_1 \ge 0$ and $\| f_1 \|_X = 1$}. 
We conclude that condition \ref{H3} holds as in the proof of  Lemma~\ref{lem:H3abstract-StrongCbis}. 

\smallskip
We next assume (2). As observed in Remark~\ref{rem:H3abstract-compactStrong}-(2), the property \eqref{lem:H3abstract-StrongCompact-c1} holds and thus also the condition \ref{H3} for the primal problem.
We claim that the same  locally uniform strong compactness property \eqref{lem:H3abstract-StrongCompact-c1} holds
for the dual problem at order $N+1$ and thus condition {\bf(H3)} holds for the dual problem. 
We may indeed use Remark~\ref{rem:H3abstract-compactStrong}-(2)  since then $\RR_{\BB^*}(\alpha)$ and $\AA^* \RR_{\BB^*}(\alpha)$ are  bounded   uniformly in $\alpha \ge \kappa_0$ and 
$$
(\AA^* \RR_{\BB^*}(\alpha))^{N+1} = \AA^* \WW(\alpha)^*  \RR_{\BB^*}(\alpha) \in \KK(Y), \quad \forall \, \alpha \ge \kappa_0,
$$
as a product of two bounded operator with a compact operator.   
\end{proof}

 \begin{rem}\label{rem:H3abstract-compactStrongAlternative} Instead of  \ref{HS1} in Lemma~\ref{lem:H3abstract-StrongC}, one can assume  
that  there exists a splitting  $\LL = \AA + \BB$ and  $N \ge 1$ such that $\BB-\alpha$ is invertible for any  $\alpha \ge \kappa_0$ and  
$$ 
\RR_\BB(\alpha) := (\alpha-\BB)^{-1}, \quad \check\VV (\alpha) :=  \sum_{i=0}^{N-1}  (\AA\RR_\BB(\alpha))^i, \quad \check\WW(\alpha) := (\AA\RR_\BB(\alpha))^N
$$ 
 are  respectively bounded in $\BBB(X)$ uniformly  with respect to $\alpha \ge \kappa_0$ and  strongly compact locally uniformly on $\alpha \ge \kappa_0$. 
 Starting indeed again from \eqref{eq:existenceStat-H3abstract-StrongC} and defining $ h_n: = (\lambda_n-\BB) \hat f_n $, we may write 
 $$
 h_n = \AA \RR_\BB(\lambda_n) h_n +  \eps_n.
$$
Observing that $ \|h_n\|_X \geq \|\RR_\BB(\lambda_n)\|_{\BBB(X)}^{-1} \geq c > 0$ by assumption, we deduce that $ \hat h_n := h_n / \|h_n\|_X$ satisfies 
$$
 \hat h_n = \check w_n + \check v_n,  \quad \check w_n := \check \WW(\lambda_n) \hat h_n ,  \quad \check v_n := \check \VV (\lambda_n) \hat \eps_n,
$$
with $\|  \hat h_n \| =1$ and $\hat \eps_n := \eps_n / \|h_n\|_X \to 0$. Similarly as in the proof of  Lemma~\ref{lem:H3abstract-StrongC}, we conclude  to the existence of subsequence $(\hat h_{n_k})$ and $h_1\in X_+\setminus\{0\}$ such that $\hat h_{n_k} \to h_1$ strongly in $X$.
Defining $f_1 := \RR_\BB(\lambda_1) h_1 / \|\RR_\BB(\lambda_1) h_1\| $,  we have again $\hat f_{n_k} \to f_1$ strongly in $X$ and next that condition \ref{H3} holds. 
  \end{rem}
   
As we see now, strong compactness is not really necessary.

\begin{lem}\label{lem:H3abstractWeakX1} We assume \ref{H1}-\ref{H2}-\ref{HS1} and there exists  $N \ge 1$ such that  
$$
  \WW(\alpha) : X \to \XX_1 \subset X \  \hbox{ is 
  uniformly bounded in  } \ \alpha \ge \kappa_0
$$
and, denoting $\XX_0 := X$,  we assume that for any $R_1 \ge R_0 > 0$ the set 
\beqn\label{eq:H3abstract-defC}
\CC = \CC_{R_0,R_1}  := \{ g \in X_+; \ \| g \|_{\XX_0} \ge R_0, \ \| g \|_{\XX_1}  \le  R_1 \}
\eeqn
is relatively sequentially compact for the weak topology $\sigma(X,Y)$ and $0 \notin \overline{\CC}$, 
where the closure is taken in the sense of the weak topology $\sigma(X,Y)$. 
{\Blue We further assume that either {\bf (1)}  $\WW(\alpha)$ is positive or {\bf (2)} there exists a constant $C_{\XX_1} > 0$ such that $\|g_+\|_{\XX_1} \le C_{\XX_1} \|g\|_{\XX_1}$ for any $g \in \XX_1$.}
Then condition \ref{H3} holds.
\end{lem}

\begin{rem}\label{rem:H3abstractWeakX1}  {\Blue When $\XX_1$ is a Banach lattice, the condition {\bf (2)} holds because of \eqref{eq:propprieteDUcone2}.}
When $\XX_1 \subset \XX_0$ with strongly compact embedding the above set $\CC$ clearly satisfies the required conditions. In particular, Lemma~\ref{lem:H3abstractWeakX1} somehow generalizes the result stated in Remark~\ref{rem:H3abstract-compactStrong}-(1).
\end{rem}

%

\begin{proof}[Proof of  Lemma~\ref{lem:H3abstractWeakX1}] 
We go back to the proof of Lemma~\ref{lem:H3abstract-StrongC} and we start with 
\eqref{eq:1hatfn=iterate}. We recall that $ \| \hat f_{n}  \|_{\XX_0} = 1$ and $ \| v_n \|_{\XX_0} \to 0$ from  \eqref{eq:lambda1approx}.  We first assume {\bf (1)} so that $w_n \ge 0$ because $\WW(\lambda_n)$ is a positive operator. 
We also observe that 
$$
 \| w_n\|_{\XX_1} \le C_\WW \| \hat f_n \|_{\XX_0} = C_\WW
$$
and
$$
 \| w_n\|_{\XX_0} \ge 1 - \| v_n \|_{\XX_0} \ge 1/2
$$
for any  $n\ge n_*$, with $n_* \ge 1$ large enough, so that $w_n \in \CC := \CC_{1/2,C_\WW}$ for any $n \ge n_*$. By the compactness properties of $\CC$, there exist a subsequence $(w_{n_k})$ and $f_1 \in X_+ \backslash \{ 0 \}$ such that 
$ w_{n_k} \wto f_1$ weakly $\sigma(X,Y)$. Since $v_n \to 0$ strongly in $X$, we deduce that $ \hat f_{n_k} \wto f_1$ weakly $\sigma(X,Y)$ and that  ends the proof of  {\bf(H3)}.

\Blue 
\smallskip
We now  assume that {\bf (2)} holds. We define $w'_n := (w_n)_+$ and $v'_n := v_n - (w_n)_-$. 
We observe that $(w_n)_- = (\hat f_n - v_n)_- \le (v_n)_+$ and thus $\| (w_n)_- \|_{\XX_0} \le \| v_n \|_{\XX_0}$ because of the compatibility of norm and order structures in the Banach space $\XX_0 = X$. 
We immediately deduce that 
$$
\hat f_n = w'_n + v'_n, \quad w'_n \ge 0,  \quad \| v'_n \|_{\XX_0} \to 0. 
$$
Arguing as above and using the additional property of  $\XX_1$, we have
$$
 \| w'_n\|_{\XX_1} \le  C_{\XX_1}  \| w_n\|_{\XX_1} \le C_{\XX_1} C_\WW
\quad\hbox{and}\quad
 \| w'_n\|_{\XX_0} \ge   1/2
$$
for any  $n\ge n_*$, with $n_* \ge 1$ large enough. We conclude the proof in the same way as previously. 
\end{proof}

\medskip
We present a typical concrete application of the preceding result.  

\begin{lem}\label{lem:H3Lp}  
We assume $X = L^p(E,\EE,\mu)$,  $p \in [1,\infty)$,  \ref{H1}-\ref{H2}-\ref{HS1} and 
\beqn\label{eq:lem=H3Lp}
  \WW(\alpha) : X \to \XX_1   \  \hbox{ is uniformly bounded in  } \ \alpha \ge \kappa_0, 
\eeqn
for some $N \ge 1$ and a  subspace $\XX_1 \subset X$ such that  $\{ g^p; \ g \ge 0, \ \| g \|_{\XX_1} \le R_1  \}$ is a weakly compact subset of $L^1(E)$ for any $R_1 > 0$. 
Then condition \ref{H3} holds.
\end{lem}

\begin{rem}\label{rem:H3Lp}  
(1) A typical example in the above statement is $\XX_1 := L^q  \cap L^p_m $  
for some exponent $q > p$ and some weight function $m : E \to [1,\infty)$ such that $m(x) \to \infty$ as $x \to \infty$. 

\smallskip
(2) The same result holds under the condition that  if $(u_n)$ is a nonnegative and bounded sequence in $L^p$ then the nonnegative sequence $w_n := \WW(\lambda_n)   u_n$ is such that $w_n^p$ is weakly compact in $L^1$. 

\smallskip
\Blue (3) We will  use Lemma~\ref{lem:H3Lp} only for the transport equation with an abstract kernel term in Section~\ref{ssec:Transport:KR}. 
\end{rem}

\begin{proof}[Proof of Lemma~\ref{lem:H3Lp}.]
For $0 < R_0 < R_1$, we define $\CC$ by \eqref{eq:H3abstract-defC} with $\XX_0 := L^p$. From 
the weak compactness property made on $\XX_1$, we observe that 
$$
\alpha(R) := \sup_{g \in \CC}  \| g {\bf 1}_{E_R^c} \|_{L^p} \to 0, \ \hbox{as} \ R \to \infty, 
$$
and 
$$
\beta(M) := \sup_{g \in \CC}  \|g {\bf 1}_{g\ge M} \|_{L^p} \to 0, \ \hbox{as} \ M \to \infty.
$$
For $g \in \CC$, we may then write 
$$
R_0 \le \| g \|_{L^p} 
\le  \| g \wedge M {\bf 1}_{E_R} \|_{L^p} + \| g {\bf 1}_{E_R^c} \|_{L^p} + \| g {\bf 1}_{g \ge M} \|_{L^p} 
$$
and thus 
$$
M^{1-1/p} \| g  {\bf 1}_{E_R} \|^{1/p}_{L^1} \ge  \| g \wedge M {\bf 1}_{E_R} \|_{L^p} \ge R_0 - \alpha (R) - \beta(M)  \ge R_0/2,
 $$
 for some $R, M > 0$ large enough. On the one hand, from the reflexivity of $L^p$ or the Dunford-Pettis theorem, the set $\CC$ is 
 relatively sequentially compact for the weak topology $\sigma(L^p,L^{p'})$. On the other hand, because ${\bf 1}_{E_R} \in L^{p'}$
 the last estimate implies that any element $g^* \in  \overline{\CC}$, 
where the closure is taken in the sense of the weak topology $\sigma(L^p,L^{p'})$, satisfies 
$$
 \langle g^*, {\bf 1}_{E_R} \rangle =  \| g^*  {\bf 1}_{E_R} \|_{L^1} \ge M^{1-p} (R_0/2)^p > 0,
 $$
 and in particular $0 \notin \overline{\CC}$.  We deduce that \ref{H3} holds as a consequence of  Lemma~\ref{lem:H3abstractWeakX1}.
 \end{proof}

We present a second kind of result where some weak compactness is involved. 

\begin{lem}\label{lem:H3abstractWeakX0} 
We assume \ref{H1}-\ref{H2}-\ref{HS1} and there exists  $N \ge 1$ such that  
\beqn\label{eq:H3abstractWeakX0} 
  \WW(\alpha) : \XX_{0} \to X \subset \XX_0  \  \hbox{ is  uniformly bounded in  } \ \alpha \ge \kappa_0
\eeqn
and, denoting $\XX_1 := X$, the set $\CC$ defined by \eqref{eq:H3abstract-defC} satisfies the same properties as the ones stated in Lemma~\ref{lem:H3abstractWeakX1}.
Then condition \ref{H3} holds.
\end{lem}

 \begin{rem}\label{rem:H3Lp&H3XX1}  
If we replace the norm $\| \cdot \|_{\XX_0}$ by a seminorm $\| f \|_{\XX_0} := \langle |f|, \varphi_0 \rangle$, $\varphi_0 \in Y_+$,   we define $\CC$ accordingly by \eqref{eq:H3abstract-defC} and  we
 assume that $X=Y'$ with $Y$ separable, then $\CC$ satisfies the same compactness properties  as required in the statement of Lemma~\ref{lem:H3abstractWeakX1}. If we further assume that \eqref{eq:H3abstractWeakX0} holds, where $\XX_0$ is endowed with the above seminorm, we may repeat the proof below in order to obtain that \ref{H3} holds in that situation (see also Lemma~\ref{lem:H3Y'}
 and its proof for a slightly generalized situation).  
\end{rem}

\begin{proof}[Proof of Lemma~\ref{lem:H3abstractWeakX0}.]
We start here again with  \eqref{eq:1hatfn=iterate}. 
We have
$$
1=  \| \hat f_n \|_{\XX_1} \le C_\WW \| \hat f_n \|_{\XX_0} + \| v_n \|_{\XX_1}, 
$$
and thus 
$$
 \|  \hat f_n \|_{\XX_0} \ge C_\WW^{-1} (1 - \| v_n \|_{\XX_1}) \ge (2C_\WW)^{-1}  
$$
for any  $n\ge n_*$, with $n_* \ge 1$ large enough, so that $\hat f_n  \in  \CC := \CC_{ (2C_\WW)^{-1}  , 1}$, for $n\ge n_*$.
By the compactness properties of $\CC$, there exist a subsequence  $(\hat f_{n_k})$ and $f_1 \in X_+ \backslash \{ 0 \}$ such that 
$\hat f_{n_k} \wto f_1$ weakly $\sigma(X,Y)$.  
\end{proof}

We present  a variant of Lemma~\ref{lem:H3Lp} which is also a concrete consequence of Lemma~\ref{lem:H3abstractWeakX1}
and  Lemma~\ref{lem:H3abstractWeakX0}.

\begin{cor}\label{cor:H3Lp} 
We assume \ref{H1}-\ref{H2}-\ref{HS1} in $X = L^{p_0}_{m_0}$, $1 \le p_0 < \infty$, 
together with the facts that  $\RR_\BB(\alpha)$ is   bounded in $\BBB(L^{p_0}_{m_0})$ and $\BBB(L^{p_1}_{m_1})$  uniformly in $\alpha \ge \kappa_0$, 
$\AA \in \BBB(L^{p_0}_{m_0})$  and $(\RR_\BB(\alpha) \AA)^{N}$ is   bounded in $\BBB(L^{p_0}_{m_0},L^{p_1}_{m_1})$  uniformly in $\alpha \ge \kappa_0$ for some $N \ge 1$, with $p_1 > p_0$ and $m_1$ such that $m_0/m_1 \in L^\vartheta$, $1/\vartheta := 1/p_0 - 1/p_1$. Then condition \ref{H3} holds for both the primal and the dual problems.
\end{cor}

\begin{proof}[Proof of Corollary~\ref{cor:H3Lp}.]
On the one hand, we have 
\bean
&& \RR_\BB(\alpha)  + \dots + (\RR_\BB(\alpha) \AA)^{N-1} \RR_\BB(\alpha)  \hbox{ is bounded in } \BBB(X) \hbox{ uniformly in  } \alpha \ge \kappa_0, 
\\ 
&& \WW(\alpha) := (\RR_\BB(\alpha) \AA)^{N}  \hbox{ is bounded in } \BBB(X,\XX_1) \hbox{ uniformly in  } \alpha \ge \kappa_0, 
\eean
with $\XX_1 := L^{p_1}_{m_1} \subset X$ and thus   $\{ (gm_0)^{p_0}; \ g \ge 0, \ \| g \|_{\XX_1} \le R_1  \}$ is a weakly compact subset of $L^1(E)$ for any $R_1 > 0$.
Condition \ref{H3} holds for the direct problem thanks to  Lemma~\ref{lem:H3Lp}.

On the other hand, we set $Y := X' = L^{q_0}_{\nu_0}$, $q_0 := p'_0$, $\nu_0 := m^{-1}_0$, and we first observe that 
$$
\RR_{\BB^*}(\alpha)  + \dots + (\RR_{\BB^*}(\alpha) \AA^*)^{N-1} \RR_{\BB^*}(\alpha)  \hbox{ is bounded in } \BBB(Y) \hbox{ uniformly in  } \alpha \ge \kappa_0. 
$$
We next observe that 
$$
(\AA^* \RR_{\BB^*}(\alpha))^{N+1} = \AA^* \WW(\alpha)^*  \RR_{\BB^*}(\alpha) 
 \hbox{ is bounded in } \BBB(\YY_0,Y) \hbox{ uniformly in  } \alpha \ge \kappa_0, 
$$
with  $\YY_0 := L^{q_1}_{\nu_1}$, $q_1 := p'_1$, $\nu_1 := m^{-1}_1$. Because  $\{ (g\nu_1)^{q_1}; \ g \ge 0, \ \| g \|_{Y} \le R_1  \}$ is a weakly compact subset of $L^1(E)$ for any $R_1 > 0$,
we have from the proof of  Lemma~\ref{lem:H3Lp} that the set $\CC$ defined by \eqref{eq:H3abstract-defC} for the norms of $\YY_0$ and $\YY_1:=Y$ satisfies the weak compactness property 
required in the statement of Lemma~\ref{lem:H3abstractWeakX1}. We may thus apply Lemma~\ref{lem:H3abstractWeakX0} and we deduce that condition \ref{H3} holds for the dual problem.
\end{proof}

Another concrete consequence of Lemma~\ref{lem:H3abstractWeakX1} and Lemma~\ref{lem:H3abstractWeakX0} is the following. 

\begin{lem}\label{lem:H3M1} We assume $E$ not compact, $X = M^1_{m_i}(E)$ for a continuous weight function $m_i$ on $E$, $i=0$ or $i=1$,  \ref{H1}-\ref{H2}-\ref{HS1} and there exists $N \ge 1$ such that 
$ (\RR_\BB(\alpha) \AA)^N : M^1_{m_0}(E)  \to M^1_{m_1}(E)$ uniformly in $\alpha \ge \kappa_0$ for another continuous weight function $m_{1-i}$ on $E$ such that $m_1(x)/m_0(x) \to \infty$ as $x \to \infty$. 
Then condition \ref{H3} holds.
\end{lem}

\begin{proof}[Proof of Lemma~\ref{lem:H3M1}.]
We define $\XX_i := M^1_{m_i}(E)$ and we consider the set 
$\CC$ defined by \eqref{eq:H3abstract-defC} which is clearly compact for the weak $*$ $\sigma(M_{m_1}^1, C_{m^{-1}_1,0})$ topology.
When  $X=M^1_{m_0}$, the result follows from  Lemma~\ref{lem:H3abstractWeakX1} while when $X=M^1_{m_1}$, the result is a  consequence of Lemma~\ref{lem:H3abstractWeakX0}. 
\end{proof}

We may slightly improve the preceding results by considering a more abstract framework and a somehow more general boundedness condition. 

\Black

\begin{lem}\label{lem:H3Y'} We assume $X = Y'$, $Y$ separable,  \ref{H1}-\ref{H2}-\ref{HS1} and there exist $N \ge 1$, $\gamma\in[0,1)$ and $\varphi\in Y_+\setminus\{0\}$ such that for any $\alpha \ge \kappa_0$, 
there holds
\begin{equation}\label{as:weakly-compact-1}
\| \WW(\alpha) f\|_X\leq \gamma\|f\|_X+\langle f,\varphi\rangle_{X,Y},
\end{equation}
for all $f\in X_+$, or there holds
\begin{equation}\label{as:weakly-compact-2}
\| \WW(\alpha)  f\|_X\leq \gamma\| f \|_X+\langle \WW(\alpha) f,\varphi\rangle_{X,Y},
\end{equation}
for all  $f \in X_+$.
Then condition \ref{H3} holds true, and the limit $f_1$ satisfies $\langle f_1,\varphi\rangle_{X,Y}\geq1-\gamma > 0$.
\end{lem}

The case $X=M^1_{m_1}(E)$  in Lemma~\ref{lem:H3M1} corresponds here to the first situation where \eqref{as:weakly-compact-1} holds with $X:=M^1_{m_1}(E)$, $\gamma:=0$, $Y := C_{m^{-1}_0,0}(E)$
and {\Cyan $\varphi :={  \big(\sup_{\alpha\geq\kappa_0}\|\WW(\alpha)\|_{M^1_{m_0}\to M^1_{m_1}}\big)} m_0$.
\Blue We will  use Lemma~\ref{lem:H3Y'} only for the  mutation-selection model with singular mutation kernel  
 in Section~\ref{ssec:appl5:singular}.}

\smallskip

\begin{proof}[Proof of Lemma~\ref{lem:H3Y'}]
Starting  with  \eqref{eq:1hatfn=iterate} and using   
\eqref{as:weakly-compact-1},  we have
\bean
\| \hat f_{n}\|_{X} 
&\le&   \|  \WW(\lambda_n)\hat f_n \|_{X} +  \|\VV(\lambda_n)\varepsilon_n  \|_{X} 
\\
&\le&   \gamma\| \hat f_n \|_X+\langle \hat f_n,\varphi\rangle_{X,Y} +  \| v_n  \|_{X} , 
\eean
so that 
$$
\langle \hat f_n,\varphi\rangle_{X,Y} \ge    1- \gamma - \|v_n  \|_{X} .
$$
By compactness, there are $f_1\geq0$ and a subsequence $(\hat f_{n'})$ such that $\hat f_{n'}\rightharpoonup f_1$ weak $*$ $\sigma(X,Y)$.
Passing to the limit as $n' \to \infty$ in the above estimate, we find 
\beqn\label{as:weakly-compact-3}
\langle f_1 ,\varphi\rangle_{X,Y}  = \lim_{n' \to \infty} \langle \hat f_{n'},\varphi\rangle_{X,Y} \ge 1- \gamma , 
\eeqn
and in particular  $f_1\neq0$. 

\smallskip
Under the assumption  \eqref{as:weakly-compact-2}, modifying  slightly the previous argument,  we have
\bean
\| \hat f_{n}\|_{X} 
\le   \gamma\| \hat f_n \|_X+ \langle w_n ,\varphi\rangle_{X,Y} +  \|  v_n  \|_{X} , 
\eean
which, together with 
\[
\langle \hat f_n,\varphi \rangle_{X,Y} =\langle w_n,\varphi \rangle_{X,Y} +\langle v_n,\varphi \rangle_{X,Y} ,
\]
implies
$$
\langle \hat f_n,\varphi\rangle_{X,Y} \ge 1- \gamma - \| v_n  \|_{X} +\langle v_n,\varphi \rangle_{X,Y} . 
$$
By compactness again, there are $f_1\geq0$ and a subsequence $(\hat f_{n'})$ such that $\hat f_{n'}\rightharpoonup f_1$ weak $*$ $\sigma(X,Y)$, and passing to the
limit $n' \to \infty$ in the above estimate, we conclude again to \eqref{as:weakly-compact-3}. 
\end{proof}

 
Let us comment on Lemma~\ref{lem:H3Y'} and in particular the condition~\eqref{as:weakly-compact-2}.
 
 \smallskip
In the case when $X=L^\infty(E,\EE,\mu)=\big(L^1(E,\EE,\mu)\big)'$, 
we  can relate condition~\eqref{as:weakly-compact-2} to the assumption that there exist $f_0\in X_+$ and $\varphi \in Y_+\setminus\{0\}$ such that 
\begin{equation}\label{as:A4Harris}
\|S_\LL(t)f_0\|_X\leq \langle S_\LL(t)f_0,\varphi\rangle, \quad \forall \, t \ge 0.
\end{equation}
This last condition is reminiscent from conditions that appear in probabilistic inspired methods for the ergodicity of semigroups,
see the condition {\it(1b)} in Theorem~\ref{theo:BCGM-Harris} but also Assumption~(A2) in~\cite{MR3449390}, both in the vein of~\cite[Condition~$\mathcal Z$]{MR1988460}.
Assume indeed~\eqref{as:A4Harris}, let $\eta>\kappa_1-\kappa_0>0$ and consider the trivial decomposition $\LL=\AA+\BB=\eta+(\LL-\eta)$.
Then set $\kappa_\BB :=\kappa_1-\eta<\kappa_0$, so that for any $\alpha>\kappa_\BB $, $\BB-\alpha=\LL-(\eta+\alpha)$ is invertible since $\eta+\alpha>\eta+\kappa_\BB =\kappa_1$.
We thus have for any $\alpha>\kappa_\BB $
\[
\WW(\alpha):=\eta(\alpha-\BB)^{-1}=\eta\int_0^\infty e^{-(\eta+\alpha)t}S_\LL(t)\,dt
\]
and~\eqref{as:A4Harris} then ensures that
\[\|\WW(\alpha)f_0\|_X\leq \langle \WW(\alpha)f_0,\varphi\rangle.\]
We recover~\eqref{as:weakly-compact-2} with $\gamma = 0$ and the difference that $f_0$ is fixed here. 
 
\medskip
As a Corollary of Lemma~\ref{lem:H3M1} or Lemma~\ref{lem:H3Y'} and anticipating on the material of part~\ref{sec:DynamicalExistenceKR}, we present now a situation very classical in stochastic processes theory.

\Black
\begin{cor}\label{cor:H3M1} We consider a positive semigroup $S = S_\LL$ defined on a Radon space $ X = M^1_\psi(E)$ for some  positive weight functions $\psi$ on $E$, in  particular \ref{H1} holds. We also assume that \ref{H2} holds for some $\kappa_0 \in \R$. We finally assume the Lyapunov condition 
\beqn\label{eq:LyapunovH3M1}
\LL^* \psi \le \kappa_\BB \psi  + M \chi, 
\eeqn
with $ \kappa_\BB < \kappa_0$, $M \ge 0$ and   $\chi \in C_{\psi,0}(E)$, $0 \le \chi \le \psi$. Then condition \ref{H3} holds true.
\end{cor}

Let us emphasize that we may assume some regularity on $\psi$ by considering $\psi \in D(\LL^*)$ so that \eqref{eq:LyapunovH3M1} makes sense in $X$ or just understand \eqref{eq:LyapunovH3M1} in the weak sense: 
$$
\langle \LL f,  \psi \rangle \le \kappa_\BB \langle f,\psi \rangle + M  \langle f, \chi \rangle, \quad \forall \, f \in D(\LL) \cap X_+.
$$

\begin{proof}[Proof of Corollary~\ref{cor:H3M1}.]
We introduce the splitting $\LL = \AA + \BB$ where $\AA$ is the bounded multiplicator  operator $\AA := M\chi/\psi$. As a bounded perturbation of $\LL$, the operator $\BB$
is the generator of a semigroup $S_\BB$.
Defining $\widetilde S_t  := S_\LL(t) e^{-Mt} \ge 0$ and $\AA^c := M(1-\chi/\psi) \ge 0$, we have the Duhamel formula 
$$
S_\BB = \widetilde S  +  \widetilde S \AA^c * S_\BB
$$
and iterating infinitely many times, we deduce the Dyson-Philips formula 
$$
S_\BB = \sum_{k=0}^\infty ( \widetilde S \AA^c)^{(*k)} *  \widetilde S. 
$$
That implies that $S_\BB \ge 0$ as a combination of positive operators. 
Alternatively,   from the very definition of $\BB$, we have $\kappa-\BB\leq(M+\kappa)-\LL$ for  any $\kappa \in \R$. 
Choosing $\kappa>\max(\omega(S_\LL),\omega(S_\BB))$ and using the direct implication in Lemma~\ref{lem:sec2-EquivPositiveS}, 
we have  $\RR_\BB(\kappa)\geq\RR_\LL(M+\kappa)\geq0$. 
Using the reciprocal implication in Lemma~\ref{lem:sec2-EquivPositiveS},  we obtain again that $S_\BB\geq0$.
  
Now, for $0 \le f_0 \in D(\BB)$ and setting $f_t := S_\BB(t) f_0$, we may compute 
$$
\frac{d}{dt} \langle f_t, \psi \rangle =  \langle \BB f_t, \psi \rangle  \le  \kappa_\BB \langle f_t, \psi \rangle, 
$$
so that 
$$
\| S_\BB(t) f_0 \|_{M^1_{\psi}} \le e^{ \kappa_\BB t} \|  f_0 \|_{M^1_{\psi}}.
$$
Using \eqref{eq:Exist1-DefRepresentationRR}, we immediately and classically deduce 
$$
\| \RR_\BB (\alpha) \|_{\BBB(M^1_{\psi})} \le \frac{1}{\alpha -  \kappa_\BB}, \quad \forall \, \alpha > \kappa_\BB, 
$$
so that $\RR_\BB (\alpha)$ is bounded in $\BBB(M^1_\psi)$ and $\RR_\BB (\alpha) \AA$ is bounded in $\BBB(M^1_\chi,M^1_{\psi})$ uniformly for $\alpha \ge \kappa_0$.
We apply Lemma~\ref{lem:H3M1} or Lemma~\ref{lem:H3Y'} (\eqref{as:weakly-compact-1} with $N=1$, $\gamma = 0$ and $\varphi =   \frac{M}{\alpha-\kappa_\BB} \chi$)  in order to conclude.
\end{proof}

In the proof of Corollary~\ref{cor:H3M1}, we may alternatively use the trivial splitting $\LL=\tilde\AA+\tilde\BB= \eta +(\LL-\eta)$
for some $\eta>\kappa_1-\kappa_0$, so that $\alpha-\tilde\BB$ is invertible for any $\alpha\geq\kappa_0$,  and reformulate the Lyapunov condition 
\[(\alpha-\tilde\BB^*)\psi\geq (\alpha+\eta-\kappa_\BB ) \psi -M\chi,
\]
for any $\alpha \ge \kappa_0$. Observing that $\widetilde\WW(\alpha):= \tilde\AA \RR_{\tilde\BB}(\alpha)  = \eta (\alpha-\tilde\BB)^{-1}$, we deduce 
\[
\widetilde\WW^*(\alpha) \psi  \leq \frac{\eta}{\eta+\alpha-\kappa_\BB} \psi +\frac{M}{\eta+\alpha-\kappa_\BB}\widetilde\WW^*(\alpha)\chi.
\]
We equivalently have
\[\
\| \widetilde \WW(\alpha) f \|_{M^1_\psi}\leq\gamma\| f \|_{M^1_\psi}+\langle \widetilde\WW(\alpha) f ,\varphi \rangle,
\]
uniformly for any $\alpha \ge \kappa_0$,  with  $\gamma:=\frac{\eta}{\eta+\kappa_0 -\kappa_\BB }<1$ and $\varphi :=\frac{M}{r+\kappa_0 -\kappa_\BB }\chi$, 
which is nothing but condition~\eqref{as:weakly-compact-2}. 
We conclude again thanks to  Lemma~\ref{lem:H3Y'}.

\medskip
{\Blue We emphasize again that in almost all our applications, we will be able to use  Lemma~\ref{lem:H3abstract-StrongC}-(2) for establishing \ref{H3}, except in very irregular situations (we will use Lemma~\ref{lem:H3Lp} in Section~\ref{ssec:Transport:KR}  and  Lemma~\ref{lem:H3Y'} in Section~\ref{ssec:appl5:singular}) and in weak dissipative frameworks (we will directly show \ref{H3} by following Lions' proof in Section~\ref{subsection:diffusionWeakPotential} and we will use Theorem~\ref{theo:KRexistTER} in Section~\ref{subsec-diffusionwith drift}).
In Section~\ref{ssec:GF:variability} we do not use exactly Lemma~\ref{lem:H3abstract-StrongC} neither, but we adapt the proof for combining two different splittings.}

\medskip
We finally come to  the existence of a solution to the first eigenvalue problem and the first eigentriplet problem.

\begin{theo}\label{theo:exist1-KRexistence}
Under conditions \ref{H1}-\ref{H2}-\ref{H3}, the first eigenvalue problem \eqref{eq:1stEVP} has a  solution $(\lambda_1,f_1)$ with $\lambda_1$ satisfying \eqref{eq:lambda1inab}. 
When furthermore \ref{H3} holds for the dual problem, then 
the first eigentriplet problem \eqref{eq:triplet1}-\eqref{eq:triplet2} admits a solution
$(\lambda_1, f_1, \phi_1) \in \R \times X \times Y$.
\end{theo} 
 
Theorem~\ref{theo:exist1-KRexistence} generalizes some known versions of the existence part of the Krein-Rutman Theorem
where either $\LL$ is assumed additionally to be the generator of a semigroup or to have strongly power compact resolvent or even where some additional conditions are made on the positive cone $X_+$. 
As mentioned in the introduction, some possible references for these previous results are 
Krein-Rutman \cite{MR0027128}, Greiner in \cite[Cor~1.2]{MR763356} and in \cite[C-IV Thm.~2.1]{MR839450}  
and Webb \cite[Prop.~2.5]{MR902796}, 
see also \cite[Thm.~2]{MR923493}, 
\cite[Thm.~5.3]{MR3489637}, 
 \cite{PL2}, 
\cite[Thm.~2.1]{Bansaye2022}, 
the textbook  \cite[Thm.~12.15]{MR3616245}  and the references therein.

\begin{proof}[Proof of Theorem~\ref{theo:exist1-KRexistence}.]
We first assume \ref{H1}-\ref{H2}-\ref{H3}. 
Because of  Lemma~\ref{lem:Existe1-Spectre2}, there exists a sequence $(\hat f_n)$ of $X$ such that \eqref{eq:lambda1approx} holds, and in particular 
\beqn\label{eq:lambda1approxBIS}
 \langle \lambda_n \hat f_n,\phi \rangle  - \langle  \hat f_n, \LL^* \phi \rangle = \langle   \eps_n, \phi \rangle, 
\eeqn
for any $\phi \in D(\LL^*)$.  Because of condition \ref{H3}, we may pass to the limit $n' \to \infty$ in equation  \eqref{eq:lambda1approxBIS} and we deduce that $(\lambda_1,f_1)$ satisfies \eqref{eq:1stEVP}  and \eqref{eq:lambda1inab}. 

\smallskip
We now additionally assume that  \ref{H3} holds for the dual problem. 
As recalled during the proof  of Lemma~\ref{lem:Existe1-H2bis} and by definition of $\lambda_1$, we have  $(\lambda_1,+\infty) \subset \rho(\LL) = \rho(\LL^*)$ and $\lambda_1 \in \Sigma(\LL) = \Sigma(\LL^*)$. Taking $\lambda_n \searrow \lambda_1$, we argue as in the proof of Lemma~\ref{lem:Existe1-Spectre2} and we get
$$ \ \exists \, \widehat \phi_n \ge 0,\  \lambda_n \widehat \phi_n - \LL^* \widehat \phi_n  \to 0, \ \| \widehat \phi_n \| = 1.
$$
Thanks to  {\bf (H3)}  for the dual problem, we deduce that there exist a subsequence $(\widehat \phi_{n_k})$ and $\phi_1 \in X'$, $\| \phi_1 \| = 1$ such that $\widehat \phi_{n_k} \to \phi_1$. We thus conclude that $\phi_1$ is a solution to the dual problem  \eqref{eq:triplet2} (for the same eigenvalue $\lambda_1$). 
\end{proof}

Let us conclude this section by some remarks. 

\begin{rem}\label{rem:existenceMM-dualPb} 

(1)  - As seen above, the condition  \ref{H1}-\ref{H2} for the primal and the dual problems are equivalent, and thus one only has to check \ref{H1}-\ref{H2}-\ref{H3} for the primal problem and \ref{H3} for the dual problem in order to solve the first eigentriplet problem. 
It is worth emphasizing that condition {\bf (H3)} on the dual problem is not a consequence of the condition \ref{H3} on the primal problem. However, as presented in Lemma~\ref{lem:H3abstract-StrongCbis},  Lemma~\ref{lem:H3abstract-StrongC} and Corollary~\ref{cor:H3Lp}, there exist 
several natural situations where both  conditions \ref{H3} for the primal and the dual problems hold together. 
 
\smallskip
(2) - Alternatively, one may also assume 
 \ref{H1}-\ref{H2}-\ref{H3} for the dual problem, and then use a more classical fixed point theorem for proving the existence of a steady state for the rescaled semigroup by using 
 for instance the Markov--Kakutani fixed point theorem \cite{MR1568507} as in \cite[Thm.~5.1]{MR4265692}, the Tychonov fixed point theorem as in \cite{MR2053942} or \cite[Thm.~1.2]{MR2114413} 
 or a Birkhoff-Von Neumann type Theorem as in  \cite[Thm.~6.1]{MR4534707}. For these last techniques, we also refer to Section~\ref{sec:DynamicalExistenceKR}, where such a dynamical approach is adapted to the present context.  One may also use the usual Doblin-Harris theory, see for instance \cite{MR2857021,MR4534707} and the references therein, and  Sections~\ref{ssec:Transport:renewal} and~\ref{ssec:appl5:singular} for applications of this approach. 
\end{rem}

{\Blue 
Let us finally emphasize that Theorem~\ref{theo:exist1-KRexistence}  provides an answer to problem \ref{S1} about the existence of a first eigentriplet. 
 In particular, for a  positive semigroup $S_\LL$ on $X$ which satisfies \ref{H1}, \ref{H2} and which enjoys the splitting structure \eqref{eq:intro:RLRB-RLARBN}, \eqref{eq:intro-regRL}, \eqref{eq:intro:stopDyson}, \eqref{eq:intro-regSL} and thus \ref{HS1}, we may use  Theorem~\ref{theo:exist1-KRexistence} and Remark~\ref{rem:existenceMM-dualPb} for establishing Theorem~\ref{theo:main-intro}-(1), at least in three situations: (1)~when $\XX_1 \subset \XX_0 := X$ is compact (see Lemma~\ref{lem:H3abstract-StrongC} and Remark~\ref{rem:H3abstract-compactStrong}), (2) when the conditions of   Lemma~\ref{lem:H3abstractWeakX1},     Lemma~\ref{lem:H3abstractWeakX0} or  Lemma~\ref{lem:H3Y'} are complete for both the primal and the dual problem
 and (3) when the assumptions of  Corollary~\ref{cor:H3Lp} are met. 
 }

\subsection{Discussion}
\label{subsect:Exist1-discussion}
\

\smallskip

We discuss now the existence results presented in the preceding section.

\medskip\quad
For further references, let us first recall that  when $X$ is a Hilbert space and $\LL$ is self-adjoint, the first eigenvalue may be simply obtained thanks to the variational problem
\beqn\label{eq:lambda1=varSA}
\lambda_1 = \sup_{f \in X_+ \backslash \{0 \}} \frac{\langle \LL f, f \rangle}{\| f \|^2} .
\eeqn

\smallskip\quad
We now explain how Theorem~\ref{theo:exist1-KRexistence} is a generalization of the classical Krein-Rutman theorem stated in Theorem \ref{theo:intro-KRresult}. 
We thus consider a Banach lattice $X$ such that $X_{++} := \hbox{\rm int}X_+ \not=\emptyset$ and  an operator $\LL$ such that,  for $\kappa_1 \in \R$ and any $\kappa > \kappa_1$, $\RR := (\kappa-\LL)^{-1} : X \to X$ is compact and $\RR : X_+\backslash\{0\} \to X_{++}$, in particular {\bf (H1)} holds true. 
As a first step, we  recall the  following very classical technical lemma of the Krein-Rutman theory. 
 
\begin{lem}
\label{lem:exists-lemfonda2} Assume $X_{++} := \hbox{\rm int} \, X_+\not=\emptyset$. 
For $g \in X_+$ and  $f \in X_{++}$, 
 there exists $C \ge 0$ such that $g \le C f$.
\end{lem}

\begin{proof}[Proof of Lemma~\ref{lem:exists-lemfonda2}.]
We argue by contradiction. Otherwise, for any $n \ge 1$,
we would have $f-g/n \in X_+^c \subset X_{++}^c$ 
 and that last set is closed. Passing to the limit, we get $f \in X_{++}^c$, 
which is in contradiction with the assumption $f \in X_{++}$. 
\end{proof}

For a given $g_0 \in X_+\backslash\{0\}$, we set $f_0 := \RR g_0 \in X_{++}$. From Lemma~\ref{lem:exists-lemfonda2}, there exists $C_0 \ge 0$ such that 
$ (\kappa-\LL) f_0 = g_0 \le C_0f_0$. That implies that Lemma~\ref{lem:Existe1-Spectre2bis}-{\bf (ii)} holds with $\kappa_0 := \kappa-C_0$, and thus \ref{H2}  also holds. One may then define $\mu_1 := \kappa-\lambda_1$, 
with 
$$
 \lambda_1 := \inf \{ \lambda \in \R; \, (\lambda'-\LL)^{-1} \in \BBB(X), \ \forall \, \lambda' \in [\lambda,\kappa] \} \ge \kappa_0.
$$ 

\smallskip
We recall that because of Lemma~\ref{lem:Existe1-Spectre2} (or its proof), there exist $(\lambda_n)$, $(\hat f_n)$ and $(\eps_n)$ such that \eqref{eq:lambda1approx} holds, namely 
$$
 \lambda_n \searrow \lambda_1, \  \hat f_n \ge 0,\  \eps_n := \lambda_n \hat f_n - \LL \hat f_n \ge 0, \ \| \hat f_n \| = 1, \ \| \eps_n \| \to 0.
$$
We may rewrite the equation as 
$$
\hat f_n = \RR [\eps_n + (\kappa-\lambda_n)\hat f_n], 
$$
so that $(\hat f_n)$ belongs to a compact set of $X$ because of the compactness assumption made on $\RR$, so that  \ref{H3}  holds true. 

\smallskip
Because of Theorem~\ref{theo:exist1-KRexistence}, we deduce that there exists $f_1 \in X_+$ such that $\| f_1 \|=1$ and $\LL f_1 = \lambda_1$. 
That implies $f_1 = \mu_1 \RR f_1$, and thus that the existence part of Theorem~\ref{theo:intro-KRresult} is a consequence of Theorem~\ref{theo:exist1-KRexistence} for an operator $\RR$ which is the positive resolvent of an operator $\LL$.

\medskip\quad
We would like to emphasize on the fact that our definition of the first eigenvalue by \eqref{eq:exist1-defI}-\eqref{eq:exist1-deflambda1} 
bears some strong similarity with the definition of the first eigenvalue  for elliptic operators in non divergence form as presented in  \cite{MR1258192}. 
Indeed, if $\lambda \in \II$, then 
$$
\exists \,  f \in X_+ \backslash \{ 0 \}, \quad \LL  f \le \lambda f.
$$
Assuming now that $X$ is a space of functions (defined on a set $E$) and that $f (x) > 0 $ for any $x \in E$, we deduce that 
$$
\lambda \ge \sup_E \frac{ \LL  f}{f}, 
$$
and thus $\lambda_1$ is characterized by 
$$
\lambda_1 = \inf_{f > 0} \sup_E \frac{ \LL  f  }{ f},
$$
which is nothing but \cite[(1.13)]{MR1258192} (with a change of sign because of a different sign convention). We thus see that our formulation
is a generalization at a more abstract level and for  resolvent positive operators of that classical $\min$-$\max$ approach for elliptic operators. 
Some more or less classical references on that subject are \cite{MR361998,MR425380},  \cite{MR1226957}, \cite{MR924555},  \cite{MR3340379} and \cite{Berestycki2016}.
In particular in  \cite{MR3340379}, two generalized principal eigenvalues 
$$
\lambda'_1 := \sup \{ \kappa_0 \in \R; \, \exists \, g_0 \in \CC_0 \   \LL g_0 \ge \kappa_0 g_0 \}
$$
and 
$$
\lambda_1'' := \inf \{ \kappa_1 \in \R; \, \exists \, g_1 \in \CC_1 \    \LL g_1 \le \kappa_1  g_1 \}
$$
are defined for appropriate cones $\CC_i \subset X_+ \backslash \{0 \}$ for problems with lack of compactness. The links between the three quantities $\lambda_1$,  $\lambda'_1$ and
 $\lambda''_1$ are discussed as well as the possible non existence of a related principal eigenfunction $f_1$. The non existence of associated principal eigenfunction 
 should not be a surprise since it is the case when one considers $\LL = \Delta$ in $X= L^2(\R^d)$ where $\LL \psi = \LL^* \psi =\lambda''_1 \psi$ with $0 < \psi = 1 \notin X = X'$ and $\lambda''_1 = 0$,
  but no associated principal eigenfunction exists in $X$. We also refer to \cite{PL2} where some examples of such a situation are discussed.

 \Black


\medskip\quad
For its own interest and further discussions, we finally state and prove a slightly less general variant of 
Theorem~\ref{theo:PLL}.

\begin{theo}\label{theo:PLLbis}
Consider a Banach lattice $X$ and a linear and bounded operator $\RR : X \to X$ such that 

\smallskip
(i) $\RR : X_+ \to X_+$;

\smallskip
(ii)  $\exists \, g_2 \in X_+ \backslash\{0\}$, $\exists \, C_2 > 0$ such that 
$\RR g_2 \le C_2 g_2$.

\smallskip
We define 
$$
K_2 := \{ g \in X_+; \, \exists a > 0, \, g \le a g_2 \}, 
$$
and next 
$$
A(g) := \inf \{ a > 0; \, g \le a g_2 \}, \ \hbox{ if } \ g \in K_2,
$$
as well as 
$$
\JJ := \{ \mu \ge 0; \ \exists h \in K_2, \ h \ge \mu \RR h + g_2 \}.  
$$

\smallskip
We further assume

\smallskip
(iii) $\mu_1 := \sup \JJ < +\infty$.

\smallskip
(iv) Any upper bounded  and increasing sequences $(g^n)$ of $X$ is convergent in the weak sense $\sigma(X,Y)$. 
More precisely, if $g_n \le g_{n+1} \le \bar g \in X$ for any $n \ge 1$, there exists $g \in X$, $g \le \bar g$, 
such that $g_n \wto g$.

\smallskip
(v) Any sequence $(g^n)$ of normalized almost first eigenvectors is relatively  compact. More precisely, for any sequence  $(g^n)$  of $K_2$ such that $A(g^n) = 1$,  $g^n = \mu^n \RR g^n+ \eps^n$, $\mu^n \nearrow \mu_1$ and $\eps^n \to 0$, there exists $g \in K_2$ and a subsequence $(g^{n_k})$ such that $g^{n_k} \to g$ and $A(g) = 1$. 

\smallskip
Then there exists $f_1 \in X_2$ such that   $f_1 = \mu_1 \RR f_1$ and $A(f_1) = 1$. 
\end{theo}

\begin{proof}[Proof of Theorem~\ref{theo:PLLbis}.]
We split the proof into three steps. 

\smallskip
{\sl Step 1. } 
We first establish that for any $\mu \in \JJ$, there exists $\tilde g = \tilde g_\mu \in K_2$ such that 
\beqn\label{eq:PL2:pbmuexist}
\tilde g = \mu \RR \tilde g+ g_2. 
\eeqn
We set $\tilde g_0 = 0$, $\tilde g_1 = g_2$, and we define $(\tilde g_n)$ recursively by $\tilde g_{n+1} = \mu \RR \tilde g_n + g_2$, for any $n\ge1$. 
We claim that 
$$
0 \le \tilde g_n \le \tilde g_{n+1} \le h, \ \hbox{ for any } \ n \ge 0,
$$
where $h$ enters in the definition of $\mu \in \JJ$. That is obviously true at order $n=0$. Assuming that last inequality 
is proved at order $n -1$ for $n\ge1$, we compute 
$$
\tilde g_{n+1} = \mu \RR \tilde g_n + g_2 \ge \mu \RR\tilde g_{n-1} + g_2 = \tilde g_n
$$
and 
$$
\tilde g_{n+1} = \mu \RR\tilde g_n + g_2 \le \mu \RR h + g_2 \le h,
$$
which proves the same inequality at order $n$, and thus for any $n \ge 0$. Using the convergence property (iv) of upper bounded increasing sequences in $X$,
we deduce that there exists $\tilde g \in X_2$ such that $\tilde g_n \to \tilde g$ and thus \eqref{eq:PL2:pbmuexist} holds. 

\smallskip
{\sl Step 2. } We obviously have $0 \in \JJ$ and $\JJ$ is an interval because if $(\mu,h)$ satisfies the condition $\mu \in \JJ$ then so do $(\mu',h)$ for any $\mu' \in [0,\mu]$. 
We finally claim that $\JJ$ is open. Take indeed $\mu \in \JJ$ and $\tilde g \in K_2$ such that \eqref{eq:PL2:pbmuexist} holds, what is possible due to Step~1. By definition, there would
exist $a > 0$ such that $\tilde g \le a g_2$. Choosing $0 < \eps \le 1/(2aC_2)$ and $M \ge 2$, we compute 
\bean
(M\tilde g) - (\mu +\eps)\RR (M \tilde g) 
&=& M g_2 - M\eps \RR\tilde g  
\\
&\ge& M g_2 - M\eps a  \RR g_2 \ge M (1 - \eps aC_2)  g_2 \ge g_2,
\eean
so that $\mu+\eps \in \JJ$. 

\smallskip
{\sl Step 3. } 
We first establish by contradiction that $A(\tilde g_\mu) \nearrow \infty$ when $\mu \nearrow \mu_1$. If it was not the case, there exists $a \in (0,\infty)$ and a sequence $(\mu^n)$ such that 
$A(\tilde g_{\mu^n}) \le a$ as $\mu^n \nearrow \mu_1$.  
Choosing $0 < \eps \le 1/(2aC_2)$ and $M \ge 2$ as in Step 2, the same computation   gives 
\bean
(M\tilde g_{\mu^n}) - (\mu +\eps)\RR (M \tilde g_{\mu^n})  \ge g_2,
\eean
so that $\mu^n +\eps \in \JJ$. That means that $\mu^n+\eps \le \mu_1$, and a contradiction with the fact $\mu^n \nearrow \mu_1$. 
We next consider $\mu^n \nearrow \mu_1$ and we define 
$$
a^n := A(\tilde g_{\mu^n}), \quad \hat g^n := \frac{\tilde g_{\mu^n} }{ a^n}, 
\quad \eps_n  := \frac{g_2 }{ a^n},  \quad \hat g^n = \mu^n \RR \hat g^n + \eps_n . 
$$
We observe that $\eps_n \to 0$ and $A( \hat g^n ) = 1$. Because of the compactness assumption {\it (v)}, we deduce that 
there exists $f_1 \in K_2$ and a subsequence $(\hat g^{n_k})$ such that $\hat g^{n_k} \to f_1$ and $A(f_1) = 1$. 
We conclude by passing to the limit in the above 
almost first eigenvalue equations. 
\end{proof}

We  may compare  Theorem~\ref{theo:PLLbis} with the results presented in the previous section.   When $\LL$ satisfies condition \ref{H1}, we may set $\RR := \RR_{\LL}(\kappa_1)$ so that $\RR \in \BBB(X)$ and $\RR$ satisfies (i). 
In that case, Theorem~\ref{theo:PLLbis} claims the existence of $f_1 \in K_2$ such that  $\LL f_1 = \lambda_1 f_1$, with $\lambda_1 := \kappa_1 - \mu_1$. 
The condition (ii) on $\RR$ translates as $\LL g_2 \le (\kappa_1 - 1/C_2) g_2$ which may be seen as an equivalent of condition \ref{H1} (when working in the space $X_{2} := K_2 - K_2$ with norm $\| g \|_2 := A(|g|)$ and $\LL$ generates a semigroup $S$.
The hypothesis (iii) is nothing but \ref{H2} and the hypothesis (iv) is very natural: it holds in the space $L^p(E)$ and $M^1(E)$ without additional condition on $\RR$ and it holds in a space of continuous functions when some  additional uniform continuity assumption is made on the range of $\RR$. Assumption (v) has to be compared with condition \ref{H3}. It is worth emphasizing that when $X \subset L^p(E)$ and $g_2 > 0$ a.e., we simply have $A(g) = \| g/g_2 \|_{L^\infty}$ for any $g \in X_+$. 
As a conclusion, although Theorem~\ref{theo:exist1-KRexistence} and Theorem~\ref{theo:PLLbis} bear some similarities none seems to be a consequence of the other. We believe that Theorem~\ref{theo:exist1-KRexistence} is more flexible since it does not impose to work with the normalization associated to the  seminorm $ g \mapsto A(|g|)$ of $L^\infty$-type.   In a sense, above Lions’ proof is dual to the probabilistic approaches in \cite{MR2834717,MR3102473,BCG2019,Bansaye2022}, these last ones being more in the spirit of Lemma~\ref{lem:H3M1}, Lemma~\ref{lem:H3Y'} and Corollary~\ref{cor:H3M1}.  
It is also worth emphasizing  the very similarity between Step 3 in the proof of Theorem~\ref{theo:PLLbis} and the proof of Lemma~\ref{lem:Existe1-Spectre2} and, on the other hand,  that Theorem~\ref{theo:intro-KRresult} is a particular case of Theorem~\ref{theo:PLLbis} by essentially exploiting the fundamental Lemma~\ref{lem:exists-lemfonda2} as
 shown in \cite{PL2}.
 We finally point out that when $Y = X'$, the condition (iv) is equivalent to a property of Banach lattices known as \emph{order continuous norm}, see for instance~\cite[Definition~2.4.1]{MR1128093}, as a consequence of~\cite[Thm.~2.4.2~(iii)]{MR1128093} along with the fact that weakly convergent increasing sequences in Banach lattices are automatically norm convergent, see {\it e.g.} \cite[Prop.~1.4.1]{MR1128093}.

  
%


%
%


\bigskip\bigskip
\section{Existence through a dynamical approach}
\label{sec:DynamicalExistenceKR}


In this part, we develop a  dynamical approach for proving the existence part of the Krein-Rutman Theorem.
We thus always consider a positive semigroup $S = S_\LL$ on a Banach lattice $X$ {\Blue  and we  give an alternative proof of the existence of $(\lambda_1,f_1)$ under conditions implying \ref{H1}-\ref{H2} and some structure conditions  in the spirit of those in the previous Section which imply \ref{H3}, see in particular  \ref{HS2}  and \ref{HS3}.} 
Above all, we 
are able to extend the existence part of the Krein-Rutman Theorem to a more general framework, namely to the case when $\LL$ only enjoys a suitable weakly dissipative condition.

 \medskip
\subsection{About dissipativity}
\label{subsect-Exist2-Dissip}
\

Let us start by recalling some classical definitions and results. 
We say that an operator $\LL$ defined in a Banach space $X$  is dissipative  if there is some number $\kappa \in \R$  
such that 
\beqn\label{eq:DefDissip}
\forall \, f \in D(\LL), \ \exists \, f^* \in J_{f}, \quad \Re e \langle f^*, \LL f \rangle \le \kappa \| f \|^2,
\eeqn
where we define the associated dual set $J_f \subset X'$ of $f$ by 
\beqn\label{eq:DefApplDual1}
J_{f} := \{ \varphi \in X'; \,\, \langle  \varphi , f \rangle = \| f \| = \| \varphi \|_{X'} \}.
\eeqn
In that situation and in order to be more precise, we should say that  $\LL-\kappa$ is dissipative. 
It is worth emphasizing that $J_f \not= \emptyset$ thanks to the corollary \eqref{eq:f*fge0} of the Hahn-Banach dominated extension theorem.
{\Blue
At least formally, denoting $f_t := S(t) f$, for $f \in D(\LL)$, we deduce from \eqref{eq:DefDissip} that 
$$
\frac12\frac{d}{dt} \| f_t \| ^2 = \Re e \langle (f_t)^*, \LL f_t \rangle \le \kappa \|  f_t \|^2, 
$$
and together with the Gr{\"o}nwall lemma, we deduce 
$$
\|  S(t) f  \| \le e^{\kappa t} \|  f \|, \quad \forall \, t \ge 0, 
$$
which is nothing but \eqref{eq:SMwt}.} 
Quite similarly, when 
\beqn\label{eq:exist2-subsupEF}
\exists \, \psi \in Y_+ \backslash \{ 0 \}, \ \exists \kappa \in \R, \quad \pm \LL^* \psi \le \kappa \psi,
\eeqn
 we may compute 
$$
\pm \frac{d }{ dt} \langle f_t,\psi\rangle = \pm \langle \LL f_t,\psi\rangle =\pm  \langle  f_t,\LL^*\psi\rangle \le \kappa  \langle  f_t,\psi\rangle,
$$
and together with the Gr{\"o}nwall lemma, we get 
\beqn\label{eq:exist2-subsupEF2}
\pm \langle S_t f ,\psi\rangle \le \pm e^{\pm\kappa t}  \langle   f ,\psi\rangle, \quad \forall \, t \ge 0. 
\eeqn

\smallskip
Two important more accurate versions of the previous ones are presented now. They will be of main importance in the sequel. On the one hand, we may assume that $\LL$ satisfies a Lyapunov type condition, namely
there exists $\psi_i \in Y_+$ and $\kappa \in \R$ such that 
\beqn\label{eq:exist2-DefLyapCond}
\LL^* \psi_2 \le \kappa \psi_2 + \psi_0,  
\eeqn
 with $\psi_2 > 0$ and $\psi_0/\psi_2 \to 0$ at infinity.  
For $f_t = S_\LL(t) f$, $f \in D(\LL) \cap X_+$, a similar computation as above gives
\bean
\frac{d }{ dt} \langle f_t , \psi_2 \rangle
=   \langle f_t , \LL^*\psi_2 \rangle
\le  \kappa \langle f_t , \psi_2 \rangle   + \langle f_t , \psi_0 \rangle.
\eean

Denoting $[f]_i := \langle |f|, \psi_i \rangle$ and using the Gr{\"o}nwall lemma, we classically deduce 
\beqn\label{eq:exist2-estimLyap}
[ S(t) f ]_2 \le  e^{\kappa  t}  [ f ]_2 +   \int_0^t e^{\kappa  (t-s)} [S(s) f]_0 \, ds, \quad \forall \, t \ge 0.
\eeqn
The Lyapunov condition \eqref{eq:exist2-DefLyapCond} is particularly relevant and useful in a Radon measures space framework $X = M^1_{\psi_2}(E)$ for some weight function $\psi_2$ on $E$.

\smallskip
On the other hand, we may generalize the above  Lyapunov condition by assuming the structure condition 

\begin{enumerate}[label={\bf(HS2)},itemindent=14mm,leftmargin=0mm,itemsep=1mm]
\item\label{HS2}  there exist  a splitting  $\LL = \AA + \BB$ and $\kappa_\BB \in \R$ such that $\AA$ is  $\BB$-bounded, that means 
$$
\exists \, C \ge 0, \ \forall \, f \in X, \quad \| \AA f \| \le C ( \| f \| + \| \BB f \|),
$$
the operator $\BB$ generates a semigroup $S_\BB$ and 
\beqn\label{eq:exist2-HS2}
\| ( S_{\BB} \AA)^{(*\ell)} *S_{\BB} (t)  \|_{\BBB(X)} = \OO(e^{\alpha t}), \qquad   \forall \, t > 0, 
\eeqn
for any $\ell  \ge 0$ and $\alpha > \kappa_\BB$.
\end{enumerate}
 Here and below, for two functions $U : \R_+ \to \BBB(\XX_0,\XX_1)$ and $V : \R_+ \to \BBB(\XX_1,\XX_2)$, we define the convolution function
$$
(V * U)(t) := \int_0^t V(t-s) U(s) \, ds, 
$$
when the integral is well-defined. For $U : \R_+ \to \BBB(\XX)$, we also recursively define $U^{(*0)} = I$ and $U^{(*(\ell+1))} = U^{(*\ell)} * U$. 
Using this convolution notation, the Duhamel formula writes
$$
S_\LL = S_\BB + S_\BB \AA * S_\LL, 
$$
and iterating this formula, for any $N \ge 1$, we get the following iterated Duhamel formula
\beqn\label{eq:itratedDuhamel}
S_\LL = S_\BB + \dots+ (S_\BB \AA)^{*(N-1)} * S_\BB + (S_\BB \AA)^{(*N)} * S_\LL. 
\eeqn
When $S_\LL$ is well defined in another space $X_0 \supset X$ and the last iterated convolution term  enjoys the regularity property $\| (S_\BB \AA)^{(*N)} (t) \|_{\BBB(X_0,X)} =
\OO(e^{\alpha t})$ for all $t > 0$ and $\alpha > \kappa_\BB$, we deduce from the above iterated Duhamel formula, the estimate 
\beqn\label{eq:exist2-estimDuhameN}
\| S(t) f \| \le C_0 e^{\alpha t} \| f \| + C_1  \int_0^t e^{\alpha (t-s)} \| S(s) f \|_0 \, ds, \quad \forall \, t \ge 0,  \ \alpha > \kappa_\BB, \ \forall \, f \in X,
\eeqn
for some constants $C_i \ge 1$ and where $\| \cdot \|_0$ stands for the norm in $X_0$. We may observe that the estimate \eqref{eq:exist2-estimLyap} in the case of a Lyapunov condition is a particular case of \eqref{eq:exist2-estimDuhameN} corresponding to the norms $\| \cdot \|= [\cdot]_2$ and $\| \cdot \|_0= [\cdot]_0$. More specifically,  in a Radon measures space framework, the splitting condition  \ref{HS2}  is obtained by introducing the bounded operator $\AA f := f \psi_0$ and the generator $\BB := \LL - \AA$. Because of \eqref{eq:exist2-DefLyapCond}, we have $\BB^* \psi_2 \le \kappa \psi_2$, and arguing as for establishing  \eqref{eq:exist2-estimLyap}, we have $[ S_\BB(t) f ]_2 \le  e^{\kappa t}  [ f ]_2$ for any $t \ge 0$ and $f \in X$. That last growth condition is equivalent to assuming that $\BB-\kappa$ is dissipative for the norm $[\cdot]_2$, so that we have established that $\LL$ enjoys the  splitting condition  \ref{HS2}.

 \medskip
\subsection{Existence in the dissipative case}
\label{subsect-Exist2-DissipExist}
\

In this section, we give an existence result for a positive semigroup $S_\LL$ on a Banach lattice $X$ satisfying a kind of regularity/compactness assumption in the spirit of the structure condition \ref{HS2} discussed above. {\Blue Here and below, for an element $\psi \in Y_+$, we define the (semi)norm
$$
 [f]_\psi := \langle \psi, |f| \rangle, \quad \forall \, f \in X.
$$
We will sometime use the shorthand $[\cdot]_{i}$ for $[\cdot]_{\psi_i}$.
}

\begin{theo} \label{theo:KRexistBIS} 
On a Banach lattice $X=Y'$, with $Y$ separable Banach lattice, consider a positive semigroup $S= S_\LL$ satisfying the growth bound \eqref{eq:SMwt}, 
and set $\kappa_1 := \omega' + \log M$ for some $\omega' > \omega(S_\LL)$.  

We assume  

\smallskip
\quad {\bf (1)} $\exists \,  \phi_0 \in  Y_{+} \backslash \{0 \}$, $\exists \, \kappa_0 \in \R$ such that  $[S (t) f]_{\phi_0} \ge e^{\kappa_0 t} [f]_{\phi_0}$ for any $t \ge 0$ and $f \in X_+$; 
%

\smallskip
\quad {\bf (2)} there exist $\kappa,  C_0, C_1 \in \R$ with $\kappa <  \kappa_0  $, $C_0 \ge 1$ and $C_1 \ge 0$,  such that 
\beqn\label{eq:theo:KRexistBIS:borneSt}
\| S(t) f \| \le C_0 e^{\kappa t} \| f \| + C_1 \int_0^t e^{\kappa (t-s)} [S(s) f]_{\phi_0} \, ds, \quad \forall \, t \ge 0, \ \forall \, f \in X. 
\eeqn

Then there exist $\lambda_1 \in [\kappa_0,\kappa_1]$ and $f_1 \in X_+\backslash \{ 0 \}$ such that $\LL f_1 = \lambda_1 f_1$. 
\end{theo}

Let us mention that this result shares similarities with~\cite[Cor.~2.7]{Martinez1993} and \cite[Thm.~4.2]{MR2834717},  see also \cite{MR0859722,MR0901175} for earlier works in that direction.

\begin{rem}\label{rem:theo:KRexistBIS}

  \smallskip
 (1)  Assumption {\bf (2)}  in the statement of Theorem~\ref{theo:KRexistBIS} holds when there exist $V,W$ such that 
\beqn\label{eq:theo:KRexistBIS-representation}
S = V + W* S, \quad W \ge 0, 
\eeqn
and there exist $\kappa, C_V, C_W \in \R$, $\kappa < \kappa_0$, $C_V \ge 1$, $C_W > 0$ such that 
\beqn\label{eq:theo:KRexistTER1}
\| V (t) \|_{\BBB(X)} \le C_V e^{\kappa  t} , \quad  \| W (t) \|_{\BBB(\XX_0,X)} \le C_W e^{\kappa  t}.
\eeqn

\smallskip
(2)  Under the structural condition \ref{HS2} together with some regularization effect on the semigroup of the type 
$$
\| ( S_{\BB} \AA)^{(*N)}   (t)  \|_{\BBB(\XX_0,X)} = \OO(e^{\kappa t}), \qquad   \forall \, t > 0,  \ \kappa \in ( \kappa_\BB,\kappa_0),
$$
we recover the above condition~\eqref{eq:theo:KRexistBIS-representation}-\eqref{eq:theo:KRexistTER1} with 
\beqn\label{eq:V=&W=} 
V := S_\BB + \dots+ (S_\BB \AA)^{*(N-1)} * S_\BB, \quad W := ( S_{\BB} \AA)^{(*N)} , 
\eeqn
because of the iterated Duhamel formula \eqref{eq:itratedDuhamel}. In that case,  the representation formula \eqref{eq:Exist1-DefRepresentationRR} holds true for any $z > \lambda_1$ from  Lemma~\ref{lem:Exist1-RkSG}-{\bf (ii)} and we easily compute
$$
\RR_\LL(z) = \VV(z) + \WW(z) \RR_\LL(z), \quad \forall \, z > \lambda_1, 
$$
with 
$$
\VV(z) := \int_0^\infty e^{-\lambda t} V(t) dt, \quad
\WW(z) := \int_0^\infty e^{-\lambda t} W(t) dt,  \quad \forall \, z > \kappa.
$$
We observe that $\WW$ satisfies  \eqref{eq:H3abstractWeakX0}  in  Lemma~\ref{lem:H3abstractWeakX0} if $W$ satisfies \eqref{eq:theo:KRexistTER1}
and the set $\CC$ defined by \eqref{eq:H3abstract-defC}  satisfies the same compactness properties  as required in the statement of Lemma~\ref{lem:H3abstractWeakX1}.
We may thus apply Lemma~\ref{lem:H3abstractWeakX0} (see also Remark~\ref{rem:H3Lp&H3XX1}) and deduce that \ref{H3} holds for the primal problem. 
We finally obtain the same conclusion as in Theorem~\ref{theo:KRexistBIS} thanks to Theorem~\ref{theo:exist1-KRexistence}.

\smallskip 
(3) Under the same structural condition \ref{HS2} as above, but assuming now that 
$$
\| W   (t)  \|_{\BBB(X,\XX_1)} = \OO(e^{\kappa t}), \qquad   \forall \, t > 0,  \ \kappa \in ( \kappa_\BB,\kappa_0),
$$
with $W := ( S_{\BB} \AA)^{(*N)}$ and $\XX_1 \subset X$ with strongly compact embedding, we observe that $S$ does not necessary satisfies the assumptions of Theorem~\ref{theo:KRexistBIS}, but it rather satisfies the 
assumptions of Lemma~\ref{lem:H3abstract-StrongCbis} with $K_T := (W * S)(T)$ and $T > 0$ large enough.
In that situation, we also obtain the same conclusion as in Theorem~\ref{theo:KRexistBIS} thanks to  Lemma~\ref{lem:H3abstract-StrongCbis} and Theorem~\ref{theo:exist1-KRexistence}.

 \smallskip
 {\Cyan (4) The condition \ref{H1} follows from Lemma~\ref{lem:Exist1-RkSG} and the above assumption {\bf (1)} implies \ref{H2} thanks to Lemma~\ref{lem:Existe1-Spectre2bis}. 
 }

\end{rem}

\begin{proof}[Proof of Theorem~\ref{theo:KRexistBIS}.]
We split   the proof into two steps.
{\Cyan During the proof, we will use the shorthand $[\cdot]_0 := [\cdot]_{\phi_0}$.}

\smallskip
{\sl Step 1.}   
We define the set 
$$
\CC := \{ f \in X_+, \ [f ]_0 = 1, \ \| f \| \le R \}, 
$$
for a convenient constant $R > 0$ to be fixed later. 
For any fixed $t > 0$, we next define the nonlinear weakly $\sigma(X,Y)$ continuous mapping
$$
\Phi_t : \CC \to X, \quad f \mapsto \frac{S_tf }{ [S_t f]_0   }.
$$
Thanks to assumption {\bf (1)}, we may observe that it is well defined because 
\beqn\label{eq:St0kappa2}
[S_t f]_0   \ge e^{\kappa_0 t}  [f]_0 = e^{\kappa_0 t}  > 0.  
\eeqn
For any $f \in \CC$, we thus immediately have $\Phi_t f \ge 0$ and $[\Phi_t f]_0  =  1$. 
On the other hand, from assumption {\bf (1)} again and the semigroup property, we have 
\beqn\label{eq:St0kappa2BIS}
[S (t) f]_0 \ge e^{\kappa_0 (t-s)} [S(s) f]_0. 
\eeqn
{\Cyan Then, for $f \in \CC$ and $t \ge 0$, taking advantage of \eqref{eq:theo:KRexistBIS:borneSt} and \eqref{eq:St0kappa2BIS},} we compute 
\bean
\| \Phi_t f \| 
&\le&  C_0  e^{-\alpha t} \| f \|   + C_1 \int_0^t  e^{-\alpha (t-s)}  \, ds
\\
&\le& C_0  e^{-\alpha t}  R   + \frac{C_1 }{ \alpha}, 
\eean
where we have set $\Cyan\alpha := \kappa_0 - \kappa > 0$.
Fixing $T_0$ such that $C_0 e^{-\alpha T_0} \le 1/2$ and next $R \ge 2 C_1/\alpha$, we have thus 
$\Phi_{T_0} : \CC \to \CC$. Thanks to the Tykonov fixed point Theorem, there exists $f_{T_0} \in \CC$ such that $\Phi_{T_0} f_{T_0} = f_{T_0}$.
In other words, we have established the existence of $f_{T_0} \in X$ such that 
\beqn\label{eq:periodicEqfT0}
f_{T_0} \ge 0, \quad [f_{T_0}]_0 = 1, \quad S_{T_0} f_{T_0} = e^{\lambda_1 T_0} f_{T_0}, 
\eeqn
with $\lambda_1 := (1/T_0)\log  [S_{T_0} f_{T_0}]_0 \in [\kappa_0,\kappa_1]$.

 \medskip
\noindent 
{\sl   Step 2.}  Rewriting equation \eqref{eq:periodicEqfT0} as 
\[0=e^{-\lambda_1T_0}S_{T_0}f_{T_0}-f_{T_0}=(\LL-\lambda_1)\int_{0}^{T_0}e^{-\lambda_1t}S_tf_{T_0}dt
\]
and defining 
$$
f_1:=\int_{0}^{T_0}e^{-\lambda_1t}S_tf_{T_0}dt, 
$$
we get that  $f_1 \in X_+\backslash \{0\}$ satisfies $\LL f_1=\lambda_1f_1$.
\end{proof}

  We present now a second  proof based on  a large times dynamical argument
which is classical in the  mean ergodicity theory of Von {N}eumann and Birkhoff introduced in \cite{vonNeumann,Birkhoff}
and which will be adapted in the weak dissipativity case in Section~\ref{subsec-Exist2-weakD2} below. 

\begin{proof}[Alternative Step 2.]
We define $\widetilde S_t := S_t e^{-\lambda_1t}$, so that $f_{T_0}$ becomes a periodic state for $\widetilde S_t $ from \eqref{eq:periodicEqfT0}, namely   
$$
\widetilde S_t f_{T_0} = \widetilde S_{t-kT_0} f_{T_0}, \quad k := [t/T_0], \quad \forall \, t > 0. 
$$
Using \eqref{eq:St0kappa2} and the above relation, we have 
\bean
[\widetilde S_t f_{T_0}]_0 &=& [\widetilde S_{t-kT_0} f_{T_0}]_0
\\
&\ge& e^{(\kappa_0 - \lambda_1)(t-kT_0)}  [  f_{T_0}]_0 \ge e^{(\kappa_0 - \lambda_1) T_0} =: r_* > 0, 
\eean
for any $t \ge 0$. On the other hand, thanks to the growth bound \eqref{eq:SMwt},  we have
\bean
\| \widetilde S_t f_{T_0} \| &=& \| \widetilde S_{t-kT_0} f_{T_0} \|
\\
&\le&  M e^{(\kappa - \lambda_1)(t-kT_0)} \|  f_{T_0} \|   \le  M e^{(\kappa - \lambda_1) T_0}  R =: R^* < \infty, 
\eean
for any $t \ge 0$. We finally define 
$$
u_T := \frac1T \int_0^T \widetilde S_t f_{T_0} \, dt.
$$
 From the previous estimates, both sequences $(\widetilde S_t f_{T_0})$ and $(u_T)$ are   bounded in 
$$
{\mathbb K} := \{ f \in X; \ f \ge 0, \ [f]_0 \ge r_*, \ \| f \| \le R^* \}. 
$$
 By compactness,  there exists a
  subsequence $(u_{T_k})$ and $f_1 \in {\mathbb K}$ such that
  $u_{T_k} \wto f_1$ in a weak sense 
  as
  $k\to\infty$. For any fixed $t > 0$, we observe that 
  \bean 
  \widetilde S_t f_1 - f_1 
 &=& \lim_{k \to \infty} \Bigl\{ \frac{1 }{ T_k}
  \int_0^{T_k} \widetilde S_t \widetilde S_s f_{T_0} ds - \frac{1 }{ T_k} \int_0^{T_k} \widetilde S_s f_{T_0}  ds  \Bigr\}
  \\
  &=& \lim_{k \to \infty} \Bigl\{ \frac{1 }{ T_k} \int_{T_k}^{T_k+t }\widetilde S_s f_{T_0} ds - \frac{1 }{ T_k} \int_0^t \widetilde S_s f_{T_0} ds \Bigr\} = 0, 
  \eean 
  where we have used that $(\widetilde S_s f_{T_0})$ is uniformly bounded in the last line. 
  As a consequence, $f_1$ is a stationary state for the rescaled semigroup $\widetilde S_t$, and thus an eigenfunction associated to the eigenvalue $\lambda_1$ for the operator $\LL$.  
\end{proof}

\subsection{About weak dissipativity}
\label{subsect:AboutWeakDissip}

\

In this section, we recall some definitions and results about the weak dissipativity. 
We say that the generator $\BB$ of a semigroup $S_\BB$ is weakly dissipative in a Banach space $X_i$ if there exist a second Banach space $X_{i-1} \supset X_i$ and some numbers $\kappa \in \R$ and $\sigma > 0$
such that 
$$
\forall \, f \in D(\BB_{|X_i}), \ \exists \, f^* \in J_{f,X_i}, \quad \langle f^*, \BB f \rangle \le \kappa \| f \|^2_{X_i} - \sigma \| f \|^2_{X_{i-1}}, 
$$
where we define the associated dual set $J_{f,X_i} \subset X_i'$ of $f$ (for the norm $\| \cdot \|_{X_i}$) by 
\beqn\label{eq:DefApplDual}
J_{f,X_i} := \{ \varphi \in X_i'; \,\, \langle  \varphi , f \rangle = \| f \|_{X_i}^2 = \| \varphi \|_{X_i'}^2 \}.
\eeqn
By translation, we may assume that $\kappa=0$, an hypothesis we will always make in the sequel of this section. 
We will furthermore assume the splitting structure $\LL = \AA + \BB$ with $\AA$ a $\BB$-bounded operator and $\BB$ a weakly dissipative generator. 

\smallskip
More precisely, we assume that there exists one more Banach lattice $X_0 \supset X_1 \supset X_{2} := X$, with norm denoted by $\| \cdot \|_k := \| \cdot \|_{X_k}$, 
such that $\BB$ generates a semigroup and is weakly dissipative in each $X_k$:  for any $k =1,2$ 
\beqn\label{eq:defWeakDissip}
\forall \, f \in D(\BB_{|X_k}), \ \exists \, f^* \in J_{f,X_k}, \quad \langle f^*, \BB f \rangle_{X_k',X_k} \le - \sigma \| f \|^2_{k-1}. 
\eeqn
This classically implies (or we can take the next inequality as a definition of the weak dissipativity condition) that 
\beqn\label{eq:SBBk}
\frac{d }{ dt} \| S_\BB(t) f \|_k + \sigma \| S_\BB(t) f \|_{k-1} \le 0, \quad \forall \, t \ge 0, \ \forall \, f \in X_k, \ \forall \, k =1,2.
\eeqn
We assume that $X_{k}$ is dense into $X_{k-1}$ for $k=1,2$ and that $X_1$ is an interpolated space between $X_0$ and $X_2$ in the sense that 
there exists a continuous and strictly decreasing  function $\eta : (0,1] \to [0,\infty)$, $\eta(\eps) \to \infty$ when $\eps \to 0$, $\eta(1) = 0$, such that 
\beqn\label{eq:interpolXk}
\qquad\| f \|_1 \le \eps \| f \|_{2} + \eta(\eps) \| f \|_{0}, \quad \forall \, \eps \in (0,1], \ \forall \, f \in X_2.
\eeqn
From \eqref{eq:SBBk} with $k=2$, we deduce 
\beqn\label{eq:ContractionSB}
\| S_\BB(t) f \|_{2} \le \|  f \|_{2}, \quad \forall \, t \ge 0, \ \forall \, f \in X_{2}.
\eeqn
Next, for $k=1$, gathering the weak dissipativity condition \eqref{eq:SBBk}, the interpolation condition  \eqref{eq:interpolXk} and  the non expansion property \eqref{eq:ContractionSB} in the space $X_2$, we get
\bean
\frac{d }{ dt} \| S_\BB(t) f \|_1 + \frac{\sigma}{\eta(\eps)} \  \| S_\BB(t) f \|_1 
&\le& \frac{\sigma\eps}{\eta(\eps)} \  \| S_\BB(t) f \|_{2}
\\
&\le&\frac{\sigma\eps}{\eta(\eps)} \  \|   f \|_{2}, 
\eean
for any $ t \ge 0$, $\eps \in (0,1)$ and $f \in X_{2}$. We deduce 
\bean
\frac{d }{ dt} \Bigl( \| S_\BB(t) f \|_1 e^{\frac{\sigma}{\eta(\eps)} t} \Bigr)
 \le \frac{\sigma\eps}{\eta(\eps)} e^{\frac{\sigma}{\eta(\eps)} t}  \|   f \|_{2}, 
\eean
and thanks to the  Gr{\"o}nwall lemma, we obtain
\beqn\label{eq:decaySBk}  
\| S_\BB(t) f \|_1 \le \Theta(t)  \|   f \|_{2}, 
\eeqn
for any $ t \ge 0$  and $f \in X_{2}$, with 
\beqn\label{eq:decaySBkTheta}  
 \Theta(t) := \inf_{\eps \in (0,1)} \bigl( e^{- \frac{\sigma}{\eta(\eps)} t} + \eps \bigr) \to 0
 \ \hbox{ as } \ t \to +\infty.
\eeqn
On the other hand, using the representation formula 
$$
\RR_\BB(z) f = \int_0^\infty e^{-zt} S_\BB(t) f \, dt, \quad \forall \, z \in \Delta_0, \ \forall \, f \in X_2, 
$$
together with  \eqref{eq:SBBk}, we get 
$$
\sigma \| \RR_\BB(z) f \|_{1} \le \int_0^\infty \sigma \| S_\BB(t) f \|_{1} \, dt \le \| f \|_2, 
$$
for any $z \in \overline{ \Delta}_0$ and $f \in X_2$.  We next assume that
\beqn\label{eq:decayASBk}
\Theta(t)^{-1} \| \AA S_\BB(t) f \|_{1} + \int_0^\infty \| \AA S_\BB(t) f \|_{1} dt \lesssim  \|   f \|_{1}, 
\eeqn
that  there exist $\alpha > 1$, $N \ge 1$, $C \ge 1$ such that 
\beqn\label{eq:estimARBk}
\sup_{x+iy \in \Delta_0} \| \AA \RR_\BB^{1+\eps_1}(x+iy) \dots \AA \RR_\BB^{1+\eps_N} (x+iy) f \|_2 \le \frac{C }{ \langle y \rangle^\alpha} \| f \|_{2}, 
\eeqn
for any  $\eps \in \{0,1\}^N$, $\eps_1 + \dots + \eps_N \le 1$, and that
\beqn\label{eq:estimARBN}
\sup_{z \in \Delta_0}  \|  (\RR_\BB (z) \AA)^N f \|_{\XX_1}  \le  \| f \|_{1}, 
\eeqn
with $\XX_1$ compactly imbedded in $X_1$. The necessity to add $(\eps_i)$ in \eqref{eq:estimARBk} is probably purely technical and not restrictive for applications. In examples, 
we can take $N=2N'$, when 
\beqn\label{eq:estimARBkBIS}
\sup_{x+iy \in \Delta_0} \| (\AA \RR_\BB)^{N'} (x+iy) f \|_3 \le \frac{C }{ \langle y \rangle^\alpha} \| f \|_{2},  
\eeqn
for some convenient space $X_3$ such that $\AA: X_1 \to X_3$ and $\sup_{z \in \Delta_0} \| \RR_\BB (z) \|_{\BBB(X_3,X_2)} < \infty$. 
 At the level of the semigroup, \eqref{eq:estimARBkBIS} is typically a consequence of  
$$
\| (\AA S_{\BB} )^{(*N'')} (t) \|_{\BBB(X_2,X^\zeta_3)} \in L^1(\R_+), 
$$
with $\zeta > 0$, where $X_3^\zeta := \{f \in X_3, \ \LL^\zeta f \in X_3 \}$ stands for the (possibly fractional) domain for the operator defined in $X_3$. 
However, \eqref{eq:estimARBk} is a bit more general than that last estimate. We refer to \cite{MR3489637,MischErratum,MR3850019,MR3859527} for precise definition, examples and discussion.
 For further references, we observe that \eqref{eq:decaySBk} and \eqref{eq:decayASBk} together imply 
\bean
 \frac1T\int_0^T   \| (S_\BB * \AA S_\BB)(t) f \|_1 \, dt
&\le& \frac1T \int_0^T \!\int_0^t \| S_\BB (t-s) \AA S_\BB(s) f \|_1 \, ds dt
 \\
&\le& \frac1T \int_0^T \!\int_0^T \| S_\BB (u) \|_{\BBB(X_{2},X_1)} \|  \AA S_\BB(s) f \|_1 \, du ds
 \\
&\lesssim& \frac1T \int_0^T  \Theta (u) \, du  \,  \|  f \|_{2}.   
\eean
Arguing in a similar way for any $\ell \ge 1$, we establish 
\beqn\label{eq:estimCesaroSBASB}
\frac1T\int_0^T   \| (S_\BB * (\AA S_\BB)^{(*\ell)})(t) f \|_1 \, dt 
\lesssim  \frac1T \int_0^T  \!\!\Theta \, du  \,  \|  f \|_{2} \to 0 \ \hbox{ as } \ T \to \infty.
\eeqn

\smallskip
For synthesizing and for further references, let us now bring out some possible general framework for semigroup enjoying weak dissipativity.  We introduce the following  structure condition 
on a semigroup $S_\LL$ and its generator $\LL$ by assuming

\begin{enumerate}[label={\bf(HS3)},itemindent=14mm,leftmargin=0mm,itemsep=1mm]
\item\label{HS3}   
there exist  a splitting  $\LL = \AA + \BB$, some Banach lattices $X_2 \subset X_1$, an integer $N \ge 1$ and some decaying functions $\Theta_i : \R_+ \to \R_+$ with $\Theta_1(t) \to 0$ as $t \to \infty$, $\Theta_2 \in L^1(\R_+)$ such that   $\AA$ is  
positive, $\BB$ generates a positive  semigroup $S_\BB$  and the following estimates hold 
\bear\label{eq:exist2-HS3}
&&\| ( S_{\BB} \AA)^{(*\ell)} *S_{\BB}    \|_{\BBB(X_2,X_1)} = \OO(\Theta_1 ), \quad \forall \, \ell \in \{0, \dots, N-1\}, 
\\ \label{eq:exist2-HS3bis}
&&\| ( S_{\BB} \AA)^{(*N)}   \|_{\BBB(X_1,X_2)} = \OO(\Theta_2 ).
\eear
\end{enumerate}

\smallskip 
We now particularize our discussion  to a Radon measures framework. 
We assume that  there exist some weight functions $\psi_i$ on $E$,  $\psi_0 \lesssim \psi_1 \le \psi_2$, with  $ \psi_2(x)/\psi_1(x) \to \infty$ as $x \to \infty$ so that $M^1_{\psi_2} \subset\subset M^1_{\psi_1}$ (compact imbedding for the weak convergence), 
a function $\chi \in C_c(E)$, $0 \le \chi \le 1$, and a constant $M \ge 0$ such that 
\begin{itemize}
\item[(i)]  $\LL^* \psi_1 \le - \psi_0 +  M\chi$; 
\item[(ii)]  $\LL^* \psi_2 \le   M \chi$; 
\item[(iii)] $\psi_1 \le \eps \psi_2 + \eta(\eps) \psi_0$ for any $\eps > 0$,
\end{itemize}
for a  function $\eta : (0,1] \to (0,\infty)$ such that  $\eta(1) = 0$, $\eta(\eps) \to \infty$ when $\eps \to 0$, and 
\beqn\label{eq:decaySBkThetaBIS}  
t \mapsto \Theta(t) := \inf_{\eps \in (0,1)} \bigl( e^{- \frac{t}{\eta(\eps)} } + \eps \bigr) \in L^1(0,\infty).  
\eeqn
It is worth emphasizing that  from the very definition, we have automatically that $\Theta$ is positive and decreasing, $\Theta(0) = 1$  and $\Theta(t) \to 0$ as $t \to \infty$. 
Arguing similarly as we did during the proof   of Corollary~\ref{cor:H3M1} and the end of Section~\ref{subsect-Exist2-Dissip}, we introduce the splitting  
$$
  \AA := M \chi, \quad \BB := \LL - \AA,
$$
and we establish that $S_\BB$ is a positive semigroup on $X = M^1_{\psi_2}(E)$. More precisely, for $0 \le f_0 \in D(\BB)$ in the domain of $S_\BB$   and denoting $f_t := S_\BB(t) f_0$, we may compute 
$$
\frac{d}{ dt} \int f_t \, \psi_2 \le  \int f_t \, \BB^* \psi_2 \le 0
$$
and similarly 
$$
 \frac{d}{ dt} \int f_t \, \psi_1 \le \int f_t \, \BB^* \psi_1 \le  -  \int f_t \,  \psi_0. 
$$
Integrating both differential inequalities, we deduce $ S_{\BB} \in L^\infty_t (\BBB(M^1_{\psi_i}))$, $ i=1,2$ and
$$
 \int_0^\infty \|   S_\BB(t) f_0 \|_{M^1_{\psi_0}} dt \le \| f_0 \|_{M^1_{\psi_1}}, \quad 
  \forall \, f_0 \in M^1_{\psi_1}. 
  $$
 We may make a slight (but important) improvement of the previous estimate  by proceeding similarly as we did for proving \eqref{eq:decaySBk}. Using the same notations as in the above computation, we indeed have 
$$
\frac{d }{ dt} \int  f_t \, \psi_1 + \frac{1 }{ \eta(\eps)} \int f_t \, \psi_1  \le \frac{\eps }{ \eta(\eps)}  \int f_t \, \psi_2  \le \frac{\eps}{ \eta(\eps)}  \int f_0 \, \psi_2,
$$
 where we have used  (i) and (iii) in the first inequality and  the previous $L^\infty_t (\BBB(M^1_{\psi_2}))$ bound in the second inequality. Integrating in time, we deduce  
$$
\| S_\BB (t) f \|_{M^1_{\psi_1}} \le \Theta(t) \|   f \|_{M^1_{\psi_2}}, \quad \forall \, t  > 0.
$$
Taking $X_i := M^1_{\psi_i}$ and $N=1$, we see that $\LL$ then satisfies \ref{HS3} with $\Theta_i = \Theta$. 
\Black


 \medskip
\subsection{First existence result  in the weakly dissipative case} 
\label{subsec-Exist2-weakD1}

\
We first come back to the proof of Theorem~\ref{theo:exist1-KRexistence} and explain what goes wrong when we try to adapt it to a weak dissipativity context. 
More precisely, we assume that $S_\LL$ is a positive semigroup (so that \ref{H1} holds) satisfying  $\LL^* \phi_0 \ge 0$ for some $\phi_0 \in X' \backslash \{0 \}$ (so that \ref{H2} holds) 
 and the splitting structure \ref{HS3} for some bounded operator $\AA$ and some weakly dissipative operator $\BB$, in the sense that \eqref{eq:defWeakDissip} holds. 
In such a situation, we may define 
$$
\lambda_1 := \inf \{ \lambda \in \R; \ \RR_\LL(\kappa) \in \BBB(X), \ \forall \, \kappa \ge \lambda \} \ge 0, 
$$
and there exist sequences $(\lambda_n)$ of $\R$ and $(\hat f_n)$ of $X_+$ such that 
$$
\lambda_n \searrow \lambda_1 \ge 0, \quad \| \hat f_n \|=1, \quad \eps_n: = \lambda_n \hat f_n - \LL \hat f_n \to 0 \hbox{ in } X, 
$$
thanks to Lemma~\ref{lem:Existe1-Spectre2}. In the simplest situation, we may further assume that $\RR_\BB(\kappa): X_1 \to X_0$ is uniformly bounded in $\kappa \ge \lambda_1$ and 
$\AA: X_0 \to X_1$ with $X = X_1 \subset X_0$. The issue is that even in that case, we may write 
$$
\hat f_n = \RR_\BB(\lambda_n) \AA \hat f_n + \RR_\BB(\lambda_n)\eps_n,
$$
but  it is not clear how  to conclude that $(\hat f_n)$ belongs to a compact set in $X$ because it is not clear that $\RR_\BB(\lambda_n)\eps_n \to 0$ in $X$. 

\smallskip
The next result aim precisely to establish that last convergence under suitable quite strong (although natural and true in some examples) assumptions on the operator $\LL$. The proof is adapted from \cite[Sec.~6.3]{MR4265692} and mixes some dynamical argument together with the stationary approach developed in Section~\ref{subsec-Exist1-KRtheorem}.

  \begin{theo} \label{th:ExistKRweakDissip} Consider a positive semigroup $S_\LL$ in a Banach lattice  $X = X_2 \subset X_1 \subset X_0$  such that its generator $\LL$ satisfies

{\bf (1)}  there exists $\phi_0 \in D(\LL^*)$, $\phi_0 \ge 0$, $\phi_0 \not= 0$,  such that $\LL^* \phi_0 \ge 0$;

{\bf (2)} $\LL = \AA + \BB$ with $\AA$ and $\BB$ satisfying  \eqref{eq:decaySBk}, \eqref{eq:decayASBk}, \eqref{eq:estimARBk} and \eqref{eq:estimARBN}.

\noindent
  Then,  there exist $\lambda_1 \ge 0$ and $f_1 \in X_1$ such that 
\beqn\label{eq:0=LfBIS}
 {\| f_1 \|}_{X_1} = 1,  \quad  f_1 \ge 0,  \quad \LL f_1 = \lambda_1 f_1 .
\eeqn
 \end{theo}

\begin{proof}[Proof of Theorem~\ref{th:ExistKRweakDissip}.]
We split the proof into four steps. 

 \smallskip
{\sl Step 1. } 
We know from Lemma~\ref{lem:Exist1-RkSG} and  Lemma~\ref{lem:Existe1-Spectre2bis}-{\bf (i)}   that \ref{H1} and \ref{H2} hold. We may then define $\lambda_1 \ge 0$ with the help of \eqref{eq:exist1-deflambda1}. 
If $\lambda_1 > 0$, we see that $\VV(\alpha)$ defined in \eqref{eq:exist1-defVVWW} is  bounded in $\BB(X)$ uniformly on $\alpha \ge \kappa_0 := \lambda_1/2$ because of  \eqref{eq:decaySBk} and \eqref{eq:decayASBk},
and that  $\WW(\alpha)$ also defined in \eqref{eq:exist1-defVVWW} satisfies \eqref{lem:H3abstract-StrongCompact-c1}
because of  \eqref{eq:decayASBk} and Remark~\ref{rem:H3abstract-compactStrong}-{\bf (1)}. 
Using Lemma~\ref{lem:H3abstract-StrongC}, we get that \ref{H3} holds, and we conclude thanks to Theorem~\ref{theo:exist1-KRexistence} in that case. 

\smallskip
In the sequel, we always assume $\lambda_1 = 0$.

 \smallskip
{\sl Step 2. } 
Let us fix $f_0 \in D(\LL)$ such that $f_0 \ge 0$ and $C_0 :=  \langle f_0, \phi_0 \rangle  > 0$, which exists by definition of $\phi_0$.  Denoting $f(t) := S_\LL(t) f_0$, we have 
$$
\frac{d }{ dt} \langle f(t) , \phi_0 \rangle =  \langle \LL f(t) , \phi_0 \rangle = \langle  f(t) , \LL^* \phi_0 \rangle \ge 0, 
$$
which in turns implies
$$
\langle f(t) , \phi_0 \rangle \ge C_0,  \quad \forall \, t \ge 0. 
$$

 \smallskip\noindent
{\sl Step 3.} 
\Black
We claim that $ \| \RR_\LL(0) \|_{\BBB(X_2,X_1)} = +\infty$. That in particular implies  $\| \RR_\LL(0) \|_{\BBB(X)} = +\infty$ and thus $0 \in \Sigma(\LL)$.
We assume by contradiction that  $\KK_{2,1} := \| \RR_\LL(0) \|_{\BBB(X_2,X_1)}  < +\infty$. 
First, because $S_\LL$ is positive, we have 
$$
|\RR_\LL(z)f|  \le  \int_0^\infty e^{-t \Re e z} S_\LL(t) |f| \, dt = |\RR_\LL(\Re e z) |f|,
$$
from which we deduce
$$
\| \RR_\LL(z) \|_{\BBB(X_2,X_1)}  \le \| \RR_\LL(\Re e z) \|_{\BBB(X_2,X_1)} , \quad \forall \, z \in \Delta_{0}.
$$
As a consequence, we have 
\beqn\label{hyp:bddRL}
\sup_{y \in \R} \| \RR_\LL(iy) \|_{\BBB(X_2,X_1)}  \le \KK_{2,1}.
\eeqn
 We write the representation formulas (taken from \cite[(2.21)]{MR3489637}) 
$$
S_\LL(t)f =  \TT_0 (t) + \lim_{M \to \infty} \TT_{1,M}(t)
$$
with 
$$
\TT_0(t) :=  \sum_{\ell=0}^{N-1}  S_\BB *  (\AA S_\BB)^{(*\ell)} (t)f 
$$
and 
$$
\TT_{1,M}(t)  :=   \frac{ i }{  2\pi}
 \int_{a - {\rm i}M}^{a + {\rm i}M}\!\!\!\!  {\rm e}^{zt} \,  
 \RR_\LL(z) \, (\AA \RR_\BB(z))^{N} f \, dz,
$$
for any $f \in D(\LL)$,  $t \ge 0$ and $a > 0$.  On the one hand, from \eqref{eq:estimCesaroSBASB}, we have the Ces\`aro mean convergence 
\beqn\label{eq:KRwEstim1}
\frac1T \int_0^T \TT_0(t) \, dt \to 0 \ \hbox{ in } X_1, \ \hbox{ as } \ T \to \infty. 
\eeqn
On the other hand, we estimate the contribution of $\TT_{1,M}$. Integrating by part, we have 
$$
\TT_{1,M} (t) = \frac{1 }{ t}
\frac{ i}{  2\pi}
 \int_{a - {\rm i}M}^{a + {\rm i}M}\!\!\!\!  {\rm e}^{zt} \,  \frac{d }{ dz} \big[\RR_\LL(z) \, (\AA \RR_\BB(z))^{N} \big] f\, dz,
$$ 
with
$$
\frac{d }{ dz} \big[\RR_\LL(z) \, (\AA \RR_\BB(z))^{N} \big] = \sum_{\eps \in \N^{N+1}, \, |\eps|=1} \RR_\LL(z)^{1+\eps_0} \AA \RR^{1+\eps_1}_\BB(z)
\dots  \AA \RR^{1+\eps_N}_\BB(z).
$$
Together with condition \eqref{eq:estimARBk} and estimate \eqref{hyp:bddRL}, we get
\bean
&&\Big\| \frac{d }{ dz} \big[\RR_\LL(z) \, (\AA \RR_\BB(z))^{N} \big] f \Big\|_1
\\
&&\quad\le ( \KK_{2,1} + \KK_{2,1}^2) N  \sup_{\eps \in \N^{N}, \, |\eps|\le1} \| \AA \RR^{1+\eps_1}_\BB(z)
\dots  \AA \RR^{1+\eps_N}_\BB(z) f \|_2
\\
&&\quad\le   \frac{C_1 }{ \langle y \rangle^\alpha} \| f \|_{2},
\eean 
uniformly for any $z = x+iy \in \Delta_{0}$, for some constant $C_1 > 0$. We deduce 
\beqn\label{eq:KRwEstim2}
 \| \lim_{M \to \infty} \TT_{1,M} (t)\|_1  \le  \frac{1 }{ t}
\frac{1  }{  2\pi} \int_\R   \frac{C_1 }{\langle y \rangle^\alpha} \, dy \| f \|_{2} \to 0,
\eeqn
as $t\to\infty$. Gathering \eqref{eq:KRwEstim1} and  \eqref{eq:KRwEstim2}, we conclude in particular that 
$$
\frac1T \int_0^T S_\LL(t) f_0 \, dt \to 0 \ \hbox{ in } X_1, \ \hbox{ as } \ T \to \infty, 
$$
which is in contradiction with the estimate of Step~2.

 \smallskip\noindent
{\sl  Step 4. Conclusion.}  
Taking advantage of the convenient blow up of $\RR_\LL(\lambda)$ as $\lambda \searrow 0$ established in the previous step, we may now argue similarly as in the proof of Theorem~\ref{theo:exist1-KRexistence}. 
More precisely, from Step~3,   there exists a sequence $(\lambda_n)$ such that $\lambda_n \to 0$ and
$$ \| \RR_{\LL} (\lambda_n) \|_{\BBB(X_2,X_1)} \to \infty. 
$$
That means that there exist $(f_n)$ and $(g_n)$ such that 
$$
 \| f_n \|_{X_1} \to \infty,  \quad  \|   g_n \|_{X_2} = 1, \quad   f_n = \RR_\LL (\lambda_n) \,  g_n,
$$
or equivalently that there exist $(\hat f_n)$ and $(\eps_n)$ (by defining $\hat f_n :=  f_{n\pm}/ \|   f_{n\pm} \|_{X_1}$, $\eps_n := g_{n\pm}/ \|  f_{n\pm} \|_{X_1}$)
satisfying 
\beqn\label{eq:gn=LfnBIS}
 \| \hat f_n \|_{X_1} = 1,  \quad \hat f_n \ge 0, \quad  \| \eps_n \|_{X_2} \to 0 , \quad \eps_n = (\lambda_n - \LL) \, \hat f_n .
\eeqn
As in the proof of   Lemma~\ref{lem:H3abstract-StrongC}, we deduce that \eqref{eq:1hatfn=iterate} holds, that is
\beqn\label{eq:RLN}
\hat f_n = \sum_{\ell=0}^{N-1}  (\RR_\BB(\lambda_n) \AA)^\ell \RR_\BB(\lambda_n) \eps_n +   (\RR_\BB (\lambda_n) \AA )^N \hat f_n.
\eeqn
Using the uniform boundedness
$$
(\RR_\BB(\lambda_n) \AA)^\ell \RR_\BB (\lambda_n)  \in \BB(X_2,X_1), \quad (\RR_\BB(\lambda_n) \AA)^N \in \BB(X_1,\XX_1), \,\, \XX_1 \subset\subset X_1, 
$$
we deduce that $(\hat f_n)$ belongs to a compact set of $X_1$, or in other words, that there exist a subsequence of $(\hat f_n)$ (not relabeled) and 
 $f_1 \in X_1$ such that $\hat f_n \to f_1$ in $X_1$. We may pass to the limit in \eqref{eq:gn=LfnBIS}, and we get \eqref{eq:0=LfBIS}.
 \end{proof}

\subsection{Second existence result  in the weakly dissipative case} 
\label{subsec-Exist2-weakD2}
  
 \
 Using a pure dynamical approach adapted from the second proof of Theorem~\ref{theo:KRexistBIS}  and from \cite[Thm.~6.1]{MR4534707}, 
 we establish a second existence result which is 
   less demanding in terms of conditions on the semigroup $S_\LL$.

\begin{theo}[ ]\label{theo:KRexistTER} 
Consider a positive semigroup $S=S_\LL$ on a Banach lattice   $X = Y'$ for a  separable Banach lattice $Y$. 
We assume 
 \begin{itemize}
\item[(i)] \Cyan there exists  $\phi_0 \in Y_{+}$ such that $[S_t f]_{\phi_0} \ge [ f]_{\phi_0}$ for any $ f \in X_+$  
and such that $f \mapsto  [f ]_{\phi_0}$ is a norm on $X$. \Black
 We then denotes  
$\XX_0$ the vector space $X$ endowed with the norm $[\cdot]_{\phi_0}$; 

 
  \item[(ii)] there exist $v \in L^\infty(\R_+;\BBB(X))$ and $0 \le w \in L^1(\R_+; \BBB(\XX_0,X))$ such that 
\beqn\label{eq:borneKRexistONE}
S = v + w * S, 
\eeqn
and we set 
\beqn\label{eq:borneKRexistBIS}
M := \sup_{t \ge 0} \| v(t) \|_{\BBB(X)} < \infty, \quad 
 \Theta(t) := \| w (t)\|_{\BBB(\XX_0,X)} \in L^1(\R_+). 
\eeqn
\end{itemize}
Then there exists a pair $(\lambda_1,f_1) \in \R_+ \times X_+ \backslash \{ 0 \}$ such that $\LL f_1 = \lambda_1f_1$.  
\end{theo}

\begin{rem}\label{rem:KRexistTER}
(1) When $S_\LL$ satisfies \ref{HS3} then \eqref{eq:borneKRexistONE} holds with 
\beqn\label{eq:KRexistTER-defvw}
v := \sum_{\ell=0}^{N-1}  S_\BB *  (\AA S_\BB)^{(*\ell)}, \quad w :=  ( S_{\BB} \AA)^{(*N)} . 
\eeqn

(2)  By definition of this norm  $[\cdot]_{\phi_0}$, we see that  $\XX_0$ is a weighted $L^1$ space or a weighted Radon measures space. 
In many applications, when both $\XX_0$ and $X$  are Radon measures spaces, one can choose $N=1$. \Black
On the other hand, when $X$ is for instance a (possibly weighted) $L^p$ space with $p>1$,  one must take $N\ge2$ in most of the applications.  
In condition (ii), the first bound is not really demanding and almost automatic in view of the estimates exhibited in Section~\ref{subsect:AboutWeakDissip}. 
The second bound is a kind of regularity estimate reminiscent of the enlarging and shrinkage technique developed in \cite{MR2197542,MR3779780,MR3488535}. 
\end{rem}
 
\begin{proof}[Proof of Theorem~\ref{theo:KRexistTER}]
 We split the proof into three steps. {\Cyan During the proof, we will use the shorthand $[\cdot]_0 := [\cdot]_{\phi_0}$.}

\smallskip
  {\sl Step 1. }   
  We define 
  $$
  R := \max(2\| \Theta \|_{L^1},\| g_0 \|), 
  $$
  for some $g_0 \in X_+$ such that $[g_0]_0 = 1$, and next the nonempty convex and compact (in the weak $* \,  \sigma(X,Y)$ sense) set 
  $$
  \CC := \{ f \in X_+; \, [f]_0 = 1, \ \| f \| \le R \}, 
  $$
  as well as the increasing function
  $$
  \lambda(t) := \inf_{f \in \CC} [S(t) f]_0, \quad \forall \, t \ge 0.
  $$
 We have the alternative 
  
  \smallskip
  $\bullet$ {\bf (1)}  $\sup\lambda > 2M$,

  \smallskip
  $\bullet$ {\bf (2)}  $\sup \lambda \le 2M$.

\medskip 
{\sl Step 2.}  We assume that the first  term {\bf (1)} of the alternative holds true, or in other words, there exists $T_0 > 0$ such that 
\beqn\label{eq:KRexistTER-alternative2}
\forall \, f \in \CC, \quad [S_{T_0} f]_0 \ge 2M.
\eeqn
We define as before 
$$
\Phi_{T_0} f := \frac{S_{T_0}  f }{ [S_{T_0} f ]_0}, \quad \forall \, f \in \CC.
$$
By construction, for any $f \in \CC$, we have $\Phi_{T_0}  f \ge 0$ and $[\Phi_{T_0}  f]_0 = 1$. On the other hand, using the splitting structure   \eqref{eq:borneKRexistONE} and the estimates  \eqref{eq:borneKRexistBIS}, we have 
$$
\| S(t) f \| \le M \| f \| + \int_0^t \Theta(t-s) [S(s) f ]_{0} \, ds. 
$$
From hypothesis (i) and the semigroup property, we also have 
$$
[S_t f ]_0 \ge [S_s f ]_0, \quad \forall \, t \ge s \ge 0.
$$
The two above estimates together imply 
\bean
\| \Phi_{T_0}  f \| 
&\le& \frac{M \|f \| }{  [S_{T_0} f]_0}  + \int_0^{T_0} \Theta(T_0-s) \frac{[S_{s} f]_0 }{  [S_{T_0} f]_0} \, ds  
\\
&\le& \frac{1 }{ 2}   \|   f \| + \| \Theta \|_{L^1} \le R, 
\eean
for any $f \in \CC$. 
We have thus proved $\Phi_{T_0} :   \CC \to \CC $. Thanks to the Tykonov fixed point Theorem, there exists $f_{T_0} \in \CC$ such that $\Phi_{T_0} f_{T_0} = f_{T_0}$. In other words, we have built a pair of ``almost eigenvalue and eigenfunction''
$$
f_{T_0} \ge 0, \quad [f_{T_0}]_0 = 1, \quad S_{T_0} f_{T_0} =  e^{\lambda_1 T_0} f_{T_0},
$$
with $e^{\lambda_1 T_0} = [S_{T_0} f]_0$ and thus $\lambda_1 \in [0,\kappa_1]$. We conclude to the existence of $f_1 \in \CC$ such that $\LL f_1 = \lambda_1 f_1$ really similarly as in Step 2 of the Second proof of Theorem~\ref{theo:KRexistBIS}.

\medskip
{\sl Step 3. } We assume that the  second term  {\bf (2)}  of the alternative holds true. In that case, for any $n \ge 1$, there exists $f_{n} \in \CC$ such that 
$[S(n) f_n]_0 \le 2M$. By compactness, there exists $f_0 \in \CC$ and a subsequence $(f_{n_k})$ such that $f_{n_k} \wto f_0 \in \CC$
and
$$
\forall \, t \ge 0 , \ \forall \, k\, (n_k\ge t), \quad [S(t) f_{n_k}]_0 \le  [S(n_k) f_{n_k}]_0 \le 2M,
$$
so that 
\beqn\label{eq:boundStf0}
\forall \, t \ge 0, \quad [S(t) f_0 ]_0 \le   2M. 
\eeqn
Using this particular initial datum, we argue similarly as in \cite[proof of Thm.~6.1]{MR4534707}, and we conclude to the existence of a stationary state. 
More precisely, we come back to the splitting structure \eqref{eq:borneKRexistONE} of the semigroup $S$ and 
  we introduce the associated Cesàro means
\beqn\label{def:UTVTWT}
  U_T := \frac{1 }{ T} \int_0^T S(t) \, dt, \quad V_T := \frac{1 }{ T}
  \int_0^T v(t)   \, dt, \quad K_T := \frac{1 }{ T} \int_0^T (w *S)(t)  \, dt,
\eeqn
  for any $T > 0$.
  We obviously have
  $$
  \| V_T \|_{\BBB(X)} \le \frac{1 }{ T} \int_0^T  \| v(t) \|_{\BBB(X)}  \, dt \le M.
  $$
  On the other hand, we have 
\bean
\int_0^T (w *S)(t)  \, dt =   \int_0^T  \int_s^{T} w (t-s) dt S (s)  \, ds
\le    \int_0^{T}  w (\tau) d\tau \int_0^T S (s)  \, ds, 
\eean
 thanks to  the Fubini theorem and the  positivity of the two operators involved in this integral
  formula. We deduce 
  \bean \| K_T f_0 \|
  &\le& \Bigl\|   \int_0^{T} w(\tau) d\tau \frac{1 }{ T}  \int_0^T S (s)  \, f_0 \, ds  \Bigr\|
  \\
&\le&  \int_0^{\infty} \bigl\|  w (\tau) \bigl\|_{\BBB(\XX_0,X)} d\tau \, \Bigl[ \frac{1 }{T}  \int_0^T  S (s)  f_0 \, ds  \Bigr]_0
=  \| \Theta \|_{L^1} \,  [ U_T f_0 ]_0,
  \eean 
  thanks to assumption (ii), 
  so  that $K_Tf_0$ is uniformly bounded in
  $X$ thanks to \eqref{eq:boundStf0} and the elementary estimate $ [ U_T f_0 ]_0 \le  [ S_T f_0 ]_0$. 
  We then deduce that $U_T = V_T + K_T$ 
  satisfies 
  $$
  \| U_T f_0\| \le M \|   f_0 \| +  2M \| \Theta \|_{L^1} 
  \quad \hbox{and} \quad 1 \le [S_T f_0]_0 \le 2M,
  $$
  for any $T > 0$. 
By compactness, there exists $T_k \to +\infty$ and $f_1\in X_+$ such that $ U_{T_k} f \wto f_1$   weakly $*$ in $X$. Thanks to the second inequality, we have   $[f_1]_0 \ge 1$. 
We then argue thanks to the usual mean ergodic theorem trick. 
  For any fixed $s > 0$, we observe that 
\bean  
S(s) f_1 -f_1 
  &=& 
  \lim_{k \to \infty} \Bigl\{ \frac{1 }{ T_k}
  \int_0^{T_k} S(s) S(t) f_0 dt - \frac{1 }{ T_k} \int_0^{T_k} S(t)f_0 \, dt \Bigr\}
  \\
  &=& \lim_{k \to \infty} \Bigl\{ \frac{1 }{ T_k} \int_{T_k}^{T_k+s}  S(t) f_0 dt - \frac{1}{T_k} \int_0^s S(t) f_0 \, dt \Bigr\} 
  \eean 
 weakly $*$ in $X$. By the lower semicontinuous property of the norm $[\cdot]_0$, we deduce 
\bean  
[S(s) f_1 -f_1]_0 
   &\le& \liminf_{k \to \infty} \Bigl\{ \frac{1 }{ T_k} \int_{T_k}^{T_k+s}
  [S(t) f_0]_0 dt + \frac{1 }{ T_k} \int_0^s [S(t) f_0]_0 \, dt \Bigr\} = 0, 
  \eean 
 so that $f_1$ is a stationary solution, and thus $f_1$ is an eigenfunction associated to the eigenvalue $\lambda_1 = 0$.
 \end{proof}

{\Blue  Let us emphasize on the fact that Theorem~\ref{theo:KRexistTER} does not implies \ref{S1}, but only provides an existence result of the direct eigenvalue problem. 
In this more involved weakly dissipative framework, in order to tackle the eigentriplet problem, we may proceed in the following way. Assuming for instance that  both  $S_\LL$ and $S_{\LL^*}$ satisfy the conditions of Theorem~\ref{theo:KRexistTER} (an alternative argument combining Theorem~\ref{theo:KRexistTER} with results of Section~\ref{sec:ExistenceKR} will be used in Section~\ref{subsec-diffusionwith drift}), we have

\begin{enumerate}[label={\bf(C1)},itemindent=13mm,leftmargin=0mm,itemsep=1mm]
\item\label{C1} the first primal and dual eigenvalue problems have a solution: there exist $\lambda_1,\lambda_1^* \in \R$, 
  $f_1 \in X_+\cap D(\LL){\setminus\{0\}}$, $\phi_1 \in Y_+ \cap D(\LL^*){ \setminus\{0\}}$ such that 
\beqn\label{eq:Irred-FirstEVpb}
\LL f_1 = \lambda_1 f_1 {\qquad\text{and}\qquad} 
 \LL^* \phi_1 = \lambda_1^*\phi_1.
\eeqn
\end{enumerate}
This conclusion is  weaker than the conclusion~\ref{S1} obtained in Section~\ref{sec:ExistenceKR}, which additionally ensures that $\lambda_1^*=\lambda_1$ and $\Sigma(\LL)\subset\{z\in\C,\ \Re e(z)\leq\lambda_1\}$.
It is however the fundamental existence step. The additional mild geometric properties of the conclusion \ref{S1} may be obtained thanks to the material developed in the next Sections~\ref{sec:Irreducibility} and~\ref{sec:geo2}.
Let us emphasize that in the previous section, $\lambda_1$ is defined thanks to \eqref{eq:exist1-deflambda1}, which because of $\II = \II_2$ in Lemma~\ref{lem:Existe1-H2bis}, automatically implies the inclusion  $\Sigma(\LL)\subset\{z\in\C,\ \Re e(z)\leq\lambda_1\}$.  Here, the construction of $\lambda_1$ is different, and the last inclusion is no longer guaranteed.
}

%

\bigskip

\bigskip 
\section{Irreducibility and geometry of the first eigenvalue} 
\label{sec:Irreducibility}


%


In this section, we start the study of the geometric and asymptotic behavior parts of the Krein-Rutman theorem. 
{\Blue More precisely, we prove the geometric part \ref{S2} and the mean ergodicity part \ref{E1} 
under the partial conclusion \ref{C1} of \ref{S1} and some additional conditions.
In particular we introduce the structure condition \ref{X1} on the Banach lattice, the variant condition \ref{H1'} of \ref{H1} and the strong positivity condition \ref{H4} which is mainly based on a strong maximum principle
 (for a positive semigroup, that last condition being nothing but the  irreducibility).}

%
%

  \medskip
\subsection{More about positivity}  
 
\label{subsec:MorePositive}
  
 For further references, we introduce several notions which are strongly related to the positivity property for semigroups. 
 
 \medskip
{\bf  The signum operator ${\rm sign}$.} 
In a real Banach lattice $X$,  we say that $\sign f \in \BBB(X,X'')$  is a signum  operator for $f \in X$, if it satisfies the following properties
\bear \label{eq:prop_signff}
&& (\hbox{sign} \, f  ) \, f =    |f|, 
\\ \label{eq:prop_signfg}
&& (\hbox{sign} \, f  ) \, g    \le |g|, \quad \forall \, g \in X.
\eear
In the sequel, we will always assume that such an operator exists. 
We refer to \cite[Sec.~C.I \& C.II]{MR839450} for a general introduction to the topic. In practice, we will only need a weak formulation of the  ${\rm sign}$ operator (see below) which may be 
defined only  in some subspace $\XX \subset X$. We  always additionally assume that the signum  operator  satisfies 
\bean
&& (\hbox{sign} \,(-f ) ) \, (-g) = (\sign f) g, \quad \forall \, g \in X,  
\\
&& (\hbox{sign} \, f  ) \, g =  g, \quad \forall \, g \in X, \ \hbox{if} \ f \in X_+, 
\eean
We also define 
$$
\hbox{sign}_+ f := \frac12 \bigl( I + \hbox{sign} f \bigr).
$$
$\bullet$ When $X$ is a space of functions, the sign operator $\sign f$ associated to $f \in X$ corresponds to the multiplication by the function $\sign f := {\bf 1}_{f > 0} - {\bf 1}_{f < 0}$. 
When $X := L^p(E)$, we obviously see that $\sign f \in \BBB(L^p(E))$ for any $f \in L^p(E)$. On the other hand, when $X := C_0(E)$, we only have $\sign f \in \BBB(C_0(E);\MM^\infty(E))$, where $\MM^\infty(E)$ denotes the space of uniformly bounded measurable functions, so that $\MM^\infty(E) \subset (C_0(E))''$. 
In a space of bounded measures $X=M^1(E)$, we may define the sign operator by means of the Radon-Nikodym theorem. For a given $f \in M^1(E)$, using Hahn decomposition,  there exists indeed a measurable function $\alpha: E \to \{-1,1\}$ such that $f  = \alpha |f| $, and we then define $(\sign f) g = \alpha g$ for any $g \in M^1(E)$. 
\Black

\smallskip
$\bullet$
When  $X$ is  $\sigma$-order complete, in the sense that any increasing and upper bounded sequence has a supremum (a common least upper bound), 
the operator $\sign$ exists and is more regular, namely  $\sign f \in \BBB(X)$ for any $f \in X$, see \cite{MR637005} and also \cite[Sec.~C.I.8]{MR839450}. We recover in particular that  $\sign f \in \BBB(L^p(E))$ for any $f \in L^p(E)$. %

 \medskip
{\bf  Weak maximum principle and Kato's inequality.} We introduce now two definitions formulated on an operator $\LL$ which are almost equivalent to the positivity property of the  semigroup $S$ when $\LL$ is the generator of $S$. 

\smallskip
$\bullet$ We  say that the operator $\LL$ satisfies the {\it weak maximum principle} when 
\beqn\label{eq:PM*}
\kappa \in \R,  \ f \in D(\LL) \hbox{ and } ( \kappa - \LL) f \ge 0 \quad\hbox{imply}\quad f \ge 0; 
\eeqn

\smallskip
$\bullet$ We say that the operator $\LL$  satisfies {\it Kato's inequality} when   
\beqn\label{eq:KatoIneq}
\forall \, f \in D(\LL), \quad \LL |f| \ge (\sign f)  \LL  f. 
\eeqn
Since $|f|$ does not necessarily belong to $D(\LL)$, the correct way to understand Kato's inequality is 
\beqn\label{eq:KatoIneqDualsense}
\forall \, f \in D(\LL), \,\,\,  \forall \, \psi \in D(\LL^*) \cap Y_+,  \quad  \langle |f|, \LL^* \psi \rangle  \ge \langle  (\sign f)  \LL  f,\psi \rangle.
\eeqn
We immediately see from the definitions that \eqref{eq:KatoIneq} is equivalent to assume
\beqn\label{eq:Kato+Ineq}
\forall \, f \in D(\LL), \quad \LL f_+  \ge (\sign_+  f) \LL  f. 
\eeqn

\begin{rem}\label{rem:sec4-EquivPositiveS}
We complement Lemma~\ref{lem:sec2-EquivPositiveS}, by claiming that for a semigroup $S = S_\LL$ on a Banach lattice $X$, there is equivalence between
the fact that $S$ is positive and  $\kappa - \LL$  satisfies the  weak maximum principle for any $\kappa > \omega(\LL)$, what is straightforward using that these properties are equivalent 
to the fact that $\RR_\LL(\kappa) \ge 0$ for any $\kappa > \omega(\LL)$. These properties also imply that 
Kato's inequality holds true, see \cite{MR637005,MR763347}, \cite[Prop.~1.1]{MR790308}, \cite[Rk.~3.10]{MR683043} and the textbook \cite[Thms  C.II.2.4, C.II.2.6 and Rk. C-II.3.12]{MR839450}. 
When for instance $f,|f| \in D(\LL)$, we may indeed compute 
$$
 (\hbox{\rm sign} \, f  ) \LL f =  (\hbox{\rm sign} \, f  )  \lim_{t \to 0} \frac{S_t f - f }{ t} \le \lim_{t \to 0} \frac{S_t |f| - |f| }{ t} = \LL |f|, 
 $$
 where we have used the very definition of the generator $\LL$ and the properties \eqref{eq:prop_signff}-\eqref{eq:prop_signfg} of  $\hbox{\rm sign} \, f $ in the inequality.

\end{rem}

\medskip
 We end this section by introducing other notions of positivity which strengthen the previously defined positivity condition.

\medskip
{\bf Strict order. } 
We may define a first  stronger order $>$ (or $<$) on $X$ by writing for $f \in X$ 
$$
f >  0  \quad\hbox{if}\quad f \in X_+ \backslash \{0\}
$$
and similarly a stronger order $>$ (or $<$) on $X'$ by writing for $\phi \in X'$ 
$$
\phi >  0 \quad\hbox{if}\quad \phi  \in X'_+ \backslash \{0\}.
$$
We may next define the strict (and stronger) order $\gg$ (or $\ll$) on $X$ by writing for $f \in X$ 
\begin{equation}\label{eq:strict-order}
f \gg  0 \hbox{ or } f \in X_{++} \quad\hbox{iff}\quad
\forall \, \psi \in X'_+ \backslash\{0\}, \,\,\, \langle \psi,f \rangle > 0 , 
\end{equation}
and similarly the strict order $\gg$ (or $\ll$) on $X'$ by writing for $\phi \in X'$ 
\begin{equation}\label{eq:dual-strict-order}
\phi \gg  0 \hbox{ or } \phi \in X'_{++}  \quad\hbox{iff}\quad
\forall \, g \in X_+ \backslash\{0\}, \,\,\, \langle \phi,g \rangle > 0 .
\end{equation}

On the two Banach lattices $X$ and $Y$, we  thus have three positivity notions with $\gg$ (associated to $X_{++}$ and $Y_{++}$) stronger than $>$ (associated to $X_{+}\backslash \{0\}$ and $Y_+ \backslash \{0\}$)  which itself is stronger than $\ge$ (associated to $X_{+}$ and $Y_+$). {\Cyan We also emphasize that $\psi \in Y_{++}$ 
 if, and only if, $[\cdot]_\psi$ is a norm on~$X$.} 

\smallskip
Let us comment on the notion of strict positivity. 

\begin{rem}\label{rem:IrredStrictPositivity}
When $X = Y'$ for instance, there are two possible strict positivity notions on $X$ given by~\eqref{eq:strict-order} for the space $X$ (namely $\phi \in X$ is strictly positive when $\langle \xi, \psi \rangle > 0$ for any  
$ \xi\in X'_+ \backslash\{0\}$) 
and by~\eqref{eq:dual-strict-order} for the space $Y$ (namely $\phi \in X$ is strictly positive when $\langle \psi,g \rangle > 0$ for any  $g\in Y_+ \backslash\{0\}$).
They clearly coincide when $X$ is reflexive, but in general the first one is stronger than the second one. In that situation,   we will always consider that $X$ is endowed with the weakest ``dual'' strict order~\eqref{eq:dual-strict-order}.
\end{rem}

\begin{exples}\label{ex:Irred:StrictOrderNorm}
In the space $C_0(E)$, the strict order   is defined by $f \gg 0$ iff $f(x) > 0$ for any $x \in E$. 
In a space $L^p(E,\EE,\mu)$, $1 \le p \le \infty$, the strict order  is defined by $f \gg 0$ iff $f(x) > 0$ for $\mu$-a.e. $x \in E$. 
In the space $M^1(E)=C_0(E)'$, the strict order is defined by duality by $f \gg 0$ iff  
$\langle f, \varphi \rangle >  0$ for any $\varphi \in C_0(E)$, $\varphi \ge 0$, $\varphi \not\equiv 0$. 
\end{exples}

\begin{rem}
In a Banach lattice $X$ such that $\hbox{\rm int}\, X_+ \not=\emptyset$, the common definition of the  strict order is $X_{++} := \hbox{\rm int}\, X_+$. 
In particular, in the case when $E$ is compact and $X = C_0(E) = C(E)$, we have $\hbox{\rm int}\, X_+ \not=\emptyset$ and the definition of $X_{++}$ introduced in Examples~\ref{ex:Irred:StrictOrderNorm}
coincides with $\hbox{\rm int}\, X_+$. 
In all the other examples considered, we have  $\hbox{\rm int}\,X_+ = \emptyset$, and thus our definition of the strict order does not coincide with the one defined through the set $\hbox{\rm int}\,X_+$. 
\end{rem}

\begin{rem}\label{rem:Ideal}
Another notion of strict order can be defined through the notions of ideals and quasi-interior points as briefly explained now, see \cite{MR839450} or \cite[Chap.~10]{MR3616245} and the references therein for details. 
Defining the segment $[g_1,g_2]$ and the set $I_f$ for $g_1,g_2 \in X$ and $f \in X_+ \backslash \{0\}$ by 
$$
[g_1,g_2] := \{ g \in X; \ g_1 \le g \le g_2 \}, \quad I_f := \bigcup_{k \ge 0} [-k f, kf]=\Span[0,f], 
$$
one shows that $I_f$ is an ideal in the sense that  $g \in I_f$ implies $|g| \in I_f$ and $0 \le g \le f$ implies $g \in I_f$. 
 We say that $f$ is an {\it order unit} if $I_f = X$. 
When $\hbox{\rm int}\, X_+ \not=\emptyset$, we find that  $f$ is an {\it order unit} iff $f \in \hbox{\rm int}\, X_+$ from Lemma~\ref{lem:exists-lemfonda2}, so that we recover the notion of strict positivity defined above. 
On the other hand, we say that $f$ is a {\it quasi-interior point} if $\bar I_f = X$.
It can be shown that $f$ is a {\it quasi-interior point} iff $f$ is strictly positive in the sense of the direct strict order \eqref{eq:strict-order},  see for instance~\cite[Thm. II.6.3]{MR0423039}. 
These two notions of strict positivity and quasi interior point thus coincide when $X$ is reflexive or when $X = L^p(E,\EE,\mu)$, $1 \le p < \infty$, see also \cite[Examples~10.16]{MR3616245} when $\mu$ is a $\sigma$-finite diffuse (or atomless) measure. 
On the other hand, it is important to point out again that the ``dual'' strict order~\eqref{eq:dual-strict-order} considered here is a weaker notion than the {\it quasi-interior point} notion. 
For instance, in $X=C_0(E)'=M^1(E)$, there is no quasi-interior point but $X_{++}\neq\emptyset$.

%

\end{rem}

 
We finally point out the following result. For a semigroup $S = S_\LL$ in a Banach lattice, under the mild assumption that there exists a strictly positive subeigenvector for the dual problem, namely
$$
\exists \, \phi \in X'_{++}, \ \exists \, b \in \R, \quad \LL^* \phi \le b \, \phi, 
$$
then Kato's inequality \eqref{eq:KatoIneq} implies that $S$ is positive, see \cite[Thm.~1.6]{MR790308}.

  \medskip
\subsection{Irreducibility  and strong maximum principle}
\label{subsec:irreducibility}
We present some other material involving the strict positivity.  When satisfied by {\Blue  $ \EEE   \subset X$ or  $\EEE \subset  Y$}, we will in particular make use of the property  
\beqn\label{eq:Irred:StrictOrderNorm}
\Blue (\xi \in \EEE, \ \xi_+ \gg0) \quad \hbox{implies} \quad \xi \gg0 \quad \hbox{(or equivalently $\xi_-=0$)}.
\eeqn
\Blue The property \eqref{eq:Irred:StrictOrderNorm} does not always hold in $\EEE = X$, but it holds  in many usual cases. 
\smallskip

%
%
%

 \Black


\begin{lem}\label{lem:Irred:f+>0&f>0}
The property~\eqref{eq:Irred:StrictOrderNorm} holds true in a space $X$  endowed with the direct strict order~\eqref{eq:strict-order}, in particular in $X=  \Blue L^p_m$, $p\in[1,\infty)$, and $X= \Blue C_{0,m}$. 
\Blue The property~\eqref{eq:Irred:StrictOrderNorm} also holds  in a set $\EEE \subset X \subset \Mloc^1$ endowed with the dual strict order~\eqref{eq:dual-strict-order} when $\EEE \subset \Lloc^1$, and in particular  in the space $L^\infty_{m^{-1}}=(L^1_m)'$. 
\end{lem}

\begin{proof}[Proof of Lemma~\ref{lem:Irred:f+>0&f>0}] 
We start recalling the proof of \eqref{eq:Irred:StrictOrderNorm} in a general  space $X$  endowed with the strict order~\eqref{eq:strict-order}. 
Consider an element $f$ of the Banach lattice $X$ and assume that $f_+ \gg0$.
The vectors $f_+$ and $f_-$ are disjoint, {\it i.e.} $f_+\land f_-=0$, see for instance~\cite[Thm.~1.1.1 iv)]{MR1128093}.
On the one hand, since $f_+\gg0$, we have 
 $$(n f_+) \land f_- \to f_-$$
with respect to the norm  as $n \to \infty$, see~\cite[Thm.~II.6.3]{MR0423039}.
On the other hand, we have
 $$ 0 \leq (n f_+) \land f_- \leq  (n f_+) \land (n f_-) =  n (f_+ \land f_-) =  0, 
 $$
for every integer $n \geq 1$, where the last equality follows from the fact that $f_+$ and $f_-$ are disjoint.
We deduce by passing to the limit $n\to\infty$ that $f_- = 0$. Thus, $f = f_+ \gg0$.

\smallskip
We now establish  \eqref{eq:Irred:StrictOrderNorm} in the second case. Take thus $f \in \EEE \subset \Lloc^1 \cap X$, with $X \subset \Mloc^1$ endowed with the dual order, so that  $\langle f_+ , \varphi \rangle > 0$ 
for any $\varphi \in C_c$, $\varphi \ge 0$, $\varphi \not\equiv0$. Because $f$ (and thus $f_+$) is a measurable function, we have 
$f_+(x) = \max(f(x),0) > 0$ for a.e. $x \in E$, so that $f > 0$ a.e. on $E$ and finally $f_- =0$ a.e. on $E$.
\end{proof}

 We give now  a counter-example in the Radon measures space case.

\begin{exple}\label{cex:Irred:f+>0&f>0dual}
Consider $M^1([0,1])= C([0,1])'$ endowed with the dual strict order~\eqref{eq:dual-strict-order}.
  Let $(q_n)$ be an enumeration of the rational numbers in $[0,1]$ and let $r$ be an irrational number in $[0,1]$.
  The functional $\phi$ given by
  $\langle\phi, f\rangle := \sum_{n=1}^\infty 1/2^n f(q_n)  -  f(r)$
  satisfies $\phi_+ \gg0$, but $\phi$ itself is not positive.
  \end{exple}

 \smallskip

For an operator $A \in \BBB(X)$, we have yet formalized a positivity condition in section~\ref{subsec:BanachLattice}, by

(P1) $A \ge 0$ if $A : X_+ \to X_+$.
 
\smallskip
Other possible definition of  positivity may be

(P2)  $A : X_+\backslash\{0\} \to X_+\backslash\{0\}$;

(P3)  $A : X_{++} \to X_{++}$.

\smallskip
We now define a stronger notion of positivity, named as strong positivity condition, 
 by 
 
(P4) $A > 0$ if $A : X_+\backslash\{0\} \to   X_{++}$.

\smallskip
We list without proof some elementary properties about these different notions and also refer to Section~\ref{sec:quantified-positivity}
for further discussion.
We have  (P2) $\Rightarrow$ (P1),   (P3) $\Rightarrow$ (P1)  as well as
(P4) $\Rightarrow$ ((P3), (P2)).
We also have 
$A : X_+ \to X_+$ iff $A^* : Y_+ \to Y_+$;
$A :  X_{++}  \to  X_{++} $ iff $A^* :  Y_{++}  \to  Y_{++} $;
$A : X_+\backslash\{0\} \to  X_{++} $ iff $A^* : Y_+\backslash\{0\} \to  Y_{++} $.

\medskip
We say that $\lambda-\LL$ satisfies the strong maximum principle if 
\Black
\beqn\label{eq:Irred-StrongMP}
 f \in X_+ \cap D(\LL), \ (\lambda-\LL) f \ge 0 \quad\hbox{imply} \quad  f \gg  0  \ \hbox{or} \ f = 0.
\eeqn
  It is worth emphasizing that if $\lambda-\LL$ satisfies the strong maximum principle for some $\lambda \in \R$ then  $\lambda'-\LL$ satisfies the strong maximum principle for any $\lambda' \le \lambda$.  

\smallskip
We say that a positive semigroup $S$ is irreducible if 
\beqn\label{eq:Irred-Sirreducible}
\forall \, f \in X_+\backslash \{0\}, \forall \, \phi \in Y_+\backslash \{0\}, \ \exists \tau >0 \quad
\langle S_\tau f, \phi \rangle > 0.
\eeqn
 
A semigroup $S$ is classically said to be 
irreducible and aperiodic if the above positivity condition holds for all sufficiently large times, namely
\beqn\label{eq:Irred-eventualirreducibilityBIS}
\forall \, f \in X_+\backslash \{0\}, \forall \, \phi \in Y_+\backslash \{0\}, \ \exists\, T>0, \,\forall \tau \geq T \quad
\langle S_\tau f, \phi \rangle > 0.
\eeqn
\Black
Other notions of strong positivity for the semigroup $S$ are 
\bear\label{eq:Irred-STstrictPositive1}
&&\exists \, T > 0, \quad S_T : X_+\backslash\{0\} \to  X_{++}, 
\\ \label{eq:Irred-STstrictPositive2}
&&\exists\, T > 0, \quad \int_0^T S(t) dt : X_+\backslash\{0\} \to  X_{++}.
\eear

We summarize some possible implications between the previous positivity notions. 

\begin{lem}\label{lem:Irred-S>0impliesR>0}
For a positive semigroup $S$, the following hold: 

{\bf (1)} The {\Cyan given time} 
strong positivity condition \eqref{eq:Irred-STstrictPositive1} implies the {\Cyan integral in time strong positivity} condition \eqref{eq:Irred-STstrictPositive2};  

{\bf (2)}  The integral strong positivity condition \eqref{eq:Irred-STstrictPositive2} implies the  irreducibility condition \eqref{eq:Irred-Sirreducible}, but the reverse implication is false.
Similarly, the irreducibility and aperiodicity condition \eqref{eq:Irred-eventualirreducibilityBIS} implies the  irreducibility condition \eqref{eq:Irred-Sirreducible}, but the reverse implication is false; 
 
{\bf (3)} 
  The irreducibility condition \eqref{eq:Irred-Sirreducible} is equivalent to  the fact that $\RR_\LL(\lambda) : X_+\backslash\{0\} \to  X_{++}$, for any $\lambda > \lambda_1$,
as well as to   the fact that $\lambda-\LL$ satisfies the strong maximum principle \eqref{eq:Irred-StrongMP} for any $\lambda \in\R$.

 \end{lem}

The result is very classic, at least for strongly positive semigroup, see e.g. \cite[Definition~C.3.1]{MR839450} or \cite[Prop.~14.10]{MR3616245}. For the sake of completeness, we however present a short proof. 

\medskip
\begin{proof}[Proof of Lemma~\ref{lem:Irred-S>0impliesR>0}]
We prove {\bf (1)}. We assume \eqref{eq:Irred-STstrictPositive1} and we fix $g \in X_+\backslash \{0\}$,  $\phi \in Y_+\backslash \{0\}$, so that $ \langle S(T) g, \phi\rangle > 0$. 
Observing that the  function $t \mapsto  \langle S(t) g, \phi\rangle$ is continuous, there exists  $\eps > 0$ such that $\langle S(t) g, \phi\rangle > 0$ for any $t \in [T-\eps,T]$, so that 
$$
\Bigl\langle \int_0^T S(t) dt g, \phi \Bigr\rangle = \int_0^T \langle S(t) g, \phi\rangle dt  > 0. 
$$
Because $\phi \in Y_+\backslash \{0\}$ may be chosen arbitrary, we deduce \eqref{eq:Irred-STstrictPositive2}. 

\smallskip

We prove {\bf (2)}.  We assume now \eqref{eq:Irred-STstrictPositive2} and we fix $g \in X_+\backslash \{0\}$,  $\phi \in Y_+\backslash \{0\}$, so that 
$$
\int_0^T \langle S(t) g, \phi\rangle dt = \Bigl\langle \int_0^T S(t) dt g, \phi \Bigr\rangle > 0, 
$$
by assumption. We get \eqref{eq:Irred-Sirreducible} by observing that the  function $t \mapsto  \langle S(t) g, \phi\rangle$ must be positive somewhere on $[0,T]$. 
For the reverse implication we refer to \cite{MR3928121,GabrielMartin}, where  is studied an 
example of growth-fragmentation operator associated to mitosis satisfying the irreducibility condition \eqref{eq:Irred-Sirreducible} but not the  integral strong positivity condition \eqref{eq:Irred-STstrictPositive2} nor the irreducibility and aperiodicity condition~\eqref{eq:Irred-eventualirreducibilityBIS},
see also Section~\ref{part:application3:GF}.

\smallskip
We prove {\bf (3)}. 
We finally assume  \eqref{eq:Irred-Sirreducible}. From the above continuity argument, for any $g \in X_+\backslash \{0\}$,  $\phi \in Y_+\backslash \{0\}$ there exist $\tau > \eps > 0$ such that $\langle S(t) g, \phi\rangle > 0$ for any $t \in [\tau-\eps,\tau+\eps]$. 
As a consequence and thanks to the representation formula \eqref{eq:Exist1-DefRepresentationRR} for any fixed $\lambda > \lambda_1$ which holds thanks to Lemma~\ref{lem:Exist1-RkSG}-{\bf (ii)}, we have 
$$
\langle \RR_\LL(\lambda) g,\phi \rangle = \Bigl\langle \int_0^\infty e^{-\lambda t } S(t) dt g, \phi \Bigr\rangle  > 0.  
$$
 Because $\phi \in Y_+\backslash \{0\}$ is arbitrary, we have established that $\RR_\LL(\lambda) g \in X_{++}$ for any $g \in X_+\backslash \{0\}$. 
  In other words, when $\lambda > \lambda_1$ and $f \in X_+ \cap D(\LL)$ satisfy $g := (\lambda-\LL) f \ge 0$, we deduce that $f = \RR_\LL(\lambda) g \in X_{++}$, what is the strong maximum principle. This one is obviously equivalent to the strong positivity property $\RR_\LL(\lambda) : X_+\backslash\{0\} \to  X_{++}$. On the other way round, writing the above identity as
 $$
 \int_0^\infty e^{-\lambda t } \bigl\langle S(t) g, \phi \bigr\rangle dt =  \langle \RR_\LL(\lambda) g,\phi \rangle, 
$$
we see that the strong maximum principle implies that the RHS term is positive for any $g \in X_+\backslash \{0\}$,  $\phi \in Y_+\backslash \{0\}$. As a consequence, the  LHS term is positive  and there exists $\tau > 0$ such that 
$\bigl\langle S(\tau) g, \phi \bigr\rangle > 0$, which is nothing but the irreducibility condition \eqref{eq:Irred-Sirreducible}. 
\end{proof}

\medskip
We present two other elementary results about the  strong maximum principle. 

 \begin{lem}\label{lem:Irred-KStrongMax&dual} Consider $\LL$ satisfying \ref{H1}  and $\lambda \in \R$. Then the following assertions are equivalent 

(1)  $\lambda-\LL$ satisfies the strong maximum principle for any $f \in D(\LL) \cap X_+$;

(2)  $\lambda-\LL$ satisfies the strong maximum principle for any $f \in D(\LL^k) \cap X_+$ for some $k \ge 1$;

(3)  $\lambda-\LL^*$ satisfies the strong maximum principle for any $\phi \in D(\LL^*) \cap Y_+$; 

(4)  $\lambda-\LL^*$ satisfies the strong maximum principle for  any $\phi \in D((\LL^*)^\ell) \cap Y_+$ for some $\ell \ge 1$.

\end{lem}

\begin{proof}[Proof of Lemma~\ref{lem:Irred-KStrongMax&dual}.]
Assume that  $\lambda-\LL$ satisfies the strong maximum principle for some $\lambda \in \R$ and $k \ge 1$ and consider $\phi \in D(\LL^*) \cap Y_+ \backslash \{ 0 \}$ such that $ (\lambda - \LL^*) \phi \ge 0$.
For any $\kappa > \max(\lambda,\lambda_1)$ 
and $g \in  D(\LL^{k-1}) \cap X_+ \backslash \{0\}$, thanks to \ref{H1} and the strong maximum principle,   there exists $f \in D(\LL^k) \cap X_{++}$  such that 
$
(\kappa-\LL) f = g.
$
As a consequence, we have 
\bean 
\langle \phi, g \rangle &=& \langle \phi, (\kappa-\LL) f   \rangle 
\\
&=&  \langle  (\kappa-\LL^*) \phi, f   \rangle \ge (\kappa - \lambda) \, \langle  \phi, f   \rangle  > 0. 
\eean
Since $g \in D(\LL^{k-1}) \cap X_+$ is arbitrary and $D(\LL^{k-1}) \cap X_+$ is dense in $X_+$, we deduce that $\phi \gg 0$.  We have proved that $\lambda-\LL^*$ satisfies the strong maximum principle.
The other implications can be proved similarly.
\end{proof} 

\begin{rem}\label{rem:Irred-KStrongMax&dual} 
  (1) In many applications, we start proving the strong maximum principle on smooth enough functions (belonging to the iterated domain) for which pointwise arguments may be used.

(2)
We may replace the condition (1) by assuming that  $\lambda-\LL$ satisfies the strong maximum principle for $f \in \CC \cap X_+$ for a subspace $\CC \subset D(\LL)$ such that $(\lambda-\LL)^{-1} \in \BBB(\CC)$ and $\CC$ is dense in $X$. 
\end{rem}

The strong maximum principle can be seen as a consequence of the weak maximum principle together with the existence of
a family of strictly positive barrier functions. We give now a typical result which can be applied (or modified in order to be applied) in many situations.

\begin{lem}\label{lem:Irred-barrier+W->S} We assume that 

\quad (i) the operator $\lambda-\LL$ satisfies the weak maximum principle and Kato's inequalities; 

  \quad (ii) there exists a subset $\GGG \subset  X_{++} \cap \{ g \in D(\LL); \, (\LL - \lambda) g \ge 0 \}$ such that $\forall \, f \in D(\LL) \cap X_+ \backslash \{0\}$, $\exists \, g \in \GGG $ such that $(g-f)_+ \in D(\LL)$.  
   
\smallskip\noindent
Then $\lambda-\LL$ satisfies the strong maximum principle. 
\end{lem}

\begin{proof}[Proof of Lemma~\ref{lem:Irred-barrier+W->S}.]
We consider $f \in D(\LL) \cap X_+ \backslash \{0\}$ such that $(\lambda-\LL) f \ge 0$ and choose $ g \in \GGG $ such that $h := (g-f)_+ \in D(\LL)$. 
We remark that from Kato's inequality \eqref{eq:Kato+Ineq}, we have
$$
(\LL - \lambda)h \ge \hbox{\rm sign}_+ (g-f)
( \LL- \lambda) (g-f) \ge 0.
$$
As a consequence of the weak maximum principle, we have $h \le 0$. That implies $h=0$, so that $g-f \le 0$ and finally $f \gg 0$.
\end{proof}

The above  barrier functions technique is also useful for obtaining the condition \ref{H2} (possibly in a constructive way).  

\begin{lem}\label{lem:Irred-barrier+W->H2} For an operator $\LL$, we assume that 

\quad (i) the condition \ref{H1} holds with a constant $\kappa_1 \in \R$;  

\quad (ii) the hypothesis (ii) in Lemma~\ref{lem:Irred-barrier+W->S} holds  with $\lambda = \kappa_1$; 

\quad (iii) there exists $h_0 \in X_+\backslash \{ 0 \}$ such that for any $g \in \GGG$ there exists $\eps >0$ such that $g \ge \eps h_0$.  

\smallskip\noindent
Then the property \ref{H2}  holds true. 
\end{lem}

\begin{proof}[Proof of Lemma~\ref{lem:Irred-barrier+W->H2}] 
Thanks to assumption (i), we may  define $f_0 \in D(\LL) \cap X_+ \backslash \{0\}$ as the solution to the equation  $(\kappa_1-\LL) f_0 = h_0$.
From the proof of Lemma~\ref{lem:Irred-barrier+W->S} and condition (iii), there exists $g \in \GGG$ and next $\eps > 0$ such that $f_0 \ge g \ge \eps h_0$.
Coming back to the equation, we have 
$$
\LL f_0 = \kappa_1 f_0 - h_0 \ge  (\kappa_1 - \eps^{-1}) f_0,
$$
so that condition \ref{H2} holds true with $\kappa_0 := \kappa_1 - \eps^{-1}$ thanks to Lemma~\ref{lem:Existe1-Spectre2bis}-{\bf (ii)}.
\end{proof}


\subsection{The geometry of the first eigenvalue problem}
  \label{subsec:geo1}
  We come back on and state a result about the geometry of  the first eigenvalue. 
\Blue  
For that purpose, we introduce a general framework for the  couple of Banach lattices we will use in the sequel: 

 \begin{enumerate}[label={\bf(X1)},itemindent=13mm,leftmargin=0mm,itemsep=1mm]
\item\label{X1} $X_{++} \not=\emptyset$, $Y_{++} \not=\emptyset$ and the  signum operator is well define in both $X$ and $Y$. 
\end{enumerate}

\Black

  
%
   
  \smallskip
{\Blue  We consider an operator   $\LL$ on $X$  which satisfies the conclusion {\bf (C1)} about the existence of  solutions $(\lambda_1, f_1)$ and $(\lambda_1^*, \phi_1)$ to the direct and dual first eigenvalue problems.}
 We next assume
  
\begin{enumerate}[label={\bf(H1$'$)},itemindent=14mm,leftmargin=0mm,itemsep=1mm]
\item\label{H1'} the  {\it weak   maximum principle}
\beqn\label{eq:Irred-weakPM}
\lambda > \lambda_1, \ f \in D(\LL), \ (\lambda-\LL) f \ge 0 \quad\hbox{imply} \quad f \ge 0
\eeqn
and  its {\it Kato's inequalities} counterpart
\beqn\label{eq:Irred-Katoineq}
(\hbox{sign} f) \LL f  \le \LL |f|,
\qquad
(\hbox{sign}_+ f) \LL f    \le \LL f_+.
\eeqn
\end{enumerate}
{\Blue We recall that these two conditions are automatically satisfied for both $\LL$ and $\LL^*$ when $\LL$ is the generator of a positive semigroup, see Remark~\ref{rem:sec4-EquivPositiveS}.}

  \smallskip\Blue
  Then we define the set of eigenvectors 
  $$
  \EEE := \{ f \in X\backslash \{ 0 \}; \ \exists \, \lambda \in \C, \  \LL f = \lambda f \}, 
\quad
 \EEE^* := \{ \phi \in  Y\backslash \{ 0 \}; \ \exists \, \lambda \in \C, \  \LL^* \phi = \lambda \phi \}, 
  $$
and we additionally assume that 

\begin{enumerate}[label={\bf(H4)},itemindent=13mm,leftmargin=0mm,itemsep=1mm]
 \item\label{H4}$\left\{\begin{array}{l}
 \text{\bf (i)}\quad  \LL \ \text{and} \ \LL^* \ \text{satisfy the {\it strong maximum principle} \eqref{eq:Irred-StrongMP};}\vspace{2mm}\\
\text{\bf (ii)}\quad \text{the property~\eqref{eq:Irred:StrictOrderNorm} holds  in both sets $\EEE$ and $\EEE^*$.}
\end{array}\right.$
 \end{enumerate}
 
 Observe that \ref{H4}{\bf-(ii)} is verified when $\EEE$ and $\EEE^*$ are subsets of $\Lloc^1$, due to Lemma~\ref{lem:Irred:f+>0&f>0},
 so in particular when $D(\LL^k),D((\LL^*)^k)\subset \Lloc^1$ for some $k\in\N$.
 
 \Black

\smallskip
We may then state our main result in this section,  where we recall that $N(A)$ denotes the null space associated to the operator $A$.

 \begin{theo}\label{theo:KRgeometry1} 
We  assume that $X$ and $Y$ are Banach lattices satisfying  \ref{X1}. We consider an unbounded operator $\LL$ on $X$ which satisfies the conclusion \ref{C1} {\Cyan about the existence part 
of the first eigenvalue problem}. We also assume that $\LL$ and $\LL^*$ both satisfy the  weak maximum principle and Kato's inequalities  \ref{H1'} as well as the {\Cyan strong positivity condition 
\ref{H4}. }

   Then the following hold
 
 i) $f_1 \gg 0$, $\phi_1 \gg 0$ and $\Blue\lambda_1=\lambda_1^*$ is the unique eigenvalue associated to a positive eigenvector.  We next make the normalization choice  
\beqn\label{eq:normalize-f1phi1}
 \| \phi_1 \| = 1, \quad \langle \phi_1, f_1 \rangle = 1.
\eeqn
 
 ii)   $\lambda_1$ is algebraically simple:
\bear \label{eq:theoKRgeometry1-L}
&&N((\LL-\lambda_1)^k) = \hbox{\rm Span}(f_1), \quad \forall \, k \ge 1,
\\ \label{eq:theoKRgeometry1-L*}
&&N((\LL^*-\lambda_1)^k) = \hbox{\rm Span}(\phi_1), \quad \forall \, k \ge 1.
\eear
 In particular $f_1$ (resp. $\phi_1$) is the unique positive and normalized eigenvector of $\LL$ (resp. of $\LL^*$) associated to $\lambda_1$. 
 Finally, the projection  on the first eigenspace (associated to $\lambda_1$) is given by
$$
\Pi\, f := \langle f, \phi_1 \rangle f_1. 
$$

%
%

 \end{theo}
 

 \begin{rem}\label{rem:KRgeometry1}
 
(1)
It is worth emphasizing again that   \eqref{eq:Irred-weakPM} and  \eqref{eq:Irred-Katoineq} for both $\LL$ and $\LL^*$ are   true when $\LL$ is the generator of a
positive semigroup (see Lemma~\ref{lem:sec2-EquivPositiveS} and Remark~\ref{rem:sec4-EquivPositiveS}) 
   and that \ref{H4}{\bf-(i)} is true when $S_\LL$ enjoys additional strong positivity (or irreducibility) condition as formulated in 
\eqref{eq:Irred-Sirreducible}, \eqref{eq:Irred-eventualirreducibilityBIS},  \eqref{eq:Irred-STstrictPositive1} or \eqref{eq:Irred-STstrictPositive2}. 
As a consequence, the conclusion of Theorem~\ref{theo:KRgeometry1} holds true when $\LL$ is the generator of a positive semigroup which satisfies the hypotheses of the existence part of the Krein-Rutman Theorem~\ref{theo:exist1-KRexistence} and one of the additional above strict positivity conditions.

\smallskip
(2) Theorem~\ref{theo:KRgeometry1} has to be compared with the seminal  Krein and Rutman Theorem~\ref{theo:intro-KRresult} (\cite{MR0027128}), to the many results gathered in \cite[Part C-III]{MR839450} (see in particular \cite[Prop.~C.3.5]{MR839450}, \cite[Thm.~C.3.8]{MR839450} and the original paper \cite{MR618205}) and to the more recent contributions   \cite[Thm.~5.3]{MR3489637}, \cite[Thm. 14.15]{MR3616245} and \cite[Thm.~5.1]{MR4265692}. Probably many of the conclusions of Theorem~\ref{theo:KRgeometry1} are very similar (or even included) in the material of \cite[Part C-III]{MR839450}. However, our assumptions slightly different since we do not make explicit reference to a positive semigroup but rather refer to the weak and strong maximum principles. 

\smallskip
(3) Our proof is quite direct and elementary and uses similar arguments as those used during the proof of \cite[Thm.~4.3]{MR3489637} and \cite[Thm.~5.1]{MR4265692}.
We learnt this kind of technique in the (less abstract and general) proof of the uniqueness part of \cite[Lem.~2.1]{MR2114128}. 

\smallskip
(4) From ii), 
we deduce that $\LL$ decomposes according to $X = X_0 \oplus X_1$ with  $X_1 := \hbox{\rm Span} \, f_1$ and $X_0 := (\hbox{\rm Span} \, \phi_1)^\perp = \{ f \in X; \, \langle f , \phi_1 \rangle = 0 \}$
in the sense of \cite[\S\,III.5.6]{MR0203473}. More precisely, $X = X_0 \oplus X_1$ is a topological direct sum, $\LL : X_0 \cap D(\LL) \to X_0$ and $\LL : X_1 \to X_1$. 
 
 \smallskip
(5) One can observe from the proof below that the conclusion (i) in Theorem~\ref{theo:KRgeometry1} holds under the sole assumptions {\Cyan\ref{X1}} for $X$ and $Y$, \ref{C1} and {\Cyan \ref{H4}{\bf-(i)}} for $\LL$ and $\LL^*$. 
The conclusion \eqref{eq:theoKRgeometry1-L}  holds when assuming furthermore that \eqref{eq:Irred:StrictOrderNorm} holds in {\Cyan $\EEE$}  and $\LL$ satisfies  \ref{H1'}. The  similar conclusion  \eqref{eq:theoKRgeometry1-L*} holds when assuming furthermore that \eqref{eq:Irred:StrictOrderNorm} holds in {\Cyan $\EEE^*$} and $\LL^*$ satisfies  \ref{H1'}. 

 \smallskip
(6)  
We finally recall that under condition \ref{H1}, the strong maximum principle {\Cyan \eqref{eq:Irred-StrongMP}} for $\LL$ is equivalent to the strong maximum principle {\Cyan \eqref{eq:Irred-StrongMP}}   for $\LL^*$ (see Lemma~\ref{lem:Irred-KStrongMax&dual}).
When furthermore  condition \ref{H2} holds and $\lambda_1$ in \ref{C1} is defined by \eqref{eq:exist1-deflambda1}, the weak maximum principle \eqref{eq:Irred-weakPM} for $\LL$ is equivalent to the weak maximum principle \eqref{eq:Irred-weakPM} for $\LL^*$ (see the proof of Lemma~\ref{lem:Existe1-H2bis}).

 \end{rem}

The proof of Theorem~\ref{theo:KRgeometry1} is split into the following Lemma~\ref{lem:f1phi1>0}, Lemma~\ref{lem:PositiveEigenvector},  Lemma~\ref{lem:Uniquenessf1} and  Lemma~\ref{lem:lambda1simple}. 

\begin{lem}\label{lem:f1phi1>0}  Under assumptions   \ref{X1}, \ref{C1}  and {\Blue \ref{H4}\bf-(i)}, we have  
\beqn\label{eq:Irred-f1>0phi1>0}
f_1 \gg 0 ,\quad \phi_1 \gg 0\quad\text{and}\quad\lambda_1=\lambda_1^*.
\eeqn
\end{lem}

\begin{proof}[Proof of Lemma~\ref{lem:f1phi1>0}.]
Applying the strong maximum principle to the two eigenfunctions $f_1 \in X \backslash\{0\}$ and $\phi_1 \in Y \backslash\{0\}$, we get the strict positivity of $f_1$ and $\phi_1$.
In particular $\langle f_1,\phi_1\rangle\neq0$.
Next, we write
$$
\lambda_1 \langle f_1,\phi_1 \rangle =   \langle \LL f_1,\phi_1 \rangle =   \langle f_1, \LL^* \phi_1 \rangle = 
\lambda_1^*  \langle f_1,\phi_1 \rangle, 
$$
which yields $\lambda_1=\lambda_1^*$.
\end{proof}
 

 \begin{lem}\label{lem:PositiveEigenvector} 
 Under assumptions   \ref{X1}, \ref{C1}  and {\Blue\ref{H4}{\bf-(i)}} for $\LL^*$ (resp.  $\LL$),  $\lambda_1$ is the unique eigenvalue associated to a positive eigenvector for $\LL$  (resp. for $\LL^*$).
\end{lem}

\begin{proof}[Proof of Lemma~\ref{lem:PositiveEigenvector}.]
 Consider $\lambda \in \C$  and $f \in X_+ \backslash \{0\}$ such that $\LL f = \lambda f$ and observe that from Lemma~\ref{lem:f1phi1>0}, we have $\phi_1 \gg 0$. 
We compute 
$$
0 = \langle (\lambda-\LL) f, \phi_1 \rangle =  \langle f, (\lambda-\LL^*) \phi_1 \rangle = (\lambda-\lambda_1) \langle f, \phi_1 \rangle,
$$
and thus $\lambda=\lambda_1$ since $ \langle f, \phi_1 \rangle > 0$.   The same proof applies to the dual problem.
\end{proof}

\begin{lem}\label{lem:Uniquenessf1} 
We assume again  \ref{X1}, \ref{C1}  and {\Blue\ref{H4}{\bf-(i)}}. 
Under the additional condition \ref{H1'} for $\LL$ and~\eqref{eq:Irred:StrictOrderNorm} in $\Blue\EEE$ (resp. \ref{H1'}  for $\LL^*$ and~\eqref{eq:Irred:StrictOrderNorm}   in $\Blue\EEE^*$), we have  $N(\LL-\lambda_1) = \hbox{\rm Span}(f_1)$ (resp. $N(\LL^*-\lambda_1) = \hbox{\rm Span}(\phi_1)$).  In particular,  $f_1$ (resp. $\phi_1$) is unique 
(because of the positivity and  the normalization conditions). 
\end{lem}

\begin{proof}[Proof of Lemma~\ref{lem:Uniquenessf1}.]
Consider a   eigenfunction $f \in X \backslash\{0\}$ associated to the eigenvalue  $\lambda_1$. 
First, we observe from Kato's inequality that
\bean
\lambda_1 |f| =   \lambda_1 \hbox{sign} (f)   f 
 = \hbox{sign} (f)   \LL f 
\le \LL |f|. 
\eean
That inequality is in fact an equality, otherwise we would have 
\bean
\lambda_1 \langle |f|,\phi_1 \rangle \not=   \langle \LL |f| ,\phi_1 \rangle  = \langle  |f| , \LL^* \phi_1 \rangle  = \lambda_1 \langle |f|,\phi_1 \rangle,
\eean
and a contradiction. As a consequence, $|f|$ is a solution to the eigenvalue problem $\lambda_1 |f|  = \LL |f|$,  so that
$$\lambda_1 f_\pm = \LL f_\pm,
$$
by writing $f_\pm = (|f|\pm f)/2$. 
The strong maximum principle assumption implies $f_\pm \gg 0$ or $f_\pm= 0$, and thus $\Blue 0 \ll f_+ \in \EEE$ or $\Blue 0 \ll f_- \in \EEE$ since $f\neq0$. 
Without loss of generality we may assume $f_+ \gg 0$.
From \eqref{eq:Irred:StrictOrderNorm} in $\EEE$, we then deduce $f \gg 0$. We introduce the normalized eigenfunctions $\tilde f := r f$ and $\tilde f_1 = r_1 f_1$ with 
\beqn\label{eq:Irred:proofuniquenessf1}
\langle \tilde f,\phi_1 \rangle = \langle \tilde f_1,\phi_1 \rangle = 1.
\eeqn
Now, thanks to Kato's inequality again, we write 
\bean
\lambda_1  (\tilde f-\tilde f_1)_+ = \hbox{sign}_+(\tilde f-\tilde f_1) \LL (\tilde f-\tilde f_1)  \le  \LL (\tilde f-\tilde f_1)_+ , 
\eean
and for the same reason  as above that last inequality is in fact an inequality. 
The strong maximum principle implies that either $(\tilde f-\tilde f_1)_+ = 0$, which also reads $\tilde f \le \tilde f_1$,  or $(\tilde f-\tilde f_1)_+ \gg 0$, which implies that $\tilde f \gg \tilde f_1$ by using again \eqref{eq:Irred:StrictOrderNorm} in $\Blue\EEE$. Because of the identity \eqref{eq:Irred:proofuniquenessf1} and the fact that $\phi_1 \in X'_+ \backslash \{0\}$, 
the second case in the above alternative is not possible. 
Repeating the same argument with $(\tilde f_1-\tilde f)_+$ we get that $\tilde f_1 \le \tilde f$ and we conclude with  $\tilde f = \tilde f_1$. 
 The same proof applies to the dual problem. 
\end{proof}

\begin{rem}\label{rem:Uniuqness-supersolution} Under the same hypotheses as in Lemma~\ref{lem:Uniquenessf1}, we have $\psi \in \hbox{\rm span} (\phi_1)$ if $\psi \in Y_+$ satisfies $\LL^*\psi \ge \lambda_1 \psi$ and  
$g \in \hbox{\rm span} (f_1)$ if $g \in X_+$ satisfies $\LL g  \ge \lambda_1 g$. In the second case, we indeed  cannot have $\LL^*g - \lambda_1 g \in X_+\backslash \{ 0 \}$, since this would implies 
$$
\langle \LL g - \lambda_1 g, \phi_1 \rangle > 0,
$$
and this would be in contradiction with the fact that 
$$
\langle \LL g - \lambda_1  g , \phi_1 \rangle = 
\langle g, \LL^* \phi_1 -  \lambda_1 \phi_1 \rangle = 0.
$$
We thus must have $\LL g  - \lambda_1g= 0$ and we conclude thanks to Lemma~\ref{lem:Uniquenessf1}.  The same proof applies to the dual problem. 
\end{rem}
\Black

 \begin{lem}\label{lem:lambda1simple}   Under the same assumptions as in Lemma~\ref{lem:Uniquenessf1}, 
$\lambda_1$ is algebraically simple  for $\LL$ (resp. for $\LL^*$).
\end{lem}

\begin{proof}[Proof of Lemma~\ref{lem:lambda1simple}.]
We use an induction argument. We have already proved that $N((\LL-\lambda_1)^k) = \hbox{\rm Span}(f_1)$ for $k=1$. 
Assume then the result proved for any $\ell$, $1 \le \ell \le k$, and consider $f \in N((\LL-\lambda_1)^{k+1})$. 
That means that $(\LL - \lambda_1) f \in N((\LL-\lambda_1)^k)$, and thus $(\LL - \lambda_1) f = r f_1$, with $r \in \R$, thanks to the induction hypothesis. 
If $r = 0$, then $f \in N(\LL-\lambda_1) = \hbox{\rm Span}(f_1)$. Otherwise, $r\not=0$,
and then  
$$
 \lambda_1 \langle f,  \phi_1 \rangle = \langle f, \LL^* \phi_1 \rangle = \langle \LL f, \phi_1 \rangle = \langle \lambda_1 f + r f_1, \phi_1 \rangle, 
 $$
 which in turn implies $r\langle f_1, \phi_1 \rangle = 0$ and a contradiction.
That concludes the proof. 
\end{proof}

\subsection{Mean ergodicity}
\label{subsec:Irred-MeanErgo}

We deduce from the above analysis a first  classical and general but rough information about the long-time behaviour of the trajectories associated to a semigroup.
 
More precisely, assuming  the existence and uniqueness of the first eigentriplet $(\lambda_1,f_1,\phi_1)$ for the generator $\LL$ of a semigroup $S$ and introducing the rescaled semigroup $\widetilde S_t := e^{-\lambda_1 t} \, S(t)$, we wish to establish the following mean ergodic property

\begin{enumerate}[label={\bf(E1)},itemindent=13mm,leftmargin=0mm,itemsep=1mm]
\item\label{E1}  for any $f \in X$, there holds
\beqn\label{eq:theo:MeanErgodicity}
\frac1T \int_0^T \widetilde S_t f  dt \to \langle f, \phi_1 \rangle f_1,  \ \hbox{ as } \ T \to \infty, 
\eeqn
in a sense to be specified.
\end{enumerate}

\smallskip

We start with a general result, taken from~\cite[Thm.~V.4.5]{MR1721989}, which states that, under the conclusions of Theorem~\ref{theo:KRgeometry1}, \ref{E1} holds for the strong topology if the semigroup $(\widetilde S_t)$ is bounded.

\begin{theo}\label{theo:MeanErgodicity-generalEN}
Consider a positive semigroup $S$ on a Banach lattice $X$ and 
assume that its generator $\LL$ satisfies the conclusions of Theorem~\ref{theo:KRgeometry1} about the existence and uniqueness 
of the first eigentriplet $(\lambda_1,f_1,\phi_1)$.
Assume furthermore that $(\widetilde S_t)_{t\geq0}$ is bounded.
Then, the above mean ergodic property \ref{E1} holds  for the strong  topology.
\end{theo}

\begin{proof}[Proof of Theorem~\ref{theo:MeanErgodicity-generalEN}]
Following the proof of~\cite[Thm.~V.4.5]{MR1721989}, we consider the subspace
\[X_0 := \Span f_1 \oplus \Span\{f-\widetilde S(t)f\,:\, f\in X, t\geq0\}\]
of $X$ and we take $\phi\in Y$ which vanishes on $X_0$.
Since $\phi$ vanishes on each element of the form $f- \widetilde S(t)f$, this implies that $\widetilde S^*(t)\phi=\phi$ for all $t\geq0$.
We deduce that $\LL^*\phi=\lambda_1 \phi$, and consequently $\phi\in\Span\phi_1$ due to the point {\it ii)} in Theorem~\ref{theo:KRgeometry1}.
Since we also have $\langle \phi,f_1\rangle=0$, we deduce that $\phi=0$ and therefore $\overline{X_0}=X$.
 We observe now 
\[\bigg(\int_0^T\widetilde S(s)ds\bigg)\big(I-\widetilde S(t)\big)=\big(I-\widetilde S(T)\big)\int_0^t \widetilde S(s)ds\]
for all $t,T\geq0$, which is an immediate consequence of the semigroup property. 
The above relation and the boundedness assumption on $(\widetilde S_T)_{T\geq0}$ imply that the convergence~\eqref{eq:theo:MeanErgodicity} 
 holds for 
$f = g - \widetilde S(t) g$ with $g \in X$, $t \ge 0$, and thus for any $f\in X_0$.
Finally, since $X_0$ is dense in $X$ and using again the fact that $(\widetilde S_t)_{t\geq0}$ is bounded, we can readily extend the validity of~\eqref{eq:theo:MeanErgodicity} to any $f\in X$.
\end{proof}


We will now give weaker versions of Theorem~\ref{theo:MeanErgodicity-generalEN} with proofs which are based on compactness arguments.
The motivation for providing such alternative proofs that require stronger assumptions is that, unlike the proof of Theorem~\ref{theo:MeanErgodicity-generalEN}, the methods can be adapted to derive stronger ergodicity results, namely without averaging in time, see Section~\ref{subsec:ergodicity}.

\begin{theo}\label{theo:MeanErgodicity}
Consider a positive semigroup $S$ on a Banach lattice $X$ and 
assume that its generator $\LL$ satisfies the conclusions of Theorem~\ref{theo:KRgeometry1} about the existence and uniqueness 
of the first eigentriplet $(\lambda_1,f_1,\phi_1)$. 
With the above notations, we  assume furthermore that 

{\bf (1)} $(\widetilde S_t)$ is bounded; 

{\bf (2)} $B_X$ is weakly compact for a topology which makes $f \mapsto \langle f, \phi_1 \rangle$ continuous.

Then, the above mean ergodic property \ref{E1} holds  for the  topology introduced in {\bf (2)}.
\end{theo}
 


\begin{proof}[Proof of Theorem~\ref{theo:MeanErgodicity}] Fix $f \in X$ and define 
$$
u_T := \frac1T \int_0^T \widetilde S_t f \, dt.
$$
 From {\bf (1)} {\Cyan and denoting by $M$ the uniform bound of $(\| \widetilde S_t \|)$}, we have 
 $$
 \| u_T \| \le  \frac1T \int_0^T \| \widetilde S_t f \| \, dt \le M \, \| f \|, \quad \forall \, T > 0.
 $$
We also compute 
$$
\langle u_T,\phi_1 \rangle =  \frac1T \int_0^T  \langle \widetilde S_t f, \phi_1 \rangle \, dt =  \langle   f, \phi_1 \rangle, 
\quad \forall \, T > 0.
 $$
Thanks to assumption {\bf (2)}, we deduce that there exists $f^* \in X$ and a sequence $(T_k)$ such that 
$$
u_{T_k} \to f^* \quad\hbox{and}\quad \langle f^*,\phi_1\rangle =  \langle   f, \phi_1 \rangle.
$$
Because $( \widetilde S_t f)$ is bounded, we may use the usual ergodicity trick as in the second  proof of Theorem~\ref{theo:KRexistBIS}
 and for any $t > 0$, we have 
$$
\widetilde S_t f^* - f^* = \lim_{k \to \infty} \frac1{T_k}\Bigl\{ \int_{T_k}^{T_k+t} \widetilde S_s f ds - \int_{0}^{t} \widetilde S_s f ds \Bigr\} = 0.
$$
We have established $( \LL - \lambda_1) f^* = 0$, so that $f^* \in \hbox{\rm Span}(f_1)$ and more precisely $f^* = \langle f, \phi_1 \rangle  f_1$. 
By uniqueness of the limit, it is the whole family $(u_T)$ which converges to $f^*$.
\end{proof}
  
\medskip
We present a variant of the previous result in which we see  that in a very general framework (including all the applications we present in the second part of this work) the above hypotheses {\bf (1)} and {\bf (2)} are not needed (or more precisely are automatically satisfied).

\begin{theo}\label{theo:MeanErgodicityVariante1}
(1) Consider a Banach lattice $X \subset \Lloc^1(E,\EE,\mu)$ and  $Y\subset \Lloc^1(E,\EE,\mu)$ (so that in particular $\phi_1\in \Lloc^1$ and $L^1_{\phi_1}$ is well-defined) and a positive semigroup $S$ on $X$ such that its generator $\LL$ satisfies 
the conclusions of Theorem~\ref{theo:KRgeometry1} about the existence, positivity and uniqueness 
of the first eigentriplet $(\lambda_1,f_1,\phi_1)$. 
Then the mean ergodic convergence \ref{E1} holds for the weak topology of $L^1_{\phi_1}$.

\smallskip
(2) Assuming additionally  that $S$ is strongly continuous and that 
\beqn\label{eq1:rem:MeanErgodicityVariante1}
\XX^k := (D(\LL^k), \| \cdot \|_{\XX^k}) \subset \Lloc^1 \hbox{ with strong compact embedding for some } k \ge 1, 
 \eeqn
 where
$$
\| f \|_{\XX^k} := \| f \|_{L^1_{\phi_1}} + \dots +  \| \LL^kf \|_{L^1_{\phi_1}} , \quad \forall \, f \in D(\LL^k), 
$$
then the mean ergodic convergence \ref{E1} 
holds for the strong topology of $L^1_{\phi_1}$.
 
 \end{theo}

\begin{proof}[Proof of Theorem~\ref{theo:MeanErgodicityVariante1}.]
{\sl Step 1.} We first recall a very classical result about conservative semigroups. Denoting $\widetilde S_t := e^{-\lambda_1 t} \, S(t)$, we observe that this rescaled semigroup satisfies 

\begin{itemize}

\item[\bf (i)] $\widetilde S_t \ge 0$; 

\item[\bf (ii)] $\widetilde S_t  f_1 = f_1$ for any $t\ge 0$; 

\item[\bf (iii)] $\langle \widetilde S_t  g, \phi_1 \rangle  = \langle  g, \phi_1 \rangle $ for any $g \in X$ and  $t\ge 0$.

\end{itemize}

We denote $[f]_1 := \langle |f|,\phi_1\rangle$ which is a norm on $X$ (we use here that $\phi_1 \gg 0$) and $\widetilde S_t$ is obviously a contraction for this one. Indeed,  for any $f \in X$,  there holds
$$
|\widetilde S_t f| =  |\widetilde S_t f_+  - \widetilde S_t f_-| \le \widetilde S_t f_+  + \widetilde S_t f_-= \widetilde S_t |f|,
$$
using  {\bf (i)} in the inequality, and next 
\beqn\label{eq:theo:MeanErgodicityVariante1-1}
[\widetilde S_t f]_1 = \langle |\widetilde S_t f|, \phi_1 \rangle \le \langle \widetilde S_t |f|, \phi_1 \rangle = [f]_1, 
\eeqn
using {\bf (iii)} in the last equality.  Abusing notations, we also denote by $\XX$ the completion of $X$ for the $L^1_{\phi_1}$ norm (so that we may identify $\XX$ to a closed subspace of $L^1_{\phi_1}$). 
We may then extend $\widetilde S_t$ to $\XX$ by uniform continuity and this extension still satisfies the properties   {\bf (i)}-{\bf (ii)}-{\bf (iii)} on $\XX$.
Consider now $f \in X$ such that $H(f/f_1) f_1 \in X$ for some convex function $H : \R \to \R$, where we use here that $X \subset \Lloc^1$, and thus in particular  $f_1 > 0$ a.e. on $E$, 
 in order to give a sense to the term $H(f/f_1) f_1$. 
From {\bf (ii)}, we have 
$$
\ell [(\widetilde S_t  f)/f_1] f_1 = \widetilde S_t  [ \ell(f/f_1) f_1], 
$$ 
for any real affine function $\ell$. Next from {\bf (i)}  and \eqref{eq;AfveeAgLEAfg}, we have 
$$
H [(\widetilde S_t  f)/f_1] f_1 \le \widetilde S_t  [H(f/f_1) f_1], 
$$ 
because of $H = \sup_{\ell \le H} \ell$ and the supremum can be taken on a numerable set of affine functions. 
Thanks to  {\bf (iii)}, we conclude that 
 \beqn\label{eq:SGMarkovEntropy}
\langle H [(\widetilde S_t f)/f_1] f_1 , \phi_1 \rangle \le \langle H [f/f_1] f_1 , \phi_1 \rangle , \quad \forall \, t \ge 0.
\eeqn

\smallskip
\noindent
{\sl Step 2.} We  normalize $\langle f_1,\phi_1 \rangle = 1$. 
For $f \in \XX \subset L^1_{\phi_1}$ so that $f\phi_1 = (f/f_1) f_1 \phi_1 \in L^1$, the de la Vallée Poussin theorem tells us that there exists an even and  convex function $H : \R \to \R_+$ such that $H(s)/s \to +\infty$ as $s \to \infty$
and $H(f/f_1) f_1 \phi_1\in L^1$. Using the notations of the proof of Theorem~\ref{theo:MeanErgodicity}, the Jensen inequality and the above estimate \eqref{eq:SGMarkovEntropy}, we deduce 
$$
\int_E H(u_T/f_1) f_1 \phi_1 d\mu \le \frac1T \int_0^T \int_E H[(\widetilde S_t f)/f_1] f_1 \phi_1 d\mu dt \le 
  \int_E H(f/f_1) f_1 \phi_1 d\mu, 
$$
for any $T > 0$. Now, for any $A \in \EE$ and $T, K > 0$, we have 
\bean
\int_E  u_T {\bf 1}_A \phi_1 d\mu  
&=& \int_E  \frac{u_T}{ f_1}  {\bf 1}_{ \frac{|u_T| }{ f_1} > K}  {\bf 1}_A f_1 \phi_1 d\mu   + 
 \int_E  \frac{u_T}{ f_1}  {\bf 1}_{ \frac{|u_T| }{ f_1} \le K}  {\bf 1}_A f_1 \phi_1 d\mu 
 \\
&\le& \frac{K }{ H(K)}  \int_E  H(u_T/ f_1) f_1 \phi_1 d\mu   + 
K  \int_E    {\bf 1}_A f_1 \phi_1 d\mu
 \\
&\le& \frac{K }{ H(K)}  \int_E H(f/f_1) f_1 \phi_1 d\mu   + 
K  \int_E    {\bf 1}_A f_1 \phi_1 d\mu,
\eean
from what we immediately deduce, by using the Dunford-Pettis theorem, that $(u_T)$ belongs to a weak compact set of $L^1_{\phi_1}$. We conclude that \eqref{eq:theo:MeanErgodicity} holds for the weak convergence in $L^1_{\phi_1}$ as in the proof of Theorem~\ref{theo:MeanErgodicity}.

\smallskip
\noindent
{\sl Step 3.} We now additionally assume that \eqref{eq1:rem:MeanErgodicityVariante1} holds
with strong compact embedding for some $k \ge 1$. 
Taking  $f \in D(\LL^k)$, we compute 
$$
\langle |\LL^j(\widetilde S_t f)|,\phi_1 \rangle = \langle |\widetilde S_t ( \LL^j f)|,\phi_1 \rangle  
\le \langle | \LL^j f|,\phi_1 \rangle , 
$$
for any $j \le k$ and any $t \ge 0$, and thus the same bound holds for $(u_T)$. From \eqref{eq1:rem:MeanErgodicityVariante1}, we deduce that up to the extraction of a subsequence, $(u_{T})$
converges a.e. on $E$. Together with the weak convergence in  $L^1_{\phi_1}$ yet established, we classically deduce  that the whole family $(u_T)$ converges for the strong topology  in  $L^1_{\phi_1}$.
We conclude that the same holds for any $f \in X$ by taking advantage
of the fact that $D(\LL^k)$ is dense in $X$ for the strong topology of $X$, and thus for the strong topology of $\XX$,  and of the estimate of contraction \eqref{eq:theo:MeanErgodicityVariante1-1}.
\end{proof}

\begin{rem}\label{rem:MeanErgodicityVariante1}

(1) A similar conclusion holds as in Theorem~\ref{theo:MeanErgodicityVariante1} when we assume $X \subset \Mloc^1$,  $D(\LL^k) \subset \Lloc^1$ and $D(\LL^{*k} )\subset \Lloc^1$ for some $k \ge 1$
instead of $X, Y \subset \Lloc^1$. 
For $f \in D(\LL^k) \subset \Lloc^1$, we may indeed repeat the proof of  Theorem~\ref{theo:MeanErgodicityVariante1} and we obtain the same conclusion. 
We next define $\XX$ as the closure of $D(\LL^k)$  for the norm $[\cdot]_{\phi_1}$. We conclude that \eqref{eq:theo:MeanErgodicity} holds weakly in $L^1_{\phi_1}$ for any $f \in \XX$ by a density argument. 

\smallskip 
(2)  The proof of Theorem~\ref{theo:MeanErgodicityVariante1} is based on so-called General Relative Entropy (GRE) techniques as developed for instance in \cite{Loskot1991}, 
\cite{MR2162224} and \cite{Bernard2022}. These ones are known to be useful for some classes of PDE and for stochastic semigroups in order to establish uniform in time estimates and longtime convergence results.
 \end{rem}

The main interest of the two previous results {\Blue stated in Theorem~\ref{theo:MeanErgodicityVariante1}} is that they do not ask any new information on the {\Blue rescaled semigroup $\widetilde S$} but they are just based on the eigentriplet stationary problem. The shortcoming is that they are formulated in terms of the norm
$[\cdot]_{\phi_1}$ instead of the norm of $X$.
{\Blue We point out that under the same conditions as in Theorem~\ref{theo:MeanErgodicityVariante1}-(1), we may directly apply Theorem~\ref{theo:MeanErgodicity} and obtain a slightly stronger conclusion.}
We present a second variant of  Theorem~\ref{theo:MeanErgodicity} which is well adapted to the  splitting framework developed in Sections~\ref{sec:ExistenceKR} and \ref{sec:DynamicalExistenceKR} 
and is precisely formulated in a weak or strong topology of a space $X_0 \supset X$.

\begin{theo}\label{theo:Irred-MeanErgoBIS} 

Consider a positive semigroup $S=S_\LL$  such that  $\LL$ satisfies the conclusions of Theorem~\ref{theo:KRgeometry1} about the existence and uniqueness 
of the first eigentriplet $(\lambda_1,f_1,\phi_1)$. Assume furthermore that $S$ satisfies the  splitting structure introduced in
\ref{HS2} in Section~\ref{subsect-Exist2-Dissip} or  \ref{HS3}  in  Section~\ref{subsect:AboutWeakDissip}, or more precisely, there exist two families of operators $(V(t))$ and $(W(t))$ such that 
$$
S = V + W * S, 
$$
a real number $\kappa \le \lambda_1$ and     some decaying functions $\Theta_i : \R_+ \to \R_+$ with $\Theta_1(t) \to 0$ as $t \to \infty$, $\Theta_2 \in L^1(\R_+)$ such that  the following estimates hold 
\bean
&&\| V(t) e^{-\kappa t}   \|_{\BBB(X)} = \OO(1 ),  \quad \| V(t) e^{-\kappa t}   \|_{\BBB(X,X_0)} = \OO(\Theta_1 ),  
\\ 
&&\| W(t) e^{-\kappa t}    \|_{\BBB(\XX_0,\XX_1)} = \OO(\Theta_2 ), 
\eean
with $\XX_1 \subset   X_0 \subset \XX_0$, where $\XX_0$ is the space $X$ endowed with the norm $[g]_{\phi_1} := \langle |g|, \phi_1 \rangle$. 

\smallskip
(1) Assume furthermore that $\XX_1 \subset X_0$ with compact embedding for the weak or the strong topology in $X_0$ and this topology makes $f \mapsto \langle f, \phi_1 \rangle$ continuous.  
Then the mean ergodic convergence \ref{E1} holds true for the above strong or weak topology. 

\smallskip

(2) Assume furthermore that $X \subset \Lloc^1$, $S$ is strongly continuous, and that the space $\XX^k$ defined by \eqref{eq1:rem:MeanErgodicityVariante1} is strongly compact embedded in $\Lloc^1$  for some $k \ge 1$.
Then the mean ergodic convergence \ref{E1} holds true for the  strong topology of $X_0$. 
\end{theo}
 
\begin{proof}[Proof of Theorem~\ref{theo:Irred-MeanErgoBIS}.]
We define 
$$
\widetilde V (t) := V(t)  e^{-\lambda_1 t}, \quad 
\widetilde W (t) := W(t)  e^{-\lambda_1 t},
$$
so that 
$$
\widetilde S = \widetilde V +\widetilde W * \widetilde S,
$$
and
\bear\label{eq:exist2-HS3BIST}
&&M := \sup_{t \ge 0} \| \widetilde V(t) \|_{\BBB(X)} < \infty, \quad 
 \| \widetilde V \|_{\BBB(X,X_0)} \lesssim  \widetilde\Theta_1 \in C_0(\R_+),
\\ \nonumber 
&&
 \widetilde\Theta_2(t) := \| \widetilde W (t)\|_{\BBB(\XX_0,\XX_1)} \in L^1(\R_+). 
\eear

\smallskip
{\sl Step 1.} We furthermore assume (1) and that the weak topology of $X_0$ makes $f \mapsto \langle f, \phi_1 \rangle$ continuous. 
We denote by $\TTT$ the weak or the strong topology $X_0$ (depending of the assumption made on the embedding $\XX_1 \subset X_0$).
For $f_0 \in X$, we split
$$
f(t) := \widetilde S_t f_0 = v(t) + k(t), \quad v(t) :=  \widetilde V (t) f_0, \quad k(t) :=  (\widetilde W * \widetilde S)(t) f_0,
$$
and we observe that $\| v(t) \|_{X_0} \to 0$ as $t \to \infty$ from the second estimate in \eqref{eq:exist2-HS3BIST}. On the other hand, we have 
\bean
\sup_{t \ge 0} \| k(t) \|_{\XX_1} \le \| \widetilde \Theta_2 \|_{L^1} \sup_{t \ge 0} \| \widetilde S_t  f_0 \|_{\XX_0} \le \| \widetilde  \Theta_2 \|_{L^1} \| f_0\|_{\XX_0},
\eean
from \eqref{eq:theo:MeanErgodicityVariante1-1}.  In particular, $k(t)$ belongs to a compact set of $\TTT$, so that 
$(f(t))_{t \ge 0}$ also belongs to a compact set for the same topology $\TTT$. The same argument used on the Cesàro function $(u_T)$ defined 
during the proof of Theorem~\ref{theo:MeanErgodicity} implies that  there exist $f^* \in X$ and a sequence $(T_k)$ such that 
$$
u_{T_k} \to f^* \ \hbox{in the sense of }\TTT\quad\hbox{and}\quad \langle f^*,\phi_1\rangle =  \langle   f, \phi_1 \rangle,
$$
the last identity following from the assumption that $f \mapsto \langle f, \phi_1 \rangle$ continuous for $\TTT$.
We may then conclude as in the proof of Theorem~\ref{theo:MeanErgodicity}.

\smallskip
{\sl Step 2.} 
We furthermore assume (2), and by linearity we may assume $f_0 \in X$, $\langle f_0, \phi_1 \rangle = 0$. 
We recall that \eqref{eq:theo:MeanErgodicity} holds for the strong topology of $L^1_{\phi_1}$ from Theorem~\ref{theo:MeanErgodicityVariante1}
and that  $\| v(t) \|_{X_0} \to 0$ as $t \to \infty$ from Step 1.
Arguing as in Step~3 of the proof of Theorem~\ref{theo:KRexistTER}, we have 
\bean
K(T) := \frac1T\int_0^T (\widetilde W *\widetilde S)(t)  \, dt 
&=&  \frac1T \int_0^T \widetilde W(s) \int_0^{T-s} \widetilde S(u) \, du ds  
\\
&=&  \int_0^T \widetilde W(s) \frac{T-s }{ T} \widetilde U(T-s)    \, ds, 
\eean
where $\widetilde U_T :=U_T e^{-\lambda_1 T}$, $U_t$ is defined by \eqref{def:UTVTWT}, so that  $ u_T = \widetilde U_T f_0$ and   $ [ u_T ]_{\phi_1} \to 0$ as $T\to\infty$ from  Theorem~\ref{theo:MeanErgodicityVariante1}.
As a consequence, we have 
\bean
\| K(T) f_0 \|_{\XX_1}
&\le& \int_0^{T/2} \widetilde \Theta_2(s)  [ \widetilde U(T-s) f_0 ]_{\phi_1}  \, ds +  \int_{T/2}^T \widetilde\Theta_2(s)  [ \widetilde U(T-s) f_0]_{\phi_1}  \, ds
\\
&\le& \| \widetilde \Theta_2 \|_{L^1} \sup_{t \ge T/2} [ \widetilde U(t) f_0]_{\phi_1}  +  \int_{T/2}^\infty \widetilde \Theta_2(s) ds  \sup_{t \ge0 } [ \widetilde U(t) f_0]_{\phi_1}  \to 0, 
\eean
as $T\to \infty$. Altogether, we have established that $\| u_T \|_{X_0} \to 0$ as $T \to \infty$.
\end{proof}

{\Blue
Let us finally summarize the results obtained in this section concerning our spectral analysis program. 
For this purpose, we emphasize that the spectrum set of $\LL$ may be sensitive to the norm considered for defining it. 
We still denote by $\Sigma(\LL)$ the spectrum set defined in $X$ (and its norm) in Section~\ref{subsec:BanachLattice}. 
We may similarly define the spectrum of $\LL$  in the space $L^1({\phi_1}) \supset X$ and we denote it  by $\Sigma'(\LL)$ in the following discussion.  

\smallskip
On the one hand, Theorem~\ref{theo:MeanErgodicity-generalEN},  Theorem~\ref{theo:MeanErgodicity} or Theorem~\ref{theo:Irred-MeanErgoBIS} (when $\XX_1 = X$ and $\XX_0 = L^1({\phi_1})$) implies that $\Sigma(\LL) \subset \{ z \in \C;\ \Re e z \le \lambda_1 \}$ in $X$, because $(\| \widetilde S(t) \|_{\BBB(X)})$ is bounded in these cases and we may use the resolvent identity \eqref{eq:Exist1-DefRepresentationRR}. 
For a similar reason, Theorems~\ref{theo:MeanErgodicityVariante1}   implies that $\Sigma'(\LL)\subset\{z\in\C;\ \Re ez\leq\lambda_1\}$ in $L^1({\phi_1})$,  because $(\| \widetilde S(t) \|_{\BBB(L^1({\phi_1}))})$ is bounded in this case. 
As a consequence, all these results imply that \ref{S1} holds in either $X$ or $L^1({\phi_1})$. 

\smallskip
On the other hand, the conclusion 
\ref{C1} together with Theorem~\ref{theo:KRgeometry1} provide the first answer \ref{S2} about the geometry of the spectrum problem. The mean ergodicity convergence \ref{E1} is next deduced in the many variants Theorem~\ref{theo:MeanErgodicity-generalEN},   Theorem~\ref{theo:MeanErgodicity}, Theorem~\ref{theo:MeanErgodicityVariante1}, Theorem~\ref{theo:Irred-MeanErgoBIS}. 
Altogether, these results thus provide a proof for Theorem~\ref{theo:main-intro}-(2). In the weak dissipativity framework, we may also observe that \ref{C1} and the above mentioned theorems imply that \ref{S1}, \ref{S2} and \ref{E1} hold at least in $L^1(\phi_1)$, and for the  weak  convergence. 


 }

\bigskip
%
%

\bigskip
\section{The geometry of the boundary  point spectrum} 
\label{sec:geo2}


{\Blue We establish the sharper geometric parts \ref{S31} and  \ref{S32} as well as the evolution parts \ref{E2} and \ref{E31} of the Krein-Rutman theorem. For that purpose, we introduce in particular the partial conclusion \ref{C2} of \ref{S1}-\ref{S2} established in the previous section, the structure conditions \ref{X2} and \ref{X3} on the Banach lattice $X$  and the aperiodicity like conditions  \ref{H5} and  \ref{H5'}.}

 \medskip 
We summarize the results established up to now by assuming that the main conclusions  in the previous sections are achieved, namely 

\begin{enumerate}[label={\bf(C2)},itemindent=13mm,leftmargin=0mm,itemsep=1mm]
\item\label{C2} the  first eigentriplet  problem {\Blue\eqref{eq:triplet1}-\eqref{eq:triplet2}} 
has a unique solution $(\lambda_1, f_1, \phi_1)$, and  furthermore,  $f_1 \gg 0$ and $\phi_1 \gg 0$. 
In that situation, we make the usual normalization \eqref{eq:normalize-f1phi1}.
\end{enumerate}
{\Blue Compared to~\ref{S1}-\ref{S2}, the (geometric and algebraic) simplicity of $f_1$ and $\phi_1$ is not required in~\ref{C2} nor the localization of the spectrum 
$\Sigma(\LL)\subset\{z\in\C,\ \Re e(z)\leq\lambda_1\}$.
}

\medskip

In this section, we aim to describe one step further the geometry of the spectrum and more precisely to get   some information {\Blue about the principal point spectrum
\[
\Sigma^1_P(\LL) :=  \Sigma_P(\LL) \cap \overline{ \Delta}_{\lambda_1}.
\]
Here and for further references below, we recall that we define the sets 
$$
\Sigma_d(\LL) \subset  \Sigma_P(\LL) \subset  \Sigma(\LL),
$$
where the point spectrum set $\Sigma_P(\LL)$ is the set of eigenvalues, namely $\lambda \in \Sigma_P(\LL)$ if $N(\LL - \lambda) \not=\{0\}$, 
and the discrete spectrum set $\Sigma_d(\LL)$ is the set of eigenvalues which are isolated and have finite algebraic multiplicity. 
It is worth pointing out that in general the  principal point spectrum $\Sigma^1_P(\LL)$ differs from the boundary point spectrum $\Sigma^+_P(\LL)$, which is defined in Section~\ref{ssec:main-result} as
$$
  \Sigma^+_P(\LL) :=  \Sigma_P (\LL) \cap \Sigma_+ (\LL).
$$
Because $\lambda_1 \in \Sigma(\LL)$ and thus $\lambda_1 \le s(\LL)$, we always have $ \Sigma^+_P(\LL) \subset \Sigma^1_P(\LL)$. 
However, these two sets coincide in the case when $\lambda_1=s(\LL)$, which is the conclusion~\eqref{eq:triplet3} in~\ref{S1}.   For further reference, we introduce the translated operator $\widetilde \LL := \LL - \lambda_1$ and the associated principal point spectrum
$$
\Sigma^1_P (\widetilde \LL) := \Sigma_P (\widetilde \LL) \cap \overline{ \Delta}_{0} = \Sigma^1_P ( \LL) - \lambda_1.
$$
Analyzing these dominant parts of the spectrum will be made} possible by first introducing a suitable and usual complexification framework and next by assuming a stronger positivity property on $\LL$ or on
 the associated semigroup.

%

 \medskip
\subsection{Complexification and the sign operator}
\label{subsec:SP-defStrongP} \

 We present some materials, most of them being very  classical, about the sign operator in a complex Banach lattice and we refer to \cite{MR0423039,MR839450} for more details. 

  \medskip
{\bf Complexification.} 
 The complexification space $X_\C$ associated to a real Banach lattice $X$ is defined by $X_\C := X + i X$ so that $f \in X_\C$ if $f = g + ih$, $g,h \in X$. 
In general, we just write $X$ without mentioning the field, although when we need to specify it, we write $X_\C$ or $X_\R$.
We extend on $X_\C$ the order defined on $X_\R$ by writing 
$$
f = g + i h  \ge 0 \quad\hbox{if}\quad g\ge 0 \hbox{ and } h= 0.
$$
The conjugate $\bar f$ of a complex vector $f=g+ih$ is   defined by $\bar f=g-ih$.
We then define the modulus
\beqn\label{def:modulus}
|f| :=  \sup_{\theta \in [0,2\pi]}  (g \cos \theta + h \sin\theta ),    
\eeqn
which indeed exists for such a family of vectors. 
One checks the usual modulus properties:  
$$
 |f| \ge 0, \quad  |f| = 0 \ \hbox{ iff } \ f=0, \quad |\lambda f| = |\lambda| \, |f|, \quad |f + g| \le |f|  + |g|,
$$
for any $f, g \in X$ and $\lambda \in \C$. We finally define the norm on $X_\C$ by 
$$
\| f \| := \| |g+ih| \|_{X_\R}, 
$$
and we observe that $X_\C$ has  a complex Banach lattice structure.  
We extend the definition of $A\in\BB(X_\R)$ to $X_\C$ by setting
$$
A (g+ih) = Ag + i Ah, \quad \forall \, g+ih \in X_\C.
$$

\medskip
{\bf  The operator ${\rm sign}$.} We  extend the $\sign$ operator  defined in Section~\ref{subsec:MorePositive} to the present complex Banach lattice framework. 
Instead of dealing with the most general case, we will use some regularity assumption on the Banach lattice $X$ which is suitable for our purpose and that we present below. 
Similarly as in Remark~\ref{rem:Ideal}, for $f \in X$, we define 
$$ 
X_f :=  \bigcup_n \{ g \in X; \ |g| \le n |f| \}, 
$$
and next, similarly as in Theorem~\ref{theo:PLLbis}, we define
$$
A_f[g] := \inf \{ C > 0; \ |g| \le C |f| \}, \quad \forall \, g \in X_f.
$$

  \smallskip
We summarize the regularity conditions we need  on the Banach lattice $X$ by assuming : 

\begin{enumerate}[label={\bf(X2)},itemindent=13mm,leftmargin=0mm,itemsep=1mm]
\item\label{X2} For any $f \in X$ such that $|f| \in X_{++}$, there exists a sign operator $\sign f \in \BBB(X)$, with the following properties
\bear\label{eq:irred:signCprop1}
&&
\sign f \circ \sign \bar f = I, \quad (\sign f ) f  =   |f|, 
\\
&& \label{eq:irred:signCprop2}
(\sign f  ) \, g =  (\sign \,( uf ) ) \, (ug),
\quad|(\sign f  ) \, g|    \le |g|,  \quad \forall \, g \in X, \ \forall \, u \in \Sp^1,
\eear
\end{enumerate}

and furthermore 

\begin{enumerate}[label={\bf(X3)},itemindent=13mm,leftmargin=0mm,itemsep=1mm]
\item\label{X3} for any $f \in X$ such that $|f| \in X_{++}$,   the inclusion $X_f\subset X$ is dense for the strong, the weak, or the weak-$*$ topology, and for all $f\in X$ and $g\in X_f$
\beqn\label{eq:lem:Geo-QbQ>0}
\big(\, g\in X_\R \ \text{and}\ |g| \le C |f| \,\big) \ \Leftrightarrow \ A_f[g - i r |f|] \le\sqrt{C^2+r^2},  \, \forall \, r \in \R.
\eeqn
\end{enumerate}


\medskip

For a space of functions, the $\sign$ operator is defined as the multiplication by (abusing notations)
\beqn\label{eq:defSignCf}
\sign f := \bar f /|f|, \quad \forall \, f \in X, \,\, |f| \in X_{++}.
\eeqn

\begin{lem}\label{lem:Ef=X}
With \eqref{eq:defSignCf}, the properties \ref{X2} and \ref{X3} hold when $X =  L^p(E,\EE,\mu)$ or $X = C_0(E)$. 
\end{lem}

\begin{proof}[Proof of Lemma~\ref{lem:Ef=X}.]
For $f \in X$, $|f| \in X_{++}$, we just  indicate the proof of $\overline X_f = X$, the other algebraic properties being clear from the definition \eqref{eq:defSignCf}. When $f \in L^p$ such that $|f|> 0$ $\mu$-a.e. and $0 \le g \in L^p$, we set $g_n := g  \wedge (n |f|)$. We have $0 \le g_n \le g$ and $g_n \to g$ strongly $L^p$ if $p < \infty$ and weakly-$*$ $L^\infty$ if $p=\infty$. The general case $g \in L^p$ is dealt in the usual way by introducing positive and negative parts  and next real and imaginary part. That concludes the proof of $\overline X_f = L^p$. The proof of $\overline X_f = C_0(E)$ is similar.
\end{proof}

A sign operator satisfying \ref{X2} and \ref{X3} can actually be built by using Kakutani's theorem in general Banach lattices whenever $|f|$ is a quasi-interior point, see for instance~\cite[Chapter 14.3]{MR3616245}.
In $X=L^\infty(E,\EE,\mu)$, being a quasi-interior point is more demanding than belonging to $X_{++}$, and our framework is thus more general in that case.
In $X=M^1(E)$, the situation is even worst since there is no quasi-interior point, so the approach via Kakutani's theorem does not provide any sign operator.
However, we can associate to $f\in M^1(E)$ such that $|f|\gg0$ a sign operator by means of the Radon-Nikodym theorem.
Denoting $\alpha:E\to\Ss^1$ the measurable function such that $f=\alpha|f|$, the multiplication by $\bar \alpha/|\alpha|$ defines a sign operator $\sign f\in\BB(X)$, or in other words (abusing notations)
\beqn\label{eq:defSignCf-M1}
\sign f := \bar \alpha/|\alpha|, \quad \forall \, f=\alpha|f| \in M^1, \,\, |f| \in M^1_{++}.
\eeqn

\begin{lem}\label{lem:Ef=M1}
With the definition~\eqref{eq:defSignCf-M1}, $X =  M^1(E)$ enjoys the properties \ref{X2} and \ref{X3}.
\end{lem}

\begin{proof}[Proof of Lemma~\ref{lem:Ef=M1}.] 
As for Lemma~\ref{lem:Ef=X}, we only sketch the proof of the density property $\overline X_f=X$, which holds here for the weak-$*$ topology, the other algebraic properties being clear from the definition \eqref{eq:defSignCf-M1}.
Without loss of generality, we may take $f \in X_{++}$,
meaning that $f(\OO)>0$ for any open set $\OO\subset E$.
For $\eps>0$ and $\varphi\in C_0(E)$, we can find a partition $E_1$, \dots, $E_n$ of $E$ and some elements $x_1$, \dots, $x_n$ of $E$ such that 
for any $  i\in\{1,\cdots,n\}$: 
\[
f(E_i)>0, \quad x_i \in E_i \quad \text{and} \quad  \sup_{x\in E_i}|\varphi(x)-\varphi(x_i)|<\eps.
 \]
For $g\in X$ and $\eps > 0$, defining $g_\eps$ by
\[g_\eps := \sum_{i=1}^n \frac{g(E_i)}{f(E_i)}f_{|E_i} \in X_f, \]
we have  
\[\bigg|\langle g_\eps,\varphi\rangle - \sum_{i=1}^n \varphi(x_i)g(E_i) \bigg| \leq \sum_{i=1}^n |g(E_i)| \int_{E_i}|\varphi(x)-\varphi(x_i)| \frac{f(dx)}{f(E_i)} \leq \eps\|g\|_X,
\]
as well as 
\[
\bigg|\langle g,\varphi\rangle - \sum_{i=1}^n \varphi(x_i)g(E_i) \bigg| \leq \sum_{i=1}^n \int_{E_i}|\varphi(x)-\varphi(x_i)| \,|g|(dx) \leq \eps\|g\|_X.
\]
We have established that $|\langle g_\eps- g,\varphi\rangle| \le 2 \eps \| g \|$ for any $\eps > 0$, from what  we deduce that $g$ belongs to the weak-$*$ closure of $X_f$.
\end{proof}

\Cyan 
We establish a technical result which will be useful in the next section. 

\Black

\begin{lem}\label{lem:Geo-QbQ>0}
 Assume \ref{X2}-\ref{X3}, and  $f \in X_{++}$. 
 Consider a linear operator $\QQ :X_{f} \to X_{f}$ such that $\QQ f = f$ and $A_{f} (\QQ g) \le A_{f}(g)$ for any $g \in X_{f}$.
 Then $\QQ \ge 0$. 
 \end{lem}
 
 
\begin{proof}[Proof of Lemma~\ref{lem:Geo-QbQ>0}.]
Take $0 \le g \in X_{f}$ such that $ g \le 2 C f$,  $C > 0$, and observe that 
  $$- C f \le g- Cf \le C f.
  $$
For any $r \in \R$, we compute 
\bean
A_f [ (\QQ  g )- Cf - ir f]  
&=&
A_f [ \QQ  (g - Cf - ir f)]  
\\
&\le& A_f [   g- Cf - ir f ] 
\\
&\le& \sqrt{C^2+r^2},
\eean
by using the non expansion property of $\QQ$ and the claim \eqref{eq:lem:Geo-QbQ>0}. Using again \eqref{eq:lem:Geo-QbQ>0}, 
we deduce $-Cf \le (\QQ  g )- Cf \le Cf$ and the conclusion. 
 \end{proof}

We generalize Kato's inequality \eqref{eq:KatoIneq} to the present complex framework by saying that an operator $\LL$ on 
$X$ satisfies (the complex) Kato's inequality if 
 \beqn\label{eq:Irred-KatoCineq}
\forall \, f \in D(\LL), \quad  \Re e (\hbox{sign} f) \LL f \le \LL |f|,
\eeqn
possibly in a dual sense as in \eqref{eq:KatoIneqDualsense}. As for the real  Kato's inequality, when $\LL$ is the generator of a semigroup,  
Kato's inequality \eqref{eq:Irred-KatoCineq} is a consequence of the positivity of the semigroup, and we refer to Remark~\ref{rem:sec4-EquivPositiveS} for references about this claim.

 \medskip
\subsection{On the subgroup and discrete structure of the {\Blue principal point spectrum}}
 \
 
 In this section, we establish that the {\Blue principal point spectrum} enjoys a subgroup structure under the same kind of hypotheses as considered in the previous sections.
 \Blue Such a deep result is classical and it is usually a central argument in the analysis of positive semigroups,  in particular  because it is used during the proof of Theorem~\ref{theo:NagelWebb}.
In the next sections,  we will develop alternative and more constructive approaches, but we will still make use of the material of this section or its consequences a couple of times in the applications (in particular in Section~\ref{ssec:GF:singular} and in the proof of Theorem~\ref{th:KFP-CvgceRate}). 
Here, we develop a method that is less abstract and a bit more general than the usual approach, although it remains somewhat technical.
%
\Black

\begin{lem}\label{lem:otherEV} 
Under assumptions \ref{C2},  \ref{X2} and the complex Kato's inequality \eqref{eq:Irred-KatoCineq}, {\Blue we have $\Sigma^1_P(\LL)\cap\Delta_{\lambda_1}=\emptyset$, or in other words $ \Sigma_P^1(\widetilde\LL) = \Sigma_P(\widetilde\LL) \cap i \R $.}
Furthermore, for any $\lambda \in \Sigma^1_{P}(\LL)$ the associated normalized eigenfunction $f$ satisfies $|f| = f_1$.
\end{lem}

\begin{proof}[Proof of Lemma~\ref{lem:otherEV}.]
{\Blue
Let $\lambda\in\C$ and $f\in D(\LL)\setminus\{0\}$ such that $\LL f = \lambda f$.}
By linearity of the operator sign and using \eqref{eq:irred:signCprop1} and Kato's inequality \eqref{eq:Irred-KatoCineq}, we have 
$$
 {\Blue (\Re e \lambda)} |f| = \Re e [ \lambda (\hbox{sign} f) f] = \Re e (\hbox{sign} f) (\lambda f) =  \Re e (\hbox{sign} f) \LL f \le \LL |f|.
$$
{\Blue Testing this inequality against $\phi_1\gg0$, we get
\[(\Re e \lambda)\langle\phi_1,|f|\rangle\leq\lambda_1\langle\phi_1,|f|\rangle\]
and so $\Re e\lambda\leq\lambda_1$.
This ensures that $\Sigma^1_P(\LL)\cap\Delta_{\lambda_1}=\emptyset$.
Now if $\lambda\in\Sigma^1_P(\LL)$, which imposes that $\Re e\lambda=\lambda_1$,} we deduce by
the duality argument introduced  during the proof of Lemma~\ref{lem:Uniquenessf1},  
{\Blue that the above inequality $\lambda_1 |f| \leq \LL |f|$   must be an equality. We conclude thanks to the simplicity of~$\lambda_1$.}
\end{proof}

 \begin{theo} \label{theo:Sigma+subgroup}  
Assume \ref{C2}, \ref{X2}, \ref{X3} and that complex Kato's inequality \eqref{eq:Irred-KatoCineq} holds true. 
Denoting $\widetilde \LL = \LL - \lambda_1$, 
 the set $ \Sigma_P^1(\widetilde\LL) \Blue= \Sigma_P(\widetilde\LL) \cap i \R $ is an additive subgroup of $i\R$ and $\hbox{\rm dim} N(\widetilde\LL-i\alpha)^k = 1$ for any $i\alpha \in \Sigma_P^1(\widetilde\LL)$ and $k \ge 1$.
 \end{theo} 
 
 Theorem~\ref{theo:Sigma+subgroup} is similar but more general than \cite[C-III, Cor.~2.12]{MR839450}  and \cite[Prop.~14.15]{MR3616245}. 
 Our proof is  also very similar to the proof of  \cite[ Prop.~14.15]{MR3616245}.
However, it is more direct and avoid the use of the $C(K)$ algebra and Kakutani's Theorem~\cite{MR5778} (see also \cite[Thm.~2.1.3]{MR1128093}). 

\begin{proof}[Proof of Theorem~\ref{theo:Sigma+subgroup}.]
We split the proof into three steps. 

\smallskip
{\sl Step 1. } 
 We consider $f$ associated to an eigenvalue $i \alpha \in \Sigma_P(\widetilde\LL) \backslash \{ 0 \}$, and we define 
 $$
 T(t) := (\sign f) e^{- i \alpha t} \widetilde S(t) (\sign \bar f).
 $$
Observing that $\widetilde S(t) f = e^{i \alpha t} f$, we have
 $$
 T(t) |f| = (\sign f) e^{- i \alpha t} \widetilde S(t) f=  (\sign f) f = |f| = \widetilde S(t) |f|.
 $$
On the other hand, we have
$$
|T(t) g| \le | \widetilde S(t) (\sign \bar f) g | \le \widetilde S(t) |g|, \quad \forall \, g \in X,
$$
which, by positivity of $\widetilde S(t)$, yields
\[|T(t) g| \le A_f(g) \widetilde S(t) |f| = A_f(g) |f|, \quad \forall \, g \in X_f.\]
Because $|f|=f_1\gg0$ from Lemma~\ref{lem:otherEV}, 
we can apply Lemma~\ref{lem:Geo-QbQ>0} to $|f|$ and $\QQ := T(t)$.
We deduce that $T (t) \ge 0$ on $X_{|f|}=X_f$,  and then on $X = \overline X_f$.
As a consequence, $0 \le T(t) g = |T(t) g| \le \widetilde S(t) g$ for any $g \ge 0$. 
In other words, we have 
$0 \le \widetilde S(t) - T(t)$ and then $0 \le \widetilde S(t)^* -  T(t)^*$. We must have $ \widetilde S(t)^* -  T(t)^* = 0$. Otherwise, there would exist $\psi \in Y_+ \backslash \{ 0 \}$ such that 
$(\widetilde S(t)^* -  T(t)^*)\psi \in Y_+ \backslash \{ 0 \}$, and we find a contradiction by computing 
$$
0 < \langle (\widetilde S(t)^* -  T(t)^*)\psi, f_1 \rangle = \langle \psi,  (\widetilde S(t) -  T(t)) f_1 \rangle = 0.
$$
We have thus established that $\widetilde S(t) = T(t)$.

\smallskip
{\sl Step 2. } Consider $\alpha, \beta \in \R$ and $f,g \in X\backslash \{ 0 \}$ such that $\widetilde \LL f = i \alpha f$ and $\widetilde \LL g = i \beta g$, and suppose first that $(\sign \bar f):D(\LL)\to D(\LL)$.
From Step 1 and the fact that $(\sign \bar f) \circ \sign f = I$,  for any $h \in D(\LL)$, we may compute
\bean
\widetilde\LL h &=& \lim_{t \to 0} \frac1t (\widetilde S(t) h - h)
\\
&=&  (\sign f) \lim_{t \to 0} \frac1t (e^{- i \alpha t} \widetilde S(t)  (\sign \bar f) h  -  (\sign \bar f)h) 
\\
&=&  (\sign f)  (\widetilde \LL - i \alpha)  (\sign \bar f) h ,  
\eean
or in other words $\widetilde \LL - i \alpha = (\sign \bar f) \widetilde \LL (\sign f)$. We have similarly $\widetilde \LL - i \beta =  (\sign \bar g) \widetilde \LL (\sign g)$.  Both equations together imply 
$$
\widetilde \LL - i (\alpha+\beta) = (\sign \bar f)  (\sign \bar g) \widetilde \LL  (\sign g)  (\sign f).
$$
Defining $h := (\sign \bar f) (\sign \bar g) f_1$, so that  $ (\sign g) (\sign f) h = f_1$, we get $\widetilde \LL h =  i (\alpha+\beta) h$, and finally $i(\alpha + \beta) \in\Sigma_P^1(\widetilde\LL)$, 
so that the  additive subgroup structure is established.

\smallskip
When the condition $(\sign \bar f):D(\LL)\to D(\LL)$ is not granted, we modify the above argument by using a resolvent approach.
For some $\lambda>0$ and  thanks to \eqref{eq:Exist1-DefRepresentationRR}, we compute
\begin{align*}
(\lambda-\widetilde\LL)^{-1} & = \int_0^\infty e^{-\lambda t}\widetilde S(t) \,dt \\
& = (\sign f) \int_0^\infty e^{-(\lambda+i\alpha) t}\widetilde S(t) \,dt \,  (\sign \bar f)\\
& = (\sign f) (\lambda+i\alpha-\widetilde\LL)^{-1}(\sign \bar f).
\end{align*}
Repeating the argument, we obtain 
\[ 
( \lambda + i (\alpha+\beta) - \widetilde \LL )^{-1} = (\sign \bar f)  (\sign \bar g)  (\lambda - \widetilde \LL)^{-1}  (\sign g)  (\sign f).
\]
Applying that last identity to the vector $h=(\sign \bar f) (\sign \bar g) f_1$ and using that $(\lambda - \widetilde \LL)^{-1}f_1=\lambda^{-1}f_1$, 
we deduce  $( \lambda + i (\alpha+\beta) - \widetilde \LL )^{-1}h=\lambda^{-1}h $. In other words, we have again $\widetilde \LL h =  i (\alpha+\beta) h$, and we conclude as before. 

\smallskip
{\sl Step 3. } From the fact that $(\sign f)$ is an invertible operator and  the equation 
$$ (\widetilde\LL - i \alpha)^k =  (\sign f)^{-1} (\widetilde\LL)^k (\sign f),$$
we see from Theorem~\ref{theo:KRgeometry1}-(ii) that 
$N(\widetilde \LL-i\alpha)^k = (\sign f)^{-1} N(\widetilde\LL)^k = (\sign f)^{-1} \hbox{\rm Span} f_1 $ for any $k \ge 1$,  so that its dimension is one. 
\end{proof}
 
 \medskip

Making an additional splitting structure hypothesis as yet introduced in Section~\ref{subsec-Exist1-KRtheorem}, we may significantly improve the conclusion. 
 We first recall a classical result on the spectrum of an operator  which holds under some power compactness assumption on the resolvent.

 \begin{theo}\label{theo:Voigt} 
 We assume that $\LL$ satisfies \ref{H1}-\ref{H2} and the splitting structure \ref{HS1} introduced in Section~\ref{subsec-Exist1-KRtheorem} with $\WW(z) \in \KK(X)$ for some $N \ge 1$ and any $z \in \Delta_{\kappa_0}$.
  Then $\Sigma(\LL) \cap \Delta_{\kappa_0} \subset \Sigma_d(\LL)$. 
  \end{theo}
  
 Theorem~\ref{theo:Voigt} is a consequence of \cite[Cor. 1.1]{MR595321}. 
We also refer to  \cite[proof of Thm.~3.1]{MR3489637} for a possible elementary proof. 

  \begin{proof}[A sketch of the proof of Theorem~\ref{theo:Voigt}] 
 Iterating the formula  $\RR_\LL = \RR_\BB + \RR_\BB \AA \RR_\LL$, we deduce 
 $$
\JJ(z)  R_\LL(z)  = \VV(z)
 $$
  with 
$
 \JJ :=  I -  (\AA \RR_\BB)^{N}$ and $ \VV := \RR_\BB + \dots +  \RR_\BB (\AA \RR_\BB)^{N-1}$.
 Because $\JJ$ is holomorphic on $\Delta_{\kappa_0}$, it is a compact perturbation of the identity and $\JJ(z) \to I$ when $\Re e z \to \infty$, one may use the theory of {\it degenerate-meromorphic functions} of Ribari\v{c} and Vidav \cite{MR236741}  (also established independently by Steinberg, see in particular  \cite[Cor.~1]{Steinberg1968}), and conclude that $\JJ(z)$ is invertible outside of a discrete set $\DD$ of $\Delta_{\kappa_0}$. 
 That implies that $\Sigma(\LL) \cap \Delta_{\kappa_0} = \DD$ is a discrete set of $\Delta_{\kappa_0}$. 
 On the other hand,  thanks to the Fredholm alternative  \cite{MR1554993}, one deduces that the eigenspace associated to each spectral value $\lambda \in \DD$ is non zero and finite dimensional, so that $\lambda \in \Sigma_d(\LL)$. See also \cite{MR1502771,MR0068733} for pioneering works in the subject. 
  \end{proof}

\smallskip
We end this section by a result which gives a more accurate description of the geometry of the boundary spectrum,  and is a variant of the classical results  \cite[C-III,  Thm.~3.12]{MR839450}, \cite[Thm. VI.1.12]{MR1721989}, \cite[Thm.~14.17]{MR3616245}. 

 \begin{theo} \label{theo:Sigma+subgroupBIS}
Under the assumptions of Theorem~\ref{theo:Sigma+subgroup} and Theorem~\ref{theo:Voigt}, 
{\Blue we have $ \Sigma_P^1 (\widetilde\LL) = \Sigma_P^+ (\widetilde\LL)$, this set}   is a {\bf discrete} additive subgroup of $i\R$ and any of its elements is an {\bf algebraically simple} eigenvalue. More precisely, 

- either $\Blue\Sigma_P^+(\widetilde \LL) = \{ 0 \}$ 
and the projection on the first eigenspace (associated to $\lambda_1$) writes
$$
\Pi\, f := \langle f, \phi_1 \rangle f_1; 
$$

- or $\Blue\Sigma_P^+(\widetilde \LL) = i\alpha \Z$ for some $\alpha > 0$ and there exists a sequence $(g_k,\psi_k)_{k\in\Z}$ such that 
$\LL g_k=(\lambda_1+ik\alpha)g_k$, $\LL^* \psi_k=(\lambda_1+ik\alpha) \psi_k$, and $\langle g_k, \psi_\ell\rangle=\delta_{k\ell}$.

 \end{theo}
 
 
  \begin{proof}[Proof of Theorem~\ref{theo:Sigma+subgroupBIS}] 
  {\Blue Because of Lemma~\ref{lem:H3abstract-StrongC} and \ref{C2}, we have $\lambda_1 = s(\LL)$ so that the principal point spectrum is the boundary point spectrum.}
Combining Theorem~\ref{theo:Sigma+subgroup}  and Theorem~\ref{theo:Voigt}, we immediately get that the subgroup $\SS :=  \Sigma_P(\widetilde\LL) \cap i \R $ satisfies $\SS \subset \Sigma_d(\LL)$, 
it is thus discrete and made of algebraically simple eigenvalues. The first case $\Sigma_P(\widetilde\LL) \cap i \R = \{0\}$ falls yet in the conclusions of Theorem~\ref{theo:KRgeometry1}. 
   \end{proof}
   
In the second case, where the boundary spectrum is not trivial, the existence of a projection on the boundary eigenspace $\overline\Span(g_k)_{k\in\Z}$ is ensured by the
 Jacobs--de Leeuw--Glicksberg theorem  provided that $\LL$ is the generator of a relatively compact semigroup, see for instance~\cite[Thm.~A.39 and Prop.~A.40]{MR3616245}
 and the references therein. We also refer to \cite[paragraphs III.6.4 and III.6.5]{MR0203473} for 
very classical results on the projector on the direct sum of eigenspaces associated to eigenvalues belonging to a subset of the spectrum. 
We can even give an explicit expression of this projection in terms of $(g_k)$ and $(\psi_k)$ under the form of a Fejér type sum, see Theorem~\ref{theo:periodic}.

 \medskip
\subsection{Stronger positivity}
 \

In order to go one step further and establish the triviality of the {\Blue principal point spectrum}, we need to reinforce the positivity of the semigroup or its generator. 
One possible condition is based on the following notion. 

 \medskip
{\bf The reverse strong positivity condition}
\smallskip


For $\Cyan A \in \BBB(X) \cap \BBB(X_+)$, we recall that from \eqref{eq;absAf<Aabsf}, we have 
\beqn\label{eq;absAf<AabsfBIS}
|Af| \le A|f|, \quad \forall \, f \in X, 
\eeqn
and we observe that the above inequality is an equality when $Af = u A|f|$ for some $u \in \Ss^1$. 
We focus now on  the case of equality in \eqref{eq;absAf<AabsfBIS}.

\begin{defin}\label{def:revserSPC} We say that $A$ satisfies the ``reverse strong positivity condition" for a subclass of vectors $\CC \subset X$ if for any $f \in \CC$
\beqn\label{eq:rspCondition}
  |A f| = A |f|\quad\hbox{implies}\quad \exists \, u \in \Ss^1, \ Af = u A|f|.
\eeqn
\end{defin}
We start observing that  $A > 0$ (as defined in Section~\ref{subsec:irreducibility})  implies the strict positivity for non-signed vectors in $X_\R$.

\begin{lem}\label{lem:strictpositiveOp} 
Consider an operator $A > 0$  and assume $X$ is reflexive.   
For $f \in X_\R$ such that $\pm f \notin X_+$, there holds
$$
|Af| \ll A|f|.
$$
 \end{lem}

\begin{proof}[Proof of Lemma~\ref{lem:strictpositiveOp}.]
Let us consider $f \in X_\R$ such that $f_\pm \not= 0$. We claim that $|Af| \ll A|f|$. Observing that 
$$
A f_+ = A f + Af_- \ge Af, 
$$
we deduce $Af_+ \ge (Af)_+$, and similarly $A f_- \ge (Af)_-$. 
We first consider the case  $(Af)_+ > 0$.
For $\phi \gg 0$, we have 
\bean
0 < \langle (A f)_+ , \phi \rangle = \sup_{0 \le \psi \le \phi} \langle Af , \psi \rangle =    \langle Af , \psi^* \rangle = \langle f , A^*\psi^* \rangle, 
\eean
for some $0 \le \psi^* \le \phi$, where we have used the very definition of $X_{++}$, the definition of $(Af)_+$ as an element of $X''$ and that $B_{X'}$ is compact for the weakly $*$ topology $\sigma(X',X)$. 
We deduce in particular that $\psi^* \not = 0$, so that $\psi^* > 0$ and finally $A^*\psi^* \gg 0$ because $A^* > 0$ (as an elementary consequence of the fact that $A >0$ listed  in Section~\ref{subsec:irreducibility}).
We deduce 
\bean
\langle (A f)_+ , \phi \rangle =  \langle f , A^*\psi^* \rangle < \langle f_+ , A^*\psi^* \rangle =\langle A f_+ , \psi^* \rangle \le \langle A f_+ , \phi \rangle
\eean
where for the strict inequality we have first used the assumption $f_- \not = 0$ and next elementary arguments. 
We thus have $(A f)_+ < Af_+$. Similarly, we establish $(A f)_- < Af_-$ when $(Af)_- > 0$. 

As a conclusion, in the three cases $Af = 0$, $(Af)_+ \not= 0$ and $(Af)_- \not= 0$, we have 
$$
|Af| = (Af)_+ + (Af)_- \ll Af_+ + Af_- = A|f|,
$$
which is the desired strict inequality.
\end{proof}

We believe that a similar result also holds true for complex vectors in a general Banach lattice framework. 
We do not try to prove such a statement but we rather establish the corresponding complex version
for our  examples of concrete Banach spaces in which the definition of the absolute value $|f|$ of a vector $f \in X$ is more tractable.
 
\begin{lem}\label{lem:strictpositiveOpTous} 
Consider an operator $A> 0$ on $X \subset \Lloc^1(E)$  for some locally and $\sigma$-compact metric space $E$. 
For $f \in X $ such that  $|f| \gg 0$, 
we have 
$$
|Af| = A|f| \quad\hbox{implies}\quad \exists\,u \in \Ss^1, \  f = u |f|, 
$$
  and thus \eqref{eq:rspCondition} holds. 
\end{lem}

\begin{proof}[Proof of Lemma~\ref{lem:strictpositiveOpTous}.]
We assume by contradiction that  $\forall v \in \Ss^1$, $|f| > \Re e (vf)$, in particular writing $f = g + i h$, $g,h \in X_\R$, we have $g,h \in X \backslash \{0\}$.
On the one hand, because of $A > 0$ and $A$ is linear, for any $v = e^{i \alpha } \in \Ss^1$, we have 
\bean
A|f|  \gg  A \Re e (e^{i\alpha} f)) 
=\cos\alpha \, (A  g)  - \sin\alpha \, (Ah).
\eean
On the other hand, in the Banach lattice we consider here, there exists $\beta : E \to \R$ measurable such that $|Af| = e^{i\beta} Af$ and thus 
$$
A |f| = |Af | = \Re e |Af| = \cos\beta \, (A  g)  - \sin\beta \,  (Ah), 
$$
and a contradiction. We have established that there exists $v \in \Ss^1$ such that $|f| \equiv \Re e (fv)$. 
Now, we have
$$
\sqrt{(\Re e (fv))^2 + (\Im m  (fv))^2}   = |fv| = |f| =  \Re e (fv),
$$
which in turn implies $\Im m  (fv)= 0$, since $\Re e (fv) \gg 0$.  That is here that we use the assumption $|f| \gg 0$ and not only $f \in X_+ \backslash \{ 0 \}$.  
We conclude that $|f| = fv$ and thus that  $f = u |f|$, with $u := v^{-1} \in \Ss^1$. 
\end{proof}

A similar result also holds in the Radon space of measures. For a  measurable space $(E,\EE)$, 
we call transition kernel, a mapping  $Q : E \times  \EE \to [0,\infty]$
such that 
\bean 
(i)&& \forall \, B \in \EE, \ x \mapsto Q(x,B) \ \hbox{is measurable}; 
\\
(ii)&& \forall \, x \in E, \ B \mapsto Q(x,B) \ \hbox{is a measure}. 
\eean
We recall the classical Markov-Riesz representation theorem which claims that for any 
linear and positive operator $B : C_0(E) \to C_0(E)$ there holds
$$
(B \phi)(x) = \int_E \phi(x) Q(x,dy) , \quad \forall \, \phi \in C_0(E),
$$
for a transition kernel $Q$ such that in the condition (i) above the mapping is furthermore continuous. 

\begin{lem}\label{lem:strictpositiveOpM1} 
Consider an operator $A> 0$ in $X = M^1 = M^1(E,\EE)$, for some 
Borel space $(E,\EE)$ where $E$ is a locally and $\sigma$-compact metric set. For $f \in X$ such that  $|f| \gg  0$, we have \eqref{eq:rspCondition}.
\end{lem}

\begin{proof}[Proof of Lemma~\ref{lem:strictpositiveOpM1}.]
By definition, the operator $A$ is the dual of a positive operator on $C_0(E)$. Using the representation formula recalled above for that adjoint operator, 
we get 
$$
(A f)(dy) = \int_E Q(x,dy) f(dx), \quad \forall \, f \in M^1,
$$
for a transition kernel $Q$. We deduce that 
$$
\langle A f, \phi \rangle =  \int_{E\times E} \phi(y) Q(x,dy) f(dx), 
$$
for any bounded Borel function $\phi : E \to \C$. In particular, the strict positivity $A > 0$ translates as $Q(x,\cdot) \gg 0$ in $M^1$ for any $x \in E$. 
We fix now $\phi \in C_0(E)$ such that $\phi \gg 0$ and $f \in M^1$ such that $|f| \gg 0$, and we observe that from 
the Radon-Nikodym theorem, 
there exist two measurable functions  $\alpha,\beta : E \to [0,2\pi)$ such that $f  = e^{i\alpha}  |f| $ and $ Af=  e^{i\beta} |Af| $. 
We next compute 
\bean
\langle A|f|  -  |Af|, \phi \rangle
&=&\Re e \bigl\{ \langle A|f|,\phi \rangle -  \langle Af, e^{-i\beta} \phi \rangle \bigr\}
\\
&=&  \int_{E\times E} \Re e \bigl\{  1-  e^{i(\alpha(x) - \beta(y))} \bigr\} \phi(y) Q(x,dy) |f|(dx)
\\
&=&  \int_{E\times E}   \bigl\{  1- \cos(\alpha(x)-\beta(y)) \bigr\} \phi(y) Q(x,dy) |f|(dx). 
\eean
In the case of equality $A|f|  =  |Af|$, we must have $1- \cos(\alpha(y)-\beta(x)) = 0$ for $\mu$-a.e. $x \in E$ and $|f|$-a.e. $y \in \hbox{supp}\, f=E$. 
We deduce that $\beta$ is a constant function, so that  $Af = e^{i\beta} |Af| = u A|f|$, for the constant $u = e^{i\beta} \in \Ss^1$. 
\end{proof}
\Black
  

{\bf The reverse Kato's inequality condition}
\smallskip

\smallskip\smallskip
We recall that it has been stated in section~\ref{subsec:MorePositive} that the generator $\LL$ of  a positive semigroup $S(t)$ satisfies 
Kato's inequality \eqref{eq:KatoIneq} which in a complex framework writes
$$
\forall \, f \in X, \quad  \Re e (\hbox{sign} f) \LL f \le \LL |f|. 
$$
We also observe that if $f = u |f|$ for some $u \in \Ss^1$, we have 
$$
 \Re e ( \sign f) \LL f = \sign(u^{-1}f) \LL (u^{-1}f) = \LL |f|,
$$
which is the case of equality in Kato's inequality.  

%

 \begin{defin}\label{defin:strongKato} 
 We say that $\LL$ satisfies a ``reverse Kato's inequality condition'' for a class of vectors $\CC \subset D(\LL)$ if for any $f \in \CC$ 
 the case of equality in Kato's inequality
$$
\LL |f| = \Re e (\hbox{\rm sign} f) \LL f  
$$
implies 
$$
 \exists \, u \in \C, \,\, f = u |f|.
$$
\end{defin}

In some situation, we may prove directly that the ``reverse Kato's inequality condition'' holds by reasoning at the level of the operator $\LL$, see for instance \cite[Proof of Thm. 5.1]{MR4265692}. 
That is also a consequence of the strong positivity of the semigroup as we see below.

\begin{lem}\label{lem:StrongKato&StrongPositivity} 
 Consider a semigroup $S$ and its generator $\LL$.
 On the set $\CC$ of vectors $f \in X\backslash\{0\}$ such that   
\beqn\label{eq:lemStrongKato&StrongPositivity}
\exists \, \lambda \in \C, \quad \LL f = \lambda f, \quad \LL |f|  = (\Re e \lambda) |f|,
\eeqn
 there is equivalence between:

(i) $S(t)$  satisfies the ``reverse strong positivity condition"  for some (and thus any) $t > 0$ {\Cyan as stated in Definition~\ref{def:revserSPC}; }

(ii) $\LL$ satisfies  the above ``reverse Kato's inequality condition''.

\end{lem}
 
\begin{rem}\label{rem:StrongKato&StrongPositivity}
When  $X \subset \Lloc^1$, the ``reverse Kato's inequality condition'' (ii) implies the ``reverse strong positivity condition" (i) on the class $\CC$  of vectors such that $f \in D(\LL)$, $0 \ll |f| \in D(\LL)$. 
Assume indeed that $\LL$ satisfies (ii) and consider $f \in \CC$ such that $|S_t f| = S_t |f| $ for any $t \ge 0$. By differentiating, we get 
\beqn\label{eq:signLLf=LLabsf}
(\sign f) \LL f = \LL |f|.
\eeqn
From the  ``reverse Kato's inequality condition'', we deduce that $f = u |f|$ for some $u \in \Sp^1$, so that (i) holds. 
\end{rem}

%

\begin{proof}[Proof of Lemma~\ref{lem:StrongKato&StrongPositivity}.]
In what follows, we fix $f \in X\backslash\{0\}$ such that \eqref{eq:lemStrongKato&StrongPositivity} holds, 
and we compute 
\beqn\label{eq:lem:StrongKato&StrongPositivity1}
\Re e (\sign f) \LL f = \Re e  (\sign f)(\lambda f) = ( \Re e\lambda ) |f| = \LL |f|.
\eeqn
For any $t > 0$, we also have $S_t f= e^{\lambda t} f$, $S_t |f| = e^{\Re e \lambda t} |f|$, and thus 
\beqn\label{eq:lem:StrongKato&StrongPositivity2}
 \quad |S_t f| = S_t|f|.
\eeqn
Assuming the  ``reverse Kato's inequality condition'', we deduce from \eqref{eq:lem:StrongKato&StrongPositivity1} that $f = u |f|$ for some $u \in \Ss^1$, thus
$S_t f = u S_t |f|$ for some $u \in \Ss^1$, which is the conclusion of the  ``reverse strong positivity condition" when \eqref{eq:lem:StrongKato&StrongPositivity2} holds. 

\smallskip
On the other way round, assuming  the  ``reverse strong positivity condition" for some $T> 0$, we deduce from  \eqref{eq:lem:StrongKato&StrongPositivity2} for $T >0$ that there exists $v \in \Ss^1$ such that 
$$
e^{\lambda T} f = S_T f = v S_T |f| = v e^{\Re e \lambda T} |f|.
$$
That implies that $f= u |f|$ with $u = v e^{- i (\Im m \lambda) T}$,  
which is nothing but the conclusion of the  ``reverse Kato's inequality condition'' when \eqref{eq:lem:StrongKato&StrongPositivity1} holds.
\end{proof}

We summarize the material developed above in the following main result  of the section. 

\begin{theo}\label{theo:StrongKato&StrongPositivity} Assume that $S$ is a positive semigroup on $X$ with $X \subset \Lloc^1(E)$ or 
$X=M^1(E)$ for 
some locally and $\sigma$-compact metric space $E$ and denote by $(E_k)$ a sequence of increasing compact sets such that $E=\lim E_k$.  We furthermore assume that for any $k \ge 1$ there exists $T > 0$ such that $S_T$ is strictly positive on $E_k$, in the sense that 
\beqn\label{eq:StrongKato&StrongPositivityEk}
\forall \, f \in X_+ \backslash \{0\},  \     f_{|E_k} \not\equiv 0, 
\ \forall \, \phi \in X'_+\backslash \{0\}, \   \hbox{\rm supp}\,  \phi \subset E_k, \quad \langle S_T f, \phi \rangle > 0. 
\eeqn
Then $\LL$ satisfies  the ``reverse Kato's inequality condition'' on the set $\CC$ of eigenvectors introduced in Lemma~\ref{lem:StrongKato&StrongPositivity}.
\end{theo}

\begin{proof}[Proof of Theorem~\ref{theo:StrongKato&StrongPositivity}.]
Let us consider $f \in X\backslash\{0\}$ such that   \eqref{eq:lemStrongKato&StrongPositivity} holds, so that 
$S_{t} |f| = e^{(\Re e \lambda) t} |f|$ for any $t \ge 0$. 
On the one hand, we may  fix $k \ge 1$ such that $|f| \not\equiv 0$ on $E_k$. Then for any $\ell \ge k$, there exists $T_\ell > 0$ such that \eqref{eq:StrongKato&StrongPositivityEk} holds, 
so that 
$$
e^{(\Re e \lambda) T_\ell} \langle |f|, \phi \rangle = \langle S_{T_\ell} |f|, \phi \rangle > 0,
$$
for any $\phi \in Y_+\backslash \{0\}$, $\hbox{\rm supp}\,  \phi \subset E_\ell$. 
That implies $\langle |f| , \phi \rangle > 0$ on for any $\phi \in Y_+ \backslash \{0\}$, and thus $|f| \gg 0$. 
Next, as in the proof of Lemma~\ref{lem:StrongKato&StrongPositivity}, we  observe that 
$$
 |S_{T_\ell} f| = S_{T_\ell} |f|, \quad \forall \, \ell \ge k.  
$$
Repeating the proof of Lemma~\ref{lem:strictpositiveOpTous} and Lemma~\ref{lem:strictpositiveOpM1}, 
we deduce that there exists $u_\ell \in \Sp^1$ such that $S_{T_\ell} f = u_\ell S_{T_\ell} |f|$ on $E_\ell$, or equivalently  there exists $v_\ell \in \Sp^1$ such that $  f = v_\ell   |f|$ on $E_\ell$, 
with $v_\ell := u_\ell e^{-i (\Im m \lambda) T_\ell}$. Because $E_\ell \supset E_1$, we have established that $  f = v_1   |f|$ on $E$ 
which is the conclusion of the  ``reverse Kato's inequality condition'' when 
\eqref{eq:lemStrongKato&StrongPositivity} holds.
\end{proof}
 \Black

%
%

\subsection{On the triviality  of the  principal point spectrum}
\label{subsec:bdarySpectrum}

{\Cyan We assume here
 the existence \ref{C2} of a unique solution $(\lambda_1, f_1,\phi_1) \in \R \times X_{++} \times Y_{++}$
to the  first eigenvalue problem \eqref{eq:triplet1}-\eqref{eq:triplet2}}
and that $\LL$ enjoys the weak maximum principle \eqref{eq:Irred-weakPM} and Kato's inequalities \eqref{eq:Irred-Katoineq} as formulated  in condition 
\ref{H1'} as well as  the strong maximum principle \ref{H4}.  Because we deal with complex eigenvalue, we also assume that 
 the {\bf complex}  Kato's inequality variant {   \eqref{eq:Irred-KatoCineq}} holds. 

\smallskip
  We introduce  a first  additional assumption: 
 
\begin{enumerate}[label={\bf(H5)},itemindent=13mm,leftmargin=0mm,itemsep=1mm]
\item\label{H5} the``reverse Kato's inequality condition'' (as defined in Definition~\ref{defin:strongKato}) holds true for the class $\CC$ defined in Lemma~\ref{lem:StrongKato&StrongPositivity}:  for $f \in X\backslash\{0\}$ such that   
\beqn\label{eq:lemStrongKato&StrongPositivity3}
\exists \, \lambda \in \C, \quad \LL f = \lambda f, \quad \LL |f|  = (\Re e \lambda) |f| = \Re e (\hbox{\rm sign} f) \LL f, 
\eeqn
we have 
$$
 \exists \, u \in \C, \,\, f = u |f|.
$$ 
\end{enumerate}
On the other hand, we do not need the structure assumption \ref{X3}.  

\smallskip
We are then able to make a more accurate analysis of the geometry of the spectrum.

 \begin{theo}\label{theo:KRgeometry2} 
Consider an  unbounded operator $\LL$ on a Banach lattice $X$ which satisfy {\Blue \ref{X1}, \ref{X2}},  \ref{C2},  \ref{H4}, \eqref{eq:Irred-KatoCineq} and \ref{H5}.
{\Blue Then $\lambda_1$ is the unique eigenvalue of $\LL$ with largest real part: 
$\Sigma^1_+(\LL) =  \{ \lambda_1 \}$.}
 \end{theo}

 \begin{rem}\label{rem:KRgeometry2}
 (1) 
It is worth emphasizing  again  that  {  \eqref{eq:Irred-KatoCineq}} is true when $\LL$ is the generator of a
positive semigroup and that \ref{H5} \Black 
 is true when 
$S_\LL(T)$  satisfies the ``reverse strong positivity condition"  for some $T > 0$ as a consequence of Lemma~\ref{lem:StrongKato&StrongPositivity}, see also Theorem~\ref{theo:StrongKato&StrongPositivity}.

\smallskip
(2) During the proof we use similar arguments as in \cite[Thm.~5.1]{MR4265692}. 

\smallskip

(3) Condition \ref{H5} is reminiscent of PDE arguments as we may find for instance in \cite[Proof of Thm.~5.1]{MR4265692} or in the discussion in \cite[4th course]{PL2} about an uniqueness argument due to L.~Tartar.


 \end{rem}

\begin{proof}[Proof of Theorem~\ref{theo:KRgeometry2}.]
Consider an eigenvalue $\lambda \in \C$ with normalized eigenvector $f \in X \backslash \{0\}$, and more precisely
$\langle |f|, \phi_1 \rangle = 1$ and $\LL f = \lambda f$. Thanks to the complex Kato's inequality \eqref{eq:Irred-KatoCineq}, we have 
\bean
(\Re e \lambda) |f| =  \Re e  \,  \hbox{sign} (f) (\lambda f) =  \Re e \,   \hbox{sign} (f) (\LL f ) \le \LL |f|. 
\eean
We consider two cases: 

\smallskip
When the above inequality is not an equality, we have
\bean
(\Re e \lambda) \langle |f|,\phi_1 \rangle <   \langle \LL |f| ,\phi_1 \rangle  = \langle  |f| , \LL^* \phi_1 \rangle  = \lambda_1 \langle |f|,\phi_1 \rangle,
\eean
and thus $\Re e \lambda < \lambda_1$.

\smallskip
When on the contrary the above inequality is an equality,  then $|f|$ is a positive eigenvector associated to the eigenvalue $\Re e \lambda$. 
Because of \ref{H4}, we have $|f| \in X_{++}$ and {\Blue  we get  $\Re e \lambda =\lambda_1$ thanks to   Lemma~\ref{lem:PositiveEigenvector}}. 
The condition   \ref{C2} implies $|f| \in \hbox{\rm Span}(f_1)$. 
On the other hand, $f$ satisfies \eqref{eq:lemStrongKato&StrongPositivity3} and thus $f \in \hbox{\rm Span}(f_1)$ from assumption \ref{H5}, in particular $\lambda = \lambda_1$. 
 \end{proof}

When $\LL$ is the generator of a positive and irreducible semigroup $S$, we may introduce 
the  alternative  assumption: 
 
\begin{enumerate}[label={\bf(H5$'$)},itemindent=14mm,leftmargin=0mm,itemsep=1mm]
\item\label{H5'} the semigroup $S$ is aperiodic as defined in~\eqref{eq:Irred-eventualirreducibilityBIS}, namely
$$
\forall \, f \in X_+\backslash \{0\}, \forall \, \phi \in Y_+\backslash \{0\}, \ \exists\, T>0, \,\forall \tau \geq T \quad
\langle S_\tau f, \phi \rangle > 0.
$$
\end{enumerate}
 
  \begin{theo}\label{theo:lambda1largestBIS} {\Blue 
Let $X$ be a Banach lattice satisfying conditions \ref{X1}, \ref{X2} and  in which the property~\eqref{eq:Irred:StrictOrderNorm} holds true with $\EEE =X$.
Consider a positive and  irreducible semigroup $S$ on $X$ which satisfies the aperiodicity condition \ref{H5'} and such that its generator $\LL$ satisfies \ref{C2}.
 Then $\lambda_1$ is the unique eigenvalue of $\LL$ with largest real part:
$\Sigma^1_P(\LL) =  \{ \lambda_1 \}$. }


%

\end{theo}
 
 
\begin{rem}\label{rem:H5'} 
{\Blue  (1) Because of the restriction made on the Banach lattice structure, this framework excludes the case when $X$ is a (weighted) signed measures space but includes the cases when $X$ is a (weighted)  Lebesgue space.  } 

(2) It is worth pointing out that since \ref{H5'} is stronger than \ref{H4}{\bf-(i)}, see the points {\bf(2)} and {\bf(3)} in Lemma~\ref{lem:Irred-S>0impliesR>0},
we can use Theorem~\ref{theo:KRgeometry1} and replace in Theorem~\ref{theo:lambda1largestBIS} the assumption that \ref{C2} is satisfied by the assumption that \ref{C1}  and \ref{H1'} for both $\LL$ and $\LL^*$ are satisfied. 

%


\end{rem}

\begin{proof}[Proof of  of Theorem~\ref{theo:lambda1largestBIS}]
We introduce the notations  $\widetilde S_t := S_t e^{-\lambda_1 t}$ and  $\widetilde \LL := \LL - \lambda_1$. {\Blue The complex Kato's inequality \eqref{eq:Irred-KatoCineq} automatically holds because $\LL$ is the generator of a positive semigroup. We may thus  use Lemma~\ref{lem:otherEV} and we deduce that $ \Sigma_P^1(\widetilde\LL) = \Sigma_P(\widetilde\LL) \cap i \R $.}
Assume that $f = g + ih \in X$, $g,h \in X_\R$, is an eigenfunction associated to the eigenvalue $\lambda = \lambda_1 + i \alpha \in \C$, $\alpha > 0$, so that 
$$
\widetilde\LL (g + i h) = i \alpha (g+ih) = \frac{2\pi i }{ t_0} (g+ih),
$$
for some $t_0 > 0$. 
On the one hand, because $\alpha \not=0$, we must have $g \not=0$ and $h \not=0$, and because of
$$
\alpha \langle g, \phi_1 \rangle =   \langle \widetilde\LL h, \phi_1 \rangle =\langle h,  \widetilde\LL^*\phi_1 \rangle = 0,
$$
and  $\phi_1 \gg 0$,  we have $g_+ \not= 0$ and $g_- \not= 0$. 
As a consequence, and because of \eqref{eq:Irred:StrictOrderNorm},
there exists $\psi \in Y_+ \backslash \{ 0 \}$ such that $\langle g_+,\psi \rangle = 0$. 
On the other hand, we compute 
$$
\widetilde S_{t_0} (g+ih) = e^{i \alpha t_0} (g+ih) = g+i h,
$$
from what we deduce  $\widetilde S_{t_0} g  =  g$, 
because $S_t$ is real. On the other hand, because  $  S_t$ is positive, we have 
$$
g_+ = (\widetilde S_{t_0} g)_+ \le \widetilde S_{t_0} g_+, 
$$
and next 
$$
\langle \phi_1, g_+ \rangle \le \langle  \phi_1, \widetilde S_{t_0} g_+ \rangle = \langle  \widetilde S^*_{t_0}  \phi_1 , g_+\rangle 
= \langle   \phi_1, g_+ \rangle,
$$
so that both inequalities are equalities (remind again that $\phi_1 \gg 0$), and thus
$$
\widetilde S_{t_0} g_+  =  g_+. 
$$
We conclude that
$$
\langle \widetilde S_{k t_0} g_+,\psi \rangle = \langle  g_+,\psi \rangle = 0, \quad \forall \, k \ge 0, 
$$
what is in contradiction with  \ref{H5'}. 
\Blue We have established that $\Sigma_P^1(\LL) = \{\lambda_1\}$.  
\end{proof}

\medskip

We end this section with a third situation where the triviality of the {\Blue principal point spectrum} is an  immediate consequence of Theorem~\ref{theo:Sigma+subgroup} and Theorem~\ref{theo:Voigt}. 

\begin{theo}\label{theo:lambda1largestTER}  (1) We make the same assumptions as in Theorem~\ref{theo:Sigma+subgroup} and also that there exists $M>0$ large enough such that $\lambda-\LL$ is invertible in $\BBB(X)$ for any $\lambda \in \C$, $\Re e \lambda = \lambda_1$, $|\lambda| \ge M$. 
Then  $\lambda_1$ is the unique eigenvalue with largest real part: {\Blue $\Sigma^1_P(\LL) =  \{ \lambda_1 \}$. }

\smallskip
(2) We furthermore assume that the hypothesis of Theorem~\ref{theo:Voigt} are met and that  $\lambda-\LL$ is invertible in $\BBB(X)$ for any   for any $\lambda \in \C$, $\Re e \lambda \ge \lambda_1-\eps$, $|\lambda| \ge M$.
Then a (non constructive) spectral gap \ref{S33} holds.  
\end{theo}
\Black

{\Blue
\begin{proof}[Proof of Theorem~\ref{theo:lambda1largestTER}]
(1) Combining the property of $\Sigma_P^1(\widetilde\LL)$ established in Theorem~\ref{theo:Sigma+subgroup} and the above invertibility condition, we see that $\Sigma_P^1(\widetilde\LL)$ is a bounded additive subgroup of $i\R$, so that $\Sigma_P^1(\widetilde\LL) = \{0\}$, what is precisely the first claim.  

(2) From Theorem~\ref{theo:Voigt}, we know that  $\Sigma(\LL) \cap \Delta_{\lambda_1 - \eps} \subset \Sigma_d(\LL)$ and thus together with (1), we get $\Sigma(\LL) \cap \Delta_{\kappa}  = \{ 0\}$ for $\lambda_1-\kappa > 0$ small enough,  what is precisely the second claim.  
\end{proof}
}

\medskip
  We summarize the main results established in this section as follows. 

\begin{enumerate}[label={\bf(C3)},itemindent=13mm,leftmargin=0mm,itemsep=1mm]
\item\label{C3} the first eigentriplet problem {\Cyan\eqref{eq:triplet1}-\eqref{eq:triplet2}} 
has a  solution $(\lambda_1, f_1, \phi_1)$,   furthermore this one is unique,  $f_1 \gg 0$,  $\phi_1 \gg 0$, $\lambda_1$ is algebraically simple (for both $\LL$ and $\LL^*$)
and {\Blue$\Sigma_{P}^1(\LL) =  \{\lambda_1\}$.}
\end{enumerate}

{\Blue We notice that~\ref{C3} is actually the combination of   \ref{C2}, \ref{S32} and the algebraic simplicity. 
}


 \medskip
\subsection{Ergodicity}
\label{subsec:ergodicity}
 
Thanks to the above analysis, we are able to formulate some convergence results on the trajectories associated to a semigroup.
More precisely, assuming  the existence and uniqueness of the first eigentriplet $(\lambda_1,f_1,\phi_1)$ for the generator $\LL$ of a semigroup $S$ and still denoting the rescaled semigroup $\widetilde S(t) := e^{-\lambda_1 t} \, S(t)$, we wish in particular  to establish the following ergodic property

\begin{enumerate}[label={\bf(E2)},itemindent=13mm,leftmargin=0mm,itemsep=1mm]
\item\label{E2} for any $f \in X$, there holds
\beqn\label{eq:theo:Ergodicity}
 \widetilde S_t f \to \langle f, \phi_1 \rangle f_1,  \ \hbox{ as } \ t \to \infty, 
\eeqn
in a sense to be specified.
\end{enumerate}

\smallskip

We start with a simple result which takes advantage of some dissipativity property of the semigroup formulated by a "reverse positivity condition". 
We next present some more involved results which use directly the spectral information. It is worth emphasizing that our results in this section
do not use any spectral gap property what contrasts with the results we will present in the next section.

\begin{theo}\label{theo:ErgodicityByLiapunov1} Consider a positive semigroup $S$ on a Banach lattice $X$ such that its generator $\LL$ enjoys the conclusions \ref{C2} of existence, uniqueness and strict positivity of the first eigentriplet $(\lambda_1,f_1,\phi_1)$ and let us set $\widetilde S_t := e^{-\lambda_1 t} S_t$. We denote 
$\XX$ the space $X$ endowed with the norm $[\cdot]$, with  $[f] := \langle |f|, \phi_1 \rangle$. 
Assume furthermore that 

  {\bf (1)}     for any $f \in X$, the trajectory $(\widetilde S_t  f)_{t \ge 0}$ is continuous in  $\XX$ and belongs to a compact set of a normed space $\XX_1$, with $\XX_1 \subset \XX$;
 



{\bf (2)} $(S_t)$ satisfies the reverse positivity condition for semigroups 
 \beqn\label{eq:reverseSPCsG}
|S_t f| = S_t |f|, \ \forall \, t > 0, \quad\hbox{implies}\quad \exists \, T > 0, \ \exists \, u_T  \in \Ss^1, \ S_T f = u_T S_T |f| . 
\eeqn

Then, the ergodicity property \ref{E2} holds  in the sense of the norm of $\XX_1$. 
 
\end{theo}
 
Let us comment on hypotheses made in the statement of Theorem~\ref{theo:ErgodicityByLiapunov1}. 
Hypothesis {\bf (1)} can be obtained as a consequence of a Lyapunov (or growth) condition reminiscent of the structure condition \ref{HS3} introduced in Section~\ref{subsect:AboutWeakDissip} and  an irreducibility condition.  
Typically, we assume first
$$
\|  \widetilde S(t) f \| \le M \| f \| + K \sup_{0 \le \tau \le t} [  \widetilde S(\tau) f ]_0, 
$$
with $[g]_0 := \langle |g|, \psi_0 \rangle$, $\psi_0 \in Y_{+} \backslash \{ 0 \}$, what can be established under the very general condition (ii) of Theorem~\ref{theo:KRexistTER}. Next we need to be able to prove that  $\psi_0 \le r \phi_1$ for some $r > 0$. In concrete situations, we may take $\psi_0$ with compact support and then the above  inequality is a consequence of the standard  strong maximum principle. We deduce
\bean
\|  \widetilde S(t) f \| 
&\le& M \| f \| + K r \sup_{0 \le \tau \le t} \langle | \widetilde S(\tau) f |,\phi_1\rangle   
\\
&\le&  M \| f \| + K r \sup_{0 \le \tau \le t} \langle  \widetilde S(\tau) |f |,\phi_1\rangle 
\\
&=&  M \| f \| + K r \langle   |f |,\phi_1\rangle,
\eean
so that  $(\widetilde S_t )$ is bounded. The hypothesis {\bf (1)} is in fact a bit more demanding, but also quite natural. Assume that $S_\LL$ enjoys the splitting structure introduced in section~\ref{subsect-Exist2-Dissip} and section~\ref{subsect:AboutWeakDissip}, 
so that 
\beqn\label{eq:ergodicity-splitting}
 \widetilde S = V + K,
\eeqn
with 
$$
V :=  \widetilde S_\BB + \dots+ ( \widetilde S_\BB \AA)^{*(N-1)} * \bar   S_\BB, 
\quad 
K :=    ( \widetilde S_\BB \AA)^{(*N)} *  \widetilde S, 
\quad
 \widetilde S_\BB (t) := e^{-\lambda_1 t} S_\BB(t). 
$$
In some applications, we typically have 
$$
\| V(t) f_0 \| \le \Theta(t) \| f_0 \|, \quad  \| ( \widetilde S_\BB \AA)^{(*N)} \|_{\Cyan\BBB(X,\XX_1)}  \le \Theta
$$
with $\Theta \in L^1(\R_+) \cap C_0(\R_+)$, $\XX_1 \subset X$ compact. In that situation, we deduce {\bf (1)}. 


 
\begin{proof}[Proof of Theorem~\ref{theo:ErgodicityByLiapunov1}.]
We fix $f \in X$ and without loss of generality, we may assume that $ \langle f, \phi_1 \rangle = 0$. 
We observe that  
\beqn\label{eq:ergo1-nonexpans}
 \langle |\widetilde S_t f|,\phi_1 \rangle =  \langle | \widetilde S_{t-s}  \widetilde S_s f|,\phi_1 \rangle \le \langle  \widetilde S_{t-s} |  \widetilde S_s f|,\phi_1 \rangle  =  \langle  |  \widetilde S_s f|,\phi_1 \rangle , 
\eeqn
 for any $t \ge s$. 
  We deduce that $(\widetilde S_t)$ is a dynamical system with compact trajectories in $\XX_1$ 
 and $\HH(f) := \langle |f|,\phi_1 \rangle$ is a Lyapunov functional.  
As a consequence, from the La Salle invariance principle, we have 
\beqn\label{eq:theo:Ergodicity2}
\inf_{g \in \omega_\HH} \langle |\widetilde S_t  f - g|, \phi_1 \rangle   \to 0 \ \hbox{ as } \ t \to \infty,
\eeqn
with 
\beqn\label{eq:defomegaH}
\omega_\HH := \{ g \in X; \  
\langle g, \phi_1 \rangle = 0,\ \forall \, t \in \R, \ \HH(\widetilde S_t  g) = \inf_{s > 0} \HH( \widetilde S_s f)\}. 
\eeqn
We next characterize $\omega_\HH$. Picking up $g \in \omega_\HH$,  we observe that 
$$
\langle |\widetilde S_t  g |,\phi_1 \rangle = \langle | g |,\phi_1 \rangle =  \langle | g |,  \widetilde S^*_t\phi_1 \rangle = \langle \widetilde S_t |g |,\phi_1 \rangle, \quad \forall \, t \ge 0, 
$$
so that 
$$
\langle \widetilde S_t  |g | - |\widetilde S_t  g |,\phi_1 \rangle = 0,\quad \forall \, t \ge 0. 
$$
In particular, using that $|\widetilde S_t  g | \le \widetilde S_t  |g|$, we have
\beqn\label{eq:Stabsg=absStg}
\widetilde S_t  |g| = |\widetilde S_t  g|, \quad \forall \, t \ge 0.
\eeqn

\smallskip

Because of the reverse positivity condition for semigroups  \eqref{eq:reverseSPCsG}, there exist $T > 0$ and 
$u_T \in \Ss^1$ such that 
$$
\widetilde S_T g =  u_T \widetilde S_T  |g|.
$$
As a consequence, by definition of the set $\omega_\HH$, we have 
$$
0 = \langle g,\phi_1 \rangle = \langle \widetilde S_T g ,\phi_1 \rangle =  u_T \langle \widetilde S_T  |g| ,\phi_1 \rangle 
= u_T \langle  |g| ,\phi_1 \rangle.
$$
Because $u_T \not = 0$, we conclude that $g = 0$. In other words, we have established that $\omega_\HH = \{ 0 \}$ and together with \eqref{eq:theo:Ergodicity2}, we obtain  \eqref{eq:theo:Ergodicity}. 
\end{proof}

We present a more concrete situation where the previous result can be invoked.
Although the hypotheses are somehow restrictive, it is yet useful in many applications and its proof is very simple. 

\begin{cor}\label{cor:Ergodicity2} Consider a strongly continuous and positive semigroup $S$ on a Banach lattice $X$ such that its generator $\LL$ enjoys the conclusions \ref{C2} of existence, uniqueness and strict positivity of the first eigentriplet $(\lambda_1,f_1,\phi_1)$. Assume further that  the reverse Kato's inequality condition (as defined in Definition~\ref{defin:strongKato}) holds true for the (large) class 
$$
\CC := \{ f \in D(\LL);  \  \LL |f|  =  \Re e (\hbox{\rm sign} f) \LL f \},
$$
that  $X \subset \Lloc^1(E,\EE,\mu)$ and that  the space $\XX^k$ defined in \eqref{eq1:rem:MeanErgodicityVariante1} satisfies $\XX^k \subset \Lloc^1$ with strongly compact embedding for some $k\ge1$.
Then the ergodicity property \ref{E2} holds  in the sense of strong topology of $L^1_{\phi_1}$.

\end{cor}

\begin{proof}[Proof of Corollary~\ref{cor:Ergodicity2}.]
Because of Step 3 in the proof  Theorem~\ref{theo:MeanErgodicityVariante1}, we see that condition (1) in  Theorem~\ref{theo:ErgodicityByLiapunov1} holds with $\XX_1 := \XX^k$. 
On the other hand, because of Remark~\ref{rem:StrongKato&StrongPositivity} and  the reverse Kato's inequality condition in $\CC$, we see that condition (2) also holds, so that we may apply  Theorem~\ref{theo:ErgodicityByLiapunov1} and conclude.
\end{proof}

 We present now a variant of the previous result which provides a convergence for various topologies, and relies on the (very general) assumption that the {\Blue  principal point spectrum} is trivial rather than on the reverse positivity condition.

 \begin{theo}\label{theo:ergodicity-compact-trajectories}
 Consider a positive semigroup $S$ on a Banach lattice $X$ such that its generator $\LL$ enjoys the conclusions \ref{C3} on the existence, uniqueness and strict positivity of the  first eigentriplet problem $(\lambda_1,f_1,\phi_1)$ and triviality of the boundary point spectrum.
Setting $\widetilde S_t  := e^{-\lambda_1 t} S_t$, we assume that we are in one of the following situations:
\begin{enumerate}[itemindent=*,leftmargin=*]
\item $S$ is strongly continuous and the trajectories $(\widetilde S_t  f)_{t \ge 0}$ are relatively compact for all $f\in X$, and we denote by $\TTT$ the strong topology of $X$;
\item $X=Y'$, $Y$ separable, and the trajectories $(\widetilde S_t  f)_{t \ge 0}$ are bounded for all $f\in X$, and we denote by $\TTT$ the weak $*\, \sigma(Y',Y)$ topology;
\item $X \subset \Lloc^1(E,\EE,\mu)$, and we denote by $\TTT$ the weak topology of $L^1_{\phi_1}$;
\item $X \subset \Lloc^1$, $S$ is strongly continuous, and for some $k \ge 1$ the space $\XX^k$ defined in \eqref{eq1:rem:MeanErgodicityVariante1} satisfies $\XX^k \subset \Lloc^1$ with strongly compact embedding, and we denote by $\TTT$ the strong topology of $L^1_{\phi_1}$.
\end{enumerate}
Then the ergodicity property \ref{E2} holds  in the sense of the topology $\TTT$. 
 \end{theo}
 
 \begin{rem}
The case {\it(4)} of Theorem~\ref{theo:ergodicity-compact-trajectories} enjoys some strong similarities with the main consequences of the General Relative Entropy technique developed in \cite{MR2162224}, 
see in particular  \cite[Thm.~3.2]{MR2162224}, \cite[Thm.~4.3]{MR2162224} and \cite[Thm.~5.2]{MR2162224}.
In particular, the aperiodicity condition that the boundary point spectrum is trivial may be compared with the fact that the first eigenvector $f_1$ is the unique (normalized and nonnegative) vector $f \in X$ with vanishing dissipation of entropy $\DD(f) = 0$ as defined in  \cite{MR2162224} or more generally  that $\Span(f_1)$ is the unique invariant space on which the functional $\DD$ vanishes. 
The present formulation is more abstract and probably more general. The drawback is the condition  $\XX^k \subset \Lloc^1$ with strongly compact embedding which can be avoided in \cite{MR2162224}, 
by using some weak compactness argument and the lower semicontinuity  property of $\DD$. That is explained by the fact that  our proof uses rather the La Salle invariance principle (similarly as in the proof of \cite[Thm.~3.2]{MR2114413}) instead of an entropy dissipation argument. 
\end{rem}
 
In the case when the boundary point spectrum is not trivial but a discrete set, the same method of
proof as for Theorem~\ref{theo:ergodicity-compact-trajectories} allows us to accurately describe the periodic long time behaviour of the semigroup.
 
\begin{theo}\label{theo:periodic} Consider a positive semigroup $S$ on a Banach lattice $X$ such that its generator $\LL$ enjoys the conclusions \ref{C2} on the existence and uniqueness of the  first eigentriplet problem $(\lambda_1,f_1,\phi_1)$,
and satisfies the complex Kato's inequality \eqref{eq:Irred-KatoCineq}.
Suppose furthermore that $X$ and $Y$ both enjoy the structure conditions \ref{X2} and \ref{X3}, that $\lambda_1$ is an isolated eigenvalue
and that the boundary spectrum is not trivial, {\it i.e.} $\Sigma_P^+\neq\{\lambda_1\}$.
Setting $\widetilde S_t  := e^{-\lambda_1 t} S_t$, we assume that we are in one of the situations (1), (2), (3) or (4) listed in statement of Theorem~\ref{theo:ergodicity-compact-trajectories}.
Then $\Sigma_P^+=\{\lambda_1+ik\alpha,\ k\in\Z\}$ for some $\alpha>0$, there exists a sequence $(g_k,\psi_k)_{k\in\Z}$ such that $\LL g_k=(\lambda_1+ik\alpha)g_k$, $\LL^*\psi_k=(\lambda_1+ik\alpha)\psi_k$ and $\langle g_k,\psi_k\rangle=1$,
and for all $f\in X$, in the sense of the topology $\TTT$, the projection
\[\Pi f=\lim_{n\to\infty}\frac1n\sum_{\ell=0}^{n}\sum_{k=-\ell}^\ell\langle f,\psi_k\rangle g_k\]
is well defined and
\[\widetilde S_t f-\widetilde S_t  \Pi f\to0\qquad \text{as}\ t\to+\infty.\]
\end{theo}

\begin{rem}\label{rem:periodic}
In Theorem~\ref{theo:periodic}, the assumptions that $\lambda_1$ is isolated and $\Sigma_P^+\neq\{\lambda_1\}$ might seem difficult to check in practice.
We indicate here some ways to verify them.

(i) The condition that $\lambda_1$ is an isolated eigenvalue is for instance guaranteed under the assumptions of Theorem~\ref{theo:Voigt}  or Theorem~\ref{theo:Harris-mean}.

(ii) The condition that $\Sigma_P^+$ is not restricted to $\{\lambda_1\}$ can be guaranteed by verifying that \ref{E2} does not hold.
Indeed, if $\Sigma_P^+=\{\lambda_1\}$, then Theorem~\ref{theo:ergodicity-compact-trajectories} imposes \ref{E2} to hold.
\end{rem}
 
 The result in Theorem~\ref{theo:periodic} can be compared for instance to~\cite[Thm.~14.19]{MR3616245}, although our hypotheses are slightly more general. 
Our proof is also more direct than in~\cite{MR3616245} and it additionally provides an explicit expression of the projection on the boundary eigenspace $\overline\Span(g_k)_{k\in\Z}$.
The proof of Theorems~\ref{theo:ergodicity-compact-trajectories} and~\ref{theo:periodic} relies on the theory of almost periodic functions which  dates back to H. Bohr. 
There is a large literature on the subject and we refer for instance to the book of Corduneanu~\cite{Corduneanu2009} for a comprehensive introduction.
There are several equivalent definitions of almost periodic functions and we will use the following one.
A function $f\in C_b(\R,X)$, {\it i.e.} a bounded continuous function from $\R$ to $X$, is said to be almost periodic if the set $\{f(\cdot+\tau),\ \tau\in\R\}$ is relatively compact in $C_b(\R,X)$.
The set of almost periodic functions is a sub-algebra of $C_b(\R,X)$, and also the closure of the space of periodic functions in $C_b(\R,X)$.
We start with the proof of Theorem~\ref{theo:ergodicity-compact-trajectories} and Theorem~\ref{theo:periodic} in the case when $S$ satisfies the condition~{\it(1)}.
Then we deduce the cases {\it(2)}, {\it(3)} and {\it(4)} from the case {\it(1)}.

\begin{proof}[Proof of Theorems~\ref{theo:ergodicity-compact-trajectories} and~\ref{theo:periodic} in the case (1)] {\sl Step 1.} 
Let $f\in X$.
Since the trajectory $(\widetilde S_t f)_{t\geq0}$ is relatively compact, we infer  from~\cite[Thm.~8]{Haraux1987} (with $U(\tau,t) = S_t$ and thus no periodicity condition on $U$) the existence of an almost periodic eternal solution $g$ of the rescaled semigroup $ \widetilde S$,
{\it i.e.} a function $g:\R\to X$ such that $g(t+\tau)= \widetilde S_\tau g(t)$ for all $t\in\R$ and $\tau\geq0$,
such that
\[\lim_{t\to+\infty}\|\widetilde S_t f-g(t)\|=0.\]
The end of the proof consists in characterizing the function $g$ in the situations of Theorems~\ref{theo:ergodicity-compact-trajectories} and~\ref{theo:periodic}.
For $\lambda\in\R$, we define the Bohr transformation of the almost-periodic function $g$ by
\[
c_\lambda(g)=\lim_{T\to+\infty}\frac1T\int_0^Te^{-i\lambda t}g(t)\,dt,
\]
which is known to exists, see~\cite[Thm.~3.4]{Corduneanu2009}, since $e^{-i\lambda t}g(t)$ is also almost periodic.
Since $e^{-i\lambda t}g(t)$ is besides an eternal solution of the semigroup $ e^{-i\lambda t}\widetilde S_t $ with infinitesimal generator $\LL_\lambda=\LL-\lambda_1-i\lambda$,
we have that
\[
\LL_\lambda\int_0^Te^{-i\lambda t}g(t)\,dt=g(T)-g(0). 
\]
 Dividing by $T$ the above expression, passing to the limit $T\to+\infty$ and using that $\LL_\lambda$ is a closed operator, we get 
\[
\LL_\lambda c_\lambda(g)=0.
\]
In other words, we have established 
\[
\LL c_\lambda(g)=(\lambda_1+i\lambda) c_\lambda(g)
\]
and $\lambda_1+i\lambda$ is an eigenvalue of $\LL$ if $c_\lambda(g)\neq0$.

\smallskip
 {\sl Step 2.}  
 We deduce that if the boundary spectrum is trivial, as in Theorem~\ref{theo:ergodicity-compact-trajectories}, then necessarily $c_\lambda(g)=0$ for all $\lambda\neq0$. By the uniqueness theorem, see for instance~\cite[Thm.~4.7]{Corduneanu2009}, we get that $g$ is constant.
Due to the conservation law $\langle\widetilde S_t f,\phi_1\rangle=\langle f,\phi_1\rangle$ and the simplicity of the eigenvalue~$0$, we get that $g=\langle f,\phi_1\rangle f_1$ and the result of the case {\it(1)} in Theorem~\ref{theo:ergodicity-compact-trajectories} is proved.

\smallskip
 {\sl Step 3.}  In the case of Theorem~\ref{theo:periodic}, the boundary spectrum is not trivial and we know from Theorem~\ref{theo:Sigma+subgroup} that $\Sigma_P^+(\widetilde\LL)$ is an additive subgroup of $i\R$, made of algebraically simple eigenvalues.
 Due to the assumption that $\lambda_1$ is isolated, this subgroup must be discrete and $\Sigma_P^+(\LL)$ is thus given by $\{\lambda_1+i\alpha k ,k\in\Z\}$ for some $\alpha>0$.
As a consequence, any $\lambda$ such that $c_\lambda(g)\neq0$ is necessarily of the form $\lambda=\alpha k$ for some $k\in\Z$.
By the uniqueness theorem, $g$ is then a $\alpha$-periodic function which is given, due to Fejér's theorem, by
\[g(t)=\lim_{n\to\infty}\frac1n\sum_{\ell=0}^n\sum_{k=-\ell}^\ell c_{\alpha k}(g)e^{i\alpha kt}.\]
Consider $(g_k,\psi_k)$ two positive direct and dual eigenvectors of $\LL$ associated to the eigenvalue $i\alpha k$ such that $\langle g_k,\psi_k\rangle=1$.
From the conservation laws $\langle\widetilde S_t f,\psi_k\rangle=\langle f,\psi_k\rangle e^{i\alpha kt}$ and the algebraic simplicity of the eigenvalues $i\alpha k$, we get that $c_{\alpha k}(g)=\langle f,\phi_k\rangle g_k$, and the result is proved.
\end{proof}

\begin{proof}[Proof of Theorems~\ref{theo:ergodicity-compact-trajectories} and~\ref{theo:periodic} in the case (2)]
Since $Y$ is separable, we can find a sequence $(\varphi_n)_{n\geq1}\subset Y$ which satisfies $\|\varphi_n\|=1$ and $\mathrm{span}(\varphi_n)$ is dense in $Y$.
We can then define on $X$ the norm $\|\cdot\|_*$ by setting
\beqn\label{eq:theo:ErgodicityVariante1-1}
\|f\|_*=\sum_{n=1}^\infty 2^{-n}|\langle f,\varphi_n\rangle|.
\eeqn
On bounded subsets of $X$, the topology of this norm is the same as the weak-* topology, or more explicitly 
$$
f_n \wto f \  *\sigma(Y',Y)\ \Leftrightarrow  \ (\sup \|f_n\| < \infty \ \ \hbox{ and } \  \ \|f_n - f \|_* \to 0).
$$ 
Since by assumption the trajectory $(\widetilde S_t f)$ is bounded, it is weakly-* relatively compact, and so relatively compact in $(X,\|\cdot\|_*)$.
  It is also clear that the semigroup $S$ is continuous for the weak norm $\|\cdot\|_*$. 
The normed space $(X,\|\cdot\|_*)$ is not a Banach space, but the proof of Theorem~\ref{theo:periodic} actually only requires, for applying~\cite[Thm.~8]{Haraux1987}, 
  that the closed balls of $X$ are complete metric spaces, which is the case for the distance induced by $\|\cdot\|_*$.  
Applying the case {\it(1)} of Theorems~\ref{theo:ergodicity-compact-trajectories} and~\ref{theo:periodic} then yields the result.
\end{proof}
 
\begin{proof}[Proof of Theorems~\ref{theo:ergodicity-compact-trajectories} and~\ref{theo:periodic} in the case (3)]
We consider $f \in X$ and, repeating the proof of Step 2 in Theorem~\ref{theo:MeanErgodicityVariante1}, we get
that $(S_t f)_{t \ge 0}$ belongs to a weak compact set $\GG$ of $L^1_{\phi_1}$.
We define the norm $\|\cdot\|_*$  by \eqref{eq:theo:ErgodicityVariante1-1} for a sequence $(\varphi_n)_{n\geq1}\subset C_c(E)$ which satisfies $\|\varphi_n\|_{L^\infty}=1$ and $\mathrm{span}(\varphi_n)$ is dense in $C_0(E)$. This norm induces a metric on $\GG$ which is topologically equivalent to the weak convergence on $L^1_{\phi_1}$. 
The trajectory $(\widetilde S_t f)$ is then relatively compact in $(\GG,\|\cdot\|_*)$ and the semigroup $S$ is continuous for the weak norm $\|\cdot\|_*$.
We conclude as in the proof of the case~{\it (2)}.
\end{proof}

\begin{proof}[Proof of Theorems~\ref{theo:ergodicity-compact-trajectories} and~\ref{theo:periodic} in the case (4)]
From the step 3 of the proof of Theorem~\ref{theo:MeanErgodicityVariante1}, we know that for any $f\in\XX^k$ the trajectory $(\widetilde S_t f)$ is compact for the strong topology of $L^1_{\phi_1}$.
We may then conclude similarly as in the case~{\it(1)}, using that $\XX^k$ is dense in $X$ for the norm of $L^1_{\phi_1}$.
\end{proof}

%

\subsection{A word about spectral analysis argument}
\label{subsec:stabilityKR-Spectral Analysis}

 
The aim of this section is to recall some more or less classical results which make possible to slightly improve the conclusions of the results presented in the previous section by additionally assuming some spectral gap at the level of the operator or the semigroup. More specifically, we are interested in some accurate versions of a  partial, but principal  {\it spectral mapping theorem} which asserts that 
\beqn\label{eq:PartialSpectralMapping}
\Sigma(e^{t\LL}) \cap B^c(0,e^{\kappa t})  = e^{t\Sigma(\LL) \cap \overline\Delta_\kappa},
\qquad\forall\, t\ge0, 
\eeqn
for some $\kappa < \lambda_1$,  and even more precisely,  we wish to establish the  following geometric (or exponential) asymptotic stability

\begin{enumerate}[label={\bf(E3$_1$)},itemindent=14mm,leftmargin=0mm,itemsep=1mm]
\item\label{E31} there exist some constants $\kappa < \lambda_1$ and $C \ge 1$ such that   for any $f \in X$, there holds
\beqn\label{eq:intro-expoCvgceBIS}
\| \widetilde S(t) f - \langle f , \phi_1 \rangle f_1 \| \le \Theta(t) \| f - \langle f , \phi_1 \rangle f_1 \|, \quad \forall \, t \ge 0, \ \forall \, f \in X, 
\eeqn
with the decay rate function $\Theta(t) :=  C \, e^{(\kappa-\lambda_1) t}$. 
\end{enumerate}
 
\smallskip

 In order to discuss the several results we present, we recall the splitting framework 
\beqn\label{eq:PartialSpectralMapping-2}
S = V + W * S, \qquad \| V(t) \|_{\BBB(X)} + \| W(t) \|_{\BBB(X)}  \lesssim e^{\kappa t}, 
\eeqn
for the same $\kappa \in \R$ as above. 

We start by recalling the spectral mapping theorem for the point spectrum,
\Blue
which states that
\begin{equation}\label{eq:pointspectralmapping}
\Sigma_P(e^{t\LL})\setminus\{0\}=e^{t\Sigma_P(\LL)}.
\end{equation}
We actually give a slightly more precise version of this result.
For $\xi\in\C$ and $t>0$, we define $\Lambda_t(\xi):=\{\lambda\in\C,\ e^{\lambda t}=\xi\}$.
Clearly if $\lambda_0\in\Lambda_t(\xi)$, then $\Lambda_t(\xi)=\{\lambda_0\}+\frac{2i\pi}{t}\Z$.

\begin{lem}[Spectral mapping theorem for the point spectrum]\label{lem:pointspectralmapping}
For a semigroup $(S_t)_{t\geq0}$ with infinitesimal generator $\LL$ we have, for any $\kappa\in\R$ and $t>0$,
\begin{equation}\label{eq:kappapointspectralmapping}
\Sigma_P(S_t) \cap \{ \xi \in \C; \, |\xi| \geq e^{\kappa t} \} =e^{t (\Sigma_P(\LL) \cap \{ \lambda\in\C; \; \Re e \lambda \geq \kappa \})}.
\end{equation}
More precisely

-  If $f\in X$ is an eigenvector of $\LL$ with corresponding eigenvalue $\lambda\in\C$, then $f$ is an eigenvector of $S_t$ with corresponding eigenvalue $e^{\lambda t}$,

- If $f\in X$ is an eigenvector of $S_t$ with corresponding eigenvalue $\xi\in\C$, then there exists $\lambda\in\Lambda_t(\xi)$ such that $g:=\int_0^t e^{-\lambda s}S_sf\,ds$ is an eigenvector of $\LL$ with corresponding eigenvalue $\lambda$.
\end{lem}

Clearly, taking the union over all $\kappa\in\R$ of~\eqref{eq:kappapointspectralmapping} yields~\eqref{eq:pointspectralmapping}.

\begin{proof}[Proof of Lemma~\ref{lem:pointspectralmapping}]
Fix $\kappa\in\R$, $t>0$, and define
\[\Sigma_S := \Sigma_P(S_t) \cap \{ \xi \in \C; \, |\xi| \geq e^{\kappa t} \}, \quad \Sigma_\LL :=  \Sigma_P(\LL) \cap \{ \lambda\in\C; \; \Re e \lambda \geq  \kappa \}.\]
 One inclusion is clear: if $\lambda \in \Sigma_\LL$ and $\LL f = \lambda f$, $f \not=0$, then $\xi := e^{\lambda t}$ is such that $S_t f  = \xi f$ and thus 
$\xi \in \Sigma_P(S_t)$, $|\xi|=e^{\Re e\lambda\,t} \geq e^{\kappa t}$. We have proved $e^{t\Sigma_\LL} \subset \Sigma_S$.
In the other way around, let us take $\xi \in \Sigma_S$ and $f \in X_\C\backslash \{ 0 \}$ such that $S_t f  = \xi f$.
Then let us pick up $\lambda_0\in\C$ such that $\xi=e^{\lambda_0 t}$ and $\phi\in X'$ such that $\langle\phi,f\rangle\neq0$.
For any $k\in\Z$, we have $\xi=e^{\lambda_0 t+2ik\pi}$ and so
\[0=e^{-(\lambda_0+\frac{2ik\pi}{t})t}S_tf-f=\Big(\LL-\lambda_0-\frac{2ik\pi}{t}\Big)\int_0^t e^{-(\lambda_0+\frac{2ik\pi}{t})s}S_sf\,ds.\]
If $g_k:=\int_0^t e^{-(\lambda_0+\frac{2ik\pi}{t})s}S_sf\,ds$ is non-zero for some $k\in\Z$, we deduce that $g_k$ is an eigenvector of $\LL$ with corresponding eigenvalue $\lambda := \lambda_0+\frac{2ik\pi}{t}$.
Since besides $\Re e\lambda = \frac1t \log|\xi| \geq \kappa$, we get that $\lambda \in \Sigma_\LL$ and so $\xi \in e^{t\Sigma_\LL} $.
It thus only remains to prove that there exists $k\in\Z$ such that $g_k\neq0$.
Assume by contradiction that $g_k=0$ for all $k\in\Z$.
This means that the continuous and 
periodic complex-valued function $s\mapsto e^{-\lambda_0 s}\langle\phi,S_sf\rangle$ has all its Fourier coefficients equal to zero, which is not possible since this function is not equally zero (its value at $s=0$ is $\langle\phi,f\rangle\neq0$).
 \end{proof}
 
 \Black


We next present a very classical result about the exponential stability of $f_1$ which is based on the quasi-compact semigroup framework of Voigt \cite{MR595321} 
(see also \cite[B-IV-2]{MR839450}
and  \cite[Sec.~V.3]{MR1721989})
and which is a more accurate version of Lemma~\ref{lem:H3abstract-StrongCbis} and Theorem~\ref{theo:Sigma+subgroupBIS}. 

\begin{theo}\label{theo:NagelWebb} 
Let $(S_t)_{t\geq0}$ be a positive irreducible semigroup on a Banach lattice $X$ satisfying the hypotheses of Lemma~\ref{lem:H3abstract-StrongCbis} and Theorem~\ref{theo:Sigma+subgroupBIS},
in particular \ref{H2} holds for a constant $\kappa_0 \in \R$ and there exists $T> 0$  such that the splitting 
\beqn\label{eq:KRexistTERS=VK}
S_T = V_T + K_T,  
\eeqn
holds with $\| V_T \|_{\BBB(X)} \le e^{\kappa T}$, $\kappa < \kappa_0$, and $K_T \in \KK(X)$.  
Then there exists a unique solution $(\lambda_1,f_1,\phi_1)$ to the eigentriplet and
the exponential stability \ref{E31}  holds (without constructive estimate).
 \end{theo}

%
%

%
%
%

\begin{rem} In the splitting framework \eqref{eq:PartialSpectralMapping-2} the critical hypothesis  $K_T \in \KK(X)$ may be obtained by assuming that 
$$
 \|  W (t) \|_{\BBB(X,\XX_1)}\lesssim  e^{\kappa t}, \quad \forall \, t \ge 0, \quad \XX_1 \subset X \hbox{ compact}.
$$
In fact, in many applications, we are also able to establish $\XX_1 \subset D(\LL^\beta)$, for some $\beta >0$, without too much more work.
\end{rem}

Theorem~\ref{theo:NagelWebb} is in fact nothing but \cite[Thm.~14.18]{MR3616245} (see also \cite[Sec.~2]{MR902796}, 
\cite[Thm.~V.3.7]{MR1721989}  or  \cite[C-IV,   Thm.~2.1 \& Rk.~2.2]{MR839450}). We give however a short proof of Theorem~\ref{theo:NagelWebb} since it is simpler and more direct than the ones we usual find in the literature
and in particular does not refer to subtle results about the spectrum and its essential part.

\begin{proof}[Proof of Theorem~\ref{theo:NagelWebb}] 
\noindent{\it First step.}
From  Lemma~\ref{lem:H3abstract-StrongCbis}, we already know that  \ref{H1}, \ref{H2} and \ref{H3} hold. Together with the irreducibility which is nothing but \ref{H4} from Lemma~\ref{lem:Irred-S>0impliesR>0}, we may apply Theorem~\ref{theo:KRgeometry1} 
and conclude to the existence, uniqueness, strict positivity {\Blue and algebraic simplicity} result about the eigentriplet solution $(\lambda_1,f_1,\phi_1)$, {\Blue with $\lambda_1 > \kappa$.}

\medskip

{\Blue
\noindent{\it Second step.}
For any $\kappa'>\kappa$, the set $\Sigma(S_T)\cap\{z\in\C, |z|\geq e^{\kappa'T}\}$ is made of a finite number of isolated eigenvalues with finite algebraic multiplicity.
Let us briefly present the argument and we refer to \cite{MR236741,MR595321} for more details. 
From \eqref{eq:KRexistTERS=VK}, we readily have the identity 
\[
(z-S_T)^{-1}=(z-V_T)^{-1}{(I-K_T(z-V_T)^{-1})}^{-1},
\]
as soon as the RHS term is well defined. On the one hand, the function $z \mapsto (z-V_T)^{-1}$ is holomorphic on $\{z \in \C; \, |z| > e^{\kappa t} \}$ with values in $\BBB(X)$ because of the boundedness assumption made on $V_T$, and thus  $z \mapsto K_T (z-V_T)^{-1}$ is holomorphic   with values in $\KK(X)$. Using \cite[Corollary~II]{MR236741} (see also the proof of \cite[Theorem~1.1]{MR595321}) we obtain that $z \mapsto (I-K_T(z-V_T)^{-1})^{-1}$ and thus  $(z-S_T)^{-1}=(z-V_T)^{-1}{(I-K_T(z-V_T)^{-1})}^{-1}$ are degenerate-meromorphic functions. That means that $(z-S_T)^{-1}$ is a   meromorphic function such that at any of its (isolated) poles, the coefficients of the (finite number of) singular terms in the corresponding Laurent series are finite rank operators. 
As a consequence, the Dunford integral on a small radius circle around a pole $\xi$ (which is equivalently a spectral value of $S_T$) defines a projection $\Pi_\xi$ with finite rank: in other words, $\Pi_\xi$ is the projection on the algebraic eigenspace associated to $\xi$. We have established that $\Sigma(S_T)\cap\{z\in\C, |z|> e^{\kappa T}\} \subset \Sigma_d(S_T)$, the set of  isolated eigenvalues with finite algebraic multiplicity. 
Observing that  $\Sigma(S_T)\cap\{z\in\C, |z|\geq e^{\kappa'T}\}$ is a compact subset of $\Sigma(S_T)\cap\{z\in\C, |z|> e^{\kappa T}\}$, since $S_T$ is a bounded operator, we readily conclude. 


\medskip
\noindent{\it Third step.}
Since $S_Tf_1=e^{\lambda_1 T}f_1$, we have that $e^{\lambda_1T}\in\Sigma(S_T)$.
We now prove the existence of a spectral gap, namely the existence of $\epsilon>0$ such that $\Sigma(S_T)\cap\{z\in\C, |z|\geq e^{(\lambda_1-\epsilon)T}\}=\{e^{\lambda_1 T}\}$.
Let us take $\kappa'\in(\kappa,\lambda_1)$.
From the previous step, there must be some $\epsilon>0$ such that 
$$
\Sigma(S_T)\cap\{z\in\C, |z|\geq e^{(\lambda_1-\epsilon)T}\}=\Sigma(S_T)\cap\{z\in\C, |z|\geq e^{\lambda_1T}\},
$$
and thus 
\begin{align*}
\Sigma(S_T)\cap\{z\in\C, |z|\geq e^{(\lambda_1-\epsilon)T}\} 
&=\Sigma_P(S_T)\cap\{z\in\C, |z|\geq e^{\lambda_1T}\}\\
&=e^{T (\Sigma_P(\LL) \cap \{ \lambda\in\C, \; \Re e \lambda \geq \lambda_1 \})}=e^{T\Sigma^+_P(\LL)},
\end{align*}
where we have used the previous step again for the first equality, Lemma~\ref{lem:pointspectralmapping} for the second equality, and the definition of the boundary point spectrum $\Sigma^+_P(\LL)$ in the last one.
We next claim that $\Sigma^+_P(\LL)$ is finite. Otherwise, it would exist a infinite family $\{(\lambda_k, f_k) \}$ of couples of  eigenvalue-eingenvector of $\LL$ with $\lambda_k \in\Sigma^+_P(\LL)$ and thanks to Lemma~\ref{lem:pointspectralmapping} a corresponding  infinite family $\{(e^{T\lambda_k}, f_k) \}$ of couples of  eigenvalue-eingenvector of $S_T$ with $|e^{T \lambda_k}| \ge e^{T\lambda_1  }$, what is in contradiction with the conclusion of the second step.
Since  the translated  boundary point spectrum  $\Sigma_{P}^{+} (\LL) -\lambda_1$  is a group from Theorem~\ref{theo:Sigma+subgroup}, we infer that the only possibility is $\Sigma^+_P(\LL)=\{\lambda_1\}$,
which finally gives $\Sigma(S_T)\cap\{z\in\C, |z|\geq e^{(\lambda_1-\epsilon)T}\}=\{e^{\lambda_1 T}\}$.

\medskip
\noindent{\it Fourth step.}  We prove that $e^{\lambda_1T}$ is (geometrically) simple with eigenspace spanned by $f_1$.
Let $f_2$ be an eigenvector of $S_T$ with corresponding eigenvalue $e^{\lambda_1T}$.
Since $\langle\phi_1,f_1\rangle\neq0$, we can find $\alpha\in\C$ such that $f_3:=f_2-\alpha f_1$ verifies $\big\langle\phi_1,\int_0^T e^{-\lambda_1s}S_sf_3\,ds\big\rangle=T(\langle\phi_1,f_2\rangle-\alpha\langle\phi_1,f_1\rangle)=0$.
We now prove that $f_3=0$, which will give the conclusion.
Assume by contradiction that $f_3\neq0$.
Then it is an eigenvector of $S_T$ with corresponding eigenvalue $e^{\lambda_1T}$, as a linear combination of $f_1$ and $f_2$.
By virtue of Lemma~\ref{lem:pointspectralmapping}, there exists $\lambda\in\Lambda_T(e^{\lambda_1T})=\{\lambda_1\}+\frac{2i\pi}{T}\Z$ such that $\int_0^t e^{-\lambda s}S_sf_3\,ds$ is an eigenvector of $\LL$ with corresponding eigenvalue $\lambda$.
We thus have
\[\lambda\in\Sigma_P(\LL)\cap\Lambda_T(e^{\lambda_1T})\subset\Sigma_P(\LL)\cap(\{\lambda_1\}+i\R)=\Sigma_P^+(\LL)=\{\lambda_1\},\]
where the last equality comes from the third step of the proof,
and $g:=\int_0^t e^{-\lambda_1 s}S_sf_3\,ds$ is an eigenvector of $\LL$ with corresponding eigenvalue $\lambda_1$.
Since we know from the first step of the proof that $\lambda_1$ is simple, we must have $g\in\Span(f_1)\setminus\{0\}$, which is not possible since $\langle\phi_1,g\rangle=0$ by construction of $f_3$.
We have obtained a contradiction, from which we get that $f_3=0$, $f_2=\alpha f_1$.

\medskip

\noindent{\it Fifth step.}
We conclude to the exponential stability~\ref{E31}.
We have proved in the two previous steps that}
$\Sigma(S_T)\cap\{z\in\C,\ |z|\geq e^{(\lambda_1-\eps) T}\}=\{e^{\lambda_1 T}\}$ and $e^{\lambda_1 T}$ is simple with eigenspace spanned by $f_1$.
The restriction $S_T^\perp$ of $S_T$ to the invariant subspace $ X_\perp:=\{f\in X,\ \langle\phi_1,f\rangle=0\}$, which does not contain $f_1$, thus has a spectral radius smaller than $e^{(\lambda_1-\eps) T}$.
The spectral radius formula  (see \cite[Thm.~10.13]{MR1157815}  for instance) then ensures that
\[\lim_{n\to\infty}\|S_{nT}^\perp\|^{1/n}=r(S_T^\perp)\leq e^{(\lambda_1-\eps) T}.\]
This guarantees, for any $\eta\in(0,\eps)$, the existence of a constant $C_\eta>0$ such that for all $f\in X_\perp$ and all $t\geq0$
\[\|e^{-\lambda_1 t}S_tf\|\leq C_\eta e^{-\eta t}\|f\|,\]
and the proof is complete.
\end{proof}

Let us now present a variant of another classical result known as the Gearhart-Pr\"{u}ss Theorem in  \cite{MR461206,MR743749}, see also the contributions of  Herbst \cite{MR715559} and Greiner 
 \cite[A-III.7]{MR839450} as well as the more constructive proof \cite[Thm.~V.1.11]{MR1721989} and recently \cite{MR4271979} 
 based on techniques developed in or related to \cite{MR1347423,MR1823064}. 

 \begin{theo}\label{theo:GearhartPruss} 
 Consider a positive semigroup $S$ on a Banach lattice $X$ such that its generator $\LL$ satisfies  the conclusions \ref{C2} about the existence, positivity and uniqueness 
of the first eigentriplet $(\lambda_1,f_1,\phi_1)$. 
We assume furthermore that $X$ is an Hilbert space  and that there exist $\kappa < \lambda_1$ and $R > 0$ such that 

\smallskip
\quad (i) $\sup_{z \in \Delta_\kappa \backslash B_R} \| \RR_\LL(z) \|_{\BBB(X)} < \infty$; 

\smallskip
\quad (ii) $\Sigma(\LL) \cap \Delta_\kappa \subset \Sigma_d(\LL)  \cap B_R$. 

\smallskip
Then the exponential stability \ref{E31}  holds (without constructive estimate).
%
%
%
 \end{theo}

\begin{proof}[Proof of  Theorem~\ref{theo:GearhartPruss}.]
The spectral information \ref{C2} and (ii)  together imply \ref{C3} (because of Theorem~\ref{theo:Sigma+subgroup}) and that
 there exists  $\kappa^* \in (\kappa,\lambda_1)$, such that  $\Sigma(\LL) \cap \Delta_{\kappa^*} = \{ \lambda_1 \}$.
The operator $  \LL  $ on $X_0 := (\hbox{\rm vect}\{f_1\})^\perp$ thus satisfies  $\sup_{z \in \Delta_{\kappa^*}} \| \RR_\LL(z) \|_{\BBB(X_0)} < \infty$, and we conclude
thanks to \cite[Thm.~V.1.11]{MR1721989}. The lack of constructively here only comes from the fact that our assumptions do not provide any information on the spectral gap $\lambda_1 - \kappa > 0$.
\end{proof}

\begin{rem}\label{rem:theo:GearhartPruss} Except of the Hilbert space framework, the assumptions made in Theorem~\ref{theo:GearhartPruss} are slightly weaker 
 than those of Theorem~\ref{theo:NagelWebb}, and are indeed established during the proof of Theorem~\ref{theo:NagelWebb}: such an  information at the level of the resolvent is a bit easier to establish than a similar estimate at the level of the semigroup.
In the splitting framework \eqref{eq:PartialSpectralMapping-2} and its resolvent counterpart \eqref{eq:exist1-defVVWW}, we typically only have to show 
\beqn\label{eq:rem:theo:GearhartPruss}
 \sup_{\kappa \le \Re e z \le \kappa_1 } \| \VV(z)  \|_{\BBB(X)} < \infty,  \quad
\lim_{r \to \infty} \sup_{\kappa \le \Re e z \le \kappa_1, \, |\Im m z| \ge r}   \| \WW(z) \|_{\BBB(X)}  = 0,
\eeqn
for some $\kappa < \lambda_1$, and $\WW(z) \in \KK(X)$ for any $z \in \Delta_\kappa$. 
That last claim is classical (see for instance \cite{MR3779780}) and we only briefly sketch the proof. On the one hand, from the first and the last estimates, we deduce 
that $\Sigma(\LL) \cap \Delta_{\kappa} \subset \Sigma_d(\LL)$ thanks to Theorem~\ref{theo:Voigt}. 
As in the proof of Theorem~\ref{theo:Voigt} and with the usual notations, we also have 
$$
(I-\WW(z)) \RR_\LL(z) = \VV(z), \quad \forall \, z \in \Delta_\kappa,
$$
where $I-\WW(z)$ is invertible and $\| (I-\WW(z))^{-1} \|_{\BB(X)} \le 2$ for any $z \in \C$ such that $\kappa \le \Re e z \le \kappa_1$, $ |\Im m z| \ge R$ and $R$ is large enough. 
We immediately deduce  that the condition (i) in Theorem~\ref{theo:GearhartPruss} holds. 
 
\end{rem}

We end this section by a more recent result which is similar to the Gearhart-Pr\"{u}ss Theorem but is not restricted to an Hilbert space.

 \begin{theo}\label{theo:MiscScher}   
  Consider a positive semigroup $S$ on a Banach lattice $X$ such that its generator $\LL$ satisfies  the conclusions \ref{C2} about the existence, positivity and uniqueness 
of the first eigentriplet $(\lambda_1,f_1,\phi_1)$. We further assume that $\LL = \AA + \BB$ with the associated operators $V$  and 
$W$ defined by \eqref{eq:V=&W=} satisfy \eqref{eq:PartialSpectralMapping-2} for some $\kappa < \lambda_1$ and that the resolvent counterpart 
$\WW$ defined by \eqref{eq:exist1-defVVWW} satisfies \eqref{eq:rem:theo:GearhartPruss} and more precisely
%
$$
 \sup_{\kappa \le \Re e z \le \kappa_1 } \| \langle z \rangle^{\alpha}\WW(z)  \|_{\BBB(X)} < \infty,   
 $$
with  $\alpha > 1$.
Then the exponential stability \ref{E31}  holds (without constructive estimate).
  \end{theo}

The proof of Theorem~\ref{theo:MiscScher}  is a mere adaptation  of \cite[Thm.~3.1]{MR3489637} (see also \cite{MischErratum}) and it is thus skipped. 
The needed estimates are a bit stronger than those of Remark~\ref{rem:theo:GearhartPruss}, but in the applications, they are not really more demanding. They also hold at the level of the resolvent 
instead of what is assumed in the statement of  Theorem~\ref{theo:NagelWebb}.

%
%
%
%
%
%
%
%
%
%
%
%
%
%
%
%

\smallskip
We conclude by emphasizing   again on the fact that all the above results are not constructive. 
We propose in the next part an alternative approach which is constructive. 

\medskip
 
{\Blue
Let us emphasize first  that \ref{C2} (which can be established using the material of Sections~\ref{sec:ExistenceKR}, \ref{sec:DynamicalExistenceKR} and ~\ref{sec:Irreducibility})  together with Theorem~\ref{theo:KRgeometry2}, Theorem~\ref{theo:lambda1largestBIS} or Theorem~\ref{theo:lambda1largestTER}-(1) imply the triviality of the principal point spectrum $\Sigma^1_P(\LL) = \{\lambda_1\}$, what is a weak version of the answer \ref{S32} about the geometry of the spectrum problem. 
With the additional property $\Sigma^1_P( \LL) = \Sigma^+_P( \LL)$  (which can be established directly using the material of Section~\ref{sec:ExistenceKR} or, for a positive semigroup, using  the material of Sections~\ref{sec:DynamicalExistenceKR} and ~\ref{sec:Irreducibility}, see in particular the discussion at the end of Section~\ref{sec:Irreducibility}), we conclude to  \ref{S32}. 
When furthermore $\LL$ is also the generator of a positive semigroup, Theorem~\ref{theo:ErgodicityByLiapunov1}, 
 Corollary~\ref{cor:Ergodicity2} or Theorem~\ref{theo:ergodicity-compact-trajectories} implies the ergodicity property \ref{E2} in several topologies. Altogether, in that way, we thus complete the proof of  Theorem~\ref{theo:main-intro}-(3). 
Alternatively, we may use Theorem~\ref{theo:NagelWebb} or \ref{C2} together with Theorem~\ref{theo:GearhartPruss} or Theorem~\ref{theo:MiscScher} in order to obtain the exponential stability \ref{E31}, 
 and thus also the spectral gap property \ref{S33}.
 In particular, these arguments provide a proof of  Theorem~\ref{theo:main-intro}-(4). 

%
%
%
 }

%
%

\Black

\medskip

\bigskip

\bigskip
\section{Quantitative stability} 
\label{sec:QuantitativeStabilityKR}


 {\Blue

 \medskip

In this section we establish some quantitative stability results in the spirit of the Doblin, Harris, Meyn-Tweedie theory for Markov semigroup, and in particular  the  geometric part  \ref{S32} and the evolution part  \ref{E31intro} of the Krein-Rutman theorem in a quantitative way. 
We do not use the previous positivity and other structure conditions {\bf (Hi)}, {\bf (Xi)}, {\bf (HSi)}, but we rather start with the condition \ref{C1} and introduce the additional necessary conditions, which are somehow stronger: they are mainly formulated in term of the semigroup and are always quantitative.}
  
\subsection{About quantified positivity conditions}
\label{sec:quantified-positivity}

 We briefly discuss some positivity conditions related to the strong maximum principle and barriers techniques. 
The issue is about how to quantify the strong maximum principle
$$
f \in X_+ \backslash \{ 0 \}, \ (\kappa_1 - \LL) f \ge 0 \ \hbox{ imply } \ f > 0    \hbox{ or }  f \gg 0
$$
or the related strong positivity of the associated semigroup. A possible way can be achieved with the help of a barrier functions family $\GG \subset X_+$
and a second weaker (semi)norm $[\cdot]$ used for normalization. Let us then introduce the two conditions  
\beqn\label{eq:QuantitativeStab-SMP}
\forall \, R > 0, \ \exists g_i \in \GG, \ \forall \, f \in X_+, \ [f] = 1, \ \| f \| \le R, 
\eeqn
we have 

(i) $S_T f \ge g_1$ (for some $T >0$)

or

(ii) $f \ge g_2$ if $(\kappa_1 - \LL) f \ge 0$.

\smallskip
Point (ii) is a quantified version of the strong maximum principle when $\GG \subset X_{++}$ and it is always a consequence of the positivity condition (i). 
Assume indeed that (i) holds (for some $T > 0$) and that $f$ satisfies the requirements \eqref{eq:QuantitativeStab-SMP} and  $(\kappa_1 - \LL) f \ge 0$. We then write 
$$
\frac{d}{dt}  ( e^{(\LL-\kappa_1)t} f ) = e^{(\LL-\kappa_1)t} (\LL - \kappa_1) f \le 0,
$$
so that 
$$
f \ge e^{(\LL-\kappa_1)T} f = e^{-\kappa_1 T} S_T f \ge e^{-\kappa_1 T} g_2=: g_1, 
$$
with $g_2$ given by condition (i). The reciprocal implication is not clear, see however Lemma~\ref{lem:Irred-S>0impliesR>0}-{\bf  (3)}.

\medskip
Let us now make a list of possible quantified positivity conditions of Doblin-Harris type for a linear (and continuous) operator  $A: X \to X$:

(P1$'$)  $\exists\,g_0 \in X_+\backslash\{0\} $, $\exists\, \psi_0 \in X_+\backslash\{0\} $, $\forall\, f \in X_+$, 
$A f  \ge  g_0 \langle f,\psi_0 \rangle$; 

(P2$'$) $\exists\,g_0 \in X_+\backslash\{0\} $, $\exists\, \psi_0 \in X'_{++}$, $\forall\, f \in X_+$,  $A f  \ge  g_0 \langle f,\psi_0 \rangle$; 

(P3$'$) $\exists\,g_0 \in X_{++} $, $\exists\, \psi_0 \in X'_+\backslash\{0\}$, $\forall\, f \in X_+$, $A f  \ge  g_0 \langle f,\psi_0 \rangle$; 

(P4$'$) $\exists\,g_0 \in X_{++} $, $\exists\, \psi_0 \in X'_{++}$, $\forall\, f \in X_+$,  $A f  \ge  g_0 \langle f,\psi_0 \rangle$. 

\smallskip
We summarize some elementary relations between these conditions and those listed in Section~\ref{subsec:irreducibility}. 

\begin{lem}\label{lem:positivity2} We have 
 (P2$'$) $\Rightarrow$ (P2) $\Rightarrow$ (P1), (P3$'$) $\Rightarrow$ (P3) $\Rightarrow$ (P1), (P4$'$) $\Rightarrow$ ((P4), (P3$'$), (P2$'$)) as well as
(P4) $\Rightarrow$ ((P3), (P2)).

\smallskip
We also have: 
$A$ satisfies (P2$'$) iff $A^*$ satisfies (P3$'$);
$A$ satisfies (P3$'$) iff $A^*$ satisfies (P2$'$);
$A$ satisfies (P4$'$) iff $A^*$ satisfies (P4$'$).

\smallskip
We finally have:
$A$ satisfies (P2$'$) implies $\exists\,g_0 \in X_+\backslash\{0\} $, $\exists\, \kappa > 0$, $A g_0  \ge  \kappa g_0$

\end{lem}

\begin{proof}[Proof of Lemma~\ref{lem:positivity2}.]
We assume $A f  \ge  g_0 \langle f,\psi_0 \rangle$ for any $f \in X_+$ and some $g_0 \in X_+$, $\psi_0 \in X'_+$.
 For any $\phi \in X'\backslash\{0\}$ and $f \in X_+$, we have 
$$
\langle A^*\phi,f\rangle = \langle \phi, Af\rangle \ge   \langle \phi,  g_0 \langle f,\psi_0 \rangle \rangle,
$$
which implies $A^*\phi \ge  \psi_0 \langle \phi,  g_0  \rangle$. We thus deduce that $A$ satisfies (P2$'$) (resp. (P3$'$), (P4)) implies that $A^*$ satisfies  (P3$'$) (resp. (P2$'$), (P4)).
The other implications can be established in a similar or even simpler way.
\end{proof}

We conclude this introductory section by emphasizing on the fact (as already mentioned above) that  $S_\LL$ satisfies (P$i'$) implies $\RR_\LL(\lambda)$ satisfies (P$i'$) for any $\lambda \ge \lambda_1$ and $i=1, \dots, 4$.

\medskip
\subsection{Asymptotic stability under Doblin condition}
\label{sec:compactcase}

We start with a simple situation. 
We assume the Doblin condition, namely 
\beqn\label{eq:hyp-Doblin}
\exists\, T > 0, \ \exists \, \psi_0 \gg 0, \ \exists \, g_0 > 0, \quad \forall \, f \ge 0, \ S_T f \ge g_0 \langle \psi_0, f \rangle, 
\eeqn
together with the companion positivity condition 
\beqn\label{eq:hyp-Doblin-BdBis}
\exists \, r_0 > 0, \quad \langle \phi_1, g_0 \rangle \ge r_0,
\eeqn
as well as the strong additional boundedness assumption
\beqn\label{eq:hyp-Doblin-Bd} 
\exists \, R_0  > 0, \quad \phi_1 \le R_0 \psi_0. 
\eeqn
When $\psi_0 := 1 \in X' \subset L^\infty$, the condition in \eqref{eq:hyp-Doblin-Bd} is automatically satisfied with $R_0 := \| \phi_1 \| =1$.   
Let us first emphasize that \eqref{eq:hyp-Doblin-BdBis} is a natural condition when  $S^*_\LL$  enjoys a splitting structure similar to \eqref{eq:ergodicity-splitting}. More precisely, when 
 $$
  \|  \widetilde S^* (t) \phi \| \le \Theta(t)   \|  \phi \| + \int_0^t \Theta(t-s)  [  \widetilde S^*(s) \phi ]_{g_0} \, ds, 
$$
with $\Theta \in L^1(\R_+) \cap C_0(\R_+)$, we deduce that 
$$
  1 = \| \phi_1 \| = \|  \widetilde S^*(t) \phi_1 \| \le \Theta(t) + \int_0^t \Theta(t-s)  [ \phi_1]_{g_0} \, ds, \quad \forall \, t > 0.
$$
Passing to the limit $t \to \infty$, we get  \eqref{eq:hyp-Doblin-BdBis}  with $r_0 := \| \Theta \|_{L^1}^{-1}$. Also \eqref{eq:hyp-Doblin-Bd} can be deduced from a  splitting structure condition on the dual problem. 
More precisely, we assume that $D(\LL^\infty) \subset \Lloc^1$ and  the splitting property
$\LL = \AA+\BB$ with $\AA \in \BBB(X)$, $\RR_\BB(\lambda) \in \BBB(X) \cap \BBB(X_+)$ for any $\lambda \ge \kappa$, with $\kappa < \kappa_0 \le \lambda_1$, and the additional regularity condition
\beqn\label{eq:hyp-Doblin-WL1Linfty}
(\RR_{\BB^*} (\lambda) \AA^*)^N : L^1_{g_0} \to L^\infty_{\psi_0^{-1}}, \quad \forall \, \lambda > \kappa.
\eeqn
Since the dual eigenvector $\phi_1$ satisfies
$$
(\lambda_1-\BB^*) \phi_1 = \AA^* \phi_1, \quad \lambda_1 > \kappa,
$$
and then $\phi_1 = (\RR_{\BB^*} (\lambda_1) \AA^*)^N \phi_1$, we may use estimate \eqref{eq:hyp-Doblin-WL1Linfty} and we get that \eqref{eq:hyp-Doblin-BdBis}-\eqref{eq:hyp-Doblin-Bd} holds with the normalization condition $r_0 := 1$
and $R_0 := \| (\RR_{\BB^*} (\lambda) \AA^*)^N \|_{\BBB(L^1_{g_0} ,  L^\infty_{\psi_0^{-1}})}$. 
  
 \medskip
 We are then able to formulate a first quantified stability result.

\begin{theo}\label{theo:Doblin} 
Consider a {\Cyan positive} semigroup $S$ on a Banach lattice $X$ such that its generator $\LL$ enjoys the conclusion \ref{C1} on the existence 
of {\Blue the two eigenpairs $(\lambda_1,f_1)$ and $(\lambda_1^*,\phi_1)$.}
We assume furthermore  the Doblin condition \eqref{eq:hyp-Doblin}-\eqref{eq:hyp-Doblin-Bd}-\eqref{eq:hyp-Doblin-BdBis}. 
Then {\Blue $\lambda_1=\lambda_1^*$ and} 
the exponential stability \ref{E31intro} in the norm $[\cdot]_{\psi_0}$ holds true, with constructive constants.  
\end{theo}

The proof closely  follows the usual contraction argument in the Doblin result, see for instance \cite{MeynTweedie}, \cite[Thm.~11]{MR3893207} or \cite[Thm.~2.1]{MR4534707}.
We do not explicitly assume the irreducibility of the semigroup, but the Doblin condition \eqref{eq:hyp-Doblin}-\eqref{eq:hyp-Doblin-Bd}-\eqref{eq:hyp-Doblin-BdBis} is in many aspects a strong positivity condition. In particular, our result implies the uniqueness of the first eigentriplet $(\lambda_1,f_1,\phi_1)$ and the triviality of the boundary spectrum. 

\medskip

\begin{proof}[Proof of Theorem~\ref{theo:Doblin}.]
{\Blue
The dual (equivalent) formulation of the Doblin condition~\eqref{eq:hyp-Doblin} is 
$$
\forall \, \psi \in X'_+, \quad S_T^* \psi  \ge \psi_0\langle \psi, g_0 \rangle.
$$
Applied to the first dual eigenvector $\phi_1$ it yields, by using~\eqref{eq:hyp-Doblin-BdBis},
\beqn\label{eq:stabilityKR-phi1bybelow}
\phi_1 = e^{-\lambda_1^*T} S_T^* \phi_1  \ge e^{-\lambda_1^*T}   \psi_0\langle \phi_1, g_0 \rangle =  e^{-\lambda_1^*T} r_0 \psi_0.
\eeqn
Together with~\eqref{eq:hyp-Doblin-Bd}, it implies that $[\cdot ]_{\phi_1}$ and $[\cdot ]_{\psi_0}$ are equivalent norms.
Another important consequence is that $\lambda_1=\lambda_1^*$ since we have
\[\lambda_1\langle\phi_1,f_1\rangle=\langle\phi_1,\LL f_1\rangle=\langle\LL^*\phi_1,f_1\rangle=\lambda_1^*\langle\phi_1,f_1\rangle\]
and $\langle\phi_1,f_1\rangle\neq0$ because $\phi_1\in X'_{++}$, due to~\eqref{eq:stabilityKR-phi1bybelow}, and $f_1\in X_+\setminus\{0\}$.

\smallskip

We now prove constructively the exponential stability \ref{E31intro} in the norm $[\cdot]_{\phi_1}$, which will finish the proof due to the (quantified) equivalence of the norms $[\cdot ]_{\phi_1}$ and $[\cdot ]_{\psi_0}$.
}
The two conditions \eqref{eq:hyp-Doblin} and \eqref{eq:hyp-Doblin-Bd} together imply  the modified
Doblin condition 
$$
\exists\, T > 0,  \ \exists \, g_1 > 0, \quad \forall \, f \ge 0, \ S_T f \ge g_1 \langle \phi_1, f \rangle,
$$
with $g_1 := g_0/R_0$. 
Take $f$ such that $\langle \phi_1, f \rangle = 0$, so that $\langle \phi_1, f_\pm \rangle = r := \langle \phi_1, |f| \rangle/2  \ge 0$ and thus 
$$
S_T f_\pm \ge  g_1  \langle \phi_1, f_\pm \rangle  =  r g_1 .
$$
We write 
$$
|S_Tf| \le |S_Tf_+-rg_1| + |S_Tf_- - r  g_1| = S_T|f| - 2 r g_1.
$$
We deduce 
$$
\langle \phi_1, |S_Tf| \rangle  \le \langle S_T^*\phi_1, |f| \rangle - 2 r  \langle \phi_1, g_1\rangle
= \Bigl(e^{\lambda_1 T} -  \langle \phi_1, g_1\rangle\Bigr)  \langle \phi_1, |f| \rangle.
$$
In other words, setting $\widetilde S_t  := e^{-\lambda_1 t} S_t$, we have 
\begin{equation}\label{eq:contraction_phi1}
[\widetilde S_T f]_{\phi_1} \le \gamma [f]_{\phi_1},  
\end{equation}
with $ \gamma < 1$ which depends explicitly of $r_0$, $R_0$, $T$  and the estimates on $\lambda_1$. 
We then classically deduce the exponential convergence in the $[\cdot]_{\phi_1}$ norm. 
\end{proof}

{\Blue
\begin{rem}
The existence of $(\lambda_1,f_1)$ is not required in the argument for proving the contraction property~\eqref{eq:contraction_phi1} of $e^{-\lambda_1^*T}S_T$.
This contraction can even be used for proving the existence of a solution $(\lambda_1,f_1)$ with $\lambda_1=\lambda^*_1$ to the direct problem~\eqref{eq:triplet1} by only assuming the existence of $(\lambda_1^*,\phi_1)$ to the dual problem~\eqref{eq:triplet2} together with Doblin's condition \eqref{eq:hyp-Doblin}-\eqref{eq:hyp-Doblin-Bd}-\eqref{eq:hyp-Doblin-BdBis}.
We refer to the proof of Theorem~\ref{theo:Harris} for more details about this approach.
\end{rem}
}
 
\medskip 
\subsection{Asymptotic stability under Doblin-Harris condition}
\label{subsec:stabilityStrong}

The Doblin condition \eqref{eq:hyp-Doblin}-\eqref{eq:hyp-Doblin-Bd}-\eqref{eq:hyp-Doblin-BdBis} is too much demanding for many applications. 
 In this section, we make the following somehow more general  Doblin-Harris type condition  complemented with a Lyapunov condition. 
 More precisely, we assume {\Blue the existence of a solution $(\lambda_1^*,\phi_1)$ to the dual eigenproblem~\eqref{eq:triplet2} and} that there exists $T > 0$ such that $\Blue \widetilde S_T :=e^{-\lambda_1^* T} S_T$ first satisfies the Lyapunov condition 
\beqn\label{eq:stabilityKR-Lyapunov}
\|  \widetilde S_T f \| \le \gamma_L \| f \| + K [ f ]_{\phi_1},
\eeqn
with $\gamma_L \in (0,1)$, $K \ge 0$. We next assume that $ \widetilde S_T $ satisfies the Doblin-Harris condition 
 \begin{equation}\label{eq:hyp-Harris}
\left\{
\begin{aligned}
&  \exists\,  A > K/(1-\gamma_L), \ \exists \, g_A > 0 \  \hbox{such that}
 \\
&\forall \, f \ge 0, \  \| f \| \le A [f]_{\phi_1} \ \hbox{there holds}  \ S_T f \ge g_A [f]_{\phi_1}. 
\end{aligned}
\right. 
\end{equation}
 We finally replace the positivity condition \eqref{eq:hyp-Doblin-BdBis} by
\beqn\label{eq:hyp-Harris-BdBis}
\exists \, r_A > 0, \quad \langle \phi_1, g_A \rangle \ge r_A. 
\eeqn
%
%
As we have seen several times, condition  \eqref{eq:stabilityKR-Lyapunov} is some kind of regularity hypothesis which is natural under a splitting structure on the semigroup $S_\LL$. 
We emphasize that conditions \eqref{eq:stabilityKR-Lyapunov}-\eqref{eq:hyp-Harris}-\eqref{eq:hyp-Harris-BdBis} slightly generalize the usual set of hypotheses for the Doblin-Harris theorem, see for instance \cite[Sect.~3]{MR4534707}.
We also point out that there is a connection between the condition~\eqref{eq:hyp-Harris} and the notion of {\it partial integral} or {\it partial kernel} operators, see for instance~\cite[Cor.~5.3]{Gerlach2024}.
The long term convergence of semigroups that contain a partially integral operator was studied in particular in~\cite{Pichor2000,Gerlach2013,Glueck2022}.

\begin{theo}\label{theo:Harris} 
Consider a {\Cyan positive} semigroup $S$ on a Banach lattice $X$ such that its generator $\LL$ {\Blue admits a solution $(\lambda_1^*,\phi_1)$ to the dual eigenproblem~\eqref{eq:triplet2}.}
We assume furthermore  the Doblin-Harris condition \eqref{eq:hyp-Harris} together with the Lyapunov condition \eqref{eq:stabilityKR-Lyapunov} and the positivity condition \eqref{eq:hyp-Harris-BdBis}.
Then {\Blue there exists a solution $(\lambda_1,f_1)$ with $\lambda_1=\lambda_1^*$ to the direct eigenproblem~\eqref{eq:triplet1}, and} the exponential stability {\ref{E31intro}}
holds true in the norm of $X$ with constructive constants.  
 \end{theo}

Of course, in order that Theorem~\ref{theo:Harris} really gives a constructive convergence result, we have to establish \eqref{eq:hyp-Harris},  \eqref{eq:stabilityKR-Lyapunov}
 and \eqref{eq:hyp-Harris-BdBis} in a constructive way.

\begin{proof}[Proof of Theorem~\ref{theo:Harris}.]
On the one hand, we have 
\beqn\label{eq:contractNpsi1}
[  \widetilde S_T f ]_{\phi_1} \le  \langle  \widetilde S_T  |f|,\phi_1 \rangle =\langle  |f|,   \widetilde S_T^* \phi_1 \rangle = [ f ]_{\phi_1}.
\eeqn
On the other hand, we  wish to establish the coupling property
\beqn\label{eq:couplingNpsi1}
[  \widetilde S_T f ]_{\phi_1} \le \gamma_H [ f ]_{\phi_1} \quad\hbox{if}\quad \| f \| \le {A'} [ f ]_{\phi_1} \ \hbox{and} \ \langle f,\phi_1\rangle = 0, 
\eeqn
for   some $\gamma_H \in (0,1)$ and with $A' := A/2$.
We thus consider $f \in X$, such that $\langle f,\phi_1\rangle = 0$ and $\| f \| \le {A'} [ f ]_{\phi_1}$, so that 
$$
\| f_\pm \| \le  \| f \| \le {A'}[ f ]_{\phi_1} = {A} [ f_\pm ]_{\phi_1}. 
$$
Using the Doblin-Harris condition \eqref{eq:hyp-Harris}, we deduce 
$$
 \widetilde S_T f_\pm \ge \vartheta g_{{A}}, \quad \vartheta := \tfrac12 e^{-\lambda_1 T} [ f ]_{\phi_1}.
$$
Similarly as in the proof of Theorem~\ref{theo:Doblin}, we next compute 
$$
| \widetilde S_T f| \le | \widetilde S_T f_+-  \vartheta g_{{A}}| + | \widetilde S_T f_- -  \vartheta g_{{A}}| \le  \widetilde S_T |f| - 2  \vartheta g_{{A}}
$$
and then 
\bean
[  \widetilde S_T f ]_{\phi_1}   
&\le& \langle  \widetilde S_T |f| - 2  \vartheta g_{{A}}, \phi_1 \rangle 
\\
&=&  \langle |f| ,  \widetilde S_T^* \phi_1 \rangle - 2  \vartheta  \langle  g_{{A}}, \phi_1 \rangle 
\\
&=& \bigl( 1 - e^{-\lambda_1 T}  \langle  g_{{A}}, \phi_1 \rangle\bigr)     [ f ]_{\phi_1},
\eean
which in turn implies  \eqref{eq:couplingNpsi1} with $\gamma_H := 1 - e^{-\kappa_0 T} r_{{A}}$.  

\medskip

Now, the two estimates \eqref{eq:contractNpsi1} and \eqref{eq:couplingNpsi1} together give
\beqn\label{eq:estimNpsi1}
[\widetilde S_T f ]_{\phi_1} \le  \gamma_H  [ f ]_{\phi_1} + \frac{1-\gamma_H }{ A'} \| f \|. 
\eeqn
From \eqref{eq:estimNpsi1} and the Lyapunov condition \eqref{eq:stabilityKR-Lyapunov}, we deduce that \Black
$$
U^{n+1} = M U^n
$$
with 
$$
U^n := \left(
       \begin{matrix}
        \| \widetilde S_T^n f \|
         \\
     [ \widetilde S_T^nf ]_{\phi_1}
       \end{matrix}
     \right)
     = \left(
       \begin{matrix}
        \| \widetilde S_{nT} f \|
         \\
     [ \widetilde S_{nT} f ]_{\phi_1}
       \end{matrix}
     \right)
\quad\hbox{and}\quad
     M := \left(
       \begin{matrix}
         \gamma_L & K
         \\
        \frac{1-\gamma_H}{A} & \gamma_H 
       \end{matrix}
     \right).
$$
The eigenvalues $\mu_\pm$ of $M$ are 
$$
\mu_\pm := \frac12 \bigl(T \pm \sqrt{T^2 - 4D} \bigr) , 
$$
with 
$$
 T := \hbox{\rm tr} M = \gamma_L + \gamma_H, \quad D  := \hbox{\rm det} M =  \gamma_L  \gamma_H - (1-\gamma_H) \frac{K }{ A'}.
$$
We observe that 
$$
\gamma_L  \gamma_H > D>  \gamma_L  \gamma_H - (1-\gamma_H) (1-\gamma_L) = T -1, 
$$
so that 
$$
(\gamma_H-\gamma_L)^2 = T^2 - 4 \gamma_L  \gamma_H <T^2 - 4D < T^2 - 4(T-1) = (T-2)^2
$$
and finally 
$$
\alpha := \max(|\mu_+|,|\mu_-|) < \max(\gamma_H,\gamma_L,|T-1|,1) = 1.
$$
We conclude that $\| M^n \|  \lesssim \alpha^n$, from what we immediately
{\Blue deduce the existence of $C>0$ such that
\begin{equation}\label{eq:contraction_Harris}
\|\widetilde S_{nT}f\|\leq C\alpha^n\|f\|
\end{equation}
for all $f\in X$ such that $\langle\phi_1,f\rangle=0$.
This classically yields~\ref{E31} and it only remains to prove the existence of $f_1$.

We first observe that~\eqref{eq:contraction_Harris} ensures that ${(\|\widetilde S_{nT}f_0\|)}_n$ is a Cauchy sequence for any $f_0\in X$. Indeed,  for any $p\in\N$,  $f=f_0-\widetilde S_{pT} f_0$ verifies
\[
\langle f,\phi_1\rangle = \langle f_0,\phi_1\rangle - \langle f_0, \widetilde S_{pT}^*\phi_1\rangle = \langle f_0,\phi_1\rangle - \langle f_0,\phi_1\rangle =0,
\]
and we then have
\[
\|\widetilde S_{nT}f_0-\widetilde S_{(n+p)T}f_0\|+[\widetilde S_{nT}f_0-\widetilde S_{(n+p)T}f_0]_{\phi_1}\lesssim \alpha^n \big(\|f_0-\widetilde S_{pT} f_0\|+[f_0-\widetilde S_{pT} f_0]_{\phi_1}\big).
\]
Choosing $f_0\in X_+$ such that $[f_0]_{\phi_1}=1$, we deduce that $(\widetilde S_{nT}f)$ converges to a fixed point $f_1$ of $\widetilde S_T$, which is not zero because
\[[f_1]_{\phi_1}=\lim[\widetilde S_{nT}f_0]_{\phi_1}=[f_0]_{\phi_1}=1,\]
and, thanks to \eqref{eq:contraction_Harris},  $f_1$ is the unique fixed point with normalization $[f_1]_{\phi_1}=1$.
Besides, $f_1\in X_+$ because of the positivity of $S$ and $f_0$.
This ensures that
\[
[\widetilde S_t f_1]_{\phi_1} = \langle \widetilde S_t f_1 ,\phi_1\rangle = \langle f_1 , \widetilde S_t^*  \phi_1\rangle = \langle f_1, \phi_1\rangle = 1,
\]
for any $t>0$.
Since on the other hand
\[\widetilde S_T \widetilde S_t f_1 = \widetilde S_{t+T} f_1 = \widetilde S_t \widetilde S_T f_1 = \widetilde S_t f_1,\]
we deduce from the uniqueness of the fixed point that $\widetilde S_t f_1 = f_1$, which yields that $f_1\in D(\LL)$ and $\LL f_1 = \lambda_1^* f_1$.}
\end{proof}

\begin{rem}\Blue
We have seen in Theorem~\ref{theo:Doblin} that under Doblin's condition \eqref{eq:hyp-Doblin}-\eqref{eq:hyp-Doblin-Bd}-\eqref{eq:hyp-Doblin-BdBis}, if we have two eigenpairs $(\lambda_1,f_1)$ and $(\lambda_1^*,\phi_1)$ solutions to the direct and dual eigenvalue problems \eqref{eq:triplet1} and~\eqref{eq:triplet2} respectively,
then necessarily $\lambda_1=\lambda_1^*$.
This is not the case for the conditions \eqref{eq:stabilityKR-Lyapunov}-\eqref{eq:hyp-Harris}-\eqref{eq:hyp-Harris-BdBis}, which only ensure that $\lambda_1\leq\lambda_1^*$.
\end{rem}

\subsection{Quantified isolation of the first eigenvalue}
\Black
In terms of the geometry of the spectrum, an immediate consequence of Theorem~\ref{theo:Harris} is that the conditions~\eqref{eq:hyp-Harris}-\eqref{eq:stabilityKR-Lyapunov}-\eqref{eq:hyp-Harris-BdBis} ensure the existence of a spectral gap, namely the existence of $\eps>0$ such that
\[\Sigma(\LL)\cap\Delta_{\lambda_1-\eps}=\{\lambda_1\}.\]
We still assume that the Lyapunov condition~\eqref{eq:stabilityKR-Lyapunov} holds for some $T >0$, $\gamma_L \in (0,1)$ and $K \ge 0$, but we relax~\eqref{eq:hyp-Harris} into the time-averaged condition
\begin{equation}\label{eq:hyp-Harris-mean}
\left\{
\begin{aligned}
&   \exists\,  A > K/(1-\gamma_L), \ \exists \, g_A > 0 \  \hbox{such that}
 \\
&\forall \, f \ge 0, \  \| f \| \le A [f]_{\phi_1} \ \hbox{there holds}  \ \int_0^T\! S_t f\,dt \ge g_A [f]_{\phi_1}.
\end{aligned}
\right. 
\end{equation}
It is worth emphasizing that \eqref{eq:hyp-Harris-mean} does not imply anymore the existence of a spectral gap, and there can be a non-trivial boundary spectrum, see Section~\ref{ssec:GF:singular} for an example. However, it is strong enough  for guaranteeing that $\lambda_1$ is isolated from the rest of the spectrum, in the sense that 
\begin{equation}\label{eq:isolation}
\Sigma(\LL)\cap B(\lambda_1,\eps)=\{\lambda_1\}, 
\end{equation}
for some $\eps>0$.
In particular, if not trivial, the boundary spectrum must be discrete from Theorem~\ref{theo:Sigma+subgroup} (under the additional assumptions listed in the statement of this last result).

\begin{theo}\label{theo:Harris-mean}  
Consider a {\Cyan positive} semigroup $S$ on a Banach lattice $X$ such that its generator $\LL$
{\Blue admits a solution $(\lambda_1^*,\phi_1)$ to the dual eigenproblem~\eqref{eq:triplet2}.}
We assume furthermore  the time-averaged Doblin-Harris condition~\eqref{eq:hyp-Harris-mean} together with the Lyapunov condition~\eqref{eq:stabilityKR-Lyapunov} and the positivity condition~\eqref{eq:hyp-Harris-BdBis}.
Then {\Blue there exists a solution $(\lambda_1,f_1)$ with $\lambda_1=\lambda_1^*$ to the direct eigenproblem~\eqref{eq:triplet1}, and} \eqref{eq:isolation} holds true for some constructive constant $\eps>0$,  providing a quantified information on~\ref{S31}.
\end{theo}


\begin{proof}[Proof of Theorem~\ref{theo:Harris-mean}]
{\Blue We start by noticing that by iteration of~\eqref{eq:stabilityKR-Lyapunov} we have
\[
\|\widetilde S_{nT} f\| \leq \gamma_L^n \, \|f\| + \frac{K}{1-\gamma_L} [f]_{\phi_1}, 
\]
for all integer $n$, from which we deduce 
\begin{equation}\label{eq:borne-Lyap-itere}
\|\widetilde S_t f\| \leq C \gamma_L^{\lfloor t/T\rfloor} \, \|f\| + \frac{C K}{1-\gamma_L} [f]_{\phi_1},
\end{equation}
for all $t\geq0$,
where $C=\sup_{0\leq t\leq T}\|\widetilde S_t\|$. 
In particular
\[\|\widetilde S_t\|\leq C+\frac{C K}{1-\gamma_L}\|\phi_1\|,\]
which implies that $\omega(\widetilde S)\leq0$ and consequently $\lambda-\LL$ is invertible for all $\lambda>\lambda_1^*$ with inverse $\RR_\LL(\lambda)$ given by the inversion formula~\eqref{eq:Exist1-DefRepresentationRR}.
Using this inversion formula   and~\eqref{eq:hyp-Harris-mean}, we readily deduce that }
\begin{equation}\label{eq:hyp-Harris-resolv}
\left\{
\begin{aligned}
& \exists\,  A > K/(1-\gamma_L), \ \exists \, g_A > 0 \  \hbox{such that}
 \\
&\forall \, f \ge 0, \  \| f \| \le A [f]_{\phi_1} \ \hbox{there holds}  \ \widetilde\RR(\lambda) f \ge   \widetilde g_A   [f]_{\phi_1}, \  \forall \, \lambda > \lambda_1^*, 
\end{aligned}
\right. 
\end{equation}
where $\widetilde \RR(\lambda) :=(\lambda-\lambda_1^*)\RR_\LL(\lambda)$ and $\widetilde g_A := (\lambda-\lambda_1^*) e^{- \lambda T} g_A$.
Next, we claim that the Lyapunov condition~\eqref{eq:stabilityKR-Lyapunov} ensures the existence of $\lambda>\lambda_1^*$ such that
\beqn\label{eq:stabilityKR-Lyapunov-resolv}
\|  \widetilde\RR(\lambda)  f \| \le \gamma'_L \| f \| + K' [ f ]_{\phi_1}
\eeqn
for all $f \in X$ and some $\gamma'_L < 1$ and $K'>0$.
Indeed, 
 we infer from~\eqref{eq:borne-Lyap-itere} and  
the inversion formula~\eqref{eq:Exist1-DefRepresentationRR} that
\[
\| \widetilde\RR(\lambda) f\| \leq \frac{C_1(\lambda-\lambda_1^*)}{\lambda-\lambda_1^*+\log\frac{1}{\gamma_L}} \|f\| + \frac{C_2}{1-\gamma_L} [f]_{\phi_1},
\]
for all $\lambda>\lambda_1^*$ and some $C_1,C_2>0$.
Then we only need to choose $\lambda$ close enough to $\lambda_1^*$ so that $\frac{C_1(\lambda-\lambda_1^*)}{\lambda-\lambda_1^*+\log\frac{1}{\gamma_L}}<1$ and we obtain~\eqref{eq:stabilityKR-Lyapunov-resolv}.


We have proved that $\widetilde\RR(\lambda)$ satisfies~\eqref{eq:stabilityKR-Lyapunov-resolv} and~\eqref{eq:hyp-Harris-resolv}. 
Together with the positivity condition~\eqref{eq:hyp-Harris-BdBis}, we can thus repeat the proof of Theorem~\ref{theo:Harris} for the operator $\widetilde\RR$ instead of $\widetilde S$ and we obtain
{the existence of $f_1$ such that $\widetilde \RR(\lambda)f_1=f_1$, which is equivalent to $\LL f_1=\lambda_1^*f_1$, and}
the existence of constructive constants $\alpha\in(0,1)$ and $C \ge 1$ such that
$$
\| \widetilde\RR (\lambda)^n f \| \le C \alpha^n \|  f \| , \quad \forall \, n \ge 1,  
$$
for any $f \in X$, $\langle f, \phi_1 \rangle = 0$. By the spectral radius formula {(see \cite[Thm.~10.13]{MR1157815}  for instance)}, we deduce {that
\[r(\widetilde R(\lambda)=\lim_{n\to\infty}\| \widetilde\RR (\lambda)^n \|^{1/n}\leq\alpha\]
and so}
\[
 \Sigma\big(\widetilde\RR(\lambda)\big)\cap \{z\in  \C,\, |z|>\alpha\}=\{1\}.
 \]
 The spectral mapping theorem for the resolvent {\Blue(see for instance~\cite[Thm. IV.1.13]{MR1721989}}) ensures that
 \[\Sigma\big(\widetilde\RR(\lambda)\big)\setminus\{0\}=\frac{\lambda-\lambda_1^*}{\lambda-\Sigma(\LL)},\]
 then yields~\eqref{eq:isolation} with $\eps=(\alpha^{-1}-1)(\lambda-\lambda_1^*)$.
\end{proof}

\Black

 \medskip
\subsection{The weak dissipativity case}
\label{subsec:stabilityWeak}

In this section, we consider a weak dissipative semigroup $(S_t)$ as considered in  Section~\ref{subsect:AboutWeakDissip} and in a sense we make precise now. 
We assume  that its generator $\LL$ {\Blue admits a solution $(\lambda_1^*,\phi_1)$ to the dual eigenproblem~\eqref{eq:triplet2}.}
We consider four Banach lattices $X_3  \subset X_2  \subset X_1  \subset X_0 =X$. 
We first make the same kind of Doblin-Harris type condition as in  the previous section, namely

\smallskip
{\bf Hypothesis (H)} (Doblin-Harris)
condition \eqref{eq:hyp-Harris} holds for the same time $T > 0$ and for both norms $\| \cdot \| = \| \cdot \|_{X_0}$ and $\| \cdot \| = \| \cdot \|_{X_2}$ as well as the companion positivity condition 
\eqref{eq:hyp-Harris-BdBis} holds.

\smallskip 
 Instead of the strong Lyapunov condition \eqref{eq:stabilityKR-Lyapunov}, we assume

\smallskip
{\bf Hypothesis (L)} (weak Lyapunov) there exist a constant $K \ge 0$ 
such that
\bean
\|  \widetilde S f \|_{1} + \|  \widetilde S f \|_{0}  &\le& \| f \|_{1} +  K [ f ]_{\phi_1}, \quad \forall \, f \in X_1, 
\\
\|  \widetilde S f \|_{3} + \|  \widetilde S f \|_{2}  &\le& \| f \|_{3} +  K  [ f ]_{\phi_1}, \quad \forall \, f \in X_3,
\eean
with $ \widetilde S = S_T e^{-\lambda^*_1T}$.

\smallskip
{\bf Hypothesis (I)}  (interpolation) there exists  an increasing function $\xi : \R_+ \to \R_+$, $\lambda \mapsto \xi_\lambda$, such that 
$$
\lambda \| f \|_{1}  \le   \| f \|_{0}  + \xi_\lambda \| f \|_{3} , \ \forall \, \lambda > 0, \quad  \xi_\lambda/\lambda\to0  \ \hbox{ as} \ \lambda \to 0. 
$$
 
\smallskip

\begin{theo}\label{theo:HarrisSubgeo1} 
Consider a semigroup $S$ on a Banach lattice $X$ such that its generator $\LL$ {\Blue admits a solution $(\lambda_1^*,\phi_1)$ to the dual eigenproblem~\eqref{eq:triplet2}.}
We assume furthermore 
 the three above conditions of weak confinement {\bf (L)}, 
Doblin-Harris strong irreducibility  {\bf (H)} 
and interpolation {\bf (I)}. 
Then,  there exist some constructive decay rate functions $\Theta$ and $ \widetilde \Theta$ such that  
\beqn\label{eq:theo:HarrisSubgeo1:Estim1}
\| \widetilde S^n f \|_{X_1}   \lesssim  \Theta(n) \| f \|_{X_3}, \quad \forall \, n \ge 1, 
\eeqn
and 
\beqn\label{eq:theo:HarrisSubgeo1:Estim2}
   \| \widetilde S^n f \|  \lesssim  \widetilde \Theta (n)  \| f \|_{X_3}, \quad \forall \, n \ge 1, 
\eeqn
for any $f \in X_{3}$, $\langle f , \phi_1\rangle = 0$. More precisely, the decay rate functions $\Theta$ and $ \widetilde \Theta$ are defined by 
\beqn\label{def:Theta&TildeTheta}
\Theta (t) := \inf_{\lambda} \Theta_{\zeta \lambda} (t), \quad \widetilde  \Theta (t) := t^{-1} \Theta([t/2]),
\eeqn
for a constructive constant $\zeta \in (0,1)$, the infimum being taken over all the decreasing function $\lambda : \R_+ \to \R_+$, $t \mapsto \lambda_t$, and $\Theta_\lambda$ is defined by
\beqn
\label{}
\Theta(t) :\simeq \inf_{\lambda > 0} \bigl( e^{-\lambda t} + \frac{\xi_\lambda}\lambda \bigr).
\eeqn
\end{theo}

The proof follows closely the proof of \cite[Thm.~4.8]{MR4534707}. 
We start with the following key argument of non expansive mapping result on a well chosen norm. 

\begin{prop}\label{prop:SubHarris1}   Consider a positive semigroup $(S_t)$ which satisfies both above conditions of weak confinement {\bf (L)} and  
Doblin-Harris strong irreducibility {\bf (H)}. 
There exist some equivalent norms $\Nt \cdot \Nt_{1}$ to $\| \cdot \|_{1}$ and $\Nt \cdot \Nt_{3}$ to $\| \cdot \|_{3}$ such that $\widetilde S_t $ is a non expansive mapping for the two new norms $\Nt \cdot \Nt_{1}$  and $\Nt \cdot \Nt_{3}$. 
More precisely, there exists $\alpha > 0$ such that 
\bear
\label{eq1:prop:subHarris1}
\NormT  \widetilde S f \NormT_{1}  + \alpha \|  \widetilde S f \|_{0}  &\le& \NormT  f \NormT_{1} , \quad \forall \, f \in X_{1},  \,\, \langle f, \phi_1\rangle = 0, 
\\
\label{eq2:prop:subHarris1}
\NormT  \widetilde S f \NormT_{3}  + \alpha \|  \widetilde S f \|_{2}  &\le& \NormT  f \NormT_{3} , \quad \forall \, f \in X_{3},  \,\, \langle f, \phi_1 \rangle = 0. 
\eear
\end{prop}
 
\begin{proof}[Proof of Proposition~\ref{prop:SubHarris1}.]
We define 
\beqn\label{def:NormeTriple}
\NormT  f \NormT_{1} := [ f ]_{\phi_1}  + \delta \| f \|_{0} +   \beta \| f \|_{1},  
\eeqn
with $\beta >\delta > 0$ conveniently chosen. We take $\beta := (1-\gamma_H)/K$, $\delta := (1-\gamma_H)/A$. We define 
$\NormT \cdot \NormT_{3}$ in the same way. 
In what follows, we then only establish \eqref{eq1:prop:subHarris1}, the proof of  \eqref{eq2:prop:subHarris1} being exactly the same. 

\smallskip
We fix $f \in X_{1}$, $ \langle f,\phi_1\rangle = 0$,  and we recall
\beqn\label{eq:prop:SubHarris0}
[  \widetilde S f]_{\phi_1}  \le  [f]_{\phi_1}. 
\eeqn
We also 
recall that from \eqref{eq:couplingNpsi1}, for any $A > 0$, there exists $\gamma_H = \gamma_H(A) \in (0,1)$ such that  the following coupling property holds
\beqn\label{eq:couplingNphi1}
[  \widetilde Sf ]_{\phi_1} \le \gamma_H [ f ]_{\phi_1} \quad\hbox{if}\quad \| f \|_0 \le A  [ f ]_{\phi_1}. 
\eeqn

We fix $A > K$ and we observe that the following alternative holds
\beqn\label{eq:prop:Harris1}
\|  f \|_{0} \le A [f]_{\phi_1} 
\eeqn
or
\beqn\label{eq:prop:Harris2}
\|  f \|_{0}  >  A  [f]_{\phi_1}.
\eeqn
 
\smallskip\noindent
{\sl Case 1. } Under condition \eqref{eq:prop:Harris1},  we use  \eqref{eq:couplingNphi1} and  the first estimate in {\bf (L)}, and we deduce  
\bean
 \NormT   \widetilde Sf  \NormT_{1} 
&=& 
[ \widetilde S f]_{\phi_1}   + \delta \|  \widetilde Sf \|_{0} +  \beta \|  \widetilde Sf \|_{1}
\\&\le&  \gamma_H [f]_{\phi_1}   + \beta  \| f \|_{1}  +  \beta K  [f]_{\phi_1}  -  (\beta -\delta) \|  \widetilde S f \|_{0}. 
\eean
From our choice of $\beta > 0$ we have  $\gamma_H+ \beta K = 1$, and we conclude that  \eqref{eq1:prop:subHarris1} holds with $\alpha := \beta-\delta > 0$. 

\smallskip\noindent
{\sl Case 2. } Under condition \eqref{eq:prop:Harris2},   the first Lyapunov condition in {\bf (L)} implies 
$$
\|  \widetilde Sf \|_{1} + \|  \widetilde Sf \|_{0} \le \|  f \|_{1} + \frac{K}A \|  f  \|_{0}.
$$
Together with  the non expansivity estimate \eqref{eq:prop:SubHarris0}, we get 
$$
[ \widetilde Sf]_{\phi_1}  + \beta \|  \widetilde S f \|_{1} + \beta \|  \widetilde S f \|_{0} \le [f]_{\phi_1} +  \beta \|  f \|_{1} + \delta \|  f  \|_{0},
$$
and we conclude to \eqref{eq1:prop:subHarris1} again. 
\end{proof}

The subgeometric convergence result is a straightforward consequence of Proposition~\ref{prop:SubHarris1}
and an interpolation argument.

\begin{prop}\label{prop:Harris1SubgeoConv} Assume that $S$   satisfies the hypotheses of Theorem~\ref{theo:HarrisSubgeo1}.
Then \eqref{eq:theo:HarrisSubgeo1:Estim1} and \eqref{eq:theo:HarrisSubgeo1:Estim2} hold true  with the same decay rate functions  $\Theta$ and $\widetilde  \Theta$ given by 
\eqref{def:Theta&TildeTheta}  (up to a modification of the constant  $\zeta$).
\end{prop}

\begin{proof}[Proof  of Proposition~\ref{prop:Harris1SubgeoConv}.]
We recall that we have already proven \eqref{eq1:prop:subHarris1} and \eqref{eq2:prop:subHarris1}.
From \eqref{eq1:prop:subHarris1}  and the  interpolation condition {\bf(I)}, we deduce 
$$
\NormT \widetilde S f \NormT_{1} + \lambda \alpha \| \widetilde S f \|_{1}     \le  \NormT  f  \NormT_{1} +   \xi_\lambda \alpha  \| \widetilde S f \|_{3}. 
$$
We observe next that from the very definition of the $\NormT\cdot \NormT_{1}$ norm 
$$
 \NormT \widetilde S f \NormT_{1} + \frac\alpha\lambda \| \widetilde S f \|_{1}   \ge Z_\lambda \NormT \widetilde S f \NormT_{1}, \quad Z_\lambda = 1 + \kappa\lambda \in (1,2],
$$
for some $\kappa > 0$ and any $\lambda \in (0,\lambda_0)$, $ \lambda_0 > 0$, and that from the very  definition of the $\NormT\cdot \NormT_{3}$ norm 
$$
 \alpha  \xi_\lambda  \| \widetilde S f \|_{3} \le B  \xi_\lambda   \NormT \widetilde S f \NormT_{3}, 
 $$
 for some $B > 0$. The three above estimates together imply 
 $$
Z_\lambda  \NormT \widetilde S f \NormT_{1}      \le  \NormT  f  \NormT_{1} + B  {\xi_\lambda}   \NormT \widetilde S f  \NormT_{3}. 
$$
Using the second estimate \eqref{eq2:prop:subHarris1} and repeating the same proof, we have
$$
Z_{\lambda_{n+1}}  \NormT \widetilde S^{n+1} f \NormT_{1}      \le  \NormT \widetilde S^n  f  \NormT_{1} + B  {\xi_{\lambda_{n+1}}}   \NormT   f  \NormT_{3},
$$
for any $n\ge0$ and for any $\lambda_{n+1} > 0$. The discrete Gr{\"o}nwall lemma implies 
\beqn
\label{eq:GronwallDconclusion}
\NormT \widetilde S^n  f  \NormT_{1}   \le A_n \NormT  f  \NormT_{1}  + \sum_{k=1}^{n} A_{k,n}    \xi_{\lambda_{k}}   B \NormT   f  \NormT_{3}, \quad \forall \, n \ge 0, 
\eeqn  where we have defined  
$$
A_n := \prod_{k=1}^{n} a_k, \quad A_{k,n} = A_n/A_k =  \prod_{i=k+1}^{n} a_i, \quad a_i := Z_{\lambda_i}^{-1}.
$$  
Observing that 
$$
A_{k,n} \lesssim e^{-\kappa \sum_{i=k}^n \lambda_i} \lesssim e^{\kappa (\Lambda(k) - \Lambda(n))},  
\quad\hbox{with}\quad\Lambda(t) := \int_0^t \lambda_s \, ds, 
$$
and $\lambda_s := \lambda_i$ if $s \in (i-1,i]$, we immediately conclude that the first estimate \eqref{eq:theo:HarrisSubgeo1:Estim1} holds true.  
We come back to the first inequality in \eqref{eq1:prop:subHarris1} that we iterate and sum up in order to obtain
$$
\NormT \widetilde S^n f \NormT_{1} + \alpha \sum_{k=[n/2]+1}^n \| \widetilde S^k f \|_{0}     \le  \NormT \widetilde S^{[n/2]} f  \NormT_{1} ,
$$
for any $n \ge 1$. Together with the non expansion inequality 
$$
[ \widetilde S^n f ]_{\phi_1} \le [ \widetilde S^k f ]_{\phi_1} \lesssim  \| \widetilde S^k f \|_{0}, \quad \forall \, n \ge k, 
$$
and the first estimate \eqref{eq:theo:HarrisSubgeo1:Estim1}, we deduce 
$$
\bigl( n - [n/2] - 1 \bigr) \alpha [ \widetilde S^n f ]_{\phi_1}  \lesssim \Theta([n/2]) \NormT   f  \NormT_{3},
$$
which is nothing but \eqref{eq:theo:HarrisSubgeo1:Estim2}.
\end{proof}

{\Blue
 
Let us finally emphasize that \ref{C2} together with Theorem~\ref{theo:Doblin} or Theorem~\ref{theo:Harris} provide the  exponential stability \ref{E31intro}
 with constructive constants, and thus also the  spectral gap property \ref{S33}. 
 On the other hand, \ref{C2} together with Theorem~\ref{theo:HarrisSubgeo1} provide a quantitative stability estimate  \ref{E32intro}.
 Both together, we deduce Theorem~\ref{theo:main-intro}-(5). 
 }

%
%
%
%
%
%


\bigskip
%
%

 \bigskip
\section{Parabolic equations}
\label{sec:application1:diffusion}


%
%
%
%
%
%
%
%
%

 In this part, we consider a general elliptic operator in divergence form
\beqn\label{eq:exPLL1}
\LL f := \partial_i(a_{ij} \partial_j f) + b_i \partial_i f +  \partial_i(\beta_i f) +  cf, \quad f \in H^1_0(\Omega),
\eeqn
where  $\Omega \subset \R^d$ is  a bounded domain (i.e. an open and connected set)  or $\Omega = \R^d$, and we always assume $d \ge 3$ (in order to simplify the discussions when using the Sobolev inequality).
We also always assume at least a boundedness and ellipticity condition on the $(a_{ij})$ matrix, namely 
\beqn\label{eq:StampHypaij}
a_{ij}  \in L^\infty(\Omega), \quad  \exists \nu > 0, \ \forall \xi \in \R^d, \ a_{ij} \xi_i\xi_j \ge \nu |\xi|^2, 
\eeqn 
and some conditions on the coefficients $b_i$, $\beta_j$ and $c$ which will be described below. 

\smallskip
We aim to establish the existence of $(\lambda_1,f_1,\phi_1)$ solution to the first eigentriplet problem
\beqn\label{eq:ex1PLL-triplet}
\lambda_1 \in \R, \quad  0 < f_1 \in H^1_0, \ \LL f_1 = \lambda_1 f_1,  \quad 0 < \phi_1 \in H^1_0, \ \LL^* \phi_1 = \lambda_1 \phi_1,
\eeqn
and the existence of some (constructive) rate function $\Theta$   such that the rescaled semigroup $\widetilde S$ associated to the generator $\widetilde\LL = \LL - \lambda_1$ satisfies 
\beqn\label{eq:ex1PLL-stab}
\| \widetilde S(t) f - \langle f, \phi_1 \rangle f_1 \|_{H_0} \le \Theta(t) \|   f - \langle f, \phi_1 \rangle f_1 \|_{H} , 
\eeqn
for any $t \ge 0$ and any $f \in H$, with $H \subset H_0 \subset L^2$.

 \medskip
\subsection{Diffusion with rough coefficients in a bounded domain} 
\label{subsec:diffusion-domain}
In this section, we consider the general elliptic operator in divergence form \eqref{eq:exPLL1} in the case of a bounded and smooth enough domain 
$\Omega \subset \R^d$ with general elliptic condition on $a_{ij}$ as formulated  above.
We further assume that
\beqn\label{eq:StampDivergenceHypH1}
b_i, \beta_j \in L^r(\Omega), \quad c  \in L^{r/2}(\Omega), \quad r  > d.  
\eeqn   

In that situation, the  first eigentriplet problem \eqref{eq:ex1PLL-triplet} has been considered by Chicco in  \cite{MR280858,MR0310462} and revisited
in a slightly less general framework (all the coefficients belong to $L^\infty$) in \cite{PL2}, where the conclusions \ref{C2} are established. 
%
We explain with all details the existence proof by following more or less the arguments presented in  \cite{PL2} stressing on the constructive way for obtaining the estimates, 
and next we present a proof of the geometric part and the stability part by taking advantage of the abstract material developed in the previous sections. It is worth emphasizing that our proof of the uniqueness of the first eigenfunction significantly differs from the one presented in  \cite{PL2} which is based on a dissipativity argument, probably related to the reverse Kato's inequality condition. 
The framework considered here  is the usual generalized solutions or weak solutions framework which goes back at least to Stampacchia \cite{MR192177,MR0251373}, but it is reminiscent of previous contributions by
Friedrichs \cite{MR9701,MR58828}, G\aa rding \cite{MR64979}, De Giorgi \cite{MR0093649}, Nash \cite{MR100158},  Morrey \cite{MR0120446}, Moser \cite{MR170091,MR159138,MR159139}, Ladyzhenskaya,  Solonnikov,  Ural'ceva \cite{MR0244627,MR0241822}, Oleinik, Kruzhkov \cite{MR0141892} and many others.
Lot of the functional arguments are picked up from the book of Gilbarg and Trudinger, and more specifically from  \cite[Chapter~8]{MR0473443}, and also in recent notes by Kavian~\cite{MR2187812} and Vasseur~\cite{NotesAV}. It is worth emphasizing that the present analysis does not apply directly to elliptic operators in non divergence form, although this framework is considered in \cite{PL2}. We expect that all the results developed below can be generalized to a non divergence form  framework, for example the one developed in \cite{MR1258192}, but we do not follow this line of research in the present work.


 
       
%

\smallskip
The proof of \eqref{eq:ex1PLL-triplet} and \eqref{eq:ex1PLL-stab} are straightforward consequences of the abstract results developed in the previous sections once we have been able to check that the corresponding hypotheses are fulfilled. In the sequel, we will then show how these hypotheses are met in the present context.

\smallskip
{\bf Condition \ref{H1}.} We recall that a weak (or variational) solution to the 
elliptic equation 
$$
\LL f = g \in H^{-1}(\Omega), \quad f \in H^1_0(\Omega), 
$$
is a function $f \in H^1_0(\Omega)$ such that 
\beqn\label{def:StampOperator1}
D_\LL(f,w) = \langle g,w\rangle, \quad \forall \, w \in H^1_0(\Omega), 
\eeqn
where the (negative) Dirichlet form $D_\LL$ is defined by 
$$
D_\LL(f,w) := - \int_\Omega (a_{ij} \partial_j f + \beta_i f  )\partial_i w + \int_\Omega (b_i \partial_i f w + cf w),
$$
for any $f,w \in H^1_0(\Omega)$. Most of the time, we will simply write 
\beqn\label{def:StampOperator}
\langle \LL f, w \rangle  = \langle g,w\rangle,  \quad \forall \, w \in H^1_0(\Omega), 
\eeqn
instead of \eqref{def:StampOperator1}. 
For the reader convenience, we repeat here some estimates picked up in \cite{MR0251373}. 
For $\lambda \in \R$ and  $f \in H^1_0(\Omega)$, we start with
\bean
\langle(\lambda-\LL) f, f\rangle
&=&\int_\Omega a_{ij} \partial_i f \partial_j f  + \int_\Omega (\beta_i - b_i) \partial_i f  f  + \int_\Omega (\lambda-c) f^2 
\\
&\ge& \| f \sqrt{c_-} \|_{L^2}^2+ \nu \| \nabla f \|_{L^2}^2 - \| |\beta-b| f \|_{L^2} \| \nabla f \|_{L^2}  -   \|  \sqrt{c_+} f \|_{L^2}^2 + \lambda \|  f \|_{L^2}^2
\\
&\ge& \| f \sqrt{c_-} \|_{L^2}^2+ \frac{\nu}2  \| \nabla f \|_{L^2}^2 - \frac1{2\nu} \| |\beta-b|  f \|_{L^2}^2  -   \|  \sqrt{c_+} f \|_{L^2}^2 
 +  \lambda  \|  f \|_{L^2}^2, 
\eean
using the Cauchy-Schwarz inequality and the Young inequality, and next
\bean
\langle (\lambda-\LL) f, f \rangle
&\ge& \| f \sqrt{c_-} \|_{L^2}^2+ \frac{\nu}4  \| \nabla f \|_{L^2}^2  + ( \lambda -  \frac{M}{2\nu} - M^{1/2})  \|  f \|_{L^2}^2
\\
&& +  \frac{\nu}4 C_\Omega \|  f \|_{L^{2^*}}^2 - \frac1{2\nu} \| |\beta-b| {\bf 1}_{|\beta-b| \ge M} f \|_{L^2}^2  -   \|  \sqrt{c_+} {\bf 1}_{c_+ \ge M}  f \|_{L^2}^2 
\\
&\ge&  \| f \sqrt{c_-} \|_{L^2}^2+\frac{\nu}4  \| \nabla f \|_{L^2}^2  + ( \lambda -  \frac{M}{2\nu} - M^{1/2})  \|  f \|_{L^2}^2
\\
&&  + \bigl( \frac{\nu}4 C_\Omega - \frac1{2\nu} \| |\beta-b| {\bf 1}_{|\beta-b| \ge M}  \|_{L^d}^2  
-   \|  c_+  {\bf 1}_{c_+ \ge M}  \|_{L^{d/2}}   \bigr) \|  f \|_{L^{2^*}}^2 , 
\eean
using the Sobolev inequality (with associated constant $C_\Omega)$ and the Holder inequality. Choosing $M>0$ large enough in such a way that the last term is positive, and next $\kappa_1> 0$ large enough, we deduce for instance that 
\beqn\label{eq:Stampacchia}
\langle (\lambda-\LL) f, f \rangle 
\ge \| f \sqrt{c_-} \|_{L^2}^2+ \frac{\nu}4  \| \nabla f \|_{L^2}^2  +    \|  f \|_{L^2}^2, \quad \forall \, \lambda \ge \kappa_1. 
\eeqn
%
%
Thanks to the Lax-Milgram theorem and the above coercivity estimate, we deduce  that $\lambda-\LL$ is invertible, and 
more precisely the mapping $(\lambda-\LL)^{-1} : H^{-1} \to H^1_0(\Omega)$ is well defined. We also claim that $\lambda-\LL$ enjoys a weak principle maximum, and more precisely
\beqn\label{eq:StampacchiaWeakMP}
f \in H^1_0(\Omega), \quad  (\lambda-\LL) f \ge 0 \quad\hbox{imply}\quad f \ge 0.
\eeqn
Indeed, for such a function $f \in H^1_0(\Omega)$, we take $w = f_- \in H^1_0(\Omega)$, as a test function, and elementary Sobolev space calculus together with the previous estimate yields
\bean
0 
&\le& \langle (\lambda-\LL) f, f_- \rangle = - \langle (\lambda-\LL) f_-, f_- \rangle
\\
&\le&-  \| f_- \sqrt{c_-} \|_{L^2}^2 - \frac{\nu}4  \| \nabla f_- \|_{L^2}^2  -  \|  f_- \|_{L^2}^2 \le 0,
\eean
so that $f_- = 0$ and $f \ge 0$. We thus deduce $(\lambda-\LL)^{-1} : L^2_+ \to L^2_+$,  and from J.-L.~Lions theory on parabolic equation (see for instance~\cite[Chapter~3]{MR0247243}), we next deduce that $\LL$ is the generator in $L^2$ of a positive semigroup $S_\LL$, so that \ref{H1} holds.  
It is worth emphasizing at this point that the semigroup $S$ built thanks to Lions's theory is defined by $S(t)f_0 = f$ for any $f_0 \in L^2$,  where 
$f \in \EE := C([0,\infty);L^2) \cap \Lloc^2([0,\infty);H^1_0) \cap \Hloc^1([0,\infty);H^{-1})$ is the unique (variational) solution to the equation
\beqn\label{eq:StampacchiaLions1}
(f(T),g(T))_{L^2} - (f_0,g(0))_{L^2} = \int_0^T \{\langle \partial_t g, f \rangle_{H^{-1},H^1_0}  + D_\LL(f,g) \} ds,
\eeqn
for any $T > 0$ and $g \in \EE$. 
Choosing $g = f$ in the above equation, we classically compute 
$$
\frac12\| f(t) \|_{L^2}^2 - \frac12\| f_0 \|_{L^2}^2  -   \int_0^t D_\LL(f,f) ds = 0, \quad \forall \, t >0, 
$$
which together with \eqref{eq:Stampacchia} implies 
$$
\frac{1}{ t} \int_0^t \frac\nu4 \| \nabla f \|_{L^2}^2 ds \le  - \bigl( \frac{f(t) - f_0 }{ t}, \frac{f(t) + f_0 }{ 2} \bigr)_{L^2} +  \frac{\kappa_1}{t} \int_0^t \| f \|_{L^2}^2 ds, \quad \forall \, t >0. 
 $$
When  $f_0 \in D(\LL)$,  the RHS is bounded and there thus  exists a sequence $t_n \to 0$ such that $\| \nabla f(t_n) \|_{L^2}$ is bounded. That implies $f_0 \in H^1_0(\Omega)$ and thus $D(\LL) \subset H^1_0(\Omega)$. Similarly, we may consider the dual Dirichlet form $D^*(f,g) := D_\LL(g,f)$ and build an associated positive semigroup $S^*$ through Lions's theory described above. More precisely  $S^*(t)g_0 = g$ for any $t \ge 0$ and $g_0 \in L^2$,  where $g \in \EE$ is the unique (variational) solution to the equation
$$
(g(t),f(t))_{L^2} - (g_0,f(0))_{L^2} = \int_0^t \{\langle \partial_t f,g \rangle_{H^{-1},H^1_0}  + D^*(g,f) \} ds,
$$
for any $t > 0$ and $f \in \EE$. Now, we fix $T >0$, $g_T \in L^2$ and we set  $g(t) := S^*(T-t)g_T$, so that $g$ is a solution to the backward evolution equation 
$$
- \partial_t g = \LL^* g, \quad g(T) = g_T, 
$$
with 
$$
\LL^* g :=  \partial_j(a_{ij} \partial_i g) -   \partial_i (b_i g) -  \beta_i \partial_i g +  cg.
$$
The variational formulation of this last problem is 
\beqn\label{eq:StampacchiaLions2}
(g_T,f(T))_{L^2} - (g(0),f(0))_{L^2} = \int_0^T \{ \langle \partial_t f,g \rangle_{H^{-1},H^1_0}  - D^*(g,f) \} ds,
\eeqn
for any  $f \in \EE$.  
Summing up \eqref{eq:StampacchiaLions1} and \eqref{eq:StampacchiaLions2} with $f(t)   := S(t) f_0$ for $f_0 \in L^2$ and  $g(t) := S^*(T-t)g_T$ for $g_T \in L^2$, we deduce 
$$
(S(T) f_0,g_T)_{L^2} = (S^*(T)g_T,f_0)_{L^2}.
$$
In other words, we have established that $S^* = (S_\LL)^*$ and thus that $\LL^*$ is the generator of the semigroup $S^*$.  
 
 \medskip
{\bf Condition \ref{H2}.}  Let us consider a ball $B_R$, $R > 0$, such that $B_{4R} \subset \Omega$ and next the solution 
\beqn\label{eq:StampacchiaLions2-deff0}
f_0 \in  H^1_0(\Omega), \quad (\kappa_1 - \LL ) f_0 = {\bf 1}_{B_R}, 
\eeqn
which exists from the above discussion. We next recall some classical results. On the one hand, from \cite[Sec.~3 \& Sec.~4]{MR192177} or  \cite[Thm.~8.15]{MR0473443} (see also the original papers  \cite{MR0093649,MR100158,MR170091}), 
the following  global $L^\infty$ De Gorgi-Nash-Moser type estimate 
\beqn\label{eq:DeGiorgiNMineq1}
\| f_+ \|_{L^\infty(\Omega)} \lesssim   \| f_+ \|_{L^2(\Omega)} +  \| g \|_{ L^{r/2}(\Omega)}
\eeqn
holds for any subsolution 
$$
f \in H^1_0(\Omega), \quad (\lambda -\LL) f \le g \in { L^{r/2}(\Omega)}.
$$
 The local estimate variant  \cite[Thm.~8.18]{MR0473443} (or {\it weak Harnack inequality}) 
\beqn\label{eq:DeGiorgiNMineq1bisNew}
\| f \|_{L^p(B_{2R})} \lesssim  \mathop{\inf{\vphantom{sup}}}_{B_{R}} f  +  \| g \|_{  L^{r/2}(\Omega)}, \quad \forall \, p \in [1,2^*/2), 
\eeqn \Black
also holds for a nonnegative supersolution 
$$
f \in H^1(\Omega), \quad f \ge 0 \hbox{ on } B_{4R} \subset \Omega, \quad (\lambda -\LL) f \ge g \in  L^{r/2}(\Omega), 
$$
from what  one deduces that a strong maximum principle  \cite[Thm.~8.19]{MR0473443} holds. 
More precisely, under the additional one side pointwise bound 
\beqn\label{eq:StampDivergenceHypH1bis}
  c + \Div \beta \le c_0   \quad\hbox{or}\quad  
c - \Div b \le c_0 , 
\eeqn 
for some $c_0\in \R$, we have that, for any $f \in H^1_0(\Omega)$, 
\beqn\label{eq:DeGiorgiNMineq1ter}
 \LL f \le 0 \ \hbox{in} \ \Omega,  \ f \ge 0 \ \hbox{in} \ \Omega \quad\hbox{imply}\quad
 f \equiv 0 \ \hbox{or}\ f > 0  \ \hbox{a.e. in}\ \Omega.
\eeqn  
When indeed $f \not\equiv 0$, we may choose $B_{4R} \subset \Omega$ such that $\| f \|_{L^1(B_{2R})}  > 0$ and thus   $\inf_{B_{R}} f > 0$ from \eqref{eq:DeGiorgiNMineq1bisNew} (with $g=0$) and because constants are supersolutions thanks to the first condition in \eqref{eq:StampDivergenceHypH1bis}. In the case only the second condition holds in \eqref{eq:StampDivergenceHypH1bis}, the same argument implies that $\LL^*$ satisfies the strong strong maximum principle and thus also $\LL$ thanks to Lemma~\ref{lem:Irred-KStrongMax&dual}. 
We conclude that $f$ is positive by a connectedness argument. An alternative and less demanding proof is presented in \cite[Cor.~1]{MR280858} where \eqref{eq:DeGiorgiNMineq1ter} is established without the additional assumption \eqref{eq:StampDivergenceHypH1bis}. 
%
%

On the other hand, the following H\"older regularity estimate  \cite[Théorème~7.1]{MR192177} and \cite[Thm.~8.29]{MR0473443} 
 (see also the original papers  \cite{MR0093649,MR100158,MR170091}) of De Gorgi-Nash-Moser type 
\beqn\label{eq:DeGiorgiNMineq2}
\| f \|_{C^{\alpha}(\Omega)} \le C \, \| (\lambda - \LL ) f \|_{L^\infty(\Omega)}
\eeqn
holds true for some $\alpha =\alpha(a_{ij}) \in (0,1)$ and $C >0$. 
These last two pieces of information together and the fact that $f_0 \not\equiv 0$ imply that there exists a constant $\theta > 0$ such that $f_0 \ge \theta  {\bf 1}_{B_R}$, and thus 
$$
\LL f_0  \ge (\kappa_1 - \theta^{-1}) f_0.
$$
That is condition  {\bf (i)} in Lemma~\ref{lem:Existe1-Spectre2bis}, so that condition \ref{H2}  holds thanks to Lemma~\ref{lem:Existe1-Spectre2bis}.  
Presented in that way, the   above estimate is not really constructive, but the constant $\theta := \inf_{B_R} (\kappa_1 - \LL )^{-1} {\bf 1}_{B_R} $ can also be considered as a geometric quantity associated to geometric properties of the operator and the domain. {\Blue In order to make the lower bound \ref{H2} constructive, we propose several strategies below, each one requiring additional mild assumptions about the coefficients
and none being really better than the others.}

\medskip
{\bf First constructive argument for \ref{H2}.} 
In the case when $\LL$ is self-adjoint, that corresponds to the case $a_{ij} = a_{ji}$ and $b_i + \beta_i = 0$, we classically know (that has been recalled in Section~\ref{subsect:Exist1-discussion}, see \eqref{eq:lambda1=varSA}) that 
$$
\lambda_1 = \inf_{f \in X_+ \backslash \{0 \}} \frac{\langle \LL f, f \rangle }{ \| f \|^2} 
= \inf_{f \in H_0^1, \,  \|f\|_{L^2} = 1}  \int_\OO \bigl\{  a \nabla f \cdot \nabla f + c f^2 \bigr\}, 
$$
from what and the Sobolev imbedding, we get 
$$
\lambda_1 \ge \inf_{f \in H^1_0, \, \|f\|_{L^2} = 1}   \bigl\{  (\nu C_\Omega - \| c_- {\bf 1}_{c_- \ge M} \|_{L^{d/2}}) \| f \|_{L^{2^*}}^2 -  M  \bigr\}  \ge - M, 
$$
by choosing $M$ large enough. That gives an explicit lower bound on $\lambda_1$.

\medskip 
{\bf Second constructive argument for \ref{H2}.} 
We give another constructive argument without assuming any self-adjointness property. We rather assume 
\beqn\label{eq:StampacchiaH2hyp1}
( \partial_ib_i - c)_+ \in M^1(\Omega), \quad  
b_i +  \beta_i - \partial_j a_{ij} \in M^1(\Omega).
\eeqn
We fix $h_0 \in C^2_c(\Omega)$ such that $c_0 {\bf 1}_{B_\rho} \le h_0 \le c_0 {\bf 1}_{B_{3\rho/2}}$ with $B_{8\rho} \subset \Omega$ and   $\| h_0 \|_{L^2} =1$. We next define $f_0$ as the (positive) solution to 
\beqn\label{eq:StampacchiaH2def-f0}
f_0 \in  H^1_0(\Omega), \quad (\kappa_1 - \LL ) f_0 = h_0,  
\eeqn
so that $f_0 \in C^\alpha(\Omega)$ from \eqref{eq:DeGiorgiNMineq1} and \eqref{eq:DeGiorgiNMineq2}, and similarly 
\beqn\label{eq:StampacchiaH2def-f0tilde}
\widetilde f_0 \in  H^1_0(B_{2\rho}), \quad (\kappa_1 - \LL ) \widetilde f_0 = h_0,  
\eeqn
so that $\widetilde f_0 \in C^\alpha(B_{2\rho})$ from \eqref{eq:DeGiorgiNMineq1} and \eqref{eq:DeGiorgiNMineq2}. 
We observe that $0 \le \widetilde f_0 \le f_0$ thanks to the weak maximum principle. We then compute  
$$
1 = \| h_0 \|_{L^2}^2 = \int_{B_{2\rho}} h_0  (\kappa_1 - \LL) \widetilde f_0   = \int_{B_{2\rho}} \widetilde  f_0 (\kappa_1 - \LL^*) h_0 
\le \| \widetilde f_0 \|_{L^\infty} \| (\kappa_1 - \LL^*) h_0 \|_{M^1}, 
$$
where the last term is finite because of the additional   hypothesis \eqref{eq:StampacchiaH2hyp1}. We conclude to a first constructive  lower bound  $\| \widetilde f_0 \|_{L^\infty(B_{2\rho})} \ge c_1 > 0$.
Because of the Holder continuity, we also have $\| \widetilde f_0 \|_{L^1(B_{2\rho})} \ge c_2$ with constructive constant $c_2 = c_2(c_1,\alpha,d) > 0$.  
Thanks to \eqref{eq:DeGiorgiNMineq1bisNew} (with $g=0$), we obtain
\bean
f_0 
&\ge& {\bf  1}_{B_{3\rho/2}}  \inf_{B_{3\rho/2}}   f_0 
\ge  {\bf  1}_{B_{3\rho/2}}    C_{wH} \| f_0 \|_{L^1(B_{3\rho/2})}  
\\
&\ge&
 {\bf  1}_{B_{3\rho/2}}    C_{wH} \| \widetilde f_0 \|_{L^1(B_{3\rho/2})}  
  \ge C_{wH} c_2 c_0^{-1} h_0.  
\eean
 

 Because all the inequalities are constructive and proceeding as above, we deduce that condition {\bf (ii)} in Lemma~\ref{lem:Existe1-Spectre2bis} holds and thus also  \ref{H2} with constructive constant 
 $\kappa_0 := \kappa_1 - C^{-1}_{wH} c_2^{-1} c_0 $.  Finally, because of $(\kappa_1 - \LL) f_0 = 0$ on $\Omega \backslash B_{3\rho/2}$, we may apply the Harnack inequality \cite[Cor.~8.21]{MR0473443}, and we classically deduce  there exist constructive constants $C > 0$ and $C_\varrho > 0$ for any $\varrho > 0$ such that 
\beqn\label{eq:diffuse-H2constructive-f0bis}
C_\varrho {\bf 1}_{\omega_\varrho} \le f_0 \le C, 
\eeqn
with $\omega_\varrho := \{ x \in \Omega; \, \delta(x) > \varrho\}$ and $\delta(x) := d(x,\partial\Omega)$ is the distance to the boundary function.

  \medskip

We can also get a constructive argument for  \ref{H2} by asking 
that condition {\bf (i)} in Lemma~\ref{lem:Existe1-Spectre2bis} holds. 
We may for instance verify that the dual counterpart of the above constructive argument holds
when $(c+\partial_i \beta_i)_- \in M^1$ and $ b_i + \beta_i + \partial_j a_{ji} \in M^1$. 
More precisely,  we establish in a similar way as above that   
the solution to the problem
\beqn\label{eq:diffuse-H2constructive-psi0}
\phi_0 \in  H^1_0(\Omega), \quad (\kappa_1 - \LL^* ) \phi_0 = h_0, 
\eeqn
satisfies 
\beqn\label{eq:Diffuse1-Harris-condpsi0}
\kappa_0 \phi_0 \le \LL^* \phi_0 \le \kappa_1 \phi_0, 
\eeqn
for some constructive constants $\kappa_0 \le \kappa_1$. \Black
Similarly as above again, there exist constructive constants $C > 0$ and $C_\varrho > 0$ for any $\varrho > 0$ such that 
\beqn\label{eq:diffuse-H2constructive-psi0bis}
C_\varrho {\bf 1}_{\omega_\varrho} \le \phi_0 \le C.   
\eeqn

\medskip
{\bf  Third constructive argument for \ref{H2}.} 
We write 
\beqn\label{eq:ellipticNotDivForm}
\LL f = a_{ij} \partial^2_{ij} f + \tilde b_i \partial_i f +  \tilde  cf, 
\eeqn
with $\tilde b_i := b_i + \partial_j a_{ji} + \beta_i$ and $\tilde c := c + \partial_i\beta_i$. We further assume  $\tilde b_i, \tilde c \in L^\infty$. In that case, 
we may also obtain 
an explicit lower bound on $\lambda_1$ by proceeding in the following way. 
We  define $f_0 (x) := \chi( |x|)$ with $\chi \in C^1_c(\R_+) \cap W^{2,\infty}(\R_+)$, ${\bf 1}_{[0,1/3]} \le \chi \le {\bf 1}_{[0,1]}$,  $\chi' \le 0$ on $[0,1]$, $\chi(s) := n^2(1-s)^2/2$ on $[\iota_n,1]$, $\iota_n := 1-1/(2n)$, for some $n \ge 1$ to be chosen.  As a consequence, $\chi'' = n^2$ on $[\iota_n,1]$, $|\chi'| \le n$ on $[\iota_n,1]$ and $\chi \ge 1/2$ on $[0,\iota_n]$.
 Denoting $s := |x|$, we compute 
\bean
\LL f_0 = a_{ij} \bigl\{  \chi''(s) \hat x_i \hat x_j  +  \chi'(s) \frac{\delta_{ij} - \hat x_i \hat x_j }{ s}  
 \bigr\}   +  \tilde b(x) \cdot \hat x \chi'(s) + \tilde c(x) \chi(s).
\eean
For $n$ large enough,    we get 
\bean
\LL f_0 &\ge&  n^2 \nu  - n  2 A - n B- C \ge 0
\ \hbox{ on } \  B_{1} \backslash B_{\iota_n}, 
\\
\LL f_0 &\ge& -  A \bigl\{ \| \chi '' \|_{L^\infty}  +  \| \chi'(s)/s \|_{L^\infty} \bigr\}   -  B  \| \chi'  \|_{L^\infty} - C  \ge \kappa_0 \chi
\ \hbox{ on } \ B_{\iota_n}, 
\eean\Black
with $A := \| a \|_{L^\infty(B_1)}$, $B :=  \| \tilde b \|_{L^\infty(B_1)}$, $C := \| \tilde c \|_{L^\infty(B_1)}$ and $\kappa_0 \in \R_-$. 
 As a conclusion, we  have again established 
condition {\bf (ii)} in Lemma~\ref{lem:Existe1-Spectre2bis}, so that condition \ref{H2}  holds. 

\medskip
{\bf Fourth constructive argument for \ref{H2}.} 
We present a last  situation when we are able to prove a quantitative version of condition \ref{H2}. 
We assume that $a \in C^0(\bar\Omega)$, $\Div \beta \in L^{r/2}$, as well as  $\tilde b_i \in L^r$ and $\tilde c \in L^{r/2}$ in the definition of \eqref{eq:ellipticNotDivForm}.  
We define $h_0$ and $f_0$ as in the second constructive argument for \ref{H2}, so that \eqref{eq:StampacchiaH2hyp1} holds. 
Choosing $p \in (1,2)$ defined by $1/p := 1/r+1/2> 2/r + 1/2^*$, we observe that 
\bean
\| \kappa_1 f_0 -    \tilde b_i \partial_i f_0 -  \tilde  cf_0 - h_0 \|_{L^p} 
&\lesssim& \kappa_1 \| f_0 \|_{L^2} +  \| \tilde b_i \|_{L^r}  \|   \partial_i f_0  \|_{L^2} +  \| \tilde c \|_{L^{r/2}}  \| f_0 \|_{L^{2^*}} + \| h_0 \|_{L^2}
\\
&\lesssim&   \| h_0 \|_{L^2}, 
\eean
from equation \eqref{eq:StampacchiaH2def-f0} and  the coercivity estimate \eqref{eq:Stampacchia}.
 From the Calderon-Zygmond regularity theory  \cite{MR52553} or \cite[Thm.~9.14]{MR0473443}, we also know
that 
\beqn\label{eq:diffuse:CalderpnZestim}
\| f_0 \|_{W^{2,p}(\Omega)} \lesssim \| a_{ij} \partial^2_{ij} f_0 \|_{L^p(\Omega)}.
\eeqn
Writing $ a_{ij} \partial^2_{ij} f_0 =  \kappa_1 f_0 -    \tilde b_i \partial_i f_0 -  \tilde  cf_0 - h_0$ and using the two above estimates, we deduce 
\beqn\label{eq:Zygmond}
 \| f_0 \|_{W^{2,1}(\Omega)} \lesssim \| h_0 \|_{L^2(\Omega)}. 
\eeqn
\Black
On the other hand, from \eqref{def:StampOperator} and the Poincaré inequality, we have
$$
1 = \| h_0 \|_{L^2}^2 = \langle (\kappa_1-\LL) f_0, h_0 \rangle \lesssim  \| \nabla f_0 \|_{L^2} \|  \nabla h_0 \|_{L^2}.
$$
Together with the estimate \eqref{eq:Zygmond} and the Gagliardo-Niremberg inequality 
$$
\| \nabla f \|_{L^2} \lesssim  \| D^2 f \|_{L^{1}}^{1/2}  \|  f \|_{L^\infty}^{1/2}, 
$$
we obtain a   lower bound $\| f_0 \|_{L^\infty} \ge C_0 > 0$. We then conclude as in the second constructive argument for \ref{H2}. 

%

 \medskip
{\bf Condition \ref{H3}.}
Because of Rellich-Kondrachov theorem on the compact embedding $H^1_0 \subset L^2$, the mapping  $(\lambda-\LL)^{-1} : L^2 \to L^2$ is compact for any $\lambda \ge \kappa_1$. 
As a consequence, introducing the splitting $\LL = \AA + \BB$ with $\AA := \kappa_1 - \kappa_\BB$, $\kappa_\BB \in \R$ arbitrary, 
 the operator $\RR_\BB(\lambda) = (\lambda + \kappa_1 - \kappa_\BB-\LL)^{-1}$ is  bounded uniformly on $\lambda \ge \kappa_\BB$ and it is compact for any $\lambda \ge \kappa_\BB$. 
 We deduce from  Lemma~\ref{lem:H3abstract-StrongC}-(2) that \ref{H3} holds for both the primal and the dual problems. 
 
 \medskip  
 We may thus   apply   Theorem~\ref{theo:exist1-KRexistence}  and deduce the existence of a solution   $(\lambda_1,f_1,\phi_1)$ 
to the  first eigentriplet problem 
\beqn\label{eq:ex1PLL-tripletBIS}
\lambda_1 \in \R, \quad  0 \le f_1 \in H^1_0, \ \LL f_1 = \lambda_1 f_1,  \quad 0 \le \phi_1 \in H^1_0, \ \LL^* \phi_1 = \lambda_1 \phi_1,
\eeqn
where both equations must be understood in the variational sense as a consequence of the discussion at the end of the proof of condition \ref{H1}.

 \Black

 
 \medskip
{\bf Condition \ref{H4}.}
 The strong maximum principle holds as already mentioned in the paragraph dedicated to condition \ref{H2}. 
As a consequence and thanks to Theorem~\ref{theo:KRgeometry1}, we know that the first eigentriplet problem \eqref{eq:ex1PLL-triplet} has a unique solution $(\lambda_1,f_1,\phi_1)$ which satisfies $f_1 > 0$, $\phi_1 > 0$, 
 $N( \LL-\lambda_1)^k = \hbox{\rm Span}(f_1)$ and $N( \LL^*-\lambda_1)^k = \hbox{\rm Span}(\phi_1)$ for any $k \ge 1$. 

 \medskip 
{\bf Condition \ref{H5}.} Consider $f \in D(\LL^\infty)$ such  that $0 < |f| \in D(\LL^\infty)$ and 
$$
 \LL |f|  =  \Re e (\hbox{\rm sign} f)   \LL f, 
$$
so that multiplying both term of the equation by $|f|$ and integrating, we have 
$$
\Re e  \langle  \LL f, \bar f \rangle =  \langle   \LL |f|,  |f| \rangle.
$$
 We next compute 
 $$
\Re e  \langle  \LL f, \bar f \rangle =  - \int_\Omega a_{kj} \Re e  (\partial_j f \partial_k \bar f )+
 \int_\Omega  (b_k  -   \beta_k) \Re e ( \bar f \partial_k f ) + \int_\Omega c |f|^2, 
 $$
 and  
\bean
\langle   \LL |f|,  |f| \rangle
=  - \int_\Omega a_{kj} \partial_j |f| \partial_k |f| + 
 \int_\Omega  (b_k  -   \beta_k)  \Re e ( \bar f \partial_k f ) + \int_\Omega c |f|^2,
 \eean
where in the last equality, we have used that 
$
 \partial_k |f| = \frac1{|f|} \Re e ( \bar f \partial_k f). 
$
From the three above equations, we deduce 
\bean
\int_\Omega a_{kj} [ \partial_j |f| \partial_k |f| - \Re e  (\partial_j f \partial_k \bar f )] = 0.
 \eean
 Introducing the real and complex part decomposition $f = u + iv$, and similarly as in \cite[Proof of Theorem~5.1]{MR4265692}, 
 we  next compute 
\bean
&& \partial_j |f| \partial_k |f| - \Re e  (\partial_j f \partial_k \bar f )
\\
&&= \frac1{|f|^2} \bigl[ uv (\partial_k u \partial_j v + \partial_k v \partial_j u) - u^2 \partial_j v \partial_k v - v^2 \partial_j u \partial_k u \bigr] 
\\
&&= \frac1{|f|^2} (u  \partial_j v - v  \partial_j u) (u   \partial_k v - v  \partial_k u ), 
 \eean
so that from the ellipticity condition on $a$, 
we have $u   \partial_k v - v  \partial_k u = 0$ a.e. on $\Omega$. On the other hand, 
from De Girogi-Nash-Moser regularity estimates \eqref{eq:DeGiorgiNMineq1} and \eqref{eq:DeGiorgiNMineq2}, 
 $f$ 
 has H\"older regularity. In particular both functions $u$ and $v$ are continuous. Because $|f| \not\equiv0$, 
one of the two function is not identically vanishing, say for instance $v\not\equiv 0$. There exists some points $x_0 \in \Omega$ such that $v(x_0) \not= 0$, say for instance $v(x_0) > 0$. 
Denoting by $\omega$ the connected component of the set $\{ x \in \Omega; \, v(x) > 0 \}$ containing $x_0$, we have  $\nabla (u/v) = 0$ on $\omega$. 
Hence  $u = \alpha\, v$ on $\omega$ for some $\alpha \in \R$, which implies that there exists $\sigma \in \Sp^1$ such that 
$f = \sigma |f|$ on $\omega$. 
If $\omega \not= \Omega$,  we would have $|f| = 0$ on $\partial\omega \cap \Omega \not=\emptyset$, which would be a contradiction with the fact that $|f| > 0$. 
We conclude that $\omega = \Omega$ and thus that 
$f = \sigma |f|$, which is nothing but the reverse Kato's inequality condition \ref{H5}.

\medskip 
 At this stage, we may use  Theorem~\ref{theo:KRgeometry2}, in order to get the conclusion \ref{C3} on the triviality of the principal point spectrum. 
 
\medskip
In order to go one step further and establish the asymptotic stability of $f_1$, we may use the two following approaches which are consequences respectively of Lemma~\ref{lem:Section7-BoccardoG} and Lemma~\ref{lem:Section7-GP}. 


\begin{lem}\label{lem:Section7-BoccardoG} For any $R > 0$,  the set 
$$
\KK := \{ f \in D(\LL); \, [f] \le R, \ [\LL f] \le R \}
$$
is strongly compact in $\Lloc^1(\Omega)$, where $[g] := \| g \|_{L^1_{\phi_1}}$. 
\end{lem}

\begin{proof}[Proof of Lemma~\ref{lem:Section7-BoccardoG}.]
Consider $f \in \KK$ so that $f \in H^1_0(\Omega)$ and 
$$
 \partial_i(a_{ij} \partial_j f) + b_i \partial_i f +  \partial_i(\beta_i f) +  cf = g \in L^2(\Omega).
$$
From the renormalization theory of elliptic equations and the GRE trick (see for instance \cite{MR2162224} and the references therein) for any renormalizing function $H \in C^2(\R)$, 
there holds 
\bean
 H''(u) f_1 \phi_1 a \nabla u \cdot \nabla u  &= &\hbox{div} (a \phi_1 \nabla (H(u) f_1)) -  \hbox{div} (f_1 H(u) a  \nabla \phi_1)   
\\
&&+  \, \hbox{div} ((b + \beta) H(u) f_1 \phi_1)  + g H'(u) f_1 \phi_1, 
\eean
 with $u := f/f_1$.
Considering $H \in W^{2,\infty}$ the even (and convex) function such that $H(0) = 0$ and $H'' := {\bf 1}_{[n,n+1]}$, so that in particular $|H'(s)| \le 1$, and integrating the previous equation, 
we deduce 
$$
\nu \int_{|u| \in [n,n+1]} |\nabla u|^2 f_1 \phi_1 \le \int |g| f_1 \phi_1 \le \|f_1\|_{L^\infty} R.
$$
We proceed along the line of the proof of \cite[Thm.~1]{MR1025884}. For a fixed $\omega \subset\subset \Omega$, we define $B_n := \{ x \in \omega; \, |u(x)| \in [n,n+1] \}$. 
Using that $f_1 > 0$ and $\phi_1 > 0$, there exists a constructive constant $C_{\omega,R} > 0$ such that 
$$
\int_{B_n} |\nabla u|^2 \le C^2_\omega, \quad \forall \, n \ge 0.
$$ 
From the Cauchy-Schwarz inequality, we have 
\beqn\label{BoccardoG}
\int_{B_n} |\nabla u| 
\le C_\omega {\rm meas}(B_n)^{1/2}, \quad \forall \, n \ge 0. 
\eeqn
On the other hand,   denoting by $1^* := d/(d-1)$ the Sobolev exponent, we have 
$$
\int_{B_n} |\nabla u| 
\le  C_{\omega,R} \Bigl( n^{-1^*} \int_{B_n} |u|^{1^*} \Bigr)^{1/2}.
$$
Summing up and using the   Cauchy-Schwarz inequality again, we have 
\bean
\sum_{n \ge 1} \int_{B_n} |\nabla u| 
 &\le& C_{\omega,R} \Bigl( \sum_{n \ge 1}  n^{-1^*}  \Bigr)^{1/2}  \Bigl( \sum_{n \ge 1}   \int_{B_n} |u|^{1^*} \Bigr)^{1/2}
\\
 &\le& C_{\omega,R} \Bigl( \sum_{n \ge 1}  n^{-1^*}  \Bigr)^{1/2} \| u \|_{L^{1^*}}^{1^*/2}.
\eean
Together with \eqref{BoccardoG} for $n=0$, we deduce 
$$
\| \nabla u \|_{L^1(\omega)} \le C'_{\omega,R} (1 + 
 \| \nabla u \|_{L^1(\omega)}^{1^*/2}). 
$$
Because $1^*/2 \le 3/4 < 1$ (recall that $d \ge 3$), we can kill the last term, and 
we obtain the estimate 
$$
\| \nabla (f/f_1) \|_{L^1(\omega)} \le C'' , \quad \forall \, f \in \KK, 
$$
for some constant $C'' := C''_{\omega,R} > 0$. We classically conclude thanks to the Rellich-Kondrachov theorem.
\end{proof}

 From the above lemma and Theorem~\ref{theo:ergodicity-compact-trajectories}, we deduce that $\widetilde S(t) f \to \langle f,\phi_1\rangle f_1$ in the $L^1_{\phi_1}$ norm sense as $t \to \infty$ for any $f \in L^2(\Omega)$. The alternative approach is based on the following result.

\begin{lem}\label{lem:Section7-GP} Setting $\kappa := \kappa_0 - 1$, there exist $A, \alpha, R > 0$ such that 

\smallskip
\quad {\bf (i)} $\sup_{z \in \Delta_\kappa} \langle y \rangle^\alpha \| \RR_\BB(z) \|_{\BBB(L^2;H^1_0)} + \sup_{z \in \Delta_\kappa \backslash B_R} \| \RR_\LL(z) \|_{\BBB(L^2;H^1_0)}  < \infty$, 

\smallskip
\quad  {\bf (ii)} $\Sigma(\LL) \cap \Delta_\kappa \subset \Sigma_d(\LL)  \cap B_R $,

\smallskip
where $\BB := \LL - A$ and $z = x + iy$, $x,y \in \R$. 
 \end{lem}

\begin{proof}[Proof of Lemma~\ref{lem:Section7-GP}.]
Let us  consider an a priori solution to the stationary problem
$$
f \in H^1_0, \quad z = x+iy \in \Delta_\kappa, \quad (\LL + z) f = g \in L^2. 
$$
This one satisfies 
$$
\Bigl| - \int (a \nabla f + \beta f)\cdot \nabla \bar f + \int b \cdot  \nabla  f  \bar f + (c+z) |f|^2 \Bigr|  =  \Bigl| \int g \bar f \Bigr|.
$$
Using the elliptic condition, the Cauchy-Schwarz inequality and  triangular inequalities, we get 
\bean
 \Bigl| \int g \bar f \Bigr|  
 &\ge&   
   \Bigl| \int a \nabla f \nabla \bar f + ((c+x)_+ + i y) |f|^2 \Bigr| - \Bigl| \int  b \cdot  \nabla  f  \bar f
   -   \beta f\cdot \nabla \bar f
   + (c+x)_- |f|^2 \Bigr|  
\\
 &\ge& \frac\nu2 \| \nabla f \|_{L^2}^2 + \bigl( \frac{|y|}2 - x_- \bigr)  \|  f \|_{L^2}^2  - \| (|b|+|\beta|) f \|_{L^2} \| \nabla f \|_{L^2} - \| \sqrt{c_-} f \|_{L^2}^2. 
 \eean
Using next similar arguments and those introduced in the paragraph dedicated to condition \ref{H1} and with similar definition for the constant $M := M(b,\beta,c) > 0$, 
we deduce 
\bean
 \Bigl| \int g \bar f \Bigr|  
 &\ge&   \bigl( \frac{|y|}2  - x_- -  M    \bigr) \| f \|_{L^2}^2 + \frac\nu4  \| \nabla f \|_{L^2}^2 . 
 \eean
Defining the sectorial set 
 $$
 \SS := \bigl\{ z = x+i y \in \C; \ |y| > 2 x_- + M \bigr\},
 $$
 we have established the a priori estimates
\bean
\| f \|_{L^2} &\le& \bigl( \frac{|y|}2  - x_- -  M \bigr)^{-1/2} \| g \|_{L^2}, 
\\
\| \nabla f \|_{L^2} &\le& 2 \nu^{-1/2}\bigl( \frac{|y|}2  - x_- -  M \bigr)^{-1/4} \| g \|_{L^2}, 
\eean
for any $z \in \SS$. 
%
We classically and immediately deduce that $\rho(\LL) \supset \SS$ and the resolvent estimate $\| \RR_\LL(z) \|_{\BBB(L^2,H^1_0)} \lesssim  \bigl( \frac{|y|}2  - x_- - M \bigr)^{-1/2}
+ \bigl( \frac{|y|}2  - x_- - M \bigr)^{-1/4}$ for
any $z \in \SS$, and in particular  the estimate {\bf (i)} holds true. 

\smallskip
 On the other hand,  because $\LL$ has compact resolvent as established just above or during the proof of \ref{H3} and using the Fredholm alternative, we have   $\Sigma(\LL) = \Sigma_d(\LL)$ and $\Sigma(\LL) \cap \Delta_\kappa$ is finite for any $\kappa \in \R$, what is nothing but the property  {\bf (ii)}. 
 \end{proof}
 
 \medskip
 
  \medskip 
 From the above lemma and Theorem~\ref{theo:GearhartPruss} or Theorem~\ref{theo:MiscScher}, we deduce that $\widetilde S(t) f \to \langle f,\phi_1\rangle f_1$ in the $L^2$ norm sense as $t \to \infty$ for any $f \in L^2(\Omega)$ with exponential rate.

%
%
%
%
%
%
%
%
%
%
%
%

\medskip
We may summarize our analysis in the following result.

\begin{theo}\label{theo:stampacchia1} Consider the elliptic operator \eqref{eq:exPLL1} in a bounded domain and assume that the coefficients satisfy \eqref{eq:StampHypaij}, \eqref{eq:StampDivergenceHypH1}
and \eqref{eq:StampacchiaH2hyp1}. 
 Then the conclusions {\Blue \ref{S1}, \ref{S2} and \ref{S33}  hold in $L^2$} as well as  \ref{E2} in $L^1_{\phi_1}$ norm and \ref{E31} in $L^2$ with non constructive rate. 
\end{theo}

{\Blue 
Under not optimal assumptions, the above result belongs to folklore. With the same kind of regularity assumptions the conclusion \ref{S1}-\ref{S2} is probably due to Chicco \cite{MR280858,MR0310462} and to Lions  \cite[8th  course]{PL2} for the constructive estimates. We believe that the ergodicity part \ref{E2} and  the convergence part \ref{E31} are  new in this framework.
It is however worth emphasizing again that the above approach for the long time asymptotic behaviour is definitively  not constructive.} We propose now an alternative approach which is constructive.

 \medskip
{\bf Quantitative estimate of stability.} 
 
Using the Doblin-Harris type approach presented in Section~\ref{sec:QuantitativeStabilityKR}, we are able to 
establish a rate of convergence to the principal dynamic, at least in a regular framework. We thus make some regularity assumptions on $\Omega$ and additional regularity assumptions on the coefficients. 

- For the domain, we assume that  there exists a constant $r_\Omega >0$ such that for any $x\in \Omega$ there is $y \in \Omega$ such that $x\in B(y, r_\Omega )\subset\Omega$, in particular,  for any $x\in \partial\Omega$ there is $y \in \Omega$ such that $x\in \partial B(y, r_\Omega )$, $B(y, r_\Omega ) \subset\Omega$. We also assume that $\Omega$ is $C^{1,1}$. 

- For the coefficients, we assume $a_{ij} \in C(\bar\Omega)$, $\tilde b_i, \tilde c \in L^\infty(\Omega)$, where $\tilde b_i$ and $\tilde c$ are defined in \eqref{eq:ellipticNotDivForm}.

\begin{theo}\label{theo:stampacchia2} Consider the elliptic operator \eqref{eq:exPLL1} in a bounded domain and assume that the assumptions of Theorem~\ref{theo:stampacchia1} hold together with 
the above additional regularity assumptions on the coefficients and the boundary. 
Then the conclusion  \ref{E31} holds with  constructive exponential rate. 
\end{theo}

  {\Blue 
  We are not aware of any such a result, apart from the recent work of Champagnat and Villemonais \cite{MR4066299,Champagnat2023}, which demonstrates a similar result using a probabilistic method (also of  Doblin-Harris type) under  Holder regularity conditions. 
Although the constants involved in Theorem~\ref{theo:stampacchia2} are constructive, we do not attempt to make them explicit:  because of the multiple steps required to establish the final estimate, tracing their dependency on the data would be tedious and probably of little interest. In particular, we do not claim that our method allows us to obtain optimal constants. 
}

\medskip
  The proof of Theorem~\ref{theo:stampacchia2} follows from Theorem~\ref{theo:Harris}. We split the proof into several steps.  
 
\smallskip

{\bf - Step 1. Regularity estimates.} Thanks to De Giorgi-Nash-Moser regularity technique for parabolic equations developed for instance in \cite{MR0181837} (in Russian), \cite[Thm. 1.3, Thm. 2.2]{MR226168} as well as  more recently in \cite[Lem.~2.7]{MR2187812} and \cite[Thm.~1.1]{JGuerand}, there exists $\alpha = \alpha(a_{ij}) \in (0,1)$ and for any $T_1 > T_0 > 0$ and any $\varrho\in (0,1)$, there exist constructive constants $C_i = C_i(\| f \|_{L^\infty_tL^2_x},T,\tau,r)$ such that any solution $f \in L^\infty(0,\infty;L^2(\Omega))$ to the parabolic equation 
$
\partial_t f = \LL f
$
satisfies 
\beqn\label{eq:DeGirogiParabolic}
\| f \|_{L^\infty([T_0,T_1] \times \Omega)} \le C_1, \quad
\| f \|_{C^\alpha([T_0,T_1] \times \omega_\varrho)} \le C_2,
\eeqn
with $\omega_r := \{ x \in \Omega; \, d(x,\partial\Omega) > r \}$. 
More precisely, in order to establish the second estimate in \eqref{eq:DeGirogiParabolic} with constructive constant, one may observe that  the proof of  \cite[Prop.~2.4]{JGuerand} may be repeated in order to 
get that solutions to the parabolic equation considered in the present framework fall into De Giorgi classes as defined in \cite[Definition~2.3]{JGuerand},
and thus \cite[Thm.~1.1]{JGuerand} applies. 

\smallskip 
On the other hand, in this context and because of the regularity assumptions, we may establish a more accurate regularity estimate. 
More precisely, by gathering the Sobolev inequality and the Calderon-Zygmond estimate \eqref{eq:diffuse:CalderpnZestim},  we obtain the classical constructive regularity estimate
\beqn\label{eq:diffuse1-H6-C01estim}
\| u \|_{C^{0,1}(\Omega)} \lesssim \| u \|_{W^{2,d+1}(\Omega)} \lesssim \| (\kappa_1-\LL) u \|_{L^{d+1}(\Omega)}, 
\eeqn
see for instance  Theorem~7.10, Theorem~7.25  and  Lemma~9.17 in \cite{MR0473443}. 
Iterating the same kind of arguments, we get 
 \beqn\label{eq:diffuse1-IterateCZestim}
\| u \|_{C^{0,1}(\Omega)} \le C \| (\kappa_1-\LL)^k u \|_{L^{2}(\Omega)}, 
\eeqn
with constructive constants $C$ and $k$.  
\Black



\smallskip\smallskip
  {\bf - Step 2.  Harnack estimate}. We claim that for any $T >  t_0 > 0$ and $ \varrho>0$, there exist a constant $C_H>0$ such that, for any $f_0 \in L^2$, the associated solution $f := S_\LL f_0$ satisfies 
\begin{equation}\label{eq:diff_Harnack}
    \sup_{\omega_\varrho} f_{t_0}\leq C_H\inf_{\omega_\varrho} f_T. 
\end{equation}
The proof mainly follows form   Aronson-Serrin \cite{MR244638} (see also \cite{MR159139,MR0183967,MR0196285,MR0171086,MR0232096,MR226198,MR226168,MR0229954} for similar results).
First, we know from  \cite[Thm.~3]{MR244638} that 
\begin{equation}\label{eq:diff_Serrin}
 \max_{Q^*(\rho)}  f \le C \min_{Q(\rho)} f,   
\end{equation}
for any $\rho >0$, $t > 0$ such that $Q^\ast(3\rho) \subset (0,\infty) \times \Omega$, where $Q(\rho) := [t-\rho^2,t] \times \CC(\rho)$, $Q^*(\rho) := [t-8\rho^2,t-7\rho^2] \times \CC(\rho)$
and $\CC(\rho)$ is a cube with length $\rho$. To avoid technical issues we assume that $w_\varrho$ is convex. In other case, the geometrical condition given above implies that there is $N\in \Z_+$ such that any two points $x,y\in \Omega$ can be connected by a polygonal path of at most $N$ segments, and we can argue as follows for any segment. 
 We define $D:=\sup_{a,b\in\Omega}d(a,b)$ the diameter of $\Omega$ and we choose  $r^\prime<\varrho/7$ such that 
\begin{equation*}
   7(\lfloor \frac{D}{2r^\prime} \rfloor+1)(r')^2<T-t_0.
\end{equation*}
For any $x,y\in \omega_\varrho$, we also define $N_c=\lfloor \frac{|x-y|}{r^\prime} \rfloor$.
Since $\omega_\varrho$ is convex,  $r^\prime<\varrho/7$, we have that the family of cubes $\{ \CC(x_i,2r')\}_{i=0,N_c}$ of center $x_i$ and length $2r'$ 
  for $x_i=x+\frac{(x-y)i}{N_c}$ satisfy that
$C(x_i,6r')\subset \Omega$ and $\CC(x_i,2r')\cap C(x_{i+1},2r')\neq \emptyset$  for any $i=0, \dots, N_c$. 
As a consequence, we can apply Aronson-Serrin estimate \eqref{eq:diff_Serrin} for each cube to obtain
\begin{equation*}
    \max_{C(x_i,2r')} f_{t_i}\leq C_{2r'}\min_{C(x_{i},2r')} f_{t_{i+1}},
\end{equation*}
with $t_i=t_0+7 i(2r')^2$. Taking $y_i\in C(x_i,2r')\cap C(x_{i+1},2r')$, we deduce
\begin{equation*}
    \max_{C(x_i,2r')} f_{t_i}\leq C_{2r'}\min_{C(x_i,2r')} f_{t_{i+1}}\leq C_{2r'}f_{t_{i+1}}(y_i)\leq C_{2r'}\max_{C(x_{i+1},2r')} f_{t_{i+1}}\leq C_{2r'}^2\min_{C(x_{i+1},2r')} f_{t_{i+2}}.
\end{equation*} 
By induction, we obtain
\begin{equation*}
    f_{t_0}(x)\leq\max_{C(x_1,2r')} f_{t_1}\leq C_{2r'}^{N_c}\min_{C(x_{N_c},2r')} f_{t_{N_c}}\leq C_{2r'}^{N_c} f_{t_{N_c}}(y),
\end{equation*}
with $t_{N_c}=t_0+7N_c r'^2\leq T.$ 
Note that in any case the constant $C_{2r'}$ is the same since it only depends on the length $2r'$ and the coefficient of the equation. We have thus established  \eqref{eq:diff_Harnack} 
with $C_H:=C_{2r'}^{\lfloor \frac{D}{2r^\prime} \rfloor+1}$.

\smallskip
On the other hand, we state an improved version of the already mentioned stationary Harnack inequality. 
Because of the interior ball condition the Hopf Lemma (see for instance the proof of  \cite[Lem.~3.4]{MR0473443}) claims that for any $\varrho \in (0, r_\Omega/2]$ there exists a constructive constant $\alpha > 0$ such that if  $u \in W^{2,p}(\Omega)$, $p > d$,  is such that 
$$
u \ge {\bf 1} _{ \omega_{\varrho} }, \quad (\kappa_1 - \LL^*) u \ge 0, 
$$ 
then $u$ satisfies
\beqn\label{eq:diffuse1-H6-Hopf}
u \ge \chi(x) := e^{-\alpha (2\varrho - \delta(x))^2} - e^{-\alpha (2\varrho)^2} \hbox{ on } \omega^c_\varrho.
\eeqn

Let us give two applications of the above sharp regularity and positivity estimates. 
First, recalling \eqref{eq:diffuse-H2constructive-psi0bis} and using \eqref{eq:diffuse1-H6-C01estim} and \eqref{eq:diffuse1-H6-Hopf}, 
we deduce that there exist two constructive constants $c_i \in (0,\infty)$ such that 
\beqn\label{eq:diffuse1-H6-psi0&delta}
c_0 \delta \le \phi_0 \le c_1 \delta \ \hbox{ on } \ \Omega.
\eeqn

Consider now $f_1 \in H^1_0(\Omega)$ the positive first eigenfunction with normalization $\| f_1 \|_{L^2} = 1$. Using the  estimate of regularity \eqref{eq:diffuse1-IterateCZestim} on the 
iterated equation $(\kappa_1 - \LL)^k f_1 = (\kappa_1-\lambda_1)^k f_1$, we have 
$$
\| f_1 \|_{L^\infty(\Omega)} \le \| f_1 \|_{C^{0,1}(\Omega)} \le C_1, 
$$
for some constructive constant $C_1 \in (0,\infty)$. Next using the elementary inequality 
$$
1 = \int_\Omega f_1^2 \le \| f_1 \|_{L^\infty}  \| f_1 \|_{L^1} \le C_1  \| f_1 \|_{L^1},
$$
we deduce 
\bean
|\Omega| \sup_{\omega_\varrho} f_1 
&\ge& \int_{\omega_\varrho} f_1  =  \|f_1\|_{L^1}  -  \int_{\omega^c_\varrho} f_1  
\\
&\ge& 1/C_1 - C_1 |\omega_\varrho^c| \ge 1/(2C_1),  
\eean
by choosing $\varrho \in (0,r_\Omega/2)$ small enough. Then, from the   Harnack inequality \cite[Cor.~8.21]{MR0473443}, we deduce 
$$
\inf_{\omega_\varrho} f_1 \ge C_H \sup_{\omega_\varrho} f_1 \ge C_H (2C_1|\Omega|)^{-1}.
$$
Finally, from the above Hopf lemma and the above Lipschitz continuity, we have established 
\beqn\label{eq:diffuse1-H6-f1&delta}
c_0 \delta \le f_1 \le c_1 \delta \ \hbox{ on } \ \Omega,
\eeqn
for two constructive constants $c_i \in (0,\infty)$. The same arguments on the normalized and positive first dual eigenfunction $\phi_1$ lead to the same estimate 
\beqn\label{eq:diffuse1-H6-phi1&delta}
c_0 \delta \le \phi_1 \le  c_1 \delta \ \hbox{ on } \ \Omega.
\eeqn
In particular, for any such $\varrho \in (0,r_\Omega/2)$, we have 
\beqn\label{eq:diffuse1-hyp-Harris-BdBis}
  \langle \phi_1, {\bf 1}_{\omega_\varrho} \rangle \ge r_\varrho,  
\eeqn
with constructive constant $r_\varrho$,  what is nothing but condition \eqref{eq:hyp-Harris-BdBis} in the Doblin-Harris theorem that we will use below.

\smallskip
\Black
{\bf - Step 3. Splitting of $\LL$.}
%
%
%
%
%
We introduce the splitting $\LL = \AA + \BB$, with $\AA f = \MMM {\bf 1}_{\omega_\varrho} f$, $\MMM \ge 0$ large enough and $\varrho > 0$ small enough that we fix just below. Using \eqref{eq:Stampacchia}, we observe that 
\bean
( \BB f , f)_{L^2} 
&=& (\LL f, f)_{L^2} - \MMM \| f \|_{L^2} + \MMM \int_{\omega^c_\varrho} f^2
\\
&\le& - \frac\nu4 \| \nabla f \|^2_{L^2}  + (\kappa_1 - \MMM)  \| f \|^2_{L^2} + \MMM |\omega^c_\varrho|^{4/d} \| f \|_{L^{2^*}}^2 \le \kappa_0 \| f \|^2_{L^2}, 
\eean
by choosing first $\MMM \ge \kappa_1 - \kappa_0$ and next $\varrho > 0$ small enough in order to be able to throw away the last term using the negative first term and the Sobolev inequality.
We  deduce 
\beqn\label{eq:diffuse-Omega-SBL2L2}
S_\BB(t) : L^2 \to L^2 \ \hbox{ with bound } \ \OO(e^{\kappa_0 t}). 
\eeqn
On the other hand, denoting $f_t := S_\LL(t) f$ for $f \in L^2(\Omega)$ and recalling that $\phi_0$ defined by \eqref{eq:diffuse-H2constructive-psi0} satisfies  \eqref{eq:Diffuse1-Harris-condpsi0},  
we have 
$$
\frac{d }{ dt} \int |f_t| \phi_0 \le \int  \LL |f_t| \phi_0 \le \int    |f_t| \LL^* \phi_0 \le \kappa_1  \int |f_t|   \phi_0 , 
$$
so that 
\beqn\label{eq:diffuse-Omega-SCL1psiddt}
 \int |f_t| \phi_0 \le e^{\kappa_1 t}   \int |f_0|   \phi_0. 
\eeqn
Arguing in the same way for $S_\BB$ and using \eqref{eq:diffuse1-H6-psi0&delta}, we have established 
\beqn\label{eq:diffuse-Omega-SCL1psiL1psi}
S_\LL(t), S_\BB(t) : L^1_{\delta} \to L^1_{\delta}  \ \hbox{ with bound } \ \OO(e^{\kappa_1 t}). 
\eeqn
For a  solution to the evolution equation $\partial_t f = \CC f$, $\CC = \LL$ or $\CC = \BB$, we  also classically compute 
\bean
\frac{d }{ dt} \int f^2 \phi_0 
&=& 2 \int (\CC f) f \phi_0
\\
&=& - 2 \int (\nabla f \cdot a \nabla f )   \phi_0 + \int f^2 \CC^* \phi_0.
\eean
Thanks to \eqref{eq:Diffuse1-Harris-condpsi0} again, we have
\beqn\label{eq:diffuse-Omega-SCL2psi0ddt}
\frac{d }{ dt} \int f^2 \phi_0  \le - 2 \nu \int |\nabla f|^2  \phi_0 + \kappa_1  \int f^2 \phi_0, 
\eeqn
from what we deduce 
$$
S_\LL(t), S_\BB(t) : L^2(\delta) \to L^2(\delta)  \ \hbox{ with bound } \ \OO(e^{\kappa_1 t/2}). 
$$
 
In the sequel, we will  need the following version of Nash inequality.

\begin{lem}[weighted  Nash inequality] \label{lem:weight-delta-Nash}
There exists a constructive constant $C_N$ such that 
\beqn\label{eq:Nash-Inequality}
\| f \|_{L^2(\delta)} \le C_N  \| \nabla f \|_{L^2(\delta)}^{\frac{d+1 }{ d+2}}   \| f \|_{L^1_\delta}^{\frac{1}{d+2}}, 
\quad \forall \, f \in H^1(\delta).
\eeqn
\end{lem}

\begin{proof}[Proof of Lemma~\ref{lem:weight-delta-Nash}.]
For $\eps >0$, we define
$$
 f_\eps(x) := \frac{1 }{ \delta_\eps(x)}  \int_{B(x,\eps)} f(y) \, \delta(y) dy,   \quad  \delta_\eps(x) = \delta(  B(x,\eps)) := \int_{  B(x,\eps)} \delta(y) dy,
 $$
 and $  B(x,\eps) := \{ y \in \Omega;  \,   |x-y| < \eps\}$. It is worth emphasizing that 
\beqn\label{eq:estim-delta-r}
 \eps^{d+1} \lesssim \delta_\eps(x) \lesssim \eps^d, \quad \forall \, \eps > 0.
\eeqn
For $f \in H^1(\delta)$, we compute 
\bean
 \| f -f_\eps \|_{L^2(\delta)}^2
 &=&  \int_{\Omega} \Bigl| \frac{1 }{ \delta_\eps(x)} \int_{ B(x,\eps)}  (f(y) - f(x))  \, \delta(y) dy \Bigr|^2   \delta(x) dx 
 \\
 &\le&
\int_{\Omega}\int_{\Omega} \, {\bf 1}_{|y-x| \le \eps}  |f(y) - f(x)|^2 \, \frac{\delta(y) }{ \delta_\eps(x)}  \delta(x) dx dy
\\
&\le&
  \eps^2 \, \int_0^{1/2}\int_{\Omega}\int_{\Omega} \,    |\nabla f ((1-t) x + t y)|^2  \, \frac{\delta(y) }{ \delta_\eps(x)}  \delta(x) dx dy dt 
\\
&&+  \eps^2 \, \int_{1/2}^1\int_{\Omega}\int_{\Omega}   |\nabla f ((1-t) x + t y)|^2  \, \frac{\delta(y) }{ \delta_\eps(x)}  \delta(x) dx dy dt 
\\
&\lesssim&
  \eps^2 \, \int_0^{1/2}\int_{\Omega}\int_{\Omega}    |\nabla f (z)|^2  \,   \frac{\delta(y) }{ \eps^{d+1}} dy2 \delta(z) \frac{dz}{ (1-t)^d} dt 
\\ 
&& +  \eps^2 \, \int_{1/2}^1\int_{\Omega}\int_{\Omega} |\nabla f (z)|^2  \,    \frac{2 \delta(z) }{ \eps^{d+1}}   \delta(x)  \frac{ dz}{ t^d} dt dx 
\eean
where for the last inequality we have used the first inequality in \eqref{eq:estim-delta-r}, the fact that $\delta(x) \le 2 \delta(z)$ when $0 < t < 1/2$ and 
 the fact that $\delta(y) \le 2 \delta(z)$ when $1/2 < t < 1$. Using the second inequality in \eqref{eq:estim-delta-r}, we straightforwardly obtain
$$
\| f - f_\eps \|_{L^2(\delta)}^2 \le C_1 \, \eps\, \| \nabla f \|^2_{L^2(\delta)}, \quad \forall \, \eps>0,
$$
for a constant $C_1 > 0$.
%
On the other hand, we also observe that 
$$
 \| f_\eps \|_{L^\infty} \le \frac{C_2 }{ \eps^{d+1}} \, \| f \|_{L^1_\delta}.
$$
Writing now
$$
f^2 = f(f-f_\eps) + f f_\eps
$$
and using the above two estimates, we deduce
\bean
\| f \|_{L^2_\delta}^2 
&\le& \| f \|_{L^2_\delta} \, \| f -f_\eps \|_{L^2_\delta} + \| f \|_{L^1_\delta} \, \| f_\eps \|_{L^\infty}   .
\\
&\le&  \| f \|_{L^2_\delta} \, C_1 \, \eps^{1/2} \| \nabla f \|_{L^2_\delta} +  C_2 \, \eps^{-d-1} \, \| f \|_{L^1_\delta}^2
\\
&\le&  \frac{1 }{ 2} \| f \|_{L^2_\delta}^2  + \frac{C_1}{ 2}   \, \eps \,  \| \nabla f \|_{L^2_\delta}^2 +  C_2 \, \eps^{-d-1} \, \| f \|_{L^1_\delta}^2, 
\eean
and we obtain the weighted  Nash inequality \eqref{eq:Nash-Inequality} by choosing $\eps := ( \| f \|_{L^1_\delta}^2/   \| \nabla f \|_{L^2_\delta}^2)^{1/(d+2)}$. 
\end{proof}

Defining 
$$
u := \int |f_t| \phi_0 dx e^{-2\kappa t}, \quad v := \int  f_t^2 \phi_0 dx e^{-2\kappa t}, 
$$
with $\kappa := \kappa_{1+}$, coming back to \eqref{eq:diffuse-Omega-SCL2psi0ddt} and using \eqref{eq:diffuse1-H6-psi0&delta}, the Nash inequality \eqref{eq:Nash-Inequality}  and the estimate     \eqref{eq:diffuse-Omega-SCL1psiddt}, we get
\bean
v'(t) 
&\le&   - 2 \nu c_0 \int |\nabla f_t|^2  \delta e^{-2\kappa t}
\\
&\le&   - 2 \nu c_0 C_N^{-2 {\tfrac{d+2}{d+1}}} 
\frac{ \Bigl(  \| f_t \|_{L^2(\delta)}^2  e^{-2\kappa t} \Bigr)^{\tfrac{d+2}{d+1}} }{ \Bigl( \| f_t \|^2_{L^1(\delta)} e^{-2\kappa t} \Bigr)^{\tfrac{1}{d+1}} }  
\\
&\le&   - C
\frac{   v(t)^{1+\alpha} }{  u(0)^{2\alpha}}   ,
\eean
with $C := 2 \nu C_N^{-2 {\tfrac{d+2}{d+1}}} c_0^{1 + 2 \tfrac{d+2}{d+1}} c_1^{- \tfrac{2}{d+1}} $ and $\alpha := 1/(d+1)$.  Integrating in time, we deduce
\bean
v(t)  &\le&   \frac{\alpha^{1/\alpha} }{ C^{1/\alpha}} \frac{u(0)^2}{ t^{1/\alpha}}, 
\quad \forall \, t > 0.
\eean
We have thus established that there exist constructive constants $K > 0$ and $\kappa \ge 0$ such that 
 \beqn\label{eq:Nash-Inequality-phi0}
\| S_\CC(t)  f\|_{L^2(\phi_0)} \le K \frac{e^{\kappa t} }{ t^{(d+1)/2}}   \|  f \|_{L^1(\phi_0)}, 
\quad \forall \, f \in L^1(\phi_0).
\eeqn
From that last result, the estimates \eqref{eq:diffuse1-H6-psi0&delta} and the properties of $\AA$, we deduce that for $N \ge 1$ large enough 
\beqn\label{eq:diffuse-Omega-SAN-L1psiL2psi}
(S_\BB \AA)^{(*N)} : L^1(\delta) \to L^2(\delta)  \ \hbox{ with bound } \ \OO(e^{\kappa t}). 
\eeqn
We refer to \cite[Prop.~3.9]{MR3779780}, \cite[Prop. 2.5]{MR3465438} and \cite[Lem.~2.4]{MR4265692} for details.

%

\smallskip
{\bf - Step 4. Lyapunov condition.}
We may next write 
$$
\widetilde S_\LL = V + W * \widetilde S_\LL , 
$$
with 
$$
V :=  \widetilde S_\BB + \dots + (\widetilde S_\BB \AA)^{(*(N-1))}, \quad V :=   (\widetilde S_\BB \AA)^{(*N)} . 
$$
On the one hand, using that $\AA : L^2 \to L^2$ is bounded and \eqref{eq:diffuse-Omega-SBL2L2}, we deduce that 
$$
V: L^2 \to L^2, \ \hbox{ with bound } \ \OO(e^{\kappa t}),  
$$
for any $\kappa \in (\kappa_0-\kappa_1, 0)$. On the other hand, using that   $\AA : L^2_\delta \to L^2$   is bounded as well as 
\eqref{eq:diffuse-Omega-SCL1psiL1psi} for $S_\LL$, \eqref{eq:diffuse-Omega-SAN-L1psiL2psi}, \eqref{eq:diffuse1-H6-phi1&delta}, \eqref{eq:diffuse1-H6-psi0&delta} and  \eqref{eq:diffuse-Omega-SBL2L2}, 
we deduce that 
$$
 W * \widetilde S_\LL : L^1_{\phi_1} \to L^2, \ \hbox{ with bound } \ \OO(e^{\kappa' t}),  
$$
for any $\kappa'  > \kappa_1 - \kappa_0$. We may thus fix $t = T$ large enough such that the following Lyapunov inequality  holds
\beqn\label{eq:diffuse1-stabilityKR-Lyapunov}
\| \widetilde S_\LL (T) f \|_{L^2} \le \frac12 \| f \|_{L^2} + M_T \| f \|_{L^1_{\phi_1}},
\eeqn
which is nothing but \eqref{eq:stabilityKR-Lyapunov} in the hypothesis of the Doblin-Harris theorem.

\smallskip
{\bf - Step 5. Doblin-Harris condition} Let $A >0$ and consider $0 \le f_0 \in L^2$ such that $\|f_0\|_2\leq A\langle f_0,\phi_0\rangle$. 
 We set $\widetilde f_t := e^{-\lambda_1 t} S_\LL(t) f_0$. 
From the first inequality in \eqref{eq:Diffuse1-Harris-condpsi0}, we have 
\begin{equation*}
    \frac{d}{dt}\langle \widetilde f_t,\phi_0\rangle=\langle \widetilde f_t,(\LL^\ast-\lambda_1)\phi_0\rangle\geq -(\lambda_1-\kappa_0)\langle \widetilde f_t,\phi_0\rangle,
\end{equation*}
and then, thanks to Gronwall lemma again, we obtain,
\begin{equation*}
    \langle \widetilde f_t,\phi_0\rangle\geq e^{-(\lambda_1-\kappa_0)t}\langle f_0,\phi_0\rangle.
\end{equation*}
This estimate, together with the previous step, shows that
\begin{align*}
\int_{\omega_\varrho} \widetilde f_{t_0}(x)\phi_0 dx 
      &= \int_{\Omega} \widetilde f_{t_0}(x)\phi_0dx-\int_{\omega_\varrho^c} \widetilde f_{t_0}(x)\phi_0dx
 \\  
      &\geq  e^{-(\lambda_1-\kappa_0){t_0}}\langle  f_0,\phi_0\rangle-\|\widetilde f_{t_0}\|_2\|\phi_0\|_\infty|\omega_\varrho^c|^{1/2}
 \\  
      &  \geq e^{-(\kappa_1-\kappa_0){t_0}}\langle  f_0,\phi_0\rangle -  e^{ (\kappa_1- \kappa_0) {t_0}} \| f_0\|_2\|\phi_0\|_\infty|\omega_\varrho^c|^{1/2}
 \\  
      &  \geq \Bigl( e^{-(\kappa_1-\kappa_0){t_0}} -  A e^{ (\kappa_1- \kappa_0) {t_0}}   \|\phi_0\|_\infty|\omega_\varrho^c|^{1/2} \Bigr) \langle  f_0,\phi_0\rangle . 
\end{align*}
Choosing $\varrho > 0$ small enough, we get 
\begin{align*}
\int_{\omega_\varrho} \widetilde f_{{t_0}}(x)\phi_0 dx 
    &  \geq  \gamma \langle  f_0,\phi_0\rangle, \quad
    \gamma :=  \frac12 e^{-(\lambda_1-\kappa_0){t_0}}. 
\end{align*}
As a consequence, there is $x_{t_0}^f\in \omega_\varrho$ such that 
\begin{equation*}
    \widetilde f_{t_0}(x_{t_0}^f)
    \geq \frac{1}{|\omega_\varrho|}\int_{\omega_\varrho} \widetilde f_{t_0}(x)dx
    \geq \frac{1}{|\Omega|c_1  \varrho}\int_{\omega_\varrho} \widetilde f_{t_0}(x)\phi_0dx
    \geq \frac{\gamma}{|\Omega|c_1\varrho} \langle  f_0,\phi_0\rangle .
\end{equation*}
On the other hand, from the Harnack inequality \eqref{eq:diff_Harnack} established in Step 2, we know that for any $T > t_0$, there  exits $C_H$ such that
 \begin{equation*}
     \widetilde f_{t_0}(x_{t_0}^f)\leq \sup_{\omega_\varrho} \widetilde f_{t_0}\leq C_H\inf_{\omega_\varrho}\widetilde f_{T}.
 \end{equation*}
The two last estimates together with 
\eqref{eq:diffuse1-H6-phi1&delta} and  \eqref{eq:diffuse1-H6-psi0&delta}   imply the Doblin-Harris type estimate 
\beqn\label{eq:diffuse1-eq:hyp-Harris}
     \widetilde f_T=\widetilde S(T)f_0 \geq 
      g_A \langle  f_0,\phi_1 \rangle,
\eeqn
 with $g_A := \frac{c_0 \gamma}{C_H|\Omega|c^2_1\varrho} \textbf{1}_{w_\varrho}$,
 which is nothing but \eqref{eq:hyp-Harris} in Doblin-Harris theorem. 

\smallskip
{\bf - Step 6. Conclusion.} 
Because of the constructive estimates \eqref{eq:diffuse1-hyp-Harris-BdBis}, \eqref{eq:diffuse1-stabilityKR-Lyapunov} and \eqref{eq:diffuse1-eq:hyp-Harris}, we may apply the Doblin-Harris type Theorem~\ref{theo:Harris}, 
and we conclude to the  exponential stability \ref{E31intro}
 in   the norm of $L^2(\Omega)$ with constructive constants.

 \medskip
\subsection{Diffusion in $\R^d$ with strong potential confinement}  
\label{subsection:diffusionStrongPotential}
We consider in this section  the  elliptic operator
\beqn\label{eq:LLf=Rd}
\LL f := \Delta f + b \cdot \nabla f   + c f, \quad f \in H^1(\R^d),
\eeqn
with $b  \in \Lloc^\infty(\R^d)$, $c \in  \Lloc^2(\R^d)$ and a confinement condition that we impose through the properties of the {\it potential function} $c$, which is roughly speaking $c \to - \infty$ as $|x| \to \infty$. 
More precisely, we 
assume
\beqn\label{eq:cTO-inftyBIS}
\sigma_{i+} \in L^{d/2}, \quad \hbox{\rm meas}\{  \sigma_i \ge  K \} < \infty, \   \forall \, K < 0,
\eeqn
with either   $\sigma_1 := c + |b|^2/\kappa$  for some constant  $\kappa \in (0,4)$ \Black
or either $\sigma_2 :=  c + \Div b/2$. When we assume
that 
$$
c \sim - |x|^\gamma \quad\hbox{and}\quad  b \sim x |x|^{\beta-1} \quad\hbox{as}\quad |x| \to \infty, 
$$
the condition \eqref{eq:cTO-inftyBIS} for $\sigma_1$ is reached when $\gamma > \max(0,2 \beta)$ or $\gamma = 2 \beta > 0$ and some conditions on the constants involved in the behavior of the coefficients. In that context, the condition \eqref{eq:cTO-inftyBIS} for $\sigma_2$ is more general since it is reached when $\gamma > \max(0,\beta-1)$ or $\gamma =  \beta - 1 > 0$ and some conditions on the constants involved in the behavior of the coefficients. 

A similar framework  is considered in \cite{PL2} and for the reader convenance we just briefly check that it falls in the framework developed before by slightly modifying the arguments presented in the previous section. The integrability conditions on $b$ and $c$ may be probably weaken. For the sake of clarity we do not follow this line of research but rather focus on the  new arguments which are necessary in order to deal with the unbounded domain $\Omega = \R^d$.

\smallskip
{\bf Condition \ref{H1}.} The definition of the operator is still made through the formula \eqref{def:StampOperator}. Under assumption   \eqref{eq:cTO-inftyBIS} on $\sigma_1$, denoting $\theta_1 := 1-\kappa/4$ and proceeding 
exactly as in the previous section during the proof of \eqref{eq:Stampacchia}, for any $f \in H^1(\R^d)$ and $\lambda \in \R$, we have  \Black
\bean
\langle(\lambda-\LL) f, f\rangle
&=&\int_{\R^d} |\nabla f|^2  + \int_{\R^d}  f \, b \cdot \nabla f   + \int_{\R^d}   (\lambda-c) f^2 
\\
&\ge&\theta_1\int_{\R^d} |\nabla f|^2   + \int_{\R^d}   (\lambda -\sigma_1)  f^2, 
\eean
by using successively the Cauchy-Schwarz inequality and the Young inequality. On the other hand,  
under assumption   \eqref{eq:cTO-inftyBIS} on $\sigma_2$, denoting $\theta_2 := 1$, for any $f \in H^1(\R^d)$ and $\lambda \in \R$, we  write
\bean
\langle(\lambda-\LL) f, f\rangle
&=&\theta_2 \int_{\R^d} |\nabla f|^2  +   \int_{\R^d}   (\lambda-\sigma_2) f^2,
\eean
  by performing one integration by part in the previous equation. \Black
In both cases, for and any $M > 0$, proceeding again as in the previous section  during the proof of \eqref{eq:Stampacchia}, and denoting from now on $\sigma = \sigma_i$, $\theta = \theta_i$
we have  
\bean
\langle(\lambda-\LL) f, f\rangle
&\ge&    \frac{\theta}{2} \| \nabla f \|_{L^2}^2+  \|  \sqrt{\sigma_-} f \|_{L^2}^2
 +  (\lambda -M )  \|  f \|_{L^2}^2 +  ( \frac{\theta C_S}{2} -  \|  \sigma   {\bf 1}_{\sigma \ge M} \|_{L^{d/2}}   )  \|  f \|_{L^{2^*}}^2,
\eean
by using  the Sobolev inequality (with associated constant $C_S)$ and the Holder inequality.
Taking $M>0$ large enough, and next $\kappa_1 > 0$ large enough, we finally obtain  
\beqn\label{eq:checkH1-ellipticC1}
\langle(\lambda-\LL) f, f\rangle
\ge    \frac{\theta}2  \| \nabla f \|_{L^2}^2+  \|  \sqrt{\sigma_-} f \|_{L^2}^2
 +     \|  f \|_{L^2}^2, \quad \forall \, \lambda \ge \kappa_1.
\eeqn
With the same arguments as in the previous section, we conclude that $\LL$ is the generator in $L^2$ of a positive semigroup $S_\LL$, so that \ref{H1} holds.

 \medskip
{\bf Condition \ref{H2}.} We may for instance use the third constructive argument (which is local) presented in section~\ref{subsec:diffusion-domain} and we establish 
$$
\exists \, f_0 \in H^1_0 \backslash \{0\}, \ f_0 \ge 0, \ \exists \, \kappa_0 \in \R, \quad \LL f_0 \ge \kappa_0 f_0.
$$
That is condition {\bf (ii)} in Lemma~\ref{lem:Existe1-Spectre2bis}, so that condition \ref{H2}  holds. 

 \medskip
{\bf Condition \ref{H3}.} We introduce again the splitting $\LL = \AA + \BB$ with $\AA := \kappa_1 - \kappa_0 +1$, so that from \eqref{eq:checkH1-ellipticC1}, the operator $\lambda - \BB = (\lambda - \kappa_0+1) + (\kappa_1-\LL)$ is invertible for any $\lambda \ge \kappa_\BB := \kappa_0-1$. 
We claim that the operator $(\lambda-\BB)^{-1}$ is compact for any $\lambda \ge \kappa_\BB$. 
 For that purpose, let us 
consider a sequence $(f_n)$ such that $(\lambda - \BB) f_n$ is bounded in $L^2$ and we have to prove that $(f_n)$ is relatively strongly compact. When condition \eqref{eq:cTO-inftyBIS} holds and because of the estimate \eqref{eq:checkH1-ellipticC1} and the very definition of $\BB$, we have 
\beqn\label{eq:checkH3-ellipticC1}
 \frac{\theta}2 \| \nabla f_n \|_{L^2}^2+  \|  \sqrt{\sigma_-} f_n \|_{L^2}^2
 +     \|  f_n \|_{L^2}^2 \le C, 
\eeqn
for some constant $C \in \R_+$. Because of the Rellich-Kondrachov theorem, we just have  to show that  
$$
\lim_{R \to \infty} \sup_n \int_{B_R^c} f_n^2 = 0. 
$$
But that last convergence may be established using the assumption  \eqref{eq:cTO-inftyBIS} in the following way. We write 
\bean
  \int_{B_R^c} f_n^2 
&=& \int_{B_R^c \cap \{ \sigma \ge K \}} f_n^2 +  \int_{B_R^c \cap \{ \sigma < K \}} f_n^2 
\\
&\le& \| f_n \|_{L^{2^*}}^{\frac{d-2}d}   [\hbox{\rm meas}(  B_R^c \cap \{ \sigma \ge K \} )]^{\frac2d} + \frac{1 }{ |K|}  \int \sigma_- f_n^2 ,
\eean
for any $K< 0$, by using the Holder inequality.  Using next the Sobolev inequality, the estimate \eqref{eq:checkH3-ellipticC1} and the assumption  \eqref{eq:cTO-inftyBIS}, \Black
we deduce 
\bean
\limsup_{R\to\infty}  \int_{B_R^c} f_n^2 
&\lesssim& \limsup_{R\to\infty}\inf_{K<0}  \Bigl\{  [\hbox{\rm meas}(  B_R^c \cap \{ \sigma \ge K \} )]^{\frac2d} + \frac{1 }{ |K|} \Bigr\} = 0, 
\eean
and the claim is proved. 
 As a consequence, we may apply  Lemma~\ref{lem:H3abstract-StrongC}-(2) and we deduce that \ref{H3} holds for both the primal and the dual problems.



 \medskip
{\bf Condition \ref{H4}.} As in \cite[Prop.~5.4]{MR4265692}, we establish the strong maximum principle by exhibiting  a barrier function and  using Lemma~\ref{lem:Irred-barrier+W->S}. 
An alternative argument should be to adapt the proof based on the Harnack inequality as presented in the previous section. 
Let us then consider $f \in D(\LL^k) \cap X_+ \backslash \{ 0 \}$ such that $(\lambda-\LL) f \ge 0$ with $k$ large enough ($k > d/2$ must be suitable) and $\lambda \ge \lambda_1$ large enough but fixed ($\lambda \ge \kappa_1$ is suitable). 
Using a very classically bootstrap argument based on iterated application of the Calderon-Zygmond elliptic regularity theorem and the  Morrey estimate, we have  $f \in C(\R^d)$. By assumption, there thus exist  $x_0 \in \R^d$, and two constants $\tau,r > 0$ such that $f  \ge \tau$ on  $B(x_0,r)$ and we take choose $x_0 = 0$ in order to simplify the notations.  We next fix $R > r$ and we observe that the function  
$$
 g(x) := \tau^*( g_0(|x|) - g_0(R)), \quad g_0(s) :=    \exp(\sigma r^2/2 -\sigma s^2/2) , 
$$
satisfies
\bean
(\tau^*)^{-1}(\lambda- \LL ) g 
&=& (\lambda -c) (g_0 - g_0(R)) + (d\sigma   - \sigma  b \cdot x - \sigma^2 r^2) \, g_0 
\\
&\le& [2 (|\lambda| + \| c \|_{L^\infty(B_R)}) +  \sigma  (d + \|  b \cdot x \|_{L^\infty(B_R)}) - \sigma^2 r^2] \, g_0 \le 0 
\eean
on $\OO := B(0,R) \backslash B(0,r)$   for $\sigma > 0$ large enough. We next fix $\tau^*$ such that $g=\tau$ on $\partial B(x_0,r)$. 
 We also observe that from  \eqref{eq:checkH1-ellipticC1}, $\lambda-\LL$ is coercive on $\OO$, in the sense that 
$$
\forall \, h \in H^1_0(\OO) \quad ( (\lambda-\LL) h, h)_{L^2(\OO)} \ge \| h \|_{L^2(\OO)}.
$$
In particular, $\lambda-\LL$ satisfies the weak maximum principle as explained in the proof of \eqref{eq:StampacchiaWeakMP}. Arguing as in the proof of Lemma~\ref{lem:Irred-barrier+W->S}, 
we deduce that $f \ge g > 0$ on $\OO$, what we also see directly by observing that $h :=   (g- f)_+ \in H_0^1(\OO)$,  $(\lambda-\LL) h \le 0$ and using that the weak maximum principle implies $h \le 0$, thus $h \equiv 0$ and finally $f \ge g$.  Because $R> r$ can be chosen arbitrarily large, we conclude with $f > 0$ on $\R^d$.


\Black

 \medskip
{\bf Condition \ref{H5}.}   
The reverse Kato's inequality condition is proved by using local arguments, so that it holds for the same reasons as in the previous section. Similarly, because the argument are local, the conclusion of 
Lemma~\ref{lem:Section7-BoccardoG} holds here. 

\medskip
As a consequence, using Theorem~\ref{theo:exist1-KRexistence}, Theorem~\ref{theo:KRgeometry1}, Theorem~\ref{theo:KRgeometry2} and Theorem~\ref{theo:ergodicity-compact-trajectories},
we may summarize our analysis in the following result. 

\begin{theo}\label{theo:diffusionStrongP} Consider the elliptic operator \eqref{eq:LLf=Rd} in the whole space and assume that the coefficients satisfy \eqref{eq:cTO-inftyBIS}. 
Then the conclusions  {\Blue \ref{S1}, \ref{S2} and \ref{S32}  hold in $L^2$}  as well as  \ref{E2} in $L^1_{\phi_1}$. 
\end{theo}

{\Blue 
While the conclusion \ref{S1}-\ref{S2} is essentially due to Lions \cite[8th course]{PL2}, we believe that the ergodicity property \ref{E2} is new.}
We do not present an exponential constructive estimate, which we believe is possible to prove, but would require significantly more development. 

 \medskip
\subsection{Diffusion in $\R^d$ with weak potential confinement} 
\label{subsection:diffusionWeakPotential}
We consider in this section  the same elliptic operator \eqref{eq:LLf=Rd} with now a weak confinement condition assuming that $c$ converges to a constant. With 
no loss of generality,  we may assume $c \to 0$.   More precisely, we consider the  elliptic operator
\beqn\label{eq:LLf=Rd-cTO0}
\LL f := \Delta f + b \cdot \nabla f + r c f, \ 
\eeqn
with $c \in C_0(\R^d)$, $b \in C_0(\R^d)$ and $r \in \R_+$ a parameter.  When not necessary in the discussion we will take $r=1$.   
The associated first eigenvalue problem in such a situation has been studied in \cite[8th and 9th courses]{PL2} to which we refer for more details. 
We define 
$$
\lambda_1 = \lambda_1 (r)  := \inf \{ \kappa \in \R; \ (\lambda - \LL)^{-1} \hbox{ well defined and positive for any } \lambda \ge \kappa\}.
$$
Proceeding exactly as in the proof of \ref{H1} in the preceding section, we see that  the operator $\lambda - \LL$ is invertible for any $\lambda > \| (rc+ |b|^2/4)_+ \|_{L^\infty}$, and then its inverse is positive. 
Because the proof of  \ref{H2} in the preceding section also applies here, we deduce that the infimum $\lambda_1$ of the set $\II$ of  real resolvent values is well defined with $\lambda_1 \in  (\kappa_0,\kappa_1)$, for some constructive constants $\kappa_i \in \R$. 

\smallskip
We split now the discussion into two cases. 

\smallskip
{\bf Case 1. } We start considering {\bf the case $b = 0$}. In that case,  $\LL$ is self-adjoint so that $\lambda_1$ is also characterized by 
$$
\lambda_1 = \sup_{\| f \|_{L^2} = 1} \EE(f), 
$$
with 
$$
\EE(f) := (\LL f, f) = r \int c f^2 - \int |\nabla f|^2.
$$
We make  the following elementary observations :  

\smallskip
$\bullet$ We claim that $\lambda_1 \ge 0$. Taking $f_n(x) := n^{-d/2} u(x/n)$ for some function $u \in H^1(\R^d)$, $\| u \|_{L^2} = 1$, we compute 
\bean
- \EE(f_n) 
&=& \int |\nabla f_n|^2 -  \int_{B_R} rc f_n^2 -\int_{B_R^c} rc f_n^2
\\
&\le& \frac1{n^2} \int |\nabla u|^2 + \| r c \|_{L^\infty(B_R)}  \int_{B_{R/n}} u^2 +\| r c \|_{L^\infty(B^c_R)}, 
\eean
for any $R > 0$, so that 
$$
-\lambda_1 \le \limsup (- \EE(f_n)) \le 0.
$$

\smallskip
$\bullet$ We claim that $\lambda_1 = 0$ when $c \le 0$. In that case, we have  $\EE(f) \le 0$ for any $f \in H^1(\R^d)$, and we deduce the reverse inequality $\lambda_1 \le 0$. 
In particular, as a function $\lambda_1 = \lambda_1(r)$ of $r \ge 0$, we have $\lambda_1(0)=0$. We also claim that $\lambda_1(r) \to \infty$ as $r\to\infty$ when $c_+ \not\equiv 0$.
We may indeed fix $f \in H^1(\R^d)$, $\| f \|_{L^2} = 1$, $\hbox{supp}\, f \subset \hbox{supp}\, c_+$, and we compute 
$$
\EE(f) = r \int_{\R^d} c_+ f^2 -  \int |\nabla f|^2  \to   \infty, \quad\hbox{as}\quad r \to \infty.
$$
 
\smallskip
$\bullet$ We finally observe that  $\lambda_1   : \R_+ \to \R_+$ is convex since it is defined as the supremum of linear
functions $r \to \EE(f)$ for any fixed $f \in H^1(\R^d)$. As a consequence, we have the 
 following alternative: 

\quad - $\lambda_1 \equiv 0$;

\quad - $\exists r_0 \in [0,\infty)$ such that $\lambda_1 (r) =  0$ for $r \le r_0$ and  $\lambda_1 (r) >  0$ for $r  > r_0$.

\smallskip
Concerning the value of $r_0$, it may happen that $r_0 > 0$, and that is the case when $c \in L^{d/2}$ because of the Sobolev inequality, or that 
$r_0 = 0$, and that is the case for instance when $c \ge 0$, $c(x) = |x|^{-m}$ for $x \in B_R^c$, $m \in (0,2)$, $R > 0$. To prove that last claim, we may take the same sequence $(f_n)$ as above, and we  compute 
 \bean
\EE(f_n) 
&\ge&  \int_{B_R^c} r |x|^{-m} f_n^2 - \int |\nabla f_n |^2 dx 
\\
&=& \frac{r }{ n^m} \int_{B_{R/n}^c} |x|^{-m}u^2 - \frac{1}{ n^2} \int |\nabla u |^2 dx  >0, 
\eean
for $n$ large enough (whatever is the value of $r>0$).

\smallskip
{\bf About condition \ref{H3}.} 
It is established in  \cite{PL2} that when $\lambda_1  = 0$, the condition \ref{H3} is not satisfied and  there does not exist a first eigenfunction $f_1 \in L^2(\R^d)$ to the 
operator $\LL$ defined by \eqref{eq:LLf=Rd-cTO0}. We refer to \cite[8th course]{PL2} for a proof of that result. 
On the other hand, we claim that the condition \ref{H3} is   satisfied
when $\lambda_1  > 0$. 
Consider indeed three sequences $(\lambda_n)$ of $\R$, $(f_n)$ of $H^1(\R^d)$ and $(\eps_n)$ of $L^2(\R^d)$ such that $(\lambda_n - \LL) f_n = \eps_n$, $\eps_n, f_n \ge 0$, $\| f_n \|_{L^2} = 1$, for any $n \ge 1$, 
$\lambda_n \to \lambda_1$ and $\eps_n \to 0$ in $L^2$ as $n\to\infty$. We then have 
$$
\lambda_n - \EE(f_n) = \langle (\lambda_n - \LL)f_n,f_n \rangle = \langle  \eps_n,f_n \rangle \to 0,
$$
as $n\to\infty$. By definition of $\EE$ and boundedness of $c$, we see that $(f_n)$ is bounded in $H^1$. As a consequence, up to the extraction of a subsequence, we have 
$f_n \to f_1 \ge 0$ in $\Lloc^2$ and thus next $(\lambda_1 - \LL) f_1 = 0$ in the variational sense and
$$
\int c f^2_n \to \int c f_1^2, \quad \| \nabla f_1 \|_{L^2} \le \liminf \| \nabla f_n \|_{L^2},
$$
where we have used the dominated convergence  theorem of Lebesgue and the fact that $c \to 0$ at infinity in order to get the first convergence. We finally deduce  
$$
\EE(f_1) \ge \limsup \EE(f_n) = \lambda_1 > 0, 
$$
so that $f_1 \not\equiv 0$, and  the condition \ref{H3} is verified. 

 \smallskip
As a conclusion, for a self-adjoint operator, condition \ref{H3} is automatically fulfilled by its adjoint and the conditions \ref{H4} and \ref{H5} have been proved in a general situation, including the present framework. 
 The same conclusions of existence, uniqueness and asymptotic stability of the first eigentriplet solution $(\lambda_1,f_1,\phi_1)$ as in section~\ref{subsection:diffusionStrongPotential} hold true when $\lambda_1 > 0$. 

 \smallskip
{\bf Case 2. } We consider the {\bf  general case $ b \in C_0(\R^d)$}.

$\bullet$ We claim that $\lambda_1 \ge 0$. Adapting the second constructive argument in the proof of \ref{H2} in  Section~\ref{subsec:diffusion-domain}, 
we consider $\chi \in C^1_c(\R_+) \cap W^{2,\infty}(\R_+)$ such that  ${\bf 1}_{[0,1/2]} \le \chi \le {\bf 1}_{[0,1]}$,  $\chi' \le 0$ on $[0,1]$, $\chi(s) :=  (1-s)^2/2$ on $[\eta,1]$ with $\eta \in (1/2,1)$ large enough
in such a way that 
\beqn\label{eq:conditionH3coursPL2}
\chi''(s) + (d-1) \chi'(s)/s \ge 1/2, \quad \forall \, s \in (\eta,1),
\eeqn
and define $f_0(x):=\chi(|x-x_0|/n)$ for $|x_0|$ large enough to be chosen later. We have $\supp f_0\subset B_n(x_0)$ for any $n\geq1$ and we compute
%
\[
\LL f_0(x)
=
\frac{1}{n^2} \bigl\{ \chi''(r/n)+\frac{d-1}{r/n}\chi'(r/n)\bigr\} +\frac{1}{n}b(x)\cdot\hat y\,\chi'(r/n)+c(x)\chi(r/n)
\]
where $y=x-x_0$ and $r=|y|$.
On $B_{n\eta}(x_0)$, we have
\[
\LL f_0(x) \geq -\frac{\|\chi''\|_\infty}{n^2}-\frac{d-1}{n^2} \Big\|\frac{\chi'(r)}{r}\Big\|_\infty -\frac{\|\chi'\|_\infty}{n}\sup_{B_n(x_0)} |b|-\|\chi\|_\infty \sup_{B_n(x_0)}|c|.
\]
On $B_n(x_0)\setminus B_{\eta n}(x_0)$, thanks to \eqref{eq:conditionH3coursPL2},  we have
\begin{align*}
\LL f_0(x)
&\geq\frac{1}{2n^2}-\frac{\|\chi'\|_\infty}{n}\sup_{B_n(x_0)}|b|- \|\chi\|_\infty\sup_{B_n(x_0)}|c|.
\end{align*}
Let now fix $\eps>0$ and choose first $n$ large enough so that
$$
-\frac{\|\chi''\|_\infty}{n^2}-\frac{d-1}{n^2}\Big\|\frac{\chi'(r)}{r}\Big\|_\infty\geq-\frac{\eps}{2}\inf_{(0,\eta)}\chi.
$$
Then, using that $b,c\in C_0(\R^d)$, we can take $|x_0|$ large enough so that
$$
-\frac{\|\chi'\|_\infty}{n}\sup_{B_n(x_0)}|b|-\|\chi\|_\infty\sup_{B_n(x_0)}|c|\geq-\frac{\eps}{2}\inf_{(0,\eta)}\chi
$$
and
$$
\frac{\|\chi'\|_\infty}{n}\sup_{B_n(x_0)}|b|+\sup_{B_n(x_0)}|c|\leq\frac{1}{2n^2}.
$$
Gathering the above inequalities, we obtain 
$$
\LL f_0\geq-\eps f_0,
$$
and the condition \ref{H2} is verified with $\kappa_0=-\eps$. Because  $\eps>0$ can be choose arbitrarily small, we conclude with   $\lambda_1\geq0$.

$\bullet$ We claim that $\lambda_1=0$ when $\sigma_2\leq0$.
Indeed, we have already seen that
$$
\langle\LL f,f\rangle=-\int_{\R^d}|\nabla f|^2+\int_{\R^d}\sigma_2f^2,
$$
from which we deduce that
$$
\frac{d}{dt}\|S_tf\|^2=2\langle\LL f,f\rangle\leq0.
$$
This ensures that \ref{H1} is verified with $\kappa_1=0$ and so $\lambda_1\leq0$.

\Black

$\bullet$ We claim that $\lambda_1 > 0$ when $c_+ \not\equiv 0$ and $r > 0$ is large enough.
For simplifying notations and up to translation and dilatation, we may reduce to the case  $c \ge c_0 {\bf 1}_{B(0,1)}$ with $c_0 > 0$. Adapting the second constructive argument in the proof of \ref{H2} in 
Section~\ref{subsec:diffusion-domain}, we consider $\chi \in C^1_c(\R_+) \cap W^{2,\infty}(\R_+)$, ${\bf 1}_{[0,1/2]} \le \chi \le {\bf 1}_{[0,1]}$, supp$\chi = [0,1]$, $\chi''(1) = 1$, $\chi' \le 0$ on $[0,1]$
and we set  $f_0 (x) := \chi( |x|)$.   We compute 
\bean
\LL f_0 =   \chi''(|x|)    +  \chi'(|x|) ((d-1)/|x| + b \cdot \hat x) + r c(x) \chi(|x|).
\eean
On the one hand, we fix $\eta \in (1/2,1)$, $1-\eta$ small enough, in such a way that 
$$
\| \chi' \|_{L^\infty(\eta,1)} \bigl( 2(d-1) + \| b \|_{L^\infty} \bigr) \le 1/4, \quad 1/2 \le \| \chi'' \|_{L^\infty(\eta,1)} ,
$$
and thus 
$$
\LL f_0 \ge \frac14  \ge \frac14 f_0 \quad\hbox{on}\quad B(0,\eta)^c.
$$
On the other hand, we fix $r > 0$, large enough, in such a way that 
$$
\| \chi'' \|_{L^\infty} + \| \chi' \|_{L^\infty} \bigl( 2(d-1) + \| b \|_{L^\infty} \bigr) \le \kappa (r) := \frac12 r c_0 \inf_{[0,\eta)} \chi, 
$$
and thus 
$$
\LL f_0 \ge   \kappa(r)  \ge \kappa(r) f_0 \quad\hbox{on}\quad B(0,\eta).
$$
As a conclusion, we  have  established that condition {\bf (ii)} in the statement Lemma~\ref{lem:Existe1-Spectre2bis} holds with $\kappa_0 := \min(1/4,\kappa(r))$, and that ends the constructive proof of condition \ref{H2} by using 
Lemma~\ref{lem:Existe1-Spectre2bis}. That implies in particular the claim since then $\lambda_1 \ge \kappa_0 > 0$. 
%
%
%
%
%
%

\smallskip
$\bullet$ We finally claim   that \ref{H3} holds when $\lambda_1 > 0$. 
To see that, we consider again three sequences $(\lambda_n)$ of $\R$, $(f_n)$ of $H^1(\R^d)$ and $(\eps_n)$ of $L^2(\R^d)$ such that 
$(\lambda_n - \LL) f_n = \eps_n$, $\eps_n, f_n \ge 0$,  $\| f_n \|_{L^2} = 1$, for any $n \ge 1$, $\lambda_n \searrow  \lambda_1$ and $\eps_n \to 0$ in $L^2$ as $n\to\infty$. As a consequence,  we have 
$$
\lambda_n + \int |\nabla f_n|^2 - \int f_n b \cdot \nabla f_n  -\int c f_n ^2  = ((\lambda_n - \LL)f_n,f_n) = \langle \eps_n,f_n \rangle \to 0,
$$
as $n\to\infty$. Using the  boundedness of $c$, $b$ and $\lambda_n$, we see that $(f_n)$ is bounded in $H^1$. As a consequence, up to the extraction of a subsequence, we have 
$f_n \to f_1 \ge 0$ in $\Lloc^2$.  We assume by contradiction that $f_1 \equiv 0$. We deduce that 
$$
\int c f^2_n \to 0, 
\quad 
\int f_n b \cdot \nabla f_n \to 0, 
 $$
where we have used the dominated convergence  theorem of Lebesgue and the fact that $b,c \to 0$ at infinity. 
We thus obtain 
$$
0 < \lambda_1 \le \lambda_n + \int |\nabla f_n|^2 =  \int f_n b \cdot \nabla f_n + \int c f_n ^2   + \langle \eps_n,f_n \rangle \to 0,
$$
and our contradiction. 
 So that $f_1 \not\equiv 0$, and  the condition \ref{H3} is verified. 
 
 \smallskip
 For the dual problem, from the above analysis, we know that there exist two sequences $(\phi_n)$ of $H^1(\R^d)$, $(\eps_n)$ of $L^2(\R^d)$ such that $(\lambda_n - \LL^*) \phi_n = \eps_n$, $\eps_n, \phi_n \ge 0$ and $\| \phi_n \|_{L^2} = 1$, for any $n \ge 1$, and $\eps_n \to 0$ in $L^2$ as $n\to\infty$. 
But we face the same situation as previously, since again
$$
\lambda_n + \int |\nabla \phi_n|^2 - \int \phi_n b \cdot \nabla \phi_n  -\int c \phi_n ^2  = ((\lambda_n - \LL^*)\phi_n,\phi_n) = (\eps_n, \phi_n) \to 0,
$$
and thus the same conclusion, namely $\phi_n \to \phi_1$, with $\phi_1 \in H^1(\R^d)$, $\phi_1 \ge 0$, $\phi_1 \not\equiv 0$. 

    {\Blue 
 \smallskip
The conditions \ref{H4} and \ref{H5} have been proved in a general situation, including the present framework. 
For the same reason as in  section~\ref{subsection:diffusionStrongPotential}, the same conclusions  hold true when $r > 0$ is large enough (and thus $\lambda_1 > 0$)
and more precisely we have established the following result.

\begin{theo}\label{theo:diffusionWeakP} 
Consider the elliptic operator \eqref{eq:LLf=Rd} in the whole space and assume that the coefficients satisfy $b,c \in C_0(\R^d)$, $c_+ \not\equiv 0$. For $r > 0$ large enough, 
the conclusions  {\Blue \ref{S1}, \ref{S2} and \ref{S32}  hold in $L^2$}   (with $\lambda_1 > 0$) as well as  \ref{E2} in $L^1_{\phi_1}$. 
\end{theo}

While the conclusion \ref{S1}-\ref{S2} is essentially due to Lions \cite[9th course]{PL2}, we believe that the ergodicity property \ref{E2} is new. 
}

 \medskip
\subsection{Diffusion in $\R^d$ with   drift confinement} 
\label{subsec-diffusionwith drift}
We now consider the elliptic operator
$$
\LL f := \Delta f + b \cdot \nabla f + c f, 
$$
with a drift confinement as it is the case for the Fokker-Planck operator. More precisely, and for the sake of simplicity, we assume here
\beqn\label{eq:diffusionWithDrift-defb&gamma}
b = \nabla U, \quad U(x) = \frac1\gamma\langle x \rangle^\gamma, \quad \gamma > 0. 
\eeqn
When $\gamma = 2$ and $c =x$, that operator corresponds to the classical harmonic Fokker-Planck operator which is known to be related to the standard Poincaré inequality and 
to the standard log-Sobolev inequality,  see \cite{MR889476,MR1307413,MR1704435} or more recently \cite{MR2386063,MR4265692} and the references therein. When $c = \Div b$, the operator $\LL$ is on divergence form and $\LL^*1 = 0$, so that $(0,1) \in \R \times L^\infty(\R^d)$ is a solution to the dual first eigenvalue problem. Existence of stationary solution $f_1$ (which is also the first eigenfunction) and its stability have been widely studied. We refer for instance to \cite{MR1751701,MR1856277,MR1912106,MR2381160}  as well as to \cite{MR4265692,MR3488535,MR3779780} which techniques will be adapted here. 

\smallskip
In the present situation, we impose that the contribution of $c$ has lower influence at the infinity that the drift term $b$ and we assume 
\beqn\label{eq:dissipDrift-HYPc1}  
c \in \Lloc^\infty(\R^d), \quad \exists \, C_0, R_0 > 0, \ \forall \, x \in B_{R_0}^c ,   \ |c(x)| = o( |x|^{2(\gamma-1)}) . 
\eeqn
We further assume that 
\beqn\label{eq:dissipDrift-HYPc2}
c \ge \Div b \quad\hbox{when}\quad \gamma \in (0,1].
\eeqn

\smallskip
The action of the drift term will be revealed through the choice of a convenient ``confining space". More precisely, for a weight function $m : \R^d \to [1,\infty)$, we will work in a weighted Lebesgue space. 
Our analysis is based on the following elementary computation which can be readily adapted from \cite[Lem.~2.1]{MR4265692}, \cite[Lem.~3.8]{MR3488535} and \cite[Lem.~3.8]{MR3779780}.

\begin{lem}\label{lem:diffus-bconf:identityL} For any $p \in [1,\infty)$, any weight function $m$ and 
 any smooth, rapidly decaying  function $f$, we  have 
\beqn\label{eq:ch3:ch3FP:Lid1}
\int (\LL f ) \, f |f|^{p-2} m^p = - (p-1)\int |\nabla f |^2 |f|^{p-2} m^p + \int |f|^p m^p \varphi_1,
\eeqn
with 
\beqn\label{eq:ch3:ch3FP:Lid1defm}
\varphi_1 :=  (p-1) \, \frac{  |\nabla m|^2 }{ m^2} + \frac{ \Delta m}{ m} 
+ \left(c  - \frac{1 }{ p}  \mbox{\rm div} \, b \right)   - b \cdot \frac{\nabla m}{ m}
\eeqn
as well as 
\beqn\label{eq:ch3:ch3FP:Lid2}
\int (\LL f ) \, f |f|^{p-2} m^p = - (p-1) \int |\nabla (fm) |^2  |fm|^{p-2} + \int |f|^p m^p \varphi_2,
\eeqn
with 
\beqn\label{eq:ch3:ch3FP:Lid2defm}
\varphi_2 := 2 (1 - \frac{1 }{ p}) \frac{|\nabla m|^2}{ m^2}  + ( \frac 2p -1)\frac{ \Delta m }m+  
 \left(c  - \frac{1 }{ p}   \mbox{\rm div} \, b \right)
     - b \cdot \frac{\nabla m }{ m} . 
\eeqn
\end{lem} 
In order to simplify the discussion, we restrict ourself to the exponent $p=2$ and to the exponential weight function $m =  e^{\kappaa \langle x \rangle^s}$, $s \in (0,\gamma]$, $a>0$. We thus work in the Banach lattice $X := L^2_m$. 
We observe that 
\bean
&&\frac{\nabla m  }{ m} = s \kappaa  x  \langle x \rangle^{s-2} \sim s \kappaa  |x|^{s-1}, 
\\
&&\frac{\Delta m }{ m}  =   s \kappaa d  \langle x \rangle^{s-2} +  s (s-2) \kappaa |x|^2 \langle x \rangle^{s-4} 
+  (s\kappaa)^2  |x|^2 \langle x \rangle^{2s-4}  \sim (s \kappaa)^2  |x|^{2s-2}, 
\\
&&\Div b =    d  \langle x \rangle^{\gamma-2} +  (\gamma-2)  |x|^2 \langle x \rangle^{\gamma-4}  \sim (d+\gamma-2)   |x|^{\gamma-2},   
\\
&&b \cdot \frac{\nabla m  }{ m}  =   s \kappaa  x  \langle x \rangle^{s-2} \cdot x  \langle x \rangle^{\gamma-2}    \sim s \kappaa  |x|^{s+\gamma-2},   
\eean
so that the contribution of $(c-\Div b/2)$ is always negligible at infinity, and we get 
\beqn\label{}
 \varphi_i  \sim     (s\kappaa)^2 |x|^{2s-2}   -  s \kappaa  |x|^{s+\gamma-2}. 
\eeqn
We denote 
\bean 
\kappaa'  &:=& s\kappaa > 0 \quad \mbox{if}\quad s \in (0, \gamma), \\
\kappaa'  &:=& \kappaa\gamma - 2(\kappaa\gamma)^2 >  0 \quad \mbox{if}\quad s = \gamma \hbox{ and } \kappaa \in (0, 1/(\sqrt{2}\gamma)).
\eean 

We then face to three cases : 

\smallskip
(i)  $\gamma > 1$ : taking $s \in ((2-\gamma)_+,\gamma)$, we have $\varphi_i \sim - \kappaa'  |x|^{s+\gamma-2} \to - \infty$ with $s+\gamma-2 > 0$; 

(ii) $\gamma = 1$ : taking $s = \gamma$, $\kappaa < 1/(\sqrt{2}\gamma)$, we have $\varphi_i \to - \kappaa' $; 

(iii) $\gamma \in (0,1)$ : taking $s = \gamma$, $\kappaa < 1/(\sqrt{2}\gamma)$, we have $\varphi_i \sim -    \kappaa'  |x|^{2\gamma-2} \to 0$ with $2\gamma-2 < 0$. 

\smallskip

\smallskip
{\bf Condition \ref{H1}.} In any of the above cases, we have from \eqref{eq:ch3:ch3FP:Lid1}
$$
((\lambda-\LL)f,f) = \int |\nabla f|^2 m^2 + \int (\lambda-\varphi_1) f^2 m^2, 
$$
for $\lambda \in \R$, with $\inf (\lambda -\varphi_1) > 0$  for $\lambda \ge \kappa_1$ and $\kappa_1 > 0$ large enough. We deduce that $\lambda-\LL$ is coercive for $\lambda \ge \kappa_1$.
With the same arguments as in  section~\ref{subsec:diffusion-domain}, we conclude that $\LL$ is the generator in $L^2_m$ of a positive semigroup $S_\LL$, so that \ref{H1} holds.

\smallskip
{\bf Condition \ref{H2}.} When $\gamma > 1$, the same arguments as in Section~\ref{subsection:diffusionStrongPotential} imply that condition \ref{H2} holds for some $\kappa_0 \in \R$.
When $\gamma \in (0,1]$, we have $\LL^* 1 = c-\Div b \ge 0$ from \eqref{eq:dissipDrift-HYPc2} and \ref{H2} holds with $\kappa_0 = 0$.

\smallskip
{\bf Conditions \ref{H4} and \ref{H5}.} The strong maximum principle holds here because for instance we may apply  the same barrier function argument as presented in Section~\ref{subsection:diffusionStrongPotential}.  
The reverse Kato's inequality condition is proved by using local arguments, so that it holds for the same reasons as in the previous section.

\smallskip
{\bf Condition \ref{H3}.} We define the multiplication operator $\AA$ and the elliptic operator $\BB$ by 
$$
\AA := M \chi_R, \quad \BB := \LL - \AA,
$$
for $M, R > 0$ and $\chi_R(x) := \chi(x/R)$ with $\chi \in \DD(\R^d)$, ${\bf 1}_{B_1} \le \chi \le {\bf 1}_{B_2}$. 
We fix $\kappa_\BB < \kappa_0$ in case (i), $\kappa_\BB := - \kappaa'/4$ in case (ii) and 
$\kappa_\BB:= 0$ in case (iii), and we set $\kappaa'' := \kappaa'/2$.
Choosing $M,R >0$ large enough,  from Lemma~\ref{lem:diffus-bconf:identityL} and the discussion which follows, we deduce that 
\beqn\label{eq:diffusb:estimB}
((\BB-\alpha) f, f) \le -\int |\nabla f|^2 m^2 - \kappaa'' \int   f^2 ({\bf 1}_{B_1} + {\bf 1}_{B_1^c} |x|^{s+\gamma-2}) m^2,
\eeqn
for any $\alpha \ge \kappa_\BB$ and any nice function $f$. We classically deduce that $\alpha -\BB$ is coercive and thus invertible. We discuss the three different cases. 

\smallskip
- In the first case $\gamma > 1$, so that $s+\gamma-2 > 0$, we see that the operator $\RR_\BB(\alpha)$ is compact from Rellich-Kondrachov theorem, so that also $\WW(\alpha) := \RR_\BB(\alpha) \AA$ for any $\alpha \ge \kappa_\BB$.  
  We may thus apply   Lemma~\ref{lem:H3abstract-StrongC}-(2) and we deduce that \ref{H3} holds for both the primal and the dual problems. 
 

\smallskip
- In the  case $\gamma = 1$, so that $2\gamma-2 \le 0$, the operator $\RR_\BB(\alpha)$ is not compact anymore. However, for any sequence $(f_n)$ which is bounded in $L^2_m$, we define the sequence $(g_n)$ by $g_n := \AA f_n$, and $(g_n)$ is bounded in $L^2_{\tilde m}$ with $\tilde m := e^{\tilde\kappaa \langle x \rangle^\gamma}$, $\tilde\kappaa \in (\kappaa,1/\sqrt{2\gamma})$. 
Using the dissipativity estimate \eqref{eq:diffusb:estimB}  in $L^2_{\tilde m}$, we see that $\BB-\alpha$ is dissipative in $L^2_{\tilde m}$ for any $\alpha \ge \kappa_\BB$, and more precisely the sequence $(h_n)$ defined by $h_n := \RR_\BB(\alpha) g_n$   satisfies
$$
\int |\nabla h_n|^2 m^2 + \tilde \kappaa'' \int   h_n^2 ({\bf 1}_{B_1} + {\bf 1}_{B_1^c} |x|^{2\gamma-2})  \tilde m^2 \le  \int  g_n^2  \tilde m^2.
$$
Using that $  |x|^{2\gamma-2}  \tilde m^2/m^2 \to \infty$ as $|x| \to \infty$, that implies that $(h_n)$ is relatively compact in $L^2_m$.
More precisely, the above estimates show that $\WW(\alpha) := \RR_\BB(\alpha) \AA :L^2_m \to H^1_m \cap L^2_{m^\sharp}$ with $m^\sharp := m^{1/2} \tilde m^{1/2}$ and in particular  we have established that $\WW(\alpha) := \RR_\BB(\alpha) \AA$ is a compact operator in $L^2_m$ uniformly on $\alpha \ge \kappa_\BB$ because of the Rellich-Kondrachov theorem and the fact that $m = o(m^\sharp)$.   
Since $\RR_\BB(\alpha)$ is  bounded in $ \BBB(L^2_m)$ uniformly for any $\alpha > \kappa_\BB$, the operator $\LL$ satisfies the splitting structure \ref{HS1}  and, applying  Lemma~\ref{lem:H3abstract-StrongC}-(2), we deduce that \ref{H3} holds for both the primal and the dual problems. 
\Black


\smallskip
At this stage, when $\gamma \ge 1$, 
we obtain   a solution $(\lambda_1,f_1,\phi_1)$ to the first eigentriplet problem \eqref{eq:ex1PLL-triplet} by using Theorem~\ref{theo:exist1-KRexistence}. 

\smallskip
{\bf Condition \ref{HS3}.} In the  case $\gamma \in (0,1)$, the same as in the case $\gamma=1$ holds except that  $\RR_\BB(\alpha)$ is not uniformly bound in $ \BBB(L^2_m)$ for $\alpha \ge \kappa_\BB$ because we are in the critical case $\kappa_\BB = \kappa_0$. We do not know how to adapt the stationary approach in that situation and we thus aim to use a dynamical approach through the 
%
 use of Theorem~\ref{theo:KRexistTER} with the above splitting $\LL = \AA + \BB$ and $N:=[d/4]+1$. 
We set $X = X_1 := L^2_m$ and $X_0 := L^1$.  
The proof of condition \ref{HS3} is an immediate consequence of the following estimate.  

\begin{prop}\label{prop:diffuse-KRexistTER} We define $\Theta_\zeta(t) := e^{-\zeta t^{\gamma/(2-\gamma)}}$. For $N := [d/4]+1$, there hold
\begin{itemize}
\item[(i)]  $S_{\BB}\in L^\infty_t (\BBB(X_1))$; 
\item[(ii)]  $S_{\BB} \AA \Theta_\zeta^{-1} \in L^\infty_t (\BBB(X_i))$ for $i=0,1$ and any $\zeta \in (0,\zeta^*)$; 
\item[(iii)]  $(S_{\BB} \AA)^{(*N)} \Theta_\zeta^{-1} \in L^\infty_t (\BBB(X_0,X_1))$ for any $\zeta \in (0,\zeta^*/2)$. 
\end{itemize}
\end{prop}
 
 The proof of Proposition~\ref{prop:diffuse-KRexistTER} is similar to the proofs of  \cite[Lem.~2.1]{MR4265692}, \cite[Lem.~2.2]{MR4265692}, \cite[Lem.~2.3]{MR4265692} and \cite[Lem.~2.4]{MR4265692}.
 For the sake of completeness we however present the main lines of the proof. We start with a technical result that we will use during the proof of Proposition~\ref{prop:diffuse-KRexistTER}.

\begin{lem} \label{lem:ch2:Tn1} Consider two Banach spaces $X_0$, $X_1$ and a function $u : \R_+ \to \BBB(X_0) + \BBB(X_1)$ which satisfies  
\begin{itemize}
\item[{\bf (a)}]  $u  \Theta^{-1} \in L^\infty(0,\infty; \BBB(X_0) \cap \BBB(X_1) )$; 
\item[{\bf (b)}]    $ u \wp \in L^\infty(0,\infty; \BBB(X_0,X_1))$;
   \end{itemize}
for any exponentially decaying function $\Theta = \Theta_\zeta = e^{-\zeta t^\varsigma}$, $\forall\zeta \in (0,\zeta^*)$, and for the power function $\wp :=  t^{-\alpha}$, with $\zeta^* >0$, $\varsigma \in (0,1]$ and $\alpha \ge 0$ fixed. 
Then 

\begin{itemize}
\item[{\bf (c)}]  there exists $N$ such that  $u^{(*N)} \widetilde\Theta \in L^\infty(0,\infty; \BBB(X_0,X_1))$,
   \end{itemize}
for any $ \widetilde\Theta = \Theta_{ \tilde\zeta}$, $\tilde\zeta \in (0,\zeta^*/2)$. 
\end{lem}

\begin{proof}[Proof of Lemma~\ref{lem:ch2:Tn1}.]
A similar argument is developed in \cite[Lem.~2.17]{MR3779780}, \cite[Lem.~2.4]{MR3488535}, \cite[Prop.~2.5]{MR3465438} and \cite[Lem.~2.4]{MR4265692}.

\smallskip
{\sl Step 1. } Consider two functions $v$ and $w$ which satisfy the estimate {\bf (a)}. For $\XX = X_0$ or $\XX =  X_1$, we compute
\bean
\| v * w (t) \|_{\XX \to \XX} &\le& \int_0^{t} \| v(t-s) w(s) \|_{\XX\to \XX} \, ds  
\\
&\le&  \int_0^{t} C^v_{\XX} \Theta(t-s) \, C^w_{\XX} \Theta(s)  \, ds \le C^v_{\XX} C^w_{\XX} \, t \,  \Theta(t),
\eean
with obvious notation and 
where we have used that $\Theta(t-s) \,  \Theta(s) \le  \Theta(t)$ for any $0 < s < t$.  Since for any $\zeta' \in (0,\zeta)$, there exists a constant $C$ such that 
$t \Theta_\zeta (t) \le C \, \Theta_{\zeta'} (t) $ for any $t \ge 0$, we see that the function $v * w$ satisfies the same estimate  {\bf (a)} for any $\Theta = \Theta_\zeta$, $\zeta \in (0,\zeta^*)$.

\smallskip\noindent
{\sl Step 2. } 
Consider two functions $v$ and $w$ which satisfy the estimates {\bf (a)} and {\bf (b)} with $\alpha \ge 1$.  
We compute 
\bean
\| v * w (t) \|_{X_0 \to X_1} &\le& \int_0^{t/2} \| v(t-s) w(s) \|_{X_0\to X_1} \, ds +  \int_{t/2}^t \| v(t-s) w(s) \|_{X_0\to X_1} \, ds
\\
&\le& \int_0^{t/2} C^v_{01} (t-s)^{-\alpha}   \, C^w_{0} \Theta(s) \, ds +  \int_{t/2}^t  C^v_{1} \Theta(t-s) \, C^w_{01} s^{-\alpha} \, ds
\\
&=& [ C^v_{1}  \, C^w_{01}  +  C^v_{01}   \, C^w_{0} ] \, \Theta(0) \, t^{-\alpha+1} \,  \int_0^{1/2} (1-\tau)^{-\alpha} \, d\tau,
\eean
with obvious notation and  we have used that $\Theta$ is a decaying function.  As a consequence, the function $v * w$ satisfies estimate  {\bf (b)} with an exponent $\alpha-1$ instead of $\alpha$.  

\smallskip\noindent
{\sl Step 3. } Consider two functions $v$ and $w$ which satisfy the estimates {\bf (a)} and {\bf (b)} with $\alpha  \in [0,1)$. We compute
\bean
\| v * w (t) \|_{X_0 \to X_1} 
&\le& \int_0^{t/2} \| v(t-s) w(s) \|_{X_0 \to X_1} \, ds +  \int_{t/2}^t \| v(t-s) w(s) \|_{X_0 \to X_1} \, ds
\\
&\le&  \int_0^{t/2} C^v_1 \Theta(t-s) \, C^w_{01} s^{-\alpha} \, ds +  \int_{t/2}^t  C^v_{01} (t-s)^{-\alpha}  \, C^w_0 \Theta(s)  \, ds
\\
&\le&   C^v_1  \, C^w_{01}  \, \Theta(t/2)  \int_0^{t/2}  s^{-\alpha} \, ds  +  C^v_{01}   \, C^w_0 \, \Theta(t/2)  \int_{t/2}^t    (t-s)^{-\alpha}    \, ds
\\
&=& [ C^v_1  \, C^w_{01} +  C^v_{01}   \, C^w_0] \,   \Theta(t/2)  \frac{t^{1-\alpha} }{ 1-\alpha}, 
\eean
with the same obvious notation and  we have used again that $\Theta$ is a decaying function.  

\smallskip\noindent
{\sl Step 4. } Iterating $n := [\alpha]$ times steps 1 and 2, we get that $u^{(*n)}$  still satisfies estimate  {\bf (a)} and satisfies the estimate 
 {\bf (b)} for the exponent $\alpha-[\alpha] \in (0,1)$. We then conclude that   {\bf (c)} holds with $N := n+1$ and any $\tilde\zeta \in (0,\zeta^*/2)$ by using the third step.
 \end{proof}

\begin{proof}[Proof of Proposition~\ref{prop:diffuse-KRexistTER}.]
We classically establish that $\BB$ generates a positive semigroup $S_\BB$ in both spaces $X_i$ and we thus concentrate on the announced estimates.  On the one hand, proceeding as for the proof of \eqref{eq:diffusb:estimB}, we have 
\beqn\label{eq:diffusb:estimBL1}
\int (\BB  f) (\sign f) m  \le   - a'' \int   |f| ({\bf 1}_{B_1} + {\bf 1}_{B_1^c} |x|^{s+\gamma-2}) m,
\eeqn
for any  nice function $f$ and any weight function $m=m_a$, with  $m_a(x) := e^{a \langle x \rangle^\gamma}$, $a \in (a_1,a_2)$, $0 < a_1 < a_2 < 1/(\sqrt{2}\gamma)$,
where we define $a'' := a \gamma/2 - (a\gamma)^2$.  That exactly means that $\BB$ is  weakly dissipative in $L^1_m$ as defined in  \eqref{eq:defWeakDissip}. From the discussion in Section~\ref{subsect:AboutWeakDissip}, we deduce that $S_\BB$ is a semigroup of contractions and satisfies the associated decay estimate \eqref{eq:decaySBk}, \eqref{eq:decaySBkTheta}, and more precisely
\beqn\label{eq:decaySBL1}  
\| S_\BB(t) f \|_{L^1_{m_a}} \le   \|   f \|_{L^1_{m_{a}}}, \quad
\| S_\BB(t) f \|_{L^1_{m_a}} \le \Theta_\zeta(t)  \|   f \|_{L^1_{m_{a'}}},
\eeqn
for any   $a,a' \in (a_1,a_2)$,  $a < a'$,  $\zeta \in (0,\zeta_*)$, $\zeta_{*} := (a' - a)^{(2 - 2\gamma)/(2 - \gamma)}(a'\gamma(1 - a'\gamma))^{\gamma/(2-\gamma)}$.  We refer to   \cite[Lem.~2.1]{MR4265692} for details. Using that $\AA : L^1 \to L^1_m$ is bounded, that establishes (ii) in $X_0$. 

Similarly, starting from \eqref{eq:ch3:ch3FP:Lid2} and proceeding as in the proof of \eqref{eq:diffusb:estimB}, we get 
\beqn\label{eq:diffusb:estimBL2}
( \BB f, f)_{L^2_m} \le -\int |\nabla (fm)|^2  - a'' \int   f^2 ({\bf 1}_{B_1} + {\bf 1}_{B_1^c} |x|^{s+\gamma-2}) m^2,
\eeqn
for any  nice function $f$ and any weight function $m=m_a$ as above. Throwing away the first term at the RHS and arguing as we did in $L^1_m$, we obtain that $S_\BB$ satisfies
\beqn\label{eq:decaySBL2}  
\| S_\BB(t) f \|_{L^2_{m_a}} \le   \|   f \|_{L^2_{m_{a}}}, \quad
\| S_\BB(t) f \|_{L^2_{m_a}} \le \Theta_\zeta(t)  \|   f \|_{L^2_{m_{a'}}},
\eeqn
for any $a,a' \in (a_1,a_2)$.  Using that $\AA : L^2_{m_a} \to L^2_{m_a'}$ is bounded, that establishes (i) and (ii) in $X_1$. 

\smallskip
On the other hand, throwing away the second term at the RHS in \eqref{eq:diffusb:estimBL2}, for any trajectory $f_t = S_\BB(t) f_0$, $f_0$ in the domain of $\BB$ in $L^2_{m}$, we have 
$$
\frac{1 }{ 2} \frac{d}{ dt}  \int_{\R^d} f_t^2   \, m^2 dx \le  - \int_{\R^d} |\nabla ( f_t m) |^2 d x. 
$$
Using Nash's inequality which for some constant $C_N \in (0,\infty)$ stipulates that 
$$
  \int_{\R^d} g^2 dx  \le C_N \, \left( \int_{\R^d} |
  \nabla  g |^2 d x \right)^{\frac{d}{d+2}} \, \left( \int_{\R^d} |g|
  d x \right)^{\frac{4}{d+2}}, \quad \forall \, g, 
$$
with  $g:=f_tm$ and the first estimate in \eqref{eq:decaySBL2}, we deduce 
\begin{equation}\label{eq:Nash0}
F'(t) \le    - 2 \, C'_N \, F(t)^{-4/d} \, G(t)^{1+\frac{2}{ d}} \le  - 2 \, C'_N \, F(t)^{-4/d} \, G(0)^{1+\frac{2}{d}} ,
\end{equation}
with $C'_N := C_N^{-1-2/d}$ and where for brevity of notations we have set
$$
F(t) := \|f_t \|_{L^2(m)}^2, \qquad G(t) := \| f_t \|_{L^1(m)}.
$$
Integrate the differential inequality \eqref{eq:Nash0}, we find
$$
\| S_\BB(t) f_0 \|_{L^2_m}^2 \lesssim t^{-d/4}\,\| f_{0}\|_{L^1_m}, \quad \forall \, t > 0,
$$
and using that $\AA : L^1  \to L^1_{m}$, we next obtain
\beqn\label{eq:SBX0toX1}
S_\BB(t)\AA   t^{d/4} \in L^\infty(0,\infty;\BBB(X_0,X_1)).
\eeqn
Setting with $u(t) := S_\BB(t)\AA$, we see that $u$ satisfies {\bf (a)} in Lemma~\ref{lem:ch2:Tn1} thanks to (ii) in $X_0$ and $X_1$. Furthermore, 
 $u$ satisfies {\bf (b)} in Lemma~\ref{lem:ch2:Tn1} thanks to \eqref{eq:SBX0toX1}. Using Lemma~\ref{lem:ch2:Tn1}, we conclude that condition (iii) holds. 
 \end{proof}

We come back to the proof of \ref{HS3}. Gathering (i) and (ii) in $X_1$ in Proposition~\ref{prop:diffuse-KRexistTER}, we get that $(S_\BB \AA )^{(*\ell)}  * S_{\BB}  \in L^\infty_t (\BBB(X_1))$
for any $\ell \in \{ 0, \dots , N-1 \}$, $N := [d/4]+1$. Using that $\Theta \in L^1(0,\infty)$ and (iii) in Proposition~\ref{prop:diffuse-KRexistTER}, we deduce that 
$(S_\BB  \AA)^{(*N)} \in L^1(0,\infty;\BBB(X_0,X_1))$. 

\smallskip
We may now handle the existence part of the first eigenvalue problem. On the one hand, recalling \ref{H2}, we have $\LL^* \phi_0 \ge 0$ with $\phi_0 = 1$ so that the condition (i) in Theorem~\ref{theo:KRexistTER} holds. 
On the other hand, the condition  (ii) in Theorem~\ref{theo:KRexistTER} is an immediate consequence of  \ref{HS3} as emphasized in Remark~\ref{rem:KRexistTER}-(1). As a conclusion, the hypotheses of Theorem~\ref{theo:KRexistTER} are thus met, and we deduce 
that there exists $(\lambda_1,f_1) \in \R_+ \times L^2_{m+}$ solution to the first eigenvalue problem. Because the strong maximum principle \ref{H4} holds, we have $f_1 > 0$ on $\R^d$.

\smallskip
 
In order to prove the existence of a first positive eigenfunction for the dual problem, we argue in the following way. We start observing that we have the alternative: $\lambda_1 = 0$ or  $\lambda_1 >0$. 

- In the first case, we may argue as in Remark~\ref{rem:Uniuqness-supersolution}. We indeed have in the same time $\LL^* \phi_0 \ge 0$ and $\langle \LL^* \phi_0,f_1 \rangle = \langle  \phi_0, \LL f_1 \rangle = 0$, so that 
$\LL^*\phi_0 = 0$ because $f_1 > 0$. The function $\phi_1 := \phi_0$ is thus a solution to the first dual eigenvalue problem. 

- In the second case $\lambda_1 > 0$,  we may argue as in the case $\gamma=1$ above. 
On the one hand, the operator $\RR_\BB(\alpha)$ is uniformly bounded in $L^2_m$ for any $\alpha \ge \kappa_\BB := \lambda_1/2 > 0$ and on the other hand the operator $\WW(\alpha) := \RR_\BB(\alpha) \AA : L^2_m \to H^1_m \cap L^2_{m^\sharp}$
 is  uniformly bounded for any $\alpha \ge \kappa_\BB$ with $m = o(m^\sharp)$, so that  $H^1_m \cap L^2_{m^\sharp} \subset L^2_m$ is compact. 
  We may thus apply Theorem~\ref{theo:exist1-KRexistence} and we conclude to the existence of a solution $(\lambda_1',f'_1,\phi'_1)$ to the eigentriplet problem.

\smallskip

 The conditions \ref{H4} and \ref{H5} being true in a general situation as well as the conclusions of Lemma~\ref{lem:Section7-BoccardoG},  
 as an intermediate conclusion, we have established under the general condition $\gamma > 0$ in \eqref{eq:diffusionWithDrift-defb&gamma} that
 yet the  same conclusions as in section~\ref{subsection:diffusionStrongPotential} hold true.

 \smallskip\smallskip  
 {\bf Quantitative stability.}  %
 We now establish a quantitative stability estimate using the Doblin-Harris approach presented in Section~\ref{sec:QuantitativeStabilityKR} and yet used in the case
of a bounded domain in Section~\ref{subsec:diffusion-domain}. We first consider the more difficult case $\gamma \in (0,1)$, so that $\lambda_1 \ge \kappa_0 = 0$, and then explain the modifications to be made 
in order to deal with the case $\gamma \ge 1$.  As explained just above, $\lambda_1 = 0$ corresponds to the conservative case $(\lambda_1,\phi_1) = (0,1)$ which has been considered in \cite{MR4265692}.
We thus focus on the case $\lambda_1 >0$ for which an adapted version of Theorem~\ref{theo:stampacchia1} already imply the exponential asymptotic stability  \ref{E31} in $L^2_m$ with {\bf non constructive rate}. We do not develop further this argument but rather establish a {\bf a constructive sub-exponential asymptotic stability}. 


\smallskip
{\bf Step1 - Lyapunov condition.} We take $m = e^{a |x|^\gamma}$ with $0 < a < \gamma^{-1}$. From Lemma~\ref{lem:diffus-bconf:identityL} or a direct computation, we have 
\bean
\LL^* m
&=& \Delta m + (c-\Div b) m - b \cdot \nabla m 
\\
&\le& (C \langle x \rangle^{\gamma-2} + |c| - 2 a^* \langle x \rangle^{2\gamma-2}) m 
\\
&\le& C_0 {\bf 1}_{B_{\varrho_0}}  - a^* \langle x \rangle^{2\gamma-2} m,
\eean
for three positive constants $C = C(d)$,  $C_0 = C_0(c, C,a,\gamma)$, $\varrho_0 = \varrho_0(c,C,a,\gamma)$ and with  $a^* := (a\gamma -(a\gamma)^2)/2 > 0$.
We now set $m_1 := m$ and $m_0 := a^* \langle x \rangle^{2\gamma-2} m$. We fix $T > 0$ and for $0 \le f_0 \in L^1_m$, we denote $f_t := \widetilde S_\LL(t) f_0$. 
Recording that $\lambda_1 \ge 0$ and using the above pointwise estimate, we deduce 
\beqn\label{eq:dissipDrift-Lyapunovm}
\int f_T m_1 + \int_0^T \!\!\int f_t m_0 dt  \le \int f_0 m_1 + C_1 \int_0^T\!\! \int_{B_{\varrho_0}} f_t dt.
\eeqn
Because the same kind of pointwise estimate holds for $\LL^* m_0$, we have 
$$
\int f_T m_0   \le \int f_t m_0 + C_1 \int_t^T\!\! \int_{B_{{\varrho_0}}} f_s ds
$$
and integrating in time, we get 
$$
T \int f_T m_0   \le  \int_0^T\!\!  \int f_t m_0 dt + C_0 \int_0^T s \int_{B_{\varrho_0}} f_s ds.
$$
Coming back to the first estimate, we deduce 
\beqn\label{eq:dissipDrift-Lyapunovm1m0}
\int f_T m_1 +  T \int f_T m_0   \le \int f_0 m_1 + (C_1+C_0 T) \int_0^T\!\! \int_{B_{\varrho_0}} f_t dt.
\eeqn

\smallskip
{\bf Step 2 - Pointwise estimates on $\phi_1$.} We define $\BB := \LL - C_0 \chi_{\varrho_0}$ which is the generator of a positive semigroup of contraction in $L^1_m$ because of the above discussion. 
For $\lambda > 0$, $0 \le g \in  L^1_m$ and $0 \le f \in L^1_m$ the solution to $(\lambda - \BB) f = g$, we compute 
\bean
\int g m = \int f  (\lambda - \BB^*) m 
\ge \int f  (\lambda m + m_0) \ge  \int f  m_0, 
\eean
from what we deduce 
$$
\| \RR_\BB(\lambda) f \|_{L^1_{m_0}} \le \|   f \|_{L^1_{m}} , \quad \forall \, f \in L^1_m.
$$
Now, we consider two weight functions $m_1$ and $m_3$ with $m_i := e^{a_i |x|^\gamma}$, $0 < a_1 < a_3 < \gamma^{-1}$, we denote $m_0 := a_1^* \langle x \rangle^{2\gamma-2} m_1$ and we compute
\bean
\| \AA \RR_\BB(\lambda) f \|_{L^1_{m_3}} 
&\le& C_0 \|  \RR_\BB(\lambda) f \|_{L^1_{m_3}(B_{2\varrho_0})} 
\\
&\le&  C_1 \|   \RR_\BB(\lambda) f \|_{L^1_{m_0}} \le C_1 \|   f \|_{L^1_{m_1}}.
\eean
By duality, we obtain 
\beqn\label{eq:dissipDrift-ARBLinfty}
\| \RR_{\BB^*}   (\lambda) \phi \|_{L^\infty_{m_1^{-1}}} \le  \|   \phi \|_{L^\infty_{m_0^{-1}}}
\ \hbox{ and } \
\| \RR_{\BB^*} \AA^* (\lambda) \phi \|_{L^\infty_{m_1^{-1}}} \le C_1\|   \phi \|_{L^\infty_{m_3^{-1}}}, 
\eeqn
for any $\lambda > 0$ and $\phi \in L^\infty_{m_0^{-1}}$. We also deduce from Proposition~\ref{prop:diffuse-KRexistTER} the regularization estimate $(\AA^* \RR_{\BB^*})^N : L^2_{m^{-1}} \to L^\infty$. 
Let us now consider $0 \le \phi_1 \in L^2_{m_1^{-1}}$ the first eigenvector for the dual problem built in the preceding paragraph. 
From the eigenvalue equation 
$$
\BB^* \phi_1 + \AA^* \phi_1 = \LL^* \phi_1 = \lambda_1 \phi_1,
$$
we deduce that $\phi_1 = (\RR_{\BB^*} \AA^* )  \phi_1$, and iterating 
$$
\phi_1 = (\RR_{\BB^*} \AA^* )^{N+1} \phi_1 = \RR_{\BB^*} (\AA^* \RR_{\BB^*} )^{N } \AA^* \phi_1. 
$$
From the above regularization estimate and the first estimate in \eqref{eq:dissipDrift-ARBLinfty}, we thus deduce that $\phi_1 \in L^\infty_{m_1^{-1}}$. Moreover, normalizing $\phi_1$ and using  the second estimate in \eqref{eq:dissipDrift-ARBLinfty}, we may obtain 
\beqn\label{eq:dissip-drift}
\| \phi_1 \|_{L^\infty_{{m_3}^{-1}}} = 1 \ \hbox{ and } \
\|   \phi \|_{L^\infty_{m_1^{-1}}} \le C_1.
\eeqn
We deduce
\bean
1 
&=& \max \Bigl( \sup_{B_{\varrho_2}} \frac{|\phi_1| }{ m_3}, \sup_{B^c_{\varrho_2}} \frac{|\phi_1| }{ m_3} \Bigr)
\\
&\le& \max \Bigl( \sup_{B_{\varrho_2}} \frac{|\phi_1| }{ m_3}, C_{0,1}  \sup_{B^c_{\varrho_2}} \frac{m_1 }{ m_3} \Bigr),
\eean
so that $\sup_{B_{{\varrho_2}}} \frac{|\phi_1| }{ m_3} =1$ by choosing ${\varrho_2} := \max(\varrho_0,{\varrho_1})$ with $C_{1} e^{(a_1-a_3) \varrho_1^\gamma} = 1$. 
As a consequence, there exists $x_0 \in B_{R}$ such that $\phi_1(x_0) \ge 1$. On the other hand, using standard regularity result for elliptic equation in the ball $B_{2R}$, we obtain that 
$\phi_1 \in C^{0,1}(B_R) \cap W^{2,p}(B_R)$ for any $p \in [1,\infty)$ with constructive bound. Making use next of the Harnack inequality  as at the end of Section~\ref{subsec:diffusion-domain} or using barrier functions as in  in the proof of \cite[Lem.~6.2]{MR4265692}, we classically deduce that 
\beqn\label{eq:dissip-drift-phi1lowerB}
\phi_1 \ge z_\varrho {\bf 1}_{B_\varrho}, \quad \forall \, \varrho > 0, 
\eeqn
for a constructive constant $z_\varrho >0$ (where we emphasize here and below the $\phi_1$ always denote the normalized  by \eqref{eq:dissip-drift} dual eigenvector). 

%

\smallskip
{\bf Step 3 - Doblin-Harris estimate.} We fix $T > 0$ (for instance $T:=1$) and $A > 0$ arbitrary. For  $0 \le f_0 \in L^1_m$ such that $\| f_0 \|_{L^1_m} \le A \| f_0 \|_{L^1_{\phi_1}}$, we denote $f_t := \widetilde S_\LL(t) f_0$. 
On the one hand, we have 
$$
\int f_t \phi_1 = 
\int f_0 \phi_1 \quad  \hbox{and}\quad
\int f_t m \le C_T \int f_0 m,
$$
for any $t \in [0,T]$, the second estimate being an immediate consequence of \eqref{eq:dissipDrift-Lyapunovm}. On the other hand, we define $\eps(r) := \sup_{|x| \ge r} (m(x)/\phi_1(x))$ and 
we compute 
\bean
\int_{B_\rho} f_t \phi_1 
&=& 
\int  f_t \phi_1 - \int_{B^c_\rho}  f_t \phi_1 
\ge
\int  f_t \phi_1 - \eps(\rho) \int  f_t m 
\\
&\ge&\int  f_0 \phi_1 - C_T \eps(\rho) \int  f_0 m 
\ge \Bigl(  1  - \frac{C_T }{ A} \eps (\rho)  \Bigr)    \int  f_0 \phi_1  \ge \frac1{2}   \int  f_0 \phi_1 ,
\eean
for any $t \in (0,T)$, by choosing  $\rho := \rho(T,A) >0$ large enough. 
 In particular, there exists $x_0(t) \in B_\rho$ such that 
$$
f(t,x_0(t)) \ge  \vartheta := \frac{1 }{ 2} \frac{1}{ \| \phi_1 \|_{L^1(B_\rho)}}   \int  f_0 \phi_1.
$$
Next, arguing exactly as in Section~\ref{subsec:diffusion-domain} or as in the proof of \cite[Lem.~6.2]{MR4265692}, we deduce 
\beqn\label{eq:dissip-drift-DoblinHarris}
\widetilde S_T f_0 \ge \eta_{T,A} {\bf 1}_{B_1}  [ f_0 ]_{L^1_{\phi_1}},
\eeqn
for some constructive constant $\eta_{T,A} > 0$. 

\smallskip

{ \Blue 
  In conclusion, using Theorem~\ref{theo:Harris} or Theorem~\ref{theo:HarrisSubgeo1},  we have established the following result.

\begin{theo}
Consider the elliptic operator \eqref{eq:LLf=Rd} in the whole space and assume that the coefficients $b$ and $c$ satisfy \eqref{eq:diffusionWithDrift-defb&gamma}, \eqref{eq:dissipDrift-HYPc1}, \eqref{eq:dissipDrift-HYPc2}. Then the associated semigroup satisfies the quantitative estimate    \ref{E31intro} when $\gamma \ge 1$ or \ref{E32intro} when $\gamma \in (0,1)$ in the Lebesgue spaces $L^2_m$, with the exponential weight function $m$ as defined above, in particular \ref{S1}, \ref{S2} and \ref{S33} hold true. 
\end{theo}

We are not aware of any  result of this kind in a non conservative framework ($c \not\equiv \Div b$). Our result generalizes to a non conservative framework the previous known results in a conservative framework, see for instance \cite{MR2386063,MR4265692} and the references therein.  
}

\bigskip

 \bigskip

\section{Transport equations}
\label{sec:Transport}


%
%
%
%
%
%
%
%
%
%

%

The main aim of this part is to investigate the long time asymptotic of the solutions to the transport equation
\beqn\label{eq:Transport-evolEqWithKernel}
\partial_t f  + \Div_y(a  f) = \KKK[f] - K f   \quad\hbox{in}\quad (0,\infty) \times \OO, 
\eeqn
on the function $f = f(t,y)$, $t \ge 0$, $y \in \OO$, with $\OO \subset \R^D$, $D \ge 1$, a smooth open connected set.
We assume that $a = a(y)$, $a : \OO \to \R^D$, $K= K(y)$, $K : \OO \to \R_+$  and that the collision operator $\KKK$ is linear and defined by 
\beqn\label{eq:secTeq-defKKg}
\KKK [g](y) := \int_{\OO} k \,  g_*  \, dy_* ,
\eeqn
for some kernel $k : \OO \times \OO \to \R_+$ and for any (conveniently) bounded function $g  : \OO \to \R$. Here and below, we use the common shorthands
$$
g_* := g(y_*), \quad k := k(y,y_*), \quad k_* := k(y_*,y).
$$
When $ \OO \not=\R^D$, 
the equation is complemented with a boundary condition which imposes the value of the trace $\gamma_- f$ of $f$ on the incoming subsets of the boundary  and takes the form
\beqn\label{eq:Transport-BdaryCond}
(\gamma_- f )(t,y) = 
 \RRR  [f(t,\cdot),\gpf(t,.) ] (y) \,\, \hbox{ on }  \,\,(0,\infty) \times \Sigma_-. 
\eeqn
Let us explain the meaning of the different terms involved in \eqref{eq:Transport-BdaryCond}. 
We denote by $\Sigma := \partial\OO$ the boundary set, by $d\sigma_{\! y}$ the Lebesgue measure on $\Sigma$, by $n : \Sigma \to \Sp^{D-1}$ the normal outward vector field,
we write $n=n_y=n(y)$, and by  $\Sigma_-$ the incoming,  $\Sigma_+$ the outgoing and $\Sigma_0$ the singular subsets of the boundary defined by
\bean
\Sigma_\pm := \{ y \in \Sigma; \ \pm a(y) \cdot  n_y  > 0 \}, \quad
\Sigma_0 := \{ y \in \Sigma;  \  a(y) \cdot  n_y = 0 \}. 
\eean
We denote  $\gamma f = f_{|(0,\infty) \times \Sigma}$ the trace of $f$ and $\gamma_\pm f := {\bf 1}_{(0,\infty) \times \Sigma_\pm} \, \gamma f$ the trace restrictions on the 
 incoming  and outgoing sets. We then assume that the boundary operator $\RRR$ splits into two pieces $\RRR[g,h] = \RRR_\OO[g] + \RRR_\Sigma[h]$, where
\beqn\label{eq:Transport-BdaryCond2}
\RRR_\OO [g](y) = \int_{\OO} g(y_*) r_\OO(y,dy_*),
\quad
\RRR_\Sigma [h](y) = \int_{\Sigma_+} h(y_*) r_\Sigma(y,dy_*),
\eeqn
for a domain transition kernel  $r_\OO : \Sigma_- \times \BBB_\OO \to [0,\infty]$, a boundary transition kernel $r_\Sigma : \Sigma_- \times \BBB_{\Sigma_+} \to [0,\infty]$
and  for any (conveniently) bounded functions $g  : \OO \to \R$ and $h  : \Sigma_+ \to \R$, where $\BBB_E$ stands for the set of Borel subsets of $E$.


\smallskip 
In the next sections we will first consider the trace problem for a general force field $a$ and next the well-posedness for the transport equation 
with given inflow at the boundary and with reflection condition at the boundary.  We will also revisit the characteristic method for general force field $a$. 
We will next consider the Krein-Rutman problem still for a general force field $a$, but making strong simplification assumptions on the kernel operators $\KKK$ and $\RRR$. 
We will finally explain how the classical age structured equation falls into the present framework.  
We will come back to more specific physical situations concerning the growth-fragmentation equation and the kinetic relaxation equation with more general and physically relevant 
hypothesis on the kernel in parts~\ref{part:application3:GF} and \ref{part:application4:Kequation}.

 \medskip
\subsection{The trace problem}
\label{subsec:TranspTrace}

\

In this section, we are concerned with the trace problem associated to a (mainly stationary) transport equation for a general  vector field $a :   \OO \to \R^D$ for which we only assume 
\beqn\label{eq:Transport-Hyp-a}
a \in \Wloc^{1,1}(\bar\OO),
\eeqn
where we recall that $\OO \subset \R^D$, $D \ge 1$, is a smooth open connected set. The regularity needed on the domain is  formulated in the following way: we assume that there exists $n : \OO \to \R^D$, 
$y \mapsto n(y)$ a vector field belonging to $W^{1,\infty}(\OO)$ and which coincides with the previously defined unit outgoing normal vector field on $\Sigma$ and satisfies $\| n \|_{L^\infty} = 1$. In that situation, it is well-known that the above vector field is the restriction of a vector field $a \in \Wloc^{1,1}(\R^D)$ (where we abuse notations denoting the restriction and the extension in the same way). 
We also consider the associated differential equation 
\beqn\label{eq:Transport-characteristics}
\frac{d Y }{ dt} = a(Y), \qquad Y(0) = y , 
\eeqn
and then define the characteristic flow $Y_t=Y(t,y)$, for any $y \in \OO$, which is the solution to \eqref{eq:Transport-characteristics} 
on a maximal time interval $(t_-(y) , t_+(y))$ where $t_-(y) < 0 < t_+(y)$ are defined by $t_- := - \tb $ and $t_+ :=  t_{\bf f}$, the backward exit time is defined by 
\beqn\label{eq:Transport-def-tb}
\tb(y)  := \sup\{\tau  >0 ; \ Y_{-t} (y) \in \OO, \  \forall  \,t \in [0,\tau] \} \in   (0,+\infty]
\eeqn
and the   forward exit time is defined by 
\beqn
 \label{eq:Transport-tf}
t_{\bf f}(y) := \sup\{\tau > 0; \ Y_t(y) \in \OO, \ \forall \, t \in [0,\tau] \} \in  (0,+\infty]. 
\eeqn
%
The   real number $t_{\ell t} (y) := \tb (y) + t_{\bf f}(y) \in (0,\infty]$ corresponds to the {\it``life time''} of the characteristic flow in  $\OO$ going by $y$.  
 The construction of the flow $(Y_t)$ is classical  when $a$ is a Lipschitz function
and we refer to \cite[Thm.~II.3]{MR1022305} for a more general situation which corresponds to the assumptions we will make in the present work (see also  Lemma~\ref{lem:Transport-Uniqueness} below). 

\smallskip
For a solution $g :   \OO \to \R$ to the transport equation 
\beqn\label{boundary:trace:TE1}
a \cdot \nabla_y g = G \quad\hbox{in}\quad \OO, 
\eeqn
for a given source term $G :  \OO \to \R$, we wish to define the trace $\gamma g$ of $g$ on the boundary set $  \Sigma$. Similarly,  for a solution $g :    (0,T)  \times \OO \to \R$, $T \in (0,+\infty]$,  to the transport equation 
\beqn\label{eq:transport:EvolAVECSource}
\partial_t g + a  \cdot \nabla_y g = G  \quad\hbox{in}\quad (0,T) \times \OO, 
\eeqn
for a given source term $G :  (0,T) \times \OO \to \R$, we wish to define the trace $\gamma g$ of $g$ on the boundary set $(0,T) \times  \Sigma$. It is worth emphasizing that the trace will be in fact only
defined out of the singular set $\Sigma_0$ and thus only on the boundary set $\Sigma \backslash \Sigma_0$. 

\smallskip
We start by recalling several possible definitions of the trace of a function $g$ satisfying \eqref{boundary:trace:TE1}  when 
\beqn\label{hyp:g&G}
a \in \Wloc^{1,s}(\bar\OO), \quad g \in \Lloc^p(\bar\OO), \quad G \in \Lloc^q(\bar\OO), \quad s,p,q \in [1,\infty].
\eeqn 
Here and below, we denote by $L(E)$ the Lebesgue space of measurable functions $g : E  \to \bar\R := [-\infty,+\infty]$ with typically $E = \OO$ or $E \subset \Sigma$, 
and by $L^0(E)=L^0(E,\mu) \subset L(E)$ the subset of almost everywhere finite measurable functions  on a measurable space $(E,\AA,\mu)$.

\begin{defin}
We say that a  function $g$ on $\OO$ satisfying \eqref{boundary:trace:TE1} and \eqref{hyp:g&G} admits a trace if one of the following assertions holds true:

\begin{itemize}
\item  {\bf Extension of the restriction on the boundary.} 
There exists  $\gamma g \in \Lloc^r(\Sigma\backslash\Sigma_0)$, $r \in [1,\infty]$, such that  
$$
g_{n|\Sigma\backslash\Sigma_0} \to \gamma g \quad\hbox{in}\quad \Lloc^r(\Sigma\backslash\Sigma_0)
$$
for any sequence $(g_n)$ satisfying 
\beqn\label{eq:Transport-def1trace}
g_n \in C^1_c(\bar\OO), \qquad
g_n \to g \quad\hbox{in}\quad \Lloc^p(\bar\OO), \qquad
a(y)  \cdot \nabla_y g_n \to G \quad\hbox{in}\quad  \Lloc^q(\bar\OO). 
\eeqn

\item {\bf Characteristics.}
There exists a measurable function   $\gamma g$ on $ \Sigma\backslash\Sigma_0$ such that for a.e. $y \in \OO$ satisfying
  $t_-(y) > - \infty$, there holds
\beqn\label{eq:Transport-def2trace}
 g(y) = \gamma g(Y(t_-(y),y)) + \int_{t_-(y)}^0 G(Y(t,y)) \, dt  ,  
\eeqn
and for a.e. $y \in \OO$ satisfying $t_+(y) < \infty$, there holds
\beqn\label{eq:Transport-def2traceBis}
g(y) = \gamma g (Y(t_{\bf +}(y),y)) - \int_0^{t_+(y)} G(Y(t,y)) \, dt.  
\eeqn

\item {\bf Green formula.} There exists  $\gamma g \in \Lloc^r (\Sigma \backslash \Sigma_0)$, $r \in [1,\infty]$,  such that
\beqn\label{eq:Transport-def3trace}
\int_\OO \bigl( G \, \varphi + g \, \Div(a\varphi)  \bigr) \, dy =
\int_\Sigma \gamma g \, \varphi \, a  \cdot n  \, d\sigma_{\! y},
\eeqn
for any $ \varphi \in C^1_c(\bar\OO \backslash \Sigma_0)$.  

\smallskip
\item {\bf Renormalized Green formula.} There exists  a measurable function $\gamma g$ on $\Sigma \backslash \Sigma_0$ such that 
\beqn\label{eq:Transport-def4trace}
\int_\OO \bigl( \beta'(g) \, G \, \varphi  + \beta(g) \, \Div(a\varphi)  \bigr) \, dy =
\int_\Sigma \beta(\gamma g) \, \varphi \, a \cdot n \, d\sigma_{\! y},
\eeqn
for any $ \varphi \in C^1_c(\bar\OO)$ and any  $\beta \in C^1(\R)$ such that $\beta' \in L^\infty(\R)$. 
\end{itemize}

\end{defin}

\begin{rem}  

(1) In order that the first definition makes sense, we implicitly assume that there exists at least one sequence $(g_n)$ which satisfies \eqref{eq:Transport-def1trace}. 
That last fact corresponds to the density of $C^1_c(\bar\OO)$ in the Sobolev space $\{ g \in L^p(\OO); \ a(y)  \cdot \nabla_y g \in  L^q(\OO)\}$, which is true as we will see in Lemma~\ref{lem:transport-lemma51} below under the regularity assumptions made on $a$ and $\OO$.   It is worth emphasizing that the last convergence in \eqref{eq:Transport-def1trace} may require additional integrability assumption, typically $a \in W^{1,s}(\OO)$ with $1/r \ge 1/p+1/s$. Such a definition has been introduced by Bardos in \cite{MR274925} for a $C^1$ vector field $a$.
It is also the point of view  adopted by  Cessenat in \cite{MR772106,MR777741} in the case of the neutron operator, see also \cite[chap.~XXI]{MR1295030} or Agoshkov \cite{MR839640,MR843433,MR1087296}. 
 
\smallskip
(2) In order that the second definition makes sense, we implicitly assume that the set of points $y \in \OO$
 such that the characteristic $Y_t(y)$ hits the boundary on $\Sigma_0$ has zero measure in $\OO$.
 It is indeed the case thanks to the Sard theorem under enough regularity assumption on $a$ and $\OO$, see \cite[Prop.~2.3]{MR274925}. 
 It is worth emphasizing that what we really need in order to write \eqref{eq:Transport-def2trace} and \eqref{eq:Transport-def2traceBis}
 is that  $t \mapsto G(Y(t,y)) \in L^1(t_-(y),t_+(y))$ for a.e. $y \in \OO$. We also mention that this characteristics description leads to a  layer cake formula linking the
 integral of a function on the domain to the integral of its trace on the boundary. Such a definition has been widely  used in kinetic theory for constructing DiPerna-Lions renormalized solution, see \cite{MR1014927,MR1179690,MR1301931} and the references therein. For the  classical kinetic operator this trace approach is developed by Arkeryd, Cercignani and co-authors in \cite{MR1136615,MR1132764,MR1245073,MR1629463} 
 while for more general (but still regular) vector fields, the approach has been developed in  \cite{MR274925,MR882376,MR872231,MR896904} and more recently by Arlotti et al. 
in \cite{MR2363946,MR2565274,MR2781917,MR4023819,MR4023826}. 

\smallskip
(3)  In order that the third definition makes sense, we need that $a \cdot n \, \gamma g \in \Lloc^1(\Sigma)$ and $p \ge s'$. 
In some situation, this third definition is in some sense the weakest: it makes sense also when $\gamma g \in M^1_{\rm  loc}(\Sigma \backslash \Sigma_0)$ for instance and can be relevant under the weak assumption $a, \Div a  \in \Lloc^{p'}(\bar\OO)$ as it is the case in the early works on weak solution to the Vlasov-Poisson equation in  \cite{MR1071632,MR1224079,MR1274153,MR1338454}.
 It is also easier to handle than the two first definitions because of the way it connects the function $g$ and its trace.   

\smallskip
(4) We will adopt the last definition which extends up to the boundary the renormalization technique introduced in \cite{MR1022305}. It is more general and adapted to the weak regularity assumption made on the vector field $a$ than the two first definitions and we recover the third definition by just letting $\beta(s) \to s$ when the conditions of integrability make the limit well defined. 
Such a kind of definition has been introduced in \cite{MR1776840,MR2721875} for kinetic equations and in \cite{MR2150445,MR2188582} for transport equations. 
\end{rem}

We start with a trace result in a $L^\infty$ framework. We denote by $C^1_{\rm pw}(\R) $ the space of continuous functions $\beta:\R \to\R$ with piecewise continuous derivative.

\begin{theo}\label{theo:traceLinfty1} Assume that $g \in L^\infty(\OO)$, $a \in \Wloc^{1,1}(\bar\OO) $ and $G \in \Lloc ^1(\bar \OO)$
satisfy the  transport equation \eqref{boundary:trace:TE1} in the distributional sense. Then, there exists a unique function
$$
\gamma g \in L^\infty( \Sigma \backslash \Sigma_0; \, d\sigma_{\! y}  ), 
\quad \| \gamma g \|_{L^\infty} \le \| g \|_{L^\infty}, 
$$
which satisfies the   renormalized Green formula  
\beqn\label{eq:Transport-def4traceBIS}
\int_\OO \bigl( \beta'(g) \, G \, \varphi  + \beta(g) \, \Div(a\varphi)  \bigr) \, dy =
\int_\Sigma \beta(\gamma g) \, \varphi \, a \cdot n \, d\sigma_{\! y},
\eeqn
for any $ \varphi \in C^1_c(\bar\OO)$ and any $\beta \in C^1_{\rm pw}(\R) $. 
As a consequence, renormalization and trace operations commute: 
\beqn\label{eq:Transport-traceCrenormalization}
\gamma \, \beta(g) = \beta (\gamma \, g), \quad \forall \, \beta  \in C^1_{\rm pw}(\R). 
\eeqn
\end{theo}

\begin{rem}
(1) Because of the very general  assumption \eqref{eq:Transport-Hyp-a} made on the vector field $a: \OO \to \R^D$ which is exactly the one made in the DiPerna-Lions  theory for transport equation in the whole space developed in \cite{MR1022305}, the above trace result slightly  improves 
the similar trace result established by Boyer in 
\cite[Thm.~3.1]{MR2150445}, where an additional assumption $a \cdot n \in L^\zeta(\partial\OO)$, $\zeta > 1$, is made. 

\smallskip
(2) 
 An alternative approach has been developed by Ambrosio  and co-authors by assuming weaker bound on $D a$ but stronger bound on $a$. 
More precisely, denoting by  $\MM_\infty$ the set of vector fields $a \in L^\infty(\OO)$ such that $\Div a \in M^1(\OO)$,  
it is established in \cite[Prop.~3.2]{MR2188582} that  there exists a linear and bounded mapping  $\Trace :  \MM_\infty(\OO) \to L^\infty(\partial\OO)$  such that ${\rm Tr} a = n \cdot a_{|\partial\OO}$ when $a \in C^1(\bar\OO)$. 
The proof relies on  Ambrosio's extension to a BV framework in  \cite{MR2096794} of the famous DiPerna-Lions 
improvement   \cite[Lem.~II.1]{MR1022305}
of  Freidrichs' type Lemma on the estimate of the commutator between directional derivative and convolution (see Lemma~\ref{lem:transport-lemma51} below). 
Moreover, it is also established in  \cite{MR2188582} (see in particular \cite[Thm.~4.2]{MR2188582}) that 
$$
\Trace (a \beta (g)) = \beta \Bigl( \frac{\Trace (a  g) }{ \Trace a} \Bigr)  \Trace (a) , \quad \forall \, \beta \in C^1(\R), 
$$
for any   $a \in BV(\OO) \cap L^\infty(\OO)$ and $g \in L^\infty(\OO)$ such that $ag \in \MM_\infty$. The above formula is then nothing but \eqref{eq:Transport-traceCrenormalization} when  $a \in W^{1,1}(\OO) \cap L^\infty(\OO)$.
 
\end{rem}


Before coming to the proof of Theorem~\ref{theo:traceLinfty1}, we state 
one technical but fundamental result. 
We define the mollifier $(\rho_\eps)_{\eps>0}$  by
\beqn\label{eq:Transport-mollifier}
\rho_\eps(z) = \tfrac1{\eps^d} \, \rho(z/\eps) , 
\quad 0 \le \rho \in \DD(\R^d), 
\quad \hbox{supp} \, \rho \subset B_1, \quad \int_{\R^N} \rho(z) \, dz = 1,
\eeqn
and for any $u \in \Lloc^1(\bar\OO)$, $v_\eps \in C_c(\R^D)$, supp$\, v_\eps \subset B_\eps$,  we introduce the convolution-translation function $u *_\eps v_\eps$   defined
by
\beqn\label{eq:convolution-for-trace}
(u *_\eps v_\eps)(y) 
:= \int_{\OO} u(z) \, v_\eps (y - 2\eps \, n(y) - z ) \, dz.
\eeqn
 
\begin{lem}\label{lem:transport-lemma51}
For $g \in \Lloc^p(\bar\OO)$, $p  \in [1,\infty]$, $a \in \Wloc^{1,p'}(\bar\OO)$ and $G \in \Lloc^1(\bar\OO)$  satisfying \eqref{boundary:trace:TE1}
in the distributional sense, 
the sequence $(g_\eps)$ 
defined by $g_\eps := g*_\eps \rho_\eps$ satisfies 
$$
g_\eps \in \Wloc^{1,\infty}(\bar\OO), \quad G_\eps := a \cdot \nabla g_\eps \to a \cdot \nabla g \ \hbox{in}\  \Lloc^{1}(\bar\OO), 
$$
as $\eps \to 0$, and 
\bean 
&&g_\eps \to g \ \hbox{ in } \ \Lloc^p(\bar\OO), \quad \hbox{if} \ p < \infty, 
\\
&&g_\eps \to g \ \hbox{ in } \ \Lloc^1(\bar\OO), \quad  (g_\eps) \ \hbox{bounded in} \ \Lloc^\infty(\bar\OO) , \quad \hbox{if} \ p = \infty.
\eean
 \end{lem}
 
We skip the  proof of  Lemma~\ref{lem:transport-lemma51} since it follows by just repeating the proofs of \cite[Lem.~II.1]{MR1022305}, \cite[Lem.~1]{MR1765137} or \cite[Lem.~3.1]{MR2150445}.

\begin{proof}[Proof of Theorem~\ref{theo:traceLinfty1}.]
Let us fix $\chi \in \DD(\bar\OO)$ such that $0 \le \chi \le 1$ and denote $R > 0$ a real number such that $\hbox{supp}\,\chi \subset B_R$. 
We observe that $ \chi \hbox{sign} ( a \cdot n) \in L^1(\Sigma)$. 
From Gagliardo trace theorem \cite[Teor.~1.II]{MR102739}, there exists $\psi \in W^{1,1}(\OO)$  such that $\gamma \psi =  \chi \hbox{sign} ( a \cdot n)$ and supp$\, \psi \subset B_R$. 
Denoting $T_1 : \R \to [-1,1]$ the truncation function which is odd and is defined by $T_1(\sigma) = \sigma \wedge 1$ for any $\sigma \ge 0$, we see that $\gamma T_1(\psi) = T_1(\gamma \psi) = \gamma \psi$, and thus we may assume $\psi \in L^\infty(\OO)$ up to replacing $\psi$ by $T_1(\psi)$. 
As a consequence, there exists a sequence  $(\psi_k)$ of $W^{1,\infty}(\OO)$ such that $\psi_k \to \psi$ in $W^{1,1}(\OO)$, with $(\psi_k)$ bounded in  $L^\infty(\OO)$, $\hbox{supp} \psi_k \subset B_R$, and $\gamma \psi_k \to \chi  \, \hbox{sign} (a \cdot n)$ in $  L^1(\Sigma)$, with $(\gamma \psi_k)$ bounded in $L^\infty(\Sigma)$. 

\smallskip
Let us then consider the sequences $(g_\eps)$ and  $(G_\eps)$ defined in Lemma~\ref{lem:transport-lemma51}.  The classical Green formula for Lipschitz functions writes
\bean
&&\int_\Sigma (g_{\eps|\Sigma} - g_{\eps'|\Sigma})^2 \, |a \cdot n| \, \chi \, d\sigma_{\! y} 
\\
&&= \int_\Sigma (g_{\eps|\Sigma} - g_{\eps'|\Sigma})^2 \, a \cdot n \, \psi_k \, d\sigma_{\! y} 
+ \int_\Sigma(g_{\eps|\Sigma} - g_{\eps'|\Sigma})^2 \, [|a \cdot n| \, \chi  - a \cdot n \,  \psi_k ] \, d\sigma_{\! y} \\
&&= \int_\OO [ 2 \, ( g_\eps - g_{\eps'}) \, (G_\eps - G_{\eps'})  \,  \psi_k \, dy 
+ ( g_\eps - g_{\eps'})^2  \Div (a  \psi_k) ] \, dy \\
&& \quad + \int_\Sigma (g_{\eps|\Sigma} - g_{\eps'|\Sigma})^2 \, [|a \cdot n| \, \chi  - a \cdot n \,  \psi_k ] \, d\sigma_{\! y}
\\
&&\le 4 \| \psi_k \|_{L^\infty} \| g \|_{L^\infty} \| G_\eps - G_{\eps'} \|_{L^1(B_R)}
+   \| \psi_k \|_{W^{1,\infty}}  \int_{B_R} ( |a| +|\Div a|)  ( g_\eps - g_{\eps'})^2 \, dy \\
&& \quad +  \, 2 \| g \|_{L^\infty}^2 \| (a \cdot n) \gamma \psi_k - \chi  \, |a \cdot n| \|_{L^1(\Sigma)} , 
\eean
for any $\eps>0$ and $k \ge 1$.
We deduce that  $(g_{\eps|\Sigma})$ is a Cauchy sequence in $L^2(|a\cdot n| \, \chi \, d\sigma)$. 
From the fact that  $( g_\eps)$ is  bounded in $L^\infty(\OO)$, we deduce 
that the sequence $(\gamma g_\eps)$ is also bounded in $L^\infty(\Sigma)$. As a consequence, there exists a function $\gamma g \in L^\infty(\Sigma)$ such that 
$\gamma g_\eps \to \gamma g$ in $L^2(|a\cdot n| \, \chi \, d\sigma)$. Next, we may write the 
Green formula 
\bean
 \int_\OO \bigl[ G_\eps \, \varphi + g_\eps \Div(a\varphi) \bigr] \, dy  =   \int_\Sigma  \gamma \, g_\eps \, \varphi \, a \cdot n \, d\sigma_{\! y}  , 
\eean
for any test function $\varphi \in C^1_c( \bar\OO)$, and we may pass to the limit as $\eps \to 0$. We deduce that the 
Green formula 
\beqn\label{eq:Transport-def3traceBis}
\int_\OO \bigl( G \, \varphi + g \, \Div(a\varphi)  \bigr) \, dy =
\int_\Sigma \gamma g \, \varphi \, a(y) \cdot n(y) \, d\sigma_{\! y},
\eeqn
holds for any $\varphi \in C^1_c( \bar\OO)$. That clearly uniquely defines the trace function $\gamma g$ on $\Sigma\backslash\Sigma_0$. 

\smallskip
Now, on the one hand, from the DiPerna-Lions renormalizing theory \cite[proof of Corollary~II.1]{MR1022305}, we know that 
 $\beta(g) \in L^\infty(  \OO)$ satisfies the  transport equation 
\beqn\label{eq:Transport-renormalizationTransportOO_b=0}
 a(y) \cdot \nabla_y \beta(g) = \beta'(g) G \quad\hbox{in}\quad \DD'(\OO), 
\eeqn
 for any renormalizing function $\beta \in \Lip(\R)$ and any test function $\varphi \in C^1_c(\OO)$.  Using the already established trace result, we know that
 there exists $\gamma \beta(g) \in L^\infty(\Sigma\backslash\Sigma_0)$ such that 
\beqn\label{eq:Transport-theo:traceLinfty:traceRenorm1}
 \int_\OO \bigl[ \beta'(g) G \, \varphi + \beta(g) \Div(a\varphi) \bigr] \, dy  =   \int_\Sigma  \gamma \beta(g) \, \varphi \, a \cdot n \, d\sigma_{\! y}  , 
\eeqn
for any test function $\varphi \in C^1_c( \bar\OO)$. 
 On the other hand, from the classical   Green formula for Lipschitz functions and because $ \beta(g_{\eps})_{|\Sigma} = \beta(g_{\eps|\Sigma}) $, we have 
 \bean
\int_\OO \bigl[ \beta'(g_\eps) G_\eps \, \varphi + \beta(g_\eps)  \Div(a\varphi) \bigr] \, dy  =  \int_\Sigma  \beta(g_{\eps|\Sigma})   \, \varphi \, a \cdot n \, d\sigma_{\! y} , 
\eean
for any renormalizing function $\beta \in \Lip(\R)$ and any  test function $\varphi \in C^1_c( \bar\OO)$. Using that 
$$
 \beta'(g_\eps) G_\eps \to \beta'(g) G, \quad \beta(g_\eps) \to \beta(g), \quad \beta(g_{\eps|\Sigma}) \to \beta(\gamma g)
 $$
 respectively in $\Lloc^1(\bar\OO)$ and in $\Lloc^1(\Sigma)$, and that the two last sequences are bounded in $L^\infty$,  we may pass to the  limit $\eps \to 0$ in the last 
Green formula,  and we thus get 
\bean
\int_\OO \bigl[ \beta'(g) G \, \varphi + \beta(g) \Div(a\varphi) \bigr] \, dy  =  \int_\Sigma  \beta(\gamma \, g) \, \varphi \, a \cdot n \, d\sigma_{\! y} . 
\eean
Together with \eqref{eq:Transport-theo:traceLinfty:traceRenorm1} and by uniqueness of the trace function, we conclude to $\gamma \beta(g) = \beta(\gamma g)$.
\end{proof}


Let us state several variants of the preceding  trace result. For the transport evolution equation \eqref{eq:transport:EvolAVECSource} a first possible trace result writes as follows. 

\begin{theo}\label{theo:traceEvolLinfty} Assume that $g \in L^\infty((0,T) \times \OO)$, $a \in L^1(0,T;\Wloc^{1,1}(\bar\OO)) $ and $G \in \Lloc ^1([0,T] \times \bar \OO)$
satisfy the evolution  transport equation \eqref{eq:transport:EvolAVECSource} in the distributional sense. Then,
$$
g \in C([0,T];\Lloc^1(\bar\OO))
$$
and  there exists a unique function
$$
\gamma g \in L^\infty( (0,T) \times \Sigma \backslash \Sigma_0; \, dt\otimes d\sigma_{\! y}  ), 
\quad \| \gamma g \|_{L^\infty} \le \| g \|_{L^\infty}, 
$$
which satisfies the   renormalized Green formula  
\bear\label{eq:Transport-def4traceBIS-evol}
&&\int_{t_0}^{t_1}\!\! \int_\OO \bigl( \beta'(g) \, G \, \varphi  + \beta(g) \, [\partial_t \varphi +  \Div(a\varphi) ] \bigr) \, dy dt
\\ \nonumber
&&\quad= \Bigl[ \int_\OO \beta(g(t,\cdot)) \varphi dy \Bigr]_{t_0}^{t_1} + \int_{t_0}^{t_1}\!\! 
\int_\Sigma \beta(\gamma g) \, \varphi \, a \cdot n \, d\sigma_{\! y}dt,
\eear
for any $ \varphi \in C^1_c([0,T] \times \bar\OO)$, any $t_0,t_1 \in [0,T]$ and any $\beta \in C^1_{\rm pw}(\R) $. In particular  renormalization and trace operations commute: 
\eqref{eq:Transport-traceCrenormalization} holds. 
\end{theo}

\smallskip
We skip the proof of Theorem~\ref{theo:traceEvolLinfty} which is very similar to the proof Theorem~\ref{theo:traceLinfty1} using the slight modifications that one can find in \cite[Thm.~2]{MR1765137} or
\cite[Thm.~3.1]{MR2150445}. Under the slightly more regularity assumption $a \in \Wloc^{1,1}([0,T] \times \bar\OO)$, Theorem~\ref{theo:traceEvolLinfty} is a direct corollary  of Theorem~\ref{theo:traceLinfty1}
applied to the for field $(1,a(t,y))$ on the open set $(0,T) \times \OO$.


\medskip
For some additional function $b : \OO \to \R$, another possible variant is the following trace result for the stationary transport equation
\beqn\label{eq:Transport-withb}
 a \cdot \nabla_y  g + b g  =  G  \quad \hbox{in} \quad \OO
\eeqn
 in the renormalized framework as introduced by DiPerna and  Lions in \cite{MR1022305}. 
Assuming $a \in \Wloc^{1,1}(\bar\OO)$, $b,G \in \Lloc^1(\bar\OO)$, we say that  $g \in \Lloc^1(\bar\OO)$ is a renormalized solution
to the transport equation \eqref{eq:Transport-withb} if 
\beqn\label{eq:Transport-renormalizationTransportOO}
 a \cdot \nabla_y \beta(g) + b \beta'(g) g = \beta'(g) G, 
\eeqn
 in the distributional sense for any renormalizing function $\beta \in C^1_*(\R)$ the set of $C^1(\R)$ functions such that $\beta$  admits some finite limits in $\pm\infty$ and $s \mapsto \langle s \rangle \beta'(s)$ is bounded on $\R$,  in particular $C^1_*(\R) \subset C^1_b\R)$. 
We also denote by $\beta \in C^1_{\rm pw,*}(\R)$ the $C^1$ piecewise variant of $C^1_{*}(\R)$.
We will repeatedly use the family of functions $\beta_\delta \in C^1_*(\R)$ defined by $\beta_\delta (s) := s/(1+\delta s^2)^{1/2}$ for any $\delta \in (0,1]$.  We observe that $\beta'_\delta(s) = (1+\delta s^2)^{-3/2}$, so that $s \beta'(s) \to 0$ as $s\to\pm\infty$.

\medskip
 Let us start formulating some basic facts on renormalized solutions to equation \eqref{eq:Transport-withb}. 
 
 \begin{lem}\label{lem:transport-linearity&others} Assume $a \in \Wloc^{1,1}(\bar\OO) $, $b,G \in \Lloc^1(\bar\OO)$.
 
(1) If $g \in \Lloc^1(\bar\OO)$ and $\alpha(g)$ satisfies equation \eqref{eq:Transport-renormalizationTransportOO} 
for one renormalizing function $\alpha : \R \to (-1,1)$ which is bijective and belongs to $C^1_{\rm pw,*}(\R)$
then $\beta(g)$ satisfies equation \eqref{eq:Transport-renormalizationTransportOO} 
 for any renormalizing function $\beta \in C^1_{\rm pw,*}(\R)$.

(2) If $g_1,g_2 \in \Lloc^1(\bar\OO)$ are two renormalized solutions to the transport equations
$$
a \cdot \nabla_y g_i + b g_i = G_i \in \Lloc^1(\bar\OO), 
 $$
 then $g := g_1 + g_2$ is a renormalized solution to the transport equation \eqref{eq:Transport-withb} with $G := G_1 + G_2$. 
 
 (3) If $g$  is a renormalized solution to  the transport equation \eqref{eq:Transport-withb} and $\Phi,c \in L^\infty(\OO)$ satisfy
$$
a \cdot  \nabla_y \Phi = c
 $$
in the distributional sense, then $h := g e^{-\Phi}$ satisfies 
\beqn\label{eq:transport:h=ge-Phi}
 a \cdot \nabla_y h + (b+c)h   = Ge^{-\Phi} 
\eeqn
in the renormalized sense.  
%
%
%
 \end{lem}

\begin{proof}[Proof of Lemma~\ref{lem:transport-linearity&others}.]
We briefly sketch the proof and for more details we refer to  \cite{MR1022305}, in particular to \cite[Lem.~II.2]{MR1022305}. 
It is worth mentioning that only the case $b \in L^\infty(\OO)$ is considered in  \cite{MR1022305}, but it readily extends to our framework. 
Assertion (1) is just a consequence of the chain rule $\beta'(s) = (\beta \circ \alpha^{-1})'(\alpha(s)) \alpha'(s)$ for smooth enough solutions and thus for any solution thanks to 
Lemma~\ref{lem:transport-lemma51} (see the proof of \cite[Cor.~II.1]{MR1022305}) and to a standard approximation procedure in order to deal with piecewise $C^1$ functions. 
In order to establish (2),  we consider two renormalized solutions $g_i$, a renormalized function $\beta \in C^1_*(\R)$ and we write 
 $$
a \cdot \nabla_y \beta(\beta_\delta(g_1) + \beta_\delta(g_2))  =
\beta'(\beta_\delta(g_1) + \beta_\delta(g_2)) [(G_1 - b g_1) \beta'_\delta (g_1) + (G_2 - b g_2) \beta'_\delta (g_2)], 
$$
where we have added the two renormalized formulations \eqref{eq:Transport-renormalizationTransportOO} associated to $\beta_\delta(g_i)$ and renormalized once more the resulting solution using (1). 
Letting $\delta \to 0$, we immediately obtain 
 $$
a \cdot \nabla_y \beta(g_1+g_2)  =
\beta'(g_1+g_2) [G_1 + G_2 - b(g_1 + g_2)] 
$$
in the distributional sense. For proving (3), we introduce the mollified sequence $(g_\eps)$ and $(\Phi_\eps)$ defined as in the statement of Lemma~\ref{lem:transport-lemma51} so that 
$$
a \cdot \nabla g_\eps + b g_\eps = G_\eps, \quad a \cdot \nabla \Phi_\eps = c_\eps
$$
with $G_\eps  \to G$ and $c_\eps \to c$ in $\Lloc^1(\bar\OO)$ as $\eps\to0$. The smooth function $h_\eps := g_\eps e^{-\Phi_\eps}$ satisfies 
$$
a  \cdot \nabla_y h_\eps + (b+c_\eps)h_\eps   = G_\eps e^{-\Phi_\eps} 
$$
and then 
$$
a  \cdot \nabla_y \beta(h_\eps) + \beta'(h_\eps) (b+c_\eps)h_\eps   = \beta'(h_\eps)  G_\eps e^{-\Phi_\eps} 
$$
for any $\beta \in C^1_*(\R)$. Passing to the limit $\eps\to0$, we obtain the renormalized formulation of \eqref{eq:transport:h=ge-Phi}.
\end{proof}

We now generalize the trace result to the framework of renormalized solutions. 

\begin{theo}\label{theo:traceL0}  
Assume that $a \in   \Wloc^{1,1}(\bar\OO)$,  $b,G \in \Lloc^1( \bar \OO)$ and that  $g \in \Lloc^1(\bar\OO)$ is a renormalized solution
to the transport equation \eqref{eq:Transport-withb}. Then there exists a unique function
$$
\gamma g \in L(  \Sigma \backslash \Sigma_0; \, d\sigma_{\! y} )
$$
which satisfies the renormalized Green formula
\beqn\label{eq:Transport-def4traceTER}
\int_\OO \bigl( \beta'(g) \, (G - b g) \, \varphi + \beta(g) \, \Div(a\varphi)  \bigr) \, dy =
\int_\Sigma \beta(\gamma g) \, \varphi \, a \cdot n \, d\sigma_{\! y},
\eeqn
for any $ \varphi \in C^1_c(\bar\OO)$ and any $\beta \in C^1_{\rm pw,*}(\R) $. 
\end{theo}


\begin{proof}[Proof of Theorem~\ref{theo:traceL0}.]
We fix $\beta_1 : \R \to \R$ defined by $\beta_1(s) := s (1+s^2)^{-1/2}$, 
 so that $\beta_1 \in C^1_b(\R)$ and $\beta_1 : \R \to (-1,1)$ is a bijection. Since then $\beta_1(g) \in L^\infty(\OO)$ and $\beta'(g)(G-bg) \in \Lloc^1(\bar\OO)$, we know from Theorem~\ref{theo:traceLinfty1} that $\gamma \beta_1 (g)$ is well defined in $L^\infty(\Sigma \backslash \Sigma_0)$ through the Green formula
 \bean
 \int_\OO \bigl[  \beta_1'(g) (G - b g) \, \varphi + \beta_1(g) \, \Div(a\varphi) \bigr] \, dy   = \int_\Sigma  \gamma \beta_1(g) \, \varphi \, a \cdot n \, d\sigma_{\! y}, 
\eean
for any test function $\varphi \in C^1_c( \bar\OO)$. 
 We set $\gamma g := \beta^{-1}_1 (\gamma \beta_1 (g)) \in L(  \Sigma \backslash \Sigma_0; \, d\sigma_{\! y} )$, with the convention $\beta_1^{-1}(\pm 1) = \pm \infty$. 
 For any $\beta \in C^1_{\rm pw,*}(\R)$, we then deduce 
 $$
 \gamma \beta(g) = \gamma [(\beta \circ \beta_1^{-1}) (\beta_1 (g))] =   \beta \circ \beta_1^{-1} (\gamma \beta_1 (g))  =  \beta ( \gamma g ),
 $$
where we have used the renormalization result stated in Theorem~\ref{theo:traceLinfty1} and the chain rule (1) stated in Lemma~\ref{lem:transport-linearity&others} in the second equality and the very definition of $\gamma g$ in the third equality.
In other words, the renormalized Green formula \eqref{eq:Transport-def4traceTER} holds. 
\end{proof}

\begin{rem}\label{rem:Transport-trace}
(1) We will see in Section~\ref{subsec:TranspEq-tau&characteristics} that under the same conditions as in Theorem~\ref{theo:traceL0}  the information on $\gamma g$  can be slightly  improved, in particular $\gamma g \in L^0(\Sigma \backslash \Sigma_0)$. 

(2)  Theorem~\ref{theo:traceL0} in particular holds when we assume $a \in  \Wloc^{1,p'}(\bar\OO)$, $b \in \Lloc^{p'}(\bar\OO)$, $G \in \Lloc^1( \bar \OO)$ and  $g \in   \Lloc^p(\bar\OO)$ satisfy the  transport equation \eqref{eq:Transport-renormalizationTransportOO} in the distributional sense. Indeed, in that situation one knows from the DiPerna-Lions renormalizing theory \cite[Cor.~II.2]{MR1022305} that $g$ is  also a renormalized solution
to the transport equation \eqref{eq:Transport-renormalizationTransportOO} (in the above sense). 
 
\smallskip
(3) Assuming more interior integrability on the functions $g$, $b$, $G$ and $a$, we may deduce more accurate information on $\gamma g$.
A typical example, is that 
$$
\int_{\Sigma \cap B_R} |\gamma g|^r (|a \cdot n| \wedge 1)^2 d\sigma_{\!y} < \infty,
$$
for some $r \in [1,\infty)$ and  any $R>0$, under the  additional assumption
$$
|g|^r   (|\Div a|  + |a  \cdot \nabla T_1(a \cdot n)| +|b| )\in \Lloc^1(\bar\OO), \quad |g|^{r-1} |G| \in \Lloc^1(\bar\OO). 
$$
The proof follows by choosing $\varphi := T_1(a \cdot n)  \, \chi $, $\chi \in C^1_c(\bar\OO)$, $0 \le \chi \le 1$,  and $\beta_k(s) = (|s| \wedge k)^r$ in the associated Green formula \eqref{eq:Transport-def4trace}, and then to pass to the limit $k \to \infty$. 

\smallskip
(4) Even more integrability on $\gpm g$ is available on one part of the boundary if additional integrability assumption is made on $\gmp g$ on the other part of the boundary. 
A typically example, is that 
$$
\int_{\Sigma_\pm \cap B_R} |\gpm g|^r |a \cdot n|   d\sigma_{\!y} < \infty,
$$
 under the  additional assumption
 $$
 |g|^r   (|\Div a|  + |a|  +|b| )\in \Lloc^1(\bar\OO), \quad |g|^{r-1} |G| \in \Lloc^1(\bar\OO),  \quad  |\gmp g|^r a \cdot n \in \Lloc^1(\bar\Sigma_\mp).
$$
The proof follows by choosing $\varphi    \in C^1_c(\bar\OO)$, $0 \le \varphi \le 1$,  and $\beta_k(s) = (|s| \wedge k)^r$ in the associated Green formula \eqref{eq:Transport-def4trace}, and then to pass to the limit $k \to \infty$. 

\smallskip
(5) The results stated in Lemma~\ref{lem:transport-linearity&others}, 
 in  Theorem~\ref{theo:traceL0} and in points (1), (2), (3) and (4) above may be straightforwardly adapted to the 
evolution  transport equation \eqref{eq:transport:EvolAVECSource}. We refer to \cite{MR1765137,MR1776840,MR2721875,MR2150445} where such results are established in a slightly less general framework. 
Let us emphasize again that when $a \in \Wloc^{1,1}([0,T] \times \bar\OO)$ (what it is the case in the time independent case when $a$ satisfies \eqref{eq:Transport-Hyp-a}) this extension is directly implied by Theorem~\ref{theo:traceL0} applied to the vector field $(1,a(t,y))$  in the open set $(0,T) \times \OO$.  \end{rem}

 \medskip
\subsection{Well-posedness for the transport equation with given inflow at the boundary} 
\label{subsec:transport-WPinflow}
\

\medskip
We deduce from the previous trace theorems and standard tools the well-posedness for the transport equation with several boundary conditions. 
In this section, we deal with the transport equation with given inflow at the boundary. 
We are first concerned with the stationary  transport equation
\bear\label{eq:Transport-PrimalTransport2}
\lambda g + a \cdot \nabla g + bg = G \ \hbox{ in } \OO, \quad
\gamma_-  g  = \frakg  \ \hbox{ on } \Sigma_-, 
\eear
for a real number $\lambda \in \R$ large enough, a given source term $G : \OO \to \R$ and a boundary term $\frakg : \Sigma_- \to \R$. 
As we will see, our analysis also apply to the  associated dual equation 
\bear  \label{eq:Transport-DualTransport2}
\lambda \varphi -  a \cdot \nabla \varphi + (b-\Div a) \varphi = \Phi \ \hbox{ in } \DD'(\OO), \quad
\gamma_+ \varphi  = \psi  \ \hbox{ on } \Sigma_+. 
\eear
We will also consider the related evolution equation
\begin{equation}\label{eq:Transport-evolTransport2}
\left\{
\begin{aligned}
&   \frac{\partial g }{ \partial t}  + a \cdot \nabla g + b g = G\ \hbox{ on }\  (0,T) \times \OO, 
 \\
&
\gamma_- g = \mathfrak{g} \ \hbox{ on }\  (0,T) \times \Sigma_-, 
\quad g(0,\cdot) = g_0 \ \hbox{ on } \ \OO,
\end{aligned}
\right. 
\end{equation}
with  given source term $G : (0,T) \times \OO \to \R$, boundary term $ \mathfrak{g}  :  (0,T) \times \Sigma_- \to \R$ and initial datum $g_0 :  \OO \to \R$. 

\smallskip
A possible simple framework consists in imposing the following conditions  
\beqn\label{eq:TranspEqHypWellPose1}
a \in \Wloc^{1,1}(\bar\OO),  \ \ b  \in \Lloc^1(\bar\OO), 
\eeqn
and 
\beqn\label{eq:TranspEqHypWellPose2}
b_-, \Div a  \in L^\infty(\OO), \ \ \frac{ a }{  \langle y \rangle \langle b_+ \rangle}  \in L^1 (\OO) + L^\infty(\OO).
\eeqn
The first condition on $a$ is useful for the renormalization trick and the definition of the trace, the second condition is needed for the existence results in a $L^p$ framework when $p\not=\infty$
and the last condition is used for proving the uniqueness result. 
In order to be able to apply our results to more general (and realistic) situations, we rather consider the following situation.  We assume that $a$ and $b$ satisfy \eqref{eq:TranspEqHypWellPose1}, 
and defining 
\beqn\label{eq:transport-def-varpi} 
 \varpi = \varpi_p := b  -  \frac1p  \Div a - a \cdot \frac{ \nabla m }{ m}, 
\eeqn
for some smooth enough weight  function $m : \bar\OO \to (0,\infty)$ and some exponent $p \in [1,\infty]$, we assume
\beqn\label{eq:TranspEqHypWellPoseWeight}
 \varpi_-  \in L^\infty(\OO), \ \  b, \Div a  \in L^\infty_{\langle \varpi_+\rangle^{-1}} (\OO), \ \ \frac{ a }{  \langle y \rangle}  \in L^\infty_{\langle \varpi_+\rangle^{-1}}(\OO) + L^1 (\OO)
\eeqn
 In the case $p=1$ and $p = \infty$, we will additionally assume $(\varpi_{q})_- \in L^\infty$ for any $q \in (1,\infty)$.
It is worth emphasizing that condition \eqref{eq:TranspEqHypWellPoseWeight} automatically holds when $m \equiv 1$ and $a,b$ satisfy \eqref{eq:TranspEqHypWellPose2}.
We also define the critical real number
\beqn\label{eq:transport-deflambda*pm} 
\lambda_p^* = \lambda_{p}^*(a,b,m) := \| \varpi_-\|_{L^\infty},
\eeqn
and we may observe that 
\beqn\label{eq:transport-lambda*pm&duality} 
\lambda^*_{p'} (-a, b - \Div a, m^{-1}) = \lambda^*_{p} (a, b,m), 
\eeqn
what links up the primal and the dual problems. 
In order to shorten notations, we introduce the three weight functions
\beqn\label{eq:transport-defmOOmSigma}
m_\OO :=   m \langle \varpi_+ \rangle^{1/p}, \quad
\widetilde m_\OO := m \langle \varpi \rangle^{-1/p'}, \quad
m_\Sigma := m |a \cdot n|^{1/p}.
\eeqn


\medskip
We start with a general discussion about a priori bounds, formal representation formulas and general stability results.

 \medskip
{\bf A priori estimates. }
Consider a solution $g$ to the stationary equation  \eqref{eq:Transport-PrimalTransport2}. For any renormalizing function $\beta : \R \to \R_+$ and any function $\varphi : \bar\OO \to (0,\infty)$, we (at least) formally have 
$$
\int_\OO [(\lambda+ b) g   \beta'(g)\varphi  - \beta(g)  (\Div ( a\varphi))] +
\int_{\Sigma_+} a\cdot n \beta(\gamma_+g) \varphi
=  \int_\OO \beta'(g) G\varphi + \int_{\Sigma_-} |a\cdot n| \beta(\frakg) \varphi.
$$
Choosing $\beta(s) := |s|^p$, $1 \le p < \infty$, and $\varphi := m^p$, we get in particular  
\beqn\label{eq:Transport-identity-existenceTranspLp}
 \int_\OO  |g|^p m^p  \bigl( \lambda + \varpi \bigr) + \frac1p\int_{\Sigma_+} |\gp g|^p m^p a\cdot n
=  \int_\OO G g |g|^{p-2} m^p +  \frac1p\int_{\Sigma_-} | \frakg|^p m^p  |a\cdot n|. 
\eeqn
For $p=1$ and any $\lambda > \lambda_1^*$, we get 
$$
 \int_\OO  |g| m \bigl\{ \lambda - \lambda^*_{1} + \varpi_+ \bigr\}
 +  \int_{\Sigma_+} | \gamma_+ g| m_\Sigma
\le   \int_\OO |G| m   +   \int_{\Sigma_-} | \frakg| m_\Sigma.
$$
For $p \in (1,\infty)$, we split $G=G_1+G_2$ and using the Young inequality, we have 
$$
 \int_\OO G g |g|^{p-2} m^p \le \eps_1 \int g^p m^p + \eps_2 \int g^p m^p_\OO +  \frac{1 }{ p} \frac{1}{ (p'\eps_1)^{p/p'}}  \int G_1^p m^p + \frac{1 }{ p} \frac{1 }{(p'\eps_2)^{p/p'}} \int G_2^p \widetilde m^p_\OO ,
 $$
 for any $\eps_i > 0$. For $\lambda > \lambda_p^*$, we choose  $\eps_1 :=(\lambda-\lambda_p^*)/2$ and  $\eps_2 := 1/2 $, 
 we get
\bean
&&\frac12\int_\OO  |g|^p m^p_\OO + \frac{\lambda-\lambda_p^*}{2}\int_\OO  | g|^p m^p + \frac{1 }{ p} \int_{\Sigma_+} | \gamma_+ g |^p m^p_\Sigma
\\
&&\qquad\le  \frac{2^{p-1} }{ p(p')^{p/p'}} (\lambda-\lambda_p^*)^{1-p} \int_\OO |G_1|^p m^p   
+ \frac{1}{ p(p'/2)^{p/p'}}   \int_\OO |G_2|^p \widetilde m^p_\OO 
+ \frac{1 }{ p} \int_{\Sigma_-} | \frakg|^p m^p_\Sigma.
\eean
We thus deduce 
\beqn\label{eq:transport-LpmEstimate}
\| g \|_{L^p_{m}}  
 \le  \frac{C_p }{ \lambda-\lambda_p^*}   \| G_1 \|_{L^p_m} +   \frac{C_p}{ ( \lambda-\lambda_p^*)^{1/p}}  ( \| \frakg \|_{L^p_{m_\Sigma}} +  \| G_2 \|_{L^p_{\widetilde m_\OO}})
\eeqn
and
\beqn\label{eq:transport-LpmEstimate+}
 \| g \|_{L^p_{m_\OO}} +   \| \gamma_+ g \|_{L^p_{m_\Sigma}}
 \le     \frac{C_p }{ ( \lambda-\lambda_p^*)^{1/p'}}    \| G_1 \|_{L^p_m} +   C_p   ( \| \frakg \|_{L^p_{m_\Sigma}} +  \| G_2 \|_{L^p_{\widetilde m_\OO}}), 
\eeqn
for some numerical constant $C_p \in (0,\infty)$ and any $p \in (1,\infty)$ and also for $p=1$ because of the previous estimate. 
%
Finally, for $\lambda > \lambda_\infty^*$ and $\alpha \in (0,\lambda-\lambda^*_\infty)$, we may proceed exactly as above, but throwing also away the contribution of $\varpi_+$, and we may thus first write 
\beqn
\Bigl(  \lambda  - {\| \varpi_- \|_{L^\infty}} - \frac{\alpha^{p'} }{ p'  } \Bigr)  \| g \|^p_{L^p_m(\OO)}  
\le  \frac{1 }{ p}  \bigl(  \| G/\alpha \|^p_{L^p_m(\OO)}      + \| \frakg |a \cdot n|^{1/p} \|_{L^p_m(\Sigma_-)}^p \bigr) , 
\eeqn
for any $p \in (1,\infty)$ large enough in such a way that the coefficient in front of $ \| g \|^p_{L^p_m(\OO)} $ is positive. 
Taking the power $1/p$ in both sides and passing first to the limit  $p \nearrow \infty$ and next to the limit $\alpha \nearrow \lambda-\lambda^*_\infty$, we end with
\beqn\label{eq:Transport-EstimBdyPbLinfty}
 \| g \|_{L^\infty_m(\OO)}
\le   \max  \Bigl(  \frac{1 }{ \lambda-\lambda^*_\infty} \| G  \|_{L^\infty_m(\OO)}  , \| \frakg   \|_{L^\infty_m(\Sigma_-)} \Bigr).
\eeqn

Consider now a solution to the evolution equation \eqref{eq:Transport-evolTransport2}.  
For any renormalizing function $\beta : \R \to \R_+$ and any test function $\varphi : [0,\infty)\times\bar\OO \to (0,\infty)$, we (at least) formally have 
\bear\label{eq:transport-evolRenormalization} 
&&\int_0^t\int_\OO\big(\beta'(g)G\varphi-\beta'(g)g\varphi+\beta(g)[\partial_t\varphi+\Div(a\varphi)]\big)dyds
\\
&&\nonumber
\qquad=\bigg[\int_\OO\beta(g(s,y))\varphi(s,y)dy\bigg]_0^t+\int_0^t\int_\Sigma \beta(\gamma g)\varphi\, a\cdot n\,d\sigma_y ds.
\eear
Choosing $\beta(s) := |s|^p$, $1 \le p < \infty$, and $\varphi(t,y) := m^p(y)$, we get in particular
\begin{align*}
\int_\OO |g(t)|^p m^p & + \int_0^t \int_{\Sigma_+} a\cdot n \, |\gamma_+g|^p m^p + p \int_0^t\int_\OO |g|^p m^p \varpi \\
& = \int_\OO |g_0|^p m^p + p \int_0^t\int_\OO g|g|^{p-2} G m^p + \int_0^t\int_{\Sigma_-} |\frakg|^p m^p |a\cdot n|.
\end{align*}
Using the Young inequality 
$$
p \int_\OO g|g|^{p-2} G m^p 
\le \frac{p }{ p'}  \int_\OO |g|^p m_\OO^p +  \int_\OO |G|^p \widetilde m_\OO^p
$$
and the Gronwall lemma, we then deduce 
\bear\label{eq:transportLpmAprioriInflow}
&&\| g(t) \|^p_{L^p_m} + \int_0^t e^{p\kappa (t-s)} ( \| g_s \|^p_{L^p_{m_\OO}} +  \| \gamma_+ g_s \|^p_{L^p_{m_\Sigma}} ) \, ds 
\\
&&\nonumber\qquad\le
e^{p\kappa t} \|  g_0 \|^p_{L^p_m} + \int_0^t e^{p\kappa (t-s)} ( \| G_s \|^p_{L^p_{\widetilde m_\OO}} +  \|  \frakg_s \|^p_{L^p_{m_\Sigma}} ) \, ds,
\quad \forall \, t \ge 0, 
\eear
 with $\kappa :=  \| \langle \varpi_- \rangle \|_{L^\infty}$. Passing to the limit $p\to\infty$, we also have 
\bear\label{eq:transportLinftymAprioriInflow}
 \max \bigl( \| g(t) \|_{L^\infty_m},    \| \gamma_+ g(t) \|_{L^\infty_{m}} \bigr) 
\le  e^{\kappa t}  \max \bigl(
 \|  g_0 \| _{L^\infty_m} ,  \sup_{[0,t]}    (\| G_s \|_{L^\infty_m}  +  \|  \frakg_s \|_{L^\infty_m})  \bigr) , 
\eear
for any $t \ge 0$.

\medskip

{\bf Representation formulas.} In a  smooth functions framework or still formally, one classically knows that 
 the solution $g$ to the evolution transport equation \eqref{eq:Transport-evolTransport2} is given by   
\bear\label{eq:transport-characteristic-rep1}
g(t,y) &=&  g_0(Y_{-t}(y)) e^{-\int_0^t b(Y_{s-t}(y)) ds} {\bf 1}_{t < t_\bb} + \frakg(t-t_\bb,y_\bb) e^{-\int_0^{t_\bb} b(Y_{s-t_\bb}(y)) ds}  {\bf 1}_{t > t_\bb}  
\\ \nonumber
&&+  \int_0^{t'_\bb} G(s,Y_{s-t'_\bb}(y)) e^{-\int_0^{t'_\bb-s} b(Y_{-\tau}(y)) d\tau}   du,
\eear
where we recall that  the characteristics $Y$ and the backward exit time $t_{\bf b}$ are defined in \eqref{eq:Transport-characteristics}-\eqref{eq:Transport-def-tb} and we denote $t'_\bb := \min (t,t_\bb)$. \Black
 Similarly, the solution $g$  to the stationary  transport equation
\eqref{eq:Transport-PrimalTransport2} is given by 
\beqn\label{eq:Transport-deftildegBIS}
g(y) =  \frakg(y_\ii(y))) e^{-\int_0^{t_\bb} b(Y_{-s}(y)) ds}  
+  \int_0^{t_\bb} G(Y_{-s}(y)) e^{-\int_0^{s} b(Y_{-\tau}(y)) d\tau}   du.
\eeqn
Alternatively, we may define a semigroup $S_b$ (say on $L^\infty(\OO)$) by 
\begin{equation}\label{eq:Transport-defSb}
(S_b(t) f_0) (y)
:= \left\{
\begin{aligned}
&   f_0( Y_{-t}(y))\exp(-\int_0^t  b( Y_{\tau-t}(y)) d\tau) , \ \hbox{ if } t \in   (0,t_\bb(y)), \\ 
&0 \ \hbox{ otherwise}.
\end{aligned}
\right.
\end{equation}
Given $f_0 : \OO \to \R$, the function $f(t,y) := (S_b (t) f_0)(y)$ is thus a solution to the evolution equation 
$$
\partial_t f + a \cdot \nabla f + b f = 0 \ \hbox{in} \ (0,\infty) \times \OO, 
\quad 
\gamma_- f = 0 \ \hbox{on} \ (0,\infty) \times \Sigma_-.
$$
For $G, \tilde g : \OO \to \R$, we next define
\beqn\label{eq:Transport-deftildeg} 
g := \tilde g + \int_0^\infty e^{-\lambda t} S_b(t) \tilde G \, dt,
\eeqn
with $\tilde G := G - \lambda \tilde g + a \cdot \nabla_x \tilde g - b \tilde g$. By construction, it is a solution  to the stationary  transport equation
\eqref{eq:Transport-PrimalTransport2}.
%

\medskip
{\bf Stability.} 
 We present some stability and continuity results. 
Generalizing slightly \cite[Definitions~2.6 and 3.1]{MR2721875}, we say that a sequence $(g_n)$ of $L(E)$ converges in the renormalized sense to $g$, we note $\ds{g_n \rto g}$, if for any 
$\delta \in \Delta$, $\Delta \subset (0,1]$, $0 \in \overline{\Delta}$,
there exists $\bar \beta_\delta \in L^\infty(E)$ 
such that 
\beqn\label{eq:Transport-defrto}
\beta_\delta (g_n) \wto \bar \beta_\delta \ *\sigma(L^\infty,L^1) \hbox{ as } n \to \infty \  \hbox{ and }
\beta_1(\bar \beta_\delta) \to \beta_1(g) \ \Lloc^1(\bar\OO) \hbox{ as }  \delta\to0.
\eeqn
%
%
%
%
%
We may observe that in particular $g_n \wto g$ weakly $L^1(E)$  or $g_n \to g$ a.e. in $E$ implies $\ds{g_n \rto g}$, {\Cyan and that when defined, the renormalized limit $g$ is unique}.
We refer to \cite{MR2721875} and the references therein for 
more material about the subject. 


 \begin{prop}\label{prop:transport-stability}  
Let us consider four sequences $(g_k)$ of $ \Lloc^1(\bar\OO)$, $(a_k)$  of $ \Wloc^{1,1}(\bar\OO)$, $(b_k)$ and  $(G_k)$  of $\Lloc^1( \bar \OO)$ such that 
$$
a_k \cdot \nabla \beta(g_k) + b_k   \beta'(g_k) g_k = \beta'(g_k)G_k \quad \hbox{in} \quad \DD'(\OO), 
$$
 for any $k \ge 1$ and any $\beta \in C^1_*(\R)$ and four functions   $g \in  \Lloc^1(\bar\OO)$, $a \in  \Wloc^{1,1}(\bar\OO)$,  $b,G \in \Lloc^1( \bar \OO)$ such that $a_k  \to a$ in $ \Wloc^{1,1}(\bar\OO)$ and $b_k \to b$, $G_k \to G$ in $\Lloc^1(\bar \OO)$. Let us denote by $\Sigma_0$ the boundary singular subset associated to $a$. 
  

 (1)  If $g_k \to g$ a.e. in $\OO$ then $g$ satisfies \eqref{eq:Transport-renormalizationTransportOO}
for any  $\beta \in C^1_*(\R)$   and, up to the extraction of a subsequence,
$ \gamma g_k \to  \gamma g$ a.e. on $\Sigma  \backslash \Sigma_0$.

(2)   If  $g_k \wto g$ weakly in $\Lloc^1(\bar\OO)$ then $g$ satisfies \eqref{eq:Transport-renormalizationTransportOO} and, up to the extraction of a subsequence,   
$\ds{\gamma g_k \rto \gamma g}$ on $\Sigma \backslash \Sigma_0$.
 
%
\end{prop}

\begin{rem}\label{rem:transport-stability}
Because of Remark~\ref{rem:Transport-trace}-(5) and the time independence made on $a$ and $b$ in \eqref{eq:TranspEqHypWellPoseWeight}, exactly the same stability result holds for the evolution equation  \eqref{eq:Transport-evolTransport2} as a consequence of  Proposition~\ref{prop:transport-stability}.
\end{rem}

\begin{proof}[Proof of Proposition~\ref{prop:transport-stability}.] 
We split the proof into two steps. 

\smallskip
{\sl Step 1. We establish (1).}
We fix $\beta_1$ as in the proof of   Theorem~\ref{theo:traceL0} and we write  the Green formula
 \bean
 \int_\OO \bigl[  \beta_1'(g_k) G_k \, \varphi - \beta_1'(g_k) g_k b_k \, \varphi + \beta_1(g_k) \, \Div(a_k\varphi) \bigr] \, dy   = \int_\Sigma  \beta_1( \gamma g_k) \, \varphi \, a_k \cdot n \, d\sigma_{\! y}, 
\eean
for any test function $\varphi \in C^1_c( \bar\OO)$. 
There exists  $\bar\beta_1 \in L^\infty(\Sigma \backslash \Sigma_0)$ and a subsequence $(g_{n_k})$ such that 
  $\beta_1( \gamma g_{n_k}) \wto \bar\beta_1$ weakly $\sigma(L^\infty,L^1)$. Passing to the limit in the above equation, we get 
 \bean
 \int_\OO \bigl[  \beta_1'(g) G \, \varphi  - \beta_1'(g) g b \, \varphi  + \beta_1(g) \, \Div(a\varphi) \bigr] \, dy   = \int_\Sigma  \bar\beta_1  \, \varphi \, a \cdot n \, d\sigma_{\! y}.
\eean
From Lemma~\ref{lem:transport-linearity&others} and Theorem~\ref{theo:traceL0}, we deduce that $\bar\beta_1 = \beta_1(\gamma g)$, so that $\beta_1(\gamma g_{n}) \wto \beta_1(\gamma g)$ weakly $\sigma(L^\infty,L^1)$. 
Fixing now $\beta_2 := \beta_1^2 \in C^1_*(\R)$ and repeating the same argument, we get  $\beta_2(\gamma g_{n}) \wto \beta_2(\gamma g)$ weakly $\sigma(L^\infty,L^1)$. 
We then immediately deduce that 
$$
(\beta_1(\gamma g_{n}) - \beta_1(\gamma g))^2 \ \wto \ 0 \quad \hbox{weakly} \ \sigma(L^\infty,L^1), 
$$
so that $\beta_1(\gamma g_{n}) \to \beta_1(\gamma g)$ in $\Lloc^1(\Sigma \backslash \Sigma_0)$. We conclude by using that $\beta_1$ is one-to-one.  

\smallskip
{\sl Step 2. We establish (2).}
We fix $\beta_\delta$ as defined just before the statement of  Lemma~\ref{lem:transport-linearity&others}  and we write  the Green formula
 \bean
 \int_\OO \bigl[  \beta_\delta'(g_k) G_k \, \varphi + \beta_\delta(g_k) \, \Div(a\varphi) \bigr] \, dy   = \int_\Sigma  \beta_\delta( \gamma g_k) \, \varphi \, a \cdot n \, d\sigma_{\! y}, 
\eean
for any test function $\varphi \in C^1_c( \bar\OO)$. 
There exist $\overline{\beta}_\delta, \widetilde{\beta}_\delta, \overline{\beta'}_{\!\delta} \in L^\infty(\OO)$,  $\overline{\gamma\beta}_\delta \in L^\infty(\Sigma \backslash \Sigma_0)$ and a subsequence $(g_{n_k})$ such that 
$\beta_\delta(g_{n_k}) \wto \overline{\beta}_\delta$, $g_{n_k} \beta'_\delta(g_{n_k}) \wto \widetilde{\beta}_\delta$, $\beta'_\delta(g_{n_k}) \wto \overline{\beta'}_{\!\delta}$ and  $\beta_\delta( \gamma g_{n_k}) \wto \overline{\gamma\beta}_\delta$ weakly $\sigma(L^\infty,L^1)$. 
Passing to the limit in the above equation, we get 
 \bean
a  \cdot \nabla_y  \overline{\beta}_\delta + b \widetilde{\beta}_\delta =   \overline{\beta'}_{\!\delta} G \,\ \hbox{ in }\,\ \DD'(\OO), \quad \gamma  \overline{\beta}_\delta =   \overline{\gamma\beta}_\delta 
\,\ \hbox{ on }\,\ \Sigma \backslash \Sigma_0.
\eean
From the fact that $(g_k)$ is locally uniformly integrable, we classically deduce that
$$
 \overline{\beta}_\delta, \widetilde\beta_\delta \to g \,\ \hbox{ in }\,\ \Lloc^1(\OO),
 \quad \overline{\beta'}_{\!\delta} G \to G  \,\ \hbox{ in }\,\ \Lloc^1(\OO),
$$
 as $\delta\to0$. More precisely, the two first convergences come from the elementary inequalities
 $$
\forall \, M > 0, \ \exists \, c_M > 0, \ \  | s - \beta_\delta(s)| \le |s - s \beta'_\delta(s)| \le c_M \delta + |s| {\bf 1}_{|s| > M}, 
 $$
 for any $s \in \R$, $\delta \in (0,1)$,  and the last convergence comes from the convexity inequality $\overline{\beta'}_{\!\delta} \ge \beta'_\delta(g)$ 
 and the elementary inequalities
 $$
\forall \, M > 0, \ \exists \, c_M > 0, \ \  0 \le 1 -  \beta'_\delta(s) \le c_M \delta +  {\bf 1}_{|s| > M}, \ \ \forall \, s \in \R, \ \forall \delta \in (0,1).
 $$
From Step~1, we deduce that $g$ satisfies \eqref{eq:Transport-renormalizationTransportOO} and that, up to the extraction of a subsequence,  
$\overline{\gamma\beta}_\delta =  \gamma\overline{\beta}_\delta
\to \gamma g$ a.e. on $\Sigma \backslash \Sigma_0$. Using the Cantor diagonal process, we obtain 
that there exist two sequences $(\delta_m)$ and $(g_{n_k})$ such that $\delta_m \searrow 0$ and 
$\ds{\gamma g_{n_k}\rto \gamma g}$ in the renormalized sense associated to $(\delta_m)$. 
\end{proof}

\medskip
{\bf Existence.} We establish two existence results of solutions to the transport equation \eqref{eq:Transport-PrimalTransport2}.

\begin{lem}[Existence in $L^\infty_m$]\label{lem:Transport-exist}
We assume that $a$ and $b$ satisfy \eqref{eq:TranspEqHypWellPoseWeight} with $p=\infty$ and some  weight function $m : \bar\OO \to [1,\infty)$. 
For any $\lambda >  \lambda^*_\infty$ and
any given functions $G \in L^\infty_m(\OO)$ and $\frakg \in L_m^\infty(\Sigma_-)$,  there exists $g \in L^\infty_m(\OO)$ solution
 to \eqref{eq:Transport-PrimalTransport2} in the distributional sense. This solution satisfies  \eqref{eq:Transport-PrimalTransport2} in the renormalized sense, the weak maximum principle,  namely 
 \beqn\label{eq:Transport-EstimBdyWeakMaxPrinc}
g\ge 0 \hbox{ in } \OO \  \hbox{ if } \ \frakg \ge 0 \hbox{ on } \Sigma_- \hbox{ and } G\ge 0 \hbox{ in } \OO , 
\eeqn
and the $L^\infty_m$ estimate \eqref{eq:Transport-EstimBdyPbLinfty}.
 \end{lem}

\begin{proof}[Proof of Lemma~\ref{lem:Transport-exist}.]
 The proof follows \cite[Lem.~3]{MR1776840} using \cite[Thm.~2.3]{MR274925}; we only sketch it. 
 Under the stronger regularity assumption $a,b \in C^1_b(\bar\OO)$, $G \in C^1_c(\OO)$, $\frakg := \tilde g_{|\Sigma_-}$ with $\tilde g  \in C^1_c(\bar\OO)$,
 both definitions \eqref{eq:Transport-deftildegBIS} and  \eqref{eq:Transport-deftildeg} provide a classical (and thus also renormalized) solution $g$  to \eqref{eq:Transport-PrimalTransport2}. 
 In such a situation, we may justify the computations made in the above a priori estimates paragraph
 and we conclude that $g$ satisfies the $L^\infty$ estimate \eqref{eq:Transport-EstimBdyPbLinfty}. 
In the general case for $a$, $b$, $\frakg$ and $G$, we introduce some sequences $(a^\eps)$, $(b^\eps)$, $(\frakg^\eps)$ and $(G^\eps)$ of regular and approximating functions so that we may apply the first  step above. 
In that way, we build a sequence $(g^\eps)$ of renormalized solutions to the approximated problem which is uniformly bounded and thus converges (up to the extraction of a subsequence)  in the weakly $*\sigma(L^\infty,L^1)$ sense to a function $g \in L^\infty(\OO)$ satisfying \eqref{eq:transportLinftymAprioriInflow}. 
We then  conclude by passing to the limit $\eps \to 0$ thanks to Proposition~\ref{prop:transport-stability}-(2).
\end{proof}

\smallskip 
We give a first version of an existence result in a $L^p$ framework with strong assumption on the boundary condition.

\begin{lem}[Existence in $L^p_{m_\OO}$]\label{lem:Transport-existp1} 
We assume that $a$ and $b$ satisfy \eqref{eq:TranspEqHypWellPoseWeight}  for some $p \in [1,\infty)$ and some  weight function $m : \bar\OO \to [1,\infty)$. 
For any $\lambda > \lambda^*_{p} $,  $G \in L^p_{\widetilde m_\OO}(\OO)$ and $\frakg \in L^p_{m_\Sigma}(\Sigma_-)$,  there exists   $g \in L^p_{m_\OO}(\OO)$ a renormalized solution
 to the transport equation \eqref{eq:Transport-PrimalTransport2}.
This one satisfies \eqref{eq:transport-LpmEstimate},  \eqref{eq:Transport-EstimBdyWeakMaxPrinc}  and $\gp g \in L^p_{m_\Sigma}(\Sigma_+)$.
 \end{lem}

\begin{proof}[Proof of Lemma~\ref{lem:Transport-existp1}.] 
We argue similarly as during the proof of Lemma~\ref{lem:Transport-exist}. When $\frakg = \tilde g_{|\Sigma_-}$ with $\tilde g$, $a$, $b$, $G$ smooth and
with compact support, the classical solution built above satisfies \eqref{eq:transport-LpmEstimate}, and thus
\beqn\label{eq:transport-LpmEstimateBIS}
\| g \|_{L^p_{m_\OO}(\OO)} \lesssim \| G \|_{ L^p_{\widetilde m_\OO} } + \| \frakg m \|_{L^\infty}    \| a \|^{1/p}_{L^1( {\rm supp}\frakg)}.
\eeqn
For $p > 1$, and under the general conditions  \eqref{eq:TranspEqHypWellPoseWeight} on $a$ and $b$, but still assuming $\frakg = \tilde g_{|\Sigma_-}$ and $G,\tilde g \in C^1_c(\bar\OO)$, 
we may introduce two sequences $(a_\eps)$ and $(b_\eps)$ of smooth functions approximating $a$ and $b$. 
Since the resulting solution $g_\eps$ satisfies \eqref{eq:transport-LpmEstimateBIS}, so that the sequence $(g_\eps)$ is bounded in $L^p_{m_\OO}$, we may argue with the same (compactness) argument as  in the proof of  Lemma~\ref{lem:Transport-exist}. 
We then conclude to the existence of a (renormalized) solution $g \in L^p_{m_\OO}$  to the transport equation  \eqref{eq:Transport-PrimalTransport2} satisfying \eqref{eq:transport-LpmEstimate}. 
Still for $p > 1$, but assuming $G \in L^p_{\widetilde m_\OO}$ and $\frakg \in L^p_{m_\Sigma}(\Sigma_-)$, we may introduce two sequences $(G_\eps)$ and $(\frakg_\eps)$ of smooth functions approximating $G$ and $\frakg$.
Thanks to  \eqref{eq:transport-LpmEstimate}, the associated sequence of solutions $(g_\eps)$  is bounded (and better it is a Cauchy sequence) in $L^p_{m_\OO}(\OO)$ and we conclude again to the existence of a (renormalized) solution $g \in L^p_{m_\OO}(\OO)$  to the transport equation  \eqref{eq:Transport-PrimalTransport2} satisfying \eqref{eq:transport-LpmEstimate}. Finally, in the case $p=1$ and  $\lambda > \lambda^*_{1}$, we may find $q > 1$ small enough such that $\lambda > \lambda^*_{q}$. For $G,\frakg \in L^1 \cap L^q$, the last step imply the existence of a renormalized solution $g \in L^q_{m_\OO} (\OO) $  to the transport equation  \eqref{eq:Transport-PrimalTransport2}. Renormalizing the equation, we deduce that $g$ satisfies \eqref{eq:transport-LpmEstimate} for $p=1$. 
When $G,\frakg \in L^1$, we introduce two sequences $(G_\eps)$ and $(\frakg_\eps)$ of $L^1 \cap L^q$ functions approximating $G$ and $\frakg$, and using \eqref{eq:transport-LpmEstimate} for $p=1$, we deduce that 
the resulting sequence $(g_\eps)$ is a Cauchy sequence  in $L^1_{m_\OO}(\OO)$. We easily conclude again. Finally  $\gp g \in L^p_{m_\Sigma}(\Sigma_+)$ from \eqref{eq:transport-LpmEstimate+}  (see also Remark~\ref{rem:Transport-trace}-(4)).
\end{proof}

\medskip

\medskip
{\bf Uniqueness.} 
We present now a uniqueness result. 
 
\begin{lem}[Uniqueness]\label{lem:Transport-Uniqueness}
We assume that $a$ and $b$ satisfy \eqref{eq:TranspEqHypWellPoseWeight} for some exponent $p \in [1,\infty]$ and some weight function $m : \bar\OO \to [1,\infty)$ as well as for $p=1$ and $m \equiv 1$. 
We additionally assume $\Div a  \in \Lloc^\infty(\bar\OO)$ (what is automatically true under assumption \eqref{eq:TranspEqHypWellPose2}). 
With obvious notations, for any  $\lambda > \max( \lambda^*_{p}(m),\lambda^*_1(1)),$ and  any solution $g \in L^p_{m_\OO} (\OO)$    to the transport equation
\beqn\label{eq:Transport-PrimalTransport2=0}
\lambda g + a \cdot \nabla g + bg = 0 \ \hbox{ in } \ \DD'(\OO), \quad
\gamma_-  g  = 0 \ \hbox{ on } \ \Sigma_-, 
\eeqn
 we have $g \equiv 0$.  
\end{lem}

\begin{proof}[Proof of Lemma~\ref{lem:Transport-Uniqueness}.]
We main follow the proof of \cite[Cor.~II.1]{MR1022305}. 
We fix $\beta \in W^{1,\infty}(\R)$, $\beta (0) = 0$, in such a way that $\beta(g) \in L^p_{m_\OO} \cap L^\infty$ is a solution to 
$$
(\lambda + b)g \beta'(g) + a \cdot \nabla \beta(g) = 0 \ \hbox{ in }\  \DD'(\OO), \quad
\gamma_-  \beta(g)  = 0 \ \hbox{ on } \ \Sigma_-. 
$$
For any  $\psi \in C_c(\OO)$ and  any $\lambda > \lambda^*_{m,p}$, we solve in $L^{p'}_{m^{-1} \langle\varpi_+\rangle^{1/p'}} \cap L^\infty $ the dual problem
\beqn\label{eq:Transport-BoundaryPbPhi}
\lambda \varphi - a \cdot \nabla \varphi + (b -\Div a)\varphi = \psi\ \hbox{ in } \ \DD'(\OO), \quad
\gamma_+ \varphi = 0  \ \hbox{ on }\  \Sigma_+, 
\eeqn
 thanks to Lemma~\ref{lem:Transport-exist} and Lemma~\ref{lem:Transport-existp1}, where we observe that, because of \eqref{eq:transport-lambda*pm&duality}, the necessary condition on $\lambda$ in these results is precisely  the one made here. 
For $\chi \in C^1_c(\R^D)$, ${\bf 1}_{B_1} \le \chi \le {\bf 1}_{B_2}$,  and $R > 0$, we define $\chi_R (x) := \chi(x/R)$.  Using the Green formula \eqref{eq:Transport-def3traceBis},  we have 
\bean
0 = \int_\OO ((\lambda+b)  \varphi - \psi)\beta(g) \chi_R - \int_\OO (\lambda+b) \varphi g\beta'(g) \chi_R +  \int_\OO  \varphi \beta(g) \frac{a }{ R} \cdot (\nabla \chi)_R.
\eean
Because on the one hand $\varphi \beta(g) \in L^1_{\langle \varpi_+\rangle} \cap L^\infty$ and on the other hand  ${a / R} \cdot (\nabla \chi)_R \to 0$ a.e. and is bounded in $  L^\infty_{\langle \varpi_+\rangle^{-1}} + L^1$ 
we deduce that  the last term vanishes when $R \to \infty$. 
Using also that $b \varphi g \in L^1$ thanks to \eqref{eq:TranspEqHypWellPoseWeight}, we may pass to the limit $R \to \infty$ in the above equation and we get
\bean
0 =\int_\OO ((\lambda+b)  \varphi - \psi)\beta(g)  - \int_\OO (\lambda+b) \varphi g\beta'(g)  . 
\eean
We take $\beta := \beta_\delta$ and we observe that $\langle b \rangle |\varpi| |\beta_\delta(g)  - g\beta'_\delta(g)| \le \langle b \rangle |\varphi g| \in L^1(\OO)$. We may then pass to the limit $\delta\to0$ in the last equation, and we get 
\bean
0 = - \int_\OO   \psi g , \quad \forall \psi \in C_c(\OO),
 \eean
 from which we conclude that $g\equiv0$.
 \end{proof}

%
%
%

\medskip
We come to the time dependent transport equation by formulating a general continuity result. 

\begin{prop}\label{prop:Transport-CL0}
Assume that $a \in \Wloc^{1,1}(\bar\OO) $, $b \in \Lloc^1(\bar\OO)$, $G \in  \Lloc^1([0,T] \times \bar\OO)$. Any  renormalized solution $g \in \Lloc^1([0,T] \times \bar \OO)$ to the first equation in \eqref{eq:Transport-evolTransport2}, meaning
$$
 \frac{\partial }{ \partial t}  \beta(g)  + a \cdot \nabla \beta(g) +  \beta'(g)b g = \beta'(g)G\ \hbox{ in }\  \DD'((0,T) \times \OO),
  $$ 
for any renormalizing function $\beta \in C^1_*(\R)$, satisfies $g \in C([0,T];L^0(\OO))$, meaning that  $\beta(g) \in C([0,T];\Lloc^1(\bar\OO))$ for any $\beta \in C_b(\R)$. 

 \end{prop}

\begin{proof}[Proof of Proposition~\ref{prop:Transport-CL0}.]
The proof is a  variant of the proof of \cite[Thm.~II.3]{MR1022305} and we just allude it. 
Because $\beta (g) \in L^\infty((0,T) \times \OO)$ is a solution to the transport equation with source term $ \beta'(g)G-\beta'(g)b g \in \Lloc^1([0,T] \times \bar\OO)$, 
we have $\beta(g) \in C([0,T];\DD'(\OO))$ for any  $\beta \in C^1_*(\R)$. Fixing  $\beta_0 \in C^1_*(\R)$ strictly increasing, we deduce that  $\beta_0(g), \beta_0(g)^2 \in C([0,T];\DD'(\OO))$, 
so that $\beta_0(g)  \in C([0,T];\Lloc^2(\bar\OO))$, and the conclusion. 
\end{proof}

We consider now the time dependent transport equation \eqref{eq:Transport-evolTransport2}.


\begin{prop}[Renormalized solutions]\label{prop:Transport-Rsolutions}
We assume \eqref{eq:TranspEqHypWellPoseWeight} for some $p \in [1,\infty]$ and some weight function $m$.
For any $g_0 \in L^p_m(\OO)$, $G \in L^p_{\widetilde m_\OO} ((0,T) \times \OO)$, $\frakg \in L^p_{m_\Sigma}((0,T) \times \Sigma)$, there exists a unique $g \in C([0,T];\Lloc^1(\bar\OO))$ satisfying the estimate  \eqref{eq:transportLpmAprioriInflow} or \eqref{eq:transportLinftymAprioriInflow} and being a solution to the  transport equation \eqref{eq:Transport-evolTransport2} in the renormalized sense, 
namely 
\begin{equation}\label{eq:Transport-lemRS1}
\left\{
\begin{aligned}
&  \frac{\partial \beta(g) }{ \partial t}  + a \cdot \nabla \beta(g) +  \beta'(g)b g = \beta'(g)G\ \hbox{ on }\  (0,T) \times \OO, 
 \\
&
\gamma_- \beta(g) = \beta(\mathfrak{g}) \ \hbox{ on }\  (0,T) \times \Sigma_-, 
\quad \beta(g)(0,\cdot) = \beta(g_0) \ \hbox{ on } \ \OO,
\end{aligned}
\right.
\end{equation}
for any $\beta \in C^1_{\rm pw,*}$. 
Furthermore,  $g \in C([0,T];L^p_m)$ when $p \in [1,\infty)$ and $g \ge 0$ if $g_0, G, \frakg \ge 0$. 
\end{prop}

\begin{rem}\label{rem:transport-WeakMP}
(1)  The above result extends some previous results due to Bardos in \cite[Chap.~III]{MR274925}, Boyer in \cite[Thm.~4.1]{MR2150445} and Crippa et al in \cite[Thm.~1.1]{MR3525452}
 and \cite[Thm.~1.1]{MR3212249}, where the cases $p=2$ or $p=\infty$ are considered with always the additional assumption $a \in L^\infty$ 
 (in the last paper however the present $W^{1,1}$ bound on $a$ is relaxed into a $BV$
 condition) by adapting the DiPerna-Lions theory developed in \cite[Sec.~II]{MR1022305}.
 
 \smallskip
 (2) We immediately deduce from the above result and Lemma~\ref{lem:transport-linearity&others}-(2) a weak maximum principle: $g_1 \le g_2$ if $g_i$ is renormalized solution 
 to the  transport equation \eqref{eq:Transport-evolTransport2} associated to the data $g_{0i}$, $G_i$, $\frakg_i$ such that $g_{01} \le g_{02}$, $G_{01} \le G_{02}$, $\frakg_{01} \le \frakg_{02}$.
  

\end{rem}

\begin{proof}[Proof of Proposition~\ref{prop:Transport-Rsolutions}.]
We proceed similarly as during the proof of Lemma~\ref{lem:Transport-exist}. 

\smallskip
{\sl Step 1. Characteristics.} 
We assume first $a \in C^1(\R^D)$, $g_0 \in C_c(\OO)$, $b \in C_b(\bar\OO)$, $\frakg \in C_c((0,T) \times \Sigma_0)$, $G \in C_c^1((0,T) \times \OO)$.
 We use the characteristics representation \eqref{eq:transport-characteristic-rep1}. 
We may verify that $\bar g$ both satisfies the transport equation in the renormalized sense  and the boundary conditions in  \eqref{eq:Transport-lemRS1} and we may justify the computations leading to the a priori estimates \eqref{eq:transportLpmAprioriInflow} and \eqref{eq:transportLinftymAprioriInflow}. 

\smallskip
{\sl Step 2. Existence.} 
In the general case, we define some regularized sequence $(a_\eps)$, $(g_{0,\eps})$, $(b_\eps)$ $(\mathfrak{g}_\eps)$, $(G_\eps)$ and thanks to the first step we deduce the existence of an associated function $g_\eps \in C([0,T]; L^p_m)$  satisfies both the equation \eqref{eq:Transport-lemRS1} in the renormalized sense and the a priori estimates
\eqref{eq:transportLpmAprioriInflow} or \eqref{eq:transportLinftymAprioriInflow}. 
When $p > 1$, the sequence $(g_\eps)$  is bounded in $L^\infty(0,T;L^p_m)$ and (up to the extraction of a subsequence)  we may pass to the limit $\eps \to 0$ using Proposition~\ref{prop:transport-stability}-(2) and Remark~\ref{rem:transport-stability}. We have established  the existence of a renormalized solution to the transport equation which satisfies the estimate  \eqref{eq:transportLpmAprioriInflow} or \eqref{eq:transportLinftymAprioriInflow}. When $p=1$, we may for instance proceed in the following way by first assuming $0 \le g_0 \in L^1_m$,  $0 \le G \in L^1_{m} ((0,T) \times \OO)$, $0 \le \frakg \in L^1_{m_\Sigma}((0,T) \times \Sigma)$. We may thus consider some nonnegative approximating sequences $(g_{0,\eps})$ in $L^p_m \cap L^1_m$, $G_\eps \in L^p_{\widetilde m_\OO}  \cap L^1_{m}  $, $\frakg_\eps \in L^p_{m_\Sigma} \cap L^p_{m_\Sigma} $ such that $g_{0,\eps} \nearrow g_0$, $G_\eps \nearrow G $ and $\frakg_\eps \nearrow \frakg$. 
The same construction as above implies the existence of  $0 \le g_\eps \in L^\infty(0,T; L^1_m \cap L^p_m)$ renormalized solution to the transport equation associated to these data and 
such that  $(g_\eps)$ is increasing and uniformly bounded in $L^\infty(0,T;L^1_m)$ thanks to the a priori $L^1_m$ estimate \eqref{eq:transportLpmAprioriInflow}. There thus exists $0 \le g \in L^\infty(0,T;L^1_m)$ such that $g_\eps \nearrow g$, and we get that $g$ is a renormalized solution to the transport equation by using again Proposition~\ref{prop:transport-stability}-(2) and Remark~\ref{rem:transport-stability}. We remove the nonnegative condition on $g_0$, $G$ and $ \frakg $ by introducing the positive and negative parts of each function, using the preceding step in order
to prove the existence of two solutions $0 \le g_\pm \in L^\infty(0,T;L^1_m)$ associated respectively to $(g_{0+},G_+,\frakg_+)$ and $(g_{0-},G_-,\frakg_-)$, 
 and finally defining $g := g_+ - g_-$ which is a renormalized
solution to the transport equation thanks to  Lemma~\ref{lem:transport-linearity&others}  and  Remark~\ref{rem:Transport-trace}-(5). 
 
   \smallskip
{\sl Step 3. Continuity. } From Proposition~\ref{prop:Transport-CL0}, we already know that $g  \in C([0,T]; L^0(\OO))$. 
Together with the a priori estimate \eqref{eq:transportLpmAprioriInflow} or  \eqref{eq:transportLinftymAprioriInflow}, we also have $g  \in C([0,T];\Lloc^1(\OO))$ when $p > 1$. 
When $p \in [1,\infty)$, we may improve the above continuity properties by arguing in the following way. 
 We define $\widetilde g := g m$ and we observe that it is a solution to the transport equation 
$$
\partial_t \widetilde g + a \cdot \nabla \widetilde g + \widetilde b   \widetilde g =  \widetilde G, 
\quad \gamma \widetilde g = \widetilde \frakg, \quad \widetilde g(0) = \widetilde g_0, 
$$
with  $\widetilde b := b  - a \cdot \nabla m / m$, $\widetilde G := G m$, $\widetilde \frakg  := \frakg m$ and $\widetilde g_0 := g_0 m$.  
We write the associated renormalized equation \eqref{eq:transport-evolRenormalization} for the renormalizing function  $\beta_M(s) := (|s| \wedge M)^p$, $M >0$, and 
the test function $\varphi :=  \chi_R$, with $\chi \in C^1_c(\R^d)$, ${\bf 1}_{B_1} \le \chi \le {\bf 1}_{B_2}$ and $\chi_R(y) := \chi(y/R)$. 
Observing in particular that 
$$
 \int_0^t\int_\OO \beta_M(\widetilde g)  a \cdot \nabla \chi_R \to 0 \ \hbox{as}\ R \to \infty, 
$$
because of \eqref{eq:transportLpmAprioriInflow} and \eqref{eq:TranspEqHypWellPoseWeight} by arguing as in the proof of Lemma~\ref{lem:Transport-Uniqueness}, we may pass to the limit in 
the associated renormalized equation 
as $R \to \infty$, and we obtain  
\bean
&&\bigg[\int_\OO\beta_M(\widetilde g)  dy\bigg]_0^t = \int_0^t \!\! \int_\OO \beta_M'(\widetilde g)\widetilde G dyds
- \int_0^t \!\!  \int_\Sigma \beta_M(\gamma \widetilde g) a\cdot n  \,d\sigma_y ds 
\\
&&\nonumber
\qquad \qquad 
+ \int_0^t\!\! \int_\OO  \beta_M(\widetilde g)  \, {\bf 1}_{\widetilde  g >  M} \Div a \,  dyds
- \int_0^t\!\!  \int_\OO p \beta_M(\widetilde  g) \, {\bf 1}_{\widetilde g \le  M} \varpi  dyds.
\eean  
Using again \eqref{eq:transportLpmAprioriInflow} and \eqref{eq:TranspEqHypWellPoseWeight}, we may next pass to the limit as $M\to\infty$ in the above equation, and we get 
$$
\frac{d }{ dt} \int_\OO |\widetilde g|^p   = - p \int_\OO  |\widetilde g|^p   \varpi   + \int_\OO p \widetilde G \widetilde g  |\widetilde g|^{p-2}   + \int_\Sigma |\gamma \widetilde g|^p a\cdot n \in L^1(0,T).
$$
We deduce that $t \mapsto \| g(t) \|_{L^p_m} = \| \widetilde g (t) \|_{L^p}$ is continuous. 
\Black
Consider then $t \in [0,T]$ and $t_k \to t$, so that in particular $\| g_{t_k} \|_{L^p_m} \to \| g_t \|_{L^p_m}$ as $k \to \infty$. 
On the other hand, we have yet established that $\| \beta_0(g_{t_k}) - \beta_0(g_{t}) \|_{L^1(\OO \cap B_R)} \to 0$ as $k\to\infty$ for any $R >0$. There exists thus a subsequence $(g_{t_{k'}})$ such that 
$g_{t_{k'}} \to g_t$ a.e. on $\OO$. Thanks to Brézis-Lieb theorem~\cite{MR699419}, we deduce that $g_{t_{k'}} \to g_t$ in $L^p_m$ and it is the whole sequence which converges by uniqueness of the limit. 
We have thus established $g \in C([0,T];L^p_m)$ when $p \in [1,\infty)$. 

%
%
%

\smallskip
{\sl Step 4. Uniqueness. } Because of  Lemma~\ref{lem:transport-linearity&others}  and  Remark~\ref{rem:Transport-trace}-(5), we just have to prove that $g \equiv 0$ if $g$ is a renormalized solution associated to vanishing data $g_0 = 0$, $G = 0$ and $\frakg = 0$. When $p \in [1,\infty)$, the previous step implies that 
$$
\frac{d }{ dt} \int_\OO |g|^p m^p = \int_\OO |g|^p m^p \varpi  \in L^1(0,T), \quad  \int_\OO |g(0)|^p m^p =0,
$$
and together with the Gronwall lemma, we deduce that $g = 0$. The case $p = \infty$ may be tackled thanks to a duality argument exactly as in the proof of  Lemma~\ref{lem:Transport-Uniqueness}.
\end{proof}

\begin{cor}\label{cor:Transport-semigroupSb} 
The semigroup $S_b$ defined by \eqref{eq:Transport-defSb} extends to a  positive semigroup of contractions in $L^p_m$.
\end{cor}

\begin{proof}[Proof of Corollary~\ref{cor:Transport-semigroupSb}.]
We just apply Proposition~\ref{prop:Transport-Rsolutions} with $G =  \frakg = 0$. When $p\in [1,\infty)$, we define in that way a mapping $L^p_m \to C(\R_+;L^p_m)$, $g_0 \mapsto g$, 
where $g$ denotes the unique
renormalized solution. Defining then  $S(t) g_0 := g(t)$ we have built a strongly continuous semigroup in $L^p_m$. The case $p=\infty$ is identical, except the fact that the semigroup 
is only weak $*\sigma(L^\infty_m,L^1_{m^{-1}})$ continuous. The positivity has been established in Proposition~\ref{prop:Transport-Rsolutions}  and the contraction property comes from 
the estimates \eqref{eq:transportLpmAprioriInflow} and \eqref{eq:transportLinftymAprioriInflow}.
\end{proof}

\Black
\begin{rem} 
 It is worth emphasizing that in Bardos  \cite{MR274925} the semigroup is defined by its representation formula
 for smooth data and by Hille-Yosida theory for $L^2$ data. Here we proceed in another way, by rather following \cite{MR1022305,MR1765137,MR2072842}.
\end{rem}

 \medskip
\subsection{Optimal weighted trace theorem and transport equation with reflection at the boundary} 
\label{subsec:TranspEq-tau}
\
We define the functions $\tau^\pm$ as  the solutions to 
 \beqn\label{eq:Transport-taupm}
\lambda_0  \tau^\pm  \mp  a \cdot \nabla  \tau^\pm = 1 \ \hbox{ in } \DD'(\OO), \quad
\gamma_\pm  \tau^\pm
 = 0 \ \hbox{ on } \Sigma_\pm, 
\eeqn
with $\lambda_0 := 1 + \| \Div a \|_{L^\infty}$. 

\begin{lem}\label{lem:tau+&-} Each of the two equations \eqref{eq:Transport-taupm} has a unique solution $\tau^\pm \in L^\infty(\OO)$ and 
$$
0 < \tau^\pm \le 1 \  \hbox{ a.e. in } \ \OO, \quad  0 < \gamma_\pm \tau^\pm \le 1 \  \hbox{ a.e. on } \ \Sigma_\mp.
$$
\end{lem}

\begin{proof}[Proof of Lemma~\ref{lem:tau+&-}.]
We follow a similar proof as in \cite[Prop.~5.1]{MR2150445} (see also  \cite[Sec.~5]{MR1765137}). 
We only deals with $\tau^-$ since the  case of $\tau^+$ can be handled in the same way. 
The existence of $\tau^-\in L^\infty$, its non negativity and the upperbound are consequences of Lemma~\ref{lem:Transport-exist} while the uniqueness is 
ensured by Lemma~\ref{lem:Transport-Uniqueness}. In order to prove the strict positivity we argue as follows. We first fix $A \in \OO$, $|A| \in (0,\infty)$ and we solve 
$$
\lambda_0\varphi - \Div(a\varphi) = {\bf 1}_A \ \hbox{ in } \DD'(\OO), \quad
\gamma_+ \varphi  
 = 0 \ \hbox{ on } \Sigma_+, 
$$
for which there exists a unique solution $\varphi \in L^1(\OO)$ thanks to  Lemma~\ref{lem:Transport-exist} and Lemma~\ref{lem:Transport-Uniqueness}, which 
furthermore satisfies $\varphi \ge 0$ and $\varphi \not\equiv 0$. We observe that $\tau^- \varphi \in L^1(\OO)$ satisfies 
$$
\Div (a \tau^- \varphi) = \varphi - \tau {\bf 1}_A \ \hbox{ in } \  \OO, \quad  \gamma(\tau^-\varphi) = 0  \ \hbox{ on } \ \Sigma \backslash\Sigma_0.
$$
Thanks to the Green formula, first written for $\beta_\delta(\tau^- \varphi)$ and next passing to the limit $\delta\to0$, we deduce 
$$
0 = \int_\Sigma \gamma(\tau^-\varphi) \, a \cdot n d\sigma_{\! y} = \int_\OO \Div (a \tau^- \varphi) dy = \int_\OO \varphi - \int_A \tau _- dy,  
$$
so that the last integral does not vanish. This being true for any $A \subset \OO$, we get $\tau^- > 0$ a.e. on $\OO$. For $A \subset \Sigma_+$ such that 
$$
0 < \int_A (a \cdot n)_+ d\sigma_{\! y} < \infty,
$$
we solve 
$$
\lambda_0 \varphi - \Div(a\varphi) = 0 \ \hbox{ in } \DD'(\OO), \quad
\gamma_+ \varphi  
 = {\bf 1}_A \ \hbox{ on } \Sigma_+,
$$
thanks to Lemma~\ref{lem:Transport-exist} and Lemma~\ref{lem:Transport-Uniqueness}, and we get a unique solution  $0 \le \varphi \in L^1(\OO)$ such that $\varphi \not \equiv 0$. 
The Green formula again implies 
$$
\int_A \gamma \tau^- (a \cdot n)_+ d\sigma_{\! y}  =   \int_\OO \Div (a \tau^- \varphi) dy = \int_\OO \varphi,
$$
so that the first integral does not vanish. This being true for any $A \subset \Sigma_+$,  we conclude  that $\gamma_+\tau^- > 0$ a.e. on $\Sigma_+$.
\end{proof}



\begin{lem}[Optimal weight]\label{lem:OptimalWeight}
We assume that $a$ satisfies \eqref{eq:TranspEqHypWellPose2} as well as $a \in \Wloc^{1,p'}(\bar\OO)$ for some $1 \le p < \infty$.
For any  $g \in L^p(\OO)$ satisfying \eqref{boundary:trace:TE1} in the distributional sense with $G \in  L^p(\OO)$, 
 the associated trace function $\gamma g$ defined in Theorem~\ref{theo:traceL0} satisfies 
$$
\gamma \, g \in L^p \bigl(\Sigma,| n \cdot a| \, \tau d\sigma    \bigr).
$$
\end{lem}

\begin{proof}[Proof of Lemma~\ref{lem:OptimalWeight}.]
One fixes $\beta_M(z) = (|z| \wedge M)^p$. 
From the DiPerna-Lions renormalizing theory, 
we have
$$ 
a \cdot \nabla ( \beta_M(g) \, \tau^+) = \beta_M'(g) \, G  \, \tau^++ \beta_M(g)   (\tau^+-1)
\ \ \hbox{ in } \DD'(\OO). 
$$
Because  $\beta_M(g_\eps) \, \tau^+ \in L^1 (\OO) \cap L^\infty(\OO)$ and $a/\langle y \rangle \in L^1 + L^\infty$, we may use  the Green formula \eqref{eq:Transport-def3traceBis} with $\phi \equiv 1$,  and we get
\bean
 \int_{\Sigma_-} \beta_M(\gamma g) \, \tau^+ \, | n \cdot a|   d\sigma  
 &=& \int_\OO \bigl\{ (\Div a) \beta_M(g) \, \tau^+ -  \beta_M'(g) \, G  \, \tau^+ + \beta_M(g)  (\tau^+-1) \bigr\} 
 \\
 &\lesssim& \| g \|^{p-1}_{L^p} \, \bigl\{ \| g \|_{L^{p}} + \| G \|_{L^{p}} \bigr\}.
\eean
Passing to the limit $M \to \infty$, we obtain $\gamma_- \, g \in L^p \bigl(\Sigma_-,| n \cdot a| \, \tau d\sigma    \bigr)$.
In a very same way, we prove $\gamma_+ \, g \in L^p \bigl(\Sigma_+,| n \cdot a| \, \tau d\sigma    \bigr)$.
\end{proof}


\medskip
We give now a second version of an existence result in a $L^p$ framework with optimal assumption on the boundary condition in the sense that it is reverse with respect to Lemma~\ref{lem:OptimalWeight}.
That also a posteriori justifies that Lemma~\ref{lem:OptimalWeight} provides the optimal trace result in term of weight function on the boundary.

\begin{lem}[Existence in $L^p$ - optimal assumption]\label{lem:Transport-existp2}
We make the same assumption on $a$, $b$ and $p$ as in Lemma~\ref{lem:Transport-existp1}. 
For any $\lambda > \lambda_{a,b,p}+1/p$ and any given functions $G \in L^p(\OO)$ and $\frakg \in L^p(\Sigma_-,\tau|a\cdot n| d\sigma)$,  there exists $g \in L^p(\OO)$ solution
 to \eqref{eq:Transport-PrimalTransport2}.
%
\end{lem}

\begin{proof}[Proof of Lemma~\ref{lem:Transport-existp2}.]
We only sketch the proof in the case of equation \eqref{eq:Transport-PrimalTransport2}, arguing along the lines of   Lemma~\ref{lem:Transport-exist}.
We start with an a priori estimate. Observing that 
$$
\Div (a \tau^+ |g|^p) = (\Div a) \tau^+ g^p + (\tau^+-1) g^p + p \tau^+ ( G g |g|^{p-2} - b |g|^p - \lambda |g|^p) ,
$$
we have 
\bean
\int_\OO |g|^p \bigl\{ 1 + p \tau^+ ( \lambda + b - \tfrac1p \Div  a - \tfrac1p\bigr\} = \int_{\Sigma_-} |\gamma_-g|^p \tau |a \cdot n| d\sigma + \int_\OO G g |g|^{p-2}\tau^+. 
\eean
Using the condition on $\lambda$, the property $0 \le \tau^+ \le 1$ and the Young inequality, we deduce 
$$
\frac1p\int_\OO |g|^p \le \int_{\Sigma_-} |\gamma_-g|^p \tau |a \cdot n| d\sigma + \frac1p \int_\OO |G|^p. 
$$
We conclude in a similar way as in the proof of Lemma~\ref{lem:Transport-existp1}.
\end{proof}



We  consider now the time dependent transport equation  with positive abstract kernels
\begin{equation}\label{eq:Transport-lemRS3}
\left\{
\begin{aligned}
&  \frac{\partial g }{ \partial t}  + a \cdot \nabla g + b g  = \KKK[g] + G\ \hbox{ on }\  (0,T) \times \OO, 
 \\
&
\gamma_- g = \RRR[g,\gamma_+g] + \frakg   \ \hbox{ on }\  (0,T) \times \Sigma_-, 
\quad g(0,\cdot) = \frakg_0 \ \hbox{ on } \ \OO, 
\end{aligned}
\right. 
\end{equation}
with notations introduced at the beginning of the Section. 
We will work in a weighted Lebesgue space $L^p_m$ with the same conditions on $p$, $m$, $a$ and $b$ as introduced at the beginning of Section~\ref{subsec:transport-WPinflow}. 

%
%

On the other hand, we assume
\bear
&&\label{eq:Transport-RKhyp2}
\KKK : L^p_{m_\OO}(\OO) \to L^p_{\widetilde m_\OO}(\OO) \hbox{ linear and positive}, 
\\
&&\label{eq:Transport-RKhyp3}
 \RRR : L^p_{m_\OO} (\OO) \times L^p_{m_\Sigma} (\Sigma_+) \to L^p_{m_\Sigma}(\Sigma_-) \hbox{ linear and positive in each variable},
\eear
where we recall that the weight functions $m_\OO$, $\widetilde m_\OO$ and $m_\Sigma$ have been defined in \eqref{eq:transport-defmOOmSigma}. 
More precisely, recalling that $\RRR = \RRR_\OO + \RRR_\Sigma$ with $ \RRR_\OO$ and $\RRR_\Sigma$ defined by \eqref{eq:Transport-BdaryCond2}, we assume 
\bear 
\label{eq:Transport-KKKbound}
&&\| \KKK[g] \|_{L^p_{\widetilde m_\OO}}^p \le \alpha_{_{\!\KKK}} \| g \|^p_{L^p_{m_\OO}} +  M_{_{\!\KKK}} \| g \|^p_{L^p_m}, 
\\
\label{eq:Transport-RRRbound}
&&\|\RRR_\OO [g] \|^p_{L^p_{m_\Sigma} }
 \le 
\alpha_{_{{\!\RRR}}} \| g  \|^p_{L^p_{m_\OO} }  +  M_{_{\!\RRR}}  \| g  \|^p_{L^p_{m}},
\quad
\| \RRR_\Sigma [h] \|^p_{L^p_{m_\Sigma} } \le \beta_{_{\!\RRR}}   \| h  \|^p_{L^p_{m_\Sigma}}, 
\eear
with $ \alpha_{_{\!\KKK}} , \,  \alpha_{_{{\!\RRR}}} , \, \beta_{_{{\!\RRR}}} \in [0,1]$, $M_{_{\!\KKK}}, \, M_{_{\!\RRR}}  \ge 0$ and 
\beqn\label{eq:Transport-WPKKKRRRcondAlphaBeta}
\vartheta_\OO :=   (1 - \alpha_{_{{\!\RRR}}}   -  \alpha_{_{\!\KKK}})/2 > 0, 
\quad
\vartheta_\Sigma :=    1 - \beta_{_{\!\RRR}} \ge 0.
\eeqn

Let us emphasize that when $p=1$, the assumption \eqref{eq:Transport-KKKbound} is equivalent to the Lyapunov type condition
$$
\KKK^*[m]  \le \alpha_{_{\!\KKK}} \varpi_+ m + M_{_{\!\KKK}} m.
$$


\begin{prop}\label{prop:Transport-RKsolutions}
 We assume that $a$, $b$, $\KKK$ and $\RRR$ satisfy the conditions  \eqref{eq:TranspEqHypWellPoseWeight},
  \eqref{eq:Transport-RKhyp2}, \eqref{eq:Transport-RKhyp3}, \eqref{eq:Transport-KKKbound}, \eqref{eq:Transport-RRRbound},
  and \eqref{eq:Transport-WPKKKRRRcondAlphaBeta} for some weight function $m:\bar\OO \to [1,\infty)$ and some exponent $p \in [1,\infty)$. 
We consider some data $g_0 \in L^p_m(\OO)$, $G \in L^p_m((0,T) \times \OO)$ and $\frakg  \in L^p_{m_\Sigma}((0,T) \times \Sigma_+)$ with either 

\quad (1) $\beta_{_{{\!\RRR}}} \in [0,1)$; 

or $\beta_{_{{\!\RRR}}} = 1$.
In the latter case,  we assume that $\frakg = 0$ and  we make one of the following additional structural assumption

\quad (2) there exist an exponent $p_0 \in [1,p]$ and a weight function $m_0$ such that $\KKK$ and $\RRR$ satisfy \eqref{eq:Transport-RKhyp2} in $L^{p_0}_{m_0}$
and \eqref{eq:Transport-RKhyp3} in $L^{p_0}_{m_{0\OO}} (\OO) \times L^{p_0}_{m_{0\Sigma}} (\Sigma_+) $, with obvious definitions for the weight functions $m_{0\OO}$ and $m_{0\Sigma}$,
and with  $L^p_m \subset  L^{p_0}_{m_0} $, $L^p_{m_\OO} \subset  L^{p_0}_{m_{0\OO}} $, $L^p_{\hat m_\Sigma} \subset L^{p_0}_{m_{0\Sigma}}$, where $\hat m_\Sigma := m (\tau^+ a \cdot n)^{2/p}$; 


\quad (3) $p=1$ and $\RRR_\Sigma$ is diffusive, namely $\RRR^*_\Sigma [  \tau^+  m_\Sigma ] \ge c_\Sigma m_\Sigma $ a.e. on $\Sigma_+$ with $c_\Sigma >0$. 
 

In the above three  cases, there exists a unique  solution $g \in L^\infty(0,T;L^p_m(\OO))  \cap C([0,T];L^{p_0}_{m_0}(\OO))$ 
satisfying the  transport equation \eqref{eq:Transport-lemRS3} in the renormalized sense as well as  $g \in L^{p_0}(0,T;L^{p_0}_{m_{0\OO}}(\OO))$ and $\gamma g \in L^{p_0}(0,T;L^{p_0}_{m_{0\OO}}(\Sigma))$, 
with $p_0 = p$ and $m_0 = m$ in the first and the third cases.  
\end{prop}

\begin{rem} \label{rem:Transport-RKsolutions}
(1) The above result extends some previous results initiated by   Bardos in \cite[Chap. III]{MR274925} and Beals et al in \cite[Thm.~1\&7]{MR872231},  where however only the kinetic case were considered.  
We refer to Section~\ref{part:application4:Kequation} for a discussion about that important model. 

%
%
%

\smallskip
(2) When $\beta_{_{{\!\RRR}}} = 1$, the existence part of the above result still holds (without any  additional structural assumption).

\smallskip
 (3) Similarly as observed in Remark~\ref{rem:transport-WeakMP}, a weak maximum principle holds: 
 $g_1 \le g_2$ if $g_1$ and $g_2$ are the renormalized solutions to two transport equations \eqref{eq:Transport-lemRS3} such that  (with obvious notations) 
 $b_1 \ge b_2$, $\KKK_1 \le \KKK_2$, $\RRR_1 \le \RRR_2$, $g_{01} \le g_{02}$, $G_{01} \le G_{02}$ and $\frakg_{01} \le \frakg_{02}$. 
That is an immediate consequence of the way we build the solutions $g_i$ thanks to the iterative scheme we present in Step~2 of the proof of Proposition~\ref{prop:Transport-RKsolutions}.

\smallskip
 (4) Another immediate consequence of the iterative way of building the solution, together with the fact that the characteristics representation~\eqref{eq:transport-characteristic-rep1} in  the very first step of the construction,
is the validity of the Duhamel formula
\[S_\LL = S_\BB + S_\BB\AA* S_\LL\]
if we denote by $S_\LL$ the semigroup generated by the transport equation~\eqref{eq:Transport-lemRS3} with $G=\frakg=0$, by $S_\BB$ the semigroup when additionally $\KKK=0$, and $\AA f = \KKK[f]$.
\end{rem}

\begin{proof}[Proof of Proposition~\ref{prop:Transport-RKsolutions}.]
We split the proof into five steps. 

\smallskip
{\sl Step 1. A priori estimates.} For a positive solution, we formally compute  
\begin{align}\label{eq:transportLpmApriori-RRR&KKK-diff}
\frac1p
\frac{d}{dt} \int g^p m^p 
= \frac1p \int_{\Sigma_-} (\RRR[g,\gamma_+g] + \frakg)^p m^p_\Sigma& d\sigma_{\!y} -  \frac1p \int_{\Sigma_+} (\gamma_+ g)^p m^p_\Sigma d\sigma_{\!y} 
\\
&+ \int_\OO \bigl\{ g^{p-1} (\KKK[g] + G) - g^p \varpi \bigr\} m^p.\nonumber
\end{align}

\smallskip

Using  the Young inequality and \eqref{eq:Transport-KKKbound}, we have 
\bean 
\int_\OO g^{p-1} \KKK[g]   m^p
&\le& 
\frac1{p'}  \int_\OO g^p \langle \varpi_+ \rangle m^p + \frac1{p}  \int_\OO \KKK[g]^p \langle \varpi_+ \rangle^{-p/p'} m^p
\\
&\le& 
\bigl( \frac1{p'} + \frac{\alpha_{_{\!\KKK}}}p \bigr) \int_\OO g^p \langle \varpi_+ \rangle m^p + \frac{M_{_{\!\KKK}}}{p}  \int_\OO g^p m^p.
\eean

$\bullet$
When $\frakg = 0$, using also \eqref{eq:Transport-RRRbound} and once more the Young inequality, we then have 
\bean 
\frac1p
\frac{d }{ dt} \int g^p m^p 
&\le& 
\frac1p (\beta_{_{\!\RRR}} - 1)  \int_{\Sigma_+} (\gamma_+ g)^p m^p a \cdot n  
+ \bigl( \frac{M_{_{\!\RRR}} }{ p} + \frac{M_{_{\!\KKK}}}{p} + \| \langle \varpi_- \rangle \|_{L^\infty} \bigr)   \int_\OO g^p m^p
\\
&&+ \bigl( \frac{\alpha_{_{\!\RRR}}}p + \frac1{p'} +   \frac{\alpha_{_{\!\KKK}}}p + \frac\eps{p'} - 1\bigr) \int_\OO g^p \langle \varpi_+ \rangle m^p + \frac{\eps^{-p/p'}}{p}  \int_\OO G^p m^p \langle \varpi_+ \rangle^{-p/p'}
 \eean
 for any $\eps>0$.
Making the choice $\eps := \vartheta_\OO p'/p$, we deduce
\bean 
 \frac{d }{ dt} \| g \|_{L^p_m} + \vartheta_\OO \| g  \|^p_{L^p_{m_\OO} }
+ \vartheta_\Sigma \| \gamma_+ g  \|^p_{L^p_{m_\Sigma}}
 \le 
p\kappa \| g  \|^p_{L^p_{m}}  + C_\OO \| G  \|^p_{L^p_{m}}, 
\eean
with 
$$
\kappa :=  \frac{M_{_{\!\RRR}} }{ p} + \frac{M_{_{\!\KKK}}}{p} + \| \langle \varpi_- \rangle \|_{L^\infty},
\quad
C_\OO :=  (\vartheta_\OO p'/p)^{-p'/p}.
$$
Using the Gronwall lemma, we then obtain
\bear\label{eq:transportLpmApriori-RRR&KKK}
&&
 \| g(t) \|^p_{L^p_m} + \int_0^t e^{p\kappa (t-s)} ( \vartheta_\OO \| g_s \|^p_{L^p_{m_\OO}} + \vartheta_\Sigma \| \gamma_+ g_s \|^p_{L^p_{m_\Sigma}} ) \, ds 
\\
&&\nonumber\qquad\le
e^{p\kappa t} \|  \frakg_0 \|^p_{L^p_m} + C_\OO \int_0^t e^{p\kappa (t-s)} \| G_s \|^p_{L^p_{\widetilde m_\OO}}  \, ds,
\quad \forall \, t \ge 0.
\eear

$\bullet$ 
When  $\frakg \not= 0$ and thus $ \vartheta_\Sigma  > 0$, we control the ingoing boundary term by 
$$
\int_{\Sigma_-} (\RRR[g,\gamma_+g] + \frakg)^p m_\Sigma^p
\le (1+\eps_1) \int_{\Sigma_-} \RRR[g,\gamma_+g] ^p m_\Sigma^p 
+ C_{\eps_1} \int_{\Sigma_-} \frakg^p m_\Sigma^p, \quad \forall \, \eps_1 >0, 
$$
and a very similar computation as above leads to the a priori estimate 
\bear\label{eq:transportLpmApriori-RRR&KKKbis}
&&
 \| g(t) \|^p_{L^p_m} + \int_0^t e^{p\kappa (t-s)} ( \vartheta'_\OO \| g_s \|^p_{L^p_{m_\OO}} + \vartheta'_\Sigma \| \gamma_+ g_s \|^p_{L^p_{m_\Sigma}} ) \, ds 
\\
&&\nonumber\qquad\le
e^{p\kappa t} \|  \frakg_0 \|^p_{L^p_m} +  \int_0^t e^{p\kappa (t-s)} (C_\OO\| G_s \|^p_{L^p_{\widetilde m_\OO}}+C_\Sigma \| \frakg_s \|^p_{L^p_{m_\Sigma}})  \, ds,
\eear
for any $t \ge 0$,  with 
\bean
&&
\vartheta'_\Sigma :=   1 - \beta_{_{{\!\RRR}}} (1+\eps_1),
\quad
\vartheta'_\OO :=   1 - \alpha_{_{{\!\RRR}}}  (1+\eps_1) -  \alpha_{_{\!\KKK}}- \eps\frac{p}{p'}, 
\\
&& 
\kappa :=   \frac{M_{_{\!\RRR}} }{ p} (1+\eps_1)+ \frac{M_{_{\!\KKK}}}{p} + \| \langle \varpi_- \rangle \|_{L^\infty},
\quad
C_\OO :=  \eps^{-p'/p}, \quad C_\Sigma := C_{\eps_1},
\eean
and where we have chosen  $\eps,\eps_1>0$ small enough in such a way  that $\vartheta'_\Sigma > 0$ and $\vartheta'_\OO > 0$.

\smallskip
$\bullet$ When $ \vartheta_\Sigma =0$ and thus $\frakg=0$, we further multiply the equation by $m^p \,\tau^\pm $, where $\tau^\pm$ is defined in \eqref{eq:Transport-taupm}, and integrating, we deduce 
\bean
&&\int_0^T \int_{\Sigma_\mp} \tau^\pm  a \cdot n  (\gamma g)^p m^p  \, d\sigma dt
= \Bigl[ \int_\OO  g^p m^p \, \tau^\pm \Bigr]_T^0 
\\
&&+ \int_0^T\int_\OO p g^{p-1} (\KKK[g] + G) m^p  \, \tau^\pm +  \int_0^T\int_\OO g^p \bigl( \frac{\Div( a m^p) }{ m^p} +1 - \lambda_0 \tau^\pm - p K m^p \tau^\pm \bigr).
\eean
Together with \eqref{eq:transportLpmApriori-RRR&KKK} and $\tau^\pm \in L^\infty(\OO)$, we obtain 
\beqn\label{eq:Transport-RK:Apriori2BIS}
 \int_0^T \int_{\Sigma_+}    (\gamma_+ g)^p   \tau^-  m_\Sigma^p  \le C_T \Bigl( \| g_0 \|^p_{L^p_m} + \| G \|^p_{L^p(0,T;L^p_m)} \Bigr)
\eeqn
and
\beqn\label{eq:Transport-RK:Apriori2}
\int_0^T \int_{\Sigma_-}   [\RRR_\Sigma(\gamma_+ g)]^p   \tau^+ m_\Sigma^p  \le C_T \Bigl( \| g_0 \|^p_{L^p_m} + \| G \|^p_{L^p(0,T;L^p_m)} \Bigr),
\eeqn
for  some constant $C_T \in (0,\infty)$.  In particular, when $p=1$ and $\RRR_\Sigma$ is diffusive, we have 
\bean
c_\Sigma \int_0^T \!\!\int_{\Sigma_+} (\gamma_+ g) m_\Sigma
\le  \int_0^T  \!\! \int_{\Sigma_+}  (\gamma_+ g) \RRR_\Sigma^* (  \tau^+  m_\Sigma)  
=  \int_0^T \!\! \int_{\Sigma_-} \RRR_\Sigma(\gamma_+ g)  \tau^+  m_\Sigma,
\eean
and together with  \eqref{eq:Transport-RK:Apriori2}, we deduce the additional estimate 
\beqn\label{eq:Transport-RK:Apriori3}
c_\Sigma \int_0^T \!\!\int_{\Sigma_+} (\gamma_+ g) m_\Sigma
\le
 C_T \Bigl( \| g_0 \|_{L^1_m} + \| G \|_{L^1(0,T;L^1_m)} \Bigr).
\eeqn


\smallskip
{\sl Step 2. Existence.} As a consequence of these a priori estimates, we may classically build a solution through an iterative scheme. 
For the sake of brevity, we only consider the (more interesting and more difficult) case $b_\Sigma = 1$ (so that $\vartheta_\Sigma=0$ and $\frakg=0$) and $G=0$.
 For a given $0 \le g_0 \in L^p_m(\OO)$, we define a sequence of solution $(h_n)$ starting from $h_0 \equiv 0$ thanks to the recursive definition
\begin{equation*}\label{eq:Transport-rec}
\left\{
\begin{aligned}
&   \frac{\partial h_{n+1} }{ \partial t}  + a  \cdot \nabla h_{n+1} + K h_{n+1} = \KKK[h_n]\ \hbox{ on }\  (0,T) \times \OO, 
 \\
&
\gamma_- h_{n+1} = \RRR[h_n,\gamma_+h_n]  \ \hbox{ on }\  (0,T) \times \Sigma_-, 
\quad h_{n+1}(0,\cdot) = \frakg_0 \ \hbox{ on } \ \OO.
\end{aligned}
\right.
\end{equation*} 
From Proposition~\ref{prop:Transport-Rsolutions}, there exists a unique renormalized solution $h_{n+1} \in C([0,T);L^p_m(\OO))$ to the above equation satisfying the estimate \eqref{eq:transportLpmAprioriInflow} with $g := h_{n+1}$, $G := \KKK[h_n] $ and $\frakg :=\RRR[h_n,\gamma_+h_n] \in L^p_{m_\Sigma}$. 
We observe that $0 \le h_n \le h_{n+1}$ thanks to the  weak maximum principle (see Remark~\ref{rem:transport-WeakMP})  and that $h_n$ satisfies the estimates  \eqref{eq:transportLpmApriori-RRR&KKK} and  \eqref{eq:Transport-RK:Apriori2} where $g$ is replaced by $h_n$. 
Thanks to the  monotonous convergence theorem of Beppo Levi, there exists $g$ satisfying estimates \eqref{eq:transportLpmApriori-RRR&KKK} and $h_n \to g$ in $L^p_{m_\OO}((0,T) \times \OO)$. 
We may pass to the limit in the equation satisfied by $(h_n)$ and we deduce that $g$ is a renormalized solution to 
$$
  \frac{\partial g }{ \partial t}  + a  \cdot \nabla g + b g = \KKK[g]\ \hbox{ on }\  (0,T) \times \OO.
$$ 
From  Theorem~\ref{theo:traceL0} and Remark~\ref{rem:Transport-trace}-(5), the function $g$ admits a trace $\gamma g$ and thanks to  Proposition~\ref{prop:transport-stability}, we have $\gamma h_n \to \gamma g$ a.e. on $\Sigma\backslash\Sigma_0$. 
Because of \eqref{eq:Transport-BdaryCond2} and  the   Beppo Levi   theorem again we deduce that $\RRR [h_n,\gamma_+ h_n] \to \RRR [g,\gamma_+ g]$ a.e. on $\Sigma_-$. 
Together with the fact that $\gamma_- h_n \to \gamma_- g$ a.e. on $\Sigma_-$, we have
established that the boundary condition in \eqref{eq:Transport-lemRS3} holds true. It is worth emphasizing here that $\gamma_+ g \in L^1(\Sigma_+; dr_\Sigma(y,\cdot))$ for a.e. $y \in \Sigma_-$ because of \eqref{eq:Transport-RK:Apriori2}. For  $\frakg_0 \in L^p_m(\OO)$, we  separate the positive and the negative parts $\frakg_0 = \frakg_{0+} - \frakg_{0-}$ and we obtain two renormalized solutions $g^\pm \in L^\infty(0,T;L^p_m)$ associated to $\frakg_{0\pm}$ respectively. By linearity, the function $g := g^+-g^- \in L^\infty(0,T;L^p_m)$ is a renormalized solution to the transport equation and the boundary condition is 
\bean
\gamma_- g 
&=&  \gamma_- g^+ -  \gamma_- g^- = \RR_\OO[g^+] - \RR_\OO[g^-] + \RR_\Sigma[\gamma_+ g^+] -\RR_\Sigma[\gamma_+ g^-] 
\\
&=&   \RR_\OO[g] + \RR_\Sigma[\gamma_+ g],
\eean
where the last term is indeed well defined a.e. from the fact that $\gamma_+ g^\pm \in L^1(\Sigma_+; dr_\Sigma(y,\cdot))$ for a.e. $y \in \Sigma_-$ and thus $\gamma_+ g = \gamma_+ g^+ - \gamma_+ g^-$ belongs to the same spaces. From Proposition~\ref{prop:Transport-CL0}, we already know that $g  \in C([0,T]; L^0(\OO))$ and thus using an interpolation argument 
$g  \in C([0,T]; L^{p_1}_{m_1}(\OO))$ for any $p_1 \in [1,p)$ and any weight function $m_1$ such that $m_1/m \in L^{pp_1/(p-p_1)}$ when $p>1$.

\smallskip
{\sl Step 3.} When $\beta_{_{{\!\RRR}}} < 1$ and $p \in [1,\infty)$, we have \eqref{eq:transportLpmApriori-RRR&KKKbis}, and we may just repeat the proof of Proposition~\ref{prop:Transport-Rsolutions} in order to get  $ g   \in C([0,T]; L^p_m(\OO))$ and the uniqueness of the solution.

\smallskip
{\sl Step 4. We assume $\beta_{_{{\!\RRR}}} = 1$ and the structural assumption (2).} From the estimate \eqref{eq:transportLpmApriori-RRR&KKK-diff} on a solution $g$ and the renormalized formulation of the equation, we deduce that
\bean 
\frac1{p_0}
\frac{d }{ dt} \int |g|^{p_0} m^{p_0} 
&=& \frac1{p_0} \int_{\Sigma_-} |\RRR[g,\gamma_+g]|^{p_0} m^{p_0}_{0\Sigma}   -  \frac1{p_0} \int_{\Sigma_+} |\gamma_+ g|^{p_0} m^{p_0}_{0\Sigma}
\\
&&+ \int_\OO \bigl\{ g |g|^{{p_0}-2} (\KKK[g] + G)m^{p_0} - |g|^{p_0} m^{p_0}_{0\OO} \bigr\} ,
\eean
with a RHS term in $L^1(0,T)$. As above, we thus deduce $g \in C([0,T]; L^{p_0}_{m_0}(\OO))$ and next the uniqueness of the solution.

\smallskip
{\sl Step 5. We assume $\beta_{_{{\!\RRR}}} = 1$ and the structural assumption (3).} In that case, we have $p=1$, $\gamma_+ g \in L^1_{m_\Sigma}((0,T) \times \Sigma_+)$ from \eqref{eq:Transport-RK:Apriori3}
and then $\gamma_- g \in L^1_{m_\Sigma}((0,T) \times \Sigma_-)$ from \eqref{eq:Transport-RRRbound}. We may thus justify the same computation as in Step 4 with $p_0 = 1$, and we deduce 
 $g \in C([0,T]; L^{1}_{m}(\OO))$ and next the uniqueness of the solution.
 \end{proof}

As an immediate consequence of the above analysis, we may associate to the transport equation \eqref{eq:Transport-lemRS3} a semigroup. 

\begin{cor}\label{cor:Transport-RKsolutions}
Under the assumptions of Proposition~\ref{prop:Transport-RKsolutions}, there exists a positive  semigroup $S$ on $L^p_m$ 
such that for any $\frakg_0 \in L^p_m(\OO)$, the function $t \mapsto g(t) := S(t) \frakg_0 \in C(\R_+;L^{p_0}_{m_0}(\OO)) \cap \Lloc^\infty(\R_+;L^p_m(\OO))$ is the unique renormalized solution to the transport equation \eqref{eq:Transport-lemRS3} associated to the initial datum $\frakg_0$ (and with $G= \frakg = 0$).  Furthermore the growth bound satisfies $\omega(S) \le \kappa$. \end{cor}

We end this section by formulating the counterpart of the above result for the associated stationary problem
\begin{equation}\label{eq:Transport-stat-KKK&RRR}
\left\{
\begin{aligned}
& \lambda g +   a \cdot \nabla g + b g  = \KKK[g] + G\ \hbox{ on }\   \OO, 
 \\
&
\gamma_- g = \RRR[g,\gamma_+g] + \frakg   \ \hbox{ on }\    \Sigma_-.
\end{aligned}
\right. 
\end{equation}
 
\begin{prop}\label{prop:Transport-Stat-RKsolutions}
We make exactly the same assumptions as in Proposition~\ref{prop:Transport-RKsolutions} on $a$, $b$, $\KKK$ and $\RRR$  for some weight function $m:\bar\OO \to [1,\infty)$ and some exponent $p \in [1,\infty)$ as well as either
$\beta_{_{{\!\RRR}}} < 1$ holds or $\beta_{_{{\!\RRR}}} = 1$ holds with $ \frakg=0$ and one of the additional structure assumptions (1) or (2).
There exists $\lambda^{**} \in \R$ such that for any $\lambda > \lambda^{**}$,  $G \in L^p_m( \OO)$ and $\frakg  \in L^p_{m_\Sigma}(\Sigma_+)$, there exists a unique  solution $g \in L^p_{m_\OO}(\OO) $
satisfying the  transport equation \eqref{eq:Transport-stat-KKK&RRR} in the renormalized sense and some additional a priori estimates listed during the proof.
\end{prop}

\begin{proof}[Proof of Proposition~\ref{prop:Transport-Stat-RKsolutions}.]
We just explain the main steps. We first establish an a priori estimate. We observe that any solution $g$ to the 
stationary problem \eqref{eq:Transport-stat-KKK&RRR} (at least formally) satisfies 
\beqn\label{eq:Transport-identity-existenceTranspLpKKK&RRR}
 \int_\OO  |g|^p m^p  \bigl( \lambda + \varpi \bigr) + \frac1p\int_{\Sigma_+} |\gp g|^p m^p_\Sigma
=  \int_\OO (\KKK[g] + G) g |g|^{p-2} m^p +  \frac1p\int_{\Sigma_-} |\RRR[g,\gamma_+g] + \frakg |^p m^p_\Sigma.
\eeqn
We then only consider the case $\frakg=0$. Repeating the same computations as in Step 1 of the proof of Proposition~\ref{prop:Transport-RKsolutions} and with the same notations, we get 
\beqn\label{eq:transport-Stat-RKsolutions}
p(\lambda-\kappa)   \| g  \|^p_{L^p_{m}}  +  \vartheta_\OO \| g  \|^p_{L^p_{m_\OO} }
+ \vartheta_\Sigma \| \gamma_+ g  \|^p_{L^p_{m_\Sigma}}
 \le C_\OO \| G  \|^p_{L^p_{m}}.
\eeqn
For $\lambda > \lambda^{**} := \max(\kappa,\lambda^*_p)$ and $G \ge 0$, we next consider the sequence $(h_k)$ in $L^p_{m_\OO}$ defined iteratively as the solution given by  Lemma~\ref{lem:Transport-existp1}  to 
\begin{equation*} 
\left\{
\begin{aligned}
& \lambda h_k +   a \cdot \nabla  h_k + b h_k  = \KKK[h_{k-1}] + G\ \hbox{ on }\   \OO, 
 \\
&
\gamma_- h_k = \RRR[h_{k-1},\gamma_+h_{k-1}]  \ \hbox{ on }\    \Sigma_-,
\end{aligned}
\right. 
\end{equation*}
for $k \ge 1$ and starting from $h_0 \equiv0$. We observe that $(h_k)$ is increasing and satisfies the estimate \eqref{eq:transport-Stat-RKsolutions} where $g$ is replaced by $h_k$. 
We may pass to the limit in the above equation and estimate and we obtain a renormalized solution $g \in L^p_{m_\OO}$ to the transport equation \eqref{eq:Transport-stat-KKK&RRR} and satisfying the estimate \eqref{eq:transport-Stat-RKsolutions}. By linearity, the same holds without sign condition on $G$. 
Finally, considering the three different cases as in Steps 3, 4 and 5 in Proposition~\ref{prop:Transport-RKsolutions}, we similarly show that $g \in L^{p_0}_{m_{0\OO}}$ and $\gamma_+ g \in L^{p_0}_{m_{0\Sigma}}$ for suitable exponent $p_0 \in [1,p]$ and weight function $m_0$. For two such solutions $g_i$ to \eqref{eq:Transport-stat-KKK&RRR}, the function $g := g_2 - g_1$ is also a renormalized solution to \eqref{eq:Transport-stat-KKK&RRR} for which we may justify the identity \eqref{eq:Transport-identity-existenceTranspLpKKK&RRR} with $p=p_0$, $m = m_0$, $G = 0$. We thus deduce that 
\eqref{eq:transport-Stat-RKsolutions} holds with $p=p_0$, $m = m_0$, $G = 0$, and we conclude that $g=0$, what ends the proof of the uniqueness.
\end{proof}

 \subsection{Characteristics} 
\label{subsec:TranspEq-tau&characteristics}
\

In this section we come back to the characteristics method for the evolution and the stationary transport equation. Our aim is in particular to discuss
the representation formula \eqref{eq:transport-characteristic-rep1}.

We consider a vector field   $a : \OO \to \R^D$ which extends to $\R^D$ and, denoting  its extension by the same letter $a$, we at least assume
\beqn\label{eq:Transport-ThODE-DiPLAmbrisio}
a \in \Wloc^{1,1}(\R^D), \quad   \quad  \Div a \in L^\infty(\R^D), 
\quad a/\langle y \rangle  \in  L^1(\R^D) + L^\infty(\R^D).
\eeqn
%
%
After DiPerna and Lions \cite{MR1022305,MR1648524} (see also   \cite[Def.~1]{MR2340443} or \cite[Def.~1]{MR2747466}),  we introduce the following notion of flow. 

\begin{defin}\label{defin:Transport-aeflow}
We name almost everywhere flow associated to \eqref{eq:Transport-characteristics} a measurable function $Y : \R \times \R^D \to \R^D$, $(t,y) \mapsto Y_t(y)$, such that 

(i) for a.e. $y \in \R^D$, the map $t \mapsto Y_t(y)$ is continuous and 
$$
\dot Y_t (y) = a (Y_t(y)) \ \hbox{in} \ \DD'(\R), \quad Y_0(y) = y;
$$

(ii) for a.e. $y \in \R^D$ and for any $s,t \in \R$, there holds $Y_{s+t} (y) = Y_s(Y_t(y))$; 

(iii)  there exists $C \ge 0$ such that  
\beqn\label{eq:transport-aeflow<}
\forall \, t \in \R, \quad e^{-CT} \lambda \le Y(t,\cdot)_\sharp \lambda \le e^{CT} \lambda, 
\eeqn
where $(Y(t,\cdot)_\sharp\lambda) (A) = \lambda(Y(-t, A))$,   $A \subset \R^D$,  is the pushforward
of the Lebesgue measure $\lambda$.
\end{defin}


\smallskip
From \cite[Thm.~III.1]{MR1022305},  \cite[Thm.~31 \& Remark~32]{MR2409676} and \cite[Sec.~3]{MR2340443} (see also Theorem~\ref{theo:Transport-ODE-DiPLbis} below), we know the existence and uniqueness of such an a.e. flow for $a$ satisfying \eqref{eq:Transport-ThODE-DiPLAmbrisio}. In the incompressible case ($\Div a = 0$), this one furthermore  satisfies:
\beqn\label{eq:transport-aeflow=}
\int_{\R^D} \varphi(Y_t(y)) e^{\int_0^t (\Div a) (Y_\tau(y)) d\tau} dy = 
\int_{\R^D} \varphi(y) dy,
\eeqn
  for any $t \in \R$ and any $\varphi \in L^\infty_c(\R^D)$, the space of $L^\infty$ functions with compact support.
 In the compressible case ($\Div a \not= 0$) and  when  $a$ only satisfies  \eqref{eq:Transport-ThODE-DiPLAmbrisio}, it seems not clear that   \cite[Thm.~III.2]{MR1022305}  or \cite[Thm.~31 \& Remark~32]{MR2409676} provides an a.e. flow such that \eqref{eq:transport-aeflow=} holds. In that general case, the  volume identity \eqref{eq:transport-aeflow=} must be replaced by the volume two sides estimate
 (or nearly-incompressible condition):
$$
e^{-t \| \Div a \|_\infty} \int_{\R^D} \varphi(y) dy \le \int_{\R^D} \varphi(Y_t(y))   dy \le
e^{t \| \Div a \|_\infty} \int_{\R^D} \varphi(y) dy, 
$$
for any $0 \le \varphi \in L^\infty_c(\R^D)$ and $t \in \R$, what is nothing but \eqref{eq:transport-aeflow<} with $C :=  \| \Div a \|_\infty$. It is however  quite straightforward to prove from \cite{MR1022305,MR2409676}   that  the  a.e. flows $Y$ satisfies  \eqref{eq:transport-aeflow=} when  $a$  additionally satisfies $\Div a \in C(\R^D)$.  
One possible way to construct the a.e. flow $Y$ is to define $Y = Y_t(y)$ as the unique renormalized solution in $C(\R;L)$ to the transport equation 
$$
\partial_t Y - a \cdot \nabla_y Y = 0 \ \hbox{ on } \ \R \times \R^D, 
\quad
Y_0(y) = y  \ \hbox{ on } \  \R^D, 
$$
or more explicitly 
$$
\frac{\partial }{ \partial t} \beta(Y) - a(Y) \cdot \nabla \beta(Y) =0 \ \hbox{ on }\  \R \times \R^D, 
\quad \beta(Y_0) = \beta(y) \ \hbox{ on } \  \R^D, 
$$
in the distributional sense for any $\beta \in C^1(\R^D,\R)$ such that $\beta$ and $|\nabla \beta(z)| (1+|z|)$ are uniformly bounded on $\R^D$.
In particular,  for any $g_0 \in C^1(\R^D)$ and next for any $g_0 \in L^0(\R^D)$ the function $g^\sharp(t,y) := g_0(Y_{-t}(y))$ is the unique renormalized solution to the transport equation 
$$
\partial_t g^\sharp + a \cdot \nabla g^\sharp = 0 \ \hbox{ on } \ \R \times \R^D, 
\quad
g^\sharp(0,\cdot) = g_0 \ \hbox{ on } \  \R^D.
$$

%
%
%
%
%
%


\medskip
We introduce some notations. 
We denote $y \in \YY$ if (i) and (ii) hold. In particular, $\YY$ is a measurable subset of $\R^D$ and $|\R^D \backslash \YY| = 0$. 
Because of (i), for $y \in \VV := \OO \cap \YY$, we may define the backward exit time
$$
t_{\mathbf{b}} (y) := \sup \bigl\{ \tau >0 \mid  Y_{-t}(y) \in \OO, \; \forall \, s \in [0, \tau] \bigr\} \in (0,+\infty], 
$$
the subset  $\VVb := \{ y \in \VV; \  t_{\mathbf{b}} (y) < + \infty \}$ and the associated {\it entering position} 
$$
y_{\mathbf{b}} (y) := Y_{-t_{\mathbf{b}} (y)}(y)  \  \hbox{ when }\  y \in \VVb. 
$$
%
We observe that the function $t_\bb : \VV \to (0,+\infty]$ is measurable,  $\VVb$ is a measurable subset of $\OO$ and 
\bear\label{eq:transport-tbbYs}
&&t_\bb (Y_s(y)) = t_\bb(y) + s, \quad y_\bb(Y_s(y)) = y_\bb(y), \quad \forall \, y \in  \VVb, \ \forall \, s \in [0,t_\bb(y)). 
\\
\label{eq:transport-meas-tbb=t}
&&\hbox{\rm meas}( \{ y \in \OO; \ t_\bb(y) = t \}) = 0, \quad \forall \, t > 0.
\eear
The properties \eqref{eq:transport-tbbYs} are straightforward while \eqref{eq:transport-meas-tbb=t} is a consequence of the fact that $ \{ y \in \VV; \ t_\bb(y) = t \} \subset Y_{t}(\Sigma)$ and of the 
nearly-incompressible condition \eqref{eq:transport-aeflow<}. 

\smallskip
We now introduce the following first  regularity assumption on $a$ at the boundary
\beqn\label{eq:transport-reg-at-the-boundary}
  \forall \, y_0 \in \Sigma, \  y \mapsto a(y) \cdot n(y_0) \hbox{ is continuous on } \bar\OO.
\eeqn

\smallskip
Let us present two examples. 

\ - It may happen that $\VVb = \emptyset$. For instance, choosing $\OO := \{ y \in \R^2; \ |y| < 1 \}$ the unit disk of the plane and $a(y) : = |y|y^\perp \in C^{0,1}(\R^2;\R^2)$, 
$y^\perp := (y_2,-y_1)$, we have  $\Div a \equiv 0$ and $a(y) \cdot n(y) = y^\perp \cdot y = 0$ for any $y \in \R^2$, so that the flows do not encounter the boundary set $\Sigma = \{ y \in \R^2; \ |y| = 1 \}$. 
In that situation $\VV = \OO$ and $\VVb = \emptyset$. 

\ - In the kinetic case, namely $y = (x,v) \in \OO := \Omega \times \R^d$, $\Omega \subset \R^d$ an open set with smooth boundary with unit normal outward $\nu_x$, so that $n(y) = (\nu_x,0)$, and $a(y) = (v,F(x,v))$, we have $\Div a = \Div_v F$ and $\Div_v F = 0$ when $F = E(x) + v \wedge B(x)$, and we have 
$a(y) \cdot n(y_0) = v \cdot \nu_{x_0}$ which is a smooth function on $\bar \OO \times \Sigma$.

\begin{lem}\label{lem:tbb1}
 Under the condition \eqref{eq:transport-reg-at-the-boundary}, the mapping $y_\bb : \VVb  \to \Sigma_- \cup \Sigma_0$ is measurable.

\end{lem}

\begin{proof}[Proof of Lemma~\ref{lem:tbb1}.]
From the very definitions and composition rules, we have  $y_\bb : \VVb  \to \Sigma \cap \YY$ is measurable.
Take   $y \in \VVb$, denote $y_0 := y_\bb(y)$ and consider a sequence $t_k \searrow 0$ so that $Y_{t_k}(y_0) \to y_0$ and $Y_{t_k} (y_0) \in \OO$ for any $k \ge 1$. From  \eqref{eq:transport-reg-at-the-boundary}, we deduce
$$
0 \ge \limsup_{k\to\infty} \frac{ Y_{t_k} (y_0) - y_0 }{ t_k} \cdot  n(y_0) = \lim_{k\to0}  \frac{1}{ t_k} \int_0^{t_k} a(Y_s(y_0)) \cdot n(y_0) \, ds = a(y_0) \cdot n(y_0), 
$$
which means that $y_0 \in  \Sigma_- \cup \Sigma_0$.  
\end{proof}

For further references, we introduce the following second additional mild regularity assumption on $a$ in the domain
\beqn\label{eq:transport-charac-aLinfty}
a\in \Lloc^\infty(\bar\OO), \quad  (a(y) \cdot y)_+ \lesssim \langle y \rangle^2.
\eeqn
This one may in fact replace the last boundedness condition on $a$ in \eqref{eq:Transport-ThODE-DiPLAmbrisio}. 
We establish a technical result which will be useful in the sequel.

\begin{lem}\label{lem:tbb2} Under assumptions  \eqref{eq:transport-reg-at-the-boundary} and \eqref{eq:transport-charac-aLinfty}, the following hold: 

\smallskip
(1) For any $R_0,T>0$, there exists $R_T,L_T > 0$ such that for any $y \in B_{R_0}$ there hold
\beqn\label{eq:transport-tbb2}
\sup_{[-T,T]} |Y_t(y)| \le R_T \quad \hbox{and}\quad |Y_{t_2}(y) - Y_{t_1}(y)| \le L_T |t_2-t_1|, \quad \forall \, t_1,t_2 \in [-T,T]. 
\eeqn

\smallskip
(2) For a sequence $(y_\eps)$  of $\VV$ such that $y_\eps  \to y_0  \in \Sigma_- \cap \YY$, we have $t_\bb(y_\eps) \to 0$ and $y_\bb(y_\eps) \to y_0$.

\end{lem}

%

\begin{proof}[Proof of Lemma~\ref{lem:tbb2}.]
{\sl Proof of (1).} Take $y \in \VV \cap B_{R_0}$. On the one hand, from  \eqref{eq:transport-charac-aLinfty}, we have 
$$
\langle Y_t(y) \rangle^2 - \langle y  \rangle^2 = 2 \int_0^t a(Y_s(y)) \cdot Y_s(y)) ds  \le C \int_0^t\langle Y_s(y) \rangle^2 ds, 
$$
and we conclude to $Y_t(y) \in B_{R_T}$ thanks to the Gronwall lemma. As a consequence, 
we have 
$$
|Y_{t_2} (y)- Y_{t_1} (y)| \le |t_2 - t_1| \| a \|_{L^\infty(B_{R_T})}, 
\quad \forall \, t_i \in [-T,T].
$$

\smallskip
{\sl Proof of (2).}
Assume by contradiction that $\limsup t_\bb(y_\eps) \ge \tau > 0$ and set $T := \tau +1$.   By assumption, there exists $R_0 > 0$ such that $y_\eps \in B_{R_0}$ and thus by step 1, \eqref{eq:transport-tbb2} holds uniformly in  $\eps \in (0,1]$. 
Thanks to the Ascoli Theorem and the contradiction hypothesis, there exists $\eps_k \to 0$ such that  $t_\bb(y_{\eps_k}) \ge \tau $ and a there exists $Y \in C([-T,T])$ such that $Y_{\bullet} (y_{\eps_k}) \to Y_\bullet$  in $C([-T,T])$. Next, passing to the limit in the conditions 
$$
Y_{-t} (y_{\eps_k}) \in \OO
\quad\hbox{and}\quad
Y_{-t} (y_{\eps_k})- y_{\eps_k} =-  \int_0^t a( Y_{-s} (y_{\eps_k}))   \, ds, 
$$
for any $k \ge 1$ and any $t \in [0,\tau]$, we get 
$$
Y_{-t} \in \bar\OO
\quad\hbox{and}\quad
Y_{-t}  - y_0 = - \int_0^t a(Y_{-s}) ds. 
$$
We deduce 
$$
0 \ge \limsup_{t \to 0}  \frac{Y_{-t} - y_0 }{ t}  \cdot n(y_0) = \lim_{t\to 0}  - \frac1t \int_0^t   a(Y_{-s})   \cdot n(y_0)\, ds  = -   a(y_0)   \cdot n(y_0), 
$$
which is in contradiction with the hypothesis $y_0 \in \Sigma_-$.  
Now, we may estimate 
\bean
|y_\bb(y_\eps) - y_0 | 
&\le& |y_\eps - y_0| + \int_0^{t_\bb(y_\eps)} |a(Y_s(y_\eps))| ds \to 0,
\eean
as $\eps \to 0$, as a consequence of the convergence $t_\bb(y_\eps) \to 0$, the first estimate in \eqref{eq:transport-tbb2} and the first condition in  \eqref{eq:transport-charac-aLinfty}.  
\end{proof}

 We reformulate some {\it  ``space continuity''} of solutions to the transport equation results picked up in \cite[Sec.~7]{MR2150445}. 
Defining 
 $$
 \OO_\alpha := \{y \in \OO;  \, \delta(y) > \alpha\}, \quad \Sigma_\alpha := \{ y \in \OO; \, \delta(y) = \alpha \} = \partial\OO_\alpha, 
 $$
we know from  \cite[Sec.~2]{MR2150445}, that there exists $\alpha_\OO > 0$ such that  for any $\alpha \in (0,\alpha_\OO)$,
 the mapping 
 $$
 \theta_\alpha : \Sigma \to \Sigma_\alpha, \quad  \theta_\alpha (z) := z - \alpha n(z)
 $$
 is an isomorphism  with associated jacobian function $J_\alpha$ and 
\bear\label{eq:transport-spaceC1}
&& \int_{\Sigma_\alpha} h(z') d\sigma_{\! \alpha}(z') =  \int_{\Sigma} h \circ \theta_\alpha(z) J_\alpha(z) d\sigma_{\! z},
\\ \label{eq:transport-spaceC2}
&&\int_{\OO \backslash \OO_\alpha} g(y) dy =  \int_0^\alpha\!\!\int_{\Sigma} g \circ \theta_{\alpha'}(z) J_{\alpha'}(z) d\sigma_{\! z} d\alpha', 
\eear
 for any $h \in L^1(\Sigma_\alpha)$ and $g \in L^1(\OO \backslash \OO_\alpha)$, where $ d\sigma_{\!\alpha}$ denotes
 the Lebesgue measure on $\Sigma_\alpha$ and where the jacobian function $J_\alpha$ satisfies $1/2 \le J_\alpha \le 3/2$ 
 as well as $J_\alpha \to I$ as $\alpha \to  0$. 

 
\begin{lem}\label{lem:tbb3} For any $g \in L^\infty(\OO)$ satisfying $a \cdot \nabla g \in L^1(\OO)$, we have 
\beqn\label{eq:s8-g=galpha}
 g = \gamma_\alpha g   \ \hbox{ a.e. on }\,  \Sigma_\alpha \backslash \{ a \cdot n = 0 \} \,  \hbox{ for a.e.  } \alpha \in [0,\alpha_\OO],
\eeqn
where we denote by $\gamma_\alpha g$ the trace of $g$ on $\Sigma_\alpha$, and 
\beqn\label{eq:s8-alphaTOgalphaIn0}
  \gamma_\alpha g \circ \theta_\alpha \to \gamma g \ \hbox{as} \ \alpha \to 0, \ \hbox{ a.e. on } \Sigma \backslash \Sigma_0.
\eeqn
\end{lem}

\begin{proof}[Proof of Lemma~\ref{lem:tbb3}.]
For $\varphi \in C_c(\bar\OO) \cap W^{1,1}(\OO)$ and  $\beta \in C^1(\R)$, the renormalized Green formula \eqref{eq:Transport-def4traceBIS} writes
\[
\int_{\Sigma} \varphi   \beta(\gamma g)  \, a\cdot n \,  d\sigma 
=
 \int_{\Sigma_\alpha}   \varphi   \beta(\gamma_\alpha g)  \, a\cdot n  \,  d\sigma_{\! \alpha}
+
  \int_{\OO \backslash \OO_\alpha} [\Div (a \varphi)   \beta(g) + \varphi  \beta'(g) a \cdot \nabla  g] \, dy, 
\]
 and thus 
\beqn\label{eq:tbb3-cvgce}
 \int_{\Sigma} \varphi   \beta(\gamma g)  \, a\cdot n \,  d\sigma 
= \lim_{\alpha\to0}
 \int_{\Sigma_\alpha}   \varphi   \beta(\gamma_\alpha g)  \, a\cdot n  \,  d\sigma_{\! \alpha}.
\eeqn
  Denoting $\psi := (1- (\delta(x)-\alpha)/s)_+$, $\alpha+s \in (\alpha,\alpha_\OO)$, observing that $\psi_{|\Sigma_\alpha} = 1$, recalling that $n = - \nabla \delta$ and using $\varphi\psi$ as a test function in the renormalized Green formula \eqref{eq:Transport-def4traceBIS-evol}, we similarly have
 \[
\int_{\Sigma_\alpha} \varphi   \beta(\gamma_\alpha g)  \, a\cdot n \,  d\sigma_\alpha  
=    \frac1s\int_{\OO_\alpha \backslash \OO_{\alpha+s}}   \varphi   \beta(g)  \, a\cdot n  \,  dy
+
  \int_{\OO_\alpha \backslash \OO_{\alpha+s}} [\Div (a \varphi)   \beta(g) + \varphi  \beta'(g) a \cdot \nabla  g] \psi \, dy, 
\]
so that 
 \bean
 \int_{\Sigma_\alpha} \varphi   \beta(\gamma_\alpha g)  \, a\cdot n \,  d\sigma_\alpha 
=  \lim_{s \to 0}  \frac1s\int_{\OO_\alpha \backslash \OO_{\alpha+s}}   \varphi   \beta(g)  \, a\cdot n  \,  dy.
\eean
We immediately deduce 
 \bean 
 \int_{\Sigma_\alpha} \varphi   \beta(\gamma_\alpha g)  \, a\cdot n \,  d\sigma_\alpha  = - \frac{d }{ d\alpha}  \int_{\OO_\alpha}   \varphi   \beta(g)  \, a\cdot n  \,  dy
\eean
and next 
 \bean
\int_0^{\alpha} \!\!\int_{\Sigma_{\alpha'}} \varphi   \beta(\gamma_{\alpha'} g)  \, a\cdot n \,  d\sigma_{\!\alpha'} d\alpha' =   \int_{\OO \backslash\OO_{\alpha}}   \varphi   \beta(g)  \, a\cdot n  \,  dy.
\eean
Together with \eqref{eq:transport-spaceC2} and \eqref{eq:transport-spaceC1}, we have established 
 \bean
 \int_0^\alpha\!\!\int_{\Sigma} [ \varphi   \beta(\gamma_{\alpha'} g)  \, a\cdot n] \circ \theta_{\alpha'}(z) J_{\alpha'}(z) d\sigma_{\! z} d\alpha'
 =
 \int_0^\alpha\!\!\int_{\Sigma} [ \varphi   \beta(g)  \, a\cdot n] \circ \theta_{\alpha'}(z) J_{\alpha'}(z) d\sigma  d\alpha', 
\eean
%
%
so that 
$$
   \beta(\gamma_{\alpha'} g)  \, a\cdot n =   \beta(g)  \, a\cdot n \ \hbox{ a.e. on } \, \Sigma_{\alpha'}  \hbox{ for a.e.  } \, \alpha' \in (0,\alpha),
$$
from what \eqref{eq:s8-g=galpha} follows.
On the other hand, using  \eqref{eq:transport-spaceC1} and \eqref{eq:tbb3-cvgce} together, we have 
$$
 \int_{\Sigma} \varphi   \beta(\gamma g)  \, a\cdot n \,  d\sigma 
= \lim_{\alpha\to0}
 \int_{\Sigma}   [\varphi   \beta( \gamma_\alpha g)  \, a\cdot n]  \circ \theta_\alpha J_\alpha  \,  d\sigma_{\! \alpha}, 
$$
which implies 
\beqn\label{eq:transport-tbb3-cvgceBIS}
 [\beta( \gamma_\alpha g)  \, a\cdot n]  \circ \theta_\alpha J_\alpha \wto  \beta(\gamma g)  \, a\cdot n \quad *\sigma(\Lloc^\infty(\Sigma),L^1_{\rm c}(\Sigma)). 
\eeqn
Repeating the same argument with $g :=  a \cdot n$ and using that $J_\alpha \to I$ uniformly, we get 
$$
 [\beta( a\cdot n)  \, a\cdot n]  \circ \theta_\alpha  \wto  \beta(a\cdot n)  \, a\cdot n \quad *\sigma(\Lloc^\infty(\Sigma),L^1_{\rm c}(\Sigma)), 
$$
for any $\beta \in C^1(\R)$. Choosing $\beta(s) = 1$ and $\beta(s) = s$,  we classically deduce that $[a\cdot n]  \circ \theta_\alpha   \to a\cdot n$ a.e. on $\Sigma$.
We finally conclude to \eqref{eq:s8-alphaTOgalphaIn0} by gathering  that last information with  \eqref{eq:transport-tbb3-cvgceBIS} written for $\beta(s) = s$ and $\beta(s) = s^2$.
\end{proof}

  \begin{rem} During the proof, we have in fact established that 
$$
[0,\alpha_\OO] \to L^1(\Sigma); \quad \alpha \mapsto [\gamma_\alpha g  \, a\cdot n]  \circ \theta_\alpha
$$
is continuous. 
\end{rem}

\begin{lem}\label{lem:tbb} We make the additional assumptions  \eqref{eq:transport-reg-at-the-boundary} and \eqref{eq:transport-charac-aLinfty}. 
 If $Y$ is an  almost everywhere flow associated to \eqref{eq:Transport-characteristics}, then the function $t_\bb \in L(\OO)$ is a renormalized solution to the equation 
\beqn\label{eq:eqtbb}
 a \cdot \nabla t_\bb = 1 \ \hbox{ in } \ \OO,  \quad \gamma_- t_\bb = 0 \ \hbox{ on } \ \Sigma_- \cap \YY. 
\eeqn
\end{lem}

\begin{proof}[Proof of Lemma~\ref{lem:tbb}. Step 1.]
We fix $\beta \in C^1_*(\R)$ and we recall that $\beta(t_\bb(Y_s(y))) = \beta(t_\bb(y) + s)$ for any $s \in \R$ for a.e. $y \in \OO$. 
For any $\varphi \in C^1_c(\OO)$, we may compute
 \bean
&& \int_{\R^d} \beta'(t_\bb(y) + s) \varphi(y) \, dy 
 = \frac{d }{ ds}  \int_{\R^d}\beta(t_\bb(y) + s) \varphi(y) \, dy 
\\
&&\qquad= \frac{d }{ ds}  \int_{\R^d} \beta(t_\bb(Y_s(y))) \varphi(y) \, dy 
 \\
&&\qquad= \frac{d }{ ds}  \int_{\R^d} \beta(t_\bb(y))   \varphi(Y_{-s}(y)) e^{-\int_{-s}^0  (\Div a)(Y_\tau(y)) d\tau }  \, dy
 \\
&&\qquad=   \int_{\R^d} \beta(t_\bb(y))   [ - a \cdot \nabla \varphi - (\Div a) \varphi](Y_{-s}(y))  e^{-\int_{-s}^0  (\Div a)(Y_\tau(y)) d\tau }  \, dy, 
 \eean
  where we have used   the change of variables property \eqref{eq:transport-aeflow=} in the third line.  
Taking $s= 0$, we conclude to 
 \bean
 \int_{\R^d} \beta'(t_\bb) \varphi \, dy 
 = \int_{\R^d} \beta(t_\bb)   [ - a \cdot \nabla \varphi - (\Div a) \varphi]   \, dy, 
 \eean
 which is nothing but the distributional formulation of the equation 
 $$
 a \cdot \nabla \beta(t_\bb) = \beta'(t_\bb).
 $$
 That last family of equations is the renormalized formulation of equation \eqref{eq:eqtbb} in the domain.
 
 \smallskip
{\sl Step 2.} Using lemma~\ref{lem:tbb3} with $g := \beta(t_\bb)$, we have 
$$
 \gamma_\alpha \beta(t_\bb) \circ\theta_\alpha \to \beta(\gamma t_\bb) \ \hbox{ a.e. on } \ \Sigma.
$$
Using Lemma~\ref{lem:tbb2}, we also have 
$$
\gamma_\alpha\beta(t_\bb) \circ\theta_\alpha \to 0  \ \hbox{ a.e. on } \ \Sigma_- \cap \YY.
$$
Both together, we find $\gamma_- t_\bb = 0$ on $\Sigma_- \cap \YY$. 
\end{proof}

\medskip
We establish the main result of this section.

\begin{theo}[characteristics method] \label{theo:Transport-ODE-DiPLbis}
Assume that $Y$ is an  almost everywhere flow associated to \eqref{eq:Transport-characteristics} with $a$ satisfying \eqref{eq:Transport-ThODE-DiPLAmbrisio}. 
For any $g_0 \in L^\infty(\OO)$, $\frakg \in L^\infty((0,T) \times \Sigma_-)$, $T > 0$, $b \in L^\infty(\OO)$,   the function 
\bear
\label{eq:Transport-theo-charac0}
\bar g(t,y) &:=& g_0 (Y_{-t}(y))e^{- \int_0^t b(Y_{\tau-t}(y)) \, d\tau}  {\bf 1}_{t < t_\bb(y)} 
\\ \nonumber
&&\qquad + \frakg(t-t_\bb(y),y_\bb(y)) e^{- \int_0^{t_\bb(y)} b(Y_{\tau-t_\bb(y)}(y)) \, d\tau}   {\bf 1}_{t > t_\bb(y)} 
\eear
is  the unique solution in $C([0,T];\Lloc^1(\OO)) \cap L^\infty( (0,T) \times \OO)$ to the evolution transport equation  
\begin{equation}\label{eq:Transport-theo-charac1}
\left\{
\begin{aligned}
&   \frac{\partial g }{ \partial t}  + a \cdot \nabla g + b g = 0\ \hbox{ on }\  (0,T) \times \OO, 
 \\
&
\gamma_- g = \frakg \ \hbox{ on }\  (0,T) \times \Sigma_-, 
\quad g(0,\cdot) = g_0 \ \hbox{ on } \ \OO.
\end{aligned}
\right.
\end{equation}

\end{theo}


{\sl First proof of Theorem~\ref{theo:Transport-ODE-DiPLbis}.} 
We additionally assume that \eqref{eq:transport-aeflow=},    \eqref{eq:transport-reg-at-the-boundary} and \eqref{eq:transport-charac-aLinfty}  hold, that  $\frakg \in C((0,T) \times \Sigma_-)$,
supp$\, \frakg \subset (0,T) \times (\Sigma_- \cap \YY)$  and $b \in C(\bar\OO)$,
and we mostly repeat the proof of \cite[Prop.~1]{MR1648524}. 
From the above definition, for a.e. $y \in \OO$ and any $t \in (0,\infty)$, we have
\beqn\label{eq:Transport-IdentityOfTheFlow}
\bar g (t+s, Y_s(y)) e^{\int_0^s b(Y_\tau (y)) d\tau} = \bar g (t, y), \quad \forall \, s \ge -t.
\eeqn
Let us then fix $\varphi \in \DD((0,T) \times \OO)$ and let us extend $\bar g$ and $\varphi$ by $0$ outside of $\OO$. We compute  
  \bean
 0 &=& \frac{d }{ ds} \int_0^T\!\!\int_{\R^d} \bar g(t,y) \varphi(t,y) \, dydt 
 \\
&=& \frac{d }{ ds} \int_0^T\!\!\int_{\R^d} \bar g(t+s,Y_s(y)) e^{\int_0^s b(Y_\tau(y)) d\tau }  \varphi(t,y) \, dydt 
\\
 &=&\frac{d }{ ds} \int_0^T\!\!\int_{\R^d} \bar g(t,y)  \varphi(t-s,Y_{-s}(y)) e^{\int_0^s (b- \Div a)(Y_\tau(y)) ds } \, dydt 
\\
 &=& \int_0^T\!\!\int_{\R^d} \bar g(t,y) \frac{d }{ ds} [  \varphi(t-s,Y_{-s}(y)) e^{\int_0^s (b- \Div a)(Y_\tau(y)) ds } ]  \, dydt 
\\
 &=& \int_0^T\!\!\int_{\R^d} \bar g(t,y) [- \partial_t \varphi  - a \cdot \nabla \varphi + (b-\Div a) \varphi](t-s,Y_{-s}(y)) \, dydt , 
 \eean
 where we have used the relation \eqref{eq:Transport-IdentityOfTheFlow} in the second line and the change of variables property \eqref{eq:transport-aeflow=} in the third line. 
 Taking $s = 0$, we get 
\bean
0 = \int_0^T\!\!\int_{\R^d} \bar g(t,y) [- \partial_t \varphi  - a \cdot \nabla \varphi + (b-\Div a) \varphi](t,y) \, dydt, 
 \eean
which exactly means that $\bar g$ is a solution to equation \eqref{eq:Transport-theo-charac1}  in the distributional sense.  
Now, because $t_\bb(y) > 0$ for any $y \in \YY$, we have $\bar g (0,y) = g_0(y)$. 
Take $t > 0$, and for $\alpha \in (0,\alpha_\OO)$, let us denote $A_\alpha := \{ y \in \Sigma_- \cap \YY; \ t_\bb \circ \theta_\alpha (y) < t \}$, so that 
$$
\bar g(t,\theta_\alpha(y)) = \frakg(t-t_{\bb,\alpha}(y),y_{\bb,\alpha}(y)) e^{- \int_0^{t_{\bb,\alpha}(y)} b(Y_{\tau-t_{\bb,\alpha}(y)}(\theta_\alpha(y)) \, d\tau} \ \hbox{ on } \ A_\alpha,  
$$
where we use the shorthands $t_{\bb,\alpha} := t_\bb \circ \theta_\alpha$ and $y_{\bb,\alpha} := y_\bb \circ \theta_\alpha$.
From $ \theta_{\alpha'} (y) \to y$ as $\alpha' \to 0$ when $y \in \Sigma$ and Lemma~\ref{lem:tbb2}, we deduce 
$$
\bar g(t,\theta_{\alpha'}(\cdot)) \to \frakg(t,\cdot) \ \hbox{ as } \alpha' \to 0, 
$$
on $A_\alpha$ for any fixed $\alpha > 0$.  Because $\cup_{\alpha > 0} A_\alpha = \Sigma_- \cap \YY$, the same convergence holds on $\Sigma_- \cap \YY$. 
%
On the other hand, we have $\gamma_\alpha\bar g \circ \theta_\alpha \to \gamma \bar g \ \hbox{as} \ \alpha \to 0$  a.e. on $ \Sigma \backslash \Sigma_0$,  from Lemma~\ref{lem:tbb3}. 
We deduce that $ \gamma \bar g = \frakg$  on  $\Sigma_- \cap \YY$.

\medskip
{\sl Second proof of Theorem~\ref{theo:Transport-ODE-DiPLbis}.} 
We do not make any additional assumption and we mainly repeat the proof presented in  \cite[Sec.~3]{MR2340443}. 
Consider the unique  solution $g \in C([0,T];\Lloc^1(\bar\OO)) \cap L^\infty((0,T) \times \OO)$ to the transport equation 
\eqref{eq:Transport-theo-charac1}. 
Regularizing by convolution
$$
g_\eps := g *_{t,x,\eps} \rho_{\eps}, 
$$
for a  time and space dependent mollifier sequence $\rho_\eps$  similarly as in  \eqref{eq:convolution-for-trace}, we have $g_\eps \in C^1([0,T] \times \bar \OO)$ and 
$$
 \frac{\partial g_\eps }{ \partial t}  + a \cdot \nabla g_\eps + b g_\eps  = R_\eps \ \hbox{ on }\  (0,T) \times \OO, 
 $$
 with the usual commutator $R_\eps := [a \cdot \nabla g+bg,\rho_\eps]$. We know from a classical time and space variant of Lemma~\ref{lem:transport-lemma51} (see also \cite[Lem.~II.1]{MR1022305}, \cite[Lem.~1]{MR1765137}, \cite[Lem.~3.1]{MR2150445} and the proof in \cite[Sec.~3]{MR2340443})  that $R_\eps \to 0$ in $L^1$. 
 Because $g_\eps$ is smooth and $Y$ is an almost everywhere flow, we may set 
$H^\sharp(s,y) :=H(t+s,Y_s(y))$, 
${\mathfrak B} (s,y) :=\int_0^s b(Y_\tau (y)) d\tau$ and compute  
 $$
\frac{d }{ ds} [  (g_\eps^\sharp   e^{{\mathfrak B}}) (s,y)  ] =  (R^\sharp_\eps  e^{{\mathfrak B}}) (s,y) ,
 $$
 from what we get
\bean
&&
\widetilde g_\eps(t,y) :=g_\eps(t,y) - g_\eps(0,  Y_{-t}(y)) e^{{\mathfrak B}(-t,y)} {\bf 1}_{t < t_\bb(y)}   - g_\eps  (t-t_\bb(y),y_\bb(y)) e^{{\mathfrak B}(-t_{\bb(y)},y)}  {\bf 1}_{t > t_\bb(y)}  
\\
&&\qquad =  \int_0^t (R^\sharp_\eps  e^{{\mathfrak B}}) (s,y)  ds  {\bf 1}_{t < t_\bb(y)}+ \int_{t-t_\bb(y)}^t (R^\sharp_\eps  e^{{\mathfrak B}}) (s,y)  ds {\bf 1}_{t >t_\bb(y)}, 
\eean
for a.e. $y \in \OO$ and any $t > 0$ and where we have used \eqref{eq:transport-meas-tbb=t} for getting rid of the set $\{y \in \OO;  t_\bb(y) = t\}$.
For $T,\varrho > 0$ and setting $\UU_{T,\varrho} := (0,T) \times (\OO \cap B_\varrho)$,    we deduce 
\bean
 \int_{\UU_{T,\varrho}}  | \widetilde g_\eps(t,y) | \, dy dt
&\le&  T \int_{\UU_{T,\varrho}}  | (R^\sharp_\eps  e^{{\mathfrak B}}) (s,y) | dy ds
\\
&\le & T e^{(C + \| b \|_{L^\infty})T} \int_{\UU_{T,\varrho}}  | R_\eps (s,y) | dy ds, 
\eean
from the near-incompressibility condition \eqref{eq:transport-aeflow<} of the flow. 
From Proposition~\ref{prop:transport-stability} and Remark~\ref{rem:transport-stability},
we know that 
\bean
&&g_\eps(t,\cdot) \to g(t,\cdot) \ \hbox{ in }\ \Lloc^1(\bar\OO)  \ \hbox{ as } \ \eps \to 0,  \ \hbox{ for any } \ t \in [0,T]; 
\\
&&g_{\eps|\Sigma_-}  \to \gamma_- g = {\mathfrak g}   \ \hbox{ in }\ \Lloc^1( \Sigma_-)  \ \hbox{ as } \ \eps \to 0.
\eean
Passing to the limit, we get 
\bean
 \int_{\UU_{T,\varrho}} | g (t,y) - g_0(Y_{-t}(y)) e^{{\mathfrak B}(-t,y)} {\bf 1}_{t < t_\bb(y)}    - \frakg  (t-t_\bb(y),y_\bb(y) e^{{\mathfrak B}(-t_{\bb(y)},y)}  {\bf 1}_{t > t_\bb(y)}    | dtdy  = 0, 
\eean
for any $T,\varrho > 0$, what is nothing but \eqref{eq:Transport-theo-charac0}. \qed

 \begin{cor}[Representation formula]\label{cor:Transport-Representation}
Under the assumptions of Lemma~\ref{lem:Transport-Uniqueness} in a smooth domain $\OO \not= \R^D$, for any $\lambda > \lambda_{a,b,p}$,  $G \in L^p(\OO)$, 
the unique solution $g \in L^p(\OO)$ to the stationary transport equation \eqref{eq:Transport-PrimalTransport2} (with $\frakg = 0$)
satisfies 
\bear\label{eq:Transport-cor-representation}
g (y)  = \int_0^\infty e^{-\lambda t} (S_b(t) G)(y) dt, \quad \hbox{for a.e } \ y \in \OO, 
\eear
with $S_b$ is defined by \eqref{eq:Transport-defSb} in which formula $Y$  and $t_{\bf b}$ stand for the characteristics and backward exit time defined just as above.  
\end{cor}

\begin{proof}[Proof of Corollary~\ref{cor:Transport-Representation}.]
That is nothing but \eqref{eq:Exist1-DefRepresentationRR}. 
\end{proof}

 Adapting the second proof of  Theorem~\ref{theo:Transport-ODE-DiPLbis}, we obtain a more accurate characterization of the backward exit time $t_\bb$ with more general assumptions on the vector field.

\begin{lem}\label{lem:tbbTER} Assume that $Y$ is an  almost everywhere flow associated to \eqref{eq:Transport-characteristics} with $a$ satisfying \eqref{eq:Transport-ThODE-DiPLAmbrisio}. 
The backward exit time $t_\bb$ is the unique renormalized solution in $L(\OO)$ to the backward exit time problem
\beqn\label{eq:lem:tbbTER}
a \cdot \nabla \tau= 1 \ \hbox{ in } \ \OO, 
\quad
\gamma_- \tau = 0  \ \hbox{ on } \ \Sigma_-. 
\eeqn
We also have 
\beqn\label{eq:lem:tbbTER-2}
  t_\bb = \int_0^\infty S_0(t) {\bf 1} dt, 
\eeqn
where ${\bf1}$ stands for the unit function in $\OO$ and $S_0$ is defined by \eqref{eq:Transport-defSb} with $b = 0$ and in which formula $Y$  and $t_{\bf b}$ stand for the characteristics and backward exit time defined just as above.  
\end{lem}

%

\begin{proof}[Proof of Lemma~\ref{lem:tbbTER}.] {\sl Step 1. Existence.} From  Lemma~\ref{lem:Transport-exist}, for any $\lambda > \lambda^*_\infty := 0$, there exists a (unique) solution $\tau_\lambda \in L^\infty(\OO)$ to the 
truncated backward exit time problem
$$
\lambda \tau_\lambda + a \cdot \nabla \tau_\lambda = 1 \ \hbox{ in } \ \OO, 
\quad
\gamma_- \tau_\lambda = 0  \ \hbox{ on } \ \Sigma_-. 
$$
From the weak maximum principle, we have $\tau_\lambda \ge 0$. As a consequence, for $0 < \lambda < \lambda'$, we have 
$$
\lambda (\tau_\lambda - \tau_{\lambda'}) + a \cdot \nabla (\tau_\lambda - \tau_{\lambda'}) =  (\lambda' - \lambda) \tau_{\lambda'} \ge 0 \ \hbox{ in } \ \OO, 
\quad
\gamma_-  (\tau_\lambda - \tau_{\lambda'}) = 0  \ \hbox{ on } \ \Sigma_-.
$$
From the weak maximum principle again, we deduce that $\tau_\lambda - \tau_{\lambda'}\ge 0$, and $(\tau_{\lambda_n})$ is an increasing sequence when $\lambda_n \searrow 0$. We set $\tau := \lim_{n\to\infty} \tau_{\lambda_n}$, so that $\tau \in L(\OO)$  is a nonnegative renormalized solution to the backward exit time problem.

\smallskip\noindent
{\sl   Step 2. Characterization.} By definition, for any $\beta \in C^1_*(\R)$, $\beta(0) = 0$, the function $\tau$ satisfies 
$$
a \cdot \nabla \beta(\tau) = \beta'(\tau) \ \hbox{ in } \ \OO,
\quad
\gamma_- \beta(\tau)  = 0  \ \hbox{ on } \ \Sigma_-.
$$
With the notations of  \eqref{eq:convolution-for-trace}, we define $b_\eps := \beta(\tau) *_\eps \rho_\eps$, $B_\eps :=  \beta'(\tau) *_\eps \rho_\eps$, and   thanks to Lemma~\ref{lem:transport-lemma51},
we have thus 
$$
a \cdot \nabla b_\eps = B_\eps + r_\eps \ \hbox{ in } \ \DD'(\OO),
$$
with 
$$
b_\eps \to  \beta(\tau), \quad  B_\eps \to \beta'(\tau), \quad r_\eps \to 0,
$$
respectively in $\Lloc^p(\bar\OO)$, $\Lloc^p(\bar\OO)$ and $\Lloc^1(\bar\OO)$ for any $p \in [1,\infty)$. Because then
$$
\frac{d }{ ds} b_\eps \circ Y_s = (a \cdot \nabla b_\eps) (Y_s) = (B_\eps + r_\eps) (Y_s), 
$$
and defining
\bean
\widetilde b_{t,\eps} (y) 
&:=& \Bigl\{ b_\eps(Y_{-t}(y)) + \int_{-t}^0 B_\eps (Y_{s}(y)) \, ds\Bigr\} {\bf 1}_{t < t_\bb(y)}   
\\
&& + \Bigl\{ (\gamma_- b_\eps)(y_\bb(y))) + \int_{-t_\bb(y)}^0   B_\eps (Y_{s}(y)) \, ds\Bigr\}  {\bf 1}_{t   >  t_\bb(y)}  , 
\eean
we have, as in the second proof of  Theorem~\ref{theo:Transport-ODE-DiPLbis}, 
$$
b_\eps(y) - \widetilde b_{t,\eps} (y) = R_{t,\eps} := \int_{-t}^0   r_\eps (Y_{s}(y)) \, ds  {\bf 1}_{t <  t_\bb(y)}  +  \int_{-t_\bb(y)}^0   r_\eps (Y_{s}(y)) \, ds  {\bf 1}_{t   >  t_\bb(y)}.
$$
Arguing similarly as in the second proof of  Theorem~\ref{theo:Transport-ODE-DiPLbis}, we have $\widetilde b_{t,\eps} \to \widetilde b_t$ and $R_{t,\eps} \to 0$ in $L^1(\UU_{\rho,T})$ as $\eps \to 0$, 
with 
$$
\widetilde b_t (y) 
:=  \Bigl\{ \beta(\tau(Y_{-t}(y))) + \int_{-t}^0   \beta'(\tau (Y_{s}(y))) \, ds\Bigr\} {\bf 1}_{t <  t_\bb(y)}   
+  \int_{-t_\bb(y)}^0   \beta'(\tau  (Y_{s}(y))) \, ds   {\bf 1}_{t   >  t_\bb(y)}.
$$
We deduce that 
$$
\beta(\tau(y)) =  \widetilde b_t(y), \quad \hbox{for a.e. } t > 0, \, y \in \OO.
$$
 Choosing a sequence $(\beta_n)$ of renormalizing function in $C^1_*(\R)$ such that $0 \le \beta_n(s) \nearrow s$ and $0 \le \beta'_n(s) \nearrow 1$ locally uniformly, writing the above equation for $\beta = \beta_n$ and passing to the limit $n \to \infty$, we obtain 
$$
\tau(y) =  \Bigl\{  \tau(Y_{-t}(y))  + t \Bigr\} {\bf 1}_{t <  t_\bb(y)}   + t_\bb(y)   {\bf 1}_{t   >  t_\bb(y)}, \quad \forall \, t > 0, 
$$
and in particular $\tau = t_\bb$ a.e. on $\OO$, by taking $t\to\infty$. That implies that $t_\bb$ is the unique renormalized solution to the  backward exit time problem \eqref{eq:lem:tbbTER}. From Corollary~\ref{cor:Transport-Representation}, for any $\lambda > 0$, we have 
$$
\tau_\lambda  = \int_0^\infty e^{-\lambda t} S_0(t) {\bf 1} dt  \quad \hbox{a.e on } \  \OO,
$$
and we deduce \eqref{eq:lem:tbbTER-2} by passing to the limit $\lambda \searrow 0$ in that identity. 
\end{proof}

%
%


\subsection{On the Krein-Rutman theorem for the transport equation with kernel terms}
\label{ssec:Transport:KR}

\

\smallskip
 
In this section we carry on our analysis of the  transport equation with kernel term \eqref{eq:Transport-evolEqWithKernel}-\eqref{eq:Transport-BdaryCond} for which we establish a Krein-Rutman result under strong positivity assumption on the kernel acting on the domain. 
%
%
%
%
As in section~\ref{subsec:TranspEq-tau}, we assume that $a$, $b$, $\KKK$ and $\RRR$ satisfy the conditions  \eqref{eq:TranspEqHypWellPoseWeight},  
  \eqref{eq:Transport-RKhyp2}, \eqref{eq:Transport-RKhyp3}, \eqref{eq:Transport-KKKbound},   \eqref{eq:Transport-RRRbound}
  and \eqref{eq:Transport-WPKKKRRRcondAlphaBeta} for some weight function $m:\bar\OO \to [1,\infty)$ and some exponent $p \in [1,\infty)$. 

\smallskip
On the kernel $\KKK$, we make the additional strong positivity hypothesis: 
for any $x \in \OO$, there exist   $r_x,\eta_x > 0$ such that 
\beqn\label{eq:Transport-hypKK}
\forall \, f \ge 0, \ \forall \, y \in B(x,r_x), \quad \KKK[f](y) \ge \eta_x \int_{B(x,r_x)} f_* dy_*
\eeqn
and 
\beqn\label{eq:Transport-hypab}
\exists \, x_0,  \quad a,b \in L^\infty(B(x_0,r_0)), \quad r_0 := r_{x_0}, 
\eeqn
as well as one of the two following regularity assumptions
\beqn\label{eq:Transport-hypKKbis}
\KKK\in \KK(L^p_m(\OO)) \ \hbox{ or } \ \KKK : L^p_m(\OO) \to L^{p_1}(\OO) \cap L^p_{m_1}(\OO),  
\eeqn
with $p_1 > p$ and $m_1/m \to \infty$ when $y\to\infty$.

 We thus consider the operator 
\beqn\label{eq:Transport-defLL}
\LL f = - a\cdot\nabla f - bf +  \KKK[f] = - \Div (a f) - Kf +  \KKK[f] 
\eeqn
with $K:= b - \Div a \geq0$,
which is complemented with the boundary condition 
\beqn\label{eq:Transport-defBdaryStat}
\gamma_- f  =  \RRR_\OO [f] + \RRR_\Sigma[{\gpf} ] \ \ \hbox{on}  \  \Sigma_-. 
\eeqn
More precisely, we define $\LL$ in the  Banach space $L^p_m(\OO)$ with domain 
$$
D(\LL)\subset  \{ f \in L^1(\OO); \, a \cdot \nabla f \in \Lloc^1(\bar\OO), \  \gamma_- f = \RRR[f,\gamma_+ f]  \}.
$$
Notice that because of Section~\ref{subsec:TranspTrace} the trace function is well defined.

\medskip
{\bf Example.} The nonlocal operator with a drift  
\begin{equation}\label{eq:transport-nonlocalmodel-eq1}
\partial_t f =  - a \,\partial_x f - b + \KKK[f] \  \hbox{ in } \OO, \quad \gamma_- f = 0 \  \hbox{ on } \Sigma_-, 
\end{equation}
with $\OO \subset \R$ a bounded interval, $a \in \Wloc^{1,1}(\OO)$, $a' \in L^\infty(\OO)$, $b \in L^\infty(\OO)$, and thus the boundary kernel is $\RRR \equiv 0$.
Motivated by some non-local reaction-diffusion models, this problem was recently investigated in~\cite{MR4157907,Coville2020,Li2017}.
It is also used in the study of selection-mutation models in changing environment, see the even newer works~\cite{Forien2022,Henry2023}.

%

\medskip

We start by checking that with the above assumptions, the conditions {\bf (H1)}--{\bf (H5)} presented in the abstract part hold true. 

\medskip
{\bf  Condition \ref{H1}.} 
 From Proposition~\ref{prop:Transport-Stat-RKsolutions}, we know that for any $\lambda > \lambda^{**}$ the stationary problem 
 $$
( \lambda - \LL ) g = G \hbox{ in } \OO, \quad \gamma_- g = \RRR[g,\gamma_+ g] \hbox{ on } \Sigma_-, 
$$
has a unique solution. More precisely, the associated inverse operator denoted by $\RR_\LL$ (without reference to the boundary operator $\RRR$) satisfies $\RR_\LL : L^p_m \to L^p_m$ and  $\RR_\LL G \ge 0$ if $G \ge 0$.

\medskip
{\bf Condition \ref{H2}.} We first consider the case when $\RRR_\OO \equiv 0$ and we denote by $\LL_0$ the associated generator.   
%
%
%
We fix $f_0 \in C^2_c(\OO)$, such that  $\int_{\BB_0} f_0 dy =1$, $f_0 > 0$ on $\BB_0$, supp$f_0 = \bar \BB_0$, as well as 
$$
\| f_0^{-1} \|_{L^\infty(\BB_\eps)} = \eps^{-2}, \quad \|   \nabla f_0 \|_{L^\infty(\BB_0 \backslash \BB_\eps)} = \eps, \quad \forall \, \eps \in (0,1/2), 
$$
where we denote  $\BB_\eps := B(x_0,(1-\eps)r_0)$. 
We also define $C_0 :=  \|   f_0 \|_{L^\infty(\BB_0)}$ and  $C_1 :=  \|  \nabla f_0 \|_{L^\infty(\BB_0)}$, both  may be bounded by a constant which only depends on $r_0$ and $d$. 
Because of \eqref{eq:Transport-hypKK}, we have 
\[
  \quad \KKK[f_0](y) \ge \eta_0 {\bf 1}_{\BB_0} \ge \frac{\eta_0 }{ C_0} f_0 .
\]
We observe that $f_0 \in D(\LL_0)$ and 
we compute  
\beqn\label{eq:transport-H2-def1kappa0}
\LL_0 f_0\ge  -  \| a \|_{L^\infty(\BB_0)} C_1 {\bf 1}_{\BB_0}  - \| b \|_{L^\infty(\BB_0)} C_0  {\bf 1}_{\BB_0}  +  \eta_0 {\bf 1}_{\BB_0}
\ge \kappa_0 f_0,
\eeqn
if $\kappa_0 :=  \eta_0/C_0  -  \| a \|_{L^\infty(\BB_0)} C_1/C_0   - \| b \|_{L^\infty(\BB_0)} \ge 0$. More generally, 
we have
\beqn\label{eq:transport-H2-def2kappa0}
\LL_0 f_0 \ge -  \| a \|_{L^\infty(\BB_0)} C_1 {\bf 1}_{\BB_\eps} -   \| a \|_{L^\infty} \eps {\bf 1}_{\BB_0 \backslash \BB_\eps} - \| b \|_{L^\infty(\BB_0)} f_0 +  \eta_0 {\bf 1}_{\BB_0} \ge \kappa_0 f_0, 
\eeqn
with $\kappa_0 :=  - \| a \|_{L^\infty(\BB_0)} C_1 \eps^{-2} - \| b \|_{L^\infty(\BB_0)} \in \R$ when  $ \| a \|_{L^\infty} \eps \le \eta_0$. Depending on how $\eta_0 > 0$ is large, we obtain in that way two constructive lower bounds of $\II$ thanks to Lemma~\ref{lem:Existe1-Spectre2bis}-{\bf (ii)} and we have thus established that   $\LL_0$ satisfies \ref{H2}. Because $f_0 \in D(\LL_0)$, we have 
$S_{\LL_0}(t) f_0 \ge e^{\kappa_0 t} f_0$ for any $t \ge 0$,  from Remark~\ref{rem:Existe1-Spectre2bis}-(2). 
On the other hand, we observe that $S_\LL(t) \ge S_{\LL_0}(t)$ for any $t \ge 0$, from  
the weak maximum principle mentioned in Remark~\ref{rem:Transport-RKsolutions}-(3). 
These two last observations together imply $S_\LL(t) f_0 \ge  e^{\kappa_0 t} f_0$, for any $t \ge 0$. We   deduce from Lemma~\ref{lem:Existe1-Spectre2bis}-{\bf(iv)} that \ref{H2} holds.

%

 \medskip
{\bf Condition \ref{H3}.} 
We introduce the semigroup $S_\BB$ associated to the transport equation 
$$
   \frac{\partial g }{ \partial t}  + a \cdot \nabla g + b g  = 0, 
   \quad
  \gamma_- g = \RRR [\gamma_+g],
 $$
 which is well defined thanks to Corollary~\ref{cor:Transport-RKsolutions} and satisfies  $\| S_\BB(t) g_0 \|_{L^p_m} \le e^{\kappa_\BB t} \|  g_0 \|_{L^p_m}$ for any $t \ge 0$ and $g_0 \in L^p_m$ with 
$ \kappa_\BB :=   \| \langle \varpi_- \rangle \|_{L^\infty} +  M_{_{\!\RRR}}/p$ because of the a priori estimate \eqref{eq:transportLpmApriori-RRR&KKK} particularized to the present case (in particular where we can take $\eps_1 = 0$ because the influx function is $\frakg = 0$ here).
We formulate the first hypothesis 
\beqn\label{eq:Transport-H23-hyp1} 
\eta_0 >   \| \langle \varpi_- \rangle \|_{L^\infty} C_0 + M_{_{\!\RRR}} C_0/p + 
   \| a \|_{L^\infty(\BB_0)} C_1    - \| b \|_{L^\infty(\BB_0)} C_0,
 \eeqn
with the same definitions as above for $\BB_0$, $C_0$ and $C_1$, so that $\kappa_0 > \kappa_\BB$ because of \eqref{eq:transport-H2-def1kappa0}. In a second case, we assume 
\beqn\label{eq:Transport-H23-hyp2} 
\RRR \equiv 0, \ \OO \ \hbox{is bounded  and  there exists} \ T_\OO \ \hbox{such that} \ t_{\bf b} (y) \le T_\OO \ \hbox{for a.e.} \ y \in \OO.
 \eeqn
In that case, the semigroup $S_\BB$ is explicitly given by
\[(S_\BB(t) f_0) (y)
= \left\{
\begin{aligned}
&   f_0( Y_{-t}(y))\exp(-\int_0^t  K( Y_{\tau-t}(y)) d\tau) , \ \hbox{ if } t \in   (0,t_\bb(y)) \\ 
&0 \ \hbox{ otherwise}, 
\end{aligned}
\right.\]
and in particular $S_\BB(t)f=0$ for any $f$ and any $t>T_\OO$. We immediately deduce $\kappa_\BB = - \infty$ and thus $\kappa_0 > \kappa_\BB$ because we have established that $\kappa_0 \in \R$. 
We next define $\AA f := \KKK[f]$. Using Lemma~\ref{lem:H3abstract-StrongC}  and Remark~\ref{rem:H3abstract-compactStrong}-(2)  or  Lemma~\ref{lem:H3Lp}   and Remark~\ref{rem:H3Lp}-(1)  depending on the assumption \eqref{eq:Transport-hypKKbis} made on $\KKK$, we deduce that the condition \ref{H3} holds in  both cases discussed above. Under the first condition in \eqref{eq:Transport-hypKKbis}, we conclude 
to the existence of eigenvalue triplet $(\lambda_1, f_1, \phi_1) \in \R \times L^p_m \times L^{p'}_{m^{-1}}$. 
Under the second condition in \eqref{eq:Transport-hypKKbis}, we may also get the same conclusion by by using \cite[Cor.~1 of Thm.~II.9.9]{MR0423039} when $p=1$ or by observing that the dual problem is similar to the primal problem when $p > 1$ and thus we may apply the same arguments for the dual problem as those  explained above  for the primal problem. 

\medskip

{\bf Condition \ref{H4}.} 
Let us consider $\lambda > \lambda^{**}$ and  $0 \le f \in L^p_{m_\OO}(\OO)$ a (renormalized) solution  to 
$$
\lambda f + a \cdot \nabla f + b f -  \KKK[f] = F \quad\hbox{in}\quad  \OO, 
\quad 
\gamma_- f = \RRR  [f ,\gpf ] \quad\hbox{on}  \quad \Sigma_-, 
$$
with $0 \le F \in L^p_m(\OO)$. If $f\not\equiv 0$, there exists $x_1 \in \OO$ such that $\int_{B(x_1,r_{_1})} f(z) \, dz > 0$. 
From \eqref{eq:Transport-hypKK}, we deduce
$$
\KKK[f](y) \ge \int_{B(x_1,r_1)} \kappa(y,z) f(z) \, dz > 0, \quad \forall \, y \in B(x_1,r_1).
$$
Now, we argue similarly as during the proof of Lemma~\ref{lem:Transport-Uniqueness} and in particular we use  the same notations.
For  $A \subset B(x_1,r_1)$, we define the solution $0 \le \varphi  \in L^{p'}_{m^{-1}}  \cap L^\infty$ to the equation
$$
\lambda \varphi - \Div (a\varphi) + b \varphi  =  {\bf 1}_{A} \ \ \hbox{in}\  \OO, \quad \gamma_+ \varphi  =   0 \ \ \hbox{on}  \  \Sigma_+,
$$
 thanks to Lemma~\ref{lem:Transport-exist} and Lemma~\ref{lem:Transport-existp1}, and we observe that $\varphi \not\equiv 0$ on $B(x_1,r_1)$  if $|A| > 0$. 
For the renormalizing function $\beta_\delta$ and a truncation function $\chi_R$, we compute  
\bean
0 &\ge& \int_\Sigma a\cdot n \beta_\delta(\gamma f)  \gamma \varphi \chi_R
\\
&=& 
 \int_\OO [\beta'_\delta(f) (F+\KKK[f]) \varphi  - \beta_\delta(f) {\bf 1}_A ]    \chi_R 
 \\
 &&\qquad + 
 \int_\OO (\beta_\delta(f) - f \beta'_\delta(f)) (\lambda+b)  \varphi  \chi_R
+  \int_\OO  \varphi \beta(f) \frac{a }{ R} \cdot (\nabla \chi)_R.
\eean
%
%
%
Passing first to the limit $R \to \infty$ and next to the limit $\delta \to 0$, we deduce  
\bean
0 \ge  \int_\OO [  (F + \KKK[f]) \varphi - f {\bf 1}_A  ], 
\eean
so that in particular 
\bean
\int_A f \, dy \ge  \int_{B(x_1,r_1)}  \varphi  \KKK[f] > 0.
\eean
This being true for any  $A \subset B(x_1,r_1)$, we deduce $f > 0$ a.e. on  $B(x_1,r_1)$.  By a classical 
continuity argument, we conclude that $f > 0$ a.e. on $\OO$. We have thus established \ref{H4} for $\lambda > \lambda^{**}$ from what we immediately
and classically deduce the general case $\lambda \in \R$. 
 

\medskip

{\bf Condition \ref{H5}.} Assume that $(\lambda,f) \in \C \times D(\LL)$ satisfies  
%
 \bean
  \LL |f|  =(\Re e \lambda) |f|\quad\hbox{in}\quad  \OO, \quad  \RRR  [|f| ,\gp |f| ] =   \gamma_- |f| \quad\hbox{on}  \quad \Sigma_-, 
\eean
and
$$
 \LL |f|  =   \Re e (\hbox{\rm sign} f) \LL f \quad\hbox{in}\quad  \OO, \quad  \RRR  [|f| ,\gp|f| ] = \Re e (\hbox{\rm sign} \gamma_-f) \RRR [f,\gpf]  \quad\hbox{on}  \quad \Sigma_-.
$$ 
From \ref{H4} and the first identity, we know that $|f| > 0$ a.e. on $\OO$. 
Using the second identity,  we get
$$
\KKK[|f|] =  \Re e (\hbox{\rm sign} f) \KKK[f].
$$
Writing $f = e^{i\alpha} |f|$, we deduce 
$$
\int_\OO k |f_*| (1 - \cos (\alpha-\alpha_*)) dy_* = 0 \quad \hbox{a.e. on } \OO.
$$
Using \eqref{eq:Transport-hypKK}, we deduce 
$$
\int_{B(y,r_\OO)} |f_*| (1 - \cos (\alpha-\alpha_*)) dy_* = 0, 
$$
and thus $\alpha=\alpha_*$ a.e. on $\OO \times \OO$. That means $f = u |f|$, for a constant $u = \Sp^1$,  that completes the proof of 
the fact that  $\LL$ satisfies the reverse Kato's inequality condition \ref{H5}.

\medskip
We summarize our analysis in the following result which is a straightforward consequence of the above checked conditions together with  Theorem~\ref{theo:exist1-KRexistence}, Theorem~\ref{theo:KRgeometry1}, Theorem~\ref{theo:KRgeometry2}   and Theorem~\ref{theo:ergodicity-compact-trajectories}.
We state the available result in that situation. \Black

\begin{theo}\label{th:KR-transport-kernel} We assume that $a$, $b$, $\KKK$ and $\RRR$ satisfy the conditions  \eqref{eq:TranspEqHypWellPoseWeight},  
  \eqref{eq:Transport-RKhyp2}, \eqref{eq:Transport-RKhyp3}, \eqref{eq:Transport-KKKbound},   \eqref{eq:Transport-RRRbound}
  and \eqref{eq:Transport-WPKKKRRRcondAlphaBeta} for some weight function $m:\bar\OO \to [1,\infty)$ and some exponent $p \in [1,\infty)$. 
Consider the semigroup $S_\LL$ associated to the transport equation \eqref{eq:Transport-evolEqWithKernel}-\eqref{eq:Transport-BdaryCond} through Corollary~\ref{cor:Transport-RKsolutions}.
We assume further that $\KKK$ satisfies the strong positivity conditions \eqref{eq:Transport-hypKK} together with \eqref{eq:Transport-hypab} and the first compactness property formulated in \eqref{eq:Transport-hypKKbis}. We finally assume that  \eqref{eq:Transport-H23-hyp1} holds or \eqref{eq:Transport-H23-hyp2} holds. In both cases, the conclusions {\Blue \ref{S1}, \ref{S2} and \ref{S32} hold in $L^p_m$} as well as the ergodicity \ref{E2} in $L^1_{\phi_1}$.
\end{theo}

 We are not aware of any similar result for such a  general transport equation, see however the next sections where more specific transport like equations are discussed.   
 We do not try to improve the convergence result in the general case, but rather we aim to make one step further in the following particular situation where Doblin approach may be used. 
 
 \medskip
 
{\bf Doblin condition.}
We suppose here that $\OO$ is bounded, $K\in L^\infty(\OO)$,  $\RR_\Sigma^*\1=\RR_\Sigma\1=\1$,  and $k(y,y_*)\geq k_0>0$.
We aim at establishing the Doblin condition
$$
S_\LL(T) f_0 \ge \kappa \langle f_0, \1 \rangle, 
$$
which is~\eqref{eq:hyp-Doblin} with $\psi_0=\1$ and $g_0=\kappa\1$. 
From~\eqref{eq:transportLpmApriori-RRR&KKK-diff} we have
\[\frac{d }{ dt} \int_\OO f dy =  \int_\OO f K dy \ge -\|K\|_\infty \int_\OO f dy\]
and so
\[\int_\OO f(t,y)\,dy\geq e^{-\|K\|_\infty t}\int_\OO f_0(y)\,dy.\]
Now we define, for $\varphi_0\in C^1_c(\OO)$, $\varphi_0\geq0$, $\int\varphi_0=1$, the solution $\varphi$ to the equation
\[\left\{\begin{array}{l}
\partial_t\varphi+\Div(a\varphi)=0,
\vspace{2mm}\\
\gamma_+\varphi=\RR_\Sigma^*[\gamma_-\varphi].
\end{array}\right.\]
 We have
 \[\frac{d}{dt}\int_\OO \varphi = 0,\qquad\text{and so}\quad \int_\OO\varphi(t,y)\,dy=\int_\OO\varphi_0(y)\,dy=1,\]
 and
\[
 \frac{d}{dt}\int_\OO f\varphi =\int_\OO\KKK[f]\varphi-\int_\OO Kf\varphi\ge k_0\int_\OO f - \|K\|_\infty \int_\OO f\varphi.
\]
We deduce from Grönwall's inequality that, for any fixed $T>0$,
\begin{align*}
\int_\OO f(T,y)\varphi_0(y)\,dy&\geq e^{-\|K\|_\infty t}\int_\OO f_0(y)\varphi_0(y)\,dy+k_0\int_0^T e^{-(T-t)\|K\|_\infty}\int_\OO f(t,y)\,dydt\\
&\geq k_0 T e^{-T\|K\|_\infty}\int_\OO f_0(y)\,dy =: \kappa \langle f_0,\1\rangle.
\end{align*}
This is nothing but the Doblin condition \eqref{eq:hyp-Doblin} since $\varphi_0$ is any non-negative function in $C^1_c(\OO)$ with $\int\varphi_0=1$.

\medskip

In order to verify~\eqref{eq:hyp-Doblin-BdBis} in a quantitative way, we suppose that the conditions  \eqref{eq:TranspEqHypWellPoseWeight},  \eqref{eq:Transport-RKhyp2}, \eqref{eq:Transport-RKhyp3}, \eqref{eq:Transport-KKKbound},   \eqref{eq:Transport-RRRbound}
and \eqref{eq:Transport-WPKKKRRRcondAlphaBeta} are verified with the weight function $m=\1$ and the exponent $p=1$.
Note that in this case we have $\varpi=K\geq0$.
The first condition in~\eqref{eq:TranspEqHypWellPoseWeight} then imposes that $K\in L^\infty(\OO)$, and~\eqref{eq:Transport-H23-hyp1} reads
  \begin{equation*}
\eta_0 >   \| \langle K \rangle \|_{L^\infty} C_0 + M_{_{\!\RRR}} C_0 + 
   \| a \|_{L^\infty(\BB_0)} C_1    - \| b \|_{L^\infty(\BB_0)} C_0.
\end{equation*}
We also assume that
\begin{equation}\label{as:R1}
\RR_\Sigma^*\1=\RR_\Sigma\1=\1,
\end{equation}
and
\begin{equation}\label{as:kminmax}
\forall y,y_*\in\OO,\quad k_0\leq k(y,y_*)\leq k_1
\end{equation}
for some $k_1>k_0>0$,

 \begin{theo}\label{th:KR-transport-kernel-Doblin} 
 We assume that $\OO$ is bounded and that the conditions  \eqref{eq:TranspEqHypWellPoseWeight},  
  \eqref{eq:Transport-RKhyp2}, \eqref{eq:Transport-RKhyp3}, \eqref{eq:Transport-KKKbound},   \eqref{eq:Transport-RRRbound}
  and \eqref{eq:Transport-WPKKKRRRcondAlphaBeta} are satisfied by  $a$, $b$, $\KKK$ and $\RRR$
  for the weight function $m=\1$ and the exponent $p=1$. 
We assume further that $\KKK$ satisfies the strong positivity conditions \eqref{eq:Transport-hypKK} together with \eqref{eq:Transport-hypab} and the first compactness property formulated in \eqref{eq:Transport-hypKKbis}.
We finally assume that  \eqref{eq:Transport-H23-hyp1}, \eqref{as:R1} and \eqref{as:kminmax} are satisfied.
Then the exponential convergence in~\ref{E31} holds in $L^1$ with constructive constants $C$ and $\omega$.
 \end{theo}

 
\begin{proof}[Proof of Theorem~\ref{th:KR-transport-kernel-Doblin}]
We work in $X=L^1(\OO)$ and we normalize $\phi_1$ by $\|\phi_1\|_{L^\infty}=1$.
We have proved above that~\eqref{eq:hyp-Doblin} holds true with $\psi_0=\1$ and $g_0=\kappa\1$ for some explicit $\kappa>0$, recalling that the assumption that $K\in L^\infty$ is nothing but the first condition in~\eqref{eq:TranspEqHypWellPoseWeight} when $m=\1$ and $p=1$ since $b=K+\Div a$.
Due to the normalization $\|\phi_1\|_{L^\infty}=1$, the condition~\eqref{eq:hyp-Doblin-Bd} holds with $R_0=1$.
It only remains to check the validity of~\eqref{eq:hyp-Doblin-BdBis} in order to be able to apply Theorem~\ref{theo:Doblin}.
Since we assume that the conditions  \eqref{eq:TranspEqHypWellPoseWeight}, \eqref{eq:Transport-RKhyp2}, \eqref{eq:Transport-RKhyp3}, \eqref{eq:Transport-KKKbound},   \eqref{eq:Transport-RRRbound}  and \eqref{eq:Transport-WPKKKRRRcondAlphaBeta} are satisfied for the weight function $m=\1$ and the exponent $p=1$, we have that $\RR_\BB(\lambda_1):L^1\to L^1$ with $\|\RR_\BB(\lambda_1)\|_{\BBB(L^1)}\leq \frac{1}{\kappa_0-\kappa_\BB}$.
This yields by duality that $\RR^*_\BB(\lambda_1):L^\infty\to L^\infty$ with $\|\RR^*_\BB(\lambda_1)\|_{\BBB(L^\infty)}\leq \frac{1}{\kappa_0-\kappa_\BB}$.
  Since $k$ is bounded by the constant $k_1$, we have on the other hand that $\AA^*=\KKK^*:L^1\to L^\infty$ with $\|\AA^*\|_{\BBB(L^1,L^\infty)}\leq k_1$.
  We thus get
  \[1=\|\phi_1\|_{L^\infty}\leq \frac{k_1}{\kappa_0-\kappa_\BB}\|\phi_1\|_{L^1},\]
  which yields~\eqref{eq:hyp-Doblin-BdBis} with $r_0=\kappa(\kappa_0-\kappa_\BB)/k_1$, and the proof is complete.
\end{proof}

\medskip

%

 \subsection{A word about the renewal equation}
\label{ssec:Transport:renewal}

We look at the case $\OO = (0,+\infty)$ and $a(y) = 1$, which corresponds to the equation
\begin{equation}\label{eq:renewal}
\partial_t f + \partial_y f + K f = 0 
\end{equation}
with the boundary condition
\begin{equation}\label{eq:renewal-boundary}
(\gamma_-f)(t,0) = \int_0^\infty r_\OO(y_*) f(t,y_*) dy_*.
\end{equation}
This renewal age structured model is standard in structured population dynamics, and the Krein-Rutman theorem is well-known for it, see for instance~\cite{BCG2019,Feller1941,MR763356,Gwiazda2006,Sharpe1977,Webb1984}.
The existence and uniqueness of $(\lambda_1,f_1,\phi_1)$ can even be obtained by explicit computations.
However, it is not covered by the cases considered in Section~\ref{ssec:Transport:KR} because $\KKK = 0$ here.

\medskip

The singularity of this transport equation lies in the fact that~\ref{H2} is only guaranteed by the boundary condition.
To fall into our splitting framework, we may replace the boundary condition by a singular source term $\AA f=\RRR_\OO[f](0)\delta_0$, where $\delta_0$ is the Dirac mass at the origin, and write $\LL=\AA+\BB$ with $\BB$ the generator of the free transport equation with zero flux boundary condition.
This forces working in a space of measures, as in~\cite{MR3859527,MR3850019}.
We briefly present an alternative approach,  which  is more in the spirit of~\cite{BCG2019,MR3893207} and which consists in working in the Lebesgue space $L^1$, first to solve the dual problem in $L^\infty=(L^1)'$ and next to use for instance Doblin's contraction to solve the primal problem.

\smallskip
We assume here that
\begin{align} \label{as:renewal:wellposedness}
&0 \le K,r_\OO\in \Lloc^\infty(0,\infty),\quad (r_\OO-\alpha K)_+ \in L^\infty, 
\\ \label{as:renewal:limK}
&\lim_{y\to\infty}K(y)=+\infty, \quad \liminf_{y\to+\infty}r_\OO(y)>0, 
\end{align}
for some $\alpha \in (0,1)$, and we verify the usual conditions for the direct or the dual problem. 

\smallskip
{\bf Condition \ref{H1}.} 
Under assumption~\eqref{as:renewal:wellposedness}, the age structured equation~\eqref{eq:renewal}-\eqref{eq:renewal-boundary} is well-posed in $L^1$ thanks to Proposition~\ref{prop:Transport-RKsolutions}
and we may associate to it a positive semigroup $S_\LL$ in $L^1$ with growth bound $\omega(S_\LL) \le \kappa_1 := \| (r_\OO-K)_+ \|_{L^\infty}$  thanks to Corollary~\ref{cor:Transport-RKsolutions}. We deduce that \ref{H1} holds for the primal problem and thus also for the dual problem thanks to Lemma~\ref{lem:Exist1-RkSG} and 
Lemma~\ref{lem:Existe1-H2bis}.

\medskip

\smallskip
{\bf Condition \ref{H2}.} The generator of the dual problem is
\[
\LL^*\phi=\partial_y \phi - K(y)\phi + \phi(0)r_\OO(y)
\]
with domain $D(\LL^*) \subset \Wloc^{1,\infty}(\OO)$. From the second  hypothesis in \eqref{as:renewal:limK}, there exist $y_0, \eta_0 \in (0,\infty)$ such that 
$r_\OO(y) \ge \eta_0$ for any $y \ge y_0$. We then define 
$$
\phi_0 (y) = {\bf 1}_{[0,y_0)} + \eta_0  (y_1 - y)  {\bf 1}_{[y_0,y_1)} , \quad y_1 := y_0 + 1/\eta_0,
$$
and we compute 
\bean
\LL^*\phi_0 &=& r_\OO(y) - K \ge - \| K - r_\OO \|_{L^\infty(0,y_0)} \quad\hbox{on} \ (0,y_0), 
\\
\LL^*\phi_0 &=& r_\OO(y) - \eta_0 - K \phi_0 \ge - \| K  \|_{L^\infty(y_0,y_1)} \phi_0 \quad\hbox{on} \ (y_0,y_1), 
\\
\LL^*\phi_0 &=& 0 \quad\hbox{on} \ (y_1,\infty), 
\eean
so that in the three case $\LL^* \phi_0 \ge \kappa_0 \phi_0$ with $\kappa_0 := - \max( \| K - r_\OO \|_{L^\infty(0,y_0)} , \| K  \|_{L^\infty(y_0,y_1)})$. 
Using Lemma~\ref{lem:Existe1-Spectre2bis}-{\bf (i)}, we have thus established that   $\LL$ satisfies \ref{H2} with constructive constant $\kappa_0$.

\medskip
{\bf Condition \ref{H3} on the dual problem.} We define the splitting $\LL^*=\AA^*+\BB^*$ with $\AA^*\phi :=(\RR_\OO^*\phi)(y)=\phi(0)r_\OO(y)$.
From the first hypothesis in \eqref{as:renewal:limK}, for any $\kappa_* \le 0$ there exists $y_* \in [0,\infty)$ such  that $K(y) \ge - \kappa_*$ for any $y \ge y_*$. 
Defining $m_* := e^{\kappa_* y} {\bf 1}_{[0,y_*)} + e^{\kappa_* y_*} {\bf 1}_{[y_*,\infty)}$, we compute 
$$
\BB^* m_* = \kappa_* e^{\kappa_* y} {\bf 1}_{[0,y_*)} - K m_* \le \kappa_* m_*.
$$
Together with Proposition~\ref{prop:Transport-RKsolutions} and Corollary~\ref{cor:Transport-RKsolutions}, we deduce that the operator $\BB - \kappa_*$, with 
domain $D(\BB) := \{ f \in L^1(\OO); \ \partial_y f + K f \in L^1(\OO), \ f(0) = 0 \}$, 
generates a contraction semigroup in $L^1_{m_*}(\OO)$, and thus a bounded semigroup in $L^1(\OO)$ because $m_*,m_*^{-1} \in L^\infty(\OO)$. In other words,
we have established that $\omega(S_\BB) = - \infty$. Now, we see that $\RR_{\BB^*} (\lambda) : L^\infty \to D(\BB^*) \subset \Wloc^{1,\infty}([0,\infty))$ is bounded for any $\lambda \in \R$
and thus $\AA^* \RR_{\BB^*} (\lambda) : L^\infty \to L^\infty$ is compact for any $\lambda \in \R$. We deduce from Lemma~\ref{lem:H3abstract-StrongC} and Remark~\ref{rem:H3abstract-compactStrongAlternative},
that $\LL^*$ satisfies \ref{H3}.

 \medskip
Using Lemma~\ref{lem:H3abstract-StrongC}-{\it(1)}, we  conclude to the existence of $(\lambda_1,\phi_1)$ solution to the dual eigenvalue problem.
Now, we turn to the existence, uniqueness, and exponential stability of $f_1\in L^1$, by verifying that Doblin's condition~\eqref{eq:hyp-Doblin} is satisfied.

\medskip
{\bf Doblin condition.}
Denoting $S_t := S_\LL(t)$,  we have from the characteristics method
\[S_tf(y)=f(y-t)e^{-\int_0^t K(y-s)ds}\1_{t<y}+N(t-y)e^{-\int_0^y K(s)ds}\1_{t> y}\]
with $N(t)=\int_0^\infty r_\OO(y_*)S_tf(y_*)dy_*$.
Iterating this formula and using the positivity of $S_t$ we get that for any $f\geq0$
\begin{align*}
S_tf(y)&\geq \bigg(\int_0^{t-y} r_\OO(y_*)N(t-y-y_*)e^{-\int_0^{y_*}K(s)ds}dy_*\bigg)e^{-\int_0^y K(s)ds}\1_{0<y<t}\\
&\geq \bigg(\int_0^{t-y} r_\OO(t-y-\tau)N(\tau)e^{-\int_0^{t-y-\tau}K(s)ds}d\tau\bigg)e^{-\int_0^y K(s)ds}\1_{0<y<t}.
\end{align*}
%
%

Choosing $t_0 > 2y_0$ so that $r_\OO(y)\geq \eta_0>0$ for all $y\geq t_0/2$, 
we obtain
\[S_{t_0}f(y)\geq \eta_0\,e^{-\int_0^{t_0}K(s)ds}\bigg(\int_0^{t_0/4}N(\tau)d\tau\bigg)e^{-\int_0^y K(s)ds}\1_{0<y<t_0/4}.\]
From the expression of $N(t)$, we get by duality, using that $r_\OO\geq\eta_0\1_{(y_0,\infty)}$, that
\begin{align}\label{eq:Doblin-renewal1}
S_{t_0}f(y)
&\geq \eta_0^2\,e^{-2\int_0^{t_0}K(s)ds}\bigg(\int_0^\infty f(y_*)\bigg(\int_0^{t_0/4} S_\tau^*\1_{(y_0,\infty)}(y_*)d\tau\bigg) dy_*\bigg)\1_{0<y<t_0/4}.
\end{align}
Applying $S_{t_1}$ to this inequality we deduce that for any $t_1>0$
\begin{align*}
S_{t_0+t_1}f(y)&\geq \eta_0^2\,e^{-2\int_0^{t_0}K(s)ds}\bigg(\int_0^\infty f(y_*)\bigg(\int_0^{t_0/4} S_\tau^*\1_{(y_0,\infty)}(y_*)d\tau\bigg) dy_*\bigg)S_{t_1}\1_{0<y<t_0/4}\\
&\geq \eta_0^2\,e^{-2\int_0^{t_0}K}e^{-\int_0^{t_1}K}\bigg(\int_0^\infty f(y_*)\bigg(\int_0^{t_0/4} S_\tau^*\1_{(y_0,\infty)}(y_*)d\tau\bigg) dy_*\bigg)\1_{t_1<y<t_0/4+t_1}.
\end{align*}
On the other hand, replacing $f$ by $S_{t_1}f$ in~\eqref{eq:Doblin-renewal1} we obtain
\begin{equation}\label{eq:Doblin-renewal2}
S_{t_0+t_1}f(y)\geq \eta_0^2\,e^{-2\int_0^{t_0}K(s)ds}\bigg(\int_0^\infty f(y_*)\bigg(\int_0^{t_0/4} S_{\tau+t_1}^*\1_{(y_0,\infty)}(y_*)d\tau\bigg) dy_*\bigg)\1_{0<y<t_0/4}.
\end{equation}
The fact that $S^*_t\phi(y)\geq \phi(t+y)e^{-\int_0^t K(y+s)ds}$ ensures that for $t_1>y_0$
\[S_{t_1}^*\1_{(y_0,\infty)}\geq e^{-\sup_{y\in[0,y_0]}\int_0^{t_1}K(y+s)ds}\1_{[0,y_0]}.\]
All together, we have proved that for any $t_0>4t_1>4y_0$ we have
\[S_{t_0+t_1}f(y)\geq c_0\bigg(\int_0^\infty f(y_*)\bigg(\int_0^{t_0/4} S_{\tau}^*\1(y_*)d\tau\bigg) dy_*\bigg)\1_{t_1<y<t_0/4}\]
for some explicit constant $c_0$ and all $f\geq0$.
This is Doblin's condition~\eqref{eq:hyp-Doblin} with $T=t_0+t_1$, and the functions $\psi_0=\int_0^{t_0/4} S_\tau^*\1\,d\tau$ and $g_0=c_0\1_{(t_1,t_0/4)}$.
We are now in position to prove the following result. 
 
 \medskip
 
 \begin{theo}\label{th:KRrenewal}
 Under the assumptions~\eqref{as:renewal:wellposedness} and \eqref{as:renewal:limK}, 
 the renewal equation~\eqref{eq:renewal}-\eqref{eq:renewal-boundary} enjoys the conclusions  {\Blue \ref{S1}, \ref{S2} and \ref{S32} in $L^1$ as well as  \ref{E31} with quantitative rate 
 in $L^1_{\psi_0}$.}
 \end{theo}
 
 Despite the numerous results about the renewal age-structured model, we are not aware of any previous result with a constructive rate of convergence under such general assumptions.

 \begin{proof}[Proof of Theorem~\ref{th:KRrenewal}]
 The conditions~\ref{H1}, \ref{H2} and \ref{H3} for $\LL^*$ ensure the existence of $\lambda_1\geq\kappa_0$ and $\phi_1\in L^\infty$, $\phi_1>0$, that we normalize by $\|\phi_1\|_{L^\infty}=1$.
 If we can prove that the conditions~\eqref{eq:hyp-Doblin}, \eqref{eq:hyp-Doblin-BdBis} and \eqref{eq:hyp-Doblin-Bd} are verified, then the conclusions  {\Blue \ref{S1}, \ref{S2} and \ref{E31} follow by applying Theorem~\ref{theo:Doblin}}. 
We have already proved~\eqref{eq:hyp-Doblin} with the functions $\psi_0=\int_0^{t_0/4} S_\tau^*\1\,d\tau$ and $g_0=c_0\1_{(t_1,t_0/4)}$.
For proving \eqref{eq:hyp-Doblin-BdBis}, we start by recalling that $\phi_1=\RR^*_\BB(\lambda_1)\AA^*\phi_1\in W^{1,\infty}_{loc}$ due to the informations derived on $\RR^*_\BB$ in \ref{H3}.
Consequently, there exists $y_1>0$ such that $\phi_1(y_1)>1/2$, and we deduce from $\phi_1'\leq(\lambda_1+K)\phi_1$ that
\[\phi_1(y)\geq \frac12e^{-\int_y^{y_1}(\lambda_1+K)}\]
for all $y\in(0,y_1)$.
Choosing in the proof of the Doblin condition $t_0$ such that $y_1<t_0/4$, we obtain that
\[\langle\phi_1,g_0\rangle\geq \frac{c_0}{2}\int_{y_1/2}^{y_1}e^{-\int_y^{y_1}(\lambda_1+K)}dy,\]
which gives \eqref{eq:hyp-Doblin-BdBis}.
For \eqref{eq:hyp-Doblin-Bd}, we use that
\[\phi_1=e^{-\lambda_1\tau}S^*_\tau\phi_1\leq e^{-\lambda_1\tau}S^*_\tau\1\]
for any $\tau>0$ to deduce that
\[\phi_1=\frac{4}{t_0}\int_0^{t_0}e^{-\lambda_1\tau}S^*_\tau\phi_1d\tau\leq\frac{4e^{|\lambda_1|t_0}}{t_0}\psi_0.\]
Finally, we check that the condition~\ref{H5'} is verified, so that~\ref{S32} is valid by virtue of Theorem~\ref{theo:lambda1largestBIS} and Remark~\ref{rem:H5'}.
The condition~\ref{H5'} is actually a direct consequence of the fact that~\eqref{eq:Doblin-renewal2} is verified for any $t_0>2y_0$ and $t_1>0$ together with the estimate
\[S_{\tau+t_1}^*\1_{(y_0,\infty)}(y)\geq e^{-\int_0^{t_1}K(y+s)ds}>0\]
for any $t_1>y_0$ and $\tau>0$, and all $y>0$.
 \end{proof}
 

%
%


\bigskip\bigskip
\section{The growth-fragmentation equation}
\label{part:application3:GF}


In this section, we are interested in the growth-fragmentation equation with equal mitosis kernel
\begin{equation}\label{eq:mitosis}
\partial_t f (t,x) + \partial_x \big(a(x)f (t,x)\big) +  K(x) f(t,x) = 4 K(2x) f(t,2x)
\end{equation}
and to its variant with an additional ``growth speed'' variable
\begin{equation}\label{eq:mitosis-var}
\partial_t f (t,x,v) + v \partial_x \big(a(x)f (t,x,v)\big) + K(x) f(t,x,v) = 4 \int_1^2 K(2x) \wp(v,v_*) f(t,2x,v_*) dv_*,
\end{equation}
with $x>0$ and $v\in[1,2]$.
For both equations, we assume that the total fragmentation rate $K$ is a continuous function defined on $\R_+$ such that
\beqn\label{eq:hypKmitose1}
\exists x_0>0,\quad K=0\ \text{on}\ (0,2x_0]\quad \text{and}\quad K>0\ \text{on}\ (2x_0,\infty).
\eeqn
This condition ensures that no particle of size less than $x_0$ can be produced by division, and we thus consider the equations posed on the size space $(x_0,\infty)$ with zero flux boundary condition $f(t,x_0)=0$ or $f(t,x_0,v)=0$.
The growth rate $a$ is supposed to be positive and globally Lipschitz on $[x_0,\infty)$, and we assume that
\begin{equation}\label{as:Kb1}
\lim_{x\to\infty}\frac{xK(x)}{a(x)}=+\infty.
\end{equation}
For quantifying the positivity of the first eigenvalue, we also make the technical assumption that
\beqn\label{as:Kinf}
\exists k>0,\quad \lim_{x\to\infty}e^{x^k}K(x)=+\infty.
\eeqn

\subsection{The mitosis equation with mixing growth rate}\label{ssec:GF:regular}

We are interested here in the growth-fragmentation equation~\eqref{eq:mitosis} in the case where
\begin{equation}\label{as:mixingb}
\exists x_1>x_0,\quad a(2x_1)\neq2a(x_1).
\end{equation}
As we will see below, this condition ensures some mixing property for the trajectories that guarantees the triviality of the boundary point spectrum.

\smallskip

We work in the space $X=L^1_m$ with a weight $m$ that can be
\begin{equation}\label{def:GF-m}
\text{either}\quad m(x)=x^r,\ \ r>1,\qquad\text{or}\quad m(x)=\exp\Big(\eta\int_{x_0}^x \frac{K}{a}\Big),\ \ 0<\eta<1.
\end{equation}
Note that due to assumption~\eqref{as:Kb1}, the weight $\exp\big(\eta\int_{x_0}^x K/a\big)$ is always stronger than $x^r$.

\begin{theo}\label{theo:GF-regular}
Suppose that~\eqref{eq:hypKmitose1}, \eqref{as:Kb1}, \eqref{as:Kinf} and \eqref{as:mixingb} are satisfied.
The first eigentriplet problem admits a unique solution $(\lambda_1,f_1,\phi_1)\in\R\times X_+\times X'_+$ with the normalization $\|\phi_1\|=\langle\phi_1,f_1\rangle=1$, and this triplet additionally satisfies $\lambda_1>0$, $f_1> 0$ and $\phi_1> 0$.
Besides, there are some constructive constants $C \ge 1$, $\omega>0$ such that
\[ {\| e^{-\lambda_1t} S_\LL(t) f - \langle \phi_1 , f \rangle f_1 \|}_X \leq C e^{-\omega t} {\| f - \langle \phi_1 , f \rangle f_1 \|}_X \]
for any $f\in X$ and $t\geq0$. {\Blue In other words, the conclusions \ref{S1},  \ref{S2},  \ref{S33} and  \ref{E31} hold with constructive constants in $L^1_m$.}
\end{theo}

This result is contained in the recent paper~\cite{MR4897760}.
The novelty here is that all the constants are obtained constructively, which is not clear in~\cite{MR4897760}.
We also provide what seems to us to be a more direct and comprehensive proof.
We also refer to~\cite{Bansaye2022,BernardGabriel20,MR4312814,MR3489637} where the same result is obtained under stronger assumptions.

\

Before starting the proof of Theorem~\ref{theo:GF-regular}, let us briefly justify the relevance of the chosen weight functions $m$ in~\eqref{def:GF-m}.
The dual operator associated to equation~\eqref{eq:mitosis} is given by
$$
\LL^*\phi(x) = a(x)\phi'(x) - K(x) \phi(x) + 2K(x) \phi(x/2).
$$
For $m(x)=x^r$, $r>1$, we can compute
\begin{equation}\label{eq:GF-Lyap-xr-1}
\LL^*m(x)=\bigg[r\frac{a(x)}{x}-(1-2^{1-r})K(x)\bigg]m(x),
\end{equation}
and for $m(x)=\exp\big(\eta\int_{x_0}^x K/a\big)$, $0<\eta<1$,
\begin{equation}\label{eq:GF-Lyap-exp-1}
\LL^*m(x)=\bigg[2\exp\bigg(\!\!-\eta\int_{x/2}^x\frac{K}{a}\bigg)-(1-\eta)\bigg]K(x)m(x).
\end{equation}
Assumption~\eqref{eq:hypKmitose1} then ensures that $\LL^*m\sim -\xi Km$ as $x\to+\infty$, with $\xi=1-2^{1-r}>0$ in the first case and $\xi=1-\eta>0$ in the second case.
In both cases, we deduce that
\begin{equation}\label{eq:GF-Lyap-2}
\LL^*m\leq \kappa\, m+M\1_{(x_0,R)}m
\end{equation}
for any $\kappa\geq0$, by choosing $M>0$ and $R>x_0$ large enough, and this type of Lyapunov inequality is pivotal in our analysis.

\medskip

{\bf Condition \ref{H1}.}
Equation~\eqref{eq:mitosis} is a particular case of equation~\eqref{eq:Transport-lemRS3} with $G=\frakg=\RRR=0$, $b=K+\Div a$ and $\KKK[g](x)=4 K(2x)g(2x)$.
We may then use Proposition~\ref{prop:Transport-RKsolutions}-{\it(1)} to infer the well-posedness of equation~\eqref{eq:mitosis} in $X=L^1_m(x_0,\infty)$, provided that the conditions~\eqref{eq:TranspEqHypWellPoseWeight} and \eqref{eq:Transport-KKKbound} are met, with $0<\alpha_{_{\!\KKK}}<1$, which is nothing but~\eqref{eq:Transport-WPKKKRRRcondAlphaBeta} when $\RRR\equiv0$.
To do so, we define the function
\[\varpi := K - a \frac{m'}{m},\]
which corresponds to $\varpi_1$ in~\eqref{eq:transport-def-varpi}.
When $m(x)=x^r$ with $r>1$, we have
\[\varpi(x) = K(x) - r \frac{a(x)}{x},\]
and for $m(x)=\exp\big(\eta\int_{x_0}^x K/a\big)$ with $0<\eta<1$, we have
\[\varpi(x) = (1-\eta) K(x).\]
In both cases, the fact that $a\in \Lip  $  ensures that $\varpi_q:=\varpi+(1-1/q)a'$ enjoys $(\varpi_q)_-\in L^\infty$ for any $q\in[1,\infty)$.
On the other hand,~\eqref{as:Kb1} guarantees that $K\lesssim \langle\varpi_+\rangle$ and $a/x\lesssim \langle\varpi_+\rangle$, and finally~\eqref{eq:TranspEqHypWellPoseWeight} is verified.
The condition~\eqref{eq:Transport-KKKbound} is equivalent to the Lyapunov type condition
\begin{equation}\label{eq:GF-LyapK}
\KKK^*[m] \leq (\alpha_{_{\!\KKK}}\varpi_+ + M_{_{\!\KKK}})m,
\end{equation}
where $\KKK^*[m](x)=2K(x)m(x/2)$.
For $m(x)=x^r$ with $r>1$, we compute
\[\KKK^*[m]/m = 2^{1-r}K,\]
and for $m(x)=\exp\big(\eta\int_{x_0}^x K/a\big)$ with $0<\eta<1$,
\[\frac{\KKK^*[m](x)}{m(x)} = 2\exp\bigg(\!\!-\eta\int_{x/2}^x\frac{K}{a}\bigg)K(x).\]
Using~\eqref{as:Kb1}, we obtain that~\eqref{eq:GF-LyapK} is satisfied, for any $\alpha_{_{\!\KKK}}\in(2^{1-r},1)$ in the first case, and for any $\alpha_{_{\!\KKK}}\in(0,1)$ in the second case, by choosing $M_{_{\!\KKK}}$ large enough.

We can then apply Proposition~\ref{prop:Transport-RKsolutions}-{\it(1)} for associating to equation~\eqref{eq:mitosis} a strongly continuous semigroup $S$ in $X=L^1_m(x_0,\infty)$, and \ref{H1} then follows from Lemma~\ref{lem:Exist1-RkSG}-(i).
Moreover, we readily have that $\kappa_1\leq\kappa+M$ for any couple $(\kappa,M)$ such that~\eqref{eq:GF-Lyap-2} is verified.

\medskip

{\bf Condition \ref{H2}.}
We aim at verifying~\ref{H2} for some $\kappa_0>0$.
Recalling assumption~\eqref{as:Kinf}, we pick up $\ell>k$ and we consider the function
\[\phi_0(x)=xe^{-x^\ell/n}\]
with $n$ large enough to be chosen later.
We compute
\[\frac{\LL^*\phi_0(x)}{\phi_0(x)}=\frac{a(x)}{x}\Big(1-\frac{\ell}{n}x^{\ell}\Big)+K(x)\big(e^{\frac{1-2^{-\ell}}{n}x^\ell}-1\big).\]
Choosing $R>x_0$ such that $xK(x)/a(x)\geq\frac{2\ell}{1-2^{-\ell}}$ and $K(x)\geq e^{-x^k}$ for all $x\geq R$, we get that
\begin{align*}
\frac{\LL^*\phi_0(x)}{\phi_0(x)}&\geq \frac{a(x)}{x}+K(x)\Big(e^{\frac{1-2^{-\ell}}{n}x^\ell}-1-\frac{1-2^{-\ell}}{2n}x^{\ell}\Big)\\
&\geq e^{-x^k}\big(e^{\frac{1-2^{-\ell}}{n}x^\ell}-e^{\frac{1-2^{-\ell}}{2n}x^\ell}\big)\geq e^{\frac{1-2^{-\ell}}{2n}x^\ell-x^k}\big(e^{\frac{1-2^{-\ell}}{2n}R^\ell}-1\big)
\end{align*}
on $[R,\infty)$.
Choosing then $n\geq\frac{\ell}{2}R^\ell$, we have 
\[\frac{\LL^*\phi_0(x)}{\phi_0(x)}\geq\frac{a(x)}{2x}
\]
on $(x_0,R)$. Gathering the two above estimates, we deduce the existence of an explicit $\kappa_0>0$ such that $\LL^*\phi_0\geq\kappa_0\phi_0$.  
We conclude by invoking Lemma~\ref{lem:Existe1-Spectre2bis}-(i).

\medskip

{\bf Condition \ref{H3}.}
We consider the weight function $m(x)=x^r$ for some $r>1$ or $m(x)=\exp\big(\eta\int_{x_0}^x K/a\big)$ with $0<\eta<1$ and we define the stronger weight function $m_1(x)=\exp\big(\eta_1\int_{x_0}^x K/a\big)$ for some $\eta_1\in(\eta,1)$.
We fix $\kappa_\BB\in[0,\kappa_0)$, $M>0$, and $R>x_0$ such that~\eqref{eq:GF-Lyap-2} is verified by~$m_1$ with $\kappa=\kappa_\BB$.
Using the splitting $\LL=\AA+\BB$ with $\AA f=M\1_{(x_0,R)}f$, the inequality~\eqref{eq:GF-Lyap-2} for $m_1$ reads $\BB^*m_1\leq\kappa_\BB m_1$ and this ensures (see the proof of Corollary~\ref{cor:H3M1}) that $\kappa-\BB$ is invertible in $L^1_{m_1}$ for any $\kappa>\kappa_\BB$, with positive inverse, and
\[\|(\kappa-\BB)^{-1}\|_{\BBB(L^1_{m_1})}\leq\frac{1}{\kappa-\kappa_\BB}.\]
The operator $\AA$ maps $L^1_m$ into $L^1_{m_1}$ with
\[\|\AA\|_{\BBB(L^1_m,L^1_{m_1})}\leq M\frac{m_1(R)}{m(x_0)}.\]
Besides, due to the derivative part $\partial_x(a\,\cdot)$ in the operator $\BB$, we also have that $\RR_\BB(\kappa)$ maps $L^1_{m_1}$ into $ \Wloc^{1,1}$.
Finally, we have $\RR_\BB(\kappa)\AA:L^1_m\to L^1_{m_1}\cap \Wloc^{1,1}$, and thus $\RR_\BB(\kappa)\AA\in\KKK(L^1_m)$, for any $\kappa\geq\kappa_0>\kappa_\BB$.
We deduce from Lemma~\ref{lem:H3abstract-StrongC}-{\it(2)} that the condition \ref{H3} holds for both the primal and the dual problems.

\medskip

\begin{proof}[Proof of the existence part of Theorem~\ref{theo:GF-regular}]
We deduce from Theorem~\ref{theo:exist1-KRexistence} that the conclusion \ref{S1} about the existence of a solution $(\lambda_1,f_1,\phi_1)\in \R\times X_+\times X'_+$ to the first eigentriplet problem holds true.
\end{proof}
Moreover, we have $\lambda_1\geq\kappa_0>0$ and $f_1\in  \Wloc^{1,1}\cap L^1_m$ with $m(x)=\exp\big(\eta\int_{x_0}^x K/a\big)$ for any $\eta\in(0,1)$.


\medskip

{\bf  \Cyan Additional  estimates on the eigenfunction $\phi_1$.}
We consider the weight function $m(x)=x^r$ for some $r>1$ or $m(x)=\exp\big(\eta\int_{x_0}^x K/a\big)$ with $0<\eta<1$ and we define the weaker weight function $m_0(x)=x^{r_0}$ for some $r_0\in(1,r)$.
We fix $\kappa_\BB\in[0,\kappa_0)$, $M>0$, and $R>x_0$ such that~\eqref{eq:GF-Lyap-2} is verified by $m_0$.
Using again the splitting $\LL=\AA+\BB$ with $\AA f=M\1_{(x_0,R)}f$, \eqref{eq:GF-Lyap-2} means that $\BB^*m_0\leq\kappa_\BB m_0$ and this ensures that for any $\kappa>\kappa_\BB$ the operator $\kappa-\BB^*$ is invertible in $L^\infty_{m_0}$, with positive inverse, and
\[\|(\kappa-\BB^*)^{-1}\|_{\BBB(L^\infty_{m_0^{-1}})}\leq\frac{1}{\kappa-\kappa_\BB}.\]
Because of the derivative part of $\BB^*$, we also have that $\RR_{\BB^*}(\kappa):L^\infty_{m_0^{-1}}\to \Wloc^{1,\infty}$.
Besides, the operator $\AA^*=\AA$ maps $L^\infty_{m^{-1}}$ into $L^\infty_{m_0^{-1}}$ with
\[\|\AA^*\|_{\BBB(L^\infty_{m^{-1}}, L^\infty_{m_0^{-1}})}\leq M\frac{m(R)}{m_0(x_0)}.\]
Finally we have that $\RR_{\BB^*}(\kappa)\AA:L^\infty_{m^{-1}}\to L^\infty_{m_0^{-1}}\cap \Wloc^{1,\infty}$.
Consequently $\phi_1\in L^\infty_{m_0^{-1}}\cap \Wloc^{1,\infty}$ and
\begin{equation}\label{eq:GFphi1'estim}
\|\phi_1\|_{L^\infty_{m_0^{-1}}}=\|(\lambda_1-\BB^*)^{-1}\AA^*\phi_1\|_{L^\infty_{m_0^{-1}}}\leq \frac{m(R)}{m_0(x_0)}\frac{M}{\kappa_0-\kappa_\BB}\|\phi_1\|_{L^\infty_{m^{-1}}}.
\end{equation}
We also easily deduce quantitative estimates of $\phi_1$ in $ \Wloc^{1,\infty}$ from the identity
\[\phi_1'(x) = \frac{1}{a(x)}\big[\lambda_1\phi_1(x)+K(x)\phi_1(x)-2K(x)\phi_1(x/2)\big].\]

\medskip

{\bf Condition \ref{H4}.}
The operator $\LL$ satisfies the strong maximum principle.
Let $\lambda\in\R$ and $f\in X_+\cap D(\LL)\setminus\{0\}$ such that $(\lambda-\LL)f\geq0$, {\it i.e.}
\[\lambda f(x)+(af)'(x)+K(x)f(x)\geq4K(2x)f(2x)\qquad \forall x>x_0.\]
Denoting by $\Lambda_\lambda$ a function such that $\Lambda_\lambda'(x)=\frac{\lambda+K(x)}{a(x)}$, we get that
\begin{equation}\label{eq:SMP-GF}
a(x)f(x)\geq4\int_{x_0}^x e^{\Lambda_\lambda(y)-\Lambda_\lambda(x)}K(2y)f(2y)\,dy.
\end{equation}
Since $K(2y)>0$ for all $y>x_0$, $f\in X_+\setminus\{0\}$, and $a(x)>0$ for all $x>x_0$, we deduce from~\eqref{eq:SMP-GF} that the set $\{x>x_0,\ f(x)>0\}$ is an interval of the form $(\underline x,+\infty)$.
Using again~\eqref{eq:SMP-GF} we remark that we must have $\underline x=\max(x_0,\underline x/2)$, which enforces $\underline x=x_0$ and finally $f>0$.

\medskip

\begin{proof}[Proof of the uniqueness and positivity part of Theorem~\ref{theo:GF-regular}]
We deduce from Theorem~\ref{theo:KRgeometry1} the validity of the conclusion  \ref{S2} about  uniqueness and positivity of the solution $(\lambda_1,f_1,\phi_1)$ to the first eigentriplet problem.
\end{proof}
For deriving the exponential stability, we start by verifying a quantified irreducibility and aperiodicity condition on $S$, given in the next lemma, which then allows us to prove that the Doblin-Harris condition~\eqref{eq:hyp-Harris} is met.

 \begin{lem}\label{lem:GFHarris} Assume that~\eqref{as:mixingb} is satisfied.
 Then for all $\eps>0$, $R_1>x_0$, and $R_2>x_0+\eps$, there exists $T_1>0$ such that for any $T > T_1$, there exists $c_T>0$ such that 
$$
S_T^* \phi \ge c_T {\bf 1}_{(x_0,R_1)} \int_{x_0+\eps}^{R_2} \phi\, dx, \quad \forall \, \phi \ge 0.
$$
\end{lem} 

\begin{proof}[Proof of Lemma~\ref{lem:GFHarris}]
Throughout the proof we denote by $c_t$ any positive constant that depends only on $t$.
It is proved in~\cite[Prop.~5]{MR4897760} the existence of $(x_2,x_3)\subset(x_0,\infty)$ such that for all $R_1>x_0$ there exists $T_0>0$ such that for any $T > T_0$ and any $\phi\geq0$
\begin{equation}\label{eq:Harris-CW}
S_T^* \phi \ge c_T {\bf 1}_{(x_0,R_1)} \int_{x_2}^{x_3} \phi(x) dx.
\end{equation}
We may now extend the integral to $[x_0+\eps,R_2]$.
The Duhamel formula
\[S^*_{\LL} = S^*_{\BB_0} + S^*_{\BB_0}\AA_0 * S^*_{\LL}\]
for the splitting $\LL^*=\AA_0^*+\BB_0^*$ with $\AA_0^*\phi = \KKK^*[\phi]$ and $\BB_0^*\phi = b\phi'-K\phi$, also reads
\begin{equation}\label{eq:Duhamel-GF}
S^*_t\phi(x)  = \phi(X_t(x))e^{-\int_0^tK(X_s(x))ds} + 2 \!\int_0^t\! K(X_{t-s}(x)) S^*_s\phi(X_{t-s}(x)/2) e^{-\int_0^{t-s}K(X_{s'}(x))ds'} ds,
\end{equation}
where $X_t(x)$ is the solution to the characteristic equation
\beqn\label{eq:GF:dotX=aX}
\dot X_t(x)=a(X_t(x))\quad \text{with}\quad X_0(x)=x.
\eeqn
Applying~\eqref{eq:Harris-CW} to $S^*_t\phi$, that we bound from below by the first in Duhamel's formula~\eqref{eq:Duhamel-GF}, we obtain
\[S^*_{T+t}\phi \geq c_T\1_{(x_0,R_1)}\int_{x_2}^{x_3}\phi(X_t(x))e^{-\int_0^tK(X_s(x))ds}dx\geq c_T c_t \1_{(x_0,R_1)} \int_{X_t(x_2)}^{X_t(x_3)}\phi(y)dy.\]
Choosing $t_0$ such that $X_{t_0}(x_2)=x_3$, we get that for all $T>T_0+t_0$
\[S_T^*\phi\geq c_T\1_{(x_0,R_1)}\int_{x_2}^{X_{t_0}(x_3)}\phi(x)dx.\]
Iterating this argument and using the strict positivity of $a$ we get for any $R_2>x_2$ the existence of a time $t_1$ such that for all $T>T_0+t_1$
\begin{equation}\label{eq:Harris-2}
S_T^*\phi\geq c_T\1_{(x_0,R_1)}\int_{x_2}^{R_2}\phi(x)dx.
\end{equation}
For decreasing the lower bound of the integral from $x_2$ to $x_0+\eps$, we iterate once Duhamel's formula~\eqref{eq:Duhamel-GF} to get
$$
S^*_t\phi(x) \ge 2 \int_0^t K(X_{t-s}(x)) \phi(X_s(X_{t-s}(x)/2)) e^{-\int_0^{t-s}K(X_{s'}(x))ds'-\int_0^sK(X_{s'}(x))ds'} ds
$$
and then, using~\eqref{eq:Harris-2},
\[S^*_{T+t}\phi \ge c_tc_T \1_{(x_0,R_1)} \int_0^t\! \int_{x_2}^{R_2}\! K(X_{t-s}(x))\phi(X_s(X_{t-s}(x)/2)) dx\, ds.\]
We can assume that $x_2>2x_0$ and $R_2>2x_2$.
The fact that $x_2>2x_0$ ensures, due to assumption~\eqref{eq:hypKmitose1}, that $K$ is bounded from below by a positive constant on $[x_2,X_t(R_2)]$.
We thus deduce, by using of a change of variables, that for any $t>0$
\[S^*_{T+t}\phi \ge c_tc_T \1_{(x_0,R_1)} \int_{X_t(X_t(x_2)/2)}^{R_2/2} \phi(y) dy.\]
Since $X_t(x)\to x$ when $t\to0$, we deduce for all $\zeta>0$ the existence of $t>0$ such that
\[S^*_{T+t}\phi \ge c_t c_T \1_{(x_0,R_1)} \int_{x_2/2+\zeta}^{R_2/2} \phi(y) dy.\]
Since $R_2>2x_2$, we deduce by combining the above inequality with~\eqref{eq:Harris-2} that for all $T>T_0+t_1+t$
\[S^*_{T}\phi \ge c_T \1_{(x_0,R_1)} \int_{x_2/2+\zeta}^{R_2} \phi(x) dx.\]
Let us take $\zeta=x_0$.
Since the sequence $(u_n)$ defined by $u_0=x_2$ and $u_{n+1}=u_n/2+x_0$ converges to $2x_0$,
we obtain by an iteration argument the existence of a time $t_2$ such that for all $T>T_0+t_2$
\[S^*_{T}\phi \ge c_T \1_{(x_0,R_1)} \int_{2x_0+\eps}^{R_2} \phi(x) dx.\]
Using a last time the argument with $\zeta=\eps/2$ yields the desired result.
\end{proof}

We now prove another positivity result which allows making the time $T$ independent of $R_1$ in Lemma~\ref{lem:GFHarris}.

\begin{lem}\label{lem:GFpositivitypropagation}
Let $R_1>2x_0$. Then there exists $t_0>0$ such that for any $R>R_1$ we have
\[S_{t_0}^*\1_{(x_0,R_1)}\geq c_R\1_{(x_0,R)}\]
for some $c_R>0$.
\end{lem}

\begin{proof}[Proof of Lemma~\ref{lem:GFpositivitypropagation}]
Since $a$ is Lipschitz continuous, we can find $t_0 > 0$ small enough so that $X_{t_0}(x)\leq\alpha x$ for all $x>x_0$, with $\alpha>1$ to be determined later.
Then for any $t\in(0,t_0]$ and any $x\in(x_0,\frac{R_1}{\alpha})$,
 we have by using the first term in~\eqref{eq:Duhamel-GF}
\[
S^*_t\1_{(x_0,R_1)}(x)\geq c_{t_0}\1_{(x_0,R_1)}(X_t(x))=c_{t_0}>0.
\]
Iterating once~\eqref{eq:Duhamel-GF} and keeping only the second term, we get that for any $t\in(0,t_0]$ and any $x\in(2x_0,\frac{2R_1}{\alpha^2})$
\[
S^*_t\1_{(x_0,R_1)}(x)\geq c_{t_0}\int_0^t\1_{(x_0,R_1)}(X_s(X_{t-s}(x)/2))\,ds=c_{t_0}t.
\]
Choosing $\alpha>1$ such that $\frac{R_1}{\alpha}>2x_0$ and $\frac{2R_1}{\alpha^2}>R_1$,
we deduce that for any $t\in(0,t_0]$ there exists $c_t>0$ such that
\[
S^*_t\1_{(x_0,R_1)}\geq c_t\1_{(x_0,2\alpha^{-2}R_1)}.
\]
Dividing $[0,t_0]$ into $n$ sub-intervals $[\frac{k}{n}t_0,\frac{k+1}{n}t_0]$, $0\leq k\leq n-1$, and iterating the above inequality with $t=t_0/n$, we deduce for all integer $n\geq1$ the existence of $c_n>0$ such that
\[
S^*_{t_0}\1_{(x_0,R_1)}\geq c_n\1_{(x_0,(2\alpha^{-2})^nR_1)}
\]
and the proof is complete since $2\alpha^{-2}>1$.
\end{proof}

With Lemmas~\ref{lem:GFHarris} and~\ref{lem:GFpositivitypropagation}, we are now in position to prove the convergence result in Theorem~\ref{theo:GF-regular}.

\begin{proof}[Proof of the exponential stability part of Theorem~\ref{theo:GF-regular}]
We apply Theorem~\ref{theo:Harris}.
We start by proving that~\eqref{eq:hyp-Harris-BdBis} is verified, in a quantitative way, for the function $g_0=\1_{(x_0+\eps,R_2)}$ with a suitable choice of $R_2$ and $\eps$.
Choosing $r_0\in(1,r)$ if $m(x)=x^r$ or any $r_0>1$ is $m(x)=\exp\big(\eta\int_{x_0}^x K/a\big)$ and defining $m_0(x)=x^{r_0}$ we have from~\eqref{eq:GFphi1'estim}, because of the normalization $\|\phi_1\|_{L^{\infty}_{m^{-1}}}=1$,
\[\|\phi_1\|_{L^\infty_{m_0^{-1}}}\leq C_0\]
for some explicit constant $C_0>0$.
Defining
 $$
 R_2 := \inf \{ R > 0;  \ m_0(x)/m(x) \le 1/2C_0, \ \forall \, x > R \} , 
 $$
 we have
 $$
 1 = \| \phi_1 \|_{L^\infty_{m^{-1}}} = \sup_{(x_0,\infty)} \frac{\phi_1 }{ m} = \sup_{(x_0,R_2)} \frac{\phi_1}{ m},
 $$
 because 
 $$
 \sup_{(R_2,\infty)} \frac{\phi_1 }{ m} \le  \sup_{(R_2,\infty)} \frac{\phi_1}{ m_0} \frac{m_0 }{ m} \le C_0 \, \frac{1}{ 2C_0}   < 1. 
 $$
 Together with the fact that $\phi'_1\in\Lloc^\infty$, with a quantitative estimate on $\| \phi'_1 \|_{L^\infty(x_0,R_2)}$, we see that $\phi_1$ has some quantifiable mass on $(x_0+\eps,R_2)$ for $\eps>0$ small enough, which exactly means that $\langle g_0,\phi_1\rangle$ is quantified from below.
 
 \smallskip

Now we prove that the Doblin-Harris condition~\eqref{eq:hyp-Harris} is verified.
Choosing $R_1>2x_0$ and combining Lemma~\ref{lem:GFHarris} and Lemma~\ref{lem:GFpositivitypropagation}, we have for any $\eps>0$ and $R_2>x_0+\eps$ the existence of $T>0$ such that for any $R>R_1$
\begin{equation}\label{eq:GF-Harris-dual}
S_T^* \phi \ge c_R {\bf 1}_{(x_0,R)} \int_{x_0+\eps}^{R_2} \phi \, dx, \quad \forall \, \phi \ge 0.
\end{equation}
Defining $g_0=\1_{(x_0+\eps,R_2)}$, 
we deduce by duality that for all $f\geq0$,
\[S_T f \ge c_R \langle f,\1_{(x_0,R)}\rangle g_0.\]
Let us now consider  $A >0$ and $f \in X_+$ such that $\| f \| \le A [f]_{\phi_1}$.
Since $m_0(x)/m(x)\to0$ as $x\to+\infty$ and $\|\phi_1/m_0\|_\infty\leq C_0$, we have
\bean 
 [f]_{\phi_1} &=& \int_{x_0}^R f  \frac{\phi_1 }{ m} m + \int_R^\infty  f m  \frac{\phi_1 }{ m_0} \frac{m_0 }{ m}
 \\
 &\le& \langle f,\1_{(x_0,R)} \rangle  \sup_{(x_0,R)} m + \| f \| C_0 \sup_{(R,\infty)} \frac{m_0 }{ m}
 \\
 &\le& \langle f,\1_{(x_0,R)} \rangle \, m(R) + \frac12 [f]_{\phi_1},
 \eean
 by choosing $R$ large enough.
 We deduce that $S_Tf\geq \frac{c_R}{2m(R)} [f]_{\phi_1} g_0$,
 which is equivalent to~\eqref{eq:hyp-Harris} where we recall the definition $\widetilde S_t  := S_t e^{-\lambda_1 t}$.

 \smallskip

Finally, we prove the Lyapunov condition~\eqref{eq:stabilityKR-Lyapunov}.
On the one hand, we get from~\eqref{eq:GF-Lyap-2} that
\[\frac{d}{dt}\widetilde S_t ^*m = \widetilde S_t^*(\LL^*-\lambda_1)m\leq (\kappa_\BB-\lambda_1)\widetilde S_t ^* m+M\widetilde S_t ^*(\1_{(x_0,R)}m).\]
On the other hand, arguing as in~\eqref{eq:stabilityKR-phi1bybelow}, we infer from~\eqref{eq:GF-Harris-dual} that
\[\phi_1 = e^{-\lambda_1 T}S^*_T\phi_1 \geq c_R e^{-\lambda_1 T}\langle g_0,\phi_1\rangle \1_{(x_0,R)}.\]
Combining both we deduce that
\[\frac{d}{dt}\widetilde S_t ^*m\leq (\kappa_\BB-\lambda_1)\widetilde S_t ^* m+\tilde M\phi_1\]
with $\tilde M = \frac{m(R)}{c_R}\frac{e^{\lambda_1T}}{\langle g_0,\phi_1\rangle}$, and Grönwall's inequality then yields
\[\widetilde S_t ^*m\leq e^{(\kappa_\BB-\lambda_1)t}m+\tilde M t e^{(\kappa_\BB-\lambda_1) t}\phi_1.\]
This guarantees that~\eqref{eq:stabilityKR-Lyapunov} is verified with $\gamma_L=e^{(\kappa_\BB-\lambda_1)T}\in(0,1)$ and $K=\tilde M T$.

\smallskip

We have proved that the conditions~\eqref{eq:hyp-Harris-mean}, \eqref{eq:stabilityKR-Lyapunov} and \eqref{eq:hyp-Harris-BdBis} are verified.
The conclusion of the proof then follows from Theorem~\ref{theo:Harris}.
\end{proof}

\subsection{The mitosis equation with non-mixing growth rate}\label{ssec:GF:singular}

In this section, we investigate the case when the mixing condition~\eqref{as:mixingb} is not verified.
In other words, we place ourselves under the singular condition  
\begin{equation}\label{as:nomixingb}
\forall x > x_0,\quad a(2x)=2a(x).
\end{equation}
In this case, we still have the existence of a unique eigentriplet $(\lambda_1,f_1,\phi_1)$ but the boundary point  spectrum is not reduced to $\lambda_1$.
As a consequence, the long time behavior of the semigroup does not stabilizes along $f_1$ but it exhibits periodic oscillations.

\begin{theo}\label{theo:GF-singular}
Suppose that~\eqref{eq:hypKmitose1}, \eqref{as:Kb1}, \eqref{as:Kinf} and \eqref{as:nomixingb} are satisfied.
{\Cyan The conclusions \ref{S1},  \ref{S2} and  \ref{S31} with $\lambda_1 > 0$  hold  in $X := L^1_m$.}
Besides, $\Sigma_P^+(\LL)=\{\lambda_1+ik\alpha,\,k\in\Z\}$ for some quantifiable $\alpha>0$, there exists a family $(g_k,\psi_k)_{k\in\Z}$ of corresponding primal and dual eigenvectors that verifies $\langle \psi_k,g_\ell \rangle = \delta_{k\ell}$, and for all $f\in L^1_{\phi_1}$, we have the convergence
\[ {\| e^{-\lambda_1 t} S_\LL(t) (f - \Pi f) \|}_{L^1_{\phi_1}} \to 0 \qquad\text{as}\ t\to+\infty,\]
where $\displaystyle \Pi f = \lim_{n\to\infty}\frac1n\sum_{\ell=0}^n\sum_{k=-\ell}^\ell \langle \psi_k , f \rangle g_k$.
\end{theo}

This new result complements the scarce literature on the long time behavior of equation~\eqref{eq:mitosis} in the singular case~\eqref{as:nomixingb} which, to the best of our knowledge, is limited to the references~\cite{MR3928121,GabrielMartin,Greiner1988}. 
We will actually prove that the convergence in Theorem~\ref{theo:GF-singular} also holds in other spaces,
and this will be the occasion to compare our method and results to the three above mentioned papers.

\medskip

The proof of the conclusion~\ref{C2} in Section~\ref{ssec:GF:regular} does not use the mixing   assumption~\eqref{as:mixingb}.
It thus also proves the existence, uniqueness and strict positivity of eigentriplet $(\lambda_1,f_1,\phi_1)$ under the assumptions of Theorem~\ref{theo:GF-singular}, as well as the fact that equation~\eqref{eq:mitosis} is associated with a semigroup $S$ in $X$.
For proving the long time convergence result, we start by verifying that this semigroup extends to other relevant Banach spaces.

\medskip

{\bf Well-posedness in entropic $L^p$ and $M^1$ spaces.}
The dual eigenfunction $\phi_1$ satisfies by definition $\LL^*\phi_1=\lambda_1\phi_1$ and the rescaled semigroup $\widetilde S_t =S_te^{-\lambda_1 t}$ is thus a contraction for the norm of $L^1_{\phi_1}$.
In particular $S_t$ is a bounded operator for this norm and, since $L^1_m$ is dense in $L^1_{\phi_1}$, we can uniquely extend the semigroup $S$ into a strongly continuous semigroup in $L^1_{\phi_1}$.
Similarly, due to the weak-$*$ density of $L^1_{\phi_1}$ into $M^1_{\phi_1}$, this semigroup extends uniquely into a weakly-$*$ continuous semigroup in $M^1_{\phi_1}$.
We still denote by $S$ these extensions.

The General Relative Entropy principle, see~\cite{MR2162224,Bernard2022}, ensures that the weighted $L^p$ sub-spaces of $L^1_{\phi_1}$ defined by
\[X_p:=L^p_{f_1^{1-p}\phi_1}(x_0,\infty)\quad\text{for}\ p\in[1,\infty)\quad\text{and}\quad X_\infty:=L^\infty_{f_1^{-1}}(x_0,\infty)\]
are invariant under the semigroup $S$ and the restriction to these spaces is a contraction.
Besides, Jensen's inequality yields that it is a decreasing sequence for the inclusion
\[p>q\quad\implies\quad X_p\supset X_q.\]
Since $X_\infty\subset X_p$ is dense, we can infer the strong continuity of $S$ in $X_p$ from the strong continuity in $X_1$ by writing for any $f\in X_\infty$
\[
\|\widetilde S_t f-f\|^p_{X_p}\leq\|\widetilde S_t f-f\|_{X_\infty}^{p-1}\|\widetilde S_t f-f\|_{X_1}\leq2^{p-1}\|f\|_{X_\infty}^{p-1}\|\widetilde S_t f-f\|_{X_1} \to 0, 
\]
as $t \to 0$. 
\medskip

{\bf Long-time convergence in $M^1_{\phi_1}$.}
We start by giving some useful properties of the dual semigroup $S^*$ in $X'=L^{\infty}_{m^{-1}}$.
Splitting $\LL^*$ as $\LL^*=\AA_0^*+\BB_0^*$ with $\AA_0^*\phi = \KKK^*[\phi]$, so that $\BB_0^*\phi = a\phi'-K\phi$,
the Duhamel formula
\[S^*_{\LL} = S^*_{\BB_0} + S^*_{\BB_0}\AA_0 * S^*_{\LL}\]
ensures that $\bar\varphi(t,x):=S^*_t\phi(x)$ is a fixed point of the operator $\Gamma$ defined by
\begin{align}\label{eq:Gamma}
\Gamma\varphi(t&,x):= S^*_{\BB_0}(t)\phi(x) + \big[S^*_{\BB_0}\AA_0 * \varphi(\cdot,x)\big] (t)\\ \nonumber
&=\phi(X_t(x))e^{-\int_0^tK(X_s(x))ds} + 2\!\int_0^t\! K(X_{t-s}(x))\varphi(s,X_{t-s}(x)/2) e^{-\int_0^{t-s}K(X_{s'}(x))ds'} ds,
\end{align}
where we recall that $X_t(x)$ is the solution to the characteristic equation \eqref{eq:GF:dotX=aX}. 
It turns out that $\Gamma$ has a unique fixed point in $\Lloc^\infty([0,\infty)\times(x_0,\infty))$, and that this fixed point also lies in any closed subset of $\Lloc^\infty([0,\infty)\times(x_0,\infty))$ which is left invariant by $\Gamma$.
This property is proved in~\cite{GabrielMartin} or in~\cite[Sec.~6.3]{Bansaye2022}, by building $\bar\varphi$ thanks to the Banach-Picard fixed point theorem.
It has very useful consequences, as for instance the fact that if $\phi\in C(x_0,\infty)$, then $\bar\varphi\in C([0,\infty)\times(x_0,\infty))$.
In particular, this implies that $C(x_0,\infty)\cap L^\infty_{m^{-1}}$ is invariant under the semigroup $S^*$.
Since $C(x_0,\infty)\cap L^\infty_{m^{-1}}$ is a dense subspace of $C_{0,\phi_1}(x_0,\infty)$, this ensures that $C_{0,\phi_1}$ is invariant under $S^*$ and that the duality relation
\[\langle S_tf,\phi\rangle = \langle f,S_t^*\phi\rangle\]
is valid for any $f\in M^1_{\phi_1}$ and $\phi\in C_{0,\phi_1}$.
The proof of the next result crucially relies on another application of the fact that the fixed point of $\Gamma$ belongs to any closed invariant subset.

\begin{prop}\label{prop:GF-weakperiodic}
Suppose that~\eqref{eq:hypKmitose1}, \eqref{as:Kb1}, \eqref{as:Kinf} and~\eqref{as:nomixingb} are satisfied.
Then $\Sigma_P^+(\LL)=\{\lambda_1+ik\alpha,\,k\in\Z\}$ for some $\alpha>0$, there exists a family $(g_k,\psi_k)_{k\in\Z}$ of corresponding primal and dual eigenvectors that verifies $\langle \psi_k,g_\ell \rangle = \delta_{k\ell}$, and for all $f\in M^1_{\phi_1}$ we have the convergence
\begin{equation}\label{eq:GF-periodic}
\widetilde S_t f-\widetilde S_t  \Pi f \to 0, 
\ \hbox{ as } \ t \to \infty, \quad 
 \Pi f := \lim_{n\to\infty}\frac1n\sum_{\ell=0}^n\sum_{k=-\ell}^\ell \langle \psi_k , f \rangle g_k,
\end{equation}
both convergences having to be understood in the sense of  the weak-$*$ topology. 
\end{prop}

Note that we did not specify the space in which we define the boundary point spectrum $\Sigma_P^+(\LL)$ in Proposition~\ref{prop:GF-weakperiodic}.
It is because this set is the same in all the Banach lattices we consider.
Indeed, any $g\in M^1_{\phi_1}$ such that $\LL g=\lambda g$ for some $\lambda\in\C$ with $\Re e(\lambda)=\lambda_1$ satisfies $|g|\in\Span(f_1)$, so that $g\in X=L^1_m$ for any weight $m$ as in~\eqref{def:GF-m}.

\begin{proof}[Proof of Proposition~\ref{prop:GF-weakperiodic}]
{\sl Step 1.} The rescaled semigroup $\widetilde S$ is a contraction semigroup in $M^1_{\phi_1} = (C_{0,\phi_1})'$.
This ensures in particular that for all $f\in M^1_{\phi_1}$ the trajectory $(\widetilde S_t f)_{t\geq0}$ is bounded in $M^1_{\phi_1}$. 
We can thus use Theorem~\ref{theo:ergodicity-compact-trajectories}-{\it(2)} to infer the non-triviality of the boundary spectrum, by proving that the conclusion cannot hold, see Remark~\ref{rem:periodic}-{\it(ii)}.
We start from the fact that for any $\phi\in C_{0,\phi_1}(x_0,\infty)$, the solution $S^*\phi$ to the dual mitosis equation is the unique fixed point of $\Gamma$ defined in~\eqref{eq:Gamma}, and that it belongs to any closed invariant subset of $C([0,\infty)\times(x_0,\infty))$.
For $y>x_0$ we define the set
\[
\mathcal E_y=\{x>x_0,\ x=2^ky\ \text{for some}\ k\in\Z\}
\]
and we consider a function $\phi$ such that $\phi(x)=0$ if $x\in\mathcal E_y$ and $\phi(x)>0$ if $x\not\in\mathcal E_y$.
Then, since~\eqref{as:nomixingb} ensures that $X_t(2x)=2X_t(x)$ for all $t\geq0$, the set
\[\big\{\varphi\in C([0,\infty)\times(x_0,\infty)),\ \varphi(t,x)=0\ \text{if}\ X_t(x)\in\mathcal E_y\ \text{and}\ \varphi(t,x)>0\ \text{if}\ X_t(x)\not\in\mathcal E_y\big\}\]
is invariant under $\Gamma$.
Consequently, the unique fixed point $S^*_t\phi$ belongs to this set, and we deduce that $S^*_t\phi(x)=0$ if and only if $X_t(x)\in\mathcal E_y$.
In other words, by duality, $\supp(S_t\delta_x)\subset\mathcal E_y$ if and only if $X_t(x)\in\mathcal E_y$,
and in particular $\supp(S_t\delta_x)\subset\mathcal E_{X_t(x)}$ for all $x>x_0$ and  $t\geq0$.
This prevents the convergence of $\widetilde S_t \delta_x$ toward $\langle \delta_x,\phi_1\rangle f_1$ and we infer from (the negation of) Theorem~\ref{theo:ergodicity-compact-trajectories}-{\it(2)} that the boundary point spectrum cannot be trivial.%

\smallskip
{\sl Step 2.} 
We next formulate the following simple but fundamental observation. If $(\lambda,f)$ is a solution to the eigenvalue problem in $X = M^1_m$, then $f \in D(\LL) \subset BV_{\rm loc} \subset L^1_{\rm loc}$ and it
 is a solution to the eigenvalue problem in $X = L^1_m$. 
 Symmetrically, if $(\lambda,\phi)$ is a solution to the eigenvalue problem in $Y = L^\infty_{m_1^{-1}}$, then $ \phi \in D(\LL^*) \subset \hbox{\rm Lip}_{\rm loc} \cap  L^\infty_{m_2^{-1}}$, $m_1(x)/m_2 (x)\to 0$ as $x \to\infty$, and it is a solution to the eigenvalue problem in $Y = C_{0,m^{-1}_1}$. In other words, the point spectrum and the associated eigenelements are the same in the two frameworks $(L^1_m,L^\infty_{m^{-1}})$ and $(M^1_m,C_{0,m^{-1}})$. 
Now, as a consequence of this observation and Step 1, we know that the  boundary point spectrum $\Sigma^+_P(\LL)$ is not trivial. Because we have proved that $(\kappa-\BB)^{-1}\AA$ is compact in $L^1_m$, $\kappa < \lambda_1$, we may apply Theorem~\ref{theo:Sigma+subgroupBIS} and we obtain that $\Sigma^+_P(\LL) = \{\lambda_1\}+i\alpha\Z$ for some $\alpha>0$,
and each eigenvalue is algebraically simple.
Using finally Theorem~\ref{theo:periodic} in the situation~{\it(2)}, we get the weak-$*$ convergence~\eqref{eq:GF-periodic}.
\end{proof}

This result is proved by means of entropy techniques in~\cite{GabrielMartin} for a linear growth rate $a(x)=x$, by taking advantage of the explicit formulation of the eigenvectors $g_k$ and $\psi_k$ in terms of $f_1$ and $\phi_1$ in that case.
Here we extend it to any $a$ satisfying $a(2x)=2a(x)$.
Note that arguing similarly as in~\cite{GabrielMartin}, the convergence~\eqref{eq:GF-periodic} may be strengthened into an exponential strong convergence in $M^1_m$ for $m(x)=x^r$, $r>1$, or $m(x)=\exp\big(\eta\int_{x_0}^x K/a\big)$, $0<\eta<1$, meaning that there is a spectral gap between $\Sigma_+(\LL)$ and the rest of the spectrum in these spaces.

\medskip

{\bf Long-time convergence in $X_p$.}
We prove the following result, the case $p=1$ of which corresponds to the convergence result of Theorem~\ref{theo:GF-singular}.

\begin{prop}\label{prop:GF-periodic-entropic}
Under the same assumptions as in Proposition~\ref{prop:GF-weakperiodic}, the convergence~\eqref{eq:GF-periodic} holds for the strong topology in $X_p$, $1\leq p<\infty$ for all $f\in X_p$, and the convergence of the Fejér sum in the definition of the projector $\Pi$ is also for this topology.
\end{prop}

\begin{proof}[Proof of Proposition~\ref{prop:GF-periodic-entropic}]
The case $p=1$ is an immediate consequence of Theorem~\ref{theo:periodic}, case~{\it(4)}.
The proof in the case $p>1$ is a direct adaptation of the case $p=1$.
We aim at verifying that the trajectories $(\widetilde S_t f)_{t\geq0}$ are relatively compact in $X_p$.
We have already seen that $X_\infty\subset X_p$ is dense.
Besides, the domain $D(\LL)$ of the generator $\LL-\lambda_1$ of $\widetilde S$ in $X_p$ is also dense in $X_p$, so that it suffices to check the relative compactness of $(\widetilde S_t f)_{t\geq0}$ for $f\in X_\infty\cap D(\LL)$.
For $f$ in $X_\infty\cap D(\LL)$ the bounds
\[\|\widetilde S_t f\|_{X_p}\leq\|f\|_{X_p},\quad\|\LL\widetilde S_t f\|_{X_p}=\|\widetilde S_t \LL f\|_{X_p}\leq\|\LL f\|_{X_p}\quad\text{and}\quad \|\widetilde S_t f\|_{X_\infty}\leq \|f\|_{X_\infty}\]
yield the relative compactness of $(\widetilde S_t f)_{t\geq0}$, the second bound guaranteeing uniform $ \Wloc^{1,1}$ estimates.
We can thus apply the case~{\it(1)} of Theorem~\ref{theo:periodic} to deduce the convergence~\eqref{eq:GF-periodic} in $X_p$ for the strong topology.
\end{proof}

Proposition~\ref{prop:GF-periodic-entropic} extends the result of~\cite{MR3928121} where it is proved in the case $p=2$ for $a(x)=x$ by taking advantage of the Hilbert structure of $X_2$ and of the explicit formulation of the eigenvectors $g_k$ and $\psi_k$ in terms of $f_1$ and $\phi_1$.
In this Hilbert setting it is proved that the Fourier series $\sum_{k=-n}^n\langle f,\psi_k\rangle g_k$ converges as $n$ goes to infinity, and $\Pi f$ is then given by the limit.

\medskip

{\bf About the value of $\alpha$.}
For ensuring that the boundary spectrum is discrete, we have used a compactness argument.
The period $2\pi/\alpha$ of the periodic semigroup $\widetilde S\Pi$ is thus not quantified.
It is expected to be equal to the time needed for a particle to double its size by following the flow of~$a$, namely
\begin{equation}\label{eq:period}
\frac{2\pi}{\alpha}=\int_{x}^{2x}\frac{dt}{a(t)},
\end{equation}
which is independent of the choice of $x\geq x_0$ due to the condition $a(2t)=2a(t)$.
This is known to be true in the case of a linear growth rate $a(x)=x$, see~\cite{Diekmann1984} or~\cite{MR3928121}, and also for $a$ general when the size domain is $(x_0,4x_0)$, see~\cite{Greiner1988}, where explicit computations can be carried out.
In the general case, we have not been able to prove~\eqref{eq:period}.
Yet the fact that for any $x>x_0$ the support of $S_t\delta_x$ is a subset of $\mathcal E_{X_t(x)}$ guarantees that the period cannot be too small as shown now. \Black 

\begin{prop}\label{prop:GFupperboundalpha}
We have the estimate
\beqn\label{eq:GF-upperboundalpha}
\frac{2\pi}{\alpha}\geq \ell_a := \int_{x}^{2x}\frac{dt}{a(t)}.
\eeqn
\end{prop}
\begin{proof}[Proof of Proposition~\ref{prop:GFupperboundalpha}]
Let $x>x_0$ such that $\psi_1(x)\neq0$ (actually any $x>x_0$ is suitable).
We have
\[
\widetilde S_t \delta_x-\widetilde S_t  \Pi \delta_x \to 0, 
\]
as $t\to+\infty$, 
and $\supp\widetilde S_t \delta_x\subset\mathcal E_{X_t(x)}$, so $\supp\widetilde S_t \Pi\delta_x\subset\mathcal E_{X_t(x)}$ since $\widetilde S_t \Pi\delta_x$ is periodic.
Besides,
\[\widetilde S_t \Pi\delta_x=\lim_{n\to\infty}\frac{1}{n}\sum_{\ell=0}^n\sum_{k=-\ell}^\ell\psi_k(x)e^{i\alpha kt} g_k\]
and, since $\psi_1(x)\neq0$, the period of this periodic function of time is $\frac{2\pi}{\alpha}$.
But since $\supp\widetilde S_t \Pi\delta_x\subset\mathcal E_{X_t(x)}$ and the period of the set $\mathcal E_{X_t(x)}$ is $\ell_a$, 
we have proved that \eqref{eq:GF-upperboundalpha} holds. 
\end{proof}

On the other hand, we can use Theorem~\ref{theo:Harris-mean} for deriving a quantified lower bound on $\alpha$, and thus an upper bound on the period.
We work in the space $X=L^1_m$, recalling that $\Sigma_P^+(\LL)$ is the same in this space and in $M^1_{\phi_1}$.

\begin{prop}\label{prop:GF-isolation}
There exists a constructive constant $\alpha_1>0$ such that $\Sigma(\LL)\cap B(\lambda_1,\alpha_1)=\{\lambda_1\}$.
In particular, $\alpha\geq\alpha_1$.
\end{prop}

To prove this result, we check that the conditions~\eqref{eq:stabilityKR-Lyapunov}, \eqref{eq:hyp-Harris-BdBis} and~\eqref{eq:hyp-Harris-mean} are verified, and we invoke Theorem~\ref{theo:Harris-mean}.
We start with a lemma which, together with Lemma~\ref{lem:GFpositivitypropagation}, will guarantee the validity of~\eqref{eq:hyp-Harris-mean}.

\begin{lem}\label{lem:GFmeanHarris}
 For all $\eps>0$, $R_1>x_0$, and $R_2>x_0+\eps$, there exist $T > 0$ and $c_T>0$ such that 
$$
\int_0^TS_t^* \phi \,dt\ge c_T {\bf 1}_{(x_0,R_1)} \int_{x_0+\eps}^{R_2} \phi \, dx, \quad \forall \, \phi \ge 0.
$$
\end{lem}

\begin{proof}[Proof of Lemma~\ref{lem:GFmeanHarris}]
Throughout the proof we denote by $c_t$ any positive constant that depends only on $t$.
From the Duhamel formula~\eqref{eq:Duhamel-GF}, we get by positivity that for any $\phi\geq0$
\[\int_0^TS^*_t\phi(x) dt\ge \int_0^T\phi(X_t(x))e^{-\int_0^tK(X_s(x))ds}dt.\]
We deduce that for all $x\in(x_0,R_1)$ and for $T_0$ large enough so that $X_{T_0}(x_0)>R_2$, we have
\[
\int_0^{T_0}S^*_t\phi(x) dt\ge c_{T_0} \int_{x}^{X_T(x)}\phi(y)dy\geq c_{T_0} \int_{R_1}^{R_2}\phi(y)dy,
\]
and the conclusion follows if $R_1\leq x_0+\eps$.
If not, we have with the same argument the existence of $T$ such that for all $x\in(x_0,R_1)$
\[
\int_0^{T}S^*_t\phi(x) dt\geq c_{T} \int_{\max(R_1,2x_0+\eps)}^{2R_2}\phi(y)dy.
\]
Iterating once Duhamel's formula and using that $X_t(x/2)=X_t(x)/2$ and~\eqref{eq:hypKmitose1} we get that for all $t\geq0$ and all $y\in(2x_0+\eps,R_2)$
\[
S^*_t\phi(y)\ge c_t \int_0^t K(X_{t-s}(y)) \phi(X_s(X_{t-s}(y)/2)) ds\ge c_t \phi(X_t(y/2)).
\]
This yields for all $x\in(x_0,R_1)$
\begin{align*}
\int_0^{T+t}S^*_s\phi(x)ds&\geq \int_t^{T+t}S^*_s\phi(x)ds=\int_0^{T}S^*_sS^*_t\phi(x)ds\\
&\geq c_{T}\int_{\max(R_1,2x_0+\eps)}^{2R_2}S^*_t\phi(y)dy\\
&\geq c_{T}c_t\int_{\max(R_1,2x_0+\eps)}^{2R_2}\phi(X_t(y/2))dy=c_Tc_t\int_{\max(X_t(\frac{R_1}{2}),X_t(x_0+\frac{\eps}{2}))}^{X_t(R_2)}\phi(z)dz.
\end{align*}
Choosing $t > 0$ small enough so that $\max(X_t(\frac{R_1}{2}),X_t(x_0+\frac{\eps}{2}))\leq \max(\frac{R_1+\eps}{2},x_0+\eps)$,  we get 
\[
\int_0^{T_1}S^*_s\phi(x)ds\geq c_{T_1}\int_{ \max(\frac{R_1+\eps}{2},x_0+\eps)}^{R_2}\phi(z)dz, \ \hbox{ for } \ T_1:=T+t.
\]
Iterating the argument we can build an increasing sequence of times $T_n$ such that 
\[
\int_0^{T_n}S^*_t\phi(x)dt\geq c_{T_n}\int_{ \max(u_n,x_0+\eps)}^{R_2}\phi(z)dz, \ \forall \, n \geq0, 
\]
where $(u_n)$ is defined by $u_0=R_1$ and $u_{n+1}=\frac{u_n+\eps}{2}$.
Since this sequence $(u_n)$ converges to $\eps<x_0+\eps$, we get the conclusion by taking $n$ large enough.
\end{proof}

We are now in position to prove Proposition~\ref{prop:GF-isolation}.

\begin{proof}[Proof of Proposition~\ref{prop:GF-isolation}]
Arguing similarly as in the proof of Theorem~\ref{theo:GF-regular} and using Lemma~\ref{lem:GFmeanHarris} instead of Lemma~\ref{lem:GFHarris}, we can prove that the conditions~\eqref{eq:stabilityKR-Lyapunov}, \eqref{eq:hyp-Harris-BdBis} and~\eqref{eq:hyp-Harris-mean} are verified.
Applying Theorem~\ref{theo:Harris-mean} then gives the result.
\end{proof}

\subsection{The model with variability}\label{ssec:GF:variability}

In this last part, we consider the model with variability given by the equation~\eqref{eq:mitosis-var}.
Compared to equation~\eqref{eq:mitosis}, the main consequence of introducing a variability in terms of the spectrum and the asymptotic behavior is that for equation~\eqref{eq:mitosis-var} the boundary spectrum is trivial and the first eigenfunction is exponentially stable, no matter if $a$ satisfies~\eqref{as:mixingb} or~\eqref{as:nomixingb}.

\medskip

Additionally to the assumptions~\eqref{eq:hypKmitose1}, \eqref{as:Kb1} and \eqref{as:Kinf}, we ask that
\begin{equation}\label{eq:GFV-Kinfty}
K(x)=O\Big(\exp\Big(\delta\int_{x_0}^{x/2}K/a\Big)\Big)\quad\text{as}\ x\to+\infty,
\end{equation}
for some $\delta>0$.
About the variability kernel $\wp$ we suppose that
\begin{equation}\label{eq:wp}
\int_1^2\wp(v,v_*)dv=1,\ \forall v_*\in[1,2],\quad\wp\in W^{1,\infty}([1,2]^2)\quad\text{and}\quad \wp\geq\wp_*
\end{equation}
for some $\wp_*>0$.
{\Red We  work in the space $X:=L^1_m(\OO)$ with $m=m(x)$ still given by \eqref{def:GF-m} and $\OO := (x_0,\infty) \times (1,2)$.}

\begin{theo}\label{theo:GFV-main}
Suppose that~\eqref{eq:hypKmitose1}, \eqref{as:Kb1}, \eqref{as:Kinf}, \eqref{eq:GFV-Kinfty} and~\eqref{eq:wp} are satisfied.
{\Cyan Then  the conclusions \ref{S1},  \ref{S2},  \ref{S33} and  \ref{E31} hold with constructive constants in $L^1_m$.}
\end{theo}

Yet expected, this result was known only in the case of a discrete set of variabilities~\cite{CloezdeSaporta,RatTournus}.
Theorem~\ref{theo:GFV-main} is thus new in the literature.

\medskip

Because of the assumption~\eqref{eq:wp}, we easily see that the construction of the semigroup and the proof of the conditions \ref{H1}, \ref{H2} and \ref{H4} given in Section~\ref{ssec:GF:regular} for the model without variability readily extend to the model~\eqref{eq:mitosis-var}.
We thus only have to verify~\ref{H3} and some Doblin-Harris type condition.

\medskip

{\bf Condition \ref{H3}.}
Let $\delta\in(0,1)$ such that~\eqref{eq:GFV-Kinfty} is verified, and consider the weight function $m(x)=x^r$ with $r>1$ or $m(x)=\exp\big(\eta\int_{x_0}^x K/a\big)$ with $\eta\in(0,1-\delta)$.
We also use the two other weights
$$
m_1(x)=\exp\big(\eta_1\int_{x_0}^x K/a\big), \quad 
m_2(x)=\exp\big(\eta_2\int_{x_0}^x K/a\big)
$$ for some $\eta_1\in(\eta,1-\delta)$ and $\eta_2=\eta_1+\delta$.
We combine the two different splittings $\LL=\AA+\BB$ and $\LL=\AA_0+\BB_0$, where 
$$
\AA f (x,v) = M\1_{(x_0,R)}(x)f(x,v), \quad \AA_0 f (x,v) = 4 \int_1^2 K(2x) \wp(v,v_*) f(x,v_*)\, dv_*.
$$
We prove that for any $\kappa>\kappa_\BB$ the operator
\[\CC := (\kappa-\BB_0)^{-1}\AA_0(\kappa-\BB)^{-1}\AA\]

is well defined and maps continuously $L^1_m$ into $L^1_{m_1}\cap  \Wloc^{1,1}$, in the sense that if $(f_n)$ is bounded in $L^1_m$ then the image is bounded in $L^1_{m_1}\cap W^{1,1}((x_0,R)\times[1,2])$ for all $R>0$. 
In particular, $\CC \in \KK(L^1_m)$. 
More precisely, for any $\kappa >0$, we prove 
\[L^1_m\xrightarrow[]{\quad\AA\quad}L^1_{m_2}\xrightarrow[]{(\kappa-\BB)^{-1}}L^1_{m_2}\xrightarrow[]{\quad\AA_0\quad}L^1_{m_1} \cap W^{1,1}_{v,\rm loc}\xrightarrow[]{(\kappa-\BB_0)^{-1}}L^1_{m_1}\cap  \Wloc^{1,1}\]
where $W^{1,1}_{v,loc}:=\{f\in \Lloc^1((x_0,\infty)\times[1,2]),\ \partial_v f \in \Lloc^1((x_0,\infty)\times[1,2]) \}$.

\smallskip

The results for $\AA$ and $(\kappa-\BB)^{-1}$ are proved as in the case without variability.
For the third one, the fact that $\AA_0$ maps $L^1_{m_2}$ in $L^1_{m_1}$ follows from assumption~\eqref{eq:GFV-Kinfty}, and the fact that the range is in $W^{1,1}_{v,loc}$ is a direct consequence of the assumption that $\wp \in W^{1,\infty}([1,2]^2)$.

Finally we consider $\kappa-\BB_0$ and we first verify that it is invertible in $L^1_{m_1}$ for any $\kappa>0$.
If $(\kappa-\BB_0)g=f$, then necessarily
\begin{equation}\label{eq:GF-invB0}
g(x,v)=\frac{1}{va(x)}\int_{x_0}^x e^{(\Lambda_\kappa(y)-\Lambda_\kappa(x))/v}f(y,v)\,dy,
\end{equation}
where $\Lambda_\kappa(x)=\int_{x_0}^x\frac{\kappa+K}{a}$, and consequently
\[g(x,v)m_1(x)=\frac{e^{(\eta_1\Lambda(x)-\Lambda_\kappa(x))/v}}{va(x)}\int_{x_0}^x e^{(\Lambda_\kappa(y)-\eta_1\Lambda_0(y))/v}f(y,v)m_1(y)\,dy.\]
Since
\[\Lambda_\kappa(x)-\eta_1\Lambda_0(x)=(1-\eta_1)\Lambda_{\kappa/(1-\eta_1)}(x)\]
we have for all $v\in[1,2]$
\begin{align*}
\int_{x_0}^\infty \Big( \frac{\kappa}{1-\eta_1} & + K(x) \Big) g(x,v) m_1(x) \,dx \\
& = \int_{x_0}^\infty \frac1v \Lambda'_{\kappa/(1-\eta_1)}(x) e^{-\frac{1-\eta_1}{v}\Lambda_{\kappa/(1-\eta_1)}(x)}\int_{x_0}^x e^{(1-\eta_1)\Lambda_{\kappa/(1-\eta_1)}(y,v)} f(y,v) m_1(y)\,dy dx\\
&=\int_{x_0}^\infty e^{\frac{1-\eta_1}{v}\Lambda_{\kappa/(1-\eta_1)}(y)}f(y,v)m_1(y)\int_y^\infty \frac 1v \Lambda'_{\kappa/(1-\eta_1)}(x) e^{-\frac{1-\eta_1}{v}\Lambda_{\kappa/(1-\eta_1)}(x)}dx\,dy\\
&=\frac{1}{1-\eta_1}\int_{x_0}^\infty f(y,v)m_1(y)\,dy.
\end{align*}
We deduce that 
the operator $\kappa-\BB_0$ is invertible in $L^1_{m_1}$ with $\|(\kappa-\BB)^{-1}\|\leq 1/\kappa$.
We have also proved that $(\kappa-\BB_0)^{-1}$ maps $L^1_{m_1}$ into $L^1_{Km_1}$ with $\|(\kappa-\BB_0)^{-1}\|_{\BBB(L^1_{m_1},L^1_{Km_1})}\leq\frac{1}{1-\eta_1}$.
The fact that it maps $W^{1,1}_{v,loc}$ into $ \Wloc^{1,1}$ readily follows from the formula~\eqref{eq:GF-invB0}.
We conclude to the compactness of $\CC$ and then to the validity of \ref{H3}.
Indeed, we can write~\eqref{eq:lambda1approx} as
\[\hat f_n = \RR_\BB(\lambda_n)\AA \hat f_n + \RR_\BB(\lambda_n)\eps_n,\]
but also as
\[\hat f_n = \RR_{\BB_0}(\lambda_n)\AA_0 \hat f_n + \RR_{\BB_0}(\lambda_n)\eps_n.\]
Combining both, we get
\[\hat f_n = \CC \hat f_n + \big[\RR_{\BB_0}(\lambda_n)\AA_0\RR_\BB(\lambda_n)+\RR_{\BB_0}(\lambda_n)\big]\eps_n.\]
Since $\CC$ is compact, we conclude to \ref{H3} with the same argument as in the proof of Lemma~\ref{lem:H3abstract-StrongC}. 

\medskip

From \ref{H1}, \ref{H2}, \ref{H3} and \ref{H4} we infer the conclusion \ref{S1}-\ref{S2} about existence and uniqueness of $(\lambda_1,f_1,\phi_1)$, which gives a part of Theorem~\ref{theo:GFV-main}.
For the quantitative exponential stability, we start with a lemma.

\begin{lem}\label{lem:GFV-Harris} For all $\eps>0$, $R_1>x_0$, and $R_2>x_0+\eps$, there exist $T>0$ and $c_T>0$ such that 
\begin{equation}\label{eq:GFV-Harris}
S_T^* \phi \ge c_T {\bf 1}_{(x_0,R_1)\times[1,2]} \int_1^2\!\int_{x_0+\eps}^{R_2} \phi\, dxdv, \quad \forall \, \phi \ge 0.
\end{equation}
\end{lem} 

\begin{proof}[Proof of Lemma~\ref{lem:GFV-Harris}]
Let us fix $\eps>0$, $R_1>x_0$, and $R_2>x_0+\eps$.
Throughout the proof we denote by $c_t$ any positive constant that depend on $t$, and also possibly on the ingredients of the model $g$, $K$, $\wp$ and on $\eps,R_1,R_2$, but is independent of $(x,v)\in(x_0,R_1)\times[1,2]$.

\smallskip
{\it First step.}
We prove first that 
there exists $T_1>0$ and $x_3>x_2>x_0$ such that $x_3>\max(R_2,2x_2)$ and
\begin{equation}\label{eq:GFV-Harris-step1}
S_{T_1}^* \phi \ge c_{T_1} {\bf 1}_{(x_0,R_1)\times[1,2]} \int_1^2\!\int_{x_2}^{x_3} \phi\, dxdv,
\end{equation}
for all $\phi\geq0$.
We start from the Duhamel formula
\begin{align}\label{eq:Duhamel-GFV}
S^*_t\phi(x,v)  = \phi&(X^v_t(x),v)e^{-\int_0^tK(X^v_s(x))ds} \\
&+ 2 \int_0^t \int_1^2 K(X^{v}_{s}(x)) S^*_{t-s}\phi(X^{v}_{s}(x)/2,v_*) e^{-\int_0^{s}K(X^{v}_{s'}(x))ds'} \wp(v_*,v)dv_*ds, \nonumber 
\end{align}
where $X^v_t(x)$ is the solution to the characteristic equation
\[\dot X^v_t(x) = v\, a(X^v_t(x))\quad \text{with}\quad X^v_0(x)=x.\]
Iterating twice~\eqref{eq:Duhamel-GFV}, using positivity and the fact that $K$ and $a$ are locally bounded and $\wp$ is bounded from below, we deduce that
\begin{align*}
S^*_t&\phi(x,v) \geq \\
&c_t\int_0^t\!\!\int_1^2\!\!\int_0^s\!\!\int_1^2\! K(X^v_{s'}(x)) K(X^{v_*}_{s-s'}(X^v_{s'}(x)/2)) \phi\big(X^{v_{**}}_{t-s}\big(X^{v_*}_{s-s'}\big(X^v_{s'}(x)/2\big)/2\big),v_{**}\big)\,dv_{**}ds'dv_*ds,
\end{align*}
on $(x_0,R_1)\times[1,2]$, for all $\phi\geq0$.
Let $t_0$ be such that $X^1_{t_0}(x_0)=2x_0+1$.
Then, for $t > 2t_0$, we deduce, from the fact that $K$ is locally bounded from below on $(2x_0,\infty)$, that
\[S^*_t\phi(x,v) \geq c_t \int_{2t_0}^t\int_1^2\int_{t_0}^s\int_1^2 \phi\big(X^{v_{**}}_{t-s}\big(X^{v_*}_{s-s'}\big(X^v_{s'}(x)/2\big)/2\big),v_{**}\big) \,dv_{**}ds'dv_*ds.\]
 For $t>2t_0+2$, by using the Fubini-Tonelli theorem, we thus have
\[S^*_t\phi(x,v) \geq c_t \int_1^2\int_{t-1}^t\int_{t_0}^{t_0+1}\bigg(\int_1^2 \phi\big(X^{v_{**}}_{t-s}\big(X^{v_*}_{s-s'}\big(X^v_{s'}(x)/2\big)/2\big),v_{**}\big) dv_*\bigg) \,ds'dsdv_{**}.\]
Using now a change of variables, we get
\begin{align*}
S^*_t\phi(x,v) &\geq c_t \int_1^2\int_{t-1}^t\int_{t_0}^{t_0+1} \bigg(\int_{X^{v_{**}}_{t-s}(X^1_{s-s'}(X^v_{s'}(x)/2)/2)}^{X^{v_{**}}_{t-s}(X^2_{s-s'}(X^v_{s'}(x)/2)/2)} \phi(y , v_{**})dy\bigg) \,ds'dsdv_{**}\\
& \geq c_t \int_1^2\int_{X^2_1(X^1_{t-t_0}(X^2_{t_0+1}(R_1)/2)/2)}^{X^2_{t-t_0-2}(x_0/2)/2} \phi(y , v_{**})\,dydv_{**}.
\end{align*}
Due to the strict positivity of $a$, we can choose $t=T_1$ large enough so that
\[
X^2_{T_1-t_0-2}(x_0/2)/2 > \max\big(R_2,2X^2_1(X^1_{T_1-t_0}(X^2_{t_0+1}(R_1)/2)/2)\big),
\]
and we obtain~\eqref{eq:GFV-Harris-step1} by setting $x_3=X^2_{T_1-t_0-2}(x_0/2)/2$ and $x_2=X^2_1(X^1_{T_1-t_0}(X^2_{t_0+1}(R_1)/2)/2)$, which concludes the first step of the proof.

\smallskip
{\it Second step.}
We deduce~\eqref{eq:GFV-Harris} from~\eqref{eq:GFV-Harris-step1}  as follows.
On the one hand, applying~\eqref{eq:GFV-Harris-step1} to the function $S_t^*\phi$, we obtain
\[
S_{T_1+t}^* \phi \ge c_{T_1} {\bf 1}_{(x_0,R_1)\times[1,2]} \int_1^2\!\int_{x_2}^{x_3} S^*_t\phi\, dxdv.
\]
On the other hand, iterating once the Duhamel formula~\eqref{eq:Duhamel-GFV},
we get by positivity that
\[
S^*_t\phi(x,v)\geq c_t\bigg[\phi(X^v_t(x),v) + \int_0^t \int_1^2 K(X^{v}_{s}(x)) \phi(X^{v_*}_{t-s}(X^{v}_{s}(x)/2),v_*) dv_*ds\bigg]. 
\]
We first assume $x_2>2x_0$. In that case, the term $K(X^{v}_{s}(x))$ is bounded from below uniformly in $s\in[0,t]$, $v\in[1,2]$ and $x\in[x_2,x_3]$, so that we infer from the two above inequalities that
\[S_{T_1+t}^* \phi \ge c_{T_1}c_t {\bf 1}_{(x_0,R_1)\times[1,2]} \int_1^2\!\int_{x_2}^{x_3}\bigg[\phi(X^v_t(x),v) + \int_0^t \int_1^2 \phi(X^{v_*}_{t-s}(X^{v}_{s}(x)/2),v_*) dv_*ds\bigg]dxdv.\]
By a change of variable, we have
\[\int_{x_2}^{x_3}\phi(X^v_t(x),v)dx \geq c_t \int_{X^2_t(x_2)}^{x_3}\phi(y,v)dy\]
and
\[\int_0^t\int_{x_2}^{x_3}\phi(X^{v_*}_{t-s}(X^{v}_{s}(x)/2),v_*) dxds \geq c_t \int_{X^2_t(X^2_t(x)/2)}^{x_3/2}\phi(y,v_*) dy.\]
Since $X^2_t(x)\to x$ as $t\to0$, we deduce that we can find, for any $\zeta>0$, a time $t>0$ such that
\[S_{T_1+t}^* \phi \ge c_{T_1}c_t {\bf 1}_{(x_0,R_1)\times[1,2]} \bigg[\int_1^2\!\int_{x_2+\zeta}^{x_3}\phi(y,v)dydv + \int_1^2 \int_{x_2/2+\zeta}^{x_3/2} \phi(y,v_*) dydv_*\bigg].\]
As $x_3>2x_2$, we can choose $\zeta$ small enough so that $x_3/2>x_2+\zeta$, and we get
\[S_{T_1+t}^* \phi \ge c_{T_1}c_t {\bf 1}_{(x_0,R_1)\times[1,2]} \int_1^2\!\int_{x_2/2+\zeta}^{x_3}\phi(y,v)dydv.\]
Impose additionally that $\zeta\leq x_0$.
Since the sequence $(u_n)$ defined by $u_0=x_2$ and $u_{n+1}=u_n/2+\zeta$ converges to $2\zeta\leq2x_0$,
we deduce by an iteration argument the existence of a time $T_2$ such that
\begin{equation}\label{eq:GFV-Harris-step2}
S^*_{T_2}\phi(x) \ge c_{T_2} \1_{(x_0,R_1)}(x) \int_1^2\!\int_{2x_0+\eps}^{x_3}\phi(y,v)dydv.
\end{equation}
Using a last time the argument with $\zeta\leq\eps/2$ yields~\eqref{eq:GFV-Harris}, since $x_3>R_2$.

In the case where $x_2 \leq 2x_0$, \eqref{eq:GFV-Harris-step1} directly implies~\eqref{eq:GFV-Harris-step2} with $T_2=T_1$, and only one iteration of the extension argument is enough for concluding.
\end{proof}

We are now in position to finish the proof of Theorem~\ref{theo:GFV-main}.

\begin{proof}[Proof of Theorem~\ref{theo:GFV-main}]
The proof is exactly the same as for Theorem~\ref{theo:GF-regular}, Lemma~\ref{lem:GFV-Harris} replacing Lemma~\ref{lem:GFHarris}.
The only missing information is a quantitative $\Lloc^\infty$ estimate on the derivatives $\partial_v\phi_1$ and $\partial_x\phi_1$, in order to use the same argument as in the proof of Proposition~\ref{prop:GF-isolation} for verifying~\eqref{eq:hyp-Harris-BdBis}.

\smallskip

The estimate on $\partial_x\phi_1$ follows directly from the equation $\LL^*\phi_1=\lambda_1\phi_1$, which also reads
\[
\partial_x\phi_1=\frac{1}{va(x)}\bigg[\lambda_1\phi_1+K\phi_1-2K(x)\int_1^2\phi_1(x,v_*)\wp(v_*,v)dv_*\bigg].
\]
For $\partial_v\phi_1$ we argue by duality, using that
\[\|\phi\|_{L^\infty}=\sup_{\|f\|_{L^1}=1}\langle\phi,f\rangle.\]
We start from
\[
\phi_1=(\lambda_1-\BB_0^*)^{-1}\AA_0^* \, \phi_1,
\]
which yields
\begin{align*}
\|\partial_v\phi_1\|_{L^\infty((x_0,R)\times[1,2]} & = \sup_{\|f\|_{L^1}=1} \langle \partial_v (\lambda_1-\BB_0^*)^{-1}\AA_0^* \, \phi_1 , f \rangle \\
& = \sup_{\|f\|_{L^1}=1} \langle \phi_1 , \AA_0(\lambda_1-\BB_0)^{-1} \partial_v f \rangle
\end{align*}
where the supremum can be taken over the functions $f\in C^1_c((x_0,R)\times(1,2))$.
Using an integration by parts in $v_*$, we have
\begin{align*}
\AA_0(\lambda_1-\BB_0)^{-1}\partial_v f (x,v )& = 4 \int_1^2 \frac{K(2x)}{v_*a(x)} \int_{x_0}^x e^{(\Lambda_{\lambda_1}(y)-\Lambda_{\lambda_1}(x))/v_*}\partial_v f(y,v_*) dy \,\wp(v,v_*) dv_* \\
& = - 4 \int_1^2 \int_{x_0}^x \partial_{v_*}\Big( \frac{K(2x)}{v_*a(x)} e^{(\Lambda_{\lambda_1}(y)-\Lambda_{\lambda_1}(x))/v_*}\wp(v,v_*) \Big) f(y,v_*) dy dv_*.
\end{align*}
Since $\|\phi_1\|_{L^\infty((x_0,R)\times[1,2])} \leq m(R)$,  we deduce 
\begin{align*}
\|\partial_v\phi_1&\|_{L^\infty((x_0,R)\times[1,2])} \\
& \leq 4m(R) \sup_{(y,v_*)\in(x_0,R)\times[1,2]} \int_{x_0}^R \int_1^2 \Big|\partial_{v_*}\Big( \frac{K(2x)}{v_*a(x)} e^{(\Lambda_{\lambda_1}(y)-\Lambda_{\lambda_1}(x))/v_*}\wp(v,v_*) \Big)\Big| dv dx,
\end{align*}
and this last quantity is finite due to the assumptions made on the functions $a,K$ and $\wp$.
\end{proof}

 \bigskip
 \bigskip


\section{The kinetic linear Boltzmann equation}
\label{part:application4:Kequation}

%
%
 
 \bigskip

 
\Black
In this section, we consider the kinetic linear Boltzmann type equation 
\beqn\label{eq:kinetic-eqLinearBoltzmann}
\partial_t f + v \cdot \nabla_x f - \nabla_x \Phi(x) \cdot  \nabla_v f = \KKK[f]- K f, \quad\hbox{in}\quad (0,\infty) \times \OO 
\eeqn
on the function $f = f(t,x,v)$, $t \ge 0$, $(x,v) \in \OO := \Omega \times \R^d$. 
We assume that $K = K(x,v) \ge 0$ and that the collision operator $\KKK$ is linear and defined by 
\beqn\label{eq:kinetic-hypKK}
 \KKK = r \KKK_1, \quad (\KKK_1 g)(x,v) := \int_{\R^d} k \,  g_*  \, dv_* ,
\eeqn
for a real number $r > 0$ and some collision kernel  $k : \Omega \times \R^d \times \R^d \to \R_+$.  
Here and below, we use the common shorthands
$$
g_* := g(v_*), \quad k := k(x,v,v_*), \quad k_* := k(x,v_*,v).
$$
The most classical example for the collisional operator $\CC = \KKK - K$ is the mass conservative operator
\beqn\label{eq:Kinetic-hypOpCons}
(\CC g )(v) := \int_{\R^d} |v-v_*|^\gamma \{ \MMM g_* - \MMM_* g \} \, dv_*, 
\eeqn
for some function $ \MMM \in L^1_+(\R^d)$ and some exponent $\gamma \in \R$, which includes    the relaxation operator 
\beqn\label{eq:Kinetic-DefRelax}
(\CC g)(v) := \sigma ( \MMM \rho_g - g), \quad \rho_g   := \int_{\R^d} g_*  \, dv_* . 
\eeqn
 We make the following strong positivity and boundedness assumption on the collision kernel $k$ and the function $K$. 
 There exist $\gamma \ge 0$ and $K_i > 0$ such that 
\beqn
\label{eq:kinetic-bornesK}
  \forall \, (x,v) \in \Omega \times   \R^d, \quad K_0 \le K(x,v) \langle v \rangle^{-\gamma} \le K_1.
 \eeqn
 There exists a weight function $m : \R^d \to [1,\infty)$ such that 
\beqn
\label{eq:kinetic-bornekp}
\forall \, p \in [1,\infty], \quad k \, m_*^{-1} m \in L^\infty_xL^p_vL^{p'}_{v_*}. 
\eeqn
For all $R > 0$, there exists $k_R > 0$ such that 
\beqn \label{eq:kinetic-bornekinf}
 \forall \, (x,v,v_*) \in \Omega \times B_R \times B_R, \quad k(x,v,v_*) \ge k_R.
\eeqn
It is worth emphasizing that for $\KKK$ and $K$ defined in \eqref{eq:Kinetic-hypOpCons}, the above assumptions are met when $m :=  \MMM^{-1/2} : \R^d \to [1,\infty)$ (so that in particular $\MMM > 0$ a.e.) and $ \MMM^{1/2} \langle v \rangle^\gamma \in L^1 \cap L^\infty$.  We finally assume that for some weight function $m_1 : \R^d \to [1,\infty)$ such that $m_1/m \to \infty$ at infinity, we have  
\beqn\label{eq:kinetic-ARB1}
k \, m_*^{-1} m_1 \in L^\infty_xL^2_{vv_*},  
\eeqn
what holds true for the relaxation operator when $\MMM m_1\in L^2(\R^d)$ and $m^{-1} \in L^2(\R^d)$, 
and  that for some weight function $m_0 : \R^d \to [1,\infty)$ such that $m_0/m \to 0$ at infinity, we have
\beqn\label{eq:kinetic-ARBinfty} 
k \, m_{0*}^{-1} m \in L^\infty_xL^1_{vv_*}, 
\eeqn
what holds true for the relaxation operator when $\MMM m  \in L^1(\R^d)$ and $m_0^{-1} \in L^1(\R^d)$.
 
%


\smallskip
For the space domain $\Omega$, we consider the two following cases: 
\bean
&(1)& \Omega := \T^d, \ \hbox{the torus};
\\
&(2)& \Omega := \R^d, \ \hbox{the whole space}. 
\eean

\smallskip
In case (1), and for the sake of simplicity, we will always assume that $\Phi = 0$. 
 In case (2), we will need a confinement mechanism which will be provided by the 
mean of the confinement force associated to the confinement potential $\Phi$. 
We do not consider here the case of a bounded domain with zero influx boundary condition because (1) our approach applies exactly as  for the torus case
and (2) this case has already been considered in the pioneering work by Vidav~\cite{MR230531}, where existence, uniqueness and exponential stability (with non constructive constants) have been established. 
We do not consider either the case of a bounded domain complemented with a reflection as we will consider in Section~\ref{part:application5:KFPequation} for the kinetic Fokker-Planck evolution equation,
because we have not been able to establish some crucial regularity estimates which seem to be necessary in our approach. 

\medskip
\subsection{The torus}\label{subsec:kinetic-torus}
In this section, we are first concerned with the kinetic linear Boltzmann equation in the torus 
\beqn\label{eq:Transport-eqtorus}
\partial_t f + v \cdot \nabla_x f = \KKK[f]- K f, \quad\hbox{in}\quad (0,\infty) \times \T^d \times \R^d.
\eeqn
We make the boundedness and strong positivity assumptions listed above together with the additional assumption
\beqn
\label{eq:kinetic-bornekinfty}
k \, m_*^{-1} m _1, k \, m_{1*}^{-1/2} m  \in    L^\infty_{xvv_*}, \quad \sum_{u \in \Z^d} \| m_1^{-1/2} (u + \cdot) \|_{L^\infty(\T^d)} < \infty,
\eeqn
 for some $m_1$ such that $m/m_1\in L^1\cap L^2$.  

\begin{theo}\label{theo:kinetic-torus-main}
For the kinetic equation  \eqref{eq:Transport-eqtorus} in the torus and under conditions \eqref{eq:kinetic-bornesK}-\eqref{eq:kinetic-bornekp}-\eqref{eq:kinetic-bornekinf}-\eqref{eq:kinetic-ARB1}-\eqref{eq:kinetic-ARBinfty} and \eqref{eq:kinetic-bornekinfty}
for some weight functions $m$, $m_0$, $m_1$, there exists $r^* > 0$ such that for any $r \ge r^*$,  
 the conclusions  \ref{S1},  \ref{S2},  \ref{S32}   hold in $L^2_m$ and the conclusion \ref{E31} holds in $L^1_m$. 
\end{theo}

Our result may be compared to \cite{MR230531} which establishes the same result without constructive rate and to  \cite{MR3449390} which establishes the same result using a probabilistic approach,
 both in the case of a bounded domain with zero influx boundary condition.
 It also extends to a non mass conservative situation the many results devoted to the conservative framework, see for instance the recent papers  \cite{MR3192457,MR3479064,MR4063917} and the references
therein. When $\gamma > 0$, we may  probably establish the same above result under the sole condition $r > 0$ (no need for $r$ to be large enough) by using some arguments developed in the next section.


\smallskip
We present now the proof of Theorem~\ref{theo:kinetic-torus-main} by establishing that 
the conditions presented in the abstract part are satisfied. 




\medskip
{\bf Condition \ref{H1}.}  For an exponent $p \in [1,\infty)$ and a weight function $m$ satisfying \eqref{eq:kinetic-bornekp}, we set 
$$
k_\infty := \| k m_*^{-1} m \|_{L^\infty_xL^p_vL^{p'}_{v_*}}   < \infty.
$$
Considering then a solution $f$ to the evolution equation \eqref{eq:Transport-eqtorus}, we compute
\bean
\frac{1 }{ p}  \frac{d}{ dt} \int f^p m^p 
&=&
 \int  \KKK [f]  f^{p-1} m^p - K(x,v) f^p m^p
\\
&\le&
\| \KKK[f] m\|_{L^p}  \Bigl( \int  f^p m^p  \Bigr)^{1-1/p} -   \int K(x,v) f^p m^p
\\
&\le&
r k_\infty  \int  f^p m^p   - K_0 \int \langle v \rangle^\gamma f^p m^p,
\eean
where we have used twice the Holder inequality. This differential inequality together with the Gronwall lemma provides an a priori estimates
about the growth of the $L^p_m$ norm. 
As a consequence, the same arguments as in section~\ref{subsec:TranspEq-tau} imply that   $S_\LL$ is a positive semigroup in $L^p_m$ with 
growth bound $\omega(S_\LL) = r k_\infty  - K_0$. In particular, condition \ref{H1} holds thanks to Lemma~\ref{lem:Exist1-RkSG}. 


\medskip
{\bf Condition \ref{H2}.} 
  For $f_0 := {\bf 1}_{\T^d \times B_1}$, where $B_1$ denotes the unit ball in $\R^d_v$,  
we compute 
$$
\LL f_0 = \KKK[ f_0]  - K f_0  \ge  \inf_{v \in B_1} \bigl\{  r \KKK_1[ f_0]  - K\} f_0.
$$
Using \eqref{eq:kinetic-bornesK} and the strong positivity condition \eqref{eq:kinetic-bornekinf}, we get 
\beqn\label{eq:kinetic-torusH2lowerkappa0}
 \inf_{v \in B_1} \bigl\{ \KKK[ f_0]  - K\} \ge  r k_1 - 2^{\gamma/2} K_0= : \kappa_0, 
\eeqn
which provides a constructive lower bound of the set $\II$ defined in \eqref{eq:exist1-defI} thanks to Lemma~\ref{lem:Existe1-Spectre2bis}-{\bf (ii)}. We have thus established that   $\LL$ satisfies \ref{H2}.

\medskip
{\bf Condition \ref{H3}.}  
We define the operator 
$$
\BB f := - v \cdot \nabla_x f  - K(x,v) f, 
$$
and we assume  $\kappa_\BB := - \inf K  \le -K_0 < \kappa_0$, what holds whenever $r \ge r^*$, with $r^*> 0$ large enough thanks to \eqref{eq:kinetic-torusH2lowerkappa0}. 
In the sequel, we assume $p=2$ and we work in $X = L^2_m$. We immediately deduce that $\BB-\kappa$ is dissipative for any $\kappa > \kappa_\BB$, and thus  $\RR_\BB(z)$ is bounded in $\BBB(L^2_m)$, uniformly in $z \in \Delta_\kappa$. For $\kappa > \kappa_\BB$ and $g \in L^2_m$, the function $f = \RR_\BB(\kappa) g$ satisfies
$$
v \cdot \nabla_x f + (\kappa + K) f = g \quad\hbox{in} \quad \OO,
$$
from what we deduce 
$$
(\kappa - \kappa_\BB)  \int_\OO f^2 m^2 \le  \int_\OO (\kappa + K)  f^2 m^2 = \int_\OO fg m^2, 
$$
and finally 
$$
\| f \|_{L^2_m}^2 \le \frac{1 }{ \kappa - \kappa_\BB} \| g \|_{L^2_m}^2. 
$$
Because of assumption \eqref{eq:kinetic-ARB1}, and defining $\AA := \KKK$, we immediately deduce that 
\beqn\label{eq:kinetic-ARB2}
\AA \RR_\BB (\kappa)  : L^2_m \to L^2_{m_1}.
\eeqn

On the other hand, from the classical averaging lemma \cite{MR923047}, we know that 
\beqn\label{eq:kinetic-ARB3}
\AA_\varphi \RR_\BB(\kappa) : L^2(\OO) \to H^{1/2}(\OO),
\eeqn
where for $\varphi = \varphi_1 \otimes \varphi_2 \in C^1_c(\OO) \otimes C^1_c(\R^d)$, we have defined  the mapping $\AA_\varphi : L^2(\OO) \to L^2(\OO)$ by 
$$
\AA_\varphi (f) (x) := \varphi_1(x,v)\int_{\R^d} f (x,v_*) \varphi_2(v_*) \, dv_*.
$$
By classical approximation arguments, there exists a sequence $(\varphi_n)$ such that $\varphi_n \to k$ in the space $L^\infty(\T^d;L^2_{m_1 \otimes m^{-1}} (\R^d \times \R^d))$ and such that $\varphi_n$ is a linear combination of functions of $C^1_c(\OO) \otimes C^1_c(\R^d)$. 
As a consequence of \eqref{eq:kinetic-ARB2} and \eqref{eq:kinetic-ARB3}, we deduce that $\AA \RR_\BB (\kappa) \in \KKK(L^2_m)$ and next 
$(\RR_\BB (\kappa) \AA)^2 \in \KKK(L^2_m)$ for any $\kappa > \kappa_\BB$. We may   use  Lemma~\ref{lem:H3abstract-StrongC} (and Remark~\ref{rem:H3abstract-compactStrong}-(2)) with $N=2$ in $X = L^2_m$, and deduce that \ref{H3} holds. 

%

\medskip
{\bf Condition \ref{H4}.}  
We start with a result of independent interest about strict positivity. Such an argument is reminiscent from \cite{MR1555365,MR1461954} in the study of the Boltzmann equation
and has been used for instance in \cite{MR2153518,MR4063917}.

\begin{lem}\label{lem:kinetic-torusDoblin}
For any $\varrho,\varrho_*,t >0$, there exists $c > 0$ such that 
\beqn\label{eq:kinetic-torusDoblin}
(S_\LL(t) f_0)(x,v) \ge c\, \1_{B_\varrho}(v) \int_{\T^d \times B_{\varrho_*}}\! f_0 \,dv_*dx_*,
\eeqn
for all $f_0\geq0$ and $(x,v)\in\T^d\times\R^d$.
\end{lem}

\begin{proof}[Proof of Lemma~\ref{lem:kinetic-torusDoblin}]
We observe that   the semigroup $S_\BB$ has explicit representation
$$
(S_\BB(t) f_0 ) (x,v) = f_0(x-vt,v)e^{- \int_0^t K(x-\tau v,v) d\tau}.
$$
 We next write the associated iterated Duhamel formula
  $$
 S_\LL = S_\BB + S_\BB  \KKK  * S_\BB + S_\BB  \KKK  * S_\BB  \KKK  * S_\BB + S_\BB  \KKK  *  S_\BB  \KKK  * S_\BB  \KKK  * S_\LL.
 $$
Since all the terms are nonnegative, we may through away the first terms  and the last one, and we get  
$$
 S_\LL \ge  S_\BB   *  \KKK  S_\BB *  \KKK   S_\BB. 
 $$
 On the one hand, using the explicit expression of $S_\BB$ and \eqref{eq:kinetic-bornekinf}, we have
 $$
 ( \KKK_1   S_\BB(s) f_0) (y,w) \ge k_{\varrho'} e^{-s K_{\varrho_*}}  \1_{B_{\varrho'}}(w) \int_{B_{\varrho_*}} f_0(y-w_* s,w_*) dw_* =: g(s,y,w),
 $$
for any $s  > 0$ and any $\varrho'\geq\varrho_*>0$,  with $K_\varrho := \sup_{z \in \T^d, |v| \le \varrho} K(z,v)$. 
On the other hand, for the same reasons, we have 
\begin{align*}
( \KKK_1   S_\BB * g (t)) (x,v) & = \int_0^t \int_{\R^d} k(x,v,v_*)   g (s,x-v_*(t-s),v_*) e^{-\int_0^{t-s}K(x-\tau v_*,v_*)d\tau} dv_*ds \\
& \ge k_{\varrho'} \1_{B_{\varrho}}(v) \int_0^t \int_{B_{\varrho'}}  g (s,x-v_*(t-s),v_*) e^{-(t-s)K_{\varrho'}} dv_*ds \\
& \ge k_{\varrho'}^2 e^{-t K_{\varrho'}}  \1_{B_{\varrho}}(v) \int_0^t  \int_{B_{\varrho'}}    
   \int_{B_{\varrho_*}} f_0(x-v_*(t-s)-w_* s,w_*) dw_* dv_*ds \\
& \ge k_{\varrho'}^2 e^{-t K_{\varrho'}} \1_{B_{\varrho}}(v) \int_0^{t/2} 
   \int_{B_{\varrho_*}}  \int_{B(x+w_*s,(t-s)\varrho')}     f_0(y_*,w_*) \frac{dy_* }{ (t-s)^d} dw_*  ds\\
& \ge k_{\varrho'}^2 \frac{e^{-t K_{\varrho'}}  }{ (t/2)^d}  \1_{B_{\varrho}}(v) \int_0^{t/2} 
   \int_{B_{\varrho_*}}  \int_{\T^d}     f_0(y_*,w_*) dy_* dw_*  ds\\
& \ge k_{\varrho'}^2 \frac{e^{-t K_{\varrho'}}  }{ (t/2)^d}  \1_{B_{\varrho}}(v) \frac{t }{ 2}  
   \int_{B_{\varrho_*}}  \int_{\T^d}     f_0(y_*,w_*) dy_* dw_*  , 
\end{align*}
for any $t > 0$ and $\varrho'\geq\max(\varrho,\varrho_*)$ such that $t \varrho'\ge 2$, in such a way that $B(z,(t/2)\varrho') \supset \T^d$. We then have 
\bean
 ( \KKK_1   S_\BB * \KKK_1   S_\BB (t)) (x,v) 
 &\ge& k^2_{\varrho'} \frac{e^{-t K_{\varrho'}}  }{ (t/2)^{d-1}} \1_{B_\varrho}(v)  
  \int_{\T^d}   \int_{B_{\varrho_*}}  f_{0*} \, dw_* dx_*
\eean
for any $t\geq2/\varrho'$.  We finally conclude
$$
S_\LL(t) f_0 (x,v) \ge 
r^2 k^2_{\varrho'} e^{-t K_{\varrho'}} \int_{2/\varrho'}^t  \frac{ds }{ (s/2)^{d-1}}  \, \1_{B_\varrho}(v)  
  \int_{\T^d}   \int_{B_{\varrho_*}}    f_{0*}  \, dw_* dx_*,
$$
from what we deduce \eqref{eq:kinetic-torusDoblin} by choosing $\varrho' = 8/t$.
\end{proof}

\medskip
We now consider $\lambda \ge \lambda_1$ and $0 \le f \in L^2_m$, $f\not\equiv 0$, such that 
$$
\lambda f + v \cdot \nabla_x f + K f -  \KKK [f]   \ge 0 \quad\hbox{in}\quad  \T^d \times \R^d. 
$$
We fix $\varrho_* > 0$ such that $f \not\equiv 0$ on $B_{\varrho_*}$. From \eqref{eq:Exist1-DefRepresentationRR}, we have 
$$
f \ge \int_0^{\infty} e^{-(1+\lambda)t} S_\LL(t) f dt, 
$$
and we conclude that $f > 0$ a.e. on any set $\T^d \times B_\varrho$ thanks to Lemma~\ref{lem:kinetic-torusDoblin}. We have established that the strong maximum
holds true, and thus \ref{H4}. 


\medskip
{\bf Condition \ref{H5}.}
Assume that $(\lambda,f) \in \C \times X\backslash\{0\}$ satisfies  
$$
\LL f = \lambda f, \quad
\LL |f|  = (\Re e \lambda) |f|  = \Re e (\hbox{\rm sign} f) \LL f. 
$$
From \ref{H4} and the first identity satisfied by $|f|$,  we know that $|f| > 0$ a.e. on $\T^d \times \R^d$. 
Using the second identity,  we get
$$
\KKK[|f|] =  \Re e (\hbox{\rm sign} f) \KKK[f] . 
$$
Writing $f = e^{i\alpha} |f|$, we deduce 
$$
\int_{\R^d} k(x,v,v_*) |f(x,v_*)| (1 - \cos (\alpha-\alpha_*)) dv_* = 0 \quad \hbox{a.e. on } \T^d \times \R^d, 
$$
and thus $\alpha = \alpha(x)$ thanks to \eqref{eq:kinetic-bornekinf}. Next, coming back to the first equation, we have 
\bean
\lambda |f| e^{i\alpha} 
&=& \LL ( |f| e^{i\alpha}) 
\\
&=&  e^{i\alpha} \LL |f| - |f| e^{i\alpha} i v \cdot \nabla_x \alpha
\\
&=&  e^{i\alpha} (\Re e \lambda) |f|  - |f| e^{i\alpha} i v \cdot \nabla_x \alpha.
\eean
The equation simplifies into 
$$
v \cdot \nabla_x \alpha = \Im m  \lambda,
$$
so that $\alpha(x) = \alpha$ is a constant and the reverse Kato's inequality holds. 

\medskip

Alternatively to~\ref{H5}, we readily infer from Lemma~\ref{lem:kinetic-torusDoblin} that the variant condition~\ref{H5'} is verified.

\medskip
At this stage, because of Theorem~\ref{theo:exist1-KRexistence}, Theorem~\ref{theo:KRgeometry1} and Theorem~\ref{theo:KRgeometry2} (or Theorem~\ref{theo:lambda1largestBIS}), we deduce the conclusions 
{\Blue  \ref{S1},  \ref{S2} and   \ref{S32}}
about the existence and uniqueness of the eigentriplet $(\lambda_1,f_1,\phi_1)$ which satisfies $f_1 > 0$, $\phi_1 > 0$, $\lambda_1$ is algebraically simple and   on the triviality of the boundary punctual spectrum. We now establish the exponential asymptotic stability with constructive constants.

%

\medskip
We start with a gain of uniform boundedness estimate. 

\begin{lem}\label{lem:ch3:kineticB:L1Linfty} 
There exists $N \ge 1$ such that  $(\AA \RR_\BB(\kappa))^N : L^1_m \to L^\infty_{m}$, for any $\kappa > \kappa_\BB$.
As a consequence,  $\phi_1 \in L^\infty_{m^{-1}}$.

\end{lem}

\begin{proof}[Proof of Lemma~\ref{lem:ch3:kineticB:L1Linfty}.] {\sl Step 1.} We argue similarly as in \cite[Sec.~3.1]{MR3591133}. On the one hand,
denoting  $\AA_1 = \KKK_1$, so that $\AA=r\AA_1$, we have for any $f_0 \in L^1_m$
\[(\AA_1 S_{\BB} (t) f_0) (x,v)
=  \int_{\R^d} k(x, v, v_*) f_0(x-v_*t,v_*) \, e^{- \int_0^t K(x-v_*\tau,v_*) d\tau} \, dv_*  \]
and, using estimate \eqref{eq:kinetic-bornekinfty}, we deduce that
\bean
\| m_1 \AA_1 S_\BB(t) f_0 \|_{L^1_xL^\infty_v}
&\le & \| k m_*^{-1} m_1 \|_{L^\infty_{xvv_*}} \int_{\OO} |f_0(x-v_*t,v_*)| m_* dv_* dx \, e^{ t \kappa_\BB}
\\
&\lesssim& \| f_0 \|_{L^1_m} \, e^{ t \kappa_\BB},
\eean
for any $t \ge 0$. Now, we consider  $f_0 \in L^1_xL^\infty_{vm}$, we write
\[
(\AA_1  S_{\BB} (t) f_0) (x,v) =  \sum_{u \in \Z^d} \int_{\T^d}  k( x,v,v_*) 
f_0(x-ut - v_*t,u+v_*) \, e^{- \int_0^t K(x-(u+v_*)\tau,u+v_*) d\tau} \, dv_*
\]
and using estimate \eqref{eq:kinetic-bornekinfty} again, we compute 
%
\begin{align*}
(\AA_1  S_{\BB} (t) f_0) (x,v) m(v)
&\le   \| k m_{1*}^{-1/2} m \|_{L^\infty_{xvv_*}}  \sum_{u \in \Z^d} \int_{\T^d}   m_1^{1/2}(u+v_*) f_0(x-ut - v_*t,u+v_*) \, e^{ t \kappa_\BB} \, dv_*
\\
&\lesssim e^{ t \kappa_\BB} \Bigl( \sum_{u \in \Z^d} \| m_1^{-1/2} (u + \cdot) \|_{L^\infty(\T^d)}\Bigr)  \int_{t \T^d} \|  f_0(y, \cdot) \|_{L^\infty_{m_1}} \frac{dy }{ t^d}
\\
&\lesssim e^{ t \kappa_\BB} \Bigl(  1 + \frac1{t^d}\Bigr) \| f_0  \|_{L^1_xL^\infty_{vm_1}}. 
\end{align*}
Defining $\widetilde u(t) := e^{-\kappa t} \AA_1 S_\BB(t)$, $\kappa > \kappa_\BB$, we have first established $\widetilde u : L^1_m \to L^1_xL^\infty_{vm_1}$ uniformly in time, and
thus $\widetilde u : L^1_xL^\infty_{vm_1} \to L^1_xL^\infty_{vm_1}$ uniformly in time because $L^1_xL^\infty_{vm_1} \subset L^1_m$  (we use here the fact that $m/m_1\in L^1$).
On the other hand, we have established 
that 
$t^{d} \widetilde u : L^1_xL^\infty_{vm_1} \to L^\infty_{m}$ uniformly in time. 
 Using  \cite[Prop.~2.5]{MR3465438} with $X:= L^1_xL^\infty_{vm_1}$ and $Y :=L^\infty_m$, we deduce $\widetilde u^{*(d+1)} : L^1_xL^\infty_{vm_1} \to L^\infty_{m}$ uniformly in time, and we thus conclude that 
$\widetilde u^{*N} : L^1_m \to L^\infty_{m}$ uniformly in time, for $N := d+2$. 
Using formula \eqref{eq:Exist1-DefRepresentationRR}, we deduce that $(\AA \RR_\BB)^N (z) : L^1_m \to L^\infty_m$, 
uniformly for any $z \in \Delta_\kappa$. 

\smallskip
{\sl Step~2.} In particular, $(\AA \RR_\BB)^N (z) : L^1_m \to L^2_m$ because $L^\infty_{m_1} \subset L^2_m$ (we use here the fact that $m/m_1\in L^2$).  
By duality, we deduce that  $( \RR_{\BB^*} \AA^*)^N (z) : L^2_{m^{-1}} \to L^\infty_{m^{-1}}$. Coming back to the eigenvalue equation 
$$
\AA^* \phi_1 + \BB^* \phi_1 = \lambda_1 \phi_1, 
$$
we deduce 
\beqn\label{eq:kinetic-phi1=ARBNphi1}
\phi_1 = \RR_{\BB^*}(\lambda_1)\AA^* \phi = (\RR_{\BB^*}(\lambda_1)\AA^*)^N \phi_1.
\eeqn
By construction $\phi_1 \in L^2_{m^{-1}}$ and we thus conclude that $\phi_1 \in L^\infty_{m^{-1}}$.
\end{proof}

%
%
%
%
%
%
%
%
%

\medskip
From now on, we choose the normalization convention $ \| \phi_1 \|_{L^\infty_{m^{-1}}} = 1$ and $\langle f_1,\phi_1 \rangle = 1$. 

\medskip 
 Because of \eqref{eq:kinetic-ARBinfty} and proceeding similarly as during the proof of condition \ref{H3}, we have 
 $$
 \AA \RR_\BB (\kappa) : L^1_{m_0} \to L^1_{m}, \quad \forall \, \kappa > \kappa_\BB, 
 $$
 so that 
 $$
  \RR_{\BB^*} (\kappa) \AA^* : L^\infty_{m^{-1}}  \to  L^\infty_{m_0^{-1}}  ,  \quad \forall \, \kappa > \kappa_\BB.
 $$
From the first identity in \eqref{eq:kinetic-phi1=ARBNphi1}, 
 we deduce 
 $$
 \| \phi_1 \|_{L^\infty_{m^{-1}_0}}  \le C_{01}  \| \phi_1 \|_{L^\infty_{m^{-1}}}, 
 $$
 with constructive constant $C_{01} \in (0,\infty)$. 
 We may here proceed along an already used argument. Consider $0 \le f \in L^1_m$ and assume  $\| f \|_{L^1_m} \le A [f]_{\phi_1}$. We then compute 
  \bean
 \int f \phi_1 &=& \int_{\OO_\varrho}  f \frac{\phi_1}{m} m +  \int_{\OO_\varrho^c}   f m  \frac{\phi_1 }{ m_0} \frac{m_0}{ m}
 \\
 &\le& \langle f, {\bf 1}_{\OO_\varrho}\rangle  \sup_{\OO_\varrho} m + \| f \|_{L^1_m} C_{01} \sup_{\OO_\varrho^c } \frac{m_0 }{ m}
 \\
 &\le& \langle f, {\bf 1}_{\OO_\varrho} \rangle  \sup_{\OO_\varrho} m + \frac12 [f]_{\phi_1}
 \eean
 by choosing $\varrho = \varrho(A)$ large enough, where we denote $\OO_\varrho := \T^d \times B_{\varrho}$. 
 Together with Lemma~\ref{lem:kinetic-torusDoblin}, 
 we deduce that there exists $T > 0$ and $g_A \ge 0$, $g_A \not\equiv 0$, such that 
 \beqn\label{eq:kinetic-torusHarris}
 S_\LL(T) f   \ge g_A  [f]_{\phi_1}, 
\eeqn
what is nothing but the Doblin-Harris condition~\eqref{eq:hyp-Harris}. On the other hand, from the above regularization estimate,
 we have in the same time
$$
1 = \| \phi_1 \|_{L^\infty_{m^{-1}}} \le C_0 \| \phi_1 \|_{L^1_{m^{-1}}}, \quad
 \| \phi_1 \|_{L^1_{m_0^{-1}}} \le C_1  \| \phi_1 \|_{L^1_{m^{-1}}},
 $$
 for some constructive constants $C_i \in (0,\infty)$. We may thus compute 
 \bean
 \int \phi_1 m^{-1} &\le& \int_{B_\varrho} \phi_1 m^{-1} + \sup_{B_\varrho^c} \frac{m_0 }{ m} \int \phi_1 m_0^{-1}
 \\
 &\le&  \int_{B_\varrho} \phi_1 m^{-1} + \sup_{B_\varrho^c} \frac{m_0 }{ m} C_1 \int \phi_1 m^{-1}, 
 \eean
 so that for $\varrho > 0$ large enough, we deduce 
 \begin{equation}\label{eq:kin:g0phi1lowbound}
 C_0^{-1} \le  \|\phi_1 \|_{L^1_{m^{-1}}} \le 2  \int_{B_\varrho} \phi_1 m^{-1}.
 \end{equation}
 Together with the definition of $g_A$, we deduce that the positivity condition~\eqref{eq:hyp-Doblin-BdBis} holds.

 \medskip

 Finally, as during the proof of \ref{H3} above, for $0 \le f_0 \in L^1_m$ and denoting $f := S_\LL(t) f_0$, we compute 
\bean
 \frac{d }{ dt} \int f m dvdx 
&=&  \int \KKK[f] m dvdx - \int K f m dxdv
\\
&\le& C_0 \int f m_0  dvdx + \kappa_\BB \int f m dxdv, 
\eean
with $C_0 := \| k m m_{0*}^{-1} \|_{L^\infty_{xv_*}L^1_v} < \infty$ and $m_0/ m \to 0$ as $v \to \infty$. Observing that 
$$
\int f m_0 \le \int_{B_\varrho} f \phi_1 \sup_{B_\varrho} \frac{m_0 }{ \phi_1} + \int_{B^c_\varrho} f m \sup_{B^c_\varrho} \frac{m_0 }{ m},
$$
for any $\kappa > \kappa_\BB$, we may choose $\varrho > 0$ large enough in such a way that $ \sup_{B^c_\varrho} \frac{m_0 }{ m} \le (\kappa - \kappa_\BB)/C_0$
and we deduce that 
\bean
 \frac{d }{ dt} \int f m dvdx 
 &\le& C_1 \int f \phi_1 dvdx + \kappa \int f m dxdv
\eean
with $C_1=\sup_{B_\varrho} \frac{m_0 }{ \phi_1}$.
From the Gronwall lemma, we obtain 
$$
\| f_t \|_{L^1_m} \le e^{\kappa t} \| f_0 \|_{L^1_m} + \frac{e^{\lambda_1 t} - e^{\kappa t} }{ \lambda_1 - \kappa} C_1 \int f_0 \phi_1, 
$$
from what we immediately deduce that $S_\LL$ satisfies the Lyapunov condition \eqref{eq:stabilityKR-Lyapunov} for any $t > 0$.
It remains to quantify the constant $C_1$.
The dual formulation of~\eqref{eq:kinetic-torusDoblin} applied to the dual eigenfunction $\phi_1$ with $t=1$ and $\varrho_*=\varrho$ yields
\[\phi_1=e^{-\lambda_1}S^*_\LL(1)\phi_1\geq e^{-\lambda_1}c\, \1_{\T^d\times B_{\varrho}}\int_{\T^d\times B_\varrho}\!\phi_1 \,dv_*dx_*.\]
Together with Equation~\eqref{eq:kin:g0phi1lowbound}, this provides the explicit bound $C_1\leq 2C_0e^{\lambda_1}c^{-1}\frac{\sup_{B_\varrho}m_0}{\inf_{B_\varrho} m}$.

\medskip
We have established that the three conditions \eqref{eq:hyp-Harris}, \eqref{eq:stabilityKR-Lyapunov} and \eqref{eq:hyp-Harris-BdBis} hold, so that we conclude the proof of Theorem~\ref{theo:kinetic-torus-main} by just applying Theorem~\ref{theo:Harris}.

\medskip
\subsection{The whole space case}\label{subsec:kinetic-rd}

\
In this section, we assume that $\Omega := \R^d$ and we consider the kinetic equation~\eqref{eq:kinetic-eqLinearBoltzmann} with  an additional force field confinement $F = - \nabla_x \Phi$ associated to a potential $\Phi$.
More precisely from now-on, we   assume that 
\beqn
\label{eq:kinetic-Hyp-Phi}
  \Phi(x) = |x|^\beta, \ \beta > 2, \quad K(v) = \langle v \rangle^\gamma , \ \gamma > 0, 
\eeqn
that \eqref{eq:kinetic-bornekp} holds (for $p=2$) and that there exist $\zeta, c_\zeta > 0$ such that 
\beqn\label{eq:kineticRd-borneInfKK}
\KKK[\MMM^\zeta] \ge  c_\zeta \langle v \rangle^\gamma \MMM^\zeta, \quad \MMM := e^{-|v|^2/2}.
\eeqn
Observe that condition \eqref{eq:kineticRd-borneInfKK} is satisfied when $\KKK$ is the positive part of the mass conservative operator \eqref{eq:Kinetic-hypOpCons}. 
\Black
For further references, we write  
$
\LL  := \TT  + \CCC 
$
with 
\bean
\TT  := - v \cdot \nabla_x f + \nabla_x \Phi \cdot \nabla_v , \quad 
\CCC f  := \KKK[f] - Kf 
\eean
and we define the Hamiltonian 
$$
\HH :=  \Phi(x) +\frac12 |v|^2. 
$$
In the sequel, we will only consider some weight functions $m = \omega(\HH)$ with $\omega(y) = y^\varrho$, $\varrho \ge 0$, or $\omega(y) = e^{\kappa y}$, $\kappa \in (0,1)$, so that  $\omega(\HH) \sim \omega(|v|^2)\omega(\Phi)$.  For $p \in [1,\infty)$, we further assume that 
$v \mapsto \omega^{-1}(|v|^2) \in L^{p'}$ (which imposes $\varrho > d/(2p')$ for a polynomial weight).

\begin{theo}\label{theo:kineticRd-main}
\Red Consider the the kinetic equation  \eqref{eq:kinetic-eqLinearBoltzmann} in the whole space with confinement force and under conditions \eqref{eq:kinetic-bornesK}-\eqref{eq:kinetic-bornekp}-\eqref{eq:kinetic-bornekinf}-\eqref{eq:kinetic-ARB1}-\eqref{eq:kinetic-Hyp-Phi}-\eqref{eq:kineticRd-borneInfKK} for some weight function $m = \omega(\HH)$ as discussed above. 
There exists $r^* > 0$ such that for any $r \ge r^*$, 
the conclusions 
   \ref{S1},  \ref{S2},  \ref{S32}
 about existence, uniqueness and positivity of the eigentriplet solution $(\lambda_1,f_1,\phi_1)$ hold in $L^p_m$ as well as the ergodicity \ref{E2} for the weak convergence in $L^1_{\phi_1}$.
\end{theo}

We are not aware of any result on the first eigentriplet problem for such linear Boltzmann like equation in the whole space. We may however compare our result to  \cite{MR3479064}
where the corresponding mass conservative framework is considered. 
We present the proof of Theorem~\ref{theo:kineticRd-main} in that situation by adapting the arguments presented in the previous section.

\medskip
{\bf Condition \ref{H1}.}
Let us consider a weight function $m = \omega(\HH)$ as introduced before and let us fix $p \in [1,\infty)$. 
For a solution to the evolution equation \eqref{eq:kinetic-eqLinearBoltzmann}, we classically compute 
\bean
\frac{d }{ dt} \int \frac{f^p}{ p}  m^p dvdx 
&=&  \int (\LL f) f^{p-1} m^p dvdx 
\\
&=&  \int  (f^p/p) \TT^* m^p  + \int (\KKK f) f^{p-1} m^p -  \int  f^p K m^p 
\\
&\le&  \frac1p  \int (\KKK f)^p  m^p +  \int  f^p \bigl( \frac1{p'} - K\bigr) m^p,
\eean
by using an integration by parts and the Young inequality. 
For the first term, we have 
\bean
 \int (\KKK f)^p  m^p dvdx
&\le&  c_\omega  \int  \omega(\Phi) \Bigl( \int f dv \Bigr)^p   dx 
\\
&=&  c_\omega  \int  \omega(\Phi)  \int f^p \omega(|v|^2)  dv   dx \| \omega^{-1}(|v|^2) \|_{L^{p'}}^p 
\\
&\lesssim& \int  f^p  m^p dvdx.
\eean
All together and thanks to the Gronwall lemma, we have established an a priori estimate on the evolution of the norm $\| f \|_{L^p_m}$ and we deduce as in section ~\ref{sec:Transport} that $\LL$ generates a positive semigroup on $L^p_m$.
In particular, the condition \ref{H1} is satisfied thanks to Lemma~\ref{lem:Exist1-RkSG}. 

\medskip 
{\bf Condition \ref{H2}.} We define $f_0 := e^{- \zeta \HH}$ and we compute 
\bean
\LL f_0 &=& \CCC f_0 = r e^{-\zeta \Phi} \KKK[e^{-\MMM^\zeta}] - K e^{-\zeta \Phi} \MMM^\zeta
\\
&\ge&  (rc_\zeta - K_0) \langle v \rangle^\gamma e^{-\zeta \HH} \ge 0,
\eean
for $r > 0$ large enough.  That implies that $\II$ is lower bounded by $\kappa_0 = 0$  by using  Lemma~\ref{lem:Existe1-Spectre2bis}-{\bf (ii)}, and we have thus established that   $\LL$ satisfies \ref{H2}.

%
%

\medskip
{\bf Condition \ref{H3}.} We introduce the collisionless operator
$$
\BB f := \TT f - K f
$$
and we define
$$
\BB^\sharp \phi := \frac12\TT^* \phi - K \phi.
$$
Our analysis is mainly a consequence of the following moment estimate. 

\begin{lem}\label{lem:kinetic-rd-moment}
There exist some weight functions $w\lesssim \HH$ and some real numbers $\alpha, c_\alpha, C_\alpha> 0$ such that 
\beqn\label{eq:lem:kinetic-rd-moment}
\BB^\sharp w \le  C_\alpha  w -c_\alpha   \HH^{1+\alpha}.
\eeqn
\end{lem}

\begin{proof}[Proof of Lemma~\ref{lem:kinetic-rd-moment}.]
We split the proof into two steps. 

\smallskip\noindent
{\sl Step 1.} We first assume $\gamma \le \beta-2$ 
and we define
$$
w :=   1 + \frac12 [x]^{1+\gamma/2} \cdot v + \HH,
$$
with $ [x]^\delta := x |x|^{\delta-1}$. 
We observe that 
$$
\big| [x]^{1+\gamma/2} \cdot v \big| \le \frac12 |x|^{2+\gamma} + \frac12 |v|^2 \le \frac12 \HH + \frac12,
$$
so that $w \sim  \langle \HH \rangle$. We now compute
$$
\TT^* w \le \frac{\gamma+2}{4} |x|^{\gamma/2} |v|^2   -    \frac\beta2 |x|^{\beta + \gamma/2}  
$$
and thus 
\bean
\BB^\sharp w
\le  
 C_1  | x |^{\gamma/2} |v|^2   -    \frac12 |x|^{\beta + \gamma/2}  - \frac{1}2 |v|^{2+\gamma} .
\eean
Using that 
$$
C_1 | x |^{\gamma/2} \le  C_1^2 + \frac{1 }{ 4} | v |^{\gamma},
$$
if $|x| \le |v|$ and  $$
C_1|v|^2 \le  (4C_1)^{\beta/(\beta-2)} + \frac14|x|^\beta,
$$
if $|v| \le |x|$, we obtain
\bean
\BB^\sharp w 
\le  \big(C^2_1 +  (4C_1)^{\beta/(\beta-2)} \big) \HH   -    \frac14 |x|^{\beta + \gamma/2} -   \frac{1}4  | v|^{2+\gamma}, 
\eean
from what we conclude with $\alpha := \gamma/ (2\beta)$. 
\Black

\smallskip
\noindent
{\sl Step 2.}  We now assume $\beta <\gamma +2$ and we define 
$$
 w: =   1 + \frac12 [x]^{\beta/2} \cdot v +   \HH,
  $$
 so that again $w \sim \langle\HH\rangle$. We easily compute
 \bean
\BB^\sharp w \le  C_0|x|^{\beta/2 -1}  |v|^2    -    \frac12 |x|^{\tfrac32\beta -1} -   \frac{1}2| v|^{\gamma+2} . 
 \eean
Using Young's inequality similarly as in the step~1, we  get that 
$$
\BB^\sharp w \le  C \HH  -    c  |x|^{\tfrac32 \beta -1} -   c |v|^{2+\gamma}, 
$$
which in turn implies  \eqref{eq:lem:kinetic-rd-moment} with  $\alpha := \min(\gamma/2,1/2-1/\beta)$.  
\end{proof}

 We classically deduce the following resolvent estimate. 
 
\begin{lem}\label{lem:kinetic-rd-RBL2m}
For any weight function $m_0 := \omega(\HH)$, there exists a weight function $m_1 := \omega_1(\HH)$ such that $m_1/m_0 \to \infty$ as $\HH \to \infty$ and for any $\kappa > \kappa_\BB  > - K_0$ 
there holds 
\beqn\label{eq:lem:kinetic-rd-RBL2m0L2m1}
\RR_\BB : L^2(m_0) \to L^2(m_1).
\eeqn
 
\end{lem}

\begin{proof}[Proof of Lemma~\ref{lem:kinetic-rd-RBL2m}.]
We split the proof into two steps. 

\smallskip\noindent
{\sl Step 1.}   
  We fix $\kappa_\BB \in (-K_0,\kappa_0)$ and $m_0 := \omega_0(\HH)$ with $\omega_0$ a function as defined above. We observe that 
 $$
 \BB^\sharp m_0 = -K m_0 \le - m_0,
 $$
 so that 
 $$
\int (\kappa - \BB) f ( f m_0) = \int f^2 (\kappa - \BB^\sharp) m_0 \le 0,
$$
which means that $\kappa - \BB$ is a dissipative operator in $L^2_{m_0}$. We deduce that $\RR_\BB(\kappa) : L^2_{m_0} \to L^2_{m_0}$ for any $\kappa > - 1$. 
 
\smallskip\noindent
{\sl Step 2.}    We take 
$$
m := \omega(\HH) w, \quad \omega(\HH) := \omega_0(\HH)/\HH,
$$
where $w$ is defined as in Step~1 of Lemma~\ref{lem:kinetic-rd-moment} when $\gamma \le \beta-2$ and  as in Step~2 of Lemma~\ref{lem:kinetic-rd-moment} when   $\gamma > \beta-2$. 
In any cases $m \lesssim m_0$. On the other hand, Lemma~\ref{lem:kinetic-rd-moment} and $\TT^* \omega(\HH) = 0$ imply together that 
$$
\BB^\sharp m \le C_\alpha m - c_\alpha \omega(\HH) \HH^{1+\alpha}, \quad \alpha > 0.
$$
This a priori estimate implies $\RR_\BB(C_\alpha) : L^2_{m} \to L^2_{m_1}$, with $m_1 := m_0 \HH^\alpha$. 
 For $g \in L^2_{m_0}$ and $\kappa > \kappa_\BB := - K_0$, the function $f := \RR_\BB(\kappa) g \in L^2_{m_0} \subset L^2(m)$ also satisfies 
 $$
 (C_\alpha - \BB) f =  g + (C_\alpha - \kappa) f.
 $$
We deduce $\| f \|_{L^2(m_1)} \lesssim \| g \|_{L^2(m)}+  \| f \|_{L^2(m)} \lesssim \| f \|_{L^2(m_0)}$, which is nothing but \eqref{eq:lem:kinetic-rd-RBL2m0L2m1}. 
\end{proof}

We  argue as during the proof of \ref{H3} in  Section~\ref{subsec:kinetic-torus}. 
By a localization argument and the averaging lemma, we have  $\AA \RR_\BB(\kappa) : L^2_{m_0} \to L^2(B_R \times B_R)$ with compact injection for any $R>0$. 
Together with Lemma~\ref{lem:kinetic-rd-RBL2m}, we deduce that  $\RR_\BB(\kappa) \in \KK(L^2_{m_0})$ for any $\kappa > \kappa_\BB$, and we conclude exactly as in 
 Section~\ref{subsec:kinetic-torus}.
 
\medskip
{\bf Condition \ref{H4} and \ref{H5'}.}  We recall that  it has been proved in \cite[Lem.~4.5]{MR4063917} that the semigroup $S_t$ associated to the operator $\LL$ satisfies the Doblin-Harris condition: 
for any $T > 0$ and $\varrho > 0$, there exists $\alpha > 0$ such that 
\beqn\label{eq:KineticRd-Harris}
S_T f \ge \alpha {\bf 1}_{B_\varrho} \int_{B_{\vartheta\varrho}} f dxdv, \quad \forall \, f \ge 0, 
\eeqn
for some constant $\vartheta \in (0,1)$ and where $B_\rho := \{(x,v) \in \R^{2d}; \, |x| < \rho, \, |v| < \rho \}$. Although the statement of \cite[Lem.~4.5]{MR4063917} is not written in that way, one may easily track the constants appearing in 
Lemmas 3.5, 3.6, 3.7 and 4.1 in  \cite{MR4063917} and one obtain \eqref{eq:KineticRd-Harris} with $\vartheta := 1/2$. Now, \eqref{eq:KineticRd-Harris} immediately implies \ref{H5'} which in turn implies \ref{H4} thanks to 
Lemma~\ref{lem:Irred-S>0impliesR>0}-(2)-(3).

\medskip
 Because $\LL$ is the generator of a semigroup it also satisfies the weak maximum principle and Kato's inequalities \ref{H1'}. We are then in position to apply  Theorem~\ref{theo:exist1-KRexistence}, Theorem~\ref{theo:KRgeometry1}, Theorem~\ref{theo:lambda1largestBIS}   and Theorem~\ref{theo:ergodicity-compact-trajectories}-(3) and thus complete the proof of Theorem~\ref{theo:kineticRd-main}.

\bigskip 
\bigskip

%


\section{The kinetic Fokker-Planck equation }
\label{part:application5:KFPequation}


%
%
%
%
%
%
%
%
%

In this part, we consider the  kinetic Fokker-Planck evolution equation  associated to the operator 
\beqn\label{eq:FPK-defOperator}
\LL f  := - v \cdot \nabla_x  f + \Delta_v f + b \cdot \nabla_v f + cf,
\eeqn
on functions $f : \OO \to \R$, where $\OO := \Omega \times \R^d$, $\Omega \subset \R^d$ is a domain,  $b : \OO \to \R^d$ is a given vector field
and  $c : \OO \to \R$ is a given function. 
In contrast with the previous part, collisions are typically {modeled} by a Fokker-Planck operator
$\Delta_v f + \Div_v (vf)$ (when $b=v$ and $c = d$) which takes into account a  thermal bath of (Gaussian) whitenoise  instead of the integral collisional operator $\KKK [f] - Kf $ in the linear Boltzmann equation \eqref{eq:kinetic-eqLinearBoltzmann}. 

\medskip

We will consider the case when $\Omega$ is a bounded domain \Black and the equation is complemented with a boundary condition. More precisely, we assume the classical balance between the values of the trace $\gamma f$ of $f$ on the outgoing and incoming velocities subsets of the boundary  
\beqn\label{eq:FPK-BdaryCond}
(\gamma_- f )(x,{v}) =  \RRR_x( \gpf(x,.) ) ({v}) \,\, \hbox{ on }  \,\,   \Sigma_-, 
\eeqn
where in this context we define $\Sigma_\pm^x := \{ {v} \in \R^3; \pm \, {v} \cdot \nu_x > 0 \}$ the sets of outgoing ($\Sigma_+^x$) and incoming ($\Sigma_-^x$) velocities at point $x \in \partial\Omega$, next the sets
$$
 \Sigma_\pm = \{ (x,{v}) \in \Sigma; \pm \nu_x \cdot {v} > 0 \} = \{ (x,{v}); \, x \in \partial\Omega, \, {v} \in \Sigma^x_\pm \}, 
$$
and finally the outgoing and incoming trace functions $\gamma_\pm f := {\bf 1}_{\Sigma_\pm} \, \gamma f$. Here and in the sequel, $\nu_x$ denotes the unit normal outward vector field defined on the boundary set $\partial\Omega$.
We similarly define the grazing velocity set
$$
 \Sigma_0 = \{ (x,{v}) \in \Sigma;  \nu_x \cdot {v} = 0 \}.
$$ 
The reflection operator $\RRR_x$ is local in position, but can be local or nonlocal in the velocity variable, so that it writes 
$$
(\RRR_x g )(v) := \int_{\Sigma_+^x} r(x,v,v_*) g(v_*) v_* \cdot \nu_x \, dv_* ,
$$
for a reflection kernel $r : \partial\Omega \times \R^d \times \R^d \to \R$.  Some classical general assumptions on $r$ are  
\beqn\label{eq:FPK-galRRR}
r \ge 0, \quad  \RRR_x^* 1 = 1, \quad \RRR_x  \MMM = \MMM,
\eeqn
for some positive function $\MMM = \MMM(v)$, see for instance \cite{MR1313028,MR1307620,zbMATH00994612}. The second (normalisation) condition corresponds to the fact that all the particles reaching the outgoing boundary are put back on the incoming boundary (no mass is lost) while the third (reciprocity) condition means (when $\MMM$ is a Gaussian function) that the wall is in a local equilibrium state and is not influenced by the incoming particles.
The normalization condition implies the local mass conservation
\beqn\label{eq:FPK-boundaryL1conservation} 
\int_{\Sigma^x_-} \RRR_x g |\nu \cdot v| dv = \int_{\Sigma^x_+} g \nu \cdot v dv,
\eeqn
while the three assumptions \eqref{eq:FPK-galRRR} on $r$ together  also imply
\bean
\int_{\Sigma^x_-} \bigl(\RRR_x g  \bigr)^2  \MMM^{-1} |\nu \cdot v| dv 
&\le& \int_{\Sigma^x_-} \bigl(\RRR_x (g^2/\MMM) \bigr)  \bigl(\RRR_x \MMM \bigr) \MMM^{-1} |\nu \cdot v| dv 
\\
&=&  \int_{\Sigma^x_+} g_*^2 \MMM_*^{-1} \, (\RRR^*1)  \nu \cdot v_*dv_*,
\eean
where we have used the Cauchy-Schwarz inequality (and the fact that $r \ge 0$) in the first line and the   reciprocity condition in the second line.
As a consequence, we have  
\beqn\label{eq:FPK-boundaryL2estim} 
\int_{\Sigma^x_-} \bigl(\RRR_x g  \bigr)^2  \MMM^{-1} |\nu \cdot v| dv 
\le
 \int_{\Sigma^x_+} g^2 \MMM^{-1} \,    \nu \cdot v dv ,
\eeqn
where we have used the normalization condition in that last step.  
In the sequel, we will rather consider the   possibly position dependent Maxwell boundary condition operator 
\beqn\label{eq:FPK-boundary} 
\RRR_x g =   \alpha(x) \DD_x g  + \beta(x)  \Gamma_x g, 
\eeqn
where the accommodation coefficients 
$\alpha,\beta : \partial\Omega \to [0,1]$ satisfy $\alpha(x)+\beta(x) =: \zeta(x) \le 1$,  
$\Gamma_x$ is the specular reflection operator
\beqn
\label{eq:FPK-def_Gamma}
\Gamma_x (g (x,\cdot))(v)  =  g (x , \VV_x v), \quad \VV_x v = v - 2 \nu(x) (\nu(x) \cdot v),
\eeqn
and $\DD_x$ is the diffusive operator
\beqn
\label{eq:FPK-def_D}
\DD_x (g (x,\cdot))(v) = c_\MMM \MMM(v) \widetilde g (x), \quad \widetilde g (x) = \int_{\Sigma^x_+} g(x,w) \, \nu(x) \cdot w \, d w.
\eeqn
Here the constant  $c_\MMM := (2\pi)^{1/2}$ is such that  $c_\MMM \widetilde{\MMM} = 1$ and $\MMM$ stands for the standard Maxwellian 
\beqn
\label{eq:FPK-def_M}
\MMM(v) := (2\pi)^{-d/2} \exp(-|v|^2/2),
\eeqn
or, more generally, $ \MMM = \MMM(|v|) \ge 0$ is such that 
\beqn\label{eq:FPK-galMMM}
\DDD_x^* 1 = 1, \quad \DDD_x  \MMM = \MMM,
\quad
 \langle v \rangle^{\vartheta}  \MMM \in L^1(\R^d),
\eeqn
with $ \vartheta \ge 1$ (that last condition is necessary in order that the second relation above makes sense). 
%
The boundary condition \eqref{eq:FPK-boundary} corresponds to the \emph{pure specular reflection} boundary condition when  $\beta \equiv 1$
and it corresponds to the \emph{pure diffusive} boundary condition when $\alpha \equiv 1$. When $\zeta\equiv1$, the Maxwell boundary condition operator \eqref{eq:FPK-boundary} 
satisfies \eqref{eq:FPK-galRRR}.
On the contrary, when $\zeta \not\equiv1$, the $L^2$ estimate \eqref{eq:FPK-boundaryL2estim} holds but not anymore the mass conservation \eqref{eq:FPK-boundaryL1conservation}. However, the following $L^1$ estimate 
\beqn\label{eq:FPK-boundaryL1estim} 
\int_{\Sigma^x_-} |\RRR_x g| |\nu \cdot v| dv 
\le \zeta^* \int_{\Sigma^x_+} |g|   \,    \nu \cdot v dv 
\eeqn
holds, with $0 \le \sup\zeta \le \zeta^* \le 1$. Finally, the case $\zeta \equiv 0$ corresponds to the zero inflow problem. 

\smallskip
Let us finally mention that similarly as in Part~\ref{sec:Transport}, the regularity needed on the domain $\Omega$  may be  formulated in the following way: we assume that $\Omega$ is locally on one side of $\partial\Omega$ and there exists a function $\delta = \delta_\Omega \in W^{2,\infty}(\R^d)$ such that for all $x$ in an interior neighborhood of $\partial\Omega$ one has $\delta(x) =  \hbox{dist}(x,\partial\Omega)$ and the vector field $\nu$ 
defined on $\R^d$ by  $x \mapsto \nu(x) = \nu_x :=  - \nabla_x \delta(x)$  coincides with the previously defined unit outward normal vector field on $\partial\Omega$ and satisfies $\| \nu \|_{L^\infty} = 1$. 
We also assume that  the Lebesgue measure on $\partial\Omega$ is well defined and it is denoted by $d\sigma_{\! x}$.

 \medskip
\subsection{The trace problem}
\label{subsec:FPK-trace}
\

\smallskip
We consider in this section the trace problem for a solution $g = g(x,{v})$ to the stationary Vlasov-Fokker-Planck equation
\beqn\label{eq:FPK-trace1}
\MM g :=  {v} \cdot \nabla_x  g - b \cdot \nabla_{v} g  
-  \Delta_{v} g= G \ \hbox{ in }   \OO,
\eeqn
for a given a vector field $b=b(x,{v})$, a source term $G = G(x,{v})$ and for a solution $g = g(t,x,{v})$ to the  evolution Vlasov-Fokker-Planck equation
\beqn\label{eq:FPK-trace2}
  \partial_t g +   {v} \cdot \nabla_x  g - b \cdot \nabla_{v} g  
-  \Delta_{v} g= G \ \hbox{ in }   (0,T) \times \OO,
\eeqn
for a given a vector field $b=b(t,x,{v})$, a source term $G = G(t,x,{v})$. The second problem has been considered first in \cite{MR1634851} and next in \cite[Sec.~4]{MR2721875},
where a strong (renormalized) trace function is proved to exist.  In the sequel, we recall these results and slightly extending them by considering a possible $L^2H^{-1}$ source term.
We introduce some notations. We denote 
$$
d\xi := |\nu(x) \cdot v| dvd\sigma_{\!x} \ \hbox{ and } \
 d\xi^2 := (\nu(x) \cdot \hat v)^2 dvd\sigma_{\!x} 
$$
the measures on the boundary set $\Sigma$. We denote 
by $\BB_1$ the class of renormalized functions $\beta \in
\Wloc^{2,\infty}(\R)$ such that $\beta''$ has a compact support,  by $\BB_2$ the class of functions $\beta \in
\Wloc^{2,\infty}(\R)$ such that $\beta'' \in L^\infty(\R)$ and by $ \DD_0(\bar \OO)$ the space of test functions $\varphi \in \DD(\bar
\OO)$ such that $\varphi = 0$ on $\Sigma_0$. 
We finally define the dual operator 
\bean
\MM^* \varphi :=  - v \cdot \nabla_x \varphi + \Div_v (b\varphi) - \Delta_v \varphi.
\eean

 
\begin{theo}\label{theo-FokkerPlanck-renormalization}
We consider \Cyan $g \in L^{2}_{{\rm loc},x}H^1_{{\rm loc},v}$, $b \in L^{2}_{{\rm loc},x}H^1_{{\rm loc},v}$, $b, \Div_v b \in  L^{\infty}_{{\rm loc}}$, \Black  $G \in L^{2}_{{\rm loc},x}H^{-1}_{{\rm loc},v}$  and we assume that $g$ is a solution to the stationary  Vlasov-Fokker-Planck equation \eqref{eq:FPK-trace1}.  Then there exists $\gamma g \in \Lloc^2(\Sigma,d\xi^2)$ 
such that the following Green renormalized formula
\bear\label{eq:FPK-traceL2}
&&\int\!\!\!\int_{\OO} \bigl( \beta(g) \, \MM^* \varphi   - \beta''(g) \, |\nabla_{v} \, g |^2   \varphi) \, d{v} dx 
+ \langle G  , \beta'(g)   \varphi \rangle = \\ \nonumber
&&\qquad
= \int\!\!\!\int_{\Sigma} \beta(\gamma \, g) \, \varphi \,\,  \nu(x) \cdot {v} \,\,
d{v} d\sigma_{\! x}
\eear
holds for any renormalized function $\beta \in \BB_1$ and any test functions  $\varphi \in \DD(\bar \OO)$, 
as well as  for any  renormalized function  $\beta \in \BB_2$ and any test functions
$\varphi \in \DD_0(\bar \OO)$.   It is worth emphasizing that  $\beta'(g) \varphi \in L^2_xH^1_v$ so that the duality product 
 $\langle G  , \beta'(g)   \varphi \rangle$ is well defined.

If furthermore $\gamma_\mp g  \in \Lloc^2( \Sigma; d\xi)$ then $\gamma_\pm g  \in \Lloc^2( \Sigma; d\xi)$ and \eqref{eq:FPK-traceL2}
holds  for any renormalized function $\beta \in \BB_2$ and any test functions
$\varphi \in \DD(\bar \OO)$.  
 \end{theo}

\begin{proof}[Proof of Theorem~\ref{theo-FokkerPlanck-renormalization}.]
We only allude the proof which uses very similar arguments as those presented in Section~\ref{part:application4:Kequation} and that can also be partially found in \cite{MR972541,MR2721875}.
In the one hand, considering the mollifier $(\rho_\eps)_{\eps>0}$  defined in \eqref{eq:Transport-mollifier} with $z := (x,v)$, we get that $g_\eps$ is smooth and  
satisfies 
\bean
g_\eps \to g \ \hbox{ in }L^{2}_{{\rm loc},x}H^1_{{\rm loc},v}
, \quad \MM g_\eps = G_\eps \to G \ \hbox{ in } \ L^{2}_{{\rm loc},x}H^{-1}_{{\rm loc},v},
\eean
which is nothing but a variant of   \cite[Lem.~II.1]{MR1022305}.
The function $g_\eps$ being smooth, for any $\beta \in C^2$ such that $\beta' \in C^1_b$, we may differentiate $\beta(g_\eps)$ and we get
$$
\MM \beta(g_\eps) + \beta''(g_\eps) |\nabla_v g_\eps|^2 = \beta'(g_\eps)  G_\eps  
 \ \hbox{ in } \  \OO.
$$
 We may thus pass to the limit as $\eps\to0$ 
and we obtain \eqref{eq:FPK-traceL2}.  
\end{proof}

\Blue

We formulate a stability result that we will use several times in the sequel.

\begin{prop}\label{prop:KFP-stability}  
 Let us consider three  sequences $(g_k)$,  $(b_k)$, $(G_k)$ and three functions $g$, $b$, $G$ which all satisfy the required bounds of Theorem~\ref{theo-FokkerPlanck-renormalization}. 
 We assume that 
$$
v \cdot \nabla_x g_k - b_k \cdot \nabla_v g_k  - \Delta_v g_k =  G_k \quad \hbox{in} \quad \DD'(\OO)
$$
 for any $k \ge 1$ and that  $g_k  \wto g$ weakly  in $L^{2}_{{\rm loc},x}H^{1}_{{\rm loc},v}$, $b_k \to b$  strongly
  in $\Lloc^2(\bar\OO)$, $(b_k), (\Div_v b_k)$ bounded in $L^{\infty}_{{\rm loc}}$ and  $G_k \to G$ strongly in $L^{2}_{{\rm loc},x}H^{-1}_{{\rm loc},v}$. 
Then  $g$ satisfies \eqref{eq:FPK-trace1}, so that it admits a trace $ \gamma g \in \Lloc^2 (\Gamma;  d\xi^2)$  and $\gamma g_k \wto \gamma g$ weakly in  $\Lloc^2 (\Sigma;  d\xi^2)$. 
\end{prop}

\begin{proof}[Proof of Proposition~\ref{prop:KFP-stability}.]
Using \eqref{eq:FPK-traceL2} with $\beta(s) := s^2$ and $\varphi := v \cdot \nu(x) \chi$, $0 \le \chi \in \DD(\bar\OO)$, we obtain that $(\gamma g_k)$ is bounded in $\Lloc^2 (\Sigma;  d\xi^2)$, and thus, up to the extraction of a subsequence,  $\gamma g_k \wto \bar\gamma$ weakly in  $\Lloc^2 (\Sigma;  d\xi^2)$, for some $ \bar\gamma \in \Lloc^2 (\Sigma;  d\xi^2)$. We may pass to the limit in the weak formulation \eqref{eq:FPK-traceL2} written for $g_k$, $\beta(s) = s$ and $\varphi \in  \DD_0(\bar \OO)$. We immediately obtain that $g$ satisfies \eqref{eq:FPK-traceL2}  for $\beta(s) = s$ and any $\varphi \in  \DD_0(\bar \OO)$ and thus is a solution to the  stationary  Vlasov-Fokker-Planck equation \eqref{eq:FPK-trace1}. Thanks to Theorem~\ref{theo-FokkerPlanck-renormalization}, we deduce that $g$ admits a trace $\gamma g$, that $\bar\gamma = \gamma g$ and that the whole sequence $(\gamma g_k)$ weakly converges to $\gamma g$. 
\end{proof}

\Black

 \subsection{Well-posedness problem with inflow term at the boundary}
\label{subsec:FPK-WellPinflow}
We consider the kinetic Fokker-Planck operator $\LL$ defined in \eqref{eq:FPK-defOperator} and we start revisiting the well posedness problem 
\beqn\label{eq:FPK-wellpose1}
 (\lambda - \LL) f  = \mathfrak{F}\quad\hbox{in}\quad \OO, 
 \quad
\gamma_- f   =   \mathfrak{g}  \,\, \hbox{ on } \quad \Sigma_-, 
\end{equation}
for given data $\mathfrak{F} : \OO \to \R$ and $\mathfrak{g}: \Sigma_- \to \R$. 

\smallskip
For a given weight function $m : \R^d \to [1,\infty)$, we define the measure
$d\xi_m := m^2 |\nu(x) \cdot v| \, dv d\sigma_{\! x}$ on the boundary $\Sigma$. We next define $L^2H^{1}_m= L^2H^{1}_m(\OO)$  the space associated to the Hilbert norm defined by 
$$
\| f \|_{L^2H^{1}_m}^2 := \| f \|_{L^2_m}^2 +  \| \nabla_v f \|_{L^2_m}^2,
$$
and we assume that $m$ satisfies the Poincaré type inequality 
\beqn\label{eq:FPK-HypSurLePoidsm}
\| f \frac{\nabla m }{ m} \|_{L^2_m} \lesssim \| f \|_{L^2H^{1}_m}, \quad \forall \, f \in L^2H^1_m.
\eeqn
Such a Poincaré  inequality  is classically known to be true when $m := \MMM^{-\vartheta}$, $\MMM$ is the Maxwellian \eqref{eq:FPK-def_M} and $\vartheta > 0$. 
We also define 
$$
L^2H^{-1}_m := \{ \mathfrak{F}  = g + \Div_v G; \ g, G_i \in L^2_m(\OO) \},
$$
so that when $m=1$ the space $L^2H^{-1}_m$ is nothing but the space of continuous and linear mappings on $L^2H^{1}$. For $\mathfrak{F} \in L^2H^{-1}_m$ and $f \in L^2H^{1}_m$, we may 
thus write 
$$
\langle \mathfrak{F} , fm^2 \rangle \le \| \mathfrak{F}  \|_{L^2H^{-1}_m} \| f  \|_{L^2H^{1}_m}.
$$
We finally define in this context 
$$
W_2 := \{ f \in L^2H^{1}_m ; \  \hat v \cdot \nabla_x f  \in L^2H^{-1}_m \}, 
$$
and  
$$
W_{2,\Sigma} := \{ g \in W_2; \ \gamma g \in L^2(\Sigma; d\xi_m) \}, 
$$
with $W_{2,\Sigma} \not= W_2 $ in general. 

\begin{theo}\label{theo:FPK-WellPosedness1} Let us fix a vector field $b  \in \Hloc^1(\bar\OO)$, a function $c \in L^\infty(\OO)$, a weight function $m : \R^d \to [1,\infty)$ and let us 
assume that  $b/\langle v \rangle \in L^\infty(\OO)$, that \eqref{eq:FPK-HypSurLePoidsm} holds and that 
\beqn\label{eq:FPK-deflambda*&varpi}
\lambda^* :=  \hbox{\rm ess\, sup} \, \varpi   < \infty, \quad \varpi := c + \frac{\Delta m^2 }{ 2m^2} - \frac 12 \Div b - b \cdot \frac{\nabla m }{ m}.
\eeqn
For any $\mathfrak{F}  \in L^2H^{-1}_m$, $\mathfrak{g}  \in L^2(\Sigma_-;d\xi_m)$ and $\lambda > \lambda^*$, there exists a unique   solution 
$f \in W_{2,\Sigma}$ to the Dirichlet problem \eqref{eq:FPK-wellpose1}. We have furthermore  $f \ge 0$ if $\mathfrak{F}  \ge 0$ and $\mathfrak{g}   \ge 0$. 
\end{theo}

A similar result is proved in \cite[Appendix A]{MR875086}   in the case $\Omega = \R^d$.  Also observe that \eqref{eq:FPK-deflambda*&varpi} holds with $m := \MMM^{-1/2}$ when 
$\MMM$ is the standard Maxwellian \eqref{eq:FPK-def_M} and $b(v) = \vartheta v$, with $\vartheta > 1/2$.  

\begin{proof}[Proof of Theorem~\ref{theo:FPK-WellPosedness1}.]
We split the proof into five steps. 

\smallskip
{\sl Step 1. A priori estimates.}  We argue similarly as in \cite{CM-KFP**,CM-Boltz**}.
 Multiplying the first equation in \eqref{eq:FPK-wellpose1} by $fm^2$, performing several integrations by part in the velocity variable and using the Green formula, we have 
$$
 \int_\OO (\lambda - \varpi) f^2 m^2+ \frac12 \int_\Sigma (\gamma f)^2 m^2 \nu \cdot v      + \int_\OO |\nabla_v f|^2 m^2
  = \langle \mathfrak{F} , fm^2 \rangle.
 $$ 
Fixing $\lambda > \lambda^*$, using  the Young inequality
$$
\| \mathfrak{F}  \|_{L^2H_m^{-1}} \| f   \|_{L^2H_m^{1}} \le \Bigl( \frac{1 }{ 2(\lambda-\lambda^*)} +  \frac{1 }{ 2} \Bigr) \| \mathfrak{F}  \|_{L^2H^{-1}_m}^2 + \frac{\lambda-\lambda^* }{2} \| f \|_{L_m^2}^2 
+ \frac{1 }{ 2} \| \nabla_v f \|_{L_m^2}^2
 $$
and the boundary condition on the incoming set $\Sigma_-$  in \eqref{eq:FPK-wellpose1}, we deduce
\beqn\label{eq:FPK-WP-AprioriE1}
 (\lambda - \lambda^*)   \int_\OO   f^2 m^2+   \int_{\Sigma_+} (\gamma_+ f)^2 d\xi_m   +  \int_\OO |\nabla_v f|^2 m^2
 \le \frac{1 + \lambda-\lambda^*}{ \lambda-\lambda^*}   \| \mathfrak{F}  \|_{L^2H^{-1}_m}^2 +   \int_{\Sigma_-} \mathfrak{g}^2 d\xi_m.
\eeqn
\medskip
$\bullet$ Because of the first equation  in \eqref{eq:FPK-wellpose1} and the above estimate, we find 
\beqn\label{eq:FPK-dansW2}
\hat v \cdot \nabla_x f = \frac1{\langle v \rangle} \bigl( \mathfrak{F}  - \lambda f + \Delta_v f + b \cdot \nabla_v f  + cf \bigr)   \in L^2H^{-1}_m, 
\eeqn
so that $f \in W_2$. 
 
\medskip
$\bullet$ 
Multiplying the first equation in \eqref{eq:FPK-wellpose1} by $f \psi$, $\psi := \nu(x) \cdot \tilde v m^2$ where here and below we use the notations $\hat v :=  v/\langle v \rangle$, $\tilde v := v/\langle v \rangle^2$, $\langle v \rangle^2 := 1 + |v|^2$, and using the Green formula and one integration by part in the velocity variable, we get
\bean
\frac12\int_\Sigma (\gamma f)^2 (\nu \cdot \hat v)^2m^2  
&=&   \frac12\int_\OO   f^2  \hat v \cdot D_x \nu_x \hat v \, m^2 
- \int_\OO |\nabla_v f|^2  \psi    
\\
&&+  \int_\OO f \nabla_v   f ( b \psi - \nabla_v \psi )   
+  \int_\OO f ^2 \psi   (c- \lambda)   +  \langle \mathfrak{F} ,   f \psi \rangle.
\eean
Observing that 
$$
| \langle \mathfrak{F} ,   f \psi \rangle| \le \| \mathfrak{F}  \|_{L^2H^{-1}_m} \| f \nu(x) \cdot \tilde v \|_{L^2H^{1}_m} \lesssim 
\| \mathfrak{F}  \|_{L^2H^{-1}_m} \| f \|_{L^2H^{1}_m}
$$
and 
$$
\| f \nabla_v \psi \|_{L^2} \lesssim \| f \|_{L^2H^{1}_m},
$$
recalling that $b/\langle v \rangle \in L^\infty(\OO)$ and using the Cauchy-Schwarz inequality, we deduce 
\beqn\label{eq:FPK-tracemauvaistau}
\| \gamma f \|^2_{L^2(\Sigma;d\xi^2_m)}
\le C(1+|\lambda|)  \| f \|_{L^2H^1_m}^2 + C \| \mathfrak{F}  \|_{L^2H^{-1}_m} \| f \|_{L^2H^{1}_m},
\eeqn
for some constant $C = C(b,c,m,\nu)$, with $d\xi^2_m :=  (\nu \cdot \hat v)^2m^2  dvd\sigma_{\! x}$. 

\medskip
$\bullet$
For latter reference, we establish an estimate about the behaviour of the solution near the boundary. 
We now introduce the following Lions-Perthame \cite{MR1166050} type weight function
\beqn\label{eq:KFP-defpsibord}
\psi := 2 \delta(x)^{1/2} \nu(x) \cdot \tilde v,
\eeqn
and we observe that $\psi = 0$ on $\Sigma$, $\langle v \rangle \psi \in L^\infty(\OO)$, $ \nabla_v \psi \in L^\infty(\OO)$ and 
$$
v \cdot \nabla_x \psi = \frac1{\delta(x)^{1/2}} (\hat v \cdot \nu(x))^2 + 2  \delta(x)^{1/2} \hat v \cdot D_x \nu(x) \hat v.
$$
Multiplying the first equation in \eqref{eq:FPK-wellpose1} by $f \psi$,  we have 
$$
  \frac12\int_\OO v \cdot \nabla_x f^2  \psi  - \int_\OO f \, \frac{b }{ \langle v \rangle}  \cdot \nabla_v f  \langle v \rangle \psi  + \int_\OO \nabla_v (f \psi)  \cdot \nabla_v f
+  \int_\OO (\lambda- c) f^2 \psi    =  \langle   F , f\psi \rangle.
 $$ 
Using Cauchy-Schwarz and Young inequalities, we deduce 
\beqn\label{eq:FPK-ProcheDuBord}
 \int_\OO f^2  \frac{(\hat v \cdot \nu(x))^2}{\delta(x)^{1/2}} dvdx
 \le C  (1+|\lambda|) ( \| f \|^2_{L^2H^1} +   \| F \|_{L^2H^{-1}}^2), 
\eeqn
for some constant $C = C(b,c,n)$. 

\medskip
$\bullet$  We finally state a somehow classical regularity estimate when $\FF \in L^2_m(\OO)$. 
Taking advantage of the fact that $\FF \in L^2_m$ and $f \in L^2H^1_m$ and localizing the problem by introducing the function $g := f \chi_\eps \in L^2_xH^1_v(\R^d \times \R^d)$, $\chi_\eps \in C^2_c(\OO)$, ${\bf 1}_{\OO_\eps} \le \chi_\eps \le 1$, $\OO_\eps := \{ (x,v) \in \OO; \ \delta(x) > \eps, \ |v| \le 1/\eps \}$, we have 
$$
v \cdot \nabla_x g - \Delta_v g  + \langle v \rangle^2 g = \GG \  \hbox{in} \ \DD'(\R^d \times \R^d),
$$
with 
$$
\GG := (\FF - \lambda f - c f - b \cdot \nabla_v f) \chi_\eps - 2 \nabla_v f \cdot \nabla_v \chi_\eps + \langle v \rangle^2 f \chi_\eps \in L^2(\R^d \times \R^d).
$$
From the quantitative Hormander's hypoellipticity estimate of Hérau \& Pravda-Starov \cite[Prop.~2.1]{MR2786222}, we then have 
$$
\| D^{2/3}_x g \|_{L^2} + \| D^{2}_v g \|_{L^2} \lesssim  \|\GG \|_{L^2} + \| g \|_{L^2}.
$$
Coming back to the function $f$ and using the previous estimates, we deduce 
\beqn\label{eq:FPK-interieurReg}
\| D^{2/3}_x f \|_{L^2(\OO_\eps)} + \| D^{2}_v f \|_{L^2(\OO_\eps)} \le C ( \| \LL f \|_{L^2(\OO)} + \| f \|_{L^2(\OO)} ), 
\eeqn
for a constant $C = C(\lambda,\eps) >0$. 

\smallskip
{\sl Step 2. Existence. We assume $\mathfrak{g} = 0$.}
A possible way for proving the existence is to use Lions' variant of the Lax-Milgram theorem  \cite[Chap~III,  \textsection 1]{MR0153974} as in \cite{MR0252817,MR875086} and as we proceed now. Defining the bilinear form
$\EE : L^2H^1_m (\OO) \times C_c^1(\OO \cup \Sigma_-) \to \R$,  by 
\bean
\EE(f,\varphi) 
&=& \int_\OO  (\lambda-\LL) f \varphi m^2
\\
&:=& \int_\OO (\lambda f - b \cdot \nabla_v f - cf) \varphi m^2 +  \nabla_v f \cdot \nabla_v ( \varphi m^2) -   f (  v \cdot \nabla_x \varphi) m^2,
\eean
for any $f \in L^2H^1_m (\OO) $ and $\varphi \in C^1_c(\OO \cup \Sigma_-)$, we observe that this one is coercive, namely 
\bean
\EE(\varphi,\varphi) 
&=&  \int_\OO (\lambda - \varpi) \varphi^2 m^2  + \int_\OO |\nabla_v \varphi|^2 m^2 + \frac12 \int_{\Sigma_-} (\gamma_- \varphi)^2 d\xi_m
\\
&\ge& \kappa \| \varphi \|^2_{L^2H^1_m},
\eean
for any $\varphi \in C_c^1(\OO \cup \Sigma_-)$, with $\kappa := \min (\lambda-\lambda^*,1)>0$.
From the above mentioned Lions's theorem, for any $\mathfrak{F} \in L^2H^{-1}_m$, there exists $f \in L^2H^1_m$ such that 
\beqn\label{eq:FPK-formulationEE}
\EE(f,\varphi) = \langle \mathfrak{F} , \varphi m^2 \rangle, \quad \forall \, \varphi \in C^1_c(\OO \cup \Sigma_-). 
\eeqn
In particular, $f$ satisfies the first equation in \eqref{eq:FPK-wellpose1} in the distributional sense  $\DD'(\OO)$, and thus from \eqref{eq:FPK-dansW2}, we deduce that $f \in W_2$. 
Thanks to the trace Theorem~\ref{theo-FokkerPlanck-renormalization} and the estimate \eqref{eq:FPK-tracemauvaistau}, the function $f$ admits a trace $\gamma f \in L^2(\Sigma;d\xi^2_m)$.
Using the Green formula \eqref{eq:FPK-traceL2} with $\beta = {\rm id} \in \BB_1$, 
we have 
\beqn\label{eq:FPK-wellposeWeakEq}
\int\!\!\!\int_{\OO} \bigl( f (\LL^*-\lambda) \varphi + \mathfrak{F}    \varphi) \, d{v} dx  = \int\!\!\!\int_{\Sigma}  \gamma f \, \varphi \,\,  \nu(x) \cdot {v} \, d{v} d\sigma_{\! x},  
\eeqn
for any $\varphi \in \DD(\bar\OO)$.
Particularizing to $\varphi \in \DD(\OO \cup \Sigma_-)$ and comparing with \eqref{eq:FPK-formulationEE}, we deduce that  $\gamma_- f = 0$. 

\smallskip
{\sl Step 3. Existence. The general case $\mathfrak{g} \in L^2(\Sigma_-;d\xi_m)$.} When $\mathfrak{g} \in C^2_c(\Sigma_-)$, there exists a function $\mathfrak{h} \in C^2_c(\OO \cup \Sigma_-)$ such that $\mathfrak{h}_{|\Sigma_-} = \mathfrak{g}$
and we consider the source term $G := \mathfrak{F}  + (\LL -\lambda)\mathfrak{h} \in L^2H^{-1}_m$ as well as the problem 
$$
(\lambda-\LL)g = G \hbox{ in }\OO, \quad \gamma_- g = 0  \hbox{ on }\Sigma_-.
$$
From Step 2, there exists a solution $g \in W_{2,\Sigma}$ to this problem and we set $f := g +\mathfrak{h}$, in such a way that $f \in W_{2,\Sigma}$ and satisfies
\bean
\int_\OO f (\lambda - \LL^*) \varphi 
&=& 
\int_\OO g (\lambda - \LL^*) \varphi + \int_\OO \mathfrak{h} (\lambda - \LL^*) \varphi 
\\
&=& 
\int_\OO G \varphi + \int_\OO (\lambda - \LL)  \mathfrak{h} \varphi - \int_\Sigma \mathfrak{h}_{|\Sigma} \varphi \, \nu \cdot v,
\eean
and thus
\beqn\label{eq:FPK-wellposeWeakEq2}
\int_\OO f (\lambda - \LL^*) \varphi 
=\int_\OO \mathfrak{F}  \varphi  - \int_{\Sigma_-} \mathfrak{g} \varphi \, \nu \cdot v,
\eeqn
for any $\varphi \in C^2_c(\OO \cup \Sigma_-)$. Together with \eqref{eq:FPK-wellposeWeakEq}, we get that $\gamma_- f = \mathfrak{g}$ on $\Sigma_-$. 
In order to deal with the general case $\mathfrak{g} \in L^2(\Sigma_-;d\xi_m)$, we introduce a sequence $(\mathfrak{g}^n)$ of $C^2_c(\Sigma_-)$ such that $\mathfrak{g}^n \to \mathfrak{g}$ in $L^2(\Sigma_-,d\xi_m)$ and 
we next consider the associated sequence of solutions $(f^n)$ of $W_{2,\Sigma}$ just built above. 
Using the estimates exhibited in Step~1, we get that $(f^n)$ is a Cauchy sequence in $W_2$, 
so that it converges to a limit $f \in W_{2,\Sigma}$. 
We may pass to the limit in \eqref{eq:FPK-wellposeWeakEq2}
written for the sequence $(f^n)$ and deduce that the same equation holds at the limit for $f$.

\smallskip
{\sl Step 4. Uniqueness.} Consider two weak  solutions $f_i \in W_2$ to the equation \eqref{eq:FPK-wellpose1} in the sense that 
$$
\EE(f_i,\varphi) = \langle \mathfrak{F} , \varphi m^2 \rangle, \quad \forall \, \varphi \in C^1_c(\OO \cup \Sigma_-). 
$$
In particular, the difference $f := f_2-f_1 \in W_2$ satisfies 
$$
\EE(f,\varphi) = 0, \quad \forall \, \varphi \in C^1_c(\OO \cup \Sigma_-),
$$
and from the above discussion  $\gamma_- f = 0 \in L^2(\Sigma_-;d\xi_m)$.  Thanks to the trace Theorem~\ref{theo-FokkerPlanck-renormalization}, we deduce that $\gamma f \in \Lloc^2(\Sigma;d\xi_m)$
and we may choose $\beta(s) = s^2$ in the  Green formula \eqref{eq:FPK-traceL2}: we get  
\bean
\int_\OO f^2 \{v \cdot \nabla_x \varphi - \Div_v(b \varphi) + \Delta_v \varphi + 2f  (c-\lambda) \varphi  \} -  2 |\nabla_v f|^2 \varphi  =  \int_{\Sigma_+}  (\gamma f)^2 \nu \cdot v \varphi, 
\eean
for any test function $\varphi \in C^2_c(\bar \OO)$.
Choosing  $\varphi = m^2 \chi_\varrho$,  with $\chi_\varrho(v) := \chi(v/\varrho)$, $\chi \in C^2_c(\R^d)$, ${\bf 1}_{B_1} \le \chi \le {\bf 1}_{B_2}$, we deduce
\bean
 \int_\OO f^2 m^2 \bigl\{ (\lambda - \varpi) \chi_\varrho + \frac12 b \cdot \nabla \chi_\varrho - \frac{\nabla m }{ m}  \cdot   \nabla \chi_\varrho   - \Delta \chi_\varrho \bigr\} \le 0.
\eean
Because $f \in L^2H^1_m$, we may pass to the limit  $\varrho \to \infty$ thanks to the dominated convergence theorem and we  obtain 
\bean
 \int_\OO f^2 m^2   (\lambda - \varpi) \le 0, 
\eean
and thus $f = 0$. 

\smallskip 
{\sl Step 5. Positivity.} We assume now that $\mathfrak{F}  \ge 0$ and $\mathfrak{g} \ge 0$. We proceed similarly as in the previous step by considering  
$\beta(s) = s_-^2$, $\varphi = m^2 \chi_M$. Letting $M\to\infty$, we deduce 
\bean
 \int_\OO f_-^2 m^2   (\lambda - \varpi) \le 0,
\eean
and thus $f _-= 0$. 
 \end{proof}
 
 \medskip
 Summing up, gathering  the above estimates \eqref{eq:FPK-WP-AprioriE1}, \eqref{eq:FPK-dansW2}, \eqref{eq:FPK-tracemauvaistau}, \eqref{eq:FPK-ProcheDuBord}, \eqref{eq:FPK-interieurReg}, 
 we see that there exists a constant $C > 0$ such that any function $f \in D(\LL)$ satisfies
\bear\label{eq:FPK-wellposeBilanEstim}
&&\| f \|_{L^2H^1_m} + \|  \hat v \cdot \nabla_x f  \|_{L^2H^{-1}_m} +  \|  f   \frac{\hat v \cdot \nu}{\delta^{1/4}} \|_{L^2} 
\\ \nonumber
&&\qquad\qquad + \,  \| \gamma f \|_{L^2(\Sigma;d\xi^2_m)}  +  \| \gamma_+ f \|_{L^2(\Sigma;d\xi_m)} 
\le  C (\| f  \|_{L^2} + \| \LL f \|_{L^2})
\eear
and for any $\eps > 0$ there exists a constant $C_\eps$ such that  any function $f \in D(\LL)$ satisfies
$$
\| D^{2/3}_x f \|_{L^2(\OO_\eps)} + \| D^{2}_v f \|_{L^2(\OO_\eps)}  
\le  C_\eps (\| f  \|_{L^2} + \| \LL f \|_{L^2}). 
$$

\medskip
\subsection{Well-posedness problem with reflection condition at the boundary}
\label{subsec:FPK-WellPreflection}
We consider now the well posedness problem associated to the stationary equation
\beqn\label{eq:FPK-wellposeReflection}
  (\lambda - \LL) f  = \mathfrak{F} \quad\hbox{in}\quad \OO, 
  \quad  \gamma_- f   =   \RRR \gamma_+ f  \,\, \hbox{ on } \quad \Sigma_-, 
\end{equation}
for a given datum $\mathfrak{F}  : \OO \to \R$, where the kinetic Fokker-Planck operator $\LL$ is still defined by  \eqref{eq:FPK-defOperator} and the reflexion operator $\RRR$
is described in \eqref{eq:FPK-boundary}, \eqref{eq:FPK-def_Gamma}, \eqref{eq:FPK-def_D}.

\begin{theo}\label{theo:FPK-WellPosedness2} 
Let us fix a vector field $b  \in \Hloc^1(\bar\OO)$, \Cyan $b, \Div_v b \in  L^{\infty}_{{\rm loc}}(\bar\OO)$, \Black and a function $c \in L^\infty(\OO)$  which satisfy the assumptions of Theorem~\ref{theo:FPK-WellPosedness1}
with a  given weight function $m : \R^d \to [1,\infty)$ for the pure specular reflection case $\alpha \equiv 0$
and with the weight function $m := \MMM^{-1/2}$ when $\alpha\not\equiv0$, where $\MMM$ is the Gaussian function \eqref{eq:FPK-def_M} or a more general equilibrium satisfying \eqref{eq:FPK-galMMM}.
In that last case, we furthermore assume one of the two following conditions

\begin{itemize}

\item[(i)]  $1 - \zeta + \alpha^2/2 \ge \delta_* > 0 $, \Black and we observe that  $L^2(\Sigma; d \xi_m) \subset L^1(\Sigma;d\xi)$, 

\item[(ii)] $\langle v \rangle^{2} \MMM \in L^1$,  and we observe that $L^2(\Sigma;d \xi^2_m) \subset L^1(\Sigma;d\xi)$, 

\end{itemize}

where we recall that we have defined $d \xi_m := m^2|\nu(x) \cdot v| dvd\sigma_{\! x}$ and $d \xi^2_m := m^2 (\nu(x) \cdot \hat v)^2 dvd\sigma_{\! x}$. 

%
%
%
%
%

For any $\mathfrak{F}  \in L^2H^{-1}_m$ and $\lambda > \lambda^*$, there exists at least  one solution 
$f \in W_2$ to the Dirichlet problem \eqref{eq:FPK-wellposeReflection}.
Assuming furthermore that  
\beqn\label{eq:FPK-deflambda**}
 \lambda^{**}  := \hbox{\rm ess\,sup} \, (c -  \Div b )  < \infty, 
\eeqn
and $\lambda > \lambda^{**} $, the solution $f$ is unique and $f \ge 0$ if $\mathfrak{F} \ge 0$.

\end{theo}

It is worth emphasizing that  the assumptions of Theorem~\ref{theo:FPK-WellPosedness2} hold when $b = v$ and $m := \MMM^{-1/2}$.
We also emphasize on the fact that the additional assumptions (i) or (ii) are made in order to prove the uniqueness of the solution during the proof.


\begin{proof}[Proof of Theorem~\ref{theo:FPK-WellPosedness2}.]
We split the proof into four steps. 

\smallskip
{\sl Step 1. A priori estimates.} We multiply the first equation in \eqref{eq:FPK-wellposeReflection} by $fm^2$. As in Step~1 of the proof of Theorem~\ref{theo:FPK-WellPosedness1}, we get
$$
 \int_\OO (\lambda - \varpi) f^2 m^2+ \frac12 \int_\Sigma (\gamma f)^2 m^2 \nu \cdot v  + \int_\OO |\nabla_v f|^2 m^2
  = \langle \mathfrak{F} , fm^2 \rangle.
 $$ 
Using for instance \cite[Lem.~3.1]{MR4581432}, we have 
\beqn\label{eq:FPK-defEgammalpha}
 \int_\Sigma (\gamma f)^2 m^2 \nu \cdot v \ge  \int_{\Sigma_+} [ (1-\zeta) (\gamma_+ f)^2 + \alpha  (\DD^\perp \gamma_+ f)^2 ] d\xi_m =: \EE_{\zeta,\alpha}(\gamma_+ f) \ge 0,
\eeqn
with $\DD^\perp g := g - \DD g$. Using that the contribution of the boundary is nonnegative in the first estimate, we first deduce 
$$
 (\lambda - \lambda^*)   \| f \|_{L^2_m}^2 +  \| \nabla f \|_{L^2_m}^2 \le  \| \mathfrak{F}  \|_{L^2H^{-1}_m}  \| f \|_{L^2H^{1}_m}, 
$$
 for  $\lambda > \lambda^*$, so that 
$$
\min (\lambda-\lambda^*,1)  \| f \|_{L^2H^{1}_m} \le \| \mathfrak{F}  \|_{L^2H^{-1}_m}. 
$$
From the three above estimates together, for  $\lambda > \lambda^*$, we obtain 
\beqn\label{eq:FPK-ApriorWP2}
  \int_\OO   (\lambda - \varpi)_+  f^2 m^2+  \int_\OO |\nabla_v f|^2 m^2 + \frac12 \EE_{\zeta,\alpha}(\gamma_+ f)  
 \le \frac1{\min (\lambda-\lambda^*,1)} 
  \| \mathfrak{F}  \|_{L^2H^{-1}_m}^2 .
\eeqn 
There is no difficulty for also getting the pieces of information \eqref{eq:FPK-dansW2}, \eqref{eq:FPK-tracemauvaistau}, \eqref{eq:FPK-ProcheDuBord} and \eqref{eq:FPK-interieurReg}, 
so that in particular $f \in W_2$. 
 It is worth emphasizing here that when $\langle v \rangle^2 \MMM \in L^1$, we have $L^2(d\xi^2_m) \subset L^1(\Sigma;d\xi)$ by using the Cauchy-Schwarz and \eqref{eq:FPK-tracemauvaistau}, so that in particular
the boundary condition is well defined.

Let us show now how the last conclusion also holds under condition (i) in the statement of the Theorem. 
We then assume $\vartheta = 1$ in \eqref{eq:FPK-galMMM} and we show how to establish an additional a priori estimate. We indeed know from \eqref{eq:FPK-tracemauvaistau} that 
$$
\int_{\Sigma_-} (\alpha \DD (\gamma_+ f))^2 (\nu \cdot \hat v)^2m^2    dvd\sigma_{\!x}
\le \int_{\Sigma} (\gamma f)^2 (\nu \cdot \hat v)^2m^2  dvd\sigma_{\!x}
\le C_\lambda \| \mathfrak{F}  \|_{L^2H^{-1}_m} ^2,
$$
and similarly as in \cite{MR1301931} or \cite[proof of Lem.~2.2]{MR2721875} that 
$$
1= \int_{\Sigma_-^x}  |\nu(x) \cdot v| \MMM dv = C \int_{\Sigma_-^x} (\nu(x) \cdot \hat v)^2 \MMM dv, \quad \forall \, x \in \partial\Omega,   
$$
for some constant $C \in (0,\infty)$,
so that 
\beqn\label{eq:FPK-ApriorWP2-bord}
\int_{\Sigma_-} (\alpha \DD (\gamma_+ f))^2  d\xi_m 
= C \int_{\Sigma_-} (\alpha \DD (\gamma_+ f))^2 (\nu \cdot \hat v)^2m^2  
\le C C_\lambda \| \mathfrak{F}  \|_{L^2H^{-1}_m} ^2.
\eeqn
%
%
Summing up \eqref{eq:FPK-ApriorWP2} and \eqref{eq:FPK-ApriorWP2-bord}, and using that
$$
(\gamma_+ f)^2 \le 2 (\DD^\perp\gamma_+ f)^2 + 2 (\DD\gamma_+ f)^2,
$$
we deduce that
\beqn\label{eq:FPK-borne-gammalpha}  
\int_{\Sigma_+} [ 1-\zeta  + \alpha^2/2  ] (\gamma_+ f)^2  d\xi_m \le  C_\lambda \| \mathfrak{F}  \|_{L^2H^{-1}_m} ^2.
\eeqn
Defining
 $$
  f \in W_{2,\RRR} := \{ g \in W_2;  \ \gamma_-g = \RRR \gamma_+ g \}, 
$$
we  see that $W_{2,\RRR} = W_{2,\Sigma}$ if $1-\zeta + \alpha^2/2 \ge \delta_* > 0$,  
but it is worth emphasizing that we may have $W_{2,\RRR} \not= W_{2,\Sigma}$ in the general case.

\smallskip
{\sl Step 2. Existence when $\mathfrak{F}\ge0$.} 
With the help of Theorem~\ref{theo:FPK-WellPosedness1}, we define  $f_0 = 0$ and, recursively for any $n \ge 1$, we define  $f_n \in W_{2,\Sigma}$ as  the solution of 
\beqn\label{eq:FPK-existence2-scheme}
(\lambda - \LL) f_n = 
\mathfrak{F} \ \hbox{ in }\  \OO, \quad \gamma_- f_n   =   \RRR \gamma_+ f_{n-1}  \  \hbox{ on } \ \Sigma_-.
\eeqn
It is worth emphasizing here that $\gamma_+ f_{n-1} \in L ^2(\Sigma_+; d\xi_m)$ implies $\RRR(\gamma_+ f_{n-1}) \in L^2(\Sigma_-; d\xi_m)$ because of \eqref{eq:FPK-boundaryL2estim}. We also notice that $f_n \ge 0$ because $\mathfrak{F}\ge0$.
By linearity 
$$
(\lambda - \LL) (f_{n+1}-f_n) = 0 \ \hbox{ in }\  \OO, \quad \gamma_- (f_{n+1}-f_n)   =   \RRR \gamma_+ (f_{n}-f_{n-1})  \  \hbox{ on } \  \Sigma_-, 
$$
and we thus show recursively that $f_{n+1}-f_n \ge 0$. In other words, $(f_n)$ is an increasing sequence and thus also is $(\gamma f_n)$. 
From \eqref{eq:FPK-defEgammalpha}, 
we have 
\bean
 \int_\Sigma (\gamma f_n)^2 d\xi_m
&=&  \int_{\Sigma_{+}}  (\gamma_+ f_n)^2  d\xi_m  
-   \int_{\Sigma_{-}}  (\RRR \gamma_- f_{n-1})^2 d\xi_m
\\
&\ge&  \int_{\Sigma_{+}}  (\gamma_+ f_n)^2  d\xi_m  
-   \int_{\Sigma_{-}}  (\RRR \gamma_- f_{n})^2  d\xi_m
\ge   \EE_{\zeta,\alpha}(\gamma_+ f_n), 
\eean
so that the estimate \eqref{eq:FPK-ApriorWP2} holds true for $f_n$ (instead of $f$) uniformly in $n \ge 1$.  From the monotonous convergence theorem, there exists $f \in L^2H^1_m$ satisfying \eqref{eq:FPK-ApriorWP2}, \eqref{eq:FPK-borne-gammalpha}, \eqref{eq:FPK-tracemauvaistau} and such that $f_n \nearrow f$ a.e.
 Thanks to Proposition~\ref{prop:KFP-stability}, we have $\gamma f_n \nearrow \gamma f$ a.e. on $\Sigma$, from what we deduce that $\RRR \gamma_+ f_n \to \RRR \gamma_+ f$ in $L^2(\Sigma_-;d\xi^2_m)$ thanks to the monotonous convergence theorem.  As a consequence, 
we may pass to the limit in the weak formulation of \eqref{eq:FPK-existence2-scheme}, and we get that $f$ is a solution of \eqref{eq:FPK-wellposeReflection}. 
 We may also pass to the liminf in the estimate \eqref{eq:FPK-ApriorWP2} written for $f_n$, and we thus deduce that the same estimate holds for $f$.

\smallskip
{\sl Step 3. Existence when $\mathfrak{F} \in L^2H^{-1}_m$.} 
 When $\mathfrak{F} \in L^2_m$, we may introduce the splitting $\mathfrak{F} = \mathfrak{F}_+ - \mathfrak{F}_-$, just use the previous step for $\mathfrak{F}_\pm$
and conclude by linearity of the equation. 
When $\mathfrak{F} \notin L^2_m$, we proceed  similarly as in \cite{MR2721875} and in the following way. We first assume $\zeta \le \zeta^* \in [0,1)$ and we consider the mapping $\Psi : W_{2,\Sigma} \to W_{2,\Sigma} $,  $g  \mapsto f = \Psi(g)$, where $f$ is the solution to the stationary problem
\beqn\label{eq:FP5}
\left\{
\begin{aligned}
& (\lambda - \LL) f  = \mathfrak{F} \quad\hbox{in}\quad \OO \\
& \gamma_- f   =   \RRR \gamma_+g   \,\, \hbox{ on } \quad \Sigma_-.
\end{aligned}
\right.
\end{equation}
The space $W_{2,\Sigma}$ is endowed with the norm $\| \cdot \|_{W_{2,\Sigma}}$ defined by 
$$
\| g \|_{W_{2,\Sigma}}^2=  \| g \|_{L^2_m}^2 + \| \nabla_v g \|_{L^2_m}^2 + \| \gamma_+ g \|^2_{L^2_m(d\xi_1)}.
$$
From \eqref{eq:FPK-WP-AprioriE1} and the estimate $\| \RRR g \|_{L^2(\Sigma_-;d\xi_m)} \le \zeta^* \| g \|_{L^2(\Sigma_+;d\xi_m)} $ what we obtain by repeating 
the proof of \eqref{eq:FPK-boundaryL2estim}, we deduce 
\bean
\frac{1 }{ C_\lambda} \| f \|_{L^2H^1_m}^2 + \| \gamma_+  f \|^2_{L^2(\Sigma_+;d\xi_m)} 
&\le& C_\lambda  \| \mathfrak{F}  \|_{L^2H^{-1}_m}^2 + \| \RRR  \gamma_+  g \|_{L^2(\Sigma_-;d\xi_m)} 
\\
&\le& C_\lambda  \| \mathfrak{F}  \|_{L^2H^{-1}_m}^2 + \zeta^* \|    \gamma_+  g \|_{L^2(\Sigma_+;d\xi_m)} ,
\eean
for some constant $C_\lambda >0$. 
By linearity of \eqref{eq:FP5}, we deduce that for two functions $g_1,g_2 \in W_{2,\Sigma}$, and denoting $f_i := \Psi(g_i)$, we have  
\bean
\frac{1 }{ C_\lambda} \| f_2 - f_1 \|_{L^2H^1_m}^2 + \| \gamma_+  f_2 - \gamma_+ f_1 \|^2_{L^2(\Sigma_+;d\xi_m)} 
 \le  \zeta^* \|    \gamma_+  g_2 - \gamma_+  g_1 \ \|^2_{L^2(\Sigma_+;d\xi_m)},
\eean
so that $\Psi$ is a contraction in $W_{2,\Sigma}$. By the Banach fixed point theorem, we deduce that there exists a solution $f \in W_{2,\Sigma}$ to the equation \eqref{eq:FPK-wellposeReflection} in that case. 
Finally, in order to deal with the case  $\zeta^* = 1$, we consider a sequence $(\zeta^*_n)$ of $[0,1)$ such that $\zeta^*_n \nearrow 1$ and the associated sequence $(f_n)$ of
solutions in $W_{2,\Sigma}$ associated to the equation \eqref{eq:FPK-wellposeReflection} with the   modified reflection kernel $\RRR_n g := \zeta^*_n \RRR g$. 
From \eqref{eq:FPK-ApriorWP2} and  \eqref{eq:FPK-tracemauvaistau}, that sequence  satisfies 
$$
  \| f_n \|_{L^2H^1_m}^2 +  \| \gamma f_n \|_{L^2(\Sigma;d\xi^2_m)}^2  + \EE_{1,\alpha}(\gamma_+ f_n)  
 \le C_\lambda  \| \mathfrak{F}  \|_{L^2H^{-1}_m}^2.
$$
When $\alpha \not\equiv 0$, the above estimate or \eqref{eq:FPK-borne-gammalpha} also implies that $(\gamma_+ f _n)$ belongs to a weakly compact set of $L^1(\Sigma_+;d\xi)$. 
As a consequence, there exist $f \in W_2$ and $\bar\gamma_\pm$ two functions defined  on $\Sigma_\pm$ such that,  up to the extraction of a subsequence, 
\bean
&&f_n \wto f \  \ L^2H^1_m, \quad 
\gamma_\pm f_n \wto \bar\gamma_\pm  \  \ L^2(\Sigma_\pm;d\xi^2_m), 
\\
&& 
  \gamma_+ f_n \wto \bar\gamma_+  \  \ L^1(\Sigma_+;d \xi) ,  
\quad
 \RRR\gamma_+ f_n \wto   \RRR \bar\gamma_+  \  \ L^1(\Sigma_-;d\xi),
\eean
where we have used 
\eqref{eq:FPK-boundaryL1estim}  for the last convergence. 
Using Proposition~\ref{prop:KFP-stability}, we may thus pass to the limit in the equation \eqref{eq:FPK-wellposeReflection} satisfied by $f_n$ with   modified reflection kernel $\RRR_n$
\Cyan and we get that $f$ is a solution to the first equation in \eqref{eq:FPK-wellposeReflection}, that $\bar\gamma_\pm = \gamma_\pm f$ and then that $\gamma f$ satisfies the second equation  in \eqref{eq:FPK-wellposeReflection}. \Black
 In the pure specular reflection case $\alpha \equiv 0$, only the first line of convergences holds, but that it is enough in 
 order to pass to the limit in the equations (we refer to \cite{MR1776840,MR2721875} for similar arguments). 
%
\Black

\smallskip
{\sl Step 4. Other properties.} We further assume $\lambda > \lambda^{**}$. 
We proceed similarly as in \cite{MR2072842}.  Consider two weak  solutions $f_i \in W_2$ to the equation \eqref{eq:FPK-wellposeReflection}.
In particular, the difference $f := f_2-f_1 \in W_2$ satisfies 
$$
(\lambda - \LL) f  = 0 \ \hbox{ in }\  \OO, \quad
 \gamma_- f   =   \RRR \gamma_+ f  \  \hbox{ on } \  \Sigma_-.
$$
Using the Green renormalized formula \eqref{eq:FPK-traceL2},   we have 
\bean
0  
=
\int_\OO \beta'(f) (\lambda  - c)f \varphi  + \beta''(f) |\nabla f|^2 \varphi + \beta(f) (\Div_v(b \varphi) - v \cdot \nabla_x \varphi -\Delta_v \varphi ) + \int_{\Sigma} \beta(\gamma f) \nu \cdot v \varphi.
\eean
 for any  $\beta \in C^2(\R)$, $\beta' \in C^1_b(\R)$ and  any test function $\varphi \in C^2_c(\bar \OO)$. 
 We choose  $\varphi = \varphi(v) \ge 0$, $\beta \ge 0$ and $\beta'' \ge 0$, so that 
 \bean
0  
\ge
\int_\OO \beta'(f) (\lambda  - c)f \varphi   + \beta(f) (\Div_v(b \varphi)   -\Delta_v \varphi ) + \int_{\Sigma} \beta(\gamma f) \nu \cdot v \varphi.
\eean
By an approximation argument, we may now take $\beta(s) = |s|$, and we get 
 \bean
0  
\ge
\int_\OO |f| \bigl\{ (\lambda  - c) \varphi   +  (\Div_v(b \varphi)   -\Delta_v \varphi )\bigr\}  + \int_{\Sigma} |\gamma f| \nu \cdot v \varphi.
\eean
We observe that in any cases we have $f \in L^2_m(\OO) \subset L^1(\OO)$ and 
$\gamma f \in L^1(\Sigma;d\xi)$. By an approximation argument, we may now take $\varphi = 1$ and using the $L^1$ estimate \eqref{eq:FPK-boundaryL1estim} on $\RRR$ (with $\zeta^* = 1$), we get 
 \bean
0  &\ge& \int_{\Sigma_-}|\RRR \gamma_+f | |\nu \cdot v| - \int_{\Sigma_+}| \gamma_+f | |\nu \cdot v|  
\\
&\ge&
\int_\OO |f| \bigl\{ \lambda  - c + \Div_v b  \bigr\}   \ge (\lambda-\lambda_{**}) \int_\OO |f| .
\eean
We deduce that $f=0$. The proof of the positivity property follows the same arguments but choosing $\beta(s) = s_-$.
\end{proof}

%

For latter reference, we state the counterpart of the preceding result for the kinetic Fokker-Planck evolution equation.

\begin{theo}\label{theo:FPK-WellPosedness2evol} 
Let us make the same assumptions as in Theorem~\ref{theo:FPK-WellPosedness2}.
 For any $f_0 \in L^2_m$, there exists a unique solution $f \in C([0,T);L^2_m) \cap L^2(0,T;H^1_m)$ for any $T > 0$ to the  kinetic Fokker-Planck  evolution equation
\beqn\label{eq:FPK-evoleq}
\left\{
\begin{aligned}
&\partial_t f = \LL f   \quad\hbox{in}\quad (0,\infty) \times \OO \\
& \gamma_- f   =   \RRR \gamma_+ f  \,\, \hbox{ on } \quad (0,\infty) \times \Sigma_-, 
\end{aligned}
\right.
\end{equation}
with $\LL$ defined in \eqref{eq:FPK-defOperator} and $\RRR$ defined in \eqref{eq:FPK-boundary}.
 \end{theo}

The proof of Theorem~\ref{theo:FPK-WellPosedness2evol} is skipped since it is a mere adaptation of the proof of Theorem~\ref{prop:Transport-RKsolutions} and Theorem~\ref{theo:FPK-WellPosedness2}. We refer to  \cite[Cor. 2 7, Lem. 2.8 and Cor. 2.8]{MR4253803} where similar well-posedness results are established (see also \cite{MR2721875} for the existence part).

 \medskip
\subsection{The first eigenvalue  problem in a domain with reflection at the boundary}
\label{subsec:FPK-FirstEigenValue}
\

We consider now the first eigenvalue  problem for the kinetic Fokker-Planck operator \eqref{eq:FPK-defOperator} in a domain with reflection at the boundary, namely 
\beqn\label{eq:FPK-defOperator-1stEVP}
\left\{
\begin{aligned}
&\lambda f + v \cdot \nabla_x  f - \Delta_v f - b \cdot \nabla_v f - cf = 0 \quad\hbox{in}\quad \OO \\
& \gamma_- f   =   \RRR \gamma_+ f  \,\, \hbox{ on } \quad \Sigma_-,
\end{aligned}
\right.
\end{equation}
and the associated dual problem. 
In this section, we assume that  $b$ and $c$ satisfy the assumptions of Theorem~\ref{theo:FPK-WellPosedness1} with the weight function $m := \MMM^{-1/2}$ when $\alpha\not\equiv0$ and for a given weight function $m : \R^d \to [1,\infty)$
when $\alpha \equiv 0$ and $\RRR$ is given by \eqref{eq:FPK-boundary}. We additionally assume that
\beqn\label{eq:FPK-hypbsymvGal}
\liminf_{|(x,v)| \to \infty} \varpi (x,v)= - \infty, 
\eeqn
where we recall that $\varpi$ is defined in \eqref{eq:FPK-deflambda*&varpi}. When $\MMM$ is the Gaussian function, we find
$$
\varpi = c +\frac{|v|^2 + d }{ 2} - \frac 12 \Div b - b \cdot v, 
$$
so that \eqref{eq:FPK-hypbsymvGal} holds when  $b$ is typically a bounded perturbation of the vector field $b_0(v) = \vartheta_0 v$, $\vartheta_0 > 1/2$, 
and more precisely
$$
\Div_v b \in L^\infty(\OO) \quad\hbox{and}\quad
\inf_{x \in \Omega} \liminf_{|v| \to \infty} (b \cdot   v \langle v \rangle^{-2} ) \ge \vartheta_0 > 1/2.
$$
The above condition is quite technical but can be seen as a compatibility condition between the thermalization due to the boundary and to the Fokker-Planck collisional operator. We are then able to work 
in the functional space $X := L^2_m(\OO)$. 
\Black

\begin{theo}\label{theo:FPK-1stEVP} 
Under the above conditions,  the first eigentriplet problem associated to \eqref{eq:FPK-defOperator} has a unique solution $(\lambda_1,f_1,\phi_1) \in \R \times X \times X'$ with $f_1 > 0$ and $\phi_1 > 0$. 
\end{theo}

The proof of Theorem~\ref{theo:FPK-1stEVP} follows from Theorem~\ref{theo:exist1-KRexistence},   Theorem~\ref{theo:KRgeometry1}   and Theorem~\ref{theo:KRgeometry2}  as a consequence of conditions \ref{H1}--\ref{H5}. We prove now that each of these conditions is satisfied. 
  Theorem~\ref{theo:FPK-1stEVP} generalizes \cite[Thm.~2.12]{MR4347490} where the same problem is tackled for the zero inflow condition ($\alpha=\beta=0$) with $b = v  -F(x)$ and $c=1$ by using the classical Krein-Rutman theorem \cite{MR0027128} in the space $X=C_b(\bar\OO)$.  We also refer to \cite[Thm.~6.8]{MR4756950} for a variant and somehow generalisation of  \cite{MR4347490}.

\smallskip\smallskip
{\bf Condition \ref{H1}.}  From Theorem~\ref{theo:FPK-WellPosedness2}, the operator $\LL$ satisfies \ref{H1} with 
$$
\kappa_1 := \max (\lambda^*,\lambda^{**}), 
$$
with $\lambda^*$ defined by \eqref{eq:FPK-deflambda*&varpi} and $\lambda^{**}$ defined by \eqref{eq:FPK-deflambda**}. 
For later reference, let us state more precisely the available estimates for $f$. 
On the one hand, repeating the proof of Step~1 in the proof of Theorem~\ref{theo:FPK-WellPosedness2}, we establish that 
for any $\lambda > \kappa_1$ and $\mathfrak{F}   \in L^2_m$, 
the solution $f \in W_2$ to   the Dirichlet problem \eqref{eq:FPK-wellposeReflection} satisfies
\beqn\label{eq:FPK-ApriorWP2BIS}
  \int_\OO   (\lambda - \varpi)_+  f^2 m^2+  \int_\OO |\nabla_v f|^2 m^2 + \frac12 \EE_{\zeta,\alpha}(\gamma_+ f)  
 \le \frac1{ \lambda-\lambda^* }  
  \| \mathfrak{F}  \|_{L^2_m}^2. 
\eeqn
On the other hand, adapting the proof of \eqref{eq:FPK-ProcheDuBord}, we straightforwardly obtain 
\beqn\label{eq:FPK-ProcheDuBordBIS}
 \int_\OO f^2  \frac{(\hat v \cdot \nu(x))^2}{\delta(x)^{1/2}} dvdx
 \le C  \int  \mathfrak{F}^2 m^2, 
 \eeqn
for some constant $C = C(b,c,\nu,\lambda)$. For $\eps_x,\eps_v,\varrho > 0$, let us now define 
\beqn\label{eq:KFP-defUUH3}
\UU := \{(x,v) \in \OO; \ d(x,\partial\Omega) > \eps_x, \ |v| < \varrho \},
\eeqn
and compute 
\bean
\int_{\UU^c} f^2 m^2
\le
\int  f^2 m^2 {\bf 1}_{|v| \ge \varrho} + \int  f^2 m^2 {\bf 1}_{A_x} + \int  f^2 m^2 {\bf 1}_B, 
\eean
with
$$
A_x := \{ v \in B_\varrho, \  (\hat v \cdot \nu(x))^2 \le \eps_v^2 \},
\quad
B := \{ (x,v); \ |v| \le \varrho, \  (\hat v \cdot n)^2 \ge \eps_v^2, \ d(x,\partial\Omega) \le \eps_x \}.
$$
For the second term, we have 
\bean
\int  f^2 m^2 {\bf 1}_{A_x} 
&\le& \int  |A_x|^{2/r'} \| f(x,\cdot)\|_{L^r_v}^2 dx 
\\
&\lesssim&  (\varrho^{d-1} \eps_v)^{2/r'}\| f \|^2_{L^2H^1_m},
\eean
where we have used the Holder inequality with  $r \in (1,2^*/2)$ in the first line and the Sobolev inequality in the second line. 
For the third term, we have 
\bean
\int  f^2 m^2 {\bf 1}_{B} 
\le m^2(\varrho) \frac{\eps_x^{1/2}}{\eps_v^2}
\int_{\OO} f^2_n  \frac{(\hat v \cdot \nu(x))^2}{\delta(x)^{1/2}}. 
\eean
Gathering these last estimates with \eqref{eq:FPK-ApriorWP2BIS} and \eqref{eq:FPK-ProcheDuBordBIS}, we have established that the solution $f$ to  
equation \eqref{eq:FPK-wellposeReflection} furthermore satisfies 
\beqn\label{eq:KFP-estimUUc}
\int_{\UU^c} f^2 m^2
\le C \bigl( \frac{1 }{ \langle \varrho \rangle^2} + \varrho^{d-1}  \eps_v + m^2(\varrho) \frac{\eps_x^{1/2}}{\eps_v^2} \bigr) 
  \int  \mathfrak{F}^2 m^2, 
\eeqn
for a constant $C = C(b,c,\Omega,\lambda)$ and for any  $\eps_x,\eps_v,\varrho > 0$. 

\medskip
{\bf The strong maximum principle}. Let us now consider 
a function $0 \le f \in W_2 \backslash \{0 \}$ 
which satisfies  the Dirichlet problem \eqref{eq:FPK-wellposeReflection} associated to 
$\lambda > \kappa_1$ and a source term $0 \le  \mathfrak{F}   \in L^2_m \cap L^\infty$. 
In order to simplify the discussion, we assume that the normalization $\| f \|_{L^2_m} = 1$ holds. 
For proving the strong maximum principle, we briefly explain how we may 
adapt  the arguments we have presented for the diffusive equation in Part~\ref{sec:application1:diffusion} by taking in particular advantage 
of the above  established estimates,  the regularity results established in \cite{MR3923847,MR4453413} 
 and some spreading positivity results we learnt in \cite[Cor.~A.20]{MR2562709}. We proceed in three steps. 
 
\smallskip
{\sl Step 1.} On the one hand, from \eqref{eq:KFP-estimUUc}, we may choose conveniently $\varrho^{-1}, \eps_v, \eps_x  > 0$ small enough in such a way that 
\bean
\int_{\UU^c} f^2 m^2
\le \frac12\| f \|^2_{L^2_m}, 
\eean
where $\UU$ is defined by \eqref{eq:KFP-defUUH3}. 
Because  of the normalization condition, we have 
\beqn\label{eq:FPK-H4localization}
\int_{\UU} f^2 m^2 \ge \frac12  \| f \|_{L^2_m}^2
\eeqn
and consequently $f(x_0,v_0)^2 \ge  \delta^2_0 := \| f \|_{L^2_m}^2( 2\| {\bf 1}_\UU \|_{L^2_m}^2)^{-1}$ for at least one point $(x_0,v_0) \in \UU$.
 
\smallskip
{\sl Step 2.} On the other hand, let us recall some integrability and regularity results established in  \cite{MR3923847} for a solution $g$ to the kinetic Fokker-Planck evolution equation
$$
\partial_t g + v \cdot \nabla_x g = \Delta_v g + B \cdot \nabla_v g + s \ \hbox{ in }\  \VV, 
$$
or a sub-solution 
$$
\partial_t g + v \cdot \nabla_x g \le  \Delta_v g + B \cdot \nabla_v g + s \ \hbox{ in }\ \VV, 
$$
for some bounded set $\VV \subset (0,T) \times \OO$, $s \in L^2(\VV)$ and $B \in L^\infty(\VV)$. 
For that purpose, given some $(t^*,x^*,v^*)$, we define 
$$
Q_r := \{ (t,x,v); \ t \in (t^*-r^2,t^*], \ |x-x^* - (t-t^*)v^*| < r^3, \ |v-v^*| < r \}.
$$
We claim then that there exist $2 < p < q < \infty$, $\alpha \in (0,1)$ and for any  $0 < r_1 < r_0$ there exists $C$ such that 
\beqn\label{eq:KFP-L2LpBound}
\| g \|_{L^p(Q_{r_1})} \le C \, (\| g \|_{L^2(Q_{r_0})} + \| s \|_{L^2(Q_{r_0})} )
\eeqn
for any nonnegative subsolution $g$ on $Q_{r_0}$ from  \cite[Thm.~6]{MR3923847}, 
\beqn\label{eq:KFP-LqLinftyBound}
\| g \|_{L^\infty(Q_{r_1})}  \le C \, (\| g \|_{L^2(Q_{r_0})} + \| s \|_{L^q(Q_{r_0})})
\eeqn
for any nonnegative subsolution $g$ on $Q_{r_0}$ from  \cite[Thm.~12]{MR3923847} and
\beqn\label{eq:KFP-LinftyCalphaBound}
\| g \|_{C^\alpha(Q_{r_1})} \le C \, (\| g \|_{L^2(Q_{r_0})} + \| s \|_{L^\infty(Q_{r_0})})
\eeqn
for any  solution $g$ on $Q_{r_0}$ from  \cite[Thm.~3]{MR3923847}.  As a consequence of \eqref{eq:KFP-L2LpBound} and a classical covering argument, for any bounded set $\UU \subset \bar\UU \subset \OO$,
there exist $C_0 = C_0(\UU)$ and $C_1= C_1(\UU,\lambda)$ such that 
$$
\| f \|_{L^p(\UU)} \le C_0 \, (\| f \|_{L^2(\OO)} + \|  \mathfrak{F} +  c f - \lambda f  \|_{L^2(\OO)} ) \le C_1 ( \| f \|_{L^2(\OO)} +  \|  \mathfrak{F}   \|_{L^2(\OO)}). 
$$
 Observing that for $\varrho=p/2 >1$,  we have 
$$
 v \cdot \nabla_x  f^\varrho - \Delta_v f^\varrho - b \cdot \nabla_v f^\varrho + \varrho f^{\varrho-1}( \lambda f - cf - \mathfrak{F}) = - 4\frac{(\varrho-1)}{ \varrho} |\nabla (f^{\varrho/2})|^2 \le 0,
$$
so that $f^\varrho$ is a weak sub-solution to the kinetic Fokker-Planck equation, we may repeat the argument and obtain in that way that $f \in L^{p_k}(\UU)$ for any $k \ge 1$, with $p_k := \varrho^k 2$. 
Now, choosing $k$ such that $p_k \ge q$ and using \eqref{eq:KFP-LqLinftyBound} (as well as again a classical covering argument), we get 
$$
\| f \|_{L^\infty(\UU)} \lesssim  \| f \|_{L^2(\OO)} + \|  \mathfrak{F} +  c f - \lambda f \|_{L^q(\OO)}  \lesssim \| f \|_{L^2(\OO)} +  \|  \mathfrak{F}   \|_{L^q(\OO)} . 
$$
Using finally \eqref{eq:KFP-LinftyCalphaBound}, we deduce that there exists a constant $C = C(\UU,\lambda)$ such that 
$$
\| f \|_{C^\alpha(\UU)} \lesssim  \| f \|_{L^2(\OO)} +  \|  \mathfrak{F} \|_{L^\infty(\OO)}. 
$$
Together with the conclusion of the first step, we deduce that there exists a constructive constant $r_0 > 0$ such that $f   \ge   \delta_0 {\bf 1}_{B((x_0,v_0),r_0)}$. 

 \smallskip
{\sl Step 3.}  From \cite[Cor.~A.20]{MR2562709},  we deduce that for any bounded set $\UU \subset \bar\UU \subset \OO$, there exists a constructive constant $\delta = \delta(\delta_0,r_0,\UU) > 0$
such that 
 $$
 f(x,v) \ge \delta \quad\hbox{for any} \quad  (x,v) \in   \UU, 
 $$
 where it is worth emphasizing  that the hypothesis $b, c \in C(\OO)$ made in  \cite[Cor.~A.20]{MR2562709} is not really necessary and can be replaced by $b,c \in L^\infty(\UU)$. 
Because $\UU$ may be chosen arbitrary, we have established that $f > 0$ on $\OO$ and the strong maximum principle.

\smallskip
{\bf Condition \ref{H2}.}
 For a given function  $0 \le h_0 \in C^2_c(\OO)$ normalized by $\| h_0 \|_{L^2_m} = 1$, we define $f_0 \in D(\LL)$ as the solution to 
$$
(\kappa_1 - \LL) f_0 = h_0 \ \hbox{ in } \ \OO,
\quad 
\gamma_- f_0   =   \RRR \gamma_+ f_0 \  \hbox{ on } \  \Sigma_-. 
$$
Taking advantage of the fact that $h_0$ has compact support, we compute 
$$
 1=
\int_\OO h^2_0 m^2 =  \int_\OO (\kappa_1 - \LL) f_0 \, h_0 m^2 =  \int_\OO f_0  (\kappa_1 - \LL^*) (h_0 m^2) 
\le C_1  \| f_0  \|_{L^2_m},
$$
with $C_1 := \| m^{-1} (\kappa_1 - \LL^*) (h_0 m^2) \|_{L^2}$.
On the other hand, from \eqref{eq:FPK-wellposeBilanEstim}, we have 
\beqn \label{eq:FPK-wellposeBilanEstim-f0}
 \| f _0\|_{L^2H^1_m}  +  \|  f_0   \frac{\hat v \cdot \nu}{\delta^{1/4}} \|_{L^2} \le C_2,
\eeqn
for a constant $C_2$ only depending on $\| h_0 \|_{L^2_m}$,  $\kappa_1$ and the constant $C$ which appears in \eqref{eq:FPK-wellposeBilanEstim}. 
Arguing as in \eqref{eq:FPK-H4localization}, we deduce that 
\beqn\label{eq:FPK-H4localizationBIS}
\int_{\UU} f_0^2 m^2 \ge  (2C_1)^{-1}, \quad \supp h_0 \subset \UU, 
\eeqn
with $\UU = \UU_\varrho$ defined in \eqref{eq:KFP-defUUH3} and $\varrho > 0$ small enough (chosen constructively from $C_2$ and $C_1$). 
From the above constructive strong maximum principle, we deduce that $f_0 \ge \eps {\bf 1}_\UU \ge 1/C_0 h_0$ for some $\eps,C_0 > 0$. We conclude as in the second constructive argument for (H2) in Section~\ref{subsec:diffusion-domain}. Coming back indeed to the equation, we have 
$$
\LL f_0 =  \kappa_1 f_0 - h_0 \ge  \kappa_1 f_0 - \| h_0\|_{L^\infty} {\bf 1}_\UU \ge  (\kappa_1 -  \| h_0\|_{L^\infty}C_0) f_0,
$$
so that \ref{H2} holds with $\kappa_0 := \kappa_1 -  \| h_0\|_{L^\infty}C_0$ from Lemma~\ref{lem:Existe1-Spectre2bis}-{\bf (ii)}.


\smallskip\smallskip
{\bf Condition \ref{H3}.} 
 Let us fix $\kappa < \kappa_0$ arbitrary.  
We define $\BB f := \LL f - n \chi_R(v) f$ for any $f \in W_{2,\RRR}$, with $\chi_R \in \DD(\R^d)$ such that ${\bf 1}_{B_R} \le \chi_R \le {\bf 1}_{B_{2R}}$ and for some given $n, R \ge 0$ to be specified below. 
We observe that, at least formally, 
\bean
\int f m^2 (\BB - \kappa) f 
&=&  \int_\OO (\varpi - \kappa -  n \chi_R) f^2 m^2 - \frac12 \int_\Sigma (\gamma f)^2 m^2 \nu \cdot v  - \int_\OO |\nabla_v f|^2 m^2. 
\eean
Thanks to \eqref{eq:FPK-hypbsymvGal}, 
there exists a constant $R > 0$ such that 
$$
 \sup_{v \in \R^d \backslash B_{R}} \varpi \le \kappa. 
$$
Choosing $n := \sup \varpi_+ - \kappa $, 
we deduce that $\varpi - \kappa -  n \chi_R \le 0$. On the other hand, because of \eqref{eq:FPK-defEgammalpha}, the contribution of the boundary term in the above identity is non positive. 
We thus deduce that $(\BB-\kappa)$ is dissipative in $L^2_m$. 
We now establish that the associated operator  $\BB$ has compact resolvent. For $\mathfrak{F} \in L^2_m$, we consider $f \in L^2_m$ the solution to 
\beqn\label{eq:FPK-resolventBB}
 -\BB f  = \mathfrak{F} \ \hbox{ in }\ \OO, \quad \gamma_- f   =   \RRR \gamma_+ f  \  \hbox{ on } \  \Sigma_-,
\end{equation}
which existence follows from Theorem~\ref{theo:FPK-WellPosedness2}.
From the above discussion (with $\kappa = -1$) and the same arguments as in Step 1 of the proof of Theorem~\ref{theo:FPK-WellPosedness1}, we have
\beqn\label{eq:FPK-intriorEstimBIS}
  \int f^2 \langle  \varpi \rangle_- m^2 + 2  \int |\nabla_v f|^2 m^2 \le  \int  \mathfrak{F}^2 m^2. 
\eeqn
 Together with the regularity estimate \eqref{eq:FPK-interieurReg} and the compact imbedding $H^{2/3}(\UU) \subset L^2(\UU)$, 
we conclude that  $\BB$ has compact resolvent. The operator $\AA$ on $L^2_m$ defined by $\AA f := n \chi_R(v) f$ being bounded, 
 we may apply  Lemma~\ref{lem:H3abstract-StrongC}-(2) and we deduce that \ref{H3} holds for both the primal and the dual problems. 

\medskip
{\bf Condition \ref{H4}} is nothing but the yet established strong maximum principle. 

\medskip

\medskip
{\bf A variant of condition \ref{H5}.} Consider $(f,\lambda)$ a pair of eigenfunction and eigenvalue such that $\lambda \in \Sigma_{P+}(\LL)$.
Arguing similarly as in the proof of condition \ref{H5} in Section~\ref{subsec:diffusion-domain}, we know that 
$$
\widetilde \LL f = i \vartheta  f, \ \vartheta \in \R, \quad  \widetilde \LL |f| = 0
 $$
and introducing the real and complex part decomposition $f = g + ih$, we have
\bean
\int_\OO  \frac1{|f|^2} |g  \nabla_v h - h \nabla_v g|^2  = 0, 
 \eean
 and finally $g  \nabla_v h - h \nabla_v g = 0$ a.e. on $\OO$. Because of the regularity estimate presented during the above proof of the strong maximum principle,
the functions $f$  has H\"older regularity, and thus $g$ and $h$ are continuous on $\OO$. Because $|f| \not\equiv0$, we may claim that there 
exists a  point $(x_0,v_0) \in \OO$ such that  $h(x_0,v_0) > 0$ for instance. 
Denoting by $\omega$ the connected component of $\{ (x,v) \in \OO; \, h(x,v) > 0 \}$ containing $(x_0,v_0)$, we have  $\nabla (g/h) = 0$ on $\omega$,
and thus $g = \alpha(x) h$ on $\omega$ for some continuous function $\alpha:\Omega\to\R$.  Coming back to the eigenvalue equation that we may write in the following system form
$$
\widetilde \LL g = - \vartheta  h, \quad \widetilde \LL h  =  \vartheta g, 
 $$
 we compute 
 $$
- \vartheta  h =  \widetilde \LL (\alpha h) = \alpha \widetilde \LL h  - h v \cdot \nabla_x \alpha   =
\alpha \vartheta g  - h v \cdot \nabla_x \alpha \ \hbox{ on } \ \omega, 
 $$
 so that 
 $$
- \vartheta     =
\alpha^2 \vartheta    -   v \cdot \nabla_x \alpha  \ \hbox{ on } \ \omega. 
 $$
 We deduce that $\alpha$ is a constant on   $\omega$ and finally $\vartheta = 0$. We have thus established that $\lambda = \lambda_1$.

 \medskip 
 
  At this stage, we may use  Theorem~\ref{theo:exist1-KRexistence}, Theorem~\ref{theo:KRgeometry1} and Theorem~\ref{theo:KRgeometry2}, in order to get the conclusions 
   {\Blue  \ref{S1},  \ref{S2} and \ref{S32} }
  about the existence and uniqueness of the eigentriplet $(\lambda_1,f_1,\phi_1)$ which satisfies $f_1 > 0$, $\phi_1 > 0$, $\lambda_1$ is algebraically simple and   on the triviality of the boundary punctual spectrum. 
  
  \medskip
  We briefly explain how we may deduce the stability of $f_1$ by adapting some arguments developed in \cite{MR2162224} and already mentioned. 
  On the one hand, we know from \cite[Lem.~1.1]{MR2162224} that any solution $f$ to the rescaled evolution equation \eqref{eq:FPK-evoleq} with $\LL$ replaced by $\widetilde \LL = \LL - \lambda_1$ satisfies 
  $$
  \partial_t (H(X) f_1 \phi_1) + \hbox{div}_x (v H(X) f_1 \phi_1) - \hbox{div}_v (\phi_1^2 \nabla_v (H(X) f_1/\phi_1) )  =  - H''(X) f_1 \phi_1 |\nabla_v X|^2, 
$$
for any convex function $H: \R \to \R$ and with $X := f/f_1$. After integration, we get 
\beqn\label{eq:FPK-GRE}
\frac{d }{ dt} \int_\OO H(X) f_1 \phi_1  +  \int_\Sigma \nu \cdot v H(X) f_1 \phi_1 = - \int_\OO H''(X) f_1 \phi_1 |\nabla_v X|^2, 
\eeqn
When $H(s) := |s|$, the boundary term is 
\bean
\int_\Sigma  |\gamma f|  \gamma \phi_1 \nu \cdot v
&=&
\int_{\Sigma_+}  |\gamma_+ f|  \RRR^* \gamma_- \phi_1 \nu \cdot v
-\int_{\Sigma_-} | \RRR \gamma_+ f|    \gamma_- \phi_1  |\nu \cdot v|
\\
&\ge&
\int_{\Sigma_+}   |\gamma_+ f|  \RRR^* \gamma_- \phi_1  |\nu \cdot v|
-\int_{\Sigma_-}\RRR |\gamma_+ f|    \gamma_- \phi_1 \,  |\nu \cdot v|  = 0, 
\eean
from what we deduce the non expansive property  
\beqn\label{eq:FPK-L1phi1}
  \int_\OO |f_{t_1}| \phi_1  \le  \int_\OO |f_{t_0}| \phi_1, \quad \forall \, t_1 \ge t_0 \ge 0.
\eeqn
On the other hand, from the Cauchy-Schwarz inequality, we have 
$$
(\RRR \gamma_+f)^2 \le (\RRR \gamma_+f_1) \RRR (\gamma_+f^2/\gamma_+f_1) ,
$$
so that 
$$
\int_{\Sigma_-} \frac{(\RRR \gamma_+ f )^2 }{ \RRR \gamma_+ f_1} \gamma_- \phi_1 |\nu \cdot v| \le 
\int_{\Sigma_-} \RRR (\gamma_+ f^2/ \gamma_+ f_1) \gamma_- \phi_1 |\nu \cdot v| 
$$
and finally 
$$
\int_\Sigma (\gamma f)^2 (\gamma f_1)^{-1}   \gamma \phi_1 \,  \nu \cdot v \le 0. 
$$
When  $H(s) = s^2$, the equation \eqref{eq:FPK-GRE} and the last inequality imply 
\beqn\label{eq:FPK-L2phi1}
\frac{d}{dt} \int_\OO  f_1 \phi_1 (f/f_1)^2   +  2  \int_\OO f_1 \phi_1 |\nabla_v (f/f_1)|^2 \le 0.
\eeqn

We next recall a classical compactness result. 
\begin{lem}\label{eq:FPK-GREbis} Let $(g_n)$ be a sequence of functions such that 
$$
(g_n) \hbox{ is bounded in } L^\infty(0,T;L^2_{xv,\rm loc}) \cap  L^2(0,T;L^2_{x,\rm loc}H^1_{v,\rm loc})
$$
and 
$$
\partial_t g_n + v \cdot \nabla_x g - \Delta_v g_n = G_n \hbox{ bounded in }  \Lloc^2,
$$
then $(g_n)$ belongs to a strong compact set of $\Lloc^2$. 
\end{lem}

\begin{proof}[Proof of Lemma~\ref{eq:FPK-GREbis}.]
We just sketch it. Because 
$$
\partial_t g_n + v \cdot \nabla_x g=   \Delta_v g_n + G_n \hbox{ bounded in }  L^2_{tx}H^{-1}_v,
$$
the usual averaging lemma in \cite{MR923047,MR1003433} implies that 
$$
(g_n * \rho ) \hbox{ belongs to a strong compact set of } \Lloc^2,
$$
for any $\rho \in \DD(\R^d)$. 
On the other hand, introducing a mollifiers sequence 
 $(\rho_\eps)$ and writing then 
$$
g_n = (g_n - g_n * \rho_\eps) + g_n * \rho_\eps, 
$$
we see that the first term is small uniformly in $n$ as $\eps \to 0$ and the second term is relatively compact thanks to the first step, from what we immediately conclude.
\end{proof}

\medskip
Now, for $0 \le f_0 \in L^1_{\phi_1}$, we introduce the sequence $f_{0,k} := (f_0 \wedge k) {\bf 1}_{\UU_k} \in L^2(f_1^{-1} \phi_1) \cap L^2$, 
with $\UU_k := \{ (x,v) \in \OO; \, \delta(x) > 1/k, \ |v| \le k \}$,  and the associated solution 
$f_{k} \in L^\infty(0,T; L^2) \cap L^2(0,\infty; L^2_xH^1_v)$. 
Because of \eqref{eq:FPK-L2phi1}, for any increasing sequence $(t_n)$ which converges to $\infty$ and for any function $\varphi_m \in \DD(\OO)$, ${\bf 1}_{\UU_m} \le \varphi_m \le 1$, 
the rescaled and truncated function $g_n (t) :=  f_k(t+t_n) f_1^{-1} e^{-\lambda_1 (t+t_n)} \varphi_m$ meet the hypothesis of Lemma~\ref{eq:FPK-GREbis}, 
from what we classically deduce that the sequence of $\tilde f_n (t) :=  f_k(t+t_n) f_1^{-1} e^{-\lambda_1 (t+t_n)}$ is relatively strongly compact in $\Lloc^2$. 
Repeating the proof of  Theorem~\ref{theo:MeanErgodicityVariante1} and Theorem~\ref{theo:ergodicity-compact-trajectories} (see also \cite[Thm.~3.2]{MR2162224}), we deduce that
$\tilde f_n (t) \to  \langle  f_{0,k},\phi_1 \rangle f_1$ as $t\to\infty$.  Together with the above non expansive property \eqref{eq:FPK-L1phi1}, we deduce that $f_{t} \to \langle f_0,\phi_1 \rangle f_1$ in $L^1_{\phi_1}$ as $t\to\infty$.

 \smallskip
 We summarize our convergence result in the following theorem. 
 
\begin{theo}\label{theo:FPK-Cvgce} For any $f_0 \in L^2_m$, the holds $f_{t} \to \langle f_0,\phi_1 \rangle f_1$ in $L^1_{\phi_1}$ as $t\to\infty$.  
\end{theo}

Theorem~\ref{theo:FPK-Cvgce} generalizes \cite[Thm.~2.18]{MR4347490} for the zero inflow condition and  \cite[Thms.~1.6 \& 1.7]{MR4776290} for the torus case.
It is worth emphasizing that in these papers  the longtime convergence is established with exponential rate (with constructive estimate in  \cite{MR4776290}).
In  \cite{MR4347490} the proof is based on a representation formula for the associated semigroup $S$ which is proved to have a kernel $p_t \in (L^1 \cap L^\infty \cap C^\infty)(\OO)$ for any $t > 0$
(see \cite[Thms.~2.4 \& 2.6]{MR4347490} as well as  \cite{MR2248986,MR3382587,MR4412380}).
One then classically deduces that $S_t \in \KKK(X)$ for any $t > 0$ and $X = L^p$, $p \in [1,\infty]$, or $X = C_0$ (see \cite[Thm.~2.18]{MR4347490}), and next one may apply Theorem~\ref{theo:NagelWebb}. 
We also refer to  \cite[Thm.~6.8]{MR4756950},  \cite{MR3778533} and \cite{MR3237885,MR3902464,MR3897919} for related results. 

\medskip
 
We follow now a similar approach as in \cite{MR4347490,MR4756950}.
We start with a series of technical results. Here, we make the additional assumption
\beqn\label{eq:FPK-condvarpidiese}
 \varpi^\sharp (x,v):= \sup_{1 \le p \le \infty} w_p(x,v)  \le \kappa_2 <  \infty, 
\eeqn
with  
$$
\varpi_p :=  \frac{(2 - p) }{ p}\frac{\Delta m_p}{ m_p}  +
\frac{2 }{ p'}  \frac{|\nabla m_p|^2}{ m_p^2} + c - \frac{1}{ p} \, \frac{{\rm div}({bm_p^p}) }{ m_p^p},
$$
and $m_p := \MMM^{-1+1/p}$.

\begin{lem}\label{lem:KFP-splitABsmart} For any fixed $\kappa < \kappa_0$
there exists $\varrho_x > 0$, $\varrho_v > 0$ and $\kappa_2 \in \R$ such that defining $\AA f := \xi_{\varrho_v}(v) \zeta_{\varrho_x}(x) f$ with 
$\xi_{\rho_v} \in \DD(\R^d)$, ${\bf 1}_{|v| \le \rho_v} \le \xi_{\rho_v} \le {\bf 1}_{|v| \le 2 \rho_v}$, 
$\zeta_{\rho_x} \in \DD(\Omega)$, ${\bf 1}_{\delta(x) \ge \rho_x/2} \le \zeta_{\varrho_x} \le {\bf 1}_{\delta(x) \ge \varrho_x}$, and next  $\BB := \LL - \AA$, there 
hold
\bear\label{eq:lem:KFP-splitABsmart1}
\| S_\BB(t) \|_{\BBB(L^2_m)} &\lesssim& e^{\kappa t} , \quad \forall \, t \ge 0, 
\\
\label{eq:lem:KFP-splitABsmart2}
\| S_\BB(t) \|_{\BBB(L^p_{m_p})} &\lesssim& e^{\kappa_2 t},\quad \forall \, t \ge 0, \  \forall \, p \in (2,\infty].
\eear
\end{lem}

\begin{proof}[Proof of Lemma~\ref{lem:KFP-splitABsmart}.]
We first recall from Step 1 of the proof of Theorem~\ref{theo:FPK-WellPosedness1} and \eqref{eq:FPK-defEgammalpha} that 
\bean
(\LL f , f)_{L^2_m} 
&=& - \int |\nabla f |^2m^2 -  \frac12 \int_\Sigma (\gamma f)^2 m^2 \nu \cdot v  + \int f^2 \varpi m^2
\\
&\le& - \int |\nabla f |^2m^2   + \int f^2 \varpi m^2
\eean
and, with $\psi$ defined in  \eqref{eq:KFP-defpsibord}, 
\bean
(-\LL f , f)_{L^2_\psi} 
&=& - \frac12 \int f^2 (v\cdot \nabla_x \psi)    - \int f \, \frac{b }{ \langle v \rangle}  \cdot \nabla_v f  \langle v \rangle \psi  + \int \nabla_v (f \psi)  \cdot \nabla_v f
-  \int c f^2 \psi    
\\
&\le& -  \int f^2  \frac{(\hat v \cdot \nu(x))^2}{\delta(x)^{1/2}} dvdx
+ C   \int (f^2 + |\nabla f |^2).
\eean
Defining then $\tilde m := m - \beta \psi$, with $\beta > 0$ small enough, and summing up the two previous  estimates, we get 
\bean
(\LL f , f)_{L^2_{\tilde m}} 
&\le& -  \beta \int f^2  \frac{(\hat v \cdot \nu(x))^2}{\delta(x)^{1/2}}  
- \frac12 \int |\nabla f |^2m^2   + \int f^2( \varpi m^2 +1). 
\eean
Similarly as in \eqref{eq:KFP-defUUH3}, we define 
$$
\UU := \{(x,v) \in \OO; \ \delta(x) > \varrho_x, \ |v| < \varrho_v \},
$$
and we observe that 
$$
\UU^c \subset A \cup B \cup C,
$$
with
$$
A := \{ v \in B_{\varrho_v}, \  |\hat v \cdot \nu(x)| \le \eps_v \},
\quad
B := \{ v \in B_{\varrho_v}, \  |\hat v \cdot n| \ge \eps_v, \ \delta(x) \le \varrho_x \},
$$
for some $\eps_x > 0$,  and $C := B_{\rho_v}^c$. We next repeat the proof of \eqref{eq:KFP-estimUUc}, and we get 
\bean 
\int_{\UU^c} f^2 m^2
&\lesssim &
 (\varrho_v^{d-1} \eps_v)^{2/r'} \int |\nabla_v f |^2 + m(\varrho_v)^2 \frac{\varrho_x^{1/2}}{\eps_v^2}   \int f^2  \frac{(\hat v \cdot \nu(x))^2}{\delta(x)^{1/2}}  
 + \frac1{\varpi_-(\rho_v)} \int f^2 \varpi_- m^2.
 \eean
 Observing that 
 $$
 \int f^2( \varpi m^2 +1) \le \kappa \int f^2 \tilde m^2  + C_\kappa \int_\UU f^2 m^2 + C_\kappa \int_{\UU^c} f^2 m^2
 $$
 with $C_\kappa := \sup( \varpi +2 - \kappa)_+ < \infty$, and $\AA \ge C_\kappa {\bf 1}_\UU$ for $n := C_\kappa$, altogether, we conclude with  
 $$
 (\BB f, f)_{L^2_{\tilde m}} \le \kappa \| f \|_{L^2_{\tilde m}}.
 $$
 We then classically deduce that \eqref{eq:lem:KFP-splitABsmart1} holds. 
 
 \smallskip
 Similarly as for the first estimate and in the proof of \cite[Lem. 3.8]{MR3488535}, for  any smooth, rapidly decaying  and positive function $f$, we have 
 \bean
\int (\LL f) f^{p-1} m_p^p=   - \int_\Sigma \frac{(m_p\gamma f)^p }{ p} \nu \cdot v 
 - (p-1) \int   |\nabla (m_pf)|^2 \, (m_pf)^{p-2} d x,
   + \int f^p \varpi_p m_p^p.
\eean
From Daroz\`es-Guiraud (or Jensen) inequality, we know that the first (boundary) term is nonpositive (see \cite{DGuiraud} or \cite[Rem.~6.4]{MR2721875})
and we then classically  conclude to \eqref{eq:lem:KFP-splitABsmart2}.  
\end{proof}

 \begin{lem}\label{lem:KFP-LpCalpha}  There exists a finite family $2 = p_0 < p_1 < \dots   < p_k < \infty$  and $\alpha \in (0,1)$ such that   for both $\CC = \BB,\LL$, 
 for any $T> \tau > 0$ and $\VV \subset\subset \OO$, 
\bear
&&\label{eq:lem:KFP-LpLq}
\int_\tau^T  \| \AA S_\CC(t) f_0 \|_{L^{p_j}_{m_{p_j}}} \, dt \le C_{p_{j-1}}^{p_j} \| f_0 \|_{L^{p_{j-1}}_{m_{p_{j-1}}}} , \quad j=1 \dots, k,
\\
&&\label{eq:lem:KFP-LqLinfty}
\sup_{t \in [\tau,T]}   \| \AA S_\BB(t) f_0 \|_{L^\infty} \le C_{p_k}^{\infty}  \| f_0 \|_{L^{p_k}},
\\
&&\label{eq:lem:KFP-LinftyCalpha}
\sup_{t \in [\tau,T]}   \| S_\BB(t) f_0 \|_{C^\alpha(\VV)} \le C_\infty^\alpha \| f_0 \|_{L^{\infty}}.
\eear
 
\end{lem}

\begin{proof}[Proof of Lemma~\ref{lem:KFP-LpCalpha}.]
For $0 \le f_0 \in L^2_m$, let us denote $f :=  S_\BB f_0$ which thus satisfies the PDE
 $$
\partial_t f - \BB f  =  s:= cf  \ \hbox{ in } \ \DD'((0,T) \times \OO).
$$
Let us fix two open sets $U_i$ such that $[\tau,T]  \times \hbox{\rm supp} \xi \times \hbox{\rm supp} \zeta \subset U_0 \subset\subset U_1 \subset \subset (0,T) \times \OO$. From \cite[Thm.~6]{MR3923847} and a covering lemma, there exists a constant $\bar C_0 > 0$ 
and $p_1 > 2$ such that 
$$
\| f \|_{L^{p_1}(U_0)} \le \bar C_0 \bigl( \| f \|_{L^2(U_1)} + \| s \|_{L^2(U_1)} \bigr).
$$
The estimate \eqref{eq:lem:KFP-LpLq} for $j=1$ then follows from Theorem~\ref{theo:FPK-WellPosedness2evol} (and the classical underlying energy estimate). 
On the other hand,  \cite[Thm.~12]{MR3923847} similarly implies that there exists a constant $\bar C_k > 0$ and $p_k \in (p_1,\infty)$ such that 
$$
\| f \|_{L^{\infty}(U_0)} \le \bar C_k \bigl( \| f \|_{L^2(U_1)} + \| s \|_{L^{p_k}(U_1)} \bigr),
$$
and interpolating with the previous estimate, we get 
$$
\| f \|_{L^{p_j}(U_0)} \le \bar C_{j-1} \bigl( \| f \|_{L^2(U_1)} + \| s \|_{L^{p_{j-1}}(U_1)} \bigr),  \ \forall \, j, \  2 \le j \le k-1.
$$
The growth bound \eqref{eq:lem:KFP-splitABsmart2} and the two last estimates imply \eqref{eq:lem:KFP-LqLinfty} and \eqref{eq:lem:KFP-LpLq} for any $2 \le j \le k-1$.  
Finally,  \cite[Thm.~3]{MR3923847}  similarly implies that there exists a constant $\bar C_{k+1} > 0$ and $\alpha  \in (0,1)$ such that 
$$
\| f \|_{C^\alpha(U_0)} \le \bar C_{k+1} \bigl( \| f \|_{L^2(U_1)} + \| s \|_{L^{\infty}(U_1)} \bigr), 
$$
from what we deduce \eqref{eq:lem:KFP-LinftyCalpha} in the same way.
\end{proof}

%
%
%
%

\begin{theo}\label{th:KFP-CvgceRate}
%
Under the conditions of Theorem~\ref{theo:FPK-1stEVP} and the additional assumption \eqref{eq:FPK-condvarpidiese}, 
the conclusion \ref{E31} holds  in $L^2_m$ with 
non constructive rate. 
\end{theo}

\begin{proof}[Proof of Theorem~\ref{th:KFP-CvgceRate}.]
We introduce the splitting 
$$
\AA g := M \Upsilon_\eps g, \quad \Upsilon_\eps g := \chi_\eps g , \quad \BB := \LL - \AA, 
$$
with $\chi_\eps \in C_c^2(\OO)$,  ${\bf 1}_{\UU_{2\eps}} \le \chi_{\eps} \le {\bf 1}_{\UU_{\eps}}$ and $\UU_\eps := \{ |v| \le 1/\eps, \ \delta(x) > \eps \}$. 
We next write the  iterated Duhamel formulas (with $N := k+2$)
\bean
S_\LL 
=V + W * S_\LL,
\eean
with the usual notations \eqref{eq:KRexistTER-defvw} for  $V$ and $W$  associated to the integer  $N:=k+2$ and $k\ge1$ has been introduced in  Lemma~\ref{lem:KFP-LpCalpha}.
Next for $T > 0$ large, $\tau \in (0,T)$ small and two functions (of operators) $a$ and $b$, we define the modified convolution operator
$$
\left\{
\begin{aligned}
(a*_\tau b)(t)  &:=  \int_\tau^{t-\tau} a(t - s) b(s) \, ds \ \hbox{ if }\  t \in [\tau,T-\tau] \\
(a*_\tau b)(t) &:=  0 \ \hbox{ if } \ t \in  [\tau,T-\tau]^c,
\end{aligned}
\right.
$$
 (with these notations $*_0 = *$) and by induction $a^{*_\tau 1} := a$, $a^{*_\tau k} := a^{*_\tau (k-1)}*_\tau a$ for $k \ge 2$. With these notations, we define the new splitting 
$$
S_\LL = V + K^c_1 + K^c_2 + K,
$$
with 
$$
K := \ \Upsilon_{\!\nu}  W_{\tau} *_{\tau} S_\LL , \quad K^c_1 := W*S_\LL - W_{\tau} *_{\tau} S_\LL, \quad K^c_2 := (1-\Upsilon_{\!\nu})   W_{\tau} *_{\tau} S_\LL, 
$$
where $W_{\tau} := (S_\BB \AA)^{*_\tau N}$ and $\nu > 0$.
For later references, we also define 
recursively  $\Xi_0 := S_\LL$,  $\Xi_\ell := S_\BB \AA *_\tau \Xi_{\ell-1}$ for $\ell \ge 1$, so that $K = \Upsilon_\eps \Xi_N$. 
The sequel of the proof is split into two steps. 
\medskip

{\sl Step 1.} 
On the one hand, we 
compute 
\bean
 \| \Xi_{N} (T) f_0 \|_{L^{p_1}_{m_1}} 
 &\le& \| S_\BB \|_{L^\infty(\BBB(L^{p_1}_{m_1}))} \int_\tau^{T-\tau} \Bigl\| \int_\tau^{t-\tau} \AA S_\BB (t-s) \AA \Xi_{k-1}(s) ds  f_0 \Bigr\|_{L^{p_1}_{m_1}} dt
\\
&\le& C_T  \int_\tau^{T}  \int_\tau^T \| \AA S_\BB (t) \AA \Xi_{k-1}(s) f_0 \|_{L^{p_1}_m} dt ds 
\\
&\le& C_T C_2^{p_1} \int_\tau^T \| \AA \Xi_{k-1}(s) f_0\|_{L^2_{m_1}} ds ,
\eean
and thus
\beqn\label{eq:KFP-CvgceRate1}
 \|  \Xi_{N} (T) f_0 \|_{L^{p_1}_{m_1}} \le C_T \| f_0\|_{L^{2}_m},
\eeqn
where we have used \eqref{eq:lem:KFP-splitABsmart2} in the first line, the Fubini theorem in the second line, 
\eqref{eq:lem:KFP-LpLq} with $j=1$ in the third line and several times \eqref{eq:lem:KFP-splitABsmart1} in the last line. 

\smallskip 
For $\kappa < \kappa_0$, we may choose $\eps > 0$  small enough such that  \eqref{eq:lem:KFP-splitABsmart1} holds. 
From the very definition of $\AA$ and $S_\BB$, we may thus fix $\kappa_\BB \in (\kappa,\kappa_0)$ arbitrary and next $T> 0$ large enough such that $\| V(T) \|_{\BBB(L^2_m)} \le \tfrac13 e^{\kappa_\BB T}$. 
We may next use \eqref{eq:lem:KFP-splitABsmart1} and fix $\tau > 0$ small enough such that 
$$
\| K_1^c (T) \|_{\BBB(L^2_m)} \le \tau C_T \le \tfrac13 e^{\kappa_\BB T}.
$$
Last, because of \eqref{eq:KFP-CvgceRate1}, we may fix $\nu > 0$ small enough, in such a way that 
$$ 
\| K_2^c (T) f_0 \|_{L^2_m} \le \eta(\nu) \| \Xi_N(T) f_0 \|_{L^{p_1}_{m_1}}  \le \tfrac13 e^{\kappa_\BB T} \| f_0 \|_{L^{2}_{m}}.
$$
The three last estimates together, we have established 
\beqn\label{eq:KFP-CvgceRate2}
\| (V + K_1^c  + K^c_1)(T)  \|_{\BBB(L^2_m)} \le   e^{\kappa_\BB T} .
\eeqn

\medskip

{\sl Step 2.} Performing the same kind of computations as for proving \eqref{eq:KFP-CvgceRate1} and in particular using \eqref{eq:lem:KFP-LpLq}, we get 
\bean
\int_0^T \| \AA \Xi_j (s) f_0\|_{L^{p_{j+1}}_{m_{p_{j+1}}} }  ds
&\le& \int_0^{T}  \int_\tau^{T-\tau} \| \AA S_\BB (t) \AA\Xi_{j-1}(s) f_0 \|_{L^{p_{j+1}}_{m_{p_{j+1}}} }   dt ds 
\\
&\le&   C_{p_{j}}^{p_{j+1}}   \int_0^T \| \AA  \Xi_{j-1}(s) f_0\|_{L^{p_{j}}_{m_{p_j}} } ds,
\eean
for $j=1, \dots, k$, and with $p_{k+1} := \infty$. Iterating and using \eqref{eq:lem:KFP-LpLq} with $j=0$, we get 
$$
\int_0^T \| \AA \Xi_{k} (s) f_0\|_{L^{\infty}_m} ds
\lesssim  \|  f_0\|_{L^{2}_m}.
$$
Similarly, we may write 
\bean
\sup_{[\tau,T]} \| \AA \Xi_{k+1} f_0\|_{L^{\infty}_m} 
&\le& \sup_{t \in [\tau,T]}   \int_\tau^t \| \AA S_\BB (t-s) \AA\Xi_{k}(s) f_0\|_{L^{\infty}_m} ds 
\\
&\le& \sup_{t \in [\tau,T]} \| \AA S_\BB (s) \|_{\BBB(L^\infty_m)} \int_\tau^T \| \AA \Xi_{k}(s) f_0\|_{L^{\infty}_m} ds,
\eean
thanks to \eqref{eq:lem:KFP-LqLinfty}, and 
\bean
\| K f_0 \|_{C^\alpha(\OO)} 
&\le& \int_\tau^{T-\tau}  \| S_\BB (T-s) \AA  \Xi_{k+1}(s) f_0 \|_{C^\alpha(\UU_\nu)}  \, ds
\\
&\le& C_\infty^\alpha T \sup_{[\tau,T]} \|  \AA  \Xi_{k+1}  f_0\|_{L^{\infty}_m},
\eean
thanks to \eqref{eq:lem:KFP-LinftyCalpha}.
The three last estimates together and the compact support property  supp$\chi_\nu \subset\subset \OO$ imply
\bean
\| K f_0 \|_{C^\alpha \cap L^2_{m_{p_1}}} \lesssim  \|  f_0\|_{L^2_m}, \quad \forall \, f_0 \in L^2_m,
\eean
from what we deduce that $K \in \KKK(L^2_m)$. 
We may apply Theorem~\ref{theo:NagelWebb} in order to conclude.
\end{proof}

 %
%


  \bigskip 
\section{A mutation-selection model } 
\label{sec:appl:selec-mut}


%

%
%
%
%
%
%
 

In this section, we consider the mutation-selection evolution equation associated to the mutation-selection operator
\begin{equation}\label{eq:MUTS-main}
 \LL f :=J*f -W(x)f 
\end{equation}
defined on functions $f : \R^d \to \R$, where  $J$ is a the mutation kernel, $*$ stands for the convolution operator and $W$ is a  confining  potential.

 \medskip
\subsection{Almost regular mutation kernel}
\label{ssec:appl5:regular}

We assume that the mutation kernel $J$ is a positive finite measure of $\R^d$ which is lower bounded on a neighborhood of the origin, or in other words
\beqn\label{eq:mutselec-hypJ}
0 \le J \in M^1(\R^d), \quad  J \ge J_* {\bf 1}_{B_r}, 
\eeqn
for some constants $J_*,r>0$. 
We also assume that the selection potential $W : \R^d \to \R$ is continuous and satisfies 
\beqn\label{eq:mutselec-hypW}
W(x) > W(0) = 0, \ \forall \, x \not=0, \qquad
W(x)\rightarrow+\infty \ \hbox{as}  \ |x|\rightarrow\infty.
\eeqn
We finally assume  the following compatibility condition between mutation and selection:  there exist $\beta>0$ and a bounded Borel set $A\subset\R^d$ such that
\bear\label{eq:MUTS-assum1}
&& a := \essinf_{x\in A_\beta}{\int_{x-A_\beta}\frac{J(dz)}{W(x-z)}}>1, 
\\ \label{eq:MUTS-assum2}
&&J = J_1 + J_2, \quad J_1 \in C^1_c(\R^d), \quad  \kappa_*:=\|J_2\|_1:=\int_{\R^d}\!dJ_2 <\kappa_0 := (a-1) \beta,
\eear
where we use the notation $A_\beta=A\cap \{W\geq\beta\}$. 
 In the sequel, we work in the Banach lattice  $X := L^1(\R^d)$.

\begin{theo}\label{theo:mutselec-Main}
Under the above assumptions, we have 
\begin{enumerate}
    \item The first eigentriplet problem \eqref{eq:triplet1}-\eqref{eq:triplet2} admits a unique solution
$(\lambda_1, f_1, \phi_1) \in \R \times X_+ \times X^\prime_+$ with   the normalization $\|\phi_1\|=\langle\phi_1,f_1\rangle=1$,
and this triplet additionally satisfies  $\lambda_1 \ge \kappa_0$, $0 < f_1 \in L^1_{\langle W \rangle} (\R^d) \cap L^\infty_{\langle W \rangle} (\R^d)$ and $0 < \phi_1 \in L^1_{\langle W \rangle} (\R^d) \cap L^\infty_{\langle W \rangle} (\R^d)$. 
\item Moreover, $\LL$ generates a semigroup $S_\LL$ on $X$ and for any $f_0\in X$, there holds
\begin{equation}\label{eq=mutselec-CvgceToInfty}
    \|e^{-\lambda_1 t}S_\LL(t) f_0 -\langle \phi_1,f_0\rangle f_1\|_{L^1} \leq Ce^{-\alpha t}\|f_0-\langle \phi_1,f_0\rangle f_1\|_{L^1} ,
\end{equation}
for any $t \ge 0$ and for some constructive constants $C \ge 1$, $\alpha>0$.
\end{enumerate}
{\Blue In particular, the conclusions \ref{S1},  \ref{S2},  \ref{S33} and  \ref{E31} hold with constructive constants in~$L^1$.} \end{theo}

Let us comment on the above result. 

\begin{rem}\label{rem:selection-mutation}

\

\ (1) 
Assumption \eqref{eq:MUTS-assum1} is satisfies when $W$ is small enough in a neighborhood of the origin.
It is for instance satisfied if $W^{-1} \notin L^1(B_1)$.
That is in particular the case in dimension $d=1$ when $W$ is Lipschitz, because of the condition $W(0) = 0$. 
 
\ (2) Assume $J(x) = \eps^{-d} \rho(\eps^{-1} x)$ with $\rho \in C^1_c(\R^d) \cap \PP(\R^d)$ and $\rho \ge \rho_* {\bf 1}_{B_1}$, $\rho_* >0$, so that $J=J_1$ and $J_2 = 0$, and $W = W(|x|)$.
We may observe that for $\beta > 0$ and $\eps > 0$ small enough
$$
\inf_{\beta \le W(x) < 2\beta} \int_{\beta \le  W(y) < 2\beta} \frac{J(x-y)}{ W(y)} \, dy =: a \ge  \frac{\rho_*}{ 2\beta} \hbox{\rm meas} \{  \R_+^d \cap B_1 \}  > 1,
$$
so that \eqref{eq:MUTS-assum1} holds with $A := \{ W(x) < 2\beta \}$. 

\
(3) Assumption \eqref{eq:MUTS-assum1} is similar to \cite[Condition (2.3)]{Li2017}, see also \cite[Assumption 2.6]{MR4622852} and the comparison with \cite[Assumption 2.4]{MR4622852},
as well as~\cite[Condition~(3.7)-(3.8)]{MR923493} and~\cite[p. 250, {\it Note added in proof.}]{BurgerBomze}.
On the other hand, the conditions on $J$ are relaxed here since $J$ may have singular  part in~\eqref{eq:MUTS-assum2}.

\
(4) Optimal conditions linking $J$ and $W$ for the existence of a spectral gap are still unknown.
In the recent paper \cite{MR4622852}, using variational methods in a $L^2$ framework, the authors obtain a quantified spectral gap and the associated exponential stability 
when the mutation kernel $J$ is additionally assumed to be symmetric. 
Up to our knowledge, Theorem~\ref{theo:mutselec-Main} is the very first result providing a  quantified  spectral gap for a non-symmetric mutation kernel $J$.

\
(5)
Condition~\eqref{eq:MUTS-assum1} can be compared to the condition
\[
\bar a:=\esssup_{x\in\R^d}\int_{\R^d}\frac{J(x-y)}{W(y)}dy<1,
\]
under which no first eigenfunction may exist in $X$.
First, we claim that $\lambda_1\geq0$.
Indeed, considering $\eps>0$ and $f_\eps=\1_{B_\eps}$, we have
\[\LL f_\eps\geq-\big(\inf_{B_\eps}W\big) f_\eps,\]
so that the condition {\bf(H2)} holds for $\kappa_0=-\inf_{B_\eps}W$ for any $\eps>0$.
Since $W$ is continuous and $W(0)=0$, we deduce that $\lambda_1\geq0$ by passing to the limit $\eps\to0$.
Assume now by contradiction that there exists $f_1\in X_+ \setminus\{0\}$ such that
\begin{equation}\label{eq:MS:Lf=lambdaf}
\lambda_1f_1=\LL f_1 = J * f_1 - W f_1 
\end{equation}
and define, for any $\eps>0$, the function $\varphi_\eps(x)=\frac{1}{\eps+W(x)}\in L^\infty(\R^d)$.
Testing~\eqref{eq:MS:Lf=lambdaf} against $\varphi_\eps$ we get for any $\eps\in(0,1)$
\begin{align*}
0\leq\lambda_1\langle f_1,\varphi_1\rangle\leq\lambda_1\langle f_1,\varphi_\eps\rangle&=\iint\frac{J(x-y)}{\eps+W(x)}f_1(y)\,dx\,dy-\int\frac{W(x)}{\eps+W(x)}f_1(x)\,dx\\
&\leq\bar a\int f_1-\int\frac{W(x)}{\eps+W(x)}f_1(x)\,dx, 
\end{align*}
and passing to the limit $\eps\to0$ we obtain the contradiction $0\leq\lambda_1\langle f_1,\varphi_1\rangle\leq(\bar a-1)\int f_1<0$.
However, there always exists a principal eigenvector $f_1$ in $M^1(\R^d)$, which might have an atom at the origin when $\bar a<1$, see for instance~\cite{BurgerBomze}.
\end{rem}

The proof of Theorem~\ref{theo:mutselec-Main}  follows from Theorem~\ref{theo:exist1-KRexistence},   Theorem~\ref{theo:KRgeometry1}   and Theorem~\ref{theo:KRgeometry2}  as a consequence of conditions \ref{H1}--\ref{H5} that we establish now.
Setting $D(\LL):=L^1_{\langle W\rangle}(\R^d)$, we observe that $\LL$ is an unbounded closed operator with dense domain $D(\LL)$.

\smallskip
{\bf Condition \ref{H1} and \ref{H1'}.}
We define the semigroup
$$
S_{W}(t) f(x) := e^{-W(x) t}f(x), \quad \forall \, f \in L^p, \ p \in [1,\infty],
$$
which is clearly a positive semigroup of contractions. We next define $S_\LL$ as a bounded perturbation of $S_W$. It is also positive and it satisfies the growth estimate $\| S_\LL(t) \|_{\BBB(L^p)} \le e^{\|J \|_{1} t }$, where we recall that $\| J \|_1$ stands for the $L^1$ norm or the total variation norm of $J$.  We deduce that \ref{H1} holds true with $\kappa_1:=\|J\|_1$
thanks to Lemma~\ref{lem:Exist1-RkSG}-{\bf (i)}. 
Multiplying $\LL f$ by $\sign f$, for $f\in D(\LL)$, we immediately get Kato's inequality
\[(\sign f)\LL f = (\sign f) J*f - W|f| \leq J*|f| - W |f| = \LL |f|.\]

\smallskip
{\bf Condition \ref{H2}.}   
Let us define $f_0 :=\frac{1}{W(x)}\textbf{1}_{A_\beta}$, where $A_\beta$ is introduced in condition \eqref{eq:MUTS-assum1}.
 We compute
\bean
\LL f_0
&=& J\ast\Big(\textbf{1}_{A_\beta}\frac{1}{W}\Big)-\textbf{1}_{A_\beta} \ge \Big(J\ast\Big(\textbf{1}_{A_\beta}\frac{1}{W}\Big)-1\Big) \textbf{1}_{A_\beta}
\\
&\ge&   \Big( \essinf_{x\in A_\beta} \Big[ J\ast(\textbf{1}_{A_\beta}\frac{1}{W}) \Big] -1\Big) \textbf{1}_{A_\beta}
\\
&=&   (a-1) \textbf{1}_{A_\beta} \ge (a - 1) \frac{\beta  }{ W} \textbf{1}_{A_\beta}  = \kappa_0 f_0, 
\eean
where in the second equality we have used the very definition of $a$ in assumption~\eqref{eq:MUTS-assum1}.
We conclude that \ref{H2} holds thanks to Lemma \ref{lem:Existe1-Spectre2bis}-(ii).
  
\smallskip
{\bf Condition \ref{H3}.}  We introduce the splitting
\beqn\label{eq:mutselec-AAetBB}
\LL=\AA+\BB, \quad \AA f := J_1 \ast f, \quad \BB f := J_2* f - W(x) f. 
\eeqn
Arguing as in the proof of condition \ref{H1}, we see that $\BB$ is the generator a positive semigroup in $L^p(\R^d)$, $1 \le p \le \infty$, with growth bound $\omega(S_\BB) \le \kappa_*$ and thus $(\alpha-\BB)$
is invertible for any $\alpha \ge \kappa_0 > \kappa_*$, with 
\beqn\label{eq:mutselec-B-1}
\| (\alpha-\BB)^{-1} \|_{\BBB(L^p)} \le \frac{1}{ \alpha - \kappa_*}. 
\eeqn
 Next, observing that 
$$
 (W+\alpha) h = (\alpha - \BB) h + J_2 * h ,
$$
for any $h \in \DD(\LL)$ and $\alpha \ge \kappa_0$, we deduce that 
\beqn\label{eq:mutselec-B-1bis}
(W+\alpha) (\alpha - \BB)^{-1} g = g + J_2 * ((\alpha - \BB)^{-1} g),
\eeqn
for any $g \in X$ and $\alpha \ge \kappa_0$. Together with \eqref{eq:mutselec-B-1}, we deduce 
\beqn\label{eq:mutselec-B-1ter}
\|  (\alpha - \BB)^{-1} g \|_{L^p_W} \le  \| g \|_{L^p} + \|J_2 * ((\alpha - \BB)^{-1} g) \|_{L^p} \le \frac{\alpha }{ \alpha - \kappa_*} \| g \|_{L^p},
\eeqn
for any $g \in L^p$ and $\alpha \ge \kappa_0$. 
Defining $\WW(\alpha) := (\alpha - \BB)^{-1} \AA$,  we finally deduce from \eqref{eq:mutselec-B-1bis} the identity
$$
\WW(\alpha) f   = \frac{1 }{ W+\alpha} \AA f  + \frac{1}{ W+\alpha} J_2 * ((\alpha - \BB)^{-1} \AA f),
$$
for any $f \in X$ and $\alpha \ge \kappa_0$. We may then compute 
\bean
\| \WW(\alpha) f \|_{L^\infty}  
\le \frac{1 }{ \alpha} \| \AA f \|_{L^\infty} + \frac{1 }{ \alpha} \| J_2 \|_{1}  \| (\alpha - \BB)^{-1} \AA f \|_{L^\infty}, 
\eean
and together with \eqref{eq:mutselec-B-1} for $p=\infty$ and \eqref{eq:mutselec-B-1ter}, we deduce 
\beqn\label{eq:mutselec-WW}
\| \WW(\alpha) f \|_{L^\infty}  
\le \|J_1\|_\infty\frac{1}{\alpha-\kappa_*} \| f \|_{L^1} \Black,  
\eeqn
for any $f \in X$ and $\alpha \ge \kappa_0$. Starting from the same identity, we prove in a similar way
\beqn\label{eq:mutselec-WWbis}
\| \WW(\alpha) f \|_{L^\infty_W}  
\le \| J_1 \|_\infty   \frac{\alpha }{ \alpha-\kappa_*} \| f \|_{L^1} \Black,  
\eeqn
for any $f \in X$ and $\alpha \ge \kappa_0$. 
As a conclusion and gathering \eqref{eq:mutselec-B-1}, \eqref{eq:mutselec-B-1ter}, \eqref{eq:mutselec-WW}  and \eqref{eq:mutselec-WWbis},
we have established that 
\beqn\label{eq:mutselec-WWter}
\WW(\alpha) : L^1 \to L^1_{\langle W \rangle} \cap L^\infty_{\langle W \rangle},
\eeqn
with uniform bound for any $\alpha \ge \kappa_0$. 
Observing that $L^1_{\langle W \rangle} \cap L^\infty_{\langle W \rangle} \subset L^1$ is weakly compact and using Lemma~\ref{lem:H3Lp} with $p=1$, we deduce that \ref{H3} holds.
We can actually strengthen the compactness by noticing that $\AA : L^1 \to L^1_W \cap W^{1,1}$ is bounded because of assumption \eqref{eq:MUTS-assum2}.
This ensures that $\AA : L^1 \to L^1$ is compact, from what we deduce that $\WW(\alpha) : L^1 \to L^1$ is strongly compact for all $\alpha \ge \kappa_0$.
We may thus apply  Lemma~\ref{lem:H3abstract-StrongC}-(2) to infer that condition \ref{H3} holds for both the primal and the dual problems.

\smallskip

\smallskip
{\bf Condition \ref{H4}.} Assume that $\lambda \ge \lambda_1$ \Black  and $f \in  D(\LL) =  L^1_{\langle W \rangle}  $ satisfy
\beqn\label{eq:mutselec-hyp-strongMP}
\| f  \|_{L^1} = 1, \quad f \ge 0, \quad (\lambda-\LL) f \ge 0.
\eeqn
Denoting $W_R :=\inf_{B_R^c} W$, we compute
 $$
 \int_{B_R}f 
 \geq \int_{\R^d} f-\frac{1}{W_R}\int_{B_R^c}f W
 \geq 1  - \frac{1}{W_R} \| f \|_{ L^1_{\langle W \rangle} } \ge 1/2, 
 $$
 for $R > 0$ large enough by taking advantage of the fact that $W(x)$ tend to infinity when $|x|\rightarrow\infty$.
In particular, there exists  $x_{0}^f\in B_R$  such that
$$
\int_{B_{r/2}(x_0^f)}f\geq \delta := \frac12 \Big(\frac{r}{2R}\Big)^d>0,
$$
where we recall that $r$ is defined in~\eqref{eq:mutselec-hypJ}.
We deduce that  
$$
(J\ast f)(x)\geq J_* \int_{B_{r/2}(x_0^f)}f(y)dy{\bf 1}_{B_{r/2}(x_{0}^f)}(x)\geq J_*\delta {\bf 1}_{B_{r/2}(x_{0}^f)}(x). 
$$
Using the equation \eqref{eq:mutselec-hyp-strongMP}, we obtain 
$$
f(x) \geq \dfrac{(J\ast f)(x)}{W(x)+\lambda} \geq \dfrac{ J_*\delta }{ W[R]+\lambda } {\bf 1}_{B_{r/2}(x_0^f)}(x),
$$
for $W[R]=\sup_{B_{R}}W$. With that last information and  \eqref{eq:mutselec-hypJ} again, we have now 
$$
J\ast f  \geq \frac{J_*}{2^d}\frac{J_*\delta }{W[R]+\lambda}  \, {\bf 1}_{B_{r}(x_0^f)} , 
$$
and, iterating the argument, we deduce 
$$
f 
\geq \frac{J_*^m}{2^{(m-1)d}(W[R]+\lambda)^{m-1}} \delta {\bf 1}_{B_{mr/2}(x_0^f)}
\geq \bar{\gamma}{\bf 1}_{B_R}, 
$$
with $\bar{\gamma}=\bar\gamma(R)>0$ for $m=m(R)$ large enough. 
Choosing $R$ an integer, we have proved that
\begin{equation}\label{eq:muta_pos}
        f \ge h_0 := \bar\gamma(R)\1_{B_R}+\sum_{n\geq R}\bar\gamma(n+1)\1_{B_{n+1}\setminus B_n}>0.
\end{equation}
That means that the \ref{H4} holds, with constructive lower bound.

 \medskip

{\bf Condition \ref{H5}.}
Let us consider $f \in L^1_{\langle W \rangle} \backslash \{0 \}$ and $\lambda \in \C$ such that \eqref{eq:lemStrongKato&StrongPositivity3} holds, in particular 
\begin{equation}\label{eq:inverseKato-selectionmutation}
\LL |f|  = (\Re e \lambda) |f| \quad\hbox{and}\quad
\LL |f|  =  \Re e (\hbox{\rm sign} f) \LL f. 
\end{equation}
The first equality means that $\Re e \lambda$ is an eigenvalue associated to a positive eigenfunction, and Lemma~\ref{lem:PositiveEigenvector} then enforces $\Re e \lambda=\lambda_1$.
Lemma~\ref{lem:Uniquenessf1} subsequently ensures that $|f|\in(\Span f_1)_+\setminus\{0\}$, and in particular $|f|>0$.
Throwing away the term $W|f|$ in each side of the second identity in~\eqref{eq:inverseKato-selectionmutation}, we have
$$
\Re e \frac{\bar f }{ |f|} \, (J*f) = J * |f|.
$$
Integrating this equation, we get 
$$
 \int_{\R^{2d}} J(x-y) \Re e \Big[ |f(y)| - \frac{\bar f (x) }{ |f(x)|} f(y)\Big] \, dy = 0.
$$
From the positivity condition \eqref{eq:mutselec-hypJ} on $J$, we deduce 
$$
 |f(y)| - \frac{\bar f (x) }{ |f(x)|} f(y) = \Re e \Big[ |f(y)| - \frac{\bar f (x)}{ |f(x)|} f(y)\Big] = 0, \quad \forall \, x, y \in \R^d, \ |x-y| < r, 
$$
and thus $\bar f(x)/|f(x)| = \bar u$ for any $x \in \R^d$ for a constant $u \in \C$. That ends the proof of the  reverse Kato's inequality \ref{H5}.

\medskip

\begin{proof}[Proof of theorem~\ref{theo:mutselec-Main} part (1)]
    We may use Theorem~\ref{theo:exist1-KRexistence} in order to establish the existence of a solution $(\lambda_1,f_1,\phi_1)\in (0,+\infty)\times L^1\times L^\infty$ to the first eigentriplet problem \eqref{eq:triplet1}-\eqref{eq:triplet2}.
    From Theorem~\ref{theo:KRgeometry1} and  Theorem~\ref{theo:KRgeometry2}, this solution is unique, $f_1>0$, $\phi_1>0$, $\lambda_1$ is algebraically simple for both $\LL$ and $\LL^*$ and it is the unique eigenvalue in $\Sigma_+(\LL)$. 
    
 Due to~\eqref{eq:mutselec-WWter}, we actually have $f_1\in L^1_{\langle W \rangle} \cap L^\infty_{\langle W \rangle}$.
Observing that $\LL^*$ is of the same type as $\LL$,
$$
\LL^* \phi = \check J * \phi - W \phi, \quad \check J (x) := J(-x), 
$$
and considering the dual problem as a primal problem in $L^1$, Theorem~\ref{theo:exist1-KRexistence} also provides the existence of $\lambda^*_1>0$ and $0< \phi^*_1\in L^1_{\langle W \rangle} \cap L^\infty_{\langle W \rangle}$ such that
\[
\LL^*\phi^*_1=\lambda^*_1\phi^*_1.
\]
{\Cyan Because of Theorem~\ref{theo:KRgeometry1}, 
 we have in fact $\lambda_1^*=\lambda_1$ and the simplicity of $\lambda_1$ then yields that $\Span \phi_1^*=\Span \phi_1$. This ensures that $\phi_1\in L^1_{\langle W \rangle} \cap L^\infty_{\langle W \rangle}$ and also that $\phi_1$ enjoys the explicit lower bound~\eqref{eq:muta_pos}.}
Besides, we can prove 
\[
\|\phi_1\|_{L^\infty_{\langle W\rangle}} \leq \|J_1\|_{L^1}\frac{\lambda_1}{\lambda_1-\kappa_*}\|\phi_1\|_{L^\infty} \leq \|J_1\|_{L^1}\frac{\kappa_1}{\kappa_0-\kappa_*}\|\phi_1\|_{L^\infty} 
\]
by arguing similarly as for~\eqref{eq:mutselec-WWbis}.
\end{proof}

In order to prove Theorem \ref{theo:mutselec-Main} part (2) with constructive constants we use a Doblin-Harris type argument

\begin{lem}[Lyapunov Condition]\label{lem:mutselec-Lyap}
Under the above assumptions, for any $T> 0$, there are $\gamma_L\in (0,1)$ and $K>0$ such that
 $$
\| \widetilde S_T f \|_{L^1} \le \gamma_L \| f \|_{L^1} + K \| f \|_{\phi_1}.
$$
\end{lem}

\begin{proof}[Proof of Lemma~\ref{lem:mutselec-Lyap}]
Writing $f_t=\widetilde S_t f = e^{-\lambda_1 t}S_\LL(t)f$, we have, since $\lambda_1\geq0$,
\begin{align*}
   \frac{d}{dt}\int_{\R^d}|f_t|  &\leq \|J\|_1\int_{\R^d}|f_t|-\int_{\R^d}W|f_t|\\
   &\leq \int_{B^c_{R}}(\|J\|_1-W)|f_t|+\frac{\|J\|_1}{\alpha_R}\int_{ B_{R}}|f_t|\phi_1,
\end{align*}
for any $R>0$ and some  $\alpha_R$ the bound by below of $\phi_1$ in $B_R$.
Choosing $R$ large enough so that $W(x)\geq\|J\|_1+1$ for $|x|\geq R$, we get
\[\frac{d}{dt}\int_{\R^d}|f_t| \leq - \int_{\R^d}|f_t| + \frac{\|J\|_1+1}{\alpha_R}\int_{\R^d}|f_t|\phi_1.\]
Since 
$$\int_{ \R^d}|f_t|\phi_1\leq\int_{ \R^d}\widetilde S_t|f_0|\phi_1=\int_{ \R^d}|f_0|\phi_1,$$
 we infer
$$\|\widetilde S_t f\|\leq e^{-t}\|f\|+\frac{\|J\|_1+1}{\alpha_R}(1-e^{-t})\|f\|_{\phi_1},$$
by Gr{\"o}nwall's lemma. \end{proof}
 
 \begin{lem}[Harris's condition]\label{lem:SelMut-Harris}
 Under the assumption above, there exist $\psi_0 \in X'_{++}$, $ g_0 \in X_{+}$ and $ T > 0$ such that 
\beqn\label{eq:muta_inegHarris}
S_T f \ge g_0 \langle f,\psi_{0} \rangle, \quad \forall \, f \in X_+.
\eeqn
 \end{lem}

\begin{proof}[Proof of Lemma~\ref{lem:SelMut-Harris}]
{\sl Step 1. proof of \eqref{eq:muta_inegHarris}}. 
From Duhamel's formula \eqref{eq:itratedDuhamel} we have
\begin{equation*}
S_\LL = S_\BB + \dots+ (S_\BB \AA)^{*(N-1)} * S_\BB + (S_\BB \AA)^{(*N)} * S_\LL. 
\end{equation*}

We note that
$$  (S_\BB \AA * S_\BB)f(x)=\intot S_{\BB}(t-s)\AA S_\BB(s) fds=\intot [\AA (fe^{W(x)s})]e^{-W(x)(t-s)}ds.$$
For any $R>r$, $x\in B_R$, it is satisfied that
\begin{equation*}
    \AA (fe^{Ws})(x)=\int_{\R^d}J(x-y)f(y)e^{-W(y)s}dy\geq J_\ast e^{-W[2R] s}\int_{B_{r}(x)}f(y)dy
\end{equation*}
with $W[R]$ defined as in the proof of   \eqref{eq:muta_pos}. Then we get
$$(S_\BB \AA * S_\BB)f(x) \geq {\bf 1}_{B_R}(x)J_\ast te^{-W[2R] t}\int_{B_{r}(x)}f(y)dy.$$
 
Subsequently, we obtain that
$$S_\BB \AA * (S_\BB \AA * S_\BB)f(x)\geq {\bf 1}_{B_{R-r}}(x)\intot J_\ast se^{-W[2R] t} \AA \left(\textbf{1}_{B_R}(x)\int_{B_{r}(x)}f(y)dy\right)ds,$$
with
$$\AA \left(\textbf{1}_{B_R}(x)\int_{B_{r}(x)}f(y)dy\right)=\int_{\R^d}J(x-y)\textbf{1}_{B_R}(y)\int_{B_{r}(y)}f(z)dzdy\geq J_\ast\int_{B_{r}(x)}\int_{B_{r}(y)}f(z)dzdy.$$
We claim that for all $a\geq r$, 
$$\int_{B_{r}(x)}\int_{B_{a}(y)}f(z)dzdy\geq |B_{r/4}|\int_{B_{a+r/2}(x)}f(z)dz.$$
 Indeed, we deduce 
\begin{align*}
    \int_{B_{r}(x)}\int_{B_{a}(y)}f(z)dzdy=\int_{B_{r}(x)}\int_{\R^d}\textbf{1}_{B_{a}(y)}(z)f(z)dzdy=\int_{\R^d}f(z)\int_{B_{r}(x)}\textbf{1}_{B_{a}(z)}(y)dy\,dz
\end{align*}
and, since for all $z\in B_{a+r/2}(x)$, $$B_{\frac{r}{4}}\left(\frac{z-x}{|z-x|}\frac{3r}{4}+x\right)\subset B_r(x)\cap B_a(z),$$
we have
$$\int_{B_{r}(x)}\textbf{1}_{B_{a}(z)}(y)dy\geq  |B_{r/4}|\textbf{1}_{B_{a+r/2}(x)}(z),$$
and consequently, 
\begin{align*}
    \int_{B_{r}(x)}\int_{B_{a}(y)}f(z)dzdy\geq |B_{r/4}|\int_{B_{a+r/2}(x)}f(z)dz.
\end{align*}
We have obtained   
$$
S_\BB \AA * (S_\BB \AA * S_\BB)f(x)\geq {\bf 1}_{B_{R-r}(x)}
J_\ast^2 t^2/2e^{-W[2R] t} \int_{B_{r+r/2}(x)}f(y)dy.$$
Iterating the same argument,  we arrive to
$$(S_\BB \AA)^{(*n)}*S_\BB f(x)\geq {\bf 1}_{B_{R-nr}}(x)J_\ast ^n \frac{t^n}{n!}e^{-W[2R] t} \int_{B_{r+(n-1)r/2}(x)}f(y)dy.$$
In consequence, for $R=(n+1)r$, we get
$$(S_\BB \AA)^{(*n)}*S_\BB f(x)\geq {\bf 1}_{B_r}(x)J_\ast^n \frac{t^n}{n!}e^{-W[2(n+1)r] t} \int_{B_{(n-1)r/2}(0)}f(y)dy.$$

Coming back to the Duhamel formula \eqref{eq:itratedDuhamel}, we deduce
$$S_{\LL}f(x)\geq  \1_{B_r}(x)\sum_{n=2}^{\infty} \frac{(J_\ast  t)^n}{n!}e^{-W[2(n+1)r] t}\int_{B_{(n-1)r/2}}f(y)dy,$$
from where \eqref{eq:muta_inegHarris} follows with 
$$
 \psi_0 :=\sum_{n=2}^{\infty} \frac{(J_\ast  T)^n}{n!}e^{-W[2(n+1)r] T}\1_{B_{(n-1)r/2}}
$$
and $g_0 :=\textbf{1}_{B_r}$. 
\end{proof}

\begin{proof}[Proof of Theorem~\ref{theo:mutselec-Main} part (2)] 
Let us consider  $A >0$ and $f \in X_+$ such that $\| f \| \le A [f]_{\phi_1}$. For any integer $n\geq1$, we have
\bean
 [f]_{\phi_1}
&=& \int_{B_n}f\phi_1+\int_{B_n^c}f\phi_1
\leq \alpha_n\langle f,\psi_0\rangle+\beta_n\|f\|
\\
&\le& \alpha_n\langle f,\psi_0\rangle+\beta_nA [f]_{\phi_1}, 
\eean
with $\alpha_n=\|\phi_1\|_{L^\infty}/\inf_{B_n}\psi_0$ and $\beta_n=\|\phi_1\|_{L^\infty_{\langle W\rangle}}/\inf_{B_n^c}W$. 
Choosing $n_A$ such that $\beta_{n_A} A \le 1/2$, we deduce the constructive estimate 
$$
 [f]_{\phi_1} \le 2  \alpha_{n_A} \langle f,\psi_0\rangle,
 $$
 and thus that \eqref{eq:hyp-Harris} holds with $g_A := (2\alpha_{n_A})^{-1} g_0$. Because of the constructive lower bound~\eqref{eq:muta_pos} on $\phi_1$, 
 we have 
 $$
 \langle \phi_1,g_R \rangle \ge  (2\alpha_{n_A})^{-1}  \langle h_0,g_0 \rangle =: r_A,
 $$
 which provides   \eqref{eq:hyp-Harris-BdBis} in a quantified way. The two above estimates and the Lyapunov condition established in Lemma~\ref{lem:mutselec-Lyap} ensure that 
 we may apply the Doblin-Harris Theorem~\ref{theo:Harris} and thus conclude to \eqref{eq=mutselec-CvgceToInfty} with constructive rate. 
 \end{proof}
 
 \medskip
\subsection{A singular mutation kernel}
\label{ssec:appl5:singular}

Here we consider a mutation kernel supported by a set of zero Lebesgue measure, which thus does not satisfy~\eqref{eq:mutselec-hypJ}.
The kernel $J\in M^1_+(\R^d)$ is defined for any test function $\varphi\in C_0(\R^d)$ by
\[\langle J,\varphi \rangle = \eps^{-1}\sum_{i=1}^d \int_\R \varphi(0,\cdots\!,0,x_i,0,\cdots\!,0)J_i(\eps^{-1} x_i) dx_i,\]
where $(J_i)_{1\leq i\leq d}$ is a family of $L^1$ positive kernels on $\R$ and $\eps>0$ is a variance parameter.
The operator $\LL$ then reads
\[\LL f (x) = \eps^{-1}\sum_{i=1}^d\int_\R f(x-z\mathbf e_i)J_i(\eps^{-1} z)dz-W(x)f(x),\]
where $\mathbf e_i$ is the $i$-th unit vector of the canonical basis of $\R^d$.
This model was recently considered and studied by~\cite{Velleret2023} through a probabilistic approach.
It shares similarities with a model of telomere shortening which is under study in~\cite{DoumicOlayeTomasevic}.
We show that the method developed in the first sections of the present paper allows us to handle this model, under similar yet slightly different assumptions on the $J_i$ and $W$ than in~\cite{Velleret2023}.
In particular we consider more general fitness functions $W$ than quadratic ones.
More precisely, we assume that $W$ is a continuous function that satisfies~\eqref{eq:mutselec-hypW} and
\begin{equation}\label{eq:mutselec-hypW-2}
\log W(x) = O(|x|^2)\qquad \text{as}\ |x|^2:=\sum_{i=1}^d x_i^2\to\infty.
\end{equation}
The kernels $J_i$ are supposed to be centered Gaussian distributions
\[J_i(z)=M_i G_{\sigma_i}(z) := \frac{M_i}{\sigma_i\sqrt{2\pi}}\,e^{-\frac{z^2}{2\sigma_i^2}},\]
for given masses $(M_i)_{1\leq i\leq d}\in(0,+\infty)^d$ and variances $(\sigma_i)_{1\leq i\leq d}\in(0,+\infty)^d$.
Similarly as in Section~\ref{ssec:appl5:regular}, we work in the Banach lattice $X=L^1(\R^d)$ and we may prove the  following result.

\begin{theo}\label{theo:mutselec-Main-singular}
Under the above assumptions, there exists a constructive $\eps_0>0$ small enough, 
 such that  for any $\eps \in (0,\eps_0)$ the following conclusions hold
\begin{enumerate}
    \item The first eigentriplet problem \eqref{eq:triplet1}-\eqref{eq:triplet2} admits a unique solution
$(\lambda_1, f_1, \phi_1) \in \R \times X_+ \times X^\prime_+$ with the normalization $\|\phi_1\|=\langle\phi_1,f_1\rangle=1$,
and this triplet additionally satisfies $\lambda_1 > 0$, $f_1 >0$ and $\phi_1>0$. 
\item Moreover, $\LL$ generates a semigroup $S_\LL$ on $X$ and for any $f_0\in X$, there holds
\begin{equation}\label{eq:mutselec-expergo}
    \|e^{-\lambda_1 t}S_\LL(t) f_0 -\langle \phi_1,f_0\rangle f_1\|_{L^1} \leq Ce^{-\alpha t}\|f_0-\langle \phi_1,f_0\rangle f_1\|_{L^1} ,
\end{equation}
for any $t \ge 0$ and for some constructive constants $C,\alpha>0$.
\end{enumerate}
{\Blue In other words, the conclusions \ref{S1},  \ref{S2},  \ref{S33} and  \ref{E31} hold with constructive constants in~$L^1$.}
\end{theo}

\medskip

\begin{rem}
The assumption of small variance $\eps$ in Theorem~\ref{theo:mutselec-Main-singular} replaces~\eqref{eq:MUTS-assum1}-\eqref{eq:MUTS-assum2} as a condition which guarantees the strict positivity of $\kappa_0$ in the condition \ref{H2}, and so the strict positivity of $\lambda_1$.
This property is fundamental for ensuring the existence of $f_1$ in $L^1$ and for the existence of a spectral gap. 
On the contrary, for large values of $\eps$, there cannot exist $f_1\in L^1$, as it is proved in Remark~\ref{rem:selection-mutation}-(5).
The reason is a concentration phenomenon which creates an atom at the origin for the principal eigenvector when the dispersion due to the mutations is too big.
This is already noticed in~\cite[Rk.~5.3.1]{Velleret2023}, and we refer to~\cite{Bonnefon2017,BurgerBomze,Cloez2022} for more details about the singularity of $f_1$ and the concentration phenomenon.
\end{rem}

For proving Theorem~\ref{theo:mutselec-Main-singular}, we first show that the conditions \ref{H1}, \ref{H2} and \ref{H3} are verified for the dual problem in $L^\infty=X'=(L^1)'$.
Then we check that the Doblin-Harris conditions are satisfied, thus ensuring the existence, uniqueness and exponential stability for the primal problem.

It is worth noticing that since the $J_i$ are symmetric, we have $\LL^*=\LL$ and the only difference between the primal and dual problems is the Banach lattice in which it is posed.

\smallskip
{\bf Condition \ref{H1} and \ref{H1'}.}
With the same proof as in Section~\ref{ssec:appl5:regular}, $\LL$ generates a positive semigroup $S$ in $L^1$ with $\omega(S)\leq\|J\|_1$ and satisfies Kato's inequality.
We deduce that \ref{H1} and \ref{H1'} are verified for both $\LL$ in $X$ and $\LL^*$ in $X'$ with
\[\kappa_1=\|J\|_1=\sum_{i=1}^d M_i.\]

\smallskip
{\bf Condition \ref{H2}.}
{\Cyan In view of condition \ref{H3}}, we aim at verifying \ref{H2} with $\kappa_0$ close enough to $\kappa_1$.
More precisely, we define $\rho\in(0,1]$ the ratio between the geometric and arithmetic means of the masses $M_i$, namely
\[\rho := \frac{\big(\prod_{i=1}^d M_i\big)^{1/d}}{\frac{1}{d}\sum_{i=1}^d M_i},\]
we set
\begin{equation*}\label{eq:mutselec-deftheta}
\zeta := \frac{d\prod_{i=1}^d M_i}{2\,\kappa_1^d} = \frac{d^{1-d}\rho^d}{2} \in(0,1/2],
\end{equation*}
and we prove that there exists $\eps_0$ such that if $\eps\in(0,\eps_0)$, then \ref{H2} is verified with
\[\kappa_0 = \theta\kappa_1 \qquad\text{with}\quad \theta := (1-\zeta^2)^{1/d}\in(0,1).\]

Let us fix $\eta>0$ small enough so that
\[1+(\eta\sigma_i)^2\leq \Big(\frac{2}{1+\theta}\Big)^2\]
for all $i\in\{1,\cdots,d\}$.
We then define
\[f_0(x)=\prod_{j=1}^d G_{\eps/\eta}(x_j),\]
and we compute
\begin{align*}
\frac{\LL f_0(x)}{f_0(x)} & = \sum_{i=1}^d M_i\frac{G_{\eps/\eta}*G_{\eps\sigma_i}(x_i)}{G_{\eps/\eta}(x_i)} - W(x) = \sum_{i=1}^d M_i\frac{G_{\eps\sqrt{\eta^{-2}+\sigma_i^2}}(x_i)}{G_{\eps/\eta}(x_i)} - W(x) \\
& = \sum_{i=1}^d \frac{M_i}{\sqrt{1+(\eta\sigma_i)^2}} \exp\Big(\frac{\eta^2(\eta\sigma_i)^2}{1+(\eta\sigma_i)^2}\frac{x_i^2}{2\eps^2}\Big) - W(x) \\
& \geq \frac{1+\theta}{2}\sum_{i=1}^d {M_i} \exp\Big(\frac{\eta^2(\eta\sigma_i)^2}{1+(\eta\sigma_i)^2}\frac{x_i^2}{2\eps^2}\Big) - W(x).
\end{align*}
Due to Assumptions~\eqref{eq:mutselec-hypW} and~\eqref{eq:mutselec-hypW-2} on $W$ and using Jensen's inequality, we have
\[W(x)\leq \frac{1-\theta}{2}\Big(\min_{1\leq i\leq d}M_i\Big) d\, e^{C|x|^2/d}\leq \frac{1-\theta}{2} \sum_{i=1}^d M_i e^{Cx_i^2}\]
for some $C>0$ large enough.
Choosing $\eps_0>0$ small enough so that
\[2\eps_0^2\leq \frac{\eta^2(\eta\sigma_i)^2}{(1+(\eta\sigma_i)^2)C}\]
for all $i\in\{1,\cdots,d\}$, we obtain that
\[\frac{\LL f_0(x)}{f_0(x)} \geq \theta \sum_{i=1}^d M_i e^{Cx_i^2} \geq \theta\kappa_1 = \kappa_0\]
for any $\eps\in(0,\eps_0]$.
By virtue of Lemma \ref{lem:Existe1-Spectre2bis}-(ii), this proves the announced result.

\smallskip
{\bf Condition \ref{H3} in $X'=L^\infty$.}
We use the splitting $\LL=\AA+\BB$ with $\BB \phi = -W\phi$, and we aim at proving that~\eqref{as:weakly-compact-1} holds with $N=d$ in order to apply Lemma~\ref{lem:H3Y'}.
More precisely, we want to find $\varphi\in L^1$ and $\gamma\in(0,1)$ such that for any $\alpha\geq\kappa_0$, there holds
\begin{equation}\label{eq:mutselec-H3Linf}
\big\|(\RR_\BB(\alpha)\AA)^d\phi\big\|_{L^\infty} \leq \gamma \|\phi\|_{L^\infty} + \int_{\R^d}\phi(x)\varphi(x)\,dx
\end{equation}
for all $\phi\in L^\infty_+$.
We have
\[
\RR_\BB(\alpha)\phi = \frac{\phi}{\alpha+W} \le \frac{\phi }{ \kappa_0}
\]
and, defining
\[\AA_r\phi(x) := d\,\eps^{-d}\int_{\R^d}\phi(x-y)J^\otimes(y/\eps)\,dy\qquad\text{with}\quad J^\otimes(y):=\prod_{i=1}^d J_i(y_i),\]
we have
\[ \AA^d = \AA_r + \AA_s\]
with both $\AA_r$ and $\AA_s$ positive operators.
Positivity ensures that
\[\|\AA_r\|_{\BBB(L^\infty)} = \AA_r\1 = d\prod_{i=1}^d M_i
\]
and 
\[ \|\AA_s\|_{\BBB(L^\infty)} = \AA_s\1 = \|\AA^d\| - \|\AA_r\| \leq \|J\|_1^d - d\prod_{i=1}^d M_i.\]
We deduce that for any $\alpha\geq\kappa_0$,
\begin{align*}
  \big\|(\RR_\BB(\alpha)\AA)^d\phi\big\|_{L^\infty} & \leq \kappa_0^{-d} \big\|\AA_s\phi\|_{L^\infty} + \kappa_0^{1-d} \big\|\RR_\BB(\alpha)\AA_r\phi\|_{L^\infty}\\
& \leq \frac{\kappa_1^d - d\prod_{i=1}^d M_i}{\kappa_0^d}\|\phi\|_{L^\infty} + \kappa_0^{1-d} \Big\|\frac{\AA_r\phi}{\kappa_0+W}\Big\|_{L^\infty}.
\end{align*}
For any $R>0$ we have
\begin{align*}
\frac{\AA_r\phi(x)}{\kappa_0+W(x)} & = \frac{d\,\eps^{-d}\1_{B_R}(x)}{\kappa_0+W(x)} \int_{B_R}\! \phi(x-y) J^\otimes(y/\eps)\,dy \\
& \quad + \frac{d\,\eps^{-d}\1_{B_R}(x)}{\kappa_0+W(x)} \int_{B_R^c} \!\phi(x-y) J^\otimes(y/\eps)\,dy
+ \frac{d\,\eps^{-d}\1_{B_R^c(x)}}{\kappa_0+W(x)} \int_{\R^d} \!\phi(x-y) J^\otimes(y/\eps)\,dy \\
& \leq \frac{d\,\eps^{-d}}{\kappa_0} \prod_{i=1}^d \frac{M_i}{\sigma_i\sqrt{2\pi}} \int_{B_{2R}} \!\phi(y)\,dy
+ \frac{d}{\kappa_0} \int_{B_{R/\eps_0}^c} \!\!\!\!\!\!J^\otimes(y)dy \, \|\phi\|_{L^\infty}
+ \frac{d\prod_{i=1}^d M_i}{\kappa_0+W_R} \, \|\phi\|_{L^\infty} \\
& \leq \kappa_0^{d-1}\int_{\R^d}\phi(x)\varphi_R(y)\,dy + \frac{\eta_R}{\kappa_0} \|\phi\|_{L^\infty}, 
\end{align*}
where
\[\varphi_R = \frac{d\prod_{i=1}^dM_i/\sigma_i}{\sqrt{2\pi}(\eps\kappa_0)^d} \1_{B_{2R}} \quad\text{and}\quad \eta_R = d \int_{B_{R/\eps_0}^c} \!\!\!\!\!\!J^\otimes(y)dy + \frac{d\kappa_0}{W_R}\prod_{i=1}^d M_i,\]
and with $W_R=\inf_{B_R^c}W$.
We may therefore infer that
\[\big\|(\RR_\BB(\alpha)\AA)^d\phi\big\|_{L^\infty} \leq \frac{\kappa_1^d - d\prod_{i=1}^d M_i+\eta_R}{\kappa_0^d} \|\phi\|_{L^\infty} +\langle \phi, \varphi_R\rangle.\]
Since $W(x)\to+\infty$ and $J^\otimes(x)\to0$ as $|x|\to\infty$, we can find $R$ large enough so that
\[\eta_R \leq \frac{d}{2}\prod_{i=1}^d M_i = \zeta\kappa_1^d.
\]
Recalling that $\kappa_0^d=(1-\zeta)\kappa_1^d$, we then obtain~\eqref{eq:mutselec-H3Linf} with 
\[\gamma = \frac{\kappa_1^d - \frac{d}{2}\prod_{i=1}^d M_i}{\kappa_0^d} = \frac{1-\zeta}{1-\zeta^2} = \frac{1}{1+\zeta} < 1\qquad\text{and}\qquad \varphi=\varphi_R\in L^1.\]
Invoking Lemma~\ref{lem:H3Y'}, we deduce that \ref{H3} holds true for $\LL^*=\LL$ in $X'=L^\infty$.

\

From conditions~\ref{H1}-\ref{H2}-\ref{H3}, we infer the existence of a solution to the dual problem.

\begin{lem}\label{lem:mutselec-phi1}
If $\eps<\eps_0$, where $\eps_0$ is defined in the paragraph about Condition~\ref{H2} above, then there exist $\lambda_1\geq\kappa_0$ and $\phi_1\in X'_+$, $\|\phi_1\|_{L^\infty}=1$, such that $\LL^*\phi_1=\lambda_1\phi_1$.
Moreover, $\phi_1\in L^\infty_W$ and $\langle\phi_1,\varphi\rangle\geq1-\gamma$.
\end{lem}

\begin{proof}[Proof of Lemma~\ref{lem:mutselec-phi1}]
The existence of $(\lambda_1,\phi_1)$ follows from applying Theorem~\ref{theo:exist1-KRexistence}.
The equation $\LL\phi_1=\lambda_1\phi$ readily gives that
\[\|\phi_1\|_{L^\infty_W} \leq \|J*\phi_1\|_{L^\infty}+\lambda_1\|\phi_1\|_{L^\infty}\leq \|J\|_1+\lambda_1,\]
and the estimate $\langle\phi_1,\varphi\rangle\geq1-\gamma$ comes from Lemma~\ref{lem:H3Y'}
\end{proof}

We now aim at verifying~\eqref{eq:hyp-Harris}, \eqref{eq:stabilityKR-Lyapunov} and~\eqref{eq:hyp-Harris-BdBis} in order to apply Theorem~\ref{theo:Harris}.
 
\Black

\begin{lem}[Lyapunov Condition]\label{lem:mutselec-Lyap-2}
Under the above assumptions, for any $T > 0$, there are $\gamma_L\in (0,1)$  and $K>0$ such that
 $$
\| \widetilde S_T f \|_{L^1} \le \gamma_L \| f \|_{L^1} + K \| f \|_{\phi_1}.
$$
\end{lem}

\begin{proof}[Proof of Lemma~\ref{lem:mutselec-Lyap-2}]
The proof is exactly the same as for Lemma~\ref{lem:mutselec-Lyap} in Section~\ref{ssec:appl5:regular}.
\end{proof}

\begin{lem}[Harris's condition]\label{lem:SelMut-Harris-2}
   Under the above assumptions, there exists $\psi_0 \in X'_{++}$, $ g_0 \in X_{+}$ and $ T > 0$ such that 
\beqn\label{eq:muta_inegHarris-2}
S_T f \ge \langle f,\psi_{0} \rangle g_0, \quad \forall \, f \in X_+.
\eeqn
  \end{lem}
  
\begin{proof}[Proof of Lemma~\ref{lem:SelMut-Harris-2}]
We prove the dual version of~\eqref{eq:muta_inegHarris-2}, namely
\begin{equation}\label{eq:muta_inegHarris-2-dual}
S_T\phi \ge \langle \phi,g_0 \rangle \psi_0 , \quad \forall \, \phi\in X'_+,
\end{equation}
where we have used that $S^*_T=S_{\LL^*}(T)=S_\LL(T)=S_T$, since $\LL^*=\LL$ due to the symmetry of $J$.
The iterated Duhamel formula~\eqref{eq:itratedDuhamel} and the positivity of $\AA$ and $S_\BB$ ensure that
\[S_\LL\geq (S_\BB\AA)^{(*d)}*S_\BB.\]
We start by estimating $(S_\BB \AA*S_\BB)(t)\phi$
for $\phi\geq0$. Since
\begin{align*}
 \AA S_\BB(s) \phi (x) & \geq \eps^{-1} \int_\R \phi(x-z\mathbf e_1) e^{-sW(x-z\mathbf e_1)} J_1(z/\eps) \, dz \\
 & \geq \eps^{-1} \int_{x_1-1}^{x_1+1} \phi(x-z\mathbf e_1) e^{-sW(x-z\mathbf e_1)} J_1(z/\eps) \, dz \\
 & \geq \eps^{-1} e^{-sW[|x|+1]} J_1\Big(\frac{|x_1|+1}{\eps}\Big) \int_{x_1-1}^{x_1+1} \phi(x-z\mathbf e_1) \, dz ,
 \end{align*}
 where we recall the notation $W[R]=\sup_{B_R}W$, we get
 \[ (S_\BB \AA*S_\BB)(t)\phi (x) \geq \frac{t}{\eps} e^{-tW[|x|+1]} J_1\Big(\frac{|x_1|+1}{\eps}\Big) \int_{x_1-1}^{x_1+1} \phi(x-z\mathbf e_1) \, dz .\]
 Using now the part $J_2$ of $J$ we obtain
\begin{align*}
\AA(S_\BB & \AA*S_\BB)(s)\phi (x) \\ 
&\geq \frac{s}{\eps^2} J_1\Big(\frac{|x_1|+1}{\eps}\Big) \int_\R  e^{-sW[|x-z_2\mathbf e_2|+1]} \int_{x_1-1}^{x_1+1} \phi(x-z_1\mathbf e_1-z_2\mathbf e_1) \, dz_1 \, J_2(z_2/\eps) \, dz_2 \\
& \geq \frac{s}{\eps^2} J_1\Big(\frac{|x_1|+1}{\eps}\Big) J_2\Big(\frac{|x_2|+1}{\eps}\Big) e^{-sW[|x|+2]} \int_{x_2-1}^{x_2+1}  \int_{x_1-1}^{x_1+1} \phi(x-z_1\mathbf e_1-z_2\mathbf e_1) \, dz_1 dz_2
\end{align*}
and then
\begin{align*}
((S_\BB \AA & )^{(*2)} *S_\BB)(t)\phi (x) \\
& \geq \frac{t^2}{2\eps^2} e^{-tW[|x|+2]} J_1\Big(\frac{|x_1|+1}{\eps}\Big) J_2\Big(\frac{|x_2|+1}{\eps}\Big) \int_{x_2-1}^{x_2+1}  \int_{x_1-1}^{x_1+1} \phi(x-z_1\mathbf e_1-z_2\mathbf e_1) \, dz_1 dz_2 .
\end{align*}
Iterating and using the successive $J_i$'s parts of $J$ we finally get
\begin{align*}
S_\LL(t)\phi(x) & \geq ((S_\BB\AA)^{(*d)}*S_\BB)(t)\phi(x)\\
& \geq \frac{t^d}{d!}\, \eps^{-d} e^{-tW[|x|+d]} J^\otimes\Big(\frac{|x|+1}{\eps}\Big) \int_{[-1,1]^d} \phi(y)\,dy,
\end{align*}
which yields~\eqref{eq:muta_inegHarris-2-dual}, and so~\eqref{eq:muta_inegHarris-2}, with
\[\psi_0(x) = \frac{T^d}{d!}\, \eps^{-d} e^{-TW[|x|+d]} J^\otimes\Big(\frac{|x|+1}{\eps}\Big)\]
and $g_0=\1_{[-1,1]^d}$.
\end{proof}

\begin{cor}\label{cor:mutselec-f1-convexpo}
For $\eps \in (0,\eps_0)$,  there exists $f_1\in X_+$ such that $\LL f_1=\lambda_1f_1$ with $\langle f_1,\phi_1\rangle=1$.
Moreover, the exponential convergence~\eqref{eq:mutselec-expergo} holds for some constructive constants $C \ge 1$ and $\alpha>0$.
\end{cor}

\begin{proof}[Proof of Corollary~\ref{cor:mutselec-f1-convexpo}]
Similarly as in the proof of Theorem~\ref{theo:mutselec-Main} part (2), we can infer from Lemma~\ref{lem:SelMut-Harris-2} that~\eqref{eq:hyp-Harris} holds with $g_R=C_R g_0$ where $C_R>0$ is an explicit constant.
The Lyapunov condition~\eqref{eq:stabilityKR-Lyapunov} is established in Lemma~\ref{lem:mutselec-Lyap-2}, and the positivity condition~\eqref{eq:hyp-Harris-BdBis} readily follows from the estimate $\langle\phi_1,\varphi\rangle\geq1-\gamma$ established in Lemma~\ref{lem:mutselec-phi1}.
We can thus apply Theorem~\ref{theo:Harris} which
gives the conclusion.
\end{proof}

\begin{proof}[Proof of Theorem~\ref{theo:mutselec-Main-singular}]
It only remains to prove the uniqueness and strict positivity properties.
Combining~\eqref{eq:muta_inegHarris-2} and~\eqref{eq:muta_inegHarris-2} with $\phi=g_0$, we get that
\[S_{2T}f = S_T(S_T f) \geq \langle f,\psi_0 \rangle S_T g_0 \geq \Big(\int g_0^2\Big)\langle f,\psi_0 \rangle \psi_0 = 2^d\langle f,\psi_0 \rangle \psi_0.\]
for all $f\in X_+$.
Since $\psi_0>0$, this ensures that~\eqref{eq:Irred-Sirreducible} is verified, and then {\bf(H4)} because of point {\bf(4)} in
Lemma~\ref{lem:Irred-S>0impliesR>0}.
This gives the result of uniqueness and strict positivity by using Theorem~\ref{theo:KRgeometry1}.
\end{proof}

\Black
\

\paragraph{\bf Acknowledgments}
The authors are very grateful to Philippe Cieutat for the crucial discussions and references about the almost periodic functions, to Boyan Sirakov for having pointed out relevant references about the Krein-Rutman theorem for parabolic equations
 and to Otared Kavian for the enlightening discussions in many occasions.
The authors would also like to warmly thank Jochen Glück for his relevant and useful comments, suggestions, and corrections on the first version of the paper.
CFS  acknowledge funding by the European Union's Horizon 2020 research and innovation program under the Marie Skłodowska-Curie grant agreement No 754362.
PG has been supported by the ANR project NOLO (ANR-20-CE40-0015), funded by the French Ministry of Research.
This work was partially funded by the European Union (ERC, SINGER, 101054787). Views and opinions expressed are however those of the author(s) only and do not necessarily reflect those of the European Union or the European Research Council. Neither the European Union nor the granting authority can be held responsible for them.

\bigskip\bigskip


\begin{thebibliography}{100}
\bibitem{MR1245073}
L.~Arkeryd and C.~Cercignani.
\newblock A global existence theorem for the initial-boundary value problem for
  the {B}oltzmann equation when the boundaries are not isothermal.
\newblock {\em Arch. Rational Mech. Anal.}, 125(3):271--287, 1993.
\end{thebibliography}

 \end{document}